\documentclass{amsart}

\newtheorem{theorem}{Theorem}[section]
\newtheorem{prop}[theorem]{Proposition}
\newtheorem{lemma}[theorem]{Lemma}
\newtheorem{cor}[theorem]{Corollary}

\theoremstyle{definition}
\newtheorem{definition}[theorem]{Definition}
\newtheorem{example}[theorem]{Example}
\newtheorem{problem}[theorem]{Problem}

\theoremstyle{remark}
\newtheorem{remark}{Remark} 

\newcommand{\Q}{\mathbb{Q}}
\newcommand{\R}{\mathbb{R}}
\newcommand{\Z}{\mathbb{Z}}
\newcommand{\Id}{\mathrm{id}}
\newcommand{\Inv}{\mathrm{inv}}
\newcommand{\Rank}{\mathrm{rank}}
\newcommand{\Map}{\mathrm{map}}
\newcommand{\Conv}{\mathrm{conv}}
\newcommand{\Affi}{\mathrm{affi}}
\newcommand{\Cone}{\mathrm{cone}}
\newcommand{\Convcone}{\mathrm{convcone}}
\newcommand{\Vect}{\mathrm{vect}}
\newcommand{\QVect}{\mathbb{Q}\text{-}\mathrm{vect}}
\newcommand{\Clos}{\mathrm{clos}}
\newcommand{\Stab}{\mathrm{stab}}
\newcommand{\Supp}{\mathrm{supp}}
\newcommand{\Ht}{\mathrm{height}}
\newcommand{\Ord}{\mathrm{ord}}
\newcommand{\In}{\mathrm{in}}
\newcommand{\Ps}{\mathrm{ps}}
\newcommand{\Spec}{\mathrm{Spec}}
\newcommand{\Comp}{\mathrm{comp}}

\newcommand{\Den}{\mathrm{den}}
\newcommand{\Op}{^\mathrm{op}}
\newcommand{\Mx}{^\mathrm{max}}
\newcommand{\Fc}{^\mathrm{fc}}

\newcommand{\A}{\alpha}
\newcommand{\B}{\beta}

\begin{document}

\title[New Ideas for Resolution of Singularities]
{New Ideas for\\
Resolution of Singularities\\
in Arbitrary Characteristic}

\author{Tohsuke Urabe}
\address{Department of Mathematical Sciences, Ibaraki University, Mito, 
Ibaraki, 310-8512, Japan}
\email{urabe@mx.ibaraki.ac.jp}
\urladdr{http://mathmuse.sci.ibaraki.ac.jp/urabe/}


\date{October 31, 2010.}

\dedicatory{Dedicated to Heisuke Hironaka and Shreeram S. Abhyankar}


\begin{abstract}
Let $k$ be \emph{any} algebraically closed field in \emph{any} characteristic, let $R$ be any regular local ring such that $R$ contains $k$ as a subring, the residue field of $R$ is isomorphic to $k$ as $k$-algebras and $\dim R\geq 1$, let $P$ be any parameter system of $R$ and let $z\in P$. We consider any $\phi\in R$ with $\phi\neq 0$.

In our main theorem we assume several conditions depending on $P$, $z$ and Newton polyhedrons. By our assumptions the normal fan $\Sigma$ of the Newton polyhedron $\Gamma_+(P,\phi)$ of $\phi$ over $P$ has simple structure and we can make a special regular subdivision $\Sigma^*$ of $\Sigma$ called an \emph{upward subdivision}, starting from a regular cone with dimension equal to $\dim R$ and repeating star subdivisions with center in a regular cone of dimension two. Let $X$ and $\sigma:X\rightarrow\Spec(R)$ denote the toric variety over $\Spec(R)$ and the toric morphism associated with $\Sigma^*$. $X$ is irreducible and smooth and $\sigma$ is a composition of finite blowing-ups with center in a closed irreducible smooth subscheme of codimension two. We consider any closed point $a\in X$ such that $\sigma(a)$ is the unique closed point of $\Spec(R)$, the local ring $\mathcal{O}_{X,a}$ of $X$ at $a$ and the morphism $\sigma^*: R \rightarrow\mathcal{O}_{X,a}$ of local $k$-algebras induced by $\sigma$. We show that our numerical invariant of $\sigma^*(\phi)\in \mathcal{O}_{X,a}$ measuring the badness of the singularity is strictly less than the same invariant of $\phi\in R$ and the singularity $\phi$ is strictly improved by $\sigma$.

We notice that this result opens a way toward the theory of resolution of singularities in arbitrary characteristic. We add several submain theorems to make bridges toward it and to show that our assumptions of the main theorem are not strong. 

By these results we can show that in a mathematical game with two players A and B related to the resolution of singularities of $\phi$, the player A can always win the game after finite steps. It follows ``the local uniformization theorem in arbitrary characteristic and in arbitrary dimension''.
\end{abstract}

\maketitle

\section{Introduction}
\label{intro}
 For any local ring $S$ we denote the unique maximal ideal of $S$ by $M(S)$ and the set of invertible elements of $S$ by $S^\times$. We have $S=M(S)\cup S^\times, M(S)\cap S^\times=\emptyset$, $M(S)\in\Spec(S)$ and $M(S)$ is the unique closed point of the affine scheme $\Spec(S)$. For any noetherian local ring $S$, we denote the completion of $S$ by $S^c$. $S^c$ is a noetherian local ring containing $S$ as a local subring. Let $\Z_0$ and $\Z_+$ denote the set of non-negative integers and the set of positive integers respectively. 

Let $k$ be \emph{any} algebraically closed field in \emph{any} characteristic, and let $R$ be any regular local ring such that it contains $k$ as a subring, the residue field $R/M(R)$ of $R$ is isomorphic to $k$ as $k$-algebras, $R$ is a localization of a finitely generated $k$-algebra and $\dim R\geq 1$. 

For any $\psi\in R^c$ with $\psi\neq 0$ and any parameter system $Q$ of $R^c$ by $\Gamma_+(Q,\psi)$ we denote the Newton polyhedron of $\psi$ over $Q$.

Let $P$ be any parameter system of $R$ and let $z\in P$ be any element. $P$ is also a parameter system of $R^c$. By $R'$ we denote the localization of the $k$-subalgebra of $R$ generated by $P-\{z\}$ by its maximal ideal generated by $P-\{z\}$. The ring $R'$ is a local $k$-subalgebra of $R$, it is a regular local ring of dimension $\dim R-1$, and $P-\{z\}$ is a parameter system of $R'$.

Consider any $\phi\in R$ with $\phi\neq 0$.

If the Newton polyhedron $\Gamma_+(P,\phi)$ is non-degenerate, then we can construct explicitly an embedded resolution of the singularity $\phi$ corresponding to any regular subdivision of the normal fan of $\Gamma_+(P,\phi)$ using the toric theory (Khovanskii~\cite{K77}, Oka~\cite{O97}, Cox~\cite{C10}, Fulton~\cite{F93}, Kempf et al.~\cite{KKMS}.). We can give resolution of the singularity $\phi$ explicitly even if the characteristic of $k$ is positive under the assumption of non-degeneracy.

In this article we consider the case where the Newton polyhedron $\Gamma_+(P,\phi)$ is not necessarily non-degenerate using the toric theory and  $\Gamma_+(P,\phi)$.

Since $R$ is a UFD, we have an invertible element $u\in R^\times$, a finite set $\Omega$ of irreducible elements of $R$ and a mapping $a:\Omega\rightarrow\Z_+$ satisfying $\phi=u\prod_{\omega\in\Omega}\omega^{a(\omega)}$. We take any $u, \Omega$ and $a$ satisfying these conditions. Let $$\Xi=\{\omega\in\Omega|\partial\omega/\partial w\in M(R^c)\text{, and any }x\in P-\{z\}\text{ does not divide }\omega\}.$$

We say that an element $\psi\in R$ is a \emph{main factor} of the triplet $(P,z,\phi)$, or a $z$-\emph{main factor over} $P$ of $\phi$, if $\psi=v\prod_{\omega\in\Xi}\omega^{a(\omega)}$ for some $v\in R^\times$.
Since $R$ is a UFD, the condition that $\psi\in R$ is a main factor of $(P,z,\phi)$ does not depend of the choice of $u, \Omega$ and $a$ we used for the definition.

If both $\psi\in R$ and $\psi'\in R$ are main factors of $(P,z,\phi)$, then $\psi=u\psi'$ for some $u\in R^\times$ and $\Gamma_+(P,\psi)=\Gamma_+(P,\psi')$.

Any element $\omega\in R^c$ satisfying  $\omega=z^h+\sum_{i=0}^{h-1} \omega'(i)z^i$ for some $h\in\Z_0$ and some mapping $\omega':\{0,1,\ldots,h-1\}\rightarrow M(R^{\prime c})$ is called a $z$-\emph{Weierstrass polynomial} over $P$, and the non-negative integer $h$ is called the \emph{degree} of $\phi$.

Below, we defne and use concepts of \emph{$z$-Weierstrass type, $z$-simple and $z$-removable face} for $\Gamma_+(P,\phi)$ and numerical invariants $\Inv(P,z,\phi)$ and $\Inv 2(P,z,\phi)$. (In Section~\ref{concept} we define them again in more general situation.)

By definition, $\Gamma_+(P,\phi)$ is of $z$-Weierstrass type, if and only if, we can write uniquely $\phi=u\prod_{x\in P-\{z\}} x^{a(x)}\omega$ for some $u\in R^{c\times}$, some mapping $a: P-\{z\}\rightarrow\Z_0$ and some $z$-Weierstrass polynomial $\omega\in R^c$ over $P$. 

We consider the case where $\Gamma_+(P,\phi)$ is of $z$-Weierstrass type. Note that if the $z$-Weierstrass polynomial $\omega$ just above has degree $h$, then  $\Gamma_+(P,\phi)$ has a typical vertex corresponding to the monomial $(\prod_{x\in P-\{z\}} x^{a(x)})z^h$, which we call the $z$\emph{-top vertex} of $\Gamma_+(P,\phi)$.
We say that $\Gamma_+(P,\phi)$ is $z$\emph{-simple}, if $\Gamma_+(P,\phi)$ is of $z$-Weierstrass type and moreover, any compact face of $\Gamma_+(P,\phi)$ has dimension $1$ or $0$. If $\dim R\leq 2$, then $\Gamma_+(P,\phi)$ is necessarily $z$-simple.

Let $\psi\in R$ be a main factor of $(P,z,\phi)$. It is easy to see that $\Gamma_+(P,\psi)$ is also of $z$-Weierstrass type and the $z$-top vertex of  $\Gamma_+(P,\psi)$ corresponds to the monomial $z^k$ for some $k\in\Z_0$ with $k\leq h$. We define $\Inv(P,z,\phi)=k$, since $k\in\Z_0$ does not depend on the choice of a main factor $\psi$ of $(P,z,\phi)$. The non-negative integer $\Inv(P,z,\phi) $ is the main numerical invariant measuring the badness of the singularity $\phi$ in our theory. $\Inv(P,z,\phi) \neq 1$.  $\Inv(P,z,\phi) = 0$, if and only if, $\phi$ is a product of an invertible element of $R$, elements in $P-\{z\}$ admitting non-negative multiplicities, and some elements $\omega$ of $M(R)$ of order one with $\partial\omega/\partial z\in R^\times$ admitting positive multiplicities.

Consider the case where $\Gamma_+(P,\phi)$ is of $z$-Weierstrass type and any $x\in P-\{z\}$ does not divide $\phi$. In this case the $z$-top vertex $V_z$ of $\Gamma_+(P,\phi)$ corresponds to the monomial $z^h$ for some $h\in\Z_0$. Consider any face $F$ of $\Gamma_+(P,\phi)$. We say that $F$ is $z$\emph{-removable}, if $F$ contains the $z$-top vertex $V_z$, $F$ contains a vertex of $\Gamma_+(P,\phi)$ corresponding to a monomial in a form $\prod_{x\in P-\{z\}}x^{b(x)}z^k$ with $k<h$ and after some coordinate change sending $x$ to $x$ itself for any $x\in P-\{z\}$ and sending $z$ to $z+\chi$ for some $\chi\in M(R^{\prime c})$, $F$ becomes containing no vertex corresponding to a monomial in a form $\prod_{x\in P-\{z\}}x^{b(x)}z^k$ with $k<h$.

Instead of non-degeneracy, we assume that $\Gamma_+(P,\phi)$ is $z$-simple, and the Newton polyhedron over $P$ of the main factor of $\phi$ has no $z$-removable faces. (Our assumptions depend on the element $z$ of $P$. We can gradually understand that they are necessary and not strong.) Since $\Gamma_+(P,\phi)$ is $z$-simple, the normal fan  $\Sigma$ of $\Gamma_+(P,\phi)$ has simple structure. Note that the support $|\Sigma|$ of $\Sigma$ is a regular cone with dimension equal to $\dim R$. Starting from the fan $\mathcal{F}(|\Sigma|)$ consisting of $|\Sigma|$ and its faces and repeating star subdivisions with center in a regular cone of dimension two, we construct most effectively a special regular subdivision $\Sigma^*$ of $\Sigma$ with $|\Sigma^*|=|\Sigma|$, which we call an \emph{upward subdivision} of $\Sigma$. 
Let $X$ and $\sigma:X\rightarrow \Spec(R)$ denote the toric variety over $\Spec(R)$ and the toric morphism associated with $\Sigma^*$. Note that the scheme $X$ is separated, irreducible, smooth and of finite type over $\Spec(R)$, and the morphism $\sigma$ is a composition of finite blowing-ups with center in a closed irreducible smooth subscheme of codimension two. 

In our main theorem, Theorem~\ref{main} we show that the singularity $\phi$ is improved by the morphism $\sigma$. Theorem~\ref{main} has two conclusions. The first conclusion  treats the case of $\Inv(P,z,\phi)>0$. It claims that at any closed point $a\in X$ with $\sigma(a)=M(R)$ there exists a parameter system $\bar{P}$ of $\mathcal{O}_{X,a}$ of the local ring of $X$ at $a$ and an element $\bar{z}\in\bar{P}$ such that if we consider the injective homomorphism of local $k$-algebras $\sigma^*:R\rightarrow\mathcal{O}_{X,a}$ from $R$ to $\mathcal{O}_{X,a}$ induced by $\sigma$, then $\sigma^*(x)$ has normal crossings over $\bar{P}$ for any $x\in P$, the Newton polyhedron $\Gamma_+(\bar{P}, \sigma^*(\phi))$ is of $\bar{z}$-Weierstrass type, and $\Inv(\bar{P}, \bar{z}, \sigma^*(\phi))<\Inv(P,z,\phi)$. 

Since $\Inv(\hat{P}, \hat{z}, \sigma^*(\phi))<\Inv(P,z,\phi)$, we can claim that any hypersurface singularity can be improved by a composition of finite blowing-ups which we can describe explicitly associated with the Newton polyhedron of $\phi$, if we assume several conditions related to Newton polyhedrons.

We will explain the second conclusion treating the case $\Inv(P,z,\phi)=0$ later.

Now, notice that this result opens a way toward the theory of resolution of singularities in arbitrary characteristic, since it holds even if the characteristic of the ground field $k$ is positive, $\Inv(P,z,\phi)$ is a very explicit invariant measuring the badness of a singularity $\phi$, and it strictly decreases under a composition of finite blowing-ups. Contrary, the key concept in Hironaka's resolution theory in characteristic zero is the maximal contact and it is not effective in positive characteristic cases. (Hironaka~\cite{H75}, \cite{H74}, \cite{H64}(II, Chapter III, sections 7-10), Giraud~\cite{G75}, \cite{G74}, Hauser~\cite{Ha10}.) Besides, it is known that there exist metastatic singularities (Hironaka's terminology in Hironaka~\cite{H08}) or wild singularities (Hauser's terminology in Hauser~\cite{Ha10}) in the case of positive characteristic. For these singularities in positive characteristic case, under some blowing-ups, there appears a point in the inverse image of the singular point under consideration such that the invariant measuring the badness of a singularity \emph{increases}, when we try to follow the analogy of Hironaka's resolution theory in characteristic zero.
To overcome this phenomenon some programs for resolution of singularities in characteristic positive may become very complicated. 

Contrary to that the composition of blowing-ups is explicitely given associated with the Newton polyhedron of the singularity in our theory, it is very hard to give it explicitely in Hironaka's theory.

We add several submain theorems in order to make bridges to the resolution theory in arbitrary characteristic and to show that our assumptions of the main theorem are not strong.

We continue to consider the above situation concerning $a\in X$, $R$, $\sigma^*(\phi)$, $\bar{P}$ and $\bar{z}$.
$\Gamma_+(\bar{P}, \sigma^*(\phi))$ is of $\bar{z}$-Weierstrass type, and $\Inv(\bar{P}, \bar{z}, \sigma^*(\phi))<\Inv(P,z,\phi)$.
If the Newton polyhedron over $\bar{P}$ of the main factor of $(\bar{P}, \bar{z}, \sigma^*(\phi))$ has no $\bar{z}$-removable faces, the Newton polyhedron $\Gamma_+(\bar{P},\sigma^*(\phi))$ is $\bar{z}$-simple, and $\Inv(\bar{P}, \bar{z}, \sigma^*(\phi))>0$, then we can apply the first conclusion of our main theorem again after replacing the quadruplet $(R,P, z,\phi)$ by $(\mathcal{O}_{X,a}, \bar{P},\bar{z},\sigma^*(\phi))$ and we can make $\Inv(\bar{P}, \bar{z}, \sigma^*(\phi))$ further smaller. However, these assumptions are not necessarily satisfied. Our submain theorem Theorem~\ref{make simple} claims that we can make these assumptions satisfied after some blowing-ups. However, for Theorem~\ref{make simple} another induction assumption on $\dim R$ is necessary.

Note here that $R'$ is a local $k$-subalgebra of $R$ and $\dim R'=\dim R-1<\dim R $, and any $\phi'\in R'$ with $\phi'\neq 0$ has normal crossings over $P-\{z\}$ if $\dim R\leq 2$. Therefore, we decide that we use induction on $\dim R$, and we can assume the following claim $(*)$:

\renewcommand{\descriptionlabel}[1]%
	{\hspace{\labelsep}\textrm{#1}}
\begin{description}
\item[$(*)$]
For any  $\phi'\in R'$ with $\phi'\neq 0$, there exists a composition $\sigma':X'\rightarrow\Spec(R')$ of finite blowing-ups with center in a closed irreducible smooth subscheme such that the divisor on $X'$ defined by the pull-back $\sigma^{\prime *}(\phi')\in\mathcal{O}_{X'}(X')$ of $\phi'$ by $\sigma'$ has normal crossings.
\end{description}

Claim $(*)$ is true, if $\dim R\leq 2$.

Let $\sigma':X'\rightarrow\Spec(R')$ be any composition of finite blowing-ups with centers in closed irreducible smooth subschemes. The scheme $X'$ is smooth. We consider a morphism $\Spec(R)\rightarrow\Spec(R')$ induced by the inclusion ring homomorphism $R'\rightarrow R$, the product scheme $X=X'\times_{\Spec(R')}\Spec(R)$, the projection $\sigma:X\rightarrow\Spec(R)$, and the projection $\pi:X\rightarrow X'$. We know the following (Lemma~\ref{pull back blowing-ups}.):
\begin{enumerate}
\item The morphism $\sigma$ is a composition of finite blowing-ups with center in a closed irreducible smooth subscheme. The scheme $X$ is smooth.
\item Let $\Spec(R/zR)$ denote the prime divisor on $\Spec(R)$ defined by $z\in R$. The pull-back $\sigma^*\Spec(R/zR)$ of $\Spec(R/zR)$ by $\sigma$ is a smooth prime divisor of $X$, and $\sigma^*\Spec(R/zR)\supset \sigma^{-1}(M(R))$.
\item The projection $\pi:X\rightarrow X'$ induces an isomorphism $\sigma^*\Spec (R/zR)\rightarrow X'$.
\item
For any closed point $a\in X$ and any parameter system $Q'$ of the local ring $\mathcal{O}_{X',\pi(a)}$ of $X'$ at $\pi(a)$, $\sigma(a)=M(R)$ and $\{\sigma^*(z)\}\cup\pi^*(Q') $ is a parameter system of the local ring $\mathcal{O}_{X,a}$ of $X$ at $a$, where $\sigma^*:R\rightarrow \mathcal{O}_{X,a}$ denotes the homomorphism of local $k$-algebras induced by $\sigma$ and $\pi^*:\mathcal{O}_{X',\pi(a)}\rightarrow \mathcal{O}_{X,a}$ denotes the homomorphism of local $k$-algebras induced by $\pi$.
\end{enumerate}

We explain the conclusions of Theorem~\ref{make simple}. We have three conclusions.

In the first conclusion, we assume the above $(*)$ and that $\Gamma_+(P,\phi)$ is of $z$-Weierstrass type, and $\Inv(P,z,\phi)>0$.

Then, there exists a composition $\sigma':X'\rightarrow\Spec(R')$ of finite blowing-ups with center in a closed irreducible smooth subschemes with the following properties: We consider the product scheme $X=X'\times_{\Spec(R')}\Spec(R)$, the projection $\sigma:X\rightarrow\Spec(R)$ and the projection $\pi:X\rightarrow X'$. At any closed point $a\in X$ with $\sigma(a)=M(R)$, there exists a parameter system $\bar{Q}$ of $\mathcal{O}_{X',\pi(a)}$ with the following properties. Let $\bar{P}=\{\sigma^*(z)\}\cup\pi^*(\bar{Q})$: 

\begin{enumerate}
\item
$\sigma^{\prime *}(x)$ has normal crossings over $\bar{Q}$ for any $x\in P-\{z\}$.
\item
One of the following two conditions holds:
\begin{enumerate}
\item
$\Inv(\bar{P},\sigma^*(z),\sigma^*(\phi))=\Inv(P,z,\phi)$ and there exists an element \hfill\break$\bar{z}\in M(\mathcal{O}_{X,a})$ such that  $\partial\bar{z}/\partial\sigma^*(z)\in \mathcal{O}_{X,a}^\times$,
$\bar{P}_{\bar{z}}=\{\bar{z}\}\cup\pi^*(\bar{Q})$ is a parameter system of $\mathcal{O}_{X,a} $,
$\Gamma_+(\bar{P}_{\bar{z}},\sigma^*(\phi))$ is of $\bar{z}$-Weierstrass type,
$\Gamma_+(\bar{P}_{\bar{z}},\bar{\psi})$ has no $\bar{z}$-removable faces where $\bar{\psi}\in\mathcal{O}_{X,a}$ denotes a main factor of $(\bar{P}_{\bar{z}}, \bar{z}, \sigma^*(\phi))$, and $\Inv(\bar{P}_{\bar{z}},\bar{z},\sigma^*(\phi))=\Inv(P,z,\phi)$.
\item
$\Inv(\bar{P},\sigma^*(z),\sigma^*(\phi))<\Inv(P,z,\phi)$, and $\Gamma_+(\bar{P},\sigma^*(\phi))$ is of \hfill\break$\sigma^*(z)$-Weierstrass type.
\end{enumerate}
\end{enumerate}

By the above first conclusion we can make our assumptions stronger.

In the second conclusion, we assume the above $(*)$ and that $\Gamma_+(P,\phi)$ is of $z$-Weierstrass type, $\Gamma_+(P,\psi)$ has no $z$-removable faces where $\psi\in R$ denotes a main factor of $(P,z,\phi)$, and $\Inv(P,z,\phi)>0$.

Then, there exists a composition $\sigma':X'\rightarrow\Spec(R')$ of finite blowing-ups with center in a closed irreducible smooth subschemes with the following properties: We consider the product scheme $X=X'\times_{\Spec(R')}\Spec(R)$, the projection $\sigma:X\rightarrow\Spec(R)$ and the projection $\pi:X\rightarrow X'$.  At any closed point $a\in X$ with $\sigma(a)=M(R)$, there exists a parameter system $\bar{Q}$ of $\mathcal{O}_{X',\pi(a)}$ with the following properties. Let $\bar{P}=\{\sigma^*(z)\}\cup\pi^*(\bar{Q})$: 

\begin{enumerate}
\item
$\sigma^{\prime *}(x)$ has normal crossings over $\bar{Q}$ for any $x\in P-\{z\}$.
\item
$\Gamma_+(\bar{P},\sigma^*(\phi))$ is $\sigma^*(z)$-simple.
\item
One of the following two conditions holds:
\begin{enumerate}
\item 
$\Inv(\bar{P},\sigma^*(z),\sigma^*(\phi))=\Inv(P,z,\phi)$, and
$\Gamma_+(\bar{P},\bar{\psi})$ has no $\sigma^*(z)$-removable faces where $\bar{\psi}\in\mathcal{O}_{X,a}$ denotes a main factor of $(\bar{P},\sigma^*(z),\sigma^*(\phi))$.
\item
$\Inv(\bar{P},\sigma^*(z),\sigma^*(\phi))<\Inv(P,z,\phi)$
\end{enumerate}
\end{enumerate}

By the above second conclusion we can apply the first conclusion of our main theorem Theorem~\ref{main} again.

We can apply above three claims repeatedly and replace the ring $R$, the parameter system $P$, the element $z\in P$ and the non-zero element $\phi\in R$ repeatedly. We know that after finite sequences of compositions of blowing-ups and replacements of the quadruplet $(R, P,z,\phi)$ by $(\mathcal{O}_{X,a}, \bar{P},\bar{z},\sigma^*(\phi))$, any non-zero $\phi\in R$ such that $\Gamma_+(P,\phi)$ is of $z$-Weierstrass type is reduced to a non-zero $\bar{\phi}\in R$ satisfying the same condition and $\Inv(P,z,\bar{\phi})=0$.

Consider any element $\phi\in R$ such that $\phi\neq 0$, $\Gamma_+(P,\phi)$ is of $z$-Weierstrass type and $\Inv(P, z, \phi)=0$. 

By definition $\phi$ is a product of an invertible element of $R$, elements in $P-\{z\}$ admitting non-negative multiplicities, and some elements $\omega$ of $M(R)$ of order one with $\partial\omega/\partial z\in R^\times$ admitting positive multiplicities. Since $R$ is a UFD, there exist $u\in R^\times$, a mapping $a:P-\{z\}\rightarrow\Z_0$, a finite subset $\Omega$ of $M(R)$ and a mapping $b:\Omega\rightarrow\Z_+$ satisfying the following three conditions:
\begin{enumerate}
\item
$\phi=u\prod_{x\in P-\{z\}}x^{a(x)}\prod_{\omega\in\Omega}\omega^{b(\omega)}$.
\item For any $\omega\in\Omega$, $\omega$ is of order one and $\partial\omega/\partial z\in R^\times$.
\item If $\omega=v\omega'$ for some $\omega\in\Omega $, some $\omega'\in\Omega $ and some $v\in R^\times$, then $v=1$ and $\omega=\omega'$.
\end{enumerate}
The number of elements $\sharp\Omega\in\Z_0$ in $\Omega$ does not depend on the choice of $u, a, \Omega, b$ satisfying the above conditions. We define $\Inv 2(P,z,\phi)= \sharp\Omega\in\Z_0$.

If $\Inv 2(P,z,\phi)\leq 1$, then $\phi$ has normal crossings.

We consider the case where $\phi$ does not have normal crossings. $\Inv 2(P,z,\phi)\geq 2$.

It is easy to see that there exists $\bar{z}\in M(R)$ such that $\partial\bar{z}/\partial z\in R^\times$ and $\bar{z}$ divides $\phi$. If we take any element $\bar{z}\in M(R)$ satisfying these conditions and we replace the pair $(P,z)$ by $(\{\bar{z}\}\cup(P-\{z\}), \bar{z})$, then our element $\phi$ under consideration satisfies the same assumptions as above and furthermore, $z$ divides $\phi$.

We apply the third conclusion of our submain theorem Theorem~\ref{make simple}. 

In the third conclusion, we assume the above $(*)$ and that $\Gamma_+(P,\phi)$ is of $z$-Weierstrass type, $z$ divides $\phi$, $\Inv(P,z,\phi)=0$, and $\Inv 2(P,z,\phi)\geq 2$.

Then, there exists a composition $\sigma':X'\rightarrow\Spec(R')$ of finite blowing-ups with center in a closed irreducible smooth subschemes with the following properties: We consider the product scheme $X=X'\times_{\Spec(R')}\Spec(R)$, the projection $\sigma:X\rightarrow\Spec(R)$ and the projection $\pi:X\rightarrow X'$.  At any closed point $a\in X$ with $\sigma(a)=M(R)$, there exists a parameter system $\bar{Q}$ of $\mathcal{O}_{X',\pi(a)}$ with the following properties. Let $\bar{P}=\{\sigma^*(z)\}\cup\pi^*(\bar{Q})$: 

\begin{enumerate}
\item
$\sigma^{\prime *}(x)$ has normal crossings over $\bar{Q}$ for any $x\in P-\{z\}$.
\item
$\Gamma_+(\bar{P},\sigma^*(\phi))$ is $\sigma^*(z)$-simple.
\item
$\sigma^*(z)$ divides $\sigma^*(\phi)$.
\item
$\Inv(\bar{P},\sigma^*(z),\sigma^*(\phi))=0$.
\item 
$\Inv 2(\bar{P},\sigma^*(z),\sigma^*(\phi))=\Inv 2(P,z,\phi)\geq 2$
\end{enumerate}

By the above third conclusion of our submain theorem Theorem~\ref{make simple} we can make our assumptions on $\phi$ further stronger. We can assume that $\Gamma_+(P,\phi)$ is $z$-simple  and $z$ divides $\phi$.

We apply the second conclusion of our main theorem Theorem~\ref{main}.

In the second conclusion, we assume that $\Gamma_+(P,\phi)$ is $z$-simple, $z$ divides $\phi$, $\Inv(P,z,\phi)=0$, and $\Inv 2(P,z,\phi)\geq 2$.

We consider the normal fan $\Sigma$ of $\Gamma_+(P,\phi)$ and an upward subdivision of $\Sigma^*$ of $\Sigma$. Let $X$ and $\sigma:X\rightarrow\Spec(R)$ denote the toric variety over $\Spec(R)$ and the toric morphism associated with $\Sigma^*$.

Then, at any closed point $a\in X$ with $\sigma(a)=M(R)$ there exist a parameter system $\bar{P}$ and an element $\bar{z}\in\bar{P}$ such that $\sigma^*(x)$ has normal crossings over $\bar{P}$, $\Gamma_+(\bar{P},\sigma^*(\phi))$ is of $\bar{z}$-Weiestrass type, $\Inv(\bar{P},\bar{z},\sigma^*(\phi))=0$ and $\Inv 2(\bar{P},\bar{z},\sigma^*(\phi))<\Inv 2(P,z,\phi)$.

Since $\Inv 2(\bar{P},\bar{z},\sigma^*(\phi))<\Inv 2(P,z,\phi)$, singularity $\phi$ is improved by the morphism $\sigma$.

By the above claims, after finite sequences of compositions of blowing-ups and replacements of a quadruplet, any non-zero $\phi\in R$ such that $\Gamma_+(P,\phi)$ is of $z$-Weierstrass type is reduced to a non-zero $\bar{\phi}\in R$ with normal crossings.

The assumption of $z$-Weierstrass type for $\Gamma_+(P,\phi)$ may be strong. Lemma~\ref{make Weierstrass type} shows that for any $\phi\in R$ with $\phi\neq 0$ there exists a parameter system $\bar{P}$ of $R$ and an element $\bar{z}\in \bar{P}$ such that $\Gamma_+(\bar{P}, \phi)$ is of $\bar{z}$-Weierstrass type.

By the above claims, after finite sequences of compositions of blowing-ups and replacements of a quadruplet, any non-zero $\phi\in R$ is reduced to a non-zero $\bar{\phi}\in R$ with normal crossings.

Consider a mathematical game with two players A and B.
At the start of the game a pair $(R,\phi)$ of any regular local ring $R$ with $\dim R\geq 1$ such that $R$ contains $k$ as a subring, the residue field $R/M(R)$ is isomorphic to $k$ as $k$-algebras and $R$ is a localization of a finitely generated $k$-algebra, and any non-zero element $\phi\in R$ is given. We play our game repeating the following step. Before the first step we put $(S,\psi)=(R,\phi)$: At the start of each step, player A chooses a composition $\sigma:X\rightarrow\Spec(S)$ of finite blowing-ups with center in a closed irreducible smooth subscheme. Then, player B chooses a closed point $a\in X$ with $\sigma(a)=M(S)$. We have a morphism $\sigma^*:S\rightarrow\mathcal{O}_{X,a}$ of local $k$-akgebras induced by $\sigma$. If the element $\sigma^*(\psi)\in \mathcal{O}_{X,a}$ has normal crossings, then the palyer A wins. Otherwise we proceed to the next step after replacing the pair $(S,\psi)$ by the pair $(\mathcal{O}_{X,a}, \sigma^*(\psi))$.

Note that the pair $(S,\psi)$ satisfies the same assumptions as $(R,\phi)$ throughout the game. A similar game can be found in Spivakovsky~\cite{S82}.

By our results outlined above we can conclude that player A can always win the game after finite steps for any $R$ and any non-zero element $\phi\in R$ even if the characteristic of the ground field $k$ is positive. (Corollary~\ref{resolution game}.) It follows from the valuation theory, ``the local uniformization theorem in arbitrary characteristic and in arbitrary dimension''. (Corollary~\ref{local uniformization}, Zariski~\cite{Z40}, Abhyankar~\cite{A56}, Zariski et al.~\cite{ZS}.)

The idea of watching the height of the $z$-top vertex of a Newton polyhedron of $z$-Weierstrass type can be found in Hironaka~\cite{H67} in low dimensional cases. However, he did not manipulate higher dimensional cases, because he did not apply the toric theory. See also Cossart et al.~\cite{CJS}.

Some ideas in this article are inspired by the appendix of Abhyankar~\cite{A98} and Bogomolov~\cite{B96}.

We do \emph{not} claim that the centers of blowing-ups are contained in the singular locus of the subscheme to be resolved. It may be possible to improve our theorems and to add stataments claiming that any smooth point of the hypersurface to be resolved is not modified in our process of resolution.

The author expresses thanks to Herwig Hauser for valuable discussions with him through e-mail.

We give proofs only to difficult parts of our claims. Most of our claims follow from definitions. 

\begin{center}
\textsc{Table of contents}
\end{center}

\begin{description}
\item[˜~\ref{intro}] Introduction
\item[˜~\ref{concept}] Notations and basic concepts
\item[˜~\ref{scheme}] Basic scheme theory
\item[˜~\ref{mainmain}] Main results
\item[˜~\ref{btcs}] Basic theory of convex sets
\item[˜~\ref{cones}] Convex cones and convex polyhedral cones
\item[˜~\ref{simplex}] Simplicial cones and regular cones
\item[˜~\ref{decomposition}] Fans
\item[˜~\ref{cpp}] Convex pseudo polytopes
\item[˜~\ref{std}] Star subdivisions
\item[˜~\ref{istd}] Iterated star subdivisions
\item[˜~\ref{simple}] Simpleness and semisimpleness
\item[˜~\ref{basic}] Basic subdivisions
\item[˜~\ref{upper}] Upper boundaries and lower boundaries
\item[˜~\ref{compatible}] Height, characteristic functions and compatible mappings
\item[˜~\ref{height inequalities}] The height inequalities
\item[˜~\ref{upward}] Upward subdivisions and the hard height inequalities
\item[˜~\ref{toric theory}] Schemes associated with fans
\item[˜~\ref{main proof}] Proof of the main theorem
\item[˜~\ref{submain proofs}] Proof of the submain theorems
\end{description}

The most important is ``the hard height inequality" in Section~\ref{upward}. 
It depends heavily on ``the height inequality" in Section~\ref{height inequalities}.
In Sections \ref{btcs}-\ref{compatible} we develop exact theory of convex sets.

\section{Notations and basic concepts}
\label{concept}
We arrange notations and basic concepts related to Newton polyhedrons and commutative rings.

By $\Z$, $\Q$, $\R$ and $\mathbb{C}$ we denote the ring of integers, the field of rational numbers, the field of real numbers and the field of complex numbers respectively.

The following six notations are useful:
$$\Z_0=\{t\in\Z|t\geq 0\}, \qquad \Z_+=\{t\in\Z|t>0\},$$
$$\Q_0=\{t\in\Q|t\geq 0\}, \qquad \Q_+=\{t\in\Q|t>0\},$$
$$\R_0=\{t\in\R|t\geq 0\}, \qquad \R_+=\{t\in\R|t>0\}.$$

The set of all subset of a set $Z$ is denoted by $2^Z$. The identity mapping of a set $Z$ is denoted by $\Id_Z:Z\rightarrow Z$. The number of elements of a finite set $Z$ is denoted by $\sharp Z$.

The set of all mappings from a set $X$ to a set $Y$ is denoted by $$\Map(X,Y).$$ 
Let $X$ be any set and let $R$ be any ring. If $Y$ is an abelian group (respectively, an abelian semigroup, an $R$-module, a ring) , then the set $\Map(X,Y)$ has a natural structure of an abelian group (respectively, an abelian semigroup, an $R$-module, a ring). Let $Z$ be a set and let $Y$ be a subset of $Z$. Note that the inclusion mapping $Y\rightarrow Z$ induces a injective mapping $\Map(X,Y)\rightarrow\Map(X,Z)$. Using this injective mapping we regard $\Map(X,Y)$ as a subset of $\Map(X,Z)$. If $Z$ is an abelian group (respectively, an abelian semigroup,  an $R$-module, a ring) and $Y$ is a subgroup of $Z$ (respectively, a subsemigroup, an  $R$-submodule, a subring), then  $\Map(X,Y)$ is a subgroup (respectively, a subsemigroup, an  $R$-submodule, a subring) of $\Map(X,Z)$. 

Consider any abelian group $Z$. 

Let $J$ be any finite set and let $X:J\rightarrow 2^Z$ be any mapping.
We define
\begin{equation*}
\begin{split}
\sum_{j\in J}X(j)=\{z\in Z|&z=\sum_{j\in J}x(j)\text{ for some mapping }x:J\rightarrow Z\text{ satisfying}\\
&\quad x(j)\in X(j)\text{ for any }j\in J\}\in 2^Z.
\end{split}
\end{equation*}
Note that $\sum_{j\in J}X(j)$ is a subset of $Z$, $\sum_{j\in J}X(j)=\{0\}$ if $J=\emptyset$, and $\sum_{j\in J}X(j)=\emptyset$, if and only if, $J\neq\emptyset$ and $X(j)=\emptyset$ for some $j\in J$.

We call $\sum_{j\in J}X(j)$ the \emph{sum} or the \emph{Minkowski sum} of subsets $X(j), j\in J$. 

For any $r\in\Z_+$ and for any mapping $X:\{1,2,\ldots\,r\}\rightarrow 2^Z$ we also write
$$X(1)+X(2)+\cdots+X(r)= \sum_{j\in \{1,2,\ldots\,r\}}X(j).$$

We denote $-X=\{z\in Z|z=-x$ for some $x\in X\}$ for any subset $X$ of $Z$.

Let $X$ be any set, and let $Y$ be any subset of $Z$ with $0\in Y$. For any $a\in\Map(X,Y)$ we denote
$$\Supp(a)=\{x\in X|a(x)\neq 0\},$$
and we call $\Supp(a)$ the \emph{support} of $a$. It is a subset of $X$. We denote
$$\Map'(X,Y)=\{a\in\Map(X,Y)| \Supp(a)\text{ is a finite set.}\}.$$
$\Map'(X,Y)\subset\Map(X,Y)$. If $X$ is a finite set, we have $\Map'(X,Y)=\Map(X,Y)$.

Let $V$ be any vector space of finite dimension over $\R$, and let $X$ be any subset of $V$. 

The subset $X$ is called \emph{convex}, if $X\neq\emptyset$ and for any two different points $x,y$ of $X$, the \emph{segment} $\{a\in V |a=(1-t)x+ty \text{ for some }t\in\R\text{ with }0\leq t\leq 1\}$ joining $x$ and $y$ is contained in $X$. It is called an \emph{affine space}, if $X\neq\emptyset$ and for any two different points $x,y$ of $X$, the \emph{line} $\{a\in V |a=(1-t)x+ty \in\R\text{ for some }t\in\R\}$ joining $x$ and $y$ is contained in $X$. It is called a \emph{cone}, if $0\in X$ and for any $x\in X$ and any $t\in\R_0$, we have $tx\in X$. It is called a \emph{convex cone}, if $0\in X$ and for any $x,y\in X$ and any $t,u\in\R_0$, we have $tx+uy\in X$. It is called a \emph{vector space over} $\R$, or simply a \emph{vector space}, if $0\in X$ and for any $x,y\in X$ and any $t,u\in\R$, we have $tx+uy\in X$. It is called a \emph{vector space over} $\Q$, if $0\in X$ and for any $x,y\in X$ and any $t,u\in\Q$, we have $tx+uy\in X$. It is called \emph{closed}, if it is a closed subset with respect to the natural Hausdorff topology of $V$. 

In case $X\neq\emptyset$ the minimum convex subset (respectively, minimum affine space) containing $X$ with respect to the inclusion relation is denoted by $\Conv(X)$ (respectively, $\Affi(X)$). We define $\Conv(\emptyset)=\Affi(\emptyset)=\emptyset$. 
The minimum cone (respectively, minimum convex cone, minimum vector space over $\R$, minimum vector space over $\Q$, minimum closed subset) containing $X$ with respect to the inclusion relation is denoted by
$\Cone(X)$ (respectively, $\Convcone(X), \Vect(X), \QVect(X), \Clos(X)$). 

The subset $X$ is called a \emph{convex polytope},(respectively, \emph{convex polyhedral cone}), if there exists a \emph{finite} subset $Y$ of $V$ satisfying $X=\Conv(Y)$ and $Y\neq\emptyset$ (respectively, $X=\Convcone(Y)$). The subset $X$ is called a \emph{convex pseudo polytope}, if there exist \emph{finite} subsets $Y, Z$ of $V$ satisfying $X=\Conv(Y)+\Convcone(Z)$ and $Y\neq\emptyset$. The subset $X$ is called a \emph{simplicial cone}, if $X=\Convcone(C)$ for some $\R$-basis $B$ of $V$ and a subset $C$ of $B$. The subset $X$ is called a \emph{lattice}, if there exists a $\R$-basis $B$ of $V$ such that $X=\{a\in V|a=\sum_{b\in B}\lambda(b)b\text{ for some }\lambda\in\Map(B, \Z)\}$.
Any lattice $N$ of $V$ is a free $\Z$-submodule of $V$ with $\Rank N=\dim V$. 

For any $t\in\R$ we write $$tX=\{a\in V|a=tx\text{ for some }x\in X\}.$$ We know $(-1)X=-X$, and $0X=\{0\}$ if $X\neq\emptyset$. We write
$$\Stab(X)=\{a\in V|X+\{a\}\subset X\},$$
and call it the \emph{stabilizer} of $X$ in $V$. The stabilizer of $X$ in $V$ is a subsemigroup of $V$ containing $0$. 

Let $N$ be any lattice in $V$. The subset $X$ is called a \emph{regular cone} over N, if $X=\Convcone(C)$ for some $\Z$-basis $B$ of $N$ and a subset $C$ of $B$.

Any regular cone is a simplicial cone. Any simplicial cone is a convex polyhedral cone. Any convex polyhedral cone is a convex pseudo polytope. Any convex polytope is a convex pseudo polytope. Affine spaces, vector spaces, convex polytopes, convex polyhedral cones, and convex pseudo polytopes are non-empty closed convex subsets of $V$. If $X$ is convex (respectively, a cone, a convex cone), then $\Clos(X)$ is again convex (respectively, a cone, a convex cone). 

For any subset $T$ of $\R$ and for any $a\in V$ we denote $$Ta=\{b\in V|b=ta\text{ for some }t\in T\},$$ and it is a subset of $V$.

The \emph{dual vector space} $V^*=\mathrm{Hom}_\R(V,\R)$ is a vector space over $\R$ with $\dim V^*=\dim V$. The \emph{canonical bilinear form}
$$\langle\quad,\quad\rangle:V^*\times V\rightarrow \R,$$
is defined by putting $\langle\omega, a\rangle=\omega(a)\in\R$ for any $\omega\in\mathrm{Hom}_\R(V,\R)=V^*$ and any $a\in V$. The dual vector space $V^{**}$ of $V^*$ is identified with $V$ by the natural isomorphism $V\rightarrow V^{**}$ of vector spaces over $\R$.

We consider any vector space $W$ of finite dimension over $\R$ and any homomorphism $\pi:V\rightarrow W$ of vector spaces over $\R$.
Putting
$$\pi^*(\A)=\A \pi\in\mathrm{Hom}_\R(V,\R)=V^*,$$
for any $\A\in\mathrm{Hom}_\R(W,\R)=W^*$,
we define a mapping $\pi^*:W^*\rightarrow V^*$, and we call $\pi^*$ the \emph{dual homomorphism} of $\pi$.
The dual homomorphism $\pi^*$ is a homomorphism of vector spaces over $\R$.
For any $\omega\in W^*$ and for any $a\in V$ the equality $\langle\pi^*(\omega),a\rangle=\langle\omega,\pi(a)\rangle$ holds.
The dual homomorphism $\pi^{**}$ of $\pi^*$ is equal to $\pi$.

Let $N$ be any lattice in $V$. We denote 
$$N^*=\{\omega\in V^*|\langle\omega, a\rangle\in\Z \text{ for any } a\in N\},$$
and call $N^*$ the \emph{dual lattice} of $N$. Indeed, $N^*$ is a lattice in $V^*$. The dual lattice $N^{**}$ of $N^*$ is equal to $N$. 

Let $S$ be any convex cone in $V$. We denote
$$S^\vee|V=\{\omega\in V^*|\langle\omega, a\rangle\geq 0 \text{ for any } a\in S\},$$
and call $ S^\vee|V$ the \emph{dual cone} of $S$ over $V$. Indeed, $ S^\vee|V$ is a closed convex cone in $V^*$. The dual cone $ S^\vee|V^\vee|V^*$ of $ S^\vee|V$ is equal to the closure $\Clos(S)$ of $S$ in $V$. $S^\vee|V^\vee|V^*=S$, if and only if, $S$ is closed in $V$. When we need not refer to $V$, we also write simply $S^\vee$, instead of $ S^\vee|V$.

Let $P$ be any finite set. Note that $\Map(P,\R)$ is a vector space of finite dimension over $\R$ with $\dim \Map(P,\R)=\sharp P$,  $\Map(P,\Z)$ is a lattice in $\Map(P,\R)$,  $\Map(P,\R_0)$ is a regular cone over $\Map(P,\Z)$ in $\Map(P,\R)$ with $\Vect(\Map(P,\R_0))= \Map(P,\R)$, and $\Map(P,\Z_0)= \Map(P,\Z)\cap \Map(P,\R_0)$. Let $x\in P$. Let $y\in P$. Putting
$$f^P_x(y)=
\begin{cases}
1&\text{if $y=x$},\\
0&\text{if $y\neq x$},
\end{cases}$$
we define an element $f^P_x\in \Map(P,\Z_0)$. Note that the subset $\{f^P_x|x\in P\}$ of $\Map(P,\Z_0)$ is an $\R$-basis of $\Map(P,\R)$,  it is a $\Z$-basis of $\Map(P,\Z)$, and  $\Map(P,\R_0)$\hfill\break$=\Convcone(\{f^P_x|x\in P\})$. The dual basis of $\{f^P_x|x\in P\}$ is denoted by $\{f^{P\vee}_x|x\in P\}$. For any $x,y\in P$
$$\langle f^{P\vee}_x, f^P_y\rangle=
\begin{cases}
1&\text{if $x=y$},\\
0&\text{if $x\neq y$}.
\end{cases}$$
Indeed, $\{f^{P\vee}_x|x\in P\}$ is a $\R$-basis of the dual vector space $\Map(P,\R)^*$ of $\Map(P,\R)$, it is a $\Z$-basis of the dual lattice  $\Map(P,\Z)^*$ of $\Map(P,\Z)$, and $\Map(P,\R_0)^{\vee}=\Convcone(\{f^{P\vee}_x|x\in P\})$.

A commutative ring with the identity element is called simply a \emph{ring}. The identity element and the zero element of a ring are denoted $1$ and $0$ respectively.
We assume that any ring homomorphism $\lambda$ preserves the identity elements, in other words, the equality $\lambda(1)=1$ holds.

Let $R$ be any ring.
Let $X$ be any subset of $R$ and let $S$ be any subring of $R$. The minimum ideal of $R$ containing $X$ with respect the inclusion relation is denoted by $XR$ or $RX$. The minimum subring of $R$ containing $S$ and $X$ with respect to the inclusion relation is denoted by $S[X]$. In the case where $X$ contains only one element $x$, we also write simply $xR$, $Rx$, $S[x]$, instead of $\{x\}R$, $R\{x\}$, $S[\{x\}]$ respectively.
The set of all invertible elements in $R$ is denoted by $R^{\times}$. $R^{\times}\subset R$ and $R^{\times}$ is an abelian group with respect the multiplication.
Any ring with a unique maximal ideal is called a \emph{local ring}. A ring $R$ is local, if and only if, $R- R^{\times}$ is an ideal of $R$. The unique maximal ideal of a local ring $R$ is denoted by $M(R)$. We have $R=R^\times\cup M(R), R^\times\cap M(R)=\emptyset$ and $1\neq 0$ for any local ring $R$. A subset $R'$ of a local ring $R$ is called a \emph{local subring}, if $R'$ is a subring of $R$, $R'$ is a local ring and $M(R')=M(R)\cap R'$. The completion of a noetherian local ring $R$ is denoted by $R^c$. $R^c$ is a noetherian local ring containing $R$ as a local subring. $M(R^c)=M(R)R^c$. The inclusion ring homomorphism $R\rightarrow R^c$ induces an isomorphism $R/M(R)\rightarrow R^c/M(R^c)$ of residue fields. $\dim R^c=\dim R$. The smallest Henselian local subring of $R^c$ containing $R$ as a local subring is called the \emph{Henselization} of $R$. We denote it by $R^h$. $R^h$ is a noetherian local subring of $R^c$ containing $R$ as a local subring. $M(R^h)=M(R)R^h$. The inclusion ring homomorphism $R\rightarrow R^h$ induces an isomorphism $R/M(R)\rightarrow R^h/M(R^h)$ of residue fields. $\dim R^h=\dim R$. 

Let $R$ and $R'$ be local rings and $\lambda:R\rightarrow R'$ be a ring homomorphism. Always we have $\lambda(R^\times)\subset R^{\prime\times}$. We say that $\lambda$ is a \emph{local homomorphism}, if $\lambda(M(R))\subset M(R')$. $\lambda$ is a local homomorphism, if and only if, $\lambda^{-1}(M(R'))=M(R)$ , if and only if, $\lambda^{-1}(R^{\prime\times})=R^\times$.

A noetherian local ring $R$ is called \emph{complete}, if $R=R^c$. A local ring $R$ is called \emph{regular}, if $R$ is noetherian and there exists a finite subset $P$ of $M(R)$ such that $\sharp P=\dim R$ and $PR=M(R)$.

A subset $P$ of a regular local ring $R$ is called a \emph{parameter system} of $R$, if $P$ is finite with $\sharp P=\dim R$ and $PR=M(R)$. If $P$ is a parameter system of a regular local ring $R$, then the completion $R^c$ of $R$ and the Henselization $R^h$ are also regular local rings and $P$ is a parameter system of $R^c$ and $R^h$.

Let $R$ be any regular local ring, let $\phi\in R$ be any element and let $P$ be any parameter system of $R$. We say that $\phi$ has \emph{normal crossings} over $P$, if $\phi=u\prod_{x\in P}x^{\Lambda(x)}$ for some $\Lambda\in\Map(P,\Z_0)$ and some invertible element $u\in R^\times$. We say that $\phi$ has \emph{normal crossings}, if $\phi$ has normal crossings over $Q$ for some parameter system $Q$ of $R$.

Let $k$ be any field. Let $R$ be any regular local ring such that $\dim R\geq 1$, $R$ contains $k$ as a subring, and the residue field $R/M(R)$ is isomorphic to $k$ as algebras over $k$. Let $P$ be any parameter system of $R$. $P$ is algebraically independent over $k$. Let $Q$ be any subset of $P$. Let $R'$ be the localization of $k[Q]$ by the maximal ideal $k[Q]\cap M(R)=Qk[Q]$. $R'$ is a regular local subring of $R$ containing $k$ as a subring. The residue field $R'/M(R')$ of $R'$ is isomorphic to $k$ as algebras over $k$. $Q$ is a parameter system of $R'$. $\dim R'=\sharp Q$. $R^{\prime c}$ is a local subring of $R^c$, and $R^{\prime h}$ is a local subring of $R^h$. (Matsumura~\cite{M}, Milne~\cite{Mi}.)

Let $k$ be any field. Let $A$ be any complete regular local ring such that $\dim A\geq 1$, $A$ contains $k$ as a subring, and the residue field $A/M(A)$ is isomorphic to $k$ as algebras over $k$. Let $P$ be any parameter system of $A$.

Let $\phi$ be any element of $A$. Then, there exists a unique element\hfill\break$c\in \Map(\Map(P,\Z_0),k)$ with
$$\phi=\sum_{\Lambda\in \Map(P,\Z_0)} c(\Lambda)\prod_{x\in P}x^{\Lambda(x)}.$$
The infinite sum in the right-hand side is the limit with respect to the $M(A)$-adic topology on $A$.
We take the unique element $c\in \Map(\Map(P,\Z_0), k)$ satisfying the above  equality. The element $c$ depends on $\phi$ and $P$.
Consider any $\Lambda\in\Map(P,\Z_0)$. We call $\Lambda$ the \emph{index}, $\prod_{x\in P}x^{\Lambda(x)}\in A$ a \emph{monomial} over $P$,  $c(\Lambda)\in k$ a \emph{coefficient} of $\phi$, $c(\Lambda)\prod_{x\in P}x^{\Lambda(x)}\in A$ a \emph{term} of $\phi$, and $\sum_{x\in P}\Lambda(x)\in\Z_0$ the \emph{degree} of the index $\Lambda$, of the monomial $\prod_{x\in P}x^{\Lambda(x)}$, or of the term $ c(\Lambda)\prod_{x\in P}x^{\Lambda(x)}$. Note that $0\in\Map(P,\Z_0)$. We denote $c(0)$ by $\phi(0)$ and we call $\phi(0)\in k$, \emph{the constant term} of $\phi$. $\phi-\phi(0)\in M(A)$. $\phi(0)=0\Leftrightarrow \phi\in M(A)$.
We denote
$$\Supp(P, \phi)=\Supp(c)=\{\Lambda\in \Map(P,\Z_0)| c(\Lambda)\neq 0\},$$
and we call $\Supp(P,\phi)$ the \emph{support} of $\phi$ over $P$.
It is a subset of $ \Map(P,\Z_0)$.
Note that $\phi=0\Leftrightarrow c=0\Leftrightarrow \Supp(P,\phi)=\emptyset$. 

Let $F$ be any subset of $\Map(P,\R)$. We denote
$$\Ps(P,F,\phi)=
\begin{cases}
\sum_{\Lambda\in \Supp(P, \phi)\cap F} c(\Lambda)\prod_{x\in P}x^{\Lambda(x)}
&\text{if $ \Supp(P, \phi)\cap F \neq \emptyset$},\\
0&\text{if $ \Supp(P, \phi)\cap F =\emptyset$},
\end{cases}$$
and we call $\Ps(P,F,\phi)\in A$ the \emph{partial sum} of $\phi$ over $P$ with respect to $F$.

Below, we consider the case $\phi\neq 0$ for a while. 

We define
$$\Gamma_+(P, \phi)= \Conv(\Supp(P,\phi))+ \Map(P,\R_0),$$
and  call  $\Gamma_+(P, \phi)$ the \emph{Newton polyhedron} of $\phi$ over $P$. By definition we have $\Gamma_+(P, \phi)\subset \Map(P,\R_0)\subset \Map(P,\R)$. We can show that there exists a non-empty \emph{finite} subset $Y$ of $\Supp(P,\phi)$ with $\Gamma_+(P, \phi)=\Conv(Y)+ \Map(P,\R_0)$, and $\Gamma_+(P, \phi)$ is a convex pseudo polytope in $\Map(P,\R)$ with $\Stab(\Gamma_+(P, \phi))= \Map(P,\R_0)$.
(Lemma~\ref{Newton1}, Lemma~\ref{Newton2}.)

Let $\omega\in \Map(P, \R_0)^{\vee}$ be any element. We know that $\{\langle\omega, a\rangle|a\in \Supp(P, \phi)\}$\break $\subset\R_0$, and the minimum element $\min\{\langle\omega, a\rangle|a\in \Supp(P, \phi)\}$ of $\{\langle\omega, a\rangle|a\in \Supp(P, \phi)\}$ exists.
We define
\begin{equation*}\begin{split}
\Ord(P,\omega,\phi)&=\min\{\langle\omega, a\rangle|a\in \Supp(P, \phi)\}\in\R_0,\\
\Supp(P, \omega, \phi)&=\{a\in \Supp(P, \phi)| \langle\omega, a\rangle=\Ord(P,\omega,\phi)\}\subset \Supp(P, \phi),\\
\In(P,\omega,\phi)&=
\sum_{\Lambda\in \Supp(P, \omega, \phi)} c(\Lambda)\prod_{x\in P}x^{\Lambda(x)}\in A.\\
\end{split}\end{equation*}

We consider the case $\phi=0$. We introduce a symbol $\infty$ satisfying the following conditions: for any $t\in\R$, we have $\infty>t, \infty\geq t, \infty\neq t, t<\infty, t\leq\infty, t\neq\infty, \infty+t=t+\infty=\infty$, and moreover $\infty+\infty=\infty$. 
Let $\omega\in \Map(P, \R_0)^{\vee}$ be any element. We define
\begin{equation*}\begin{split}
\Ord(P,\omega,0)&= \infty,\\
\In(P,\omega,0)&=0\in A.\\
\end{split}\end{equation*}

Let $\omega\in \Map(P, \R_0)^{\vee}$ be any element. 
In the general case including the case of $\phi=0$, we have defined $\Ord(P,\omega,\phi)\in\R_0\cup\{\infty\}$ and $\In(P,\omega,\phi)\in A$.
We call $\Ord(P,\omega,\phi)\in\R_0\cup\{\infty\}$ the \emph{order} of $\phi$ over $P$ with respect to $\omega$. By definition $\Ord(P,\omega,\phi)=\infty$ if and only if $\phi=0$.
We call $\In(P,\omega,\phi)\in A$ the \emph{initial sum} of $\phi$ over $P$ with respect to $\omega$. By definition $\In(P,\omega,\phi)=0$ if and only if $\phi=0$.

We can show that the following holds for any $\omega\in\Map(P,\R_0)^\vee$, any $\phi\in A$, any $\psi\in A$ and any $\A\in k$ with $\A\neq 0$:
\begin{enumerate}
\item
$\Ord(P,\omega,-\phi)= \Ord(P,\omega,\phi)$. $\In(P,\omega,-\phi)= -\In(P,\omega,\phi)$. \hfill\break$\Ord(P,\omega,\A\phi)= \Ord(P,\omega,\phi)$. $\In(P,\omega,\A\phi)= \A\In(P,\omega,\phi)$.
\item
$\Ord(P,\omega, \phi+\psi)\geq\min\{\Ord(P,\omega,\phi), \Ord(P,\omega,\psi)\}$.
\item
$\Ord(P,\omega, \phi+\psi)=\min\{\Ord(P,\omega,\phi), \Ord(P,\omega,\psi)\}$, if and only if, \hfill\break$\Ord(P,\omega,\phi)\neq\Ord(P,\omega,\psi)$ or $\In(P,\omega,\phi)+ \In(P,\omega,\psi)\neq 0$.
\item $\Ord(P,\omega, \phi\psi)=\Ord(P,\omega,\phi)+\Ord(P,\omega,\psi)$.\hfill\break
$\In(P,\omega, \phi\psi)=\In(P,\omega,\phi)\In(P,\omega,\psi)$.
\end{enumerate}

Let $R$ be any regular local ring.
We have $\cap_{m\in\Z_0}M(R)^m=\{0\}$.

We define a mapping
$$\Ord:R\rightarrow \Z_0\cup\{\infty\}$$
by putting
$$\Ord(\phi)=
\begin{cases}
m&\text{if $m\in\Z_0$ and $\phi\in M(R)^m-M(R)^{m+1}$},\\
\infty&\text{if $\phi=0$},
\end{cases}$$
for any $\phi\in R$. For any $\phi\in R$ with $\phi\neq 0$ and any $m\in\Z_0$, we say that $\phi$ is \emph{of order} $m$, if $\Ord(\phi)=m$.

Assume moreover, that $R$ contains $k$ as a subring, the residue field $R/M(R)$ is isomorphic to $k$ as $k$-algebras and $\dim R=\dim A$. Under these assumptions the completion $R^c$ of $R$ is a ring containing $R$ as a local subring and $R^c$ and $A$ are isomorphic as $k$-algebras. We take any isomorphism $\lambda:R^c\rightarrow A$ of $k$-algebras.
We know that $\Ord(\phi)=\Ord(P,\sum_{x\in P}f^{P\vee}_x,\lambda(\phi))$ for any $\phi\in R$.
Consider any element $x\in R$ with $\lambda(x)\in P$, any $\phi\in R$ and any $m\in\Z_0$.
$\Ord(P,f^{P\vee}_{\lambda(x)},\lambda(\phi))\geq m$, if and only if, $\phi=x^m\psi$ for some $\psi\in R$. 

Consider any $\phi\in A$ with $\phi\neq 0$.

A subset $F$ of  $\Map(P,\R)$ is a \emph{face} of $\Gamma_+(P,\phi)$, if and only if, there exists $\omega \in \Map(P,\R_0)^\vee$ such that $F=\{a \in \Gamma_+(P,\phi)|\langle\omega,a\rangle=\Ord(P,\omega,\phi)\}$. Any face of $\Gamma_+(P,\phi)$ is a non-empty closed subset of $\Gamma_+(P,\phi)$, and is a convex pseudo polytope. For any face $F$ of $\Gamma_+(P,\phi)$, we put $\dim F=\dim\Affi(F)\in\Z_0$ and we call $\dim F$ the \emph{dimension} of $F$. Any face of $\Gamma_+(P,\phi)$ with dimension zero is called a \emph{vertex} of $\Gamma_+(P,\phi)$. (See Definition~\ref{faces of cpp}.)

By $\mathcal{V}(\Gamma_+(P, \phi))$ we denote the union of all \emph{vertices} of $\Gamma_+(P, \phi)$. By definition we have
\begin{equation*}\begin{split}
\mathcal{V}(\Gamma_+(P, \phi))&=\{a\in\Gamma_+(P, \phi)|
\text{ There exists }\omega\in\Map(P,\R_0)^\vee\text{ such that for any }\\
&\qquad\qquad\quad b\in\Gamma_+(P, \phi)\text{ with }\langle\omega,b\rangle=\langle\omega,a\rangle\text{, we have }b=a\}.\\
\end{split}\end{equation*}
We call $ \mathcal{V}(\Gamma_+(P, \phi))$ the \emph{skeleton} of $\Gamma_+(P, \phi)$. The set $\mathcal{V}(\Gamma_+(P, \phi))$ is a non-empty finite subset of $\Supp(P,\phi)$, and $\Gamma_+(P, \phi)= \Conv(\mathcal{V}(\Gamma_+(P, \phi)))+ \Map(P, \R_0)$.
We denote $c(\Gamma_+(P, \phi))=\sharp\mathcal{V}(\Gamma_+(P, \phi))\in\Z_+$,
and we call $c(\Gamma_+(P, \phi))$ the \emph{characteristic number} of $\Gamma_+(P, \phi)$.

We know that $\Gamma_+(P, \phi)$ has only one vertex$\Leftrightarrow c(\Gamma_+(P, \phi))=1\Leftrightarrow\phi$ has formal normal crossings over $P$ (Lemma~\ref{Newton2}.6.), and that these equivalent conditions always hold, if $\dim A=1$.

For any $\phi\in A$ with $\phi\neq 0$ and any $\psi\in A$ with $\psi\neq 0$, $\Gamma_+(P, \phi\psi)= \Gamma_+(P, \phi)+ \Gamma_+(P, \psi)$. (Lemma~\ref{Newton2}.8.) For any $\phi\in A$ with $\phi\neq 0$ and any $u\in A^\times$, $\Gamma_+(P, \phi)= \Gamma_+(P, u\phi)$.

Consider any $\phi\in A$ with  $\phi\neq 0$ and any $z\in P$. 

Note that for any $a\in \mathcal{V}(\Gamma_+(P, \phi))$, we have $\langle f^{P\vee}_z, a\rangle\in\Z_0$. We define
\begin{equation*}\begin{split}
&\Ht(z,\Gamma_+(P, \phi)) \\
= &\max\{\langle f^{P\vee}_z, a\rangle|a\in \mathcal{V}(\Gamma_+(P, \phi))\}
-\min\{\langle f^{P\vee}_z, a\rangle|a\in \mathcal{V}(\Gamma_+(P, \phi))\}\in\Z_0,
\end{split}\end{equation*}
and we call $\Ht(z,\Gamma_+(P, \phi)) $ the \emph{height} of $\Gamma_+(P, \phi)$ with respect to $z$, or simply $z$-\emph{height} of  $\Gamma_+(P, \phi)$. It is a non-negative integer. By definition, $\Ht(z,\Gamma_+(P, \phi)) =0$ if and only if the value $\langle f^{P\vee}_z, a\rangle$ does not depend on $a\in \mathcal{V}(\Gamma_+(P, \phi))$. $\Ht(z,$\hfill\break$\Gamma_+(P, \phi))=0$, if $\dim A=1$.

Let $a\in \mathcal{V}(\Gamma_+(P, \phi))$. We say that $\{a\}$ is a $z$-\emph{top vertex} of $\Gamma_+(P,\phi)$, if $\langle f^{P\vee}_z, a\rangle= \max\{\langle f^{P\vee}_z, b\rangle|b\in \mathcal{V}(\Gamma_+(P, \phi))\}$. We say that $\{a\}$ is a $z$-\emph{bottom vertex} of $\Gamma_+(P,\phi)$, if $\langle f^{P\vee}_z, a\rangle= \min\{\langle f^{P\vee}_z, b\rangle|b\in \mathcal{V}(\Gamma_+(P, \phi))\}$.

Let $A'$ denote the completion of the subring $k[P - \{z\}]$ of $A$ by the maximal ideal $k[P - \{z\}]\cap M(A)=( P - \{z\}) k[P - \{z\}]$. The ring $A'$ is a complete regular local subring of $A$ and $M(A')=M(A)\cap A' =(P -\{z\})A'$. The set $P-\{z\}$ is a parameter system of $A'$. The completion of the subring $ A'[z]$ of $A$ by the prime ideal $z A'[z]$ coincides with $A$.

Any element $\phi$ in $A$ such that $\phi=z^h+\sum_{i=0}^{h-1} \phi'(i) z^i$ for some $h\in\Z_0$ and some mapping $\phi':\{0, 1,\ldots, h-1\}\rightarrow M(A')$ is called a $z$-\emph{Weierstrass polynomial} over $P$, and the integer $h$ is called \emph{degree} of $\phi$. For any $z$-Weierstrass polynomial $\phi$ over $P$, we denote its degree by
$$\deg(P,z,\phi)\in\Z_0.$$
Note that $\deg(P,z,\phi)=0\Leftrightarrow \phi=1$.

We say that $\Gamma_+(P,\phi)$ is \emph{of $z$-Weierstrass type}, if there exists $a\in \Gamma_+(P, \phi)$ satisfying the equality $\langle f^{P\vee}_x, a\rangle=\Ord(P, f^{P\vee}_x,\phi)$ for any $x\in P - \{z\}$.

The following two conditions are equivalent (Lemma~\ref{Newton2}.10.):
\begin{enumerate}
\item The Newton polyhedron $\Gamma_+(P,\phi)$ is of $z$-Weierstrass type.
\item There exist uniquely an invertible element $u\in A^{\times}$, a mapping $a:P-\{z\}\rightarrow\Z_0$, and a $z$-Weierstrass polynomial $\psi\in A$ over $P$ satisfying
$$\phi=u \prod_{x\in P-\{z\}}x^{a(x)} \psi$$.
\end{enumerate}

Assume that $\Gamma_+(P,\phi)$ is of $z$-Weierstrass type, we know the following \hfill\break(Lemma~\ref{Newton2}.9.$(a)$):
\begin{enumerate}
\item
If $\psi\in A$, $\omega\in A$ and $\phi=\psi\omega$, then both $\Gamma_+(P,\psi)$ and $\Gamma_+(P,\omega)$ are of $z$-Weierstrass type.
\item
$\Ht(z,\Gamma_+(P, \phi))=0\Leftrightarrow \Gamma_+(P,\phi)$ has only one vertex $\Leftrightarrow c(\Gamma_+(P, \phi))=1\Leftrightarrow \phi$ has normal crossings over $P$.
\item 
The Newton polyhedron $\Gamma_+(P,\phi)$ has a unique $z$-top vertex.
\end{enumerate}

Below, by $\{a_1\}$ we denote the unique $z$-top vertex of $\Gamma_+(P,\phi)$. 
Let $b=\Ord(P,$\hfill\break$ f^{P\vee}_z,\phi)\in\Z_0$ and let $h=\Ht(z,\Gamma_+(P, \phi)) \in\Z_0$.

\begin{enumerate}
\setcounter{enumi}{3}
\item
Consider any $a\in \Gamma_+(P, \phi)$. The equality $\langle f^{P\vee}_x, a\rangle=\Ord(P, f^{P\vee}_x,\phi)$ holds for any $x\in P - \{z\}\Leftrightarrow a-a_1\in\R_0f^P_z$.
\item
$\langle f^{P\vee}_z, a_1\rangle=b+h$.
\item
There exist uniquely an invertible element $u\in A^{\times}$ and a mapping $\phi':\{0, 1,\ldots, h-1\}\rightarrow M(A')$ satisfying
$\phi=u (\prod_{x\in P -\{z\}}x^{\langle f^{P\vee}_x, a_1\rangle}) z^b (z^h+\sum_{i=0}^{h-1} \phi'(i) z^i)$,
and $\phi'(0)\neq 0$ if $h>0$
\end{enumerate}

The concept of $z$-removable faces is important.

Assume that $\Gamma_+(P,\phi)$ is of $z$-Weierstrass type and any $x\in P-\{z\}$ does not divide $\phi$. Under this assumption we can give the definition of $z$-removable faces. 

By assumption there exists uniquely $h\in\Z_0$ such that $\{hf^P_z\}$ is the unique $z$-top vertex of $\Gamma_+(P,\phi)$. We take the unique $h\in\Z_0$ satisfying this condition.

Let $F$ be any face of $\Gamma_+(P,\phi)$. We say that $F$ is $z$-\emph{removable}, if $hf^P_z \in F$, $F\not\subset\{hf^P_z\}+\Map(P,\R_0)$ and there exist an invertible element $u\in A^{\times}$ and an element $\chi\in M(A')$ satisfying 
$$\Ps(P,F,\phi)=u (z+\chi)^h.$$
The face $F$ is $z$-removable, if and only if, $hf^P_z \in F$,  $F\not\subset\{hf^P_z\}+\Map(P,\R_0)$ and after some coordinate change sending $x$ to $x$ itself for any $x\in P-\{z\}$ and sending $z$ to $z+\chi$ for some $\chi\in M(A')$, $F$ becomes a part of $\{hf^P_z\}+\Map(P,\R_0)$. 
If $\phi\in A^\times$, then any face $F$ of $\Gamma_+(P,\phi)$ is not $z$-removable, since $h=0$ and $F\subset\Map(P,\R_0)=\{hf^P_z\}+\Map(P,\R_0)$.
If $\dim A=1$, then any face $F$ of $\Gamma_+(P,\phi)$ is not $z$-removable, since $F\subset\Gamma_+(P,\phi)= \{hf^P_z\}+\Map(P,\R_0)$.

We would like to explain the relation betwen the concept of $z$-removable faces and Hironaka's maximal contact here. We assume that the field $k$ has characteristic zero, and consider any $z$-Weierstrass polynomial $\psi\in A$ over $P$ of positive degree. We take the unique pair of a positive integer $h$ and 
a mapping $\psi':\{0,1,\ldots,h-1\}\rightarrow M(A')$
satisfying the equality $\psi=z^h+\sum_{i=0}^{h-1} \psi'(i) z^i$. Let $\hat{z}=z+(\psi'(h-1)/h)\in M(A)$ and let $\hat{P}=\{\hat{z}\}\cup(P-\{z\})$. We know that $\Ord(\psi)\leq h$, $\hat{P}$ is a parameter system of $A$ and the Newton polyhedron $\Gamma_+(\hat{P},\psi)$ is of $\hat{z}$-Weierstrass type and has no $\hat{z}$-removable faces. Now, we assume moreover that $\Ord(\psi)=h$. This condition is equivalent to that $\Ord(\psi'(i))\geq h-i$ for any $i\in\{0,1,\ldots,h-1\}$. Then, the smooth subscheme $\Spec(A/\hat{z}A)$ of $\Spec(A)$ is Hironaka's maximal contact of the subscheme $\Spec(A/\psi A)$. (Giraud~\cite{G75}.)

Note that we cannot define the element $\hat{z}=z+(\psi'(h-1)/h) \in A$, if the characteristic of the ground field $k$ is positive and the characteristic divides $h$.

The concept of $z$-simple is also important.

We say that $\Gamma_+(P,\phi)$ is \emph{$z$-simple}, if $\Gamma_+(P,\phi)$ is of $z$-Weierstrass type and any compact face $F$ of $\Gamma_+(P,\phi)$ satisfies $\dim F\leq 1$.

If $\dim A\leq 2$, then always $\Gamma_+(P,\phi)$ is $z$-simple.
If $\Gamma_+(P,\phi)$ is $z$-simple, then $\Gamma_+(P,\phi)$ is of $z$-Weierstrass type. If $\Gamma_+(P,\phi)$ is $z$-simple, $\psi\in A$, $\omega\in A$ and $\phi=\psi\omega$, then both $\Gamma_+(P,\psi)$ and $\Gamma_+(P,\omega)$ are $z$-simple. (Lemma~\ref{Newton2}.14.)

Let $R$ be any regular local ring containing the field $k$ as a subring and the residue field $R/M(R)$ is isomorphic to $k$ as $k$-algebras. 

The completion $R^c$ is a complete regular local ring containing $R$ as a local subring and the residue field $R^c/M(R^c)$ is isomorphic to $k$ as $k$-algebras. 
Consider any parameter system $P$ of $R^c$ and any $\phi\in R$ with $\phi\neq 0$. We have the Newton polyhedron $\Gamma_+(P,\phi)$ of $\phi$ over $P$, if we regard the element $\phi$ of $R$ as an element of $R^c$.

\begin{lemma}
\label{coordinate change}
Consider any parameter system $P$ of $R$, any element $z\in P$, and any $w\in M(R^c)$ with $\partial w/\partial z\in R^{c\times}$. We denote $P_w=\{w\}\cup(P-\{z\})$. Let $R'$ be the localization of $k[P-\{z\}]$ by the maximal ideal $k[P-\{z\}]\cap M(R)=(P-\{z\}) k[P-\{z\}]$.
\begin{enumerate}
\item
$P$ is also a parameter system of $R^c$, and for any $\phi\in R$ with $\phi\neq 0$, $c(\Gamma_+(P,\phi))=1\Leftrightarrow \phi$ has normal crossings over $P$.
\item $P_w$ is a parameter system of $R^c$ with $w\in P_w$ and $P_w-\{w\}=P-\{z\}\subset R$. If $w\in M(R^h)$, then $P_w$ is a parameter system of $R^h$ with $w\in P_w$. If $w\in M(R)$, then $P_w$ is a parameter system of $R$ with $w\in P_w$.
\item
There exist uniquely $u\in R^{c\times}$ and $\chi\in M(R^{\prime c})$ with $w=uz+\chi$. If $w\in M(R^h)$, then $u\in R^{h\times}$ and $\chi\in M(R^{\prime h})$.
\item
There exist uniquely $v\in R^{c\times}$ and $\omega\in M(R^{\prime c})$ with $w=v(z+\omega)$. If $w\in M(R^h)$, then $v\in R^{h\times}$ and $\omega\in M(R^{\prime h})$.
\item
Assume that  $u\in R^{c\times}$, $\chi\in M(R^{\prime c})$, $v\in R^{c\times}$, $\omega\in M(R^{\prime c})$ and $w=uz+\chi=v(z+\omega)$. We take the unique pair $v_0\in R^{\prime c}$ and $v_1\in R^c$ with $v=v_0+zv_1$. 

Then, $v_0\in R^{\prime c\times}$, $\chi=v_0\omega$, $u=v_1(z+\omega)+v_0$ and $\Gamma_+(P-\{z\},\chi)= \Gamma_+(P-\{z\},\omega)$.
\end{enumerate}

The bijection $P_w\rightarrow P$ sending $w\in P_w$ to $z\in P$ and sending any $x\in P_w-\{w\}=P-\{z\}$ to $x\in P-\{z\}$ itself induces an isomorphism $\Map(P,\R)\rightarrow\Map(P_w,\R)$ of vector spaces over $\R$.
By this isomorphism we identify $\Map(P,\R)$ and $\Map(P_w,\R)$. $f^P_z\in \Map(P,\R)$ and $f^{P_w}_{w}\in \Map(P_w,\R)$ are identified. For any $x\in P-\{z\}$, $f^P_x\in \Map(P,\R)$ and $f^{P_w}_x\in \Map(P_w,\R)$ are identified.

\begin{enumerate}\setcounter{enumi}{5}
\item If $w'\in M(R^c)$ and $w=uw'$ for some $u\in R^{c\times}$, then $\partial w'/\partial z\in R^{c\times}$ and $\Gamma_+(P_w, \phi)= \Gamma_+(P_{w'}, \phi)$ for any $\phi\in R^c$ with $\phi\neq 0$.
\item Consider any $w'\in M(R^c)$ with $\partial w'/\partial z\in R^{c\times}$ and any $\phi\in R^c$ with $\phi\neq 0$. $\Gamma_+(P_w, \phi)$ is of $w$-Weiestrass type, if and only if, $\Gamma_+(P_{w'}, \phi)$ is of $w'$-Weiestrass type. If these equivalent conditions are satisfied, then the unique $w$-top vertex of $\Gamma_+(P_w, \phi)$ and the unique $w'$-top vertex of $\Gamma_+(P_{w'}, \phi)$ are equal.
\item For any $\phi\in R^h$ with $\phi\neq 0$, there exists $\psi\in R$ satisfying $(\phi R^h)\cap R=\psi R$ and $\psi\neq 0$.
\item
Assume moreover, that $R$ is a localization of a finitely generated $k$-algebra. 
Let $Q$ be any subset of $P$. Let $S$ be the localization of $k[Q]$ by the maximal ideal $k[Q]\cap M(R)=Qk[Q]$.
\begin{enumerate}
\item
If $Q=P$, then $S^c= R^c$ and $S^h=R^h$.
\item
$S^h=S^c\cap R^h$.
\item
Let $\pi:\Map(P,\Z_0)\rightarrow\Map(Q,\Z_0)$ denote the surjective mapping induced by the inclusion mapping $Q\rightarrow P$. For any $\phi\in R^h$ and any $\Lambda\in\Map(Q,\Z_0)$, $\Ps(P, \pi^{-1}(\Lambda),\phi)\in R^h$.
\item
For any $\omega\in\Map(P,\R_0)^\vee$ and any $\phi\in R^h$, $\In(P,\omega,\phi)\in R^h$.

For any $\omega\in(\Map(P,\R_0)^\vee)^\circ$ and any $\phi\in R^c$, $\In(P,\omega,\phi)\in k[P]$.
\end{enumerate}
\end{enumerate}
\end{lemma}
\begin{remark}
In claim $9$.(d), $(\Map(P,\R_0)^\vee)^\circ$ denotes the interior of the regular cone $\Map(P,\R_0)^\vee$. (Definition~\ref{dimension}.) $(\Map(P,\R_0)^\vee)^\circ=\sum_{x\in P}\R_+f^{P\vee}_x\subset\Map(P,\R)^*$.
\end{remark}

\begin{proof}
Claim $4$ follows from Henselian Weierstrass Theorem in Hironaka~\cite{H772}.
\end{proof}

Consider any parameter system $P$ of $R$, any element $z\in P$, any $w\in M(R^c)$ with $\partial w/\partial z\in R^{c\times}$, and any $\phi\in R$ with $\phi\neq 0$. We denote $P_w=\{w\}\cup(P-\{z\})$.

Since $R$ is a UFD, we have an invertible element $u\in R^\times$, a finite set $\Omega$ of irreducible elements of $R$ and a mapping $a:\Omega\rightarrow\Z_+$ satisfying $\phi=u\prod_{\omega\in\Omega}\omega^{a(\omega)}$. We take any $u, \Omega$ and $a$ satisfying these conditions. Let $$\Xi=\{\omega\in\Omega|\partial\omega/\partial w\in M(R^c)\text{, and any }x\in P_w-\{w\}\text{ does not divide }\omega\}.$$

We say that an element $\psi\in R$ is a \emph{main factor} of the triplet $(P_w,w,\phi)$, or a $w$-\emph{main factor over} $P_w$ of $\phi$, if $\psi=v\prod_{\omega\in\Xi}\omega^{a(\omega)}$ for some $v\in R^\times$.
Since $R$ is a UFD, the condition that $\psi\in R$ is a main factor of $(P_w,w,\phi)$ does not depend of the choice of $u, \Omega$ and $a$ we used for the definition.

If  $\psi\in R$ is a main factor of $(P_w,w,\phi)$ and $\psi'\in R$ is a main factor of $(P_w,w,\phi)$, then by definition, $\psi=v\psi'$ for some $v\in R^\times$, $\Gamma_+(P,\psi)=\Gamma_+(P,\psi')$, and any $x\in P_w-\{w\}$ does not divide $\psi$.

We consider the case where $\Gamma_+(P_w,\phi)$ is of $w$-Weierstrass type. Let $\psi\in R$ be any main factor of $(P_w,w,\phi)$. The Newton polyhedron $\Gamma_+(P_w,\psi)$ does not depend on the choice of the main factor $\psi$ and it is of $w$-Weierstrass type. 
Furthermore, there exists uniquely a non-negative integer $h$ such that $\{hf^{P_w}_w\}$ is the unique $w$-top vertex of  $\Gamma_+(P_w,\psi)$. We take $h\in\Z_0$ satisfying this condition and we define
$$\Inv(P_w,w,\phi)=h\in\Z_0.$$

The non-negative integer $\Inv(P_w,w,\phi) $ is our main invariant measuring the badness of the singularity $\phi$. 

\begin{lemma}
\label{main factor}
Consider any parameter system $P$ of $R$, any element $z\in P$, any $w\in M(R^c)$ with $\partial w/\partial z\in R^{c\times}$, and any $\phi\in R$ with $\phi\neq 0$. We denote $P_w=\{w\}\cup(P-\{z\})$.

The bijection $P_w\rightarrow P$ sending $w\in P_w$ to $z\in P$ and sending any $x\in P_w-\{w\}=P-\{z\}$ to $x\in P-\{z\}$ itself induces an isomorphism $\Map(P,\R)\rightarrow\Map(P_w,\R)$ of vector spaces over $\R$.
By this isomorphism we identify $\Map(P,\R)$ and $\Map(P_w,\R)$.

\begin{enumerate}
\item
For any $\psi\in R$, $\psi$ is a main factor of $(P_w,w,\phi)$, if and only if, the following three conditions are satisfied:
\begin{enumerate}
\item $\phi=u(\prod_{x\in P_w-\{w\}}x^{a(x)})(\prod_{\omega\in\Omega}\omega^{b(\omega)})\psi$ and $\partial \omega/\partial w\in R^{c\times}$ for any $\omega\in \Omega$, for some $u\in R^\times$, some mapping $a:P_w-\{w\}\rightarrow \Z_0$, some finite subset $\Omega$ of $M(R)$ and some mapping $b:\Omega\rightarrow\Z_+$.
\item Any element of $P_w-\{w\}$ does not divide $\psi$.
\item Any $\omega\in M(R)$ with $\partial\omega/\partial w\in R^{c\times}$ does not divide $\psi$.
\end{enumerate}
\item
If both $\psi\in R$ and $\psi'\in R$ are main factors of $(P_w,w,\phi)$, then $\psi=u\psi'$ for some $u\in R^\times$ and $\Gamma_+(P,\psi)= \Gamma_+(P,\psi')$.
\item
For any $\psi\in R$ and any $w'\in M(R^c)$ with $\partial w'/\partial z\in R^{c\times}$, $\psi$ is a main factor of $(P_w,w,\phi)\Leftrightarrow \psi$ is a main factor of $(P_{w'}, w',\phi)$.
\item
Let $\psi\in R$ be any main factor of $(P_w,w,\phi)$. 
$\psi\not\in xR^c$ for any $x\in P_w-\{w\}$. If $w\in M(R)$, then $\psi\not\in wR^c$.
$\psi\not\in zR^c$.
\item The following three conditions are equivalent:
\begin{enumerate}
\item Any main factor of $(P_w,w,\phi)$ is an invertible element of $R$.
\item Some main factor of $(P_w,w,\phi)$ is an invertible element of $R$.
\item
$\phi=u(\prod_{x\in P_w-\{w\}}x^{a(x)})(\prod_{\omega\in\Omega}\omega^{b(\omega)})$ and $\partial \omega/\partial w\in R^{c\times}$ for any $\omega\in \Omega$, for some $u\in R^\times$, some mapping $a:P_w-\{w\}\rightarrow \Z_0$, some finite subset $\Omega$ of $M(R)$ and some mapping $b:\Omega\rightarrow\Z_+$.
\end{enumerate}

If $\Gamma_+(P_w,\phi)$ is of $w$-Weierstrass type, then the condition below is also equivalent to the above three conditions.
\begin{enumerate}\setcounter{enumii}{3} 
\item $\Inv(P_w,w,\phi)=0$.
\end{enumerate}
\end{enumerate}

Below, we assume that $\Gamma_+(P_w,\phi)$ is of $w$-Weierstrass type.
\begin{enumerate}
\setcounter{enumi}{5}
\item $\Inv(P_w,w,\phi) \neq 1$.
\item
Let $\psi\in R$ be any main factor of $(P_w,w,\phi)$. 
\begin{enumerate}
\item
$\Gamma_+(P_w,\psi)$ is of $w$-Weierstrass type and $\{\Inv(P_w,w,\phi)f^{P_w}_w\}$ is the unique $w$-top vertex of $\Gamma_+(P_w,\psi)$.
\item Let $w'\in M(R^c)$ be any element with $\partial w'/\partial z\in R^{c\times}$. 
$\Inv(P_w,w,\phi)=\Inv(P_{w'}, w',\phi)$.
\item There exist uniquely $u\in R^{c\times}$ and $\psi'\in R^c$ such that $\psi'$ is a $w$-Weierstrass polynomial over $P_w$ and $\psi=u\psi'$.
\end{enumerate}
We take $u\in R^{c\times}$ and $\psi'\in R^c$ satisfying the above conditions.
\begin{enumerate}\setcounter{enumii}{3}
\item
$\Gamma_+(P_w,\psi)= \Gamma_+(P_w,\psi')$.
\item
$\Inv(P_w,w,\phi)=\Ht(w,\Gamma_+(P_w,\psi))+\Ord(P_w,f^{P_w\vee}_w,\psi)= \deg(P_w,w,$\break$\psi')\geq\Ord(\psi)=\Ord(\psi')$.

If $w\in M(R)$, then $\Ord(P,f^{P_w\vee}_w,\psi)=0$ and $\Inv(P_w,w,\phi)=\Ht(w,$\break$\Gamma_+(P_w,\psi))$.
\end{enumerate}
\item
Assume both $\psi\in R$ and $\psi'\in R$ are main factors of $(P_w,w,\phi)$. We take $u\in R^{c\times}$ and $\psi''\in R^c$ such that $\psi''$ is a $w$-Weierstrass polynomial over $P_w$ and $\psi=u\psi''$. Consider any face $F$ of $\Gamma_+(P_w,\psi)= \Gamma_+(P_w,\psi') = \Gamma_+(P_w,\psi'')$. 

$F$ is a $w$-removable face of $\Gamma_+(P_w,\psi)$, if and only if, $F$ is a $w$-removable face of $\Gamma_+(P_w,\psi')$ , if and only if, $F$ is a $w$-removable face of $\Gamma_+(P_w,\psi'')$.
\item
Assume $\Inv(P_w, w, \phi)=0$.

There exist $u\in R^\times$, a mapping $a:P_w-\{w\}\rightarrow\Z_0$, a finite subset $\Omega$ of $M(R)$ and a mapping $b:\Omega\rightarrow\Z_+$ satisfying the following three conditions:
\begin{enumerate}
\item
$\phi=u\prod_{x\in P_w-\{w\}}x^{a(x)}\prod_{\omega\in\Omega}\omega^{b(\omega)}$.
\item $\partial \omega/\partial w\in R^{c\times}$ for any $\omega\in\Omega$.
\item If $\omega=v\omega'$ for some $\omega\in\Omega $, some $\omega'\in\Omega $ and some $v\in R^\times$, then $v=1$ and $\omega=\omega'$.
\end{enumerate}

Let $w'\in M(R^c)$ be any element with $\partial w'/\partial z\in R^{c\times}$.
There exist $u'\in R^\times$, a mapping $a':P-\{z\}\rightarrow\Z_0$, a finite subset $\Omega'$ of $R$ and a mapping $b':\Omega\rightarrow\Z_+$ satisfying the following three conditions:
\begin{enumerate}
\item
$\phi=u'\prod_{x\in P_{w'}-\{w'\}}x^{a' (x)}\prod_{\omega\in\Omega'}\omega^{b'(\omega)}$.
\item $\partial \omega/\partial w'\in R^{c\times}$ for any $\omega\in\Omega'$.
\item If $\omega=v\omega'$ for some $\omega\in\Omega'$, some $\omega'\in\Omega'$ and some $v\in R^\times$, then $v=1$ and $\omega=\omega'$.
\end{enumerate}

If $u,a,\Omega,b,u',a',\Omega'$ and $b'$ satisfy the above conditions, then $\sharp\Omega=\sharp\Omega'$.
\end{enumerate}
\end{lemma}

For any $w\in M(R^c)$ with $\partial w/\partial z\in R^{c\times}$ and any $\phi\in R$ such that $\phi\neq 0$, $\Gamma_+(P_w,\phi)$ is of $w$-Weierstrass type, and $\Inv(P_w, w, \phi)=0$ where $P_w=\{w\}\cup(P-\{z\})$, we take $u, a, \Omega, b$ satisfying the conditions in the above Lemma~\ref{main factor}.9 and we define
$$\Inv 2(P_w,w,\phi)= \sharp\Omega\in\Z_0.$$

\begin{lemma}
\label{inv2}
Consider any parameter system $P$ of $R$, any element $z\in P$, any $w\in M(R^c)$ with $\partial w/\partial z\in R^{c\times}$. We denote $P_w=\{w\}\cup(P-\{z\})$.

Consider any element $\phi\in R$ such that $\phi\neq 0$, $\Gamma_+(P_w,\phi)$ is of $w$-Weierstrass type, and $\Inv(P_w, w, \phi)=0$. 
\begin{enumerate}
\item
If $\Inv 2(P_w,w,\phi)\leq 1$, then $\phi$ has normal crossings.
\item
If $\Inv 2(P_w,w,\phi)\geq 1$, then there exists $\bar{z}\in M(R)$ such that $\partial \bar{z}/\partial z\in R^\times$ and $\bar{z}$ divides $\phi$.
\item
Let $w'\in M(R^c)$ be any element with $\partial w'/\partial z\in R^{c\times}$.
$\Inv(P_{w'}, w', \phi)=0$ and $\Inv 2(P_{w'}, w',\phi)=\Inv 2(P_w,w,\phi)$.
\end{enumerate}
\end{lemma}

\section{Basic scheme theory}
\label{scheme}

We develop the basic scheme theory. 

Let $X$ be any scheme, and let $\mathcal{I}$ be any ideal sheaf in the structure sheaf $\mathcal{O}_X$ of $X$, in other words, any sheaf of $\mathcal{O}_X$-modules which is a subsheaf of $\mathcal{O}_X$. The ideal sheaf $\mathcal{I}$ is called \emph{locally principal}, if for any $a\in X$ there exists $\phi\in \mathcal{O}_{X,a}$ such that $\phi$ is not a zero-divisor of $\mathcal{O}_{X,a}$ and $\mathcal{I}_a=\phi \mathcal{O}_{X,a}$, where $\mathcal{I}_a$ denotes the stalk of $\mathcal{I}$ at $a$. Note that for any scheme $Y$ and for any morphism $\gamma:Y\rightarrow X$ of schemes, the pull-back $\gamma^*\mathcal{I}$ of $\mathcal{I}$ as an ideal sheaf is defined, and $\gamma^*\mathcal{I}$ is a sheaf of $\mathcal{O}_{Y}$-modules which is a subsheaf of $\mathcal{O}_{Y}$.

Grothendieck has shown that there exists a scheme $X'$ and a morphism $\sigma:X'\rightarrow X$ satisfying the following universal mapping property:  
\begin{enumerate}
\item The ideal sheaf $\sigma^*\mathcal{I}$ is locally principal. 
\item If $Y$ is a scheme, $\gamma:Y\rightarrow X$ is a morphism, and the ideal sheaf $\gamma^*\mathcal{I}$ is locally principal, then there exists uniquely a morphism $\tau:Y\rightarrow X'$ with $\sigma\tau=\gamma$.
\end{enumerate}

By the universal mapping property we know that the pair $(X', \sigma)$ satisfying the above conditions is unique up to isomorphism of schemes over $X$. The pair $(X', \sigma)$ satisfying the above conditions is called the \emph{blowing-up with center in} an ideal sheaf $\mathcal{I}$, or the \emph{blowing-up with center in} $Z$, where $Z$ denotes the closed subscheme of $X$ defined by the ideal sheaf $\mathcal{I}$. Note that any closed subscheme of $X$ has a unique ideal sheaf in $\mathcal{O}_{X}$ defining it. If $\mathcal{I}$ is locally principal, then $\sigma$ is an isomorphism. When a closed subset $Z$ of $X$ is given, we take the unique ideal sheaf $\mathcal{I}$ in $\mathcal{O}_X$ defining the reduced scheme structure on $Z$ and we call the blowing-up with center in $\mathcal{I}$ the blowing-up with center in $Z$.

Let $(X', \sigma)$ be the blowing-up with center in $\mathcal{I}$. By $Z$ we denote the closed subscheme of $X$ defined by the ideal sheaf $\mathcal{I}$. We call the inverse image $\sigma^{-1}(Z)$ the \emph{exceptional divisor} of $\sigma$. For any closed irreducible subset $W$ of $X$ with $W\not\subset Z_{red}$, the closure in $X'$ of $\sigma^{-1}(W-Z)$ is called the \emph{strict transform} of $W$ by $\sigma$. If $X$ is separated, noetherian, irreducible and smooth, and $Z$ is irreducible and smooth, then $X'$ is also separated, noetherian, irreducible and smooth, and $\sigma^{-1}(Z)$ is irreducible and smooth.

For any ring $R$ and an ideal $I$ of $R$, we regard the affine scheme $\Spec(R/I)$ as the closed subcheme of the affine scheme $\Spec(R)$ by using the closed embedding $\Spec(R/I)\rightarrow\Spec(R)$ induced by the canonical surjective ring homomorphism $R\rightarrow R/I$ to the residue ring.

Let $X$ be any separated noethrian irreducible smooth scheme with $\dim X\geq 1$, let $D$ be any effective divisor of $X$, and let $a\in X$ be any point. We consider the natural morphism $\delta:\Spec(\mathcal{O}_{X,a})\rightarrow X$. It is dominant, and the pull-back $\delta^*D$ of $D$ by this morphism $\delta$ is defined.
We say that $D$ has \emph{normal crossings at} $a\in X$, if there exist a parameter system $P$ of the local ring $\mathcal{O}_{X,a}$ of $X$ at $a$, and an element $\Lambda\in\Map(P,\Z_0)$ such that $$\delta^*D=
\Spec(\mathcal{O}_{X,a}/\prod_{x\in P}x^{\Lambda(x)}\mathcal{O}_{X,a}).$$
We say that $D$ has \emph{normal crossings} or $D$ is a \emph{normal crossing divisor}, if it has normal crossings at any point of $X$. It follows from definition that any component of $D$ is smooth and the intersection of any finite number of components of $D$ is smooth (however, the intersection of two or more components of $D$ may be empty or reducible), if $D$ has normal crossings.

Here we give the definition of the concept of \emph{normal crossing schemes} over an algebraically closed field and introduce some notations associated with it. Let $k$ denote any algebraically closed field below in this section.

A pair
$$(X, D),$$
satisfying the following five conditions is called a  \emph{normal crossing scheme} over $k$. 
\begin{enumerate}
\item The first item $X$ is a separated noetherian irreducible smooth scheme over $k$ with $\dim X\geq 1$ such that any closed point $a\in X$ is a $k$-valued point.
\item The second item $D$ is a non-zero effective normal crossing divisor of $X$.
\end{enumerate}

We use the following notations: The set of components of $D$ is denoted by $\Comp(D)$. For any point $a\in X$ we denote $$\Comp(D)(a)=\{C\in \Comp(D)|a\in C\},$$ and $$(D)_0=\{a\in X|\sharp\Comp(D)(a)=\dim X\}.$$ For any $a\in(D)_0$ we write
$$U(X,D,a)=X-(\bigcup_{C\in \Comp(D)- \Comp(D)(a)}C).$$
We write simply $U(a)$, instead of $U(X,D,a)$, when we need not refer to the pair $(X,D)$.

\begin{enumerate}
\setcounter{enumi}{2}
\item For any subset $Q$ of $\Comp(D)$, $\bigcap_{C\in Q}C$ is irreducible.
\item For any subset $Q$ of $\Comp(D)$ with $\bigcap_{C\in Q}C\neq\emptyset$, there exists $a\in(D)_0$ such that $Q\subset\Comp(D)(a)$.
\item For any $a\in(D)_0$, $U(a)$ is an affine open subset of $X$.
\end{enumerate}

Let $(X, D)$ be a normal crossing scheme over $k$. For any $a\in(D)_0$, we consider a mapping
$$\xi_a:\Comp(D)(a)\rightarrow \mathcal{O}_X(U(a)).$$
Consider any point $a\in(D)_0$. If $\xi_a$ satisfies the following two conditions, then we call $\xi_a$ a \emph{coordinate system} of $(X, D)$ at $a$:

\begin{enumerate} 
\item For any $C\in \Comp(D)(a)$ we have
$$C\cap U(a)=
\Spec(\mathcal{O}_X(U(a))/\xi_a(C)\mathcal{O}_X(U(a))).$$
\item 
For any $k$-valued point $b$ in $U(a)$, the set $\{\xi_{a}(C)- \xi_{a}(C)(b)|C\in\Comp(D)(a)\}$ is a parameter system of the local ring $\mathcal{O}_{X,b}$ of $X$ at $b$. Here $\xi_{a}(C)(b)\in k$ denotes the value of  $\xi_{a}(C)\in \mathcal{O}_X(U(a))$ at $b$.
\end{enumerate}

If $\xi_a$ is a coordinate system of $(X, D)$ at $a$ for any $a\in(D)_0$, then we call the collection $\xi=\{\xi_a|a\in(D)_0\}$ a \emph{coordinate system} of $(X, D)$. For a coordinate system $\xi$ of $(X, D)$ we denote the element of $\xi$ corresponding to $a\in(D)_0$ by $\xi_a$.

A triplet $(X, D, \xi)$ such that the pair $(X,D)$ is a normal crossing scheme over $k$ and $\xi$ is a coordinate system of $(X,D)$ is called a \emph{coordinated normal crossing scheme} over $k$.

\begin{example}
Let $R$ be any regular local ring such that $R$ contains $k$ as a subring, the residue field $R/M(R)$ is isomorphic to $k$ as $k$-algebras, and $\dim R\geq 1$; let $P$ be any parameter system of $R$; and let $\Lambda\in\Map(P,\Z_+)$.

Note that $M(R)\in\Spec(R)$ and $M(R)$ is the unique closed point of $\Spec(R)$. Let $D=\Spec(R/\prod_{x\in P}x^{\Lambda(x)}R)$. The pair $(\Spec(R),D)$ is a normal crossing scheme over $k$. We have $(D)_0=\{M(R)\}$, $\Comp(D)=\Comp(D)(M(R))=\{\Spec(R/xR)|x\in P\}$, and $U(\Spec(R),D, M(R))=\Spec(R)$. For any $x\in P$, we put $\xi_{M(R)}(\Spec(R/xR$\break$))=x$. We obtain a mapping $\xi_{M(R)}: \Comp(D)(M(R))\rightarrow\mathcal{O}_{\Spec(R)}( U(\Spec(R),D, $\break$M(R)))$.
The mapping $\xi_{M(R)}$ is a coordinate system of $(\Spec(R),D)$ at $M(R)$, and the triplet $(\Spec(R),D,\{\xi_{M(R)}\})$ is a coordinated normal crossing scheme over $k$.

We consider the subring $k[P]$ of $R$. We denote $\bar{M}= k[P]\cap M(R)=Pk[P]$. $\bar{M}\in\Spec(k[P])$ and $\bar{M}$ is a closed point of $\Spec(k[P])$. Let $\bar{D}=\Spec(k[P]/\prod_{x\in P}x^{\Lambda(x)} $\break$k[P])$. The pair $(\Spec(k[P]), \bar{D})$ is a normal crossing scheme over $k$. We have $(\bar{D})_0=\{\bar{M}\}$, $\Comp(\bar{D})=\Comp(\bar{D})(\bar{M})=\{\Spec(k[P]/xk[P])|x\in P\}$, and $U(\Spec(k[P]), $\break$\bar{D}, \bar{M})=\Spec(k[P])$. For any $x\in P$, we put $\xi_{\bar{M}}(\Spec(k[P]/xk[P]))=x$. We obtain a mapping $\xi_{\bar{M}}: \Comp(\bar{D})(\bar{M})\rightarrow\mathcal{O}_{\Spec(k[P])}( U(\Spec(k[P]), \bar{D}, \bar{M}))$.
The mapping $\xi_{\bar{M}}$ is a coordinate system of $(\Spec(k[P]), \bar{D})$ at $\bar{M}$, and the triplet $(\Spec(k[P]), \bar{D},$\break$\{\xi_{\bar{M}}\})$ is a coordinated normal crossing scheme over $k$.
\end{example}

The four lemmas below easily follow from definitions.

\begin{lemma}
\label{ncv}
Let $(X, D)$ be a normal crossing scheme over $k$.
\begin{enumerate}
\item
The set $\Comp(D)$ is non-empty and finite.
\item
Consider any non-empty subset $Q$ of $\Comp(D)$ with $\bigcap_{C\in Q}C\neq\emptyset$. $\bigcap_{C\in Q}C$ is a closed irreducible smooth subset of $X$, and $\dim \bigcap_{C\in Q}C=\dim X-\sharp Q$. If we give the reduced scheme structure to any $C\in Q$, then the intersection scheme $\bigcap_{C\in Q}C$ is reduced and smooth.
\item
The set $(D)_0$ is a non-empty finite set of $k$-valued points of $X$.
\item
For any $a\in(D)_0$ and any $b\in(D)_0$, $\Comp(D)(a)=\Comp(D)(b)$, if and only if, $a=b$.
\item
$$X=\bigcup_{a\in(D)_0}U(a).$$
\end{enumerate}

Let $Q$ be any subset of $\Comp(D)$ with $\sharp Q\geq 2$ and $\bigcap_{C\in Q}C\neq\emptyset$. We denote $Z=\bigcap_{C\in Q}C$, and the blowing-up with center in $Z$ by $\sigma:X'\rightarrow X$. Furthermore, by $E'$ we denote the exceptional divisor of $\sigma$, and by $C'$ we denote the strict transform of $C\in\Comp(D)$ by $\sigma$ for any $C\in\Comp(D)$.

\begin{enumerate}
\setcounter{enumi}{5}
\item
The pair $(X', \sigma^*D)$ is a normal crossing scheme over $k$.
\item
$E'\in\Comp(\sigma^*D)$. $\Comp(\sigma^*D)-\{E'\}=\{C'|C\in\Comp(D)\}$.
$\sharp\Comp(\sigma^*D)=\sharp\Comp(D)+1$.
\item
For any $C\in\Comp(D)- Q$, we have $\sigma^*C=C'$.
For any $C\in Q$, we have $\sigma^*C=C'+E'$.
\item
$\sigma((\sigma^*D)_0)=(D)_0$. 
\item
For any $a\in(D)_0$ with $a\not\in Z$, we have
$\sharp\sigma^{-1}(a)\cap(\sigma^*D)_0=1$, and the unique element $a'$ in $\sigma^{-1}(a)\cap(\sigma^*D)_0$ satisfies 
$\{a'\}=
\bigcap_{C\in\Comp(D)(a)}C'$,
$\Comp(\sigma^*D)(a')=\{C'|C\in\Comp(D)(a)\}$,
and $U(X', \sigma^*D, a')=\sigma^{-1}(U (X, D, $\hfill\break$a))$.

If moreover, a mapping
$\xi_a:\Comp(D)(a)\rightarrow \mathcal{O}_X(U(X,D,a))$ is a coordinate system of $(X,D)$ at $a$, then there exists a unique coordinate system $\xi'_{a'}:\Comp(\sigma^*D)(a')\rightarrow \mathcal{O}_{X'} (U(X',\sigma^*D,a'))$ of $(X',\sigma^*D)$ at $a'$ satisfying $\sigma^*(\xi_a(C))=\xi'_{a'}(C')$ for any $C\in\Comp(D)(a)$, where 
$\sigma^*: \mathcal{O}_X(U(X,D,a))\rightarrow \mathcal{O}_{X'} (U(X',\sigma^*D,a'))$ denotes the ring homomorphism induced by $\sigma$.
\item
For any $a\in(D)_0$ with $a\in Z$, we have
$\sharp\sigma^{-1}(a)\cap(\sigma^*D)_0=\sharp Q\geq 2$, and
there exists a unique one-to-one mapping $a':Q\rightarrow\sigma^{-1}(a)\cap(\sigma^*D)_0$ such that for any $B\in Q$ we have
$\{a'(B)\}=
E'\cap\bigcap_{C\in\Comp(D)(a)-\{B\}}C'$,
$\Comp(\sigma^*D)(a'(B))=\{E'\}\cup\{C'|C\in\Comp(D)(a)-\{B\}\}$, and
\hfill\break$U(X', \sigma^*D, a'(B))= \sigma^{-1}(U (X, D, a))-B'$.

If moreover, a mapping
$\xi_a:\Comp(D)(a)\rightarrow \mathcal{O}_X(U(X,D,a))$ is a coordinate system of $(X,D)$ at $a$, then for any $B\in Q$, there exists a unique coordinate system $\xi'_{a'(B)}:\Comp(\sigma^*D)(a'(B))\rightarrow \mathcal{O}_{X'} (U(X',\sigma^*D,a'(B)))$ of $(X',\sigma^*D)$ at $a'(B)$ satisfying 
\begin{equation*}
\sigma^*(\xi_a(C))=
\begin{cases}
\xi'_{a'(B)}(E')&\text{if $C=B$},\\
\xi'_{a'(B)}(C')\xi'_{a'(B)}(E')&\text{if $C\in Q-\{B\}$},\\
\xi'_{a'(B)}(C')&\text{if $C\in\Comp(D)(a)- Q$},
\end{cases}
\end{equation*}
for any $C\in\Comp(D)(a)$, where 
$\sigma^*: \mathcal{O}_X(U(X,D,a))\rightarrow$\hfill\break $\mathcal{O}_{X'} (U(X',\sigma^*D,a'(B)))$ denotes the ring homomorphism induced by $\sigma$.
\end{enumerate}
\end{lemma}

Let $(X, D)$ be a normal crossing scheme over $k$. 

We call a non-empty closed subscheme $Z$ of $X$ such that there exists a non-empty subset $Q$ of $\Comp(D)$ satisfying $Z=\bigcap_{C\in Q}C$ a \emph{closed stratum} of $D$. We call a blowing-up whose center is a closed stratum of $D$ an \emph{admissible} blowing-up over $D$.

Let $Q$ be any subset of $\Comp(D)$ with $\sharp Q\geq 2$ and $\bigcap_{C\in Q}C\neq\emptyset$. Let $\sigma:X'\rightarrow X$ denote the admissible blowing-up with center in $\bigcap_{C\in Q}C$. We call the normal crossing scheme $(X', \sigma^*D)$ over $k$ the \emph{pull-back} of $(X,D)$ by $\sigma$. 
Let $a'\in(\sigma^*D)_0$ be any point, and let $\xi_{\sigma(a')}:\Comp(D)({\sigma(a')})\rightarrow \mathcal{O}_X(U(X,D,{\sigma(a')}))$ be a coordinate system of $(X,D)$ at ${\sigma(a')}$. We have the coordinate system $\xi'_{a'}:\Comp(D)({a'})\rightarrow \mathcal{O}_{X'} (U(X',\sigma^*D,{a'}))$  of $(X',\sigma^*D)$ at $a'$ described in Lemma~\ref{ncv}.$10$ or Lemma~\ref{ncv}.$11$. The coordinate system $\xi'_{a'}$ is called the \emph{pull-back} of $\xi_{\sigma(a')}$ at $a'$ by $\sigma$. Let $\xi=\{\xi_a|a\in(D)_0\}$ be a coordinate system of $(X,D)$. We denote 
$$\sigma^*\xi=\{\xi'_{a'}|a'\in(\sigma^*D)_0\},$$
and call $\sigma^*\xi$ the \emph{pull-back} of $\xi$ by $\sigma$. Note that triplets $(X,D,\xi)$ and $(X',\sigma^*D,\sigma^*\xi)$ are coordinated normal crossing scheme over $k$. We call the coordinated normal crossing scheme $(X',\sigma^*D,\sigma^*\xi)$ over $k$ the \emph{pull-back} of  $(X,D,\xi)$ by $\sigma$.

Let $X'$ be a scheme, and let $\sigma:X'\rightarrow X$ be a morphism. We call $\sigma$ an \emph{admissible composition of blowing-ups} over $D$, if there exist a non-negative integer $m$, $(m+1)$ of normal crossing schemes $(X(i), D(i))$ over $k$, $i\in\{0,1,\ldots,m\}$, and $m$ of morphisms $\sigma(i):X(i)\rightarrow X(i-1), i\in\{1,2,\ldots,m\}$ satisfying the following two conditions:
\begin{enumerate}
\item $X(0)=X,D(0)=D,X(m)=X'$ and $\sigma=\sigma(1)\sigma(2)\cdots\sigma(m)$.
\item For any $i\in\{1,2,\ldots\,m\}$, $\sigma(i)$ is an admissible blowing-up over $D(i-1)$ and $D(i)=\sigma(i)^*D(i-1)$.
\end{enumerate}
If moreover, the center of $\sigma(i)$ has codimension two for any $i\in\{1,2,\ldots\,m\}$, then we call $\sigma$ an admissible composition of blowing-ups \emph{with center of codimension two} over $D$.

\begin{lemma}
\label{admissible}
\mbox{\rm{1.}}
Let $(X,D)$ be a normal crossing scheme over $k$, let $X'$ be a scheme, and let $\sigma: X'\rightarrow X$ be an admissible composition of blowing-ups over $D$. Then, the pair $(X',\sigma^*D)$ is a normal crossing scheme over $k$.

\noindent \mbox{\rm{2.}}
Let $(X,D,\xi)$ be a coordinated normal crossing scheme over $k$, let $X'$ be a scheme, and let $\sigma: X'\rightarrow X$ be an admissible composition of blowing-ups over $D$. 
Assume that $m\in\Z_0$, $(m+1)$ of normal crossing schemes $(X(i), D(i))$ over $k$, $i\in\{0,1,\ldots,m\}$, and $m$ of morphisms $\sigma(i):X(i)\rightarrow X(i-1), i\in\{1,2,\ldots,m\}$ satisfy the above two conditions. We write
$\sigma^*\xi=\sigma(m)^*\sigma(m-1)^*\cdots\sigma(1)^*\xi$.
Then, the triplet $(X',\sigma^*D,\sigma^*\xi)$ is a coordinated normal crossing scheme over $k$, and the coordinate system $\sigma^*\xi$ of $(X',\sigma^*D)$ does not depend on the choice of $m\in\Z_0$, $(m+1)$ of normal crossing schemes $(X(i), D(i))$ over $k$, $i\in\{0,1,\ldots,m\}$, and $m$ of morphisms $\sigma(i):X(i)\rightarrow X(i-1), i\in\{1,2,\ldots,m\}$ satisfying the two conditions.
\end{lemma} 

Let $(X,D)$ be a normal crossing scheme over $k$, let $X'$ be a scheme, and let $\sigma: X'\rightarrow X$ be an admissible composition of blowing-ups over $D$. We call the normal crossing scheme $(X',\sigma^*D)$ over $k$ the \emph{pull-back} of $(X,D)$ by $\sigma$.

Let $(X,D,\xi)$ be a coordinated normal crossing scheme over $k$, let $X'$ be a scheme, and let $\sigma: X'\rightarrow X$ be an admissible composition of blowing-ups over $D$. 
Choosing $m\in\Z_0$, $(m+1)$ of normal crossing schemes $(X(i), D(i))$ over $k$, $i\in\{0,1,\ldots,m\}$, and $m$ of morphisms $\sigma(i):X(i)\rightarrow X(i-1), i\in\{1,2,\ldots,m\}$ satisfying the above two conditions, we define the coordinate system $\sigma^*\xi$ of $(X',\sigma^*D)$ by putting
$\sigma^*\xi=\sigma(m)^*\sigma(m-1)^*\cdots\sigma(1)^*\xi$.
The coordinate system $\sigma^*\xi$ does not depend on the choice of $m\in\Z_0$, $(m+1)$ of normal crossing schemes $(X(i), D(i)), i\in\{0,1,\ldots,m\}$, and $m$ of morphisms $\sigma(i):X(i)\rightarrow X(i-1), i\in\{1,2,\ldots,m\}$ satisfying the above two conditions. We call $\sigma^*\xi$ the \emph{pull-back} of $\xi$ by $\sigma$. For any $a'\in (\sigma^*D)_0$, we call $(\sigma^*\xi)_{a'}$ the \emph{pull-back} of $\xi_{\sigma(a')}$ by $\sigma$. We call the coordinated normal crossing scheme $(X',\sigma^*D, \sigma^*\xi)$ over $k$ the \emph{pull-back} of $(X,D, \xi)$ by $\sigma$.

\begin{lemma}
\label{ext}
Consider any scheme $X$ over $k$ and any divisors $D$ and $D'$ of $X$ such that both 
$(X, D)$ and $(X, D')$ are normal crossing schemes over $k$ and $\Supp(D)=\Supp(D')$. 
\begin{enumerate}
\item $\Comp(D)=\Comp(D')$,  $\Comp(D)(a)= \Comp(D')(a)$ for any $a\in X$, $(D)_0=(D')_0$, and $U(X, D',a)=U(X,D,a)$ for any $a\in(D)_0$.
\item If $\xi=\{\xi_a|a\in(D)_0\}$ is a coordinate system of $(X,D)$, then $\xi$ is a coordinate system of $(X,D')$.
\end{enumerate}
\end{lemma}

\begin{lemma}
\label{pull back blowing-ups}
Recall that $k$ denotes any algebraically closed field. Let $R$ be any regular local ring such that $R$ contains $k$ as a subring, the residue field  $R/M(R)$ is isomorphic to $k$ as algebras over $k$, and $\dim R\geq 2$, let $P$ be any parameter system of $R$, and let $z\in P$ be any element. 

Let $R'$ denote the localization of $k[P - \{z\}]$ by the maximal ideal $k[P - \{z\}]\cap M(A)=( P - \{z\})k[P - \{z\}]$. The ring $R'$ is a regular local subring of $R$. The set $P-\{z\}$ is a parameter system of $R'$.

Let $\sigma':X'\rightarrow\Spec(R')$ be any composition of finite blowing-ups with center in a closed irreducible smooth subscheme. The scheme $X'$ is smooth. We consider a morphism $\Spec(R)\rightarrow\Spec(R')$ induced by the inclusion ring homomorphism $R'\rightarrow R$, the product scheme $X=X'\times_{\Spec(R')}\Spec(R)$, the projection $\sigma:X\rightarrow\Spec(R)$, and the projection $\pi:X\rightarrow X'$. We know the following:
\begin{enumerate}
\item The morphism $\sigma$ is a composition of finite blowing-ups with center in a closed irreducible smooth subscheme. The scheme $X$ is smooth.
\item The pull-back $\sigma^*\Spec(R/zR)$ of the prime divisor $\Spec(R/zR)$ on $\Spec(R)$ by $\sigma$ is a smooth prime divisor on $X$, and $\sigma^*\Spec(R/zR)\supset\sigma^{-1}(M(R))$.
\item The projection $\pi:X\rightarrow X'$ induces an isomorphism $\sigma^*\Spec(R/zR)\rightarrow X'$.
\item
For any closed point $a\in X$, any $w\in M(R^h)$ with $\partial w/\partial z\in R^{h\times}$ and any parameter system $Q'$ of the local ring $\mathcal{O}_{X',\pi(a)}$ of $X'$ at $\pi(a)$, $\sigma(a)=M(R)$ and $\{\sigma^*(w)\}\cup\pi^*(Q') $ is a parameter system of the Henselization $\mathcal{O}_{X,a}^h$ of the local ring $\mathcal{O}_{X,a}$ of $X$ at $a$ with $\pi^*(Q')\subset\mathcal{O}_{X,a}$, where $\sigma^*:R^h\rightarrow \mathcal{O}_{X,a}^h$ denotes the homomorphism of local $k$-algebras induced by $\sigma$ on the Henselizations of local rings and $\pi^*:\mathcal{O}_{X',\pi(a)}\rightarrow \mathcal{O}_{X,a}\subset \mathcal{O}_{X,a}^h$ denotes the homomorphism of local $k$-algebras induced by $\pi$.

If $w\in M(R)$, then $\{\sigma^*(w)\}\cup\pi^*(Q') $ is a parameter system of the local ring $\mathcal{O}_{X,a}$ of $X$ at $a$
\item
For any closed point $a'\in X'$, $\sigma'(a')=M(R')$ and the fiber $\pi^{-1}(a')$ of $\pi$ over $a'$ is isomorphic to $\Spec(R/(P-\{z\})R)$ as $k$-schemes.
\item
For any affine open subset $U'$ of $X'$, the inverse image $\pi^{-1}(U')$ of $U'$ by $\pi$ is an affine open subset of $X$.
\end{enumerate}

Let $D=\Spec(R/\prod_{x\in P}xR)$, and let $D'=\Spec(R'/\prod_{x\in P-\{z\}}xR')$. 
\begin{enumerate}
\setcounter{enumi}{6}
\item If $\sigma'$ is an admissible composition of blowing-ups over $D'$, then $\sigma$ is an admissible composition of blowing-ups over $D$.
\item
If $\sigma'$ is an admissible composition of blowing-ups with center of codimension two over $D'$, then $\sigma$ is an admissible composition of blowing-ups with center of codimension two over $D$.
\end{enumerate}
\end{lemma}

\section{Main results}
\label{mainmain}
We state our main results. Their proofs will be given in Section~\ref{main proof} and Section~\ref{submain proofs}. 

We fix notations we use throughout this section.

Let $k$ be any \emph{algebraically closed} field, let $R$ be any regular local ring such that $R$ contains $k$ as a subring, the residue field $R/M(R)$ is isomorphic to $k$ as $k$-algebras, $R$ is a localization of a finitely generated $k$-algebra and $\dim R\geq 1$, let $P$ be any parameter system of $R$, and let $z\in P$ be any element.

By $R'$ we denote the localization of $k[P-\{z\}]$ by the maximal ideal $k[P-\{z\}]\cap M(R)=(P-\{z\}) k[P-\{z\}]$.

Furthermore, we denote
\begin{equation*}\begin{split}
PW&=\{\phi\in R|\phi\neq 0, \Gamma_+(P,\phi) \text{ is of } z \text{-Weierstrass type.}\},\\
RW&=\{\phi\in PW|\Gamma_+(P,\psi) \text{ has no } z \text{-removable faces, where }\psi\text{ denotes a}\\
&\qquad\qquad \text{main factor of } (P,z,\phi).\},\\
SW&=\{\phi\in RW|\Gamma_+(P,\phi) \text{ is } z \text{-simple.}\}.\\
\end{split}\end{equation*}

Note that $R\supset PW\supset RW\supset SW\neq \emptyset$, if $\phi\in R$, $\phi\neq 0$ and $\Gamma_+(P,\phi)$ is $z$-simple then $\Gamma_+(P,\phi)$ is of $z$-Weierstrass type, and $R=\{0\}\cup PW$ if $\dim R\leq 2$.

For the proof of our main theorem below, we apply the theory of convex sets and the toric theory. By our main theorem any element in $SW$ is reduced to an element in $PW$ with a strictly smaller value of $\Inv$ or $\Inv 2$.

\begin{theorem}
\label{main}
Assume $\dim R\geq 2$. Consider any $\phi\in R$ such that $\phi\neq 0$, $\Gamma_+(P,\phi)$ is $z$-simple, and $\phi$ satisfies one of the following two conditions:
\begin{enumerate}
\item
$\Inv(P,z,\phi)>0$, $\Gamma_+(P, \psi)$ has no $z$-removable faces, where $\psi$ denotes a main factor of $(P,z,\phi)$.
\item
$\Inv(P,z,\phi)=0$, $\Inv 2(P,z,\phi)\geq 2$ and $z$ divides $\phi$.
\end{enumerate}

Let $D=\Spec(R/\prod_{x\in P}x R)$, which is a normal crossing divisor on $\Spec(R)$. We define a coordinate system $\xi_{M(R)}:\Comp(D)\rightarrow R$ of the normal crossing scheme $(\Spec(R),D)$ at $M(R)$ by putting $\xi_{M(R)}(\Spec(R/x R))=x$ for any $x\in P$. The triplet $(\Spec(R),D, \{\xi_{M(R)}\})$ is a coordinated normal crossing scheme over $k$.

There exist a smooth scheme $X$ over $\Spec(R)$ and an admissible composition of blowing-ups $\sigma:X\rightarrow \Spec(R)$ with centers of codimension two over $D$ such that for any closed point $a\in X$ with $\sigma(a)=M(R)$, there exist a closed point $b\in X$ and a component $\bar{C}$ passing through $b$ of the pull-back $\sigma^*D$ of the divisor $D$ by $\sigma$ satisfying the following five conditions:
\begin{enumerate}
\item The number of components of the normal crossing divisor $\sigma^*D$ on $X$ passing through $b$ is equal to $\dim R=\dim X$.
\item
The point $a$ belongs to the complement $U(X,\sigma^*D,b)$ in $X$ of the union of all components of $\sigma^*D$ not passing through $b$. 
\end{enumerate}

We consider the pull-back $(\sigma^*\xi)_b$ at $b$ of the coordinate system $\xi_{M(R)}$ at $M(R)$ by $\sigma$. $(\sigma^*\xi)_b$ is a coordinate system associated with the normal crossing divisor $\sigma^*D\cap U(X,\sigma^*D,b)$ on an affine open set $U(X,\sigma^*D,b)$. $\Comp(\sigma^*D)(b)$ denotes the set of all components of $\sigma^*D$ passing through $b$. Note that for any $C\in \Comp(\sigma^*D)(b)$, $(\sigma^*\xi)_b(C)$ is a regular function over $U(X,\sigma^*D,b)$ and its value $(\sigma^*\xi)_b(C)(a)\in k$ at $a$ is defined. The local ring $\mathcal{O}_{X,a}$ of $X$ at $a$ is a regular local ring containing $k$ as a subring, the residue field $\mathcal{O}_{X,a}/M(\mathcal{O}_{X,a})$ is isomorphic to $k$ as $k$-algebras and $\mathcal{O}_{X,a}$ is a localization of a finitely generated $k$-algebra.
We denote $\bar{P}=\{(\sigma^*\xi)_b(C)- (\sigma^*\xi)_b(C)(a)|C\in \Comp(\sigma^*D)(b)\}$, which is a parameter system of $\mathcal{O}_{X,a}$, and we denote $\bar{z}= (\sigma^*\xi)_b(\bar{C})- (\sigma^*\xi)_b(\bar{C})(a)\in \bar{P}$. We consider the local $k$-algebra homomorphism $\sigma^*:R\rightarrow \mathcal{O}_{X,a}$ induced by $\sigma$.
\begin{enumerate}
\setcounter{enumi}{2}
\item $\sigma^*(\phi)\neq 0$ and the Newton polyhedron $\Gamma_+(\bar{P}, \sigma^*(\phi))$ is of $\bar{z}$-Weierstrass type.
\item If $\Inv(P,z,\phi)>0$, then $\Inv(\bar{P}, \bar{z}, \sigma^*(\phi))<\Inv(P,z,\phi)$.
\item If $\Inv(P,z,\phi)=0$, then $\Inv(\bar{P}, \bar{z}, \sigma^*(\phi))=0$ and $\Inv 2(\bar{P}, \bar{z}, \sigma^*(\phi))<\Inv 2(P,$\hfill\break$z,\phi)$.
\end{enumerate}
\end{theorem}

\begin{remark}
The smooth scheme $X$ and the admissible composition of blowing-ups $\sigma:X\rightarrow \Spec(R)$ in the above theorem are concretely constructed from the Newton polyhedron $\Gamma_+(P,\phi)$ using the toric theory.

Now, since $\Gamma_+(P,\phi)$ is $z$-simple by our assumption, the normal fan  $\Sigma$ of $\Gamma_+(P,\phi)$ has simple structure, and the support $|\Sigma|$ of $\Sigma$ is a regular cone with dimension equal to $\dim R$. Starting from the fan $\mathcal{F}(|\Sigma|)$ consisting of $|\Sigma|$ and its faces and repeating star subdivisions with center in a regular cone of dimension two, we construct most effectively a regular subdivision $\Sigma^*$ of $\Sigma$ satisfying $|\Sigma^*|=|\Sigma|$, which we call an \emph{upward subdivision} of $\Sigma$. We explain how to construct $\Sigma^*$ in Section~\ref{upward}.

Our scheme $X$ and our morphism $\sigma$ are the toric variety over $\Spec(R)$ and the toric morphism associated with an upward subdivision $\Sigma^*$ of the normal fan $\Sigma$ of $\Gamma_+(P,\phi)$. 
\end{remark}

In the theorem below we study properties of $z$-removable faces closely.

\begin{theorem}
\label{erase faces}
Assume $\dim R\geq 2$.

Consider any element $w\in M(R^c)$ with $\partial w/\partial z\in R^{c\times}$.
We denote $P_w=\{w\}\cup(P-\{z\})$. (Lemma~\ref{coordinate change}.)

The bijection $P_w\rightarrow P$ sending $w\in P_w$ to $z\in P$ and sending any $x\in P_w-\{w\}=P-\{z\}$ to $x\in P-\{z\}$ itself induces an isomorphism $\Map(P,\R)\rightarrow\Map(P_w,\R)$ of vector spaces over $\R$.
By this isomorphism we identify $\Map(P,\R)$ and $\Map(P_w,\R)$.

Consider any element $\psi\in R$ such that $\psi\neq 0$, $\Gamma_+(P,\psi)$ is of $z$-Weierstrass type and any $x\in P-\{z\}$ does not divide $\psi$. We take the unique non-negative integer $h$ such that $\{hf^P_z\}$ is the unique $z$-top vertex of $\Gamma_+(P,\psi)$.

Recall that $\Gamma_+(P,\psi)\subset\Map(P,\R)$ and $\{f^P_x|x\in P\}$ is an $\R$-basis of the vector space $\Map(P,\R)$. Let $U=\{a\in \Map(P,\R)|\langle f^{P\vee}_z,a\rangle<h\}$ and $V=\{a\in \Map(P,\R)|\langle f^{P\vee}_z,a\rangle=0\}$. We put $\rho(a)=(a-\langle f^{P\vee}_z,a\rangle f^P_z)/(h-\langle f^{P\vee}_z,a\rangle)\in V$ for any $a\in U$ and we define a mapping $\rho:U\rightarrow V$.

Note that $V$ is an $\R$-vector subspace of $\Map(P,\R)$ with $\dim V=\dim\Map(P,$\break$\R)-1$ and the set $\{f^P_x|x\in P-\{z\}\}$ is an $\R$-basis of $V$. Using  the isomorphism $\Map(P-\{z\},\R)\rightarrow V$ of vector spaces over $\R$ sending $f^{P-\{z\}}_x\in\Map(P-\{z\},\R)$ to $f^P_x\in V$ for any $x\in P-\{z\}$ we identify $\Map(P-\{z\},\R)$ and $V$. 

We identify the dual vector space $V^*$ of $V$ with the vector subspace $\{\bar{\omega}\in\Map(P,$\break$\R)^*|\langle \bar{\omega}, f^P_z\rangle=0\}$ in the dual vector space $\Map(P,\R)^*$ of $\Map(P,\R)$. Under this identification $(\Map(P,\R_0)\cap V)^\vee|V=(\Map(P,\R_0)^\vee|\Map(P,\R))\cap V^*$.
\begin{enumerate}
\item $\rho(\Gamma_+(P_w,\psi)\cap U)$ is either an empty set or a rational pseudo polytope over the lattice $\Map(P,\Z)\cap V$ in $V$ such that $\rho(\Gamma_+(P_w,\psi)\cap U)=\Conv(Y)+(\Map(P,\R_0)\cap V)$ for some non-empty finite subset $Y$ of $\Map(P,\Q)\cap V$. 
\item
There exists uniquely an element $\chi_0\in M(R^{\prime c})$ such that 
$\Gamma_+(P_{z+\chi_0},\psi)$ has no $(z+\chi_0)$-removable faces and $\Supp(P-\{z\},\chi_0)\subset \rho(\Gamma_+(P,\psi)\cap U)- \rho(\Gamma_+(P_{z+\chi_0},\psi)\cap U)$. 
\end{enumerate}

Below, we assume that $\chi_0\in M(R^{\prime c})$, $\Gamma_+(P_{z+\chi_0},\psi)$ has no $(z+\chi_0)$-removable faces and $\Supp(P-\{z\},\chi_0)\subset \rho(\Gamma_+(P,\psi)\cap U)- \rho(\Gamma_+(P_{z+\chi_0},\psi)\cap U)$ and that $v\in R^{c\times}$, $\mu\in M(R^{\prime c})$ and $w=v(z+\mu)$.
\begin{enumerate}\setcounter{enumi}{2}
\item 
$\rho(\Gamma_+(P_w,\psi)\cap U)=\Conv(\rho(\Gamma_+(P_{z+\chi_0},\psi)\cap U)\cup\Gamma_+(P-\{z\},\mu-\chi_0))$.

For any face $F$ of $\Gamma_+(P_w,\psi)$ with $hf^P_z\in F$ and $F\cap U\neq \emptyset$, $F$ is $w$-removable, if and only if, $\rho(F\cap U)\cap\rho(\Gamma_+(P_{z+\chi_0},\psi)\cap U)=\emptyset$.  
\item 
The following three conditions are equivalent:
\begin{enumerate}
\item
$\rho(\Gamma_+(P_w,\psi)\cap U)= \rho(\Gamma_+(P_{z+\chi_0},\psi)\cap U)$.
\item
$\Supp(P-\{z\},\mu-\chi_0)\subset \rho(\Gamma_+(P_{z+\chi_0},\psi)\cap U)$.
\item
$\Gamma_+(P_w,\psi)$ has no $w$-removable faces.
\end{enumerate}
\item $\rho(\Gamma_+(P_{z+\chi_0},\psi)\cap U)=\emptyset$, if and only if, $\psi=u\lambda^h$ for some $u\in R^\times$ and some $\lambda\in M(R)$.
\item $\chi_0\in M(R^{\prime h})$.
\item
Assume moreover, that either $\rho(\Gamma_+(P_{z+\chi_0},\psi)\cap U)$ has at most one vertex, or $\dim R\leq 3$. Then, there exists $w_1\in M(R)$ such that $\partial w_1/\partial z\in R^\times$, and $\Gamma_+(P_{w_1},\psi)$ has no $w_1$-removable faces. 
\end{enumerate}
\end{theorem}

We would like to solve the following problem:

\begin{problem}
\label{goal one}
Show that for any $\phi\in R$ with $\phi\neq 0$, there exists a composition $\sigma:X\rightarrow\Spec(R)$ of finite blowing-ups with center in a closed irreducible smooth subscheme such that the divisor on $X$ defined by the pull-back $\sigma^*(\phi)\in\mathcal{O}_X(X)$ of $\phi$ by $\sigma$ has normal crossings.

\end{problem}

Note here that $\dim R'=\dim R-1<\dim R$. 

We consider the case $\dim R=1$. 

Consider any $\phi\in R$ with $\phi\neq 0$.
$\phi$ has normal crossings over $P$. 

Put $X=\Spec(R) $ and we consider the identity morphism $\sigma:X\rightarrow\Spec(R)=X$. We know that $\sigma$ is a composition of blowing-ups with center in a closed irreducible smooth subscheme, and the divisor  defined by $\sigma^*(\phi)=\phi$ on $X=\Spec(R)$ has normal crossings.
We can easily solve the Problem~\ref{goal one}, if $\dim R=1$.

Therefore, we decide that we use induction on $\dim R$, and we can assume the following claim $(*)$ whenever $\dim R\geq 2$.:

\begin{description}
\item[$(*)$]
For any  $\phi'\in R'$ with $\phi'\neq 0$, there exists a composition $\sigma':X'\rightarrow\Spec(R')$ of finite blowing-ups with center in a closed irreducible smooth subscheme such that the divisor on $X'$ defined by the pull-back $\sigma^{\prime *}(\phi')\in\mathcal{O}_{X'}(X')$ of $\phi'$ by $\sigma'$ has normal crossings.
\end{description}

Claim $(*)$ is true, if $\dim R\leq 2$.

Let $\sigma':X'\rightarrow\Spec(R')$ be any composition of blowing-ups with center in a closed irreducible smooth subscheme. The scheme $X'$ is smooth. We consider a morphism $\Spec(R)\rightarrow\Spec(R')$ induced by the inclusion ring homomorphism $R'\rightarrow R$, the product scheme $X=X'\times_{\Spec(R')}\Spec(R)$, the projection $\sigma:X\rightarrow\Spec(R)$, and the projection $\pi:X\rightarrow X'$. We know the following (Lemma~\ref{pull back blowing-ups}.):
\begin{enumerate}
\item The morphism $\sigma$ is a composition of finite blowing-ups with center in a closed irreducible smooth subscheme. The scheme $X$ is smooth.
\item We consider the prime divisor $\Spec(R/zR)$ on $\Spec(R)$  defined by $z\in R$. The pull-back $\sigma^*\Spec(R/zR)$ of $\Spec(R/zR)$ by $\sigma$ is a smooth prime divisor of $X$, and $\sigma^*\Spec(R/zR)\supset \sigma^{-1}(M(R))$.
\item The projection $\pi:X\rightarrow X'$ induces an isomorphism $\sigma^*\Spec (R/zR)\rightarrow X'$.
\item
For any closed point $a\in X$, any $w\in M(R^h)$ with $\partial w/\partial z\in R^{h\times}$ and any parameter system $Q'$ of the local ring $\mathcal{O}_{X',\pi(a)}$ of $X'$ at $\pi(a)$, $\sigma(a)=M(R)$ and $\{\sigma^*(w)\}\cup\pi^*(Q') $ is a parameter system of the Henselization $\mathcal{O}_{X,a}^h$ of the local ring $\mathcal{O}_{X,a}$ of $X$ at $a$ with $\pi^*(Q')\subset\mathcal{O}_{X,a}$, where $\sigma^*:R^h\rightarrow \mathcal{O}_{X,a}^h$ denotes the homomorphism of local $k$-algebras induced by $\sigma$ on the Henselizations of local rings and $\pi^*:\mathcal{O}_{X',\pi(a)}\rightarrow \mathcal{O}_{X,a}\subset \mathcal{O}_{X,a}^h$ denotes the homomorphism of local $k$-algebras induced by $\pi$.
\end{enumerate}

The theorem below plays three roles. First, any element in $PW$ with a positive value of $\Inv$ is reduced either to an element in $RW$ with the same value of $\Inv$ or to an element in $PW$ with a strictly smaller value of $\Inv$. Second, any element in $RW$ with a positive value of $\Inv$ is reduced either to an element in $SW$ with the same value of $\Inv$, or to an element in $PW$ with a strictly smaller value of $\Inv$. Third, any element in $PW$ with the value zero of $\Inv$ is reduced to an element in $SW$ with the value zero of $\Inv$ and with the same value of $\Inv 2$.

\begin{theorem}
\label{make simple}
Assume the above $(*)$ and $\dim R\geq 2$.

Consider any element $\phi\in R$ and any $w\in M(R^h)$ such that $\phi\neq 0$, $\partial w/\partial z\in R^{h\times}$, and $\Gamma_+(P_w,\phi)$ is of $w$-Weierstrass type, where $P_w=\{w\}\cup(P-\{z\})$. By $\psi$ we denote a main factor of $(P_w,w,\phi)$.

There exists a composition $\sigma':X'\rightarrow\Spec(R')$ of finite blowing-ups with center in a closed irreducible smooth subscheme with the following properties:

We consider the product scheme $X=X'\times_{\Spec(R')}\Spec(R)$, the projection $\sigma:X\rightarrow\Spec(R)$ and the projection $\pi:X\rightarrow X'$. $X$ and $X'$ are smooth. Note that for any closed point $a\in X$ with $\sigma(a)=M(R)$, we have the homomorphism $\sigma^*:R\rightarrow \mathcal{O}_{X,a}$ of local $k$-algebras induced by $\sigma$ from $R$ to the local ring $\mathcal{O}_{X,a}$ of $X$ at $a$, the homomorphism $\pi^*:\mathcal{O}_{X',\pi(a)}\rightarrow \mathcal{O}_{X,a}$ of local $k$-algebras induced by $\pi$ from the local ring $\mathcal{O}_{X',\pi(a)}$ of $X'$ at $\pi(a)$ to $\mathcal{O}_{X,a}$ and the homomorphism $\sigma^{\prime *}:R'\rightarrow \mathcal{O}_{X',\pi(a)}$ of local $k$-algebras induced by $\sigma'$ from $R'$ to $\mathcal{O}_{X',\pi(a)}$, and $\sigma^*$ induces a homomorphism $\sigma^*:R^h\rightarrow \mathcal{O}_{X,a}^h$ of local $k$-algebras from the Henselization $R^h$ of $R$ to the Henselization $\mathcal{O}_{X,a}^h$ of $\mathcal{O}_{X,a}$.

For any closed point $a\in X$ with $\sigma(a)=M(R)$, $\sigma^*(\phi)\neq 0$, and there exists a parameter system $\bar{Q}$ of $\mathcal{O}_{X',\pi(a)}$ satisfying the following six conditions. We denote $\bar{P}_w=\{\sigma^*(w)\}\cup \pi^*(\bar{Q})$ and by $\bar{\psi}$ we denote a main factor of $(\bar{P}_w, \sigma^*(w), \sigma^*(\phi))$:
\begin{enumerate}
\item
$\sigma^{\prime *}(x)$ has normal closings over $\bar{Q}$ for any $x\in P_w-\{w\}=P-\{z\}$
\item
$\Gamma_+(\bar{P}_w,\sigma^*(\phi))$ is $\sigma^*(w)$-simple.
\item
$\Inv(\bar{P}_w,\sigma^*(w),\sigma^*(\phi)) \leq\Inv(P_w,w,\phi)$.
\item
If $\Inv(\bar{P}_w,\sigma^*(w),\sigma^*(\phi))=\Inv(P_w, w, \phi)$ and $\Gamma_+(P_w,\psi)$ has no $w$-removable faces, then 
$\Gamma_+(\bar{P}_w,\bar{\psi})$ has no $\sigma^*(w)$-removable faces.
\item
Assume that $\Inv(\bar{P}_w,\sigma^*(w),\sigma^*(\phi))=\Inv(P_w, w, \phi)$ and $w=z+\chi_0$ where $\chi_0\in M(R^{\prime h})$ is the unique element in Theorem~\ref{erase faces}.2. There exists an element $\bar{w}\in M(\mathcal{O}_{X,a})$ such that $\partial\bar{w}/\partial \sigma^*(w)\in \mathcal{O}_{X,a}^{h\times}$ and if we denote $\bar{P}_{\bar{w}}=\{\bar{w}\}\cup \pi^*(\bar{Q})$, then $\bar{P}_{\bar{w}}$ is a parameter system of $\mathcal{O}_{X,a}$ and $\Gamma_+(\bar{P}_{\bar{w}},\bar{\psi})$ has no $\bar{w}$-removable faces, and $\Inv(\bar{P}_{\bar{w}},\bar{w},\sigma^*(\phi))=\Inv(P_w, w, \phi)$. 
\item
If $\Inv(P_w, w, \phi)=0$, then $\Inv(\bar{P}_w,\sigma^*(w), \sigma^*(\phi))=0$ and  $\Inv 2(\bar{P}_w,\sigma^*(w),$\hfill\break$\sigma^*(\phi))=\Inv 2(P_w, w, \phi)$.

If $w$ divides $\phi$, then $\sigma^*(w)$ divides $\sigma^*(\phi)$.
\end{enumerate}
\end{theorem}

By the lemma below any non-zero element in $R$ is reduced to an element in $PW$.

\begin{lemma}
\label{make Weierstrass type}
Consider any $\phi\in R$ with $\phi\neq 0$. Let $h=\Ord(\phi)\in\Z_0$.
There exists a mapping $\A:P-\{z\}\rightarrow k$ such that $\bar{P}=\{z\}\cup\{x-\A(x)z|x\in P-\{z\}\}$ is a parameter system of $R$ containing $z$, $\Gamma_+(\bar{P},\phi)$ is of $z$-Weierstrass type, the unique $z$-top vertex of $\Gamma_+(\bar{P},\phi)$ is $\{hf^{\bar{P}}_z\}$, and $\Inv(\bar{P},z,\phi)\leq h$. 
\end{lemma}

\begin{cor}[Resolution game]
\label{resolution game}
Consider a mathematical game with two players A and B.
At the start of the game a pair $(R,\phi)$ of any regular local ring $R$ with $\dim R\geq 1$ such that $R$ contains $k$ as a subring, the residue field $R/M(R)$ is isomorphic to $k$ as $k$-algebras and $R$ is a localization of a finitely generated $k$-algebra, and any non-zero element $\phi\in R$ is given. We play our game repeating the following step. Before the first step we put $(S,\psi)=(R,\phi)$: At the start of each step, player A chooses a composition $\sigma:X\rightarrow\Spec(S)$ of finite blowing-ups with center in a closed irreducible smooth subscheme. Then, player B chooses a closed point $a\in X$ with $\sigma(a)=M(S)$. We have a morphism $\sigma^*:S\rightarrow\mathcal{O}_{X,a}$ of local $k$-algebras induced by $\sigma$. If the element $\sigma^*(\psi)\in \mathcal{O}_{X,a}$ has normal crossings, then the palyer A wins. Otherwise we proceed to the next step after replacing the pair $(S,\psi)$ by the pair $(\mathcal{O}_{X,a}, \sigma^*(\psi))$.

At this game, player A can always win the game after finite steps for any $R$ and any non-zero element $\phi\in R$, even if the characteristic of the ground field $k$ is positive. 
\end{cor}

\begin{remark}
Note that the pair $(S, \psi)$ satisfies the same conditions as $(R,\phi)$ throughout the game.

A similar game can be found in Spivakovsky~\cite{S82}.

By valuation theory we know that the above Corollary implies ``the local uniformization theorem in arbitrary characteristic and in arbitrary dimension''.
(Zariski et al.~\cite{ZS}, Zariski~\cite{Z40}, Abhyankar~\cite{A56}.)
\end{remark}

\begin{cor}
\label{local uniformization}
\emph{(The local uniformization theorem in arbitrary characteristic and in arbitrary dimension)}
Given any field $\Sigma$ such that $\Sigma$ contains $k$ as a subfield and $\Sigma$ is finitely generated over $k$, given any projective model $X_0$ of $\Sigma$ and given any valuation $B$ of dimension zero of $\Sigma$ containing $k$ with center $a_0$ on $X_0$, there exists a projective model $X$ of $\Sigma$ on which the center of $B$ is at a smooth point $a$ of $X$ such that the inclusion relation $\mathcal{O}_{X,a}\supset \mathcal{O}_{X_0,a_0}$ of local rings holds.
\end{cor}

\section{Basic theory of convex sets}
\label{btcs}

In this section we begin the study of convex sets to develop the toric theory. The theory of convex sets will be applied to the proof of our main theorem, Therem~\ref{main} in Section~\ref{main proof}.

Let $V$ be any vector space of finite dimension over $\R$. 

In Section~\ref{concept} we defined eight mappings
$$\Conv,\Affi,\Cone,\Convcone,\Vect,\QVect,\Clos,\Stab:2^V\rightarrow 2^V.$$

\begin{lemma}
\label{a}
Let $X$ and $Y$ be any subsets of $V$.
\begin{enumerate}
\item
\begin{equation*}
\begin{split}
\Conv(\emptyset)&=\Affi(\emptyset)=\Clos(\emptyset)=\emptyset,\\
\Cone(\emptyset)&=\Convcone(\emptyset)=\Vect(\emptyset)=\QVect(\emptyset)=\{0\},\\
\Stab(\emptyset)&=V.\\
\end{split}
\end{equation*}
\item
\begin{equation*}
\begin{split}
\Conv(X)&=\{a\in V|a=\sum_{x\in\Supp(\lambda)}\lambda(x)x
\text{ for some }\lambda\in\Map'(X,\R_0)\text{ with}\\
&\qquad\qquad\sum_{x\in\Supp(\lambda)}\lambda(x)=1\},\\
\Affi(X)= &\{a\in V|a=\sum_{x\in\Supp(\lambda)}\lambda(x)x
\text{ for some }\lambda\in\Map'(X,\R)\text{ with}\\
&\qquad\qquad\sum_{x\in\Supp(\lambda)}\lambda(x)=1\},
\end{split}\end{equation*}
\begin{equation*}
\Cone(X)=
\begin{cases}
\{a\in V|a=\lambda x\text{ for some }\lambda\in\R_0\text{ and some }x\in X\}&
\text{ if $X\neq\emptyset$},\\
\{0\}&\text{ if $X=\emptyset$},
\end{cases}\end{equation*}
\begin{equation*}\begin{split}
\Convcone(X)&=\{a\in V|a=\sum_{x\in\Supp(\lambda)}\lambda(x)x
\text{ for some }\lambda\in\Map'(X,\R_0)\},\\
\Vect(X)&=\{a\in V|a=\sum_{x\in\Supp(\lambda)}\lambda(x)x
\text{ for some }\lambda\in\Map'(X,\R)\},\\
\QVect(X)&=\{a\in V|a=\sum_{x\in\Supp(\lambda)}\lambda(x)x
\text{ for some }\lambda\in\Map'(X,\Q)\}.\\
\end{split}
\end{equation*}
\item
If $X$ is a finite set, then we have
$\Convcone(X)=\sum_{x\in X}\R_0 x$, $\Vect(X)=\sum_{x\in X}\R x$, and $\QVect(X)=\sum_{x\in X}\Q x$.
\item
For any vector space $W$ of finite dimension over $\R$ and any homomorphism $\pi:V\rightarrow W$ of vector spaces over $\R$, we have
$\pi(\Conv(X))=\Conv(\pi(X))$,
$\pi(\Affi(X))=\Affi(\pi(X))$,
$\pi(\Cone(X))=\Cone(\pi(X))$,
$\pi(\Convcone(X))=$\hfill\break$\Convcone(\pi(X))$,
$\pi(\Vect(X))=\Vect(\pi(X))$, and
$\pi(\QVect(X))=\QVect(\pi(X))$.
\item
For any $a\in V$, we have
$\Conv(X+\{a\})=\Conv(X)+\{a\}$,
$\Affi(X+\{a\})=\Affi(X)+\{a\}$.
\item
\begin{equation*}
\begin{gathered}
X\subset\Conv(X)\subset\Affi(X)\subset\Vect(X),\\
X\cup\{0\}\subset\Cone(X)\subset\Convcone(X)\subset\Vect(X),\\
\Conv(X)\subset\Convcone(X),\\
X\cup\{0\}\subset\QVect(X)\subset\Vect(X),\\
X\subset\Clos(X).\\
\end{gathered}
\end{equation*}
\item
If $X\subset Y$, then $\Conv(X)\subset\Conv(Y)$, $\Affi(X)\subset\Affi(Y)$, $\Cone(X)\subset\Cone(Y)$, $\Convcone(X)\subset\Convcone(Y)$, $\Vect(X)\subset\Vect(Y)$, $\QVect(X)\subset\QVect(Y)$, and $\Clos(X)\subset\Clos(Y)$.
\item
$\Conv(\Conv(X))=\Conv(X)$, $\Affi(\Affi(X))=\Affi(X)$, $\Cone(\Cone(X))=\Cone(X)$, $\Convcone(\Convcone(X))=\Convcone(X)$,
$\Vect(\Vect(X))=\Vect(X)$,\hfill\break$\QVect(\QVect(X))=\QVect(X)$, and $\Clos(\Clos(X))=\Clos(X)$.
\item
\begin{equation*}
\begin{split}
\Vect(\Conv(X))&=\Conv(\Vect(X))=\Vect(X),\\
\Vect(\Affi(X))&=\Affi(\Vect(X))=\Vect(X),\\
\Vect(\Cone(X))&=\Cone(\Vect(X))=\Vect(X),\\
\Vect(\Convcone(X))&=\Convcone(\Vect(X))=\Vect(X).\\
\end{split}
\end{equation*}
\item
\begin{equation*}
\begin{split}
\Convcone(\Conv(X))&=\Conv(\Convcone(X))=\Convcone(X),\\
\Convcone(\Cone(X))&=\Cone(\Convcone(X))=\Convcone(X).\\
\Cone(\Conv(X))&=\Conv(\Cone(X))=\Convcone(X),\\
\Convcone(\Affi(X))=\Cone(\Affi(X))\subset&\Affi(\Convcone(X))= \Affi(\Cone(X))=\Vect(X),\\
\end{split}
\end{equation*}
\item
$\Affi(\Conv(X))=\Conv(\Affi(X))=\Affi(X)$.
\item
\begin{equation*}
\begin{split}
\Clos(\Affi(X))&=\Affi(\Clos(X))=\Affi(X),\\
\Clos(\QVect(X))&=\Clos(\Vect(X))=\Vect(\Clos(X))=\Vect(X),\\
\Clos(\Conv(X))&=\Conv(\Clos(\Conv(X))),\\
\Clos(\Cone(X))&=\Cone(\Clos(\Cone(X))),\\
\Clos(\Convcone(X))&=\Convcone(\Clos(\Convcone(X))).\\
\end{split}
\end{equation*}
\end{enumerate}
\end{lemma}

\begin{lemma}
\label{b}
Let $X$ and $Y$ be any subsets of $V$.
\begin{enumerate}
\item
$\Conv(X)+\Conv(Y)=\Conv(X+Y)$.
\item
$\Affi(X)+\Affi(Y)=\Affi(X+Y)$.
\item
$\Convcone(X)+\Convcone(Y)=\Convcone(X\cup Y)$.
\item
$\Vect(X)+\Vect(Y)=\Vect(X\cup Y)$.
\item 
For any vector space $W$ of finite dimension over $\R$ and any homomorphism $\pi:V\rightarrow W$ of vector spaces over $\R$, we have
$\pi(X)+\pi(Y)=\pi(X+Y)$.
\end{enumerate}
\end{lemma}

In Section~\ref{concept} we defined concepts of segments, lines, convex sets, affine spaces, cones, convex cones, vector spaces and vector spaces over $\Q$.

\begin{lemma}
\label{c}
Let $S$ be any non-empty subset of $V$.
\begin{enumerate}
\item
The following four conditions are equivalent;
\begin{enumerate}
\item $S$ is convex.
\item For any $t\in\R_0$ and any $u\in\R_0$, $tS+uS=(t+u)S$.
\item $S\supset\Conv(S)$.
\item $S=\Conv(X)$ for some non-empty subset $X$ of $V$.
\end{enumerate}
\item
If $S$ is convex, then $\Clos(S)$ is also convex, and $\Affi(S)=\Affi(\Clos(S))$.
\item
The following six conditions are equivalent;
\begin{enumerate}
\item $S$ is an affine space.
\item For any $t\in\R$ and any $u\in\R$ with $t+u\neq 0$, $tS+uS=(t+u)S$.
\item $S\supset\Affi(S)$.
\item $S=\Affi(X)$ for some non-empty subset $X$ of $V$.
\item $S\neq\emptyset$ and $S+\{-a\}$ is a vector space for some $a\in S$.
\item $S\neq\emptyset$ and $S+\{-a\}$ is a vector space for any $a\in S$.
\end{enumerate}
\item
The following three conditions are equivalent;
\begin{enumerate}
\item $S$ is an affine space containing $0$.
\item $S$ is an affine space with $S=\Stab(S)$.
\item $S$ is a vector space.
\end{enumerate}
\item
Assume that $S$ is an affine space. Then, $\Stab(S)$ is a vector space, and we have $S=\Stab(S)+\{a\}$ and $\Stab(S)=S+\{-a\}$ for any $a\in S$.
\item
Any affine space is closed and convex.
\item
$0\in\Stab(S)\subset\Stab(\Affi(S))\subset\Vect(S)$. $\Stab(S)+\Stab(S)=\Stab(S)$.
\item
The following three conditions are equivalent;
\begin{enumerate}
\item $S$ is a cone.
\item $S\supset\Cone(S)$.
\item $S=\Cone(X)$ for some subset $X$ of $V$.
\end{enumerate}
\item Any cone contains $0$.
\item
If $S$ is a cone, then $\Clos(S)$ is also a cone, and $\Vect(S)=\Vect(\Clos(S))$.
\item
The following four conditions are equivalent;
\begin{enumerate}
\item $S$ is a convex cone.
\item $S$ is convex and $S$ is a cone.
\item $S\supset\Convcone(S)$.
\item $S=\Convcone(X)$ for some subset $X$ of $V$.
\end{enumerate}
\item
Any convex cone contains $0$.
\item
If $S$ is a convex cone, then $\Clos(S)$ is also a convex cone, and $\Vect(S)=\Vect(\Clos(S))$.
\item
If $S$ is a convex cone, then $S\cap(-S)$ is the maximal vector space contained in $S$ with respect to the inclusion relation.
\item
The following three conditions are equivalent;
\begin{enumerate}
\item $S$ is a vector space.
\item $S\supset\Vect(S)$.
\item $S=\Vect(X)$ for some subset $X$ of $V$.
\end{enumerate}
\item
Any vector space contains $0$.
\item
Any vector space is closed, it is an affine space containing $0$, and it is a convex cone.
\item
The following three conditions are equivalent;
\begin{enumerate}
\item $S$ is a vector space over $\Q$.
\item $S\supset\QVect(S)$.
\item $S=\QVect(X)$ for some subset $X$ of $V$.
\end{enumerate}
\item
Any vector space over $\Q$ contains $0$.
\item
The following three conditions are equivalent;
\begin{enumerate}
\item $S$ is closed.
\item $S\supset\Clos(S)$.
\item $S=\Clos(X)$ for some subset $X$ of $V$.
\end{enumerate}
\end{enumerate}
\end{lemma}

\begin{lemma}
\label{d}
Let $S$ and $T$ be any subsets of $V$.
\begin{enumerate}
\item
If $S$ and $T$ are convex, then $S+T$ is convex. If $S$ and $T$ are convex and $S\cap T\neq\emptyset$, then $S\cap T$ is convex.
\item
If $S$ and $T$ are affine spaces, then $S+T$ is an affine space. If $S$ and $T$ are affine spaces and $S\cap T\neq\emptyset$, then $S\cap T$ is an affine space.
\item
If $S$ and $T$ are cones, then $S+T$ and $S\cap T$ are cones.
\item
If $S$ and $T$ are convex cones, then $S+T$ and $S\cap T$ are convex cones.
\item
If $S$ and $T$ are vector spaces, then $S+T$ and $S\cap T$ are vector spaces.
\item
If $S$ and $T$ are vector spaces over $\Q$, then $S+T$ and $S\cap T$ are vector spaces over $\Q$.
\end{enumerate}
\end{lemma}

For any affine space $S$ of $V$, the dimension $\dim S$ of $S$ is defined. It satisfies $\dim S\in\Z_0$, $0\leq\dim S\leq \dim V$ and $\dim S=\dim \Stab(S)$.

\begin{definition}
\label{dimension}
Let $S$ be any convex subset of $V$.
We define
$$\dim S=\dim\Affi(S)\in\Z_0,$$
and we call $\dim S$ the \emph{dimension} of $S$.

We define
\begin{equation*}
\begin{split}
\partial S&=S\cap\Clos(\Affi(S)- S),\\
S^\circ&=S-\Clos(\Affi(S)- S),
\end{split}
\end{equation*}
we call $\partial S$ the \emph{boundary} of $S$, and we call $S^\circ$ the \emph{interior} of $S$.
\end{definition}

\begin{lemma}
\label{f}
\begin{enumerate}
\item
Let $S$ be a convex subset of $V$. $\dim S\in\Z_0$, and $0\leq\dim S\leq\dim V$.
If $S$ is an affine space, then the dimension of $S$ as a convex set and the dimension of $S$ as an affine space are equal.
\item
For any convex set $S$ of $V$, we have $\dim S=\dim\Affi(S)=\dim \Clos(S)$.
\item
Let $S$ and $T$ be convex subsets of $V$ with $S\subset T$. We have $\dim S\leq\dim T$.
\item
Let $S$ be a convex subset of $V$.
$\partial S\cup S^\circ=S$.
$\partial S\cap S^\circ=\emptyset$.
$S^\circ$ is a non-empty open subset of $\Affi(S)$.
If $S$ is closed, then $\partial S$ is also closed.
\item
For any convex cone $S$ of $V$, we have $\Affi(S)=\Vect(S)$, $\dim S=\dim\Vect(S)$,
$\partial S=S\cap\Clos(\Vect(S)- S)$, $S^\circ=S-\Clos(\Vect(S)- S)$ and $S^\circ$ is a non-empty open subset of $\Vect(S)$.
\end{enumerate}
\end{lemma}

\begin{remark}
Consider any convex subset $S$ and $T$ of $V$ with $S\subset T$. We have $\dim S\leq\dim T$. If $S$ and $T$ are affine spaces and $\dim S=\dim T$, then we have $S=T$. However, in general, it does not follow $S=T$ from the assumption $\dim S=\dim T$.
\end{remark}

In Section~\ref{concept} we defined concepts of convex polytopes, convex polyhedral cones, convex pseudo polytopes and simplicial cones.

\begin{lemma}
\label{g}
\begin{enumerate}
\item
Any convex polytope in $V$ is convex, compact and closed.
\item
Any convex polyhedral cone in $V$ is a closed convex cone.
\item
Any convex pseudo polytope in $V$ is convex and closed.
\item
Any vector space in $V$ is a convex polyhedral cone. Any simplicial cone in $V$ is a convex polyhedral cone.
\item
Any affine space in $V$ is a convex pseudo polytope.
Any convex polyhedral cone in $V$ is a convex pseudo polytope.
Any convex polytope in $V$ is a convex pseudo polytope. Any compact convex pseudo polytope in $V$ is a convex polytope.
\end{enumerate}
\end{lemma}

In Section~\ref{concept} we defined concepts of lattices, dual lattices and regular cones.

By definition we know that there exists a lattice $N$ of $V$.
Let $N$ be any lattice of $V$.

\begin{definition}
Let $S$ be any subset of $V$.
\begin{enumerate}
\item
We say that $S$ is a \emph{rational convex polyhedral cone} over $N$, or a convex polyhedral cone $S$ is \emph{rational} over $N$, if there exists a finite subset $X$ of $\QVect(N)$ with $S=\Convcone(X)$.
\item
We say that $S$ is a \emph{rational convex pseudo polytope} over $N$, or a convex pseudo polytope $S$ is \emph{rational} over $N$, if there exist finite subsets $X, Y$ of $\QVect(N)$ with $S=\Conv(X)+\Convcone(Y)$ and $X\neq\emptyset$.
\end{enumerate}
\end{definition}

The dual lattice $N^*$ is defined. We have
$$N^*=\{\omega\in V^*|\langle\omega, a\rangle\in\Z \text{ for any } a\in N\}\subset V^*,$$
by definition.

\begin{lemma}
\label{h}
\begin{enumerate}
\item
$N$ is a submodule of $V$, $N$ is a free module of finite rank over $\Z$ with $\Rank N=\dim V$, and $\Vect(N)=V$.
\item
Any $\Z$-basis of $N$ is a $\Q$-basis of $\QVect(N)$, and is a $\R$-basis of $V$.
\item
For any non-empty finite subset $X$ of $N$, the following three conditions are equivalent:
\begin{enumerate}
\item
$X$ is linearly independent over $\Z$.
\item
$X$ is linearly independent over $\Q$.
\item
$X$ is linearly independent over $\R$.
\end{enumerate}
\item
The dual lattice $N^*$ of $N$ is a lattice of the dual vector space $V^*$ of $V$.
The dual lattice $N^{**}$ of $N^*$ is equal to $N$.
\item
A convex polyhedral cone $S$ in $V$ is rational over $N$, if and only if, $S=\Convcone(X)$ for some finite subset $X$ of $N$.
\item
For any vector space $S$ in $V$ the following three conditions are equivalent:
\begin{enumerate}
\item
$S$ is a rational polyhedral cone over $N$.
\item
$S=\Vect(X)$ for some finite subset $X$ of $N$.
\item
$N\cap S$ is a lattice of $S$.
\end{enumerate}
\item
For any rational polyhedral cone $S$ over $N$ in $V$, $\Vect(S)$ is a rational vector space over $N$ in $V$.
\item
Any regular cone over $N$ in $V$ is a rational simplicial cone over $N$.
\item
A convex pseudo polytope $S$ in $V$ is rational over $N$, if and only if, $S=\Conv (X)+\Convcone(Y)$ for some non-empty finite subset $X$ of $\QVect(N)$ and some finite subset $Y$ of $N$.
\item
For any affine space $S$ in $V$, the following two conditions are equivalent:
\begin{enumerate}
\item
$S$ is a rational convex pseudo polytope over $N$.
\item
$S=\{x\}+\Vect(Y)$ for some $x\in\QVect(N)$ and some finite subset $Y$ of $N$.
\end{enumerate}
\item
For any rational affine space $S$ over $N$ in $V$, $\Stab(S)$ is a rational vector space over $N$ in $V$.
\item
For any rational convex pseudo polytope $S$ over $N$ in $V$, $\Affi(S)$ is a rational affine space over $N$ in $V$.
\end{enumerate}
We consider any vector space $W$ of finite dimension over $\R$ and any homomorphism $\pi:V\rightarrow W$ of vector spaces over $\R$.
The dual homomorphism $\pi^*:W^*\rightarrow V^*$ is defined, and is a homomorphism of vector spaces over $\R$.
The kernel $\pi^{-1}(0)$ of $\pi$ is a vector subspace of $V$, the image $\pi(V)$ is a vector subspace of $W$, the image $\pi^*(W^*)$ of $\pi^*$ is a vector subspace of $V^*$, and the kernel $\pi^{*-1}(0)$ of $\pi^*$ is a vector subspace of $W^*$.
\begin{enumerate}
\setcounter{enumi}{12}
\item
The following seven conditions are equivalent;
\begin{enumerate}
\item
$\pi^{-1}(0)$ is rational over $N$.
\item
$N\cap\pi^{-1}(0)$ is a lattice of $\pi^{-1}(0)$.
\item
$\pi(N)$ is a lattice of $\pi(W)$.
\item
There exists a lattice $Q$ of $W$ satisfying $Q\cap \pi(V)=\pi(N)$.
\item
$\pi^*(W^*)$ is rational over $N^*$.
\item
$N^*\cap\pi^*(W^*)$ is a lattice of $\pi^*(W^*)$.
\item
There exists a lattice $\bar{Q}$ of $W^*$ satisfying $\pi^*(\bar{Q})=N^*\cap\pi^*(W^*)$.
\end{enumerate}
\item
Assume that equivalent seven conditions in the above $13$ hold.
For any lattice $\bar{Q}$ of $W^*$ satisfying $\pi^*(\bar{Q})=N^*\cap\pi^*(W^*)$,
$\bar{Q}\cap \pi^{*-1}(0)$ is a lattice of $\pi^{*-1}(0)$.
\end{enumerate}
\end{lemma}

\section{Convex cones and convex polyhedral cones}
\label{cones}

We study convex cones and convex polyhedral cones.

Let $V$ be any vector space of finite dimension over $\R$, and let $N$ be any lattice of $V$.

\begin{definition}
\label{strongly convex}
Let $S$ be any cone in $V$. We say that $S$ is \emph{strongly convex}, if $S$ is convex and $S\cap(-S)=\{0\}$.
\end{definition}

In Section~\ref{concept} we defined the dual cone $S^\vee|V$ of any convex cone $S$ in $V$.
By definition
$$S^\vee|V=\{\omega\in V^*|\langle\omega, a\rangle\geq 0 \text{ for any } a\in S\}\subset V^*,$$
for any convex cone $S$ in $V$.

\begin{lemma}
\label{i}
Let $S$ be any convex cone in $V$.
\begin{enumerate}
\item
The dual cone $S^\vee|V$ of $S$ is a closed convex cone in the dual vector space $V^*$ of $V$.
\item
Let $W$ be any vector subspace in $V$ with $S\subset W$.
$S$ is a convex cone in $W$ and the dual cone $S^\vee|W$ of $S$ in $W^*$ is defined.

Let $\iota:W\rightarrow V$ denote the inclusion homomorphism.
The dual homomorphism $\iota^*:V^*\rightarrow W^*$ is defined, which is surjective.

$S^\vee|V=\iota^{*-1}( S^\vee|W)$.
\item
$S^\vee|V^\vee|V^*=\Clos(S)$.
\item
$S^\vee|V^\vee|V^*=S$, if and only if, $S$ is closed.
\end{enumerate}
\end{lemma}

When we need not refer to $V$, we also write simply $S^\vee$, instead of $ S^\vee|V$.
In Section~\ref{concept} we defined the concepts of simplicial cones and regular cones.

\begin{lemma}
\label{j}
Let $S$ and $T$ be any convex cones in $V$.
\begin{enumerate}
\item
If $S\subset T$, then $S^\vee\supset T^\vee$.
\item
Assume that $S$ and $T$ are closed. $S\subset T$, if and only if, $S^\vee\supset T^\vee$.
\item
$(S+T)^\vee=S^\vee\cap T^\vee$.
\item
Assume that $S$ and $T$ are closed.
$(S\cap T)^\vee=\Clos(S^\vee+T^\vee)$.
\item
Assume that $S$ is a vector space.
By $\iota:S\rightarrow V$ we denote the inclusion homomorphism.
$\iota^*:V^*\rightarrow S^*$.
$$S^\vee =\{\omega\in V^*|\langle\omega, a\rangle=0 \text{ for any } a\in S\}=\iota^{*-1}(0)\subset V^*,$$
$S^\vee$ is also a vector space, and $\dim S+\dim S^\vee =\dim V$.
\item
$\Vect(S)^\vee=S^\vee\cap (-S^\vee)$.
$\Vect(S)^\vee$ is the maximal vector space contained in $S^\vee$ with respect to the inclusion relation.
\item
Assume that $S$ is closed.
$\Vect(S^\vee)=(S\cap(-S))^\vee$ and $\dim S^\vee =\dim V-\dim(S\cap(-S))$.
$\dim S^\vee =\dim V$, if and only if, $S$ is strongly convex.
\item
$\dim S+\dim S^\vee\geq \dim V$.  $\dim S+\dim S^\vee=\dim V$, if and only if, $S$ is a vector space.
\item
We denote $n=\dim V\in\Z_0$. Let $e: \{1,2,\ldots,n\}\rightarrow V$ be any mapping such that the image $e(\{1,2,\ldots,n\})$ of $e$ is a basis of $V$, and let $I, J, K, L$ be any subset of $\{1,2,\ldots,n\}$ such that $I\cup J\cup K\cup L=\{1,2,\ldots,n\}$ and the intersection of any two of $I, J, K, L$ is empty.
We denote the dual basis of  $\{e(i)|i\in\{1,2,\ldots,n\}\}$ by $\{e^\vee(i)|i\in\{1,2,\ldots,n\}\}$.
We assume
$$\langle e^\vee(i), e(j)\rangle=
\begin{cases}
1&\text{if $i=j$},\\
0&\text{if $i\neq j$},
\end{cases}$$
for any $i\in\{1,2,\ldots,n\}$ and any $j\in\{1,2,\ldots,n\}$.

If
$$S=\sum_{i\in I}\{0\}e(i)+\sum_{j\in J}\R e(j)+\sum_{k\in K}\R_0 e(k)+\sum_{\ell\in L}\R_0(-e(\ell)),$$
then
$$S^\vee=\sum_{i\in I}\R e^\vee (i)+\sum_{j\in J}\{0\}
e^\vee (j)+\sum_{k\in K}\R_0 e^\vee (k)+\sum_{\ell\in L}\R_0(-e^\vee (\ell)).$$
\item
If $S$ is a simplicial cone with $\dim S=\dim V$, then $S^\vee$ is a simplicial cone with $\dim S^\vee=\dim V^*$.
\item If $S$ is a regular cone over $N$ with $\dim S=\dim V$, then $S^\vee$ is a regular cone over $N^*$ with $\dim S^\vee=\dim V^*$.
\end{enumerate}
\end{lemma}

\begin{remark}
Assume $V=\R^3$, $S=\{(x, y, z)\in V|x\geq 0, y\geq 0, z\geq -2\sqrt{xy}\}$ and $T=\{(x, y, z)\in V|x=0, z\geq 0\}$. $S$ and $T$ are closed convex cones in the vector space $V$ with $\dim V=3$. $(0,0,-1)\in\Clos(S+T)$ and $(0,0,-1)\not\in S+T$.
Therefore $S+T$ is not closed.
\end{remark}

\begin{lemma}
\label{one dimensional regular cone}
For any non-empty subset $S$ of $V$, the following two conditions are equivalent:
\begin{enumerate}
\item
$S$ is a regular cone over $N$ with $\dim S=1$.
\item
$S$ is a rational strongly convex polyhedral cone over $N$ with $\dim S=1$.
\end{enumerate}
\end{lemma}

\begin{definition}
\label{faces of cones}
Let $S$ be any convex polyhedral cone in $V$.
We consider the dual cone $S^\vee=S^\vee|V\subset V^*$ of $S$.
\begin{enumerate}
\item
For any $\omega\in S^\vee$, we denote
$$\Delta(\omega,S|V)=\{x\in S|\langle\omega, x\rangle=0\}\subset S.$$

When we need not refer to $V$ or to the pair $(S, V)$, we also write simply  $\Delta(\omega, S)$ or $\Delta(\omega)$, instead of $\Delta(\omega, S|V)$.
\item
Let $F$ be any subset of $S$.
We say that $F$ is a \emph{face} of $S$, if $F=\Delta(\omega, S|V)$ for some $\omega\in S^\vee.$

It is easy to see that any face $F$ of $S$ is a closed convex cone, and the dimension $\dim F\in\Z_0$ of $F$, the boundary $\partial F$ of $F$ and the interior $F^\circ$ of $F$ are defined.

Any face $F$ of $S$ with $\dim F=0$ is called a \emph{vertex} of $S$. Any vertex of $S$ is a subset of $S$ with only one element. Any face $F$ of $S$ with $\dim F=1$ is called an \emph{edge} of $S$. Any face $F$ of $S$ with $\dim F=\dim S-1$ is called a \emph{facet} of $S$. Any face $F$ of $S$ with $F\neq S$ is called a \emph{proper} face of $S$.
The subset $S\cap(-S)$ of $S$ is called the \emph{minimal face} of $S$.
\item
By $\mathcal{F}(S)$ we denote the set of all faces of $S$.

For any $i\in\Z$, the set of all faces $F$ with $\dim F=i$ is denoted by $\mathcal{F} (S)_i$, and the set of all faces $F$ with $\dim F=\dim S-i$ is denoted by $\mathcal{F}(S)^i$.
\item
Let $F$ be any face of $S$.
We denote
\begin{equation*}
\begin{split}
\Delta^\circ(F,S|V)&=\{\omega\in S^\vee|F=\Delta(\omega, S|V)\}\subset S^\vee\subset V^*,\\
\Delta(F,S|V)&=\{\omega\in S^\vee|F\subset\Delta(\omega, S|V)\}\subset S^\vee\subset V^*.\\
\end{split}
\end{equation*}

We call $\Delta^\circ(F,S|V)$ the \emph{open normal cone} of $F$, and we call $\Delta (F,S|V)$ the \emph{normal cone} of $F$.

When we need not refer to $V$ or to the pair $(S, V)$, we also write simply  $\Delta^\circ (F, S)$ or $\Delta^\circ(F)$, $\Delta(F, S)$ or $\Delta(F)$ respectively, instead of $\Delta^\circ(F, S|V)$, $\Delta(F, S|V)$.
\end{enumerate}
\end{definition}

\begin{theorem}
\label{property of cones}
Let $S$ be any convex polyhedral cone in $V$, and let $X$ be any finite subset of $V$ with $S=\Convcone(X)$.
We consider the dual cone $S^\vee=S^\vee|V\subset V^*$ of $S$.
For simplicity we denote
$s=\dim S\in\Z_0$, $L=S\cap(-S)\subset S$, $\ell=\dim L\in\Z_0$, $M=S^\vee\cap(-S^\vee)\subset S^\vee$.
\begin{enumerate}
\item
Consider any vector space $U$ of finite dimension over $\R$ with $\dim S\leq\dim U\leq\dim V$, any injective homomorphism $\nu:U\rightarrow V$ of vector spaces over $\R$ such that $S\subset\nu(U)$, and any subset $F$ of $S$.
The inverse image $\nu^{-1}(S)$ is a convex polyhedral cone in $U$.
The set $F$ is a face of $S$, if and only if, $\nu^{-1}(F)$ is a face of $\nu^{-1}(S)$.
\item
Consider any vector space $W$ of finite dimension over $\R$ with $\dim V\leq\dim W$, any injective homomorphism $\pi:V\rightarrow W$ of vector spaces over $\R$, and any subset $F$ of $S$.
The image $\pi(S)$ is a convex polyhedral cone in $W$.
The set $F$ is a face of $S$, if and only if, $\pi(F)$ is a face of $\pi(S)$.
\item
$\ell\leq s$. $\ell=s\Leftrightarrow L=S\Leftrightarrow S=\Vect(S)$.
\item
Let $F$ be any face of $S$.
\begin{enumerate}
\item
$F=\Convcone(X\cap F)=S\cap\Vect(F)$. $\Vect(F)=\Vect(X\cap F)$.
\item
$F$ is a convex polyhedral cone in $V$.
\item
If $S$ is rational over $N$, then $F$ is also rational over $N$.
If $S$ is a simplicial cone, then $F$ is also a simplicial cone.
If $S$ is a regular cone over $N$, then $F$ is also a regular cone over $N$.
\item
$L=F\cap(-F)\subset F$. $\ell\leq\dim F\leq s$.
\item
Let $G$ be any face of $S$ with $G\subset F$. We have $\dim G\leq\dim F$.  $\dim G=\dim F$, if and only if, $G=F$.
\item
Let $G$ be any subset of $F$. $G$ is a face of the convex polyhedral cone $F$, if and only if, $G$ is a face of $S$ with $G\subset F$.
\item If $x\in S$, $y\in S$ and $x+y\in F$, then $x\in F$ and $y\in F$.
\end{enumerate}
\item
$\mathcal{F}(S)$ is a finite set.
$S\in\mathcal{F}(S)_s$ and $\mathcal{F}(S)_s=\{S\}$.
$S$ contains any face of $S$.
$L\in\mathcal{F}(S)_\ell$ and $\mathcal{F}(S)_\ell=\{L\}$.
$L$ is contained in any face of $S$.
$L=\Convcone(X\cap L)=\Vect(X\cap L)$. 
For any $i\in\Z_0$, $\mathcal{F}(S)_i\neq\emptyset$ if and only if $\ell\leq i\leq s$.
\item
Let $F$ and $G$ be any face of $S$ with $F\subset G$.
We denote $f=\dim F$ and $g=\dim G$.
$\ell\leq f\leq g\leq s$.
There exist $(s-\ell+1)$ of faces $F(\ell), F(\ell+1),\ldots, F(s)$ satisfying the following three conditions:
\begin{enumerate}
\item
For any $i\in\{\ell,\ell+1,\ldots,s-1\}$, $F(i)\subset F(i+1)$.
\item
For any $i\in\{\ell,\ell+1,\ldots,s\}$, $\dim F(i)=i$.
\item
$F(\ell)=L, F(f)=F, F(g)=G, F(s)=S$.
\end{enumerate}
\item
Let $F$ be any face of $S$.
\begin{enumerate}
\item
$F=\partial F\cup F^\circ$. $\partial F\cap F^\circ=\emptyset$.
\item
$F^\circ=F\Leftrightarrow \partial F=\emptyset\Leftrightarrow F=L$.
\item
$$\partial F=\bigcup_{G\in\mathcal{F}(F)-\{F\}}G.$$
\item
$F^\circ $ is a non-empty open subset of $\Vect(F)$.
For any $a\in F^\circ$ and for any $b\in F$, $\Conv(\{a,b\})-\{b\}\subset F^\circ$.
$F^\circ$ is convex.
$\Clos(F^\circ)=F$.
\end{enumerate}
\item
Consider any $m\in\Z_+$ and any mapping $F:\{1,2,\ldots,m\}\rightarrow\mathcal{F}(S)$. 
The intersection $\cap_{i\in\{1,2,\ldots,m\}}F(i)$ is a face of $S$.
\item
Any proper face $F$ of $S$ is the intersection of all facets of $S$ containing $F$.
\item
We consider any two faces $F, G$ of $S$.
$F^\circ\cap G\neq\emptyset$, if and only if, $F\subset G$.
$F^\circ\cap G^\circ\neq\emptyset$, if and only if, $F=G$.
\item
$M=\Delta(S)=\Vect(S)^\vee$.
\item
Assume $\ell<s$. Let $F\in\mathcal{F}(S)^1$ be any facet.
\begin{enumerate}
\item
$M\subset\Delta(F)$. $M\neq\Delta(F)$.
\item
For any $\omega_F\in\Delta(F)- M$ we have $\Delta(F)=\R_0\omega_F+M$.
\item
If $S$ is rational over $N$, then $(\Delta(F)- M)\cap N^*\neq\emptyset$.
\end{enumerate}
\item
Note that $\ell<s$, if and only if, $\mathcal{F}(S)^1\neq\emptyset$.
In case $\ell<s$ we take any element $\omega_F\in\Delta(F)- M$ for any $F\in\mathcal{F}(S)^1$.
$$
S^\vee=
\Convcone(\{\omega_F|F\in\mathcal{F}(S)^1\})+M,\ 
S=
\bigcap_{F\in\mathcal{F}(S)^1}(\R_0\omega_F)^\vee\cap\Vect(S).
$$
\item
$S^\vee$ is a convex polyhedral cone in $V^*$.
If $S$ is rational over $N$, then $S^\vee$ is rational over $N^*$.
\item
Let $F$ be any face of $S$.
\begin{enumerate}
\item
$\Delta(F)$ is a face of $S^\vee$.
\item
$\Delta(F)=\Vect(F)^\vee\cap S^\vee$. $\Vect(\Delta(F))=\Vect(F)^\vee$.
\item
$\Delta^\circ(F)=\Delta(F)^\circ$. $\Delta(F)=\Clos(\Delta^\circ(F))$.
\end{enumerate}
\item
For any face $F$ of $S$, $\Delta(F,S|V)$ is a face of $S^\vee$, and $\dim\Delta(F,S|V)=\dim V-\dim F$. For any two faces $F$, $G$ of $S$ with $F\subset G$, $\Delta(F,S|V)\supset\Delta(G,S|V)$.

For any face $\bar{F}$ of $S^\vee$, $\Delta(\bar{F}, S^\vee |V^*)$ is a face of $S $, and $\dim\Delta(\bar{F}, S^\vee |V^*)=\dim V-\dim \bar{F}$. For any two faces $\bar{F}$, $\bar{G}$ of $S^\vee$ with $\bar{F}\subset \bar{G}$, $\Delta(\bar{F}, S^\vee |V^*)\supset\Delta(\bar{G}, S^\vee |V^*)$.

The mapping from $\mathcal{F}(S)$ to $\mathcal{F}(S^\vee)$ sending $F\in\mathcal{F}(S)$ to $\Delta(F,S|V)\in\mathcal{F}(S^\vee)$ and the mapping from $\mathcal{F}(S^\vee)$ to $\mathcal{F}(S)$ sending $\bar{F}\in\mathcal{F}(S^\vee)$ to\hfill\break $\Delta(\bar{F}, S^\vee |V^*)\in\mathcal{F}(S)$ are bijective mappings reversing the inclusion relation between $\mathcal{F}(S)$ and $\mathcal{F}(S^\vee)$, and they are the inverse mappings of each other.
Furthermore, if $F\in\mathcal{F}(S)$ and $\bar{F}\in\mathcal{F}(S^\vee)$ correspond to each other by them, then $\dim F+\dim\bar{F}=\dim V$.

\item
Assume $\ell<s$. Let $F\in\mathcal{F}(S)_{\ell+1}$ be any element.
\begin{enumerate}
\item
$L\subset F$. $L\neq F$.
\item
For any $t_F\in F- L$ we have $F=\R_0t_F+L$.
\end{enumerate}
\item
In case $\ell<s$ we take any element $t_F\in F- L$ for any $F\in\mathcal{F}(S)_{\ell+1}$.
Note that $\ell<s$, if and only if, $\mathcal{F}(S)_{\ell+1}\neq\emptyset$.
$$S =
\Convcone(\{t_F|F\in\mathcal{F}(S)_{\ell+1}\})+L,\ 
S^\vee=
\bigcap_{F\in\mathcal{F}(S)_{\ell+1}}(\R_0t_F)^\vee\cap\Vect(S^\vee).
$$
\item
$S$ is strongly convex$\Leftrightarrow \{0\}$ is a face of $S\Leftrightarrow S$ contains no vector subspace of $V$ of dimension positive $\Leftrightarrow S\cap(-S)=\{0\}\Leftrightarrow\dim S^\vee=\dim V$.
\item
The family $\{F^\circ|F\in\mathcal{F}(S)\}$ of subsets of $S$ gives the equivalence class decomposition of $S$.
In other words, the following three conditions hold:
\begin{enumerate}
\item
$F^\circ\neq\emptyset$ for any $F\in\mathcal{F}(S)$.
\item
If $F\in\mathcal{F}(S)$, $G\in\mathcal{F}(S)$, and $F^\circ\cap G^\circ\neq\emptyset$, then $F^\circ=G^\circ$.
\item
$$S=\bigcup_{ F\in\mathcal{F}(S)}F^\circ.$$
\end{enumerate}
\item
$\mathcal{F}(S)$ is a fan in $V$, and the support of  $\mathcal{F}(S)$ is equal to $S$.
In other words, the following four conditions hold:
\begin{enumerate}
\item
$\mathcal{F}(S)$ is a non-empty finite set whose elements are convex polyhedral cones in $V$.
\item
For any $F\in\mathcal{F}(S)$ and for any $G\in\mathcal{F}(S)$, $F\cap G$ is a face of $F$, and $F\cap G$ is a face of $G$.
\item
If $F\in\mathcal{F}(S)$ and $G$ is a face of $F$, then $G\in\mathcal{F}(S)$.
\item
$$S=\bigcup_{ F\in\mathcal{F}(S)}F.$$
\end{enumerate}
\item
Consider any vector space $W$ of finite dimension over $\R$ and any homomorphism $\pi:V\rightarrow W$ of vector spaces over $\R$.
The image $\pi(S)$ is a convex polyhedral cone in $W$, and it satisfies $\pi(S)^\circ=\pi(S^\circ)$.
\end{enumerate}
\end{theorem}

\begin{proof}
See Fulton~\cite{F93} and Cox~\cite{C10}.
\end{proof}

\begin{cor}
\label{good category}
\begin{enumerate}
\item
For any convex polyhedral cone $S$ in $V$, the dual cone $S^\vee$ is a convex polyhedral cone in $V^*$.
If moreover, $S$ is rational over $N$, then $S^\vee$ is rational over $N^*$.
\item
For any convex polyhedral cones $S$ and $T$ in $V$, $S+T$ and $S\cap T$ are convex polyhedral cones in $V$.
If moreover, $S$ and $T$ are rational over $N$, then $S+T$ and $S\cap T$ are rational over $N$.
\item
For any convex polyhedral cones $S$ and $T$ in $V$, $(S+T)^\vee=S^\vee\cap T^\vee$ and
$(S\cap T)^\vee=S^\vee+T^\vee$.
\item
For any convex polyhedral cone $S$ in $V$, any vector space $W$ of finite dimension over $\R$ and any homomorphism $\pi:V\rightarrow W$ of vector spaces over $\R$, $\pi(S)$ is a convex polyhedral cone in $W$.
If moreover, $S$ and $\pi^{-1}(0)$ are rational over $N$, then $\pi(S)$ is rational over $Q$ for any lattice $Q$ of $W$ with $\pi(N)=Q\cap\pi(V)$.
\item
For any convex polyhedral cone $S$ in $V$, any vector space $U$ of finite dimension over $\R$ and any homomorphism $\nu:U\rightarrow V$ of vector spaces over $\R$, $\nu^{-1}(S)$ is a convex polyhedral cone in $U$.
If moreover, $S$ and $\nu(U)$ are rational over $N$, then $\nu^{-1}(S)$ is rational over $K$ for any lattice $K$ of $U$ with $\nu(K)=N\cap\nu(U)$.

\end{enumerate}
\end{cor}

\begin{lemma}
\label{several cones}
Let $m\in \Z_+$ be any positive integer, and let $S$ be any mapping from the set $\{1,2,\ldots,m\}$ to the set of all convex polyhedral cones in $V$.
We denote
$$\bar{S}=\bigcap_{i\in\{1,2,\ldots,m\}}S(i)\subset V.$$
\begin{enumerate}
\item
$\bar{S}$ is a convex polyhedral cone in $V$.
$\bar{S}^\vee=\sum_{i\in\{1,2,\ldots,m\}} S(i)^\vee$.
If $S(i)$ is rational over $N$ for any $i\in\{1,2,\ldots,m\}$, then $\bar{S}$ is rational over $N$.
\item
If $\cap_{i\in\{1,2,\ldots,m\}}S(i)^\circ\neq\emptyset$, then $\bar{S}^\circ=\cap_{i\in\{1,2,\ldots,m\}}S(i)^\circ$.
\end{enumerate}

Let $\bar{F}$ be any face of $\bar{S}$.
\begin{enumerate}
\setcounter{enumi}{2}
\item
There exists uniquely a face $F(i)$ of $S(i)$ with $\bar{F}^\circ\subset F(i)^\circ$ for any $i\in\{1,2,\ldots, m\}$.
\end{enumerate}

Below, we assume that $F(i)\in\mathcal{F}(S(i))$ and $\bar{F}^\circ\subset F(i)^\circ$ for any $i\in\{1,2,\ldots,m\}$.

\begin{enumerate}
\setcounter{enumi}{3}
\item
$\bar{F}\subset F(i)$ for any $i\in\{1,2,\ldots,m\}$.
\item
$$\bar{F}=\bigcap_{i\in\{1,2,\ldots,m\}}F(i).$$
\item
$$\bar{F}^\circ=\bigcap_{i\in\{1,2,\ldots,m\}}F(i)^\circ.$$
\item
$$\Vect(\bar{F})=\bigcap_{i\in\{1,2,\ldots,m\}}\Vect(F(i)).$$
\item
$$\Delta(\bar{F},\bar{S})=\sum_{i\in\{1,2,\ldots,m\}}\Delta(F(i),S(i)).$$
\end{enumerate}

Let $G(i)$ be any face of $S(i)$ for any $i\in\{1,2,\ldots,m\}$.
\begin{enumerate}
\setcounter{enumi}{8}
\item
The intersection $\cap_{i\in\{1,2,\ldots,m\}}G(i)$ is a face of $\bar{S}$.
\item
If $\bar{F}\subset G(i)$ then $F(i)\subset G(i)$, for any $i\in\{1,2,\ldots,m\}$.
\item
If
$$\bar{F}=\bigcap_{i\in\{1,2,\ldots,m\}}G(i),$$
then $F(i)\subset G(i)$ for any $i\in\{1,2,\ldots,m\}$ and the following three conditions are equivalent:
\begin{enumerate}
\item
$$\bigcap_{i\in\{1,2,\ldots,m\}}G(i)^\circ\neq\emptyset.$$
\item
$F(i)=G(i)$ for any $i\in\{1,2,\ldots,m\}$.
\item
$$\bar{F}^\circ=\bigcap_{i\in\{1,2,\ldots,m\}}G(i)^\circ$$
\end{enumerate}
\end{enumerate}
\end{lemma}

\section{Simplicial cones and regular cones}
\label{simplex}

We study simplicial cones and regular cones.

Let $V$ be any vector space of finite dimension over $\R$, and let $N$ be any lattice of $V$.

\begin{lemma}
\label{simplex1}
Let $S$ be any simplicial cone in $V$. 
\begin{enumerate}
\item
The cone $S$ is strongly convex. The set $\{0\}$ is a face of $S$.
\item
$\sharp\mathcal{F}(S)_1=\dim S$.
$0\leq\dim S\leq\dim V$.
\item
Any face of $S$ is a simplicial cone in $V$.
For any face $F$ of $S$, $\mathcal{F}(F)\subset\mathcal{F}(S)$ and $\mathcal{F}(F)_1\subset\mathcal{F}(S)_1$.
\end{enumerate}
\end{lemma}

\begin{lemma}
\label{simplex2}
Let $S$ be any simplicial cone in $V$.
\begin{enumerate}
\item
$\mathcal{F}(F)_1\in 2^{\mathcal{F}(S)_1}$ and $\sharp \mathcal{F}(F)_1=\dim F$  for any $F\in\mathcal{F}(S)$.
$\mathcal{F}(F)_1\subset\mathcal{F}(G)_1$ for any $F\in\mathcal{F}(S)$ and any $G\in\mathcal{F}(S)$ with $F\subset G$.

$\sum_{E\in X}E\in\mathcal{F}(S)$ and $\dim(\sum_{E\in X}E)=\sharp X$ for any $X\in 2^{\mathcal{F}(S)_1}$.
$\sum_{E\in X}E\subset\sum_{E\in Y}E$ for any $X\in 2^{\mathcal{F}(S)_1}$ and any $Y\in 2^{\mathcal{F}(S)_1}$ with $X\subset Y$.

The mapping from $\mathcal{F}(S)$ to $2^{\mathcal{F}(S)_1}$ sending $F\in\mathcal{F}(S)$ to $\mathcal{F}(F)_1\in 2^{\mathcal{F}(S)_1}$ and the mapping from $2^{\mathcal{F}(S)_1}$ to $\mathcal{F}(S)$ sending $X\in 2^{\mathcal{F}(S)_1}$ to $\sum_{E\in X}E\in\mathcal{F}(S)$ are bijective mappings preserving the inclusion relation between  $\mathcal{F}(S)$ and $2^{\mathcal{F}(S)_1}$, and they are the inverse mappings of each other.

Furthermore, if $F\in\mathcal{F}(S)$ corresponds to $X\in 2^{\mathcal{F}(S)_1}$ by them, then $\dim F=\sharp X$.
The element $\{0\}\in\mathcal{F}(S)$ corresponds to $\emptyset\in 2^{\mathcal{F}(S)_1}$ by them, and $S\in\mathcal{F}(S)$ corresponds to $\mathcal{F}(S)_1\in 2^{\mathcal{F}(S)_1}$ by them.
\item
For any $X\in 2^{\mathcal{F}(S)_1}$ and any $Y\in 2^{\mathcal{F}(S)_1}$,
\begin{equation*}
\begin{split}
(\sum_{E\in X}E)\cap(\sum_{E\in Y}E)&=\sum_{E\in X\cap Y}E,\\
(\sum_{E\in X}E)+(\sum_{E\in Y}E)&=\sum_{E\in X\cup Y}E.
\end{split}
\end{equation*}
\item
For any $F\in\mathcal{F}(S)$ and any $G\in\mathcal{F}(S)$ the following claims hold:
\begin{enumerate}
\item
$F\cap G\in\mathcal{F}(S)$ and $\mathcal{F}(F\cap G)_1=\mathcal{F}(F)_1\cap\mathcal{F}(G)_1$.
\item
$F+G\in\mathcal{F}(S)$ and $\mathcal{F}(F+ G)_1=\mathcal{F}(F)_1\cup\mathcal{F}(G)_1$.
$F\subset F+G$ and $G\subset F+G$.
If $H\in \mathcal{F}(S)$ satisfies $F\subset H$ and $G\subset H$, then $F+G\subset H$.
$(F+G)^\circ=F^\circ+G^\circ$.
\item
$F\cap G=\{0\}$ and $F+G=S$, if and only if, $\mathcal{F}(F)_1\cap\mathcal{F}(G)_1=\emptyset$ and $\mathcal{F}(F)_1\cup\mathcal{F}(G)_1=\mathcal{F}(S)_1$
\end{enumerate}
\end{enumerate}
\end{lemma}

\begin{definition}
\label{opposite}
Let $S$ be any simplicial cone in $V$, and let $F$ be any face of $S$.
We denote
$$F\Op|S=\sum_{E\in\mathcal{F}(S)_1-\mathcal{F}(F)_1}E\in\mathcal{F}(S),$$
and we call $F\Op|S$ the \emph{opposite face} of $F$ over $S$.
When we need not refer to $S$, we also write simply $F\Op$, instead of $F\Op|S$.
\end{definition}

\begin{lemma}
\label{opposite2}
Let $S$ be any simplicial cone in $V$.
\begin{enumerate}
\item
For any face $F$ of $S$, the following claims hold:
\begin{enumerate}
\item
$F\Op=F\Op|S$ is a face of $S$.
$\dim F+\dim F\Op=\dim S$.
\item
$F\cap F\Op=\{0\}$ and $F+F\Op=S$. If $G\in\mathcal{F}(S)$, $F\cap G=\{0\}$ and $F+G=S$, then $G=F\Op$.
\item
$(F\Op)\Op=F$.
\item
$\mathcal{F}(F\Op)_1=\mathcal{F}(S)_1-\mathcal{F}(F)_1$.
If $G\in\mathcal{F}(S)$ and $\mathcal{F}(G)_1=\mathcal{F}(S)_1-\mathcal{F}(F)_1$,
then $G= F\Op$
\end{enumerate}
\item
$\{0\}\Op=S$. $S\Op=\{0\}$.
\item
For any $F\in\mathcal{F}(S)$ and any $G\in\mathcal{F}(S)$, the following claims hols:
\begin{enumerate}
\item
$F\subset G$, if and only if, $F\Op\supset G\Op$.
\item
$(F\cap G)\Op=F\Op+G\Op$.
\item
$(F+G)\Op=F\Op\cap G\Op$.
\end{enumerate}
\item
Consider any $F\in\mathcal{F}(S)$ and any $G\in\mathcal{F}(S)$ with $F\subset G$.
$F\in\mathcal{F}(G)$, $F\Op|G=(F\Op|S)\cap G$, $F\Op|S=(F\Op|G)+(G\Op|S)$ and $(F\Op|G)\cap(G\Op|S)=\{0\}$.
\item
The mapping from $\mathcal{F}(S)$ to itself sending $F\in \mathcal{F}(S)$ to $F\Op\in\mathcal{F}(S)$ is a bijective mapping reversing the inclusion relation. Its inverse mapping is equal to itself.
\end{enumerate}
\end{lemma}

\begin{lemma}
\label{mirror}
Let $S$ be any simplicial cone in $V$ with $\dim S=\dim V$ and let $B$ be any $\R$-basis of $V$ with $S=\Convcone(B)$. The dual basis of $B$ is denoted by $\{f^\vee|f\in B\}$. $\{f^\vee|f\in B\}$ is an $\R$-basis of $V^*$ and for any $f\in B$ and any $g\in B$
\begin{equation*}
\langle g^\vee, f\rangle=
\begin{cases}
1&\text{if $g=f$},\\
0&\text{if $g\neq f$}.
\end{cases}
\end{equation*}

$S^\vee=S^\vee|V=\Convcone(\{f^\vee|f\in B\})$. For any subset $X$ of $B$, 
$\sum_{f\in X}\R_0f$ is a face of $S$, $\sum_{f\in X}\R_0f^\vee$ is a face of $S^\vee$, and

\begin{equation*}
\begin{split}
\Delta((\sum_{E\in X}\R_0f)\Op|S,S)&=\sum_{E\in X}\R_0f^\vee,\\
\Delta((\sum_{E\in X}\R_0f^\vee)\Op| S^\vee,S^\vee)&=\sum_{E\in X}\R_0f.
\end{split}
\end{equation*}
\end{lemma}

\begin{lemma}
\label{regular1}
Let $S$ be any regular cone over $N$ in $V$.
\begin{enumerate}
\item
The cone $S$ is a rational simplicial cone over $N$ in $V$.
\item
The intersection $N\cap\Vect(S)$ is a lattice of $\Vect(S)$.
$S$ is a regular cone over $N\cap\Vect(S)$ in $\Vect(S)$.
The residue module $N/(N\cap\Vect(S))$ is a free module over $\Z$ of finite rank.
$\Rank(N/(N\cap\Vect(S)))=\dim V-\dim S$.
\item
Any face of $S$ is a regular cone over $N$ in $V$.
\item
If $\dim S=1$, then there exists uniquely an element $b_{S/N}$ of $S\cap N$ satisfying $S\cap N=\Z_0b_{S/N}$.
\end{enumerate}
\end{lemma}

\begin{definition}
\label{barycenter}
Let $S$ be any regular cone over $N$ in $V$.

If $\dim S=1$, we take the unique element $b_{S/N}$ of $S\cap N$ satisfying $S\cap N=\Z_0b_{S/N}$.

If $\dim S\neq 1$, we put
$$b_{S/N}=\sum_{E\in\mathcal{F}(S)_1} b_{E/N}\in S\cap N.$$

We call $b_{S/N}\in S\cap N$ the \emph{barycenter} of $S$ over $N$.
When we need not refer to $N$, we also write simply $b_S$, instead of $b_{S/N}$.
\end{definition}

\begin{lemma}
\label{regular2}
Let $S$ be any regular cone over $N$ in $V$.
\begin{enumerate}
\item
The set $\{b_E|E\in\mathcal{F}(S)_1\}$ is a $\Z$-basis of the lattice $N\cap\Vect(S)$ of $\Vect(S)$, and it is an $\R$-basis of the vector space $\Vect(S)$.
\begin{equation*}
\begin{split}
S&=\sum_{E\in\mathcal{F}(S)_1}\R_0b_E=\Convcone(\{b_E|E\in\mathcal{F}(S)_1\}),\\
S^\circ&=\sum_{E\in\mathcal{F}(S)_1}\R_+b_E.
\end{split}
\end{equation*}
For any $E\in\mathcal{F}(S)_1$, $E=\R_0b_E$.
$\mathcal{F}(S)_1=\{\R_0b_E| E\in\mathcal{F}(S)_1\}$.
\item
If a basis $B$ over $\Z$ of $N$ and a subset $C$ of $B$ satisfies $S=\Convcone(C)$, then $\sharp C=\dim S$ and $C=\{b_E|E\in\mathcal{F}(S)_1\}$.
\end{enumerate}
\end{lemma}

\begin{lemma}
\label{basis corresp}
\begin{enumerate}
\item
For any regular cone $S$ over $N$ in $V$ with $\dim S=\dim V$, the set $\{b_E|E\in\mathcal{F}(S)_1\}$ is a $\Z$-basis of $N$.

For any $\Z$-basis $B$ of $N$, $\Convcone(B)$ is a regular cone over $N$ in $V$ with $\dim \Convcone(B)=\dim V$.

The mapping sending any regular cone $S$ over $N$ in $V$ with $\dim S=\dim V$ to $\{b_E|E\in\mathcal{F}(S)_1\}$ and the mapping sending any $\Z$-basis $B$ of $N$ to $\Convcone(B)$ are bijective mappings between the set of all rgular cones $S$ over $N$ in $V$ with $\dim S=\dim V$ and the set of all $\Z$-bases $B$ of $N$, and they are the inverse mappings of each other.
\end{enumerate}

Below we consider any regular cone $S$ over $N$ in $V$ with $\dim S=\dim V$. Note that $\{b_E|E\in\mathcal{F}(S)_1\}$ is a $\Z$-basis of $N$ and it is an $\R$-basis of $V$. By $\{b_E^\vee|E\in\mathcal{F}(S)_1\}$ we denote the dual basis of $\{b_E|E\in\mathcal{F}(S)_1\}$, which is a $\R$-basis of $V^*$. We assume that for any $D\in\mathcal{F}(S)_1$ and any $E\in\mathcal{F}(S)_1$
\begin{equation*}
\langle b_D^\vee, b_E\rangle=
\begin{cases}
1&\text{if $D=E$},\\
0&\text{if $D\neq E$}.
\end{cases}
\end{equation*}

\begin{enumerate}
\setcounter{enumi}{1}
\item
The set $\{b_E^\vee|E\in\mathcal{F}(S)_1\}$ is a $\Z$-basis of $N^*$.
\item
$S^\vee$ is a regular cove over $N^*$ in $V^*$ with $\dim S^\vee=\dim V^*$.
$$S^\vee=\Convcone(\{b_E^\vee|E\in\mathcal{F}(S)_1\}).$$
For any $E\in\mathcal{F}(S)_1$, $\R_0 b_E^\vee\in\mathcal{F}(S^\vee)_1$ and $b_{\R_0 b_E^\vee}= b_E^\vee$.
$\mathcal{F}(S^\vee)_1=\{\R_0b_E^\vee| E\in\mathcal{F}(S)_1\}$, and
$\{b_D|D\in \mathcal{F}(S^\vee)_1\}=\{b_E^\vee|E\in\mathcal{F}(S)_1\}$.
\end{enumerate}
\end{lemma}

\begin{lemma}
\label{barycenter property}
Let $S$ be any regular cone over $N$ in $V$.
\begin{enumerate}
\item
Let $W$ be any vector space of finite dimension over $\R$, let $Q$ be any lattice of $W$, and let $T$ be any regular cone over $Q$ in $W$.
$\dim S=\dim T$, if and only if, there exists an isomorphism $\phi:\Vect(S)\rightarrow\Vect(T)$ of vector spaces over $\R$ satisfying $\phi(S)=T$ and $\phi(N\cap\Vect(S))=Q\cap\Vect(T)$.
\item
$b_S=b_{S/N}\in S^\circ\cap N$.
$\R_0b_S\cap N=\Z_0b_S$.
\item
Consider any $a\in S\cap N$ with $\R_0a\cap N=\Z_0a$.
$\phi(a)=a$ for any homomorphism $\phi:\Vect(S)\rightarrow\Vect(S)$ of vector spaces over $\R$ satisfying $\phi(S)=S$ and $\phi(N\cap\Vect(S))=N\cap\Vect(S)$, if and only if, $a=0$ or $a=b_S$.
\item
$b_S=0\Leftrightarrow\dim S=0$.
$\R_0b_S\subset S$.
$\R_0b_S=S\Leftrightarrow \R_0b_S$ is a face of $S\Leftrightarrow \dim S\leq 1$
\item
Assume $F\in\mathcal{F}(S), \dim F\geq 1, \Lambda\in\mathcal{F}(S)$, and $F\not\subset\Lambda$.
\begin{enumerate}
\item
$\Lambda+\R_0b_F$ is a regular cone over $N$ in $V$.
$\dim(\Lambda+\R_0b_F)=\dim\Lambda+1$.
\item
$\R_0b_F\in\mathcal{F}(\Lambda+\R_0b_F)_1$.
$\Lambda\in\mathcal{F}(\Lambda+\R_0b_F)^1$.
$\R_0b_F\cap\Lambda=\{0\}$.
$\R_0b_F=\Lambda\Op|(\Lambda+\R_0b_F)$.
$\Lambda=(\R_0b_F)\Op|(\Lambda+\R_0b_F)$.
\item
$\mathcal{F}(\Lambda)\subset\{\Lambda'\in\mathcal{F}(S)|F\not\subset\Lambda'\}$.
$\mathcal{F}(\Lambda)=\{\Lambda'\in\mathcal{F}(\Lambda+\R_0b_F)| \R_0b_F\not\subset\Lambda'\}$.
$\{\Lambda'+\R_0b_F|\Lambda'\in\mathcal{F}(\Lambda)\} =\{\Lambda'\in\mathcal{F}(\Lambda+\R_0b_F)| \R_0b_F\subset\Lambda'\}$.
\item
$\Lambda+\R_0b_F\subset \Lambda+F\in\mathcal{F}(S)$.
$(\Lambda+\R_0b_F)^\circ\subset (\Lambda+F)^\circ$.
\item
If $\dim F\geq 2$, then $\Lambda+\R_0b_F\not\in\mathcal{F}(S)$.
If $\dim F=1$, then $\Lambda+\R_0b_F=\Lambda+F\in\mathcal{F}(S)$.
\end{enumerate}
\end{enumerate}
\end{lemma}

\section{Fans}
\label{decomposition}

We begin the study of fans. We define notations and concepts to develop our theory.

Let $V$ be any vector space of finite dimension over $\R$, and let $N$ be any lattice of $V$.
Let $\Phi$ be any finite set whose elements are convex polyhedral cones in $V$.

We say that $\Phi$ is \emph{strongly convex}, if any $\Delta\in\Phi$ is strongly convex.
We say that $\Phi$ is \emph{simplicial}, if any $\Delta\in\Phi$ is a simplicial cone.
We say that $\Phi$ is \emph{rational} over $N$, if any $\Delta\in\Phi$ is rational over $N$.
We say that $\Phi$ is \emph{regular} over $N$, if any $\Delta\in\Phi$ is a regular cone over $N$.

We denote
\begin{equation*}
\begin{split}
|\Phi|&=\bigcup_{\Delta\in \Phi}\Delta\subset V,\\
|\Phi|^\circ&=\bigcup_{\Delta\in \Phi}\Delta^\circ\subset V.
\end{split}
\end{equation*}
We call $|\Phi|$ and $|\Phi|^\circ$ the \emph{support} of $\Phi$ and the \emph{open support} of $\Phi$ respectively.
We denote
$$\Phi\Mx=\{\Delta\in \Phi|\text{ If }\Lambda\in\Phi\text{ and }\Delta\subset\Lambda\text{, then }\Delta=\Lambda\}\subset\Phi.$$
We call any element in $\Phi\Mx$ a \emph{maximal element} of $\Phi$
and we call $\Phi\Mx$ the set of \emph{maximal elements} of $\Phi$.

Note that $\mathcal{F}(\Delta)$ is a non-empty finite set whose elements are convex polyhedral cones in $V$ for any $\Delta\in\Phi$.
We denote
$$\Phi\Fc=\bigcup_{\Delta\in\Phi}\mathcal{F}(\Delta)\subset 2^V,$$
and we call $\Phi\Fc$ the \emph{face closure} of $\Phi$.

In case $\Phi\neq\emptyset$ we define
$$\dim\Phi=\max\{\dim\Delta|\Delta\in\Phi\}\in\Z_0,$$
and we call $\dim\Phi$ the \emph{dimension} of $\Phi$.
In case  $\Phi=\emptyset$ we do not define $\dim\Phi$.
For any $i\in\Z$ we denote
$$\Phi_i=\{\Delta\in\Phi|\dim\Delta=i\},$$
\begin{equation*}
\Phi^i=
\begin{cases}
\{\Delta\in\Phi|\dim\Delta=\dim\Phi-i\}&\text{if $\Phi\neq\emptyset$},\\
\emptyset&\text{if $\Phi=\emptyset$}.
\end{cases}
\end{equation*}
$\Phi_i$ and $\Phi^i$ are subsets of $\Phi$.

Consider any subset $F$ of $V$.
We denote
\begin{equation*}
\begin{split}
\Phi\backslash F&=\{\Delta\in\Phi|\Delta\subset F\}\subset \Phi,\\
\Phi/F&=\{\Delta\in\Phi|\Delta\supset F\}\subset\Phi.
\end{split}
\end{equation*}

Consider any vector space $W$ of finite dimension and any homomorphism $\pi:V\rightarrow W$ of vector spaces over $\R$.
We denote
$$\pi_*\Phi=\{\pi(\Delta)|\Delta\in\Phi\}\subset 2^W,$$
and we call $\pi_*\Phi$ the \emph{push-down} of $\Phi$ by $\pi$.

Consider any vector space $U$ of finite dimension and any homomorphism $\nu:U\rightarrow V$ of vector spaces over $\R$.
We denote
$$\nu^*\Phi=\{\nu^{-1}(\Delta)|\Delta\in\Phi\}\subset 2^U,$$
and we call $\nu^*\Phi$ the \emph{pull-back} of $\Phi$ by $\nu$.

Consider the case $\Phi\neq\emptyset$.
We say that $\Phi$ is \emph{flat}, if $\dim\Phi=\dim\Vect(|\Phi|)$ and $\Phi\Mx=\Phi^0$. Let $\Delta\in\Phi$ be any element. We say that $\Phi$ is \emph{starry with center in} $\Delta$, if $\Phi=(\Phi/\Delta)\Fc$.

\begin{lemma}
\label{first step}
Let $\Phi$ be any finite set whose elements are convex polyhedral cones in $V$.
\begin{enumerate}
\item
$|\Phi|$ is a closed subset of $V$. If $|\Phi|\neq\emptyset$, $|\Phi|$ is a cone in $V$.
\item
$|\Phi|^\circ\subset|\Phi|$.
$\Clos(|\Phi|^\circ)=|\Phi|$.
\item
$\Phi\Mx\subset\Phi$.
$(\Phi\Mx)\Mx=\Phi\Mx$.
$|\Phi\Mx|=|\Phi|$.
\item
$\Phi\Fc$ is a finite set whose elements are convex polyhedral cones in $V$.
$\Phi\subset\Phi\Fc$.
$(\Phi\Fc)\Fc=\Phi\Fc$.
$|\Phi\Fc|=|\Phi|$.
$(\Phi\Fc)\Mx=\Phi\Mx$.
If $\Phi$ is simplicial, then $\Phi\Fc$ is also simplicial.
If $\Phi$ is rational over $N$, then $\Phi\Fc$ is also rational over $N$.
If $\Phi$ is regular over $N$, then $\Phi\Fc$ is also regular over $N$.
\item
Consider any vector space $W$ of finite dimension and any homomorphism $\pi:V\rightarrow W$ of vector spaces over $\R$. The set $\pi_*\Phi$ is a finite set whose elements are convex polyhedral cones in $W$.
$|\pi_*\Phi|=\pi(|\Phi|)$.
If $\Phi$ and $\pi^{-1}(0)$ are rational over $N$, then $\pi_*\Phi$ is rational over $Q$ for any lattice $Q$ of $W$ with $\pi(N)=Q\cap\pi(V)$.
\item
$\Id_{V*}\Phi=\Phi$.
For any vector spaces $W$, $W'$ of finite dimension and any homomorphisms $\pi:V\rightarrow W$, $\pi':W\rightarrow W'$ of vector spaces over $\R$,
$(\pi'\pi)_*\Phi=\pi'_*\pi_*\Phi$.
\item
Consider any vector space $U$ of finite dimension and any homomorphism $\nu:U\rightarrow V$ of vector spaces over $\R$. The set $\nu^*\Phi$ is a finite set whose elements are convex polyhedral cones in $U$.
$|\nu^*\Phi|=\nu^{-1}(|\Phi|)$.
If $\Phi$ and $\nu(U)$ are rational over $N$, then $\nu^*\Phi$ is rational over $K$ for any lattice $K$ of $U$ with $\nu(K)=N\cap\nu(U)$.
\item
$\Id_V^*\Phi=\Phi$.
For any vector spaces $U$, $U'$ of finite dimension and any homomorphisms $\nu:U\rightarrow V$, $\nu':U'\rightarrow U$ of vector spaces over $\R$,
$(\nu\nu')^*\Phi=\nu^{\prime*}\nu^*\Phi$.
\item
Consider any vector spaces $W$, $U$ of finite dimension and any homomorphisms $\pi:V\rightarrow W$, $\nu:U\rightarrow V$ of vector spaces over $\R$.

$\Phi=\emptyset\Leftrightarrow|\Phi|=\emptyset\Leftrightarrow0\not\in|\Phi|\Leftrightarrow|\Phi|^\circ=\emptyset\Leftrightarrow\Phi\Mx=\emptyset\Leftrightarrow\Phi\Fc=\emptyset\Leftrightarrow\pi_*\Phi=\emptyset\Leftrightarrow\nu^*\Phi=\emptyset$.

If $\Phi\neq\emptyset$, then $\dim\pi_*\Phi\leq\dim\Phi=\dim\Phi\Mx=\dim\Phi\Fc\leq\dim\Vect(|\Phi|)$.
\item
Consider any subset $\Sigma$ of $\Phi$, any $i\in\Z$, any subset $F$ of $V$, any vector spaces $W$, $U$ of finite dimension and any homomorphisms $\pi:V\rightarrow W$, $\nu:U\rightarrow V$ of vector spaces over $\R$.

$|\Sigma|\subset|\Phi|$, $|\Sigma|^\circ\subset|\Phi|^\circ$,
$\Sigma\Fc\subset\Phi\Fc$,
$\Sigma_i\subset\Phi_i$,
$\Sigma\backslash F\subset\Phi\backslash F$,
$\Sigma/F\subset\Phi/F$,
$\pi_*\Sigma\subset\pi_*\Phi$, and
$\nu^*\Sigma\subset\nu^*\Phi$.
\item
Let $F$ and $G$ be any subsets of $V$.
If $F\supset G$, $(\Phi\backslash F)\backslash G=\Phi\backslash G$.
If $F\subset G$, $(\Phi/F)/G=\Phi/G$.
$\Phi\backslash F=\Phi\Leftrightarrow |\Phi|\subset F$.
\item
Assume $\Phi\neq\emptyset$. 

$\Phi\Mx\supset\Phi^0\neq\emptyset$.
$\Phi\Mx=\Phi^0$, if and only if, $\dim\Delta=\dim\Phi$ for any $\Delta\in\Phi\Mx$.

The following four conditions are equivalent:
\begin{enumerate}
\item
$\Phi$ is flat. In other words, $\dim\Phi=\dim\Vect(|\Phi|)$ and $\Phi\Mx=\Phi^0$.
\item
$\dim\Delta=\dim\Vect(|\Phi|)$ for any $\Delta\in\Phi\Mx$.
\item
$\Vect(\Delta)= \Vect(|\Phi|)$ for any $\Delta\in\Phi\Mx$.
\item
$\Vect(\Delta)= \Vect(\Lambda)$ for any $\Delta\in\Phi\Mx$ and any $\Lambda\in\Phi\Mx$.
\end{enumerate}
\end{enumerate}
\end{lemma}

\begin{definition}
\label{cpcd}
\begin{enumerate}
\item
Any subset $\Sigma$ of $2^V$ satisfying the following three conditions is called a \emph{fan} in $V$.
\begin{enumerate}
\item
The set $\Sigma$ is a non-empty finite set whose elements are convex polyhedral cones in $V$.
\item
For any $\Delta\in\Sigma$ and any $\Lambda\in\Sigma$, $\Delta\cap\Lambda$ is a face of $\Delta$, and $\Delta\cap\Lambda$ is a face of $\Lambda$.
\item
For any $\Delta\in\Sigma$ and any face $\Lambda$ of $\Delta$, $\Lambda\in\Sigma$.
\end{enumerate}
\item
Let $\Sigma$ be any fan in $V$. We call an element $L\in\Sigma$ the \emph{mimimum element} of $\Sigma$, if $\dim L\leq\dim\Delta$ for any $\Delta\in\Sigma$.
\item
Let $\Sigma$ and $\Phi$ be any finite sets whose elements are convex polyhedral cones in $V$.

If for any $\Delta\in\Sigma$ there exists $\Lambda\in\Phi$ with $\Delta\subset\Lambda$, then we say that $\Sigma$ is a \emph{subdivision} of $\Phi$.
\item
Let $J$ be any finite set and let $\Phi$ be any mapping from $J$ to the set of all finite sets whose elements are convex polyhedral cones in $V$.
For any $j\in J$, $\Phi(j)$ is a finite set whose elements are convex polyhedral cones in $V$.
We denote
\begin{equation*}
\begin{split}
\bigcap_{j\in J}^\wedge\Phi(j)
=\{\bar{\Delta}\in 2^V|&\bar{\Delta}=\cap_{j\in J}\Delta(j)\text{ for some mapping }\Delta:J\rightarrow 2^V\\
&\quad\text{such that }\Delta(j)\in\Phi(j)\text{ for any }j\in J\}\subset 2^V,
\end{split}
\end{equation*}
and we call $\hat{\cap}_{j\in J}\Phi(j)$ the \emph{real intersection} of $\Phi(j), j\in J$.

Note that $\hat{\cap}_{j\in J}\Phi(j)$ is different from the intersection $\cap_{j\in J}\Phi(j)$ of $\Phi(j), j\in J$.

When $J=\{1,2,\ldots,m\}$ for some $m\in\Z_+$, we also denote
$$\Phi(1)\hat{\cap}\Phi(2) \hat{\cap}\cdots\hat{\cap} \Phi(m)
=\bigcap_{j\in J}^\wedge\Phi(j).$$
\item
Let $\Sigma$ be any fan in $V$, let $W$ be any vector space of finite dimension over $\R$, let $T$ be any subset of $W$, and let $\phi:| \Sigma|\rightarrow T$ be any mapping.

We say that $\phi$ is \emph{piecewise linear}, if for any $\Delta\in\Sigma$ there exists a homomorphism $\phi_\Delta:\Vect(\Delta)\rightarrow W$ of vector spaces over $\R$ such that $\phi(a)= \phi_\Delta(a)$ for any $a\in\Delta$.

\item
Let $\Sigma$ be any fan in $V$, and let $\phi:| \Sigma|\rightarrow \R$ be any piecewise linear function.

Assume that $\Sigma$ is rational over $N$. We say that $\phi$ is \emph{rational} over $N$, if for any $\Delta\in\Sigma$, there exists a linear function $\phi_\Delta:\Vect(\Delta)\rightarrow \R$ such that $\phi(a)= \phi_\Delta(a)$ for any $a\in\Delta$ and $\phi_\Delta(N\cap\Vect(\Delta))\subset\Q$.

Assume that the support $|\Sigma|$ of $\Sigma$ is convex.
We say that $\phi$ is \emph{convex} over $\Sigma$, if the following two conditions are satisfied:
\begin{enumerate}
\item
For any $a\in |\Sigma|$, any $b\in |\Sigma|$ and any $t\in\R$ with $0\leq t\leq 1$, $\phi((1-t)a+tb)\geq (1-t)\phi(a)+t\phi(b)$.
\item
If $a\in |\Sigma|$, $b\in |\Sigma|$, $t\in\R$, $0<t<1$ and $\phi((1-t)a+tb)=(1-t)\phi(a)+t\phi(b)$, then $\{a,b\}\subset\Delta$ for some $\Delta\in\Sigma$.
\end{enumerate}
\end{enumerate}
\end{definition}

\begin{example}
Let $S$ be any convex polyhedral cone in $V$.
$\mathcal{F}(S)$ is a fan in $V$.
$|\mathcal{F}(S)|=S$.
If $S$ is a simplicial cone, then $\mathcal{F}(S)$ is a simplicial fan.
If $S$ is a rational convex polyhedral cone over $N$, then $\mathcal{F}(S)$ is a rational fan over $N$.
If $S$ is a regular cone over $N$, then $\mathcal{F}(S)$ is a regular fan over $N$.
\end{example}

\begin{lemma}
\label{basic cpcd}
Let $\Sigma$ be any fan in $V$.
\begin{enumerate}
\item
$\Delta\cap\Lambda\in\Sigma$ for any $\Delta\in\Sigma$ and any  $\Lambda\in\Sigma$.
\item
Consider any $\Delta\in\Sigma$ and any $\Lambda\in\Sigma$.
The following three conditions are equivalent:
\begin{enumerate}
\item
$\Delta\subset\Lambda$.
\item
$\Delta$ is a face of $\Lambda$.
\item
$\Delta^\circ\cap\Lambda\neq\emptyset$.
\end{enumerate}
\item
Consider any $\Delta\in\Sigma$ and any $\Lambda\in\Sigma$.
The following three conditions are equivalent:
\begin{enumerate}
\item
$\Delta=\Lambda$.
\item
$\dim\Delta=\dim\Lambda$ and $\Delta\subset\Lambda$, or $\dim\Delta=\dim\Lambda$ and $\Delta\supset\Lambda$.
\item
$\Delta^\circ\cap\Lambda^\circ\neq\emptyset$.
\end{enumerate}
\item
$\Sigma=\Sigma\Fc=(\Sigma\Mx)\Fc$.
$|\Sigma|=|\Sigma|^\circ=|\Sigma\Mx|$.
For any subset $\Phi$ of $\Sigma$, $\Phi\Fc=\Sigma\backslash|\Phi|\subset\Sigma$.
For any $\Delta\in\Sigma$, $\mathcal{F}(\Delta)=\{\Delta\}\Fc=\Sigma\backslash\Delta\subset\Sigma$.
\item
The family $\{\Delta^\circ|\Delta\in\Sigma\}$ of subsets of $V$ gives an equivalence class decomposition of $|\Sigma|$, in other words, the following three conditions hold:
\begin{enumerate}
\item
$\Delta^\circ\neq\emptyset$ for any $\Delta\in\Sigma$.
\item
If $\Delta^\circ\cap\Lambda^\circ\neq\emptyset$, then $\Delta^\circ=\Lambda^\circ$ for any $\Delta\in\Sigma$ and any $\Lambda\in\Sigma$.
\item
$$|\Sigma|=\bigcup_{\Delta\in\Sigma}\Delta^\circ.$$
\end{enumerate}
\item
For any subset $\Phi$ of $\Sigma$, the following three conditions are equivalent:
\begin{enumerate}
\item
$\Phi$ is a fan in $V$.
\item
$\Phi\neq\emptyset$ and $\Sigma\backslash\Delta\subset\Phi$ for any $\Delta\in\Phi$.
\item
$\Phi\neq\emptyset$ and $\Phi=\Phi\Fc$.
\end{enumerate}
\item
For any non-empty subset $\Phi$ of $\Sigma$, $\Phi\Fc$ is a fan in $V$.
\item
Consider any subset $F$ of $V$.
If $\Sigma\backslash F\neq\emptyset$, then $\Sigma\backslash F$ is a fan.
If $\Sigma-(\Sigma/F)\neq\emptyset$, then $\Sigma-(\Sigma/F)$ is a fan.
\item
For any subset $X$ of $|\Sigma|$, the following three conditions are equivalent:
\begin{enumerate}
\item
$X$ is a closed subset of $V$.
\item
$X$ is a closed subset of $|\Sigma|$.
\item
$X\cap\Delta$ is a closed subset of $\Delta$ for any $\Delta\in\Sigma$.
\end{enumerate}
\item
For any subset $Y$ of $|\Sigma|$, the following two conditions are equivalent:
\begin{enumerate}
\item
$Y$ is an open subset of $|\Sigma|$.
\item
$Y\cap\Delta$ is an open subset of $\Delta$ for any $\Delta\in\Sigma$.
\end{enumerate}
\item
Consider any $\Delta\in\Sigma$. $\Delta\in\Sigma\Mx$, if and only if, $\Delta^\circ$ is an open subset of $|\Sigma|$.
\item
Consider any $L\in\Sigma$. 
The following five conditions are equivalent:
\begin{enumerate}
\item
$L$ is the minimum element of $\Sigma$, in other words, $\dim L\leq\dim\Delta$ for any $\Delta\in\Sigma$.
\item
$L\subset\Delta$ for any $\Delta\in\Sigma$.
\item
$L=\Delta\cap(-\Delta)$ for any $\Delta\in\Sigma$.
\item
$L=\Delta\cap(-\Delta)$ for some $\Delta\in\Sigma$.
\item
$L$ is a vector subspace over $\R$ of $V$.
\end{enumerate}

There exists a unique element $L\in\Sigma$ satisfying the above five conditions.
\end{enumerate}

Below, we assume that $L\in\Sigma$ is the minimum element of $\Sigma$.
\begin{enumerate}
\setcounter{enumi}{12}
\item
$\Sigma_{\dim L}=\{L\}$.
$\Sigma_i\neq\emptyset$, if and only if, $\dim L\leq i\leq\dim\Sigma$ for any $i\in\Z$.
$\Sigma^i\neq\emptyset$, if and only if, $0\leq i\leq\dim\Sigma-\dim L $ for any $i\in\Z$.
\item
$L=\Delta\cap(-\Delta)\in\mathcal{F}(\Delta)$,
$\Delta+L=\Delta$, $\Delta^\circ+L=\Delta^\circ$ and  $\Vect(\Delta)+L=\Vect(\Delta)$ for any  $\Delta\in\Sigma$.
\item
$\Sigma$ is strongly convex, if and only if, $L=\{0\}$, if and only if, $\{0\}\in\Sigma$.
\item
$\Sigma$ is rational over $N$, if and only if, 
any $\Delta\in\Sigma$ with $\dim\Delta\leq\dim L+1$ is rational over $N$.
\end{enumerate}
Below, we consider any vector space $W$ of finite dimension over $\R$ and any homomorphism $\pi:V\rightarrow W$ of vector spaces over $\R$ satisfying $\pi^{-1}(0)\subset L$.

\begin{enumerate}
\setcounter{enumi}{16}
\item
The push-down $\pi_*\Sigma$ is a fan in $W$, and $\pi^*\pi_*\Sigma=\Sigma$.
\item
For any $\Delta\in\Sigma$, $\pi(\Delta)\in\pi_*\Sigma$.
For any $\bar{\Delta}\in\pi_*\Sigma$, $\pi^{-1}(\bar{\Delta})\in\Sigma$.

The mapping from $\Sigma$ to $\pi_*\Sigma$ sending any $\Delta\in\Sigma$ to $\pi(\Delta)\in\pi_*\Sigma$ and the mapping from $\pi_*\Sigma$ to $\Sigma$ sending any $\bar{\Delta}\in\pi_*\Sigma$ to $\pi^{-1}(\bar{\Delta})\in\Sigma$ are bijective mappings preserving the inclusion relation between $\Sigma$ and $\pi_*\Sigma$, and they are the inverse mappings of each other.

Furthermore, if $\Delta\in\Sigma$ and $\bar{\Delta}\in\pi_*\Sigma$ correspond to each other by them, the following equalities hold:
\begin{enumerate}
\item
$\dim\bar{\Delta}=\dim\Delta-\dim\pi^{-1}(0)$.
\item
$\Vect(\bar{\Delta})=\pi(\Vect(\Delta))$,
$\pi^{-1}(\Vect(\bar{\Delta}))=\Vect(\Delta)$.
\item
$\bar{\Delta}^\circ=\pi(\Delta^\circ)$,
$\pi^{-1}(\bar{\Delta}^\circ)= \Delta^\circ$.
\item
$\mathcal{F}(\bar{\Delta})=\pi_*\mathcal{F}(\Delta)$,
$\pi^*\mathcal{F}(\bar{\Delta})=\mathcal{F}(\Delta)$.
\end{enumerate}
\end{enumerate}
\end{lemma}

\begin{lemma}
\label{starry}
Let $\Sigma$ be any simplicial fan in $V$ and let $H\in\Sigma$ be any element. Assume $\Sigma$ is starry with center in $H$, in other words, $\Sigma=(\Sigma/H)\Fc$.

Let $W=V/\Vect(H)$ denote the residue vector space of $V$ by $\Vect(H)$ and $\pi:V\rightarrow W$ denote the canonical surjective homomorphism of vector spaces over $\R$.
\begin{enumerate}
\item
The push-down $\pi_*\Sigma$ is a simplicial fan in $W$. $\dim \pi_*\Sigma=\dim\Sigma-\dim H$.
\item If $\Sigma$ is rational over $N$, then $\pi(N)$ is a lattice in $W$ and $\pi_*\Sigma$ is rational over $\pi(N)$. If $\Sigma$ is regular over $N$, then $\pi(N)$ is a lattice in $W$ and $\pi_*\Sigma$ is regular over $\pi(N)$. If $\Sigma$ is flat, then $\pi_*\Sigma$ is flat.
\end{enumerate}

Let $\bar{\Sigma}=\{\Delta\in\Sigma|\Delta\cap H=\{0\}\}\subset\Sigma$.
\begin{enumerate}
\setcounter{enumi}{2}
\item
$\bar{\Sigma}$ is a simplicial fan in $V$.  $\dim \bar{\Sigma}=\dim\Sigma-\dim H$.
If $\Sigma$ is rational over $N$, then $\bar{\Sigma}$ is rational over $N$.
If $\Sigma$ is regular over $N$, then $\bar{\Sigma}$ is regular over $N$.
\item
For any $\Delta\in\bar{\Sigma}$, $\pi(\Delta)\in \pi_*\Sigma$. For any $\bar{\Delta}\in \bar{\Sigma}$, $\pi^{-1}(\bar{\Delta})\cap|\bar{\Sigma}|\in\bar{\Sigma}$.

The mapping from $\bar{\Sigma}$ to $\pi_*\Sigma$ sending $\Delta\in\bar{\Sigma}$ to $\pi(\Delta)\in \pi_*\Sigma$ and the maping from $\pi_*\Sigma$ to $\bar{\Sigma}$ sending $\bar{\Delta}\in \bar{\Sigma}$ to $\pi^{-1}(\bar{\Delta})\cap|\bar{\Sigma}|\in\bar{\Sigma}$ are bijective mapping preserving the inclusion relation and the dimension between $\bar{\Sigma}$ and $\pi_*\Sigma$, and they are the inverse mappings of each other.
\item
$\pi(|\Sigma|)=\pi(|\bar{\Sigma}|)=|\pi_*\Sigma |$.
$\pi(\Vect(|\Sigma|))=\pi(\Vect(|\bar{\Sigma}|))=\Vect(|\pi_*\Sigma |)$.
The mapping $\pi: |\bar{\Sigma}|\rightarrow |\pi_*\Sigma |$ induced by $\pi$ is a homeomorphism.
\end{enumerate}

Consider the product vector space $W\times\Vect(H)$. Let
$$\pi_*\Sigma\hat{\times}\mathcal{F}(H)=\{\bar{\Delta}\times\Lambda|\bar{\Delta}\in \pi_*\Sigma, \Lambda\in\mathcal{F}(H)\}.$$
\begin{enumerate}
\setcounter{enumi}{5}

\item
$\pi_*\Sigma\hat{\times}\mathcal{F}(H)$ is a simplicial fan in $W\times\Vect(H)$.
$\dim \pi_*\Sigma\hat{\times}\mathcal{F}(H)=\dim\Sigma$.
\item
If $\Sigma$ is rational over $N$, then $\pi(N)\times(N\cap\Vect(H))$ is a lattice in $W\times\Vect(H)$ and $\pi_*\Sigma\hat{\times}\mathcal{F}(H)$ is rational over $\pi(N)\times(N\cap\Vect(H))$. If $\Sigma$ is regular over $N$, then $\pi(N)\times(N\cap\Vect(H))$ is a lattice in $W\times\Vect(H)$ and $\pi_*\Sigma\hat{\times}\mathcal{F}(H)$ is regular over  $\pi(N)\times(N\cap\Vect(H))$. If $\Sigma$ is flat, then $\pi_*\Sigma\hat{\times}\mathcal{F}(H)$ is flat.
\item
For any $\Delta\in\Sigma$, $\pi(\Delta)\times(\Delta\cap H)\in \pi_*\Sigma\hat{\times}\mathcal{F}(H)$.

The mapping from $\Sigma$ to $\pi_*\Sigma\hat{\times}\mathcal{F}(H)$ sending $\Delta\in\Sigma$ to $\pi(\Delta)\times(\Delta\cap H)\in \pi_*\Sigma\hat{\times}\mathcal{F}(H)$ is bijective, it preserves the dimension, and itself and its inverse mapping preserve the inclusion relation.
\item
Assume $\dim H=1$.

$\bar{\Sigma}=\Sigma-(\Sigma/H)$. $\dim(\Sigma-(\Sigma/H))=\dim\Sigma-1$.

For any $\Delta\in\Sigma/H$, $H\in\mathcal{F}(\Delta)_1$ and $H\Op|\Delta\in\Sigma-(\Sigma/H)$.
For any $\bar{\Delta}\in\Sigma-(\Sigma/H)$, $\bar{\Delta}+H\in\Sigma/H$.

The mapping from $\Sigma/H$ to $\Sigma-(\Sigma/H)$ sending $\Delta\in\Sigma/H$ to $H\Op|\Delta\in\Sigma-(\Sigma/H)$ and the mapping from $\Sigma-(\Sigma/H)$ to $\Sigma/H$ sending $\bar{\Delta}\in\Sigma-(\Sigma/H)$ to $\bar{\Delta}+H\in\Sigma/H$ are bijective mappings preserving the inclusion relation between $\Sigma/H$ and $\Sigma-(\Sigma/H)$, and they are the inverse mappings of each other.

Furthermore, if $\Delta\in\Sigma/H$ and $\bar{\Delta}\in\Sigma-(\Sigma/H)$ correspond to each other by them, then $\dim\Delta=\dim\bar{\Delta}+1$.

$|\Sigma|=|\Sigma-(\Sigma/H)|\cup|\Sigma/H|^\circ$.
$|\Sigma-(\Sigma/H)|\cap|\Sigma/H|^\circ=\emptyset$.

$\Sigma\Mx\subset\Sigma/H$.
\end{enumerate}
\end{lemma}

\begin{lemma}
\label{construction}
Let $\Phi$ be any non-empty finite set whose elements are convex polyhedral cones in $V$ satisfying the following two conditions \emph{Z:}
\begin{description}
\item[\emph{(a)}]
$\Delta\cap\Lambda$ is a face of $\Delta$, and $\Delta\cap\Lambda$ is a face of $\Lambda$ for any $\Delta\in\Phi$ and any $\Lambda\in\Phi$.
\item[\emph{(b)}]
$\Delta\cap(-\Delta)=\Lambda\cap(-\Lambda)$ for any $\Delta\in\Phi$ and any $\Lambda\in\Phi$.
\end{description}

Choosing any element $\Delta\in\Phi$, we put $L=\Delta\cap(-\Delta)\subset V$.
$L$ does not depend on the choice of $\Delta\in\Phi$.
Put $\Sigma=\Phi\Fc$.
\begin{enumerate}
\item
$\Sigma$ is a fan in $V$.
\item
$\Sigma\supset\Phi$.
$|\Sigma|=|\Phi|$.
$\Sigma\Mx=\Phi\Mx$.
\item
$L\in\Sigma$.
$L$ is the minimum element of $\Sigma$.
\item
If $\Phi$ is simplicial, then $\Sigma$ is simplicial.
If $\Phi$ is rational over $N$, then $\Sigma$ is rational over $N$.
If $\Phi$ is regular over $N$, then $\Sigma$ is regular over $N$.
If $\Phi$ is flat, then $\Sigma$ is flat.
\end{enumerate}

Assume $\dim V\geq 2$. Let $S$ be any convex polyhedral cone in $V$ with $\dim S\geq 2$, let $m\in\Z_0$, let $H$ be any mapping from $\{1,2,\ldots,m\}$ to the set of all vector subspaces of $V$ of codimension one satisfying the following three conditions:
\begin{description}
\item[\emph{(c)}]
$H(i)\neq H(j)$ for any $i\in\{1,2,\ldots,m\}$ and any $j\in\{1,2,\ldots,m\}$ with $i\neq j$.
\item[\emph{(d)}]
$\Vect(S)\not\subset H(i)$ and $H(i)\cap S^\circ\neq\emptyset$ for any $i\in\{1,2,\ldots,m\}$.
\item[\emph{(e)}]
$H(i)\cap H(j)\cap S^\circ=\emptyset$ for any $i\in\{1,2,\ldots,m\}$ and any $j\in\{1,2,\ldots,m\}$ with $i\neq j$.
\end{description}
\begin{enumerate}
\setcounter{enumi}{4}
\item
The difference $S^\circ-(\cup_{i\in\{1,2,\ldots,m\}}H(i))$ is a non-empty open set of $\Vect(S)$.
It has $(m+1)$ connected components.
The closure of any connected component of it is a convex polyhedral cone in $V$ whose dimension is equal to $\dim S$.
\end{enumerate}

Let $\bar{\Phi}$ denote the finite set whose elements are $(m+1)$ of closures of connected components of $S^\circ-(\cup_{i\in\{1,2,\ldots,m\}}H(i))$.
Let $\bar{\Sigma}=\bar{\Phi}\Fc$, and let $\bar{L}=S\cap(-S)\cap(\cap_{i\in\{1,2,\ldots,m\}}H(i))$.
\begin{enumerate}
\setcounter{enumi}{5}
\item
$\bar{\Phi}$ satisfies the above two conditions \emph{Z}.
\item 
$\bar{\Sigma}$ is a flat fan in $V$.
$\dim\bar{\Sigma}=\dim S$.
$|\bar{\Sigma}|=S$.
$\bar{\Sigma}\Mx=\bar{\Sigma}^0=\bar{\Phi}$.
$\bar{L}$ is the miminum element of $\bar{\Sigma}$.
$\{\Delta\in \bar{\Sigma}^1|\Delta^\circ\subset S^\circ\}
=\{H(i)\cap S|i\in\{1,2,\ldots,m\}\}$.
$\sharp\bar{\Sigma}^0=\sharp\{\Delta\in \bar{\Sigma}^1|\Delta^\circ\subset S^\circ\}+1=m+1$.
For any $\Delta\in\bar{\Sigma}$ with $\dim\Delta\leq\dim V-2$, $\Delta\subset\partial S$.
If $S$ is rational over $N$ and $H(i)$ is rational over $N$ for any $i\in\{1,2,\ldots,m\}$, then $\bar{\Sigma}$ is rational over $N$.
\end{enumerate}
\end{lemma}

\begin{lemma}
Let $\Sigma$, $\Phi$ and $\Psi$ be fans in $V$.
\begin{enumerate}
\item
The following three conditions are equivalent:
\begin{enumerate}
\item
$\Sigma$ is a subdivision of $\Phi$, in other words, for any $\Delta\in\Sigma$ there exists $\Lambda\in\Phi$ with $\Delta\subset\Lambda$.
\item
For any $\Delta\in\Sigma$, there exists uniquely $\Lambda\in\Phi$ with $\Delta^\circ\subset\Lambda^\circ$.
\item
$|\Sigma|\subset|\Phi|$, and if $\Delta\in\Sigma$, $\Lambda\in\Phi$ and $\Delta^\circ\cap\Lambda^\circ\neq\emptyset$ then $\Delta^\circ\subset\Lambda^\circ$.
\end{enumerate}
\item
If $\Sigma$ is a subdivision of $\Phi$ and $\Phi$ is a subdivision of $\Sigma$, then $\Sigma=\Phi$.
\item
If $\Sigma$ is a subdivision of $\Phi$ and $\Sigma$ is a subdivision of $\Psi$, then $\Sigma$ is a subdivision of $\Psi$.
\item
If $\Sigma$ is a subdivision of $\Phi$, then $|\Sigma|\subset|\Phi|$.
\item
Assume that $\Sigma$ is a subdivision of $\Phi$.
For any $\Delta\in\Sigma$ and any $\Lambda\in\Phi$ the following three conditions are equivalent:
\begin{enumerate}
\item
$\Delta^\circ\subset\Lambda^\circ$.
\item
$\Delta\subset\Lambda$, and $\Phi/\Delta\subset\Phi/\Lambda$.
\item
$\Delta^\circ\cap\Lambda^\circ\neq\emptyset$
\end{enumerate}
\item
If $\Delta\in\Sigma$, $\Lambda\in\Phi$ and $\Delta^\circ\subset\Lambda^\circ$,
then $\Delta\subset\Lambda$, and $\dim\Delta\leq\dim\Lambda$.
\item
The following three conditions are equivalent:
\begin{enumerate}
\item
$\Sigma$ is a subdivision of $\Phi$ and $|\Sigma|=|\Phi|$.
\item
$|\Sigma|=|\Phi|$ and for any $\Lambda\in\Phi$, 
$\Lambda^\circ=\cup_{\Delta\in\Sigma, \Delta^\circ\subset\Lambda^\circ}\Delta^\circ$.
\item
$|\Sigma|=|\Phi|$ and $|\Sigma-\Sigma\Mx|\supset|\Phi-\Phi\Mx|$.
\end{enumerate}
\item
Assume that $\Sigma$ is a subdivision of $\Phi$ and $|\Sigma|=|\Phi|$.
For any $\Lambda\in\Phi$, there exists $\Delta\in\Sigma$ with $\Delta^\circ\subset\Lambda^\circ$.
For any $\Lambda\in\Phi\Mx$, there exists $\Delta\in\Sigma\Mx$ with $\Delta^\circ\subset\Lambda^\circ$.
\item
Assume that $\Sigma$ is a subdivision of $\Phi$, $|\Sigma|=|\Phi|$, $\Delta\in\Sigma$, $\Lambda\in\Phi$ and $\Delta^\circ\subset\Lambda^\circ$.
$\Delta\in\Sigma\Mx$, if and only if, $\Lambda\in\Phi\Mx$ and $\dim\Delta=\dim\Lambda$.
$\dim\Sigma=\dim\Phi$.
\item
If $\Sigma$ is a subdivision of $\Phi$, $|\Sigma|=|\Phi|$, 
 and $\Phi$ is flat, then $\Sigma$ is flat.
\end{enumerate}
\end{lemma}

\begin{lemma}
\label{scpc}
Let $\Sigma$ be any fan in $V$ such that the support $|\Sigma|$ of $\Sigma$ is a convex polyhedral cone in $V$.
By $L\in\Sigma$ we denote the minimum element of $\Sigma$.
\begin{enumerate}
\item
$\Sigma$ is a subdivision of $\mathcal{F}(|\Sigma|)$ and $|\Sigma|=|\mathcal{F}(|\Sigma|)|$.
\item
$\dim\Sigma=\dim |\Sigma|=\dim\Vect(|\Sigma|)$.
$\Sigma\Mx=\Sigma^0$. $\Sigma$ is flat.
\item
For any $\Delta\in\Sigma\Mx$, $\Delta^\circ\subset|\Sigma|^\circ$.
\item
For any $\Delta\in\Sigma$, $\Delta\not\subset\partial|\Sigma|$, if and only if, $\Delta^\circ\subset|\Sigma|^\circ$.
\item
Consider any $\Lambda\in\mathcal{F}(|\Sigma|)$.
$|\Sigma\backslash\Lambda|=\Lambda$.
$\Sigma\backslash\Lambda$ is a subdivision of $\mathcal{F}(\Lambda)$ and $|\Sigma\backslash\Lambda|=|\mathcal{F}(\Lambda)|$.
$\dim(\Sigma\backslash\Lambda)=\dim\Lambda$.
$(\Sigma\backslash\Lambda)\Mx=(\Sigma\backslash\Lambda)^0=
\{\Delta\cap\Lambda|\Delta\in\Sigma^0, \dim(\Delta\cap\Lambda)=\dim\Lambda\}$.
\item
Consider any $\Delta\in\Sigma$.
Take the unique $\Lambda\in\mathcal{F}(|\Sigma|)$ with $\Delta^\circ\subset\Lambda^\circ$.
Then, $\Delta\in\Sigma\backslash\Lambda$,
$(\Sigma\backslash\Lambda)\Mx/\Delta\neq\emptyset$,
$\Delta=\cap_{\bar{\Delta}\in(\Sigma\backslash\Lambda)\Mx/\Delta}\bar{\Delta}$, and $|\Sigma/\Delta|+\Vect(\Delta)=|\Sigma|+\Vect(\Delta)= |\Sigma|+\Vect(\Lambda)$.
\item
$L\subset|\Sigma|\cap(-|\Sigma|)$.
$\dim L\leq\dim(|\Sigma|\cap(-|\Sigma|))\leq\dim|\Sigma|$.
\item
Assume $\dim\Sigma-\dim L\geq 1$.
Consider any $\Delta\in\Sigma^1$.

$\Vect(\Delta)\subset\Vect(|\Sigma|)$, and
$\dim\Vect(|\Sigma|)=\dim\Vect(\Delta)+1$.
Let $H^{\circ\prime}$ and $H^{\circ\prime\prime}$ denote the two connected components of $\Vect(|\Sigma|)-\Vect(\Delta)$.
Let $H'=\Clos(H^{\circ\prime})$, and let $H''=\Clos(H^{\circ\prime\prime})$.
$H'\cup H''=\Vect(|\Sigma|)$.
$H'\cap H''=\Vect(\Delta)$.

We consider the case where $\Delta\not\subset\partial|\Sigma|$.
$\sharp(\Sigma^0/\Delta)=2$.
Let $\Delta'$ and $\Delta''$ denote the two elements of $\Sigma^0/\Delta$.
$\{\Delta'+\Vect(\Delta), \Delta''+\Vect(\Delta)\}=\{H', H''\}$, and $\Delta'+\Vect(\Delta)\neq\Delta''+\Vect(\Delta)$.

We consider the case where $\Delta\subset\partial|\Sigma|$. $\sharp(\Sigma^0/\Delta)=1$.
Let $\Delta'$ denote the unique element of $\Sigma^0/\Delta$.
$\Delta'+\Vect(\Delta)=|\Sigma|+\Vect(\Delta)$.
$|\Sigma|+\Vect(\Delta)=H'$ or $|\Sigma|+\Vect(\Delta)=H''$.
\end{enumerate}
\end{lemma}

\begin{lemma}
\label{ris}
Let $J$ be any finite set and let $\Phi$ be any mapping from $J$ to the set of all finite sets whose elements are convex polyhedral cones in $V$.
For any $j\in J$, $\Phi(j)$ is a finite set whose elements are convex polyhedral cones in $V$.
We consider the real intersection $\hat{\cap}_{j\in J}\Phi(j)$ of $\Phi(j), j\in J$.
By definition
\begin{equation*}
\begin{split}
\bigcap_{j\in J}^\wedge\Phi(j)
=\{\bar{\Delta}\in 2^V|&\bar{\Delta}=\cap_{j\in J}\Delta(j)\text{ for some mapping }\Delta:J\rightarrow 2^V\\
&\quad\text{such that }\Delta(j)\in\Phi(j)\text{ for any }j\in J\}\subset 2^V.
\end{split}
\end{equation*}
\begin{enumerate}
\item
$\hat{\cap}_{j\in J}\Phi(j)$ is a finite set whose elements are convex polyhedral cones in $V$.
\item
$\hat{\cap}_{j\in J}\Phi(j)$ is a subdivision of $\Phi(j)$ for any $j\in J$.

Let $\Sigma$ be any finite set whose elements are convex polyhedral cones in $V$.
If $\Sigma$ is a subdivision of $\Phi(j)$ for any $j\in J$, then $\Sigma$ is a subdivision of $\hat{\cap}_{j\in J}\Phi(j)$.
\item
$$|\bigcap_{j\in J}^\wedge\Phi(j)|=\bigcap_{j\in J}|\Phi(j)|.$$
\item
If $J=\emptyset$, then $\hat{\cap}_{j\in J}\Phi(j)=\{V\}$.
$\hat{\cap}_{j\in J}\Phi(j)=\emptyset$, if and only if, $J\neq\emptyset$ and $\Phi(j)=\emptyset$ for some $j\in J$.
\item
If $\Phi(j)$ is a fan for any $j\in J$, then $\hat{\cap}_{j\in J}\Phi(j)$ is also a fan.
\item
If $\Phi(j)$ is rational over $N$ for any $j\in J$, then $\hat{\cap}_{j\in J}\Phi(j)$ is also rational over $N$.
\item
For any subsets $J',J''$ of $J$ with $J'\cup J''=J$ and $J'\cap J''=\emptyset$,
$$\bigcap_{j\in J}^\wedge\Phi(j)
=(\bigcap_{j\in J'}^\wedge\Phi(j))\hat{\cap}(\bigcap_{j\in J''}^\wedge\Phi(j)).$$
\item
For any bijective mapping $\sigma:J\rightarrow J$
$$\bigcap_{j\in J}^\wedge\Phi(\sigma(j))= \bigcap_{j\in J}^\wedge\Phi(j).$$
\end{enumerate}
\end{lemma}

\begin{lemma}
\label{several fans}
Let $m\in \Z_+$ be any positive integer, and let $\Sigma$ be any mapping from the set $\{1,2,\ldots,m\}$ to the set of all fans in $V$.
For any $i\in\{1,2,\ldots,m\}$, $\Sigma(i)$ is a fan in $V$.

We denote
$$\bar{\Sigma}=\bigcap_{i\in\{1,2,\ldots,m\}}^\wedge\Sigma(i)\subset 2^V.$$

\begin{enumerate}
\item
$\bar{\Sigma}$ is a fan in $V$.
$|\bar{\Sigma}|=\cap_{i\in\{1,2,\ldots,m\}}|\Sigma(i)|$.
If $\Sigma(i)$ is rational over $N$ for any $i\in\{1,2,\ldots,m\}$, then $\bar{\Sigma}$ is rational over $N$.
\item
$\bar{\Sigma}$ is a subdivision of $\Sigma(i)$ for any $i\in\{1,2,\ldots,m\}$.

Let $\bar{\Phi}$ be any fan in $V$.
If $\bar{\Phi}$ is a subdivision of $\Sigma(i)$ for any $i\in\{1,2,\ldots,m\}$, then $\bar{\Phi}$ is a subdivision of $\bar{\Sigma}$.
\item
Let $L(j)\in\Sigma(j)$ be the minimum element of $\Sigma(j)$ for any $j\in J$. $\cap_{j\in J}L(j)$ is the minimum element of $\bar{\Sigma}$.
\end{enumerate}

Let $\bar{\Delta}$ be any element of $\bar{\Sigma}$.
\begin{enumerate}\setcounter{enumi}{3}
\item
There exists uniquely an element $\Delta(i)\in\Sigma(i)$ with $\bar{\Delta}^\circ\subset \Delta(i)^\circ$ for any $i\in\{1,2,\ldots,m\}$.
\end{enumerate}

Below, we assume that $\Delta(i)\in\Sigma(i)$ and $\bar{\Delta}^\circ\subset \Delta(i)^\circ$ for any $i\in\{1,2,\ldots,m\}$.

\begin{enumerate}
\setcounter{enumi}{4}
\item
$\bar{\Delta}\subset \Delta(i)$ for any $i\in\{1,2,\ldots,m\}$.
\item
$$\bar{\Delta}=\bigcap_{i\in\{1,2,\ldots,m\}}\Delta(i).$$
\item
$$\bar{\Delta}^\circ=\bigcap_{i\in\{1,2,\ldots,m\}}\Delta(i)^\circ.$$
\item
$$\Vect(\bar{\Delta})=\bigcap_{i\in\{1,2,\ldots,m\}}\Vect(\Delta(i)).$$
\item
Consider any subset $\bar{\Lambda}$ of $V$.
$\bar{\Lambda}$ is a face of $\bar{\Delta}$, if and only if, there exists a mapping $\Lambda: \{1,2,\ldots,m\}\rightarrow 2^V$ such that $\bar{\Lambda}=\cap_{i\in \{1,2,\ldots,m\}}\Lambda(i)$ and $\Lambda(i)$ is a face of $\Delta(i)$ for any $i\in\{1,2,\ldots,m\}$.
\end{enumerate}

Let $\Lambda(i)$ be any element of $\Sigma(i)$ for any $i\in\{1,2,\ldots,m\}$.
\begin{enumerate}
\setcounter{enumi}{9}
\item
The intersection $\cap_{i\in\{1,2,\ldots,m\}}\Lambda(i)$ is an element of $\bar{\Sigma}$.
\item
If $\bar{\Delta}\subset \Lambda(i)$ then $\Delta(i)\subset \Lambda(i)$, for any $i\in\{1,2,\ldots,m\}$.
\item
If
$$\bar{\Delta}=\bigcap_{i\in\{1,2,\ldots,m\}}\Lambda(i),$$
then $\Delta(i)\subset \Lambda(i)$ for any $i\in\{1,2,\ldots,m\}$ and the following three conditions are equivalent:
\begin{enumerate}
\item
$$\bigcap_{i\in\{1,2,\ldots,m\}}\Lambda(i)^\circ\neq\emptyset.$$
\item
$\Delta(i)=\Lambda(i)$ for any $i\in\{1,2,\ldots,m\}$.
\item
$$\bar{\Delta}^\circ=\bigcap_{i\in\{1,2,\ldots,m\}}\Lambda(i)^\circ$$
\end{enumerate}
\end{enumerate}
\end{lemma}

\begin{lemma}
\label{pbcpcd}
Let $\Sigma$ be any fan in $V$, let $U$ be any vector space of finite dimension over $\R$, and let $\nu:U\rightarrow V$ be any homomorphism of vector spaces over $\R$.

The pull back $\nu^*\Sigma$ of $\Sigma$ by $\nu$ is a fan in $U$.
\end{lemma}

\section{Convex pseudo polytopes}
\label{cpp}

We study convex pseudo polytopes.

Let $V$ be any vector space of finite dimension over $\R$, and let $N$ be any lattice of $V$.

\begin{lemma}
\label{cpp1}
Let $S$ be any convex pseudo polytope in $V$, and let $X$ and $Y$ be any finite subset of $V$ satisfying $S=\Conv(X)+\Convcone(Y)$ and $X\neq\emptyset$.
\begin{enumerate}
\item
$\Stab(S)=\Convcone(Y)$.
$\Stab(S)$ is a convex polyhedral cone in $V$.
If $S$ is rational over $N$, then $\Stab(S)$ is also rational over $N$.
\item
$\Vect(\Stab(S))\subset\Stab(\Affi(S))$.
\item
The following four conditions are equivalent:
\begin{enumerate}
\item
$S$ is a convex polytope.
\item
$S=\Conv(X)$.
\item
$\Stab(S)=\{0\}$.
\item
$S$ is compact.
\end{enumerate}
\item
We consider the dual vector space $V^*$ of $V$ and the dual cone $\Stab(S)^\vee=\Stab(S)^\vee|V\subset V^*$ of $\Stab(S)$.

For any $\omega\in V^*$, the following three conditions are equivalent:
\begin{enumerate}
\item
$\omega\in\Stab(S)^\vee$.
\item
There exists the minimum element $\min\{\langle\omega, x\rangle|x\in S\}$ of the subset $\{\langle\omega, x\rangle|$\break$x\in S\}$ of $\R$.
\item
The subset $\{\langle\omega, x\rangle|x\in S\}$ of $\R$ is bounded below.
\end{enumerate}
\end{enumerate}
\end{lemma}

\begin{definition}
\label{faces of cpp}
Let $S$ be any convex pseudo polytope in $V$.
We consider the dual cone $\Stab(S)^\vee=\Stab(S)^\vee|V\subset V^*$ of $\Stab(S)$.
\begin{enumerate}
\item
For any $\omega\in \Stab(S)^\vee$, we denote
\begin{equation*}
\begin{split}
\Ord(\omega,S|V)=&\min\{\langle\omega,x\rangle|x\in S\}\in\R\\
\Delta(\omega,S|V)=&\{x\in S|\langle\omega, x\rangle=\Ord(\omega,S|V)\}\subset S.
\end{split}
\end{equation*}

When we need not refer to $V$ or to the pair $(S, V)$, we also write simply  $\Ord(\omega, S)$ or $\Ord(\omega)$, $\Delta(\omega, S)$ or $\Delta(\omega)$ respectively, instead of $\Ord(\omega,S|V)$, $\Delta(\omega, S|V)$.
\item
Let $F$ be any subset of $S$.
We say that $F$ is a \emph{face} of $S$, if $F=\Delta(\omega, S|V)$ for some $\omega\in \Stab(S)^\vee.$

It is easy to see that any face $F$ of $S$ is a closed convex subset of $V$, and the dimension $\dim F\in\Z_0$ of $F$, the boundary $\partial F$ of $F$, and the interior $F^\circ$ of $F$ are defined.

Any face $F$ of $S$ with $\dim F=0$ is called a \emph{vertex} of $S$.
Any vertex of $S$ is a subset of $S$ with only one element.
Any face $F$ of $S$ with $\dim F=1$ is called an \emph{edge} of $S$.
Any face $F$ of $S$ with $\dim F=\dim S-1$ is called an \emph{facet} of $S$.
Any face $F$ of $S$ with $F\neq S$ is called a \emph{proper} face of $S$.
\item
By $\mathcal{F}(S)$ we denote the set of all faces of $S$.

For any $i\in\Z$, the set of all faces $F$ with $\dim F=i$ is denoted by $\mathcal{F} (S)_i$, and the set of all faces $F$ with $\dim F=\dim S-i$ is denoted by $\mathcal{F}(S)^i$.
\item
Let $\ell=\dim(\Stab(S)\cap(-\Stab(S)))\in\Z_0$.
We denote
\begin{equation*}
\begin{split}
c(S)=&\sharp\mathcal{F}(S)_\ell\in\Z_0,\\
\mathcal{V}(S)=&\bigcup_{F\in\mathcal{F}(S)_\ell}F\subset S,
\end{split}
\end{equation*}
we call $c(S)$ the \emph{characteristic number} of $S$, and we call $\mathcal{V}(S)$ the \emph{skeleton} of $S$.
We call any face $F$ of $S$ with $\dim F=\ell$ a \emph{minimal} face of $S$.
\item
Let $F$ be any face of $S$.
We denote
\begin{equation*}
\begin{split}
\Delta^\circ(F,S|V)=&\{\omega\in \Stab(S)^\vee|F=\Delta(\omega, S|V)\}\subset \Stab(S)^\vee\subset V^*,\\
\Delta(F,S|V)=&\{\omega\in \Stab(S)^\vee|F\subset\Delta(\omega, S|V)\}\subset \Stab(S)^\vee\subset V^*.
\end{split}
\end{equation*}

We call $\Delta^\circ(F,S|V)$ the \emph{open normal cone} of $F$, and we call $\Delta (F,S|V)$ the \emph{normal cone} of $F$.

When we need not refer to $V$ or to the pair $(S, V)$, we also write simply  $\Delta^\circ (F, S)$ or $\Delta^\circ(F)$, $\Delta(F, S)$ or $\Delta(F)$ respectively, instead of $\Delta^\circ(F, S|V)$, $\Delta(F, S|V)$.
\item
We denote
$$\Sigma(S|V)=\{\Delta(F, S|V)|F\in\mathcal{F}(S)\}\subset 2^{\Stab(S)^\vee}\subset 2^{V^*},$$
and we call $\Sigma(S|V)$ the \emph{normal fan} of $S$.

When we need not refer to $V$, we also write simply $\Sigma(S)$, instead of $\Sigma(S|V)$.
\end{enumerate}
\end{definition}

Let $W$ be any vector space of finite dimension over $\R$ containing $V$ as a vector subspace over $\R$ with $\dim W=\dim V+1$, and let $z\in W- V$ be any point.

Let $\pi:V\rightarrow W$ denote the inclusion homomorphism.
Putting $\pi'(t)=tz\in W$ for any $t\in\R$ we define an injective homomorphism $\pi':\R\rightarrow W$ of vector spaces over $\R$.

For any $a\in W$, choosing the unique pair $b\in V$ and $t\in\R$ with $a=b+tz$ and putting $\rho(a)=b$ and $\rho'(a)=t$ we define homomorphisms $\rho:W\rightarrow V$ and $\rho':W\rightarrow \R$ of vector spaces over $\R$.

Putting $\iota(\omega)=\omega(1)\in\R$ for any $\omega\in \R^*=\mathrm{Hom}_{\R}(\R,\R)$, we define an isomorphism $\iota:\R^*\rightarrow \R$ of vector spaces over $\R$.
For any $\omega\in\R^*$ and any $t\in\R$ we have $\langle\omega, t\rangle=\iota(\omega)t$.
Below, using this isomorphism $\iota$ we identify $\R^*$ with $\R$.
For any $t\in\R=\R^*$ and any $u\in\R$ we have $\langle t,u\rangle=tu$.

We have eight homomorphisms of vector spaces over $\R$.
\begin{gather*}
\begin{matrix}
\pi:&V\rightarrow W,&\quad& \pi'&:\R\rightarrow W,\\
\rho:&W\rightarrow V,&\quad& \rho'&:W\rightarrow\R,\\
\pi^*:&W^*\rightarrow V^*,&\quad& \pi^{\prime*}&:W^*\rightarrow \R,\\
\rho^*:&V^*\rightarrow W^*,&\quad& \rho^{\prime*}&:\R\rightarrow W^*.
\end{matrix}
\end{gather*}
Four homomorphisms $\pi, \pi', \rho^*, \rho^{\prime*}$ are injective.
The other four $\rho, \rho', \pi^*, \pi^{\prime*}$ are surjective.
We denote
$H=V+\R_0z\subset W$ and $\zeta=\rho^{\prime*}(1)\in W^*$.

\begin{lemma}
\label{eight mappings}
\begin{enumerate}
\item
$\rho\pi=\Id_V$, $\rho'\pi'=\Id_{\R}$, $\pi\rho+\pi'\rho'=\Id_W$,
$V=\pi(V)= \rho^{\prime-1}(0)$, $\R z=\pi'(\R)=\rho^{-1}(0)$, $\pi'(1)=z$.
For any $x\in W$, $\rho'(x)=\langle\zeta,x\rangle$.
\item
$\pi^*\rho^*=\Id_{V^*}$, $\pi^{\prime*}\rho^{\prime*}=\Id_{\R}$, $\rho^*\pi^*+\rho^{\prime*}\pi^{\prime*}=\Id_{W^*}$,
$(\R z)^\vee =\rho^*(V^*)= \pi^{\prime*-1}(0)$, $V^\vee =\R\zeta=\rho^{\prime*}(\R)=\pi^{*-1}(0)$.
For any $\xi\in W^*$, $\pi^{\prime*}(\xi)=\langle\xi,z\rangle$. $\langle\zeta,z\rangle=1$.
\item
$N+\Z x$ is a lattice of $W$. $H$ is a rational convex polyhedral cone over $N+\Z x$ in $W$. 
$\dim H=\dim W=\dim V+1$.
$H=\{x\in W|\rho'(x)\geq 0\}$.
$W=H\cup(-H)=\Vect(H)=\Vect(-H)$.
$V=H\cap(-H)=\partial H=\partial(-H)$.
$z\in H^\circ=V+\R_+z=\{x\in W|\rho'(x)> 0\}$.
\item
$(N+\Z z)^*=\rho^*(N^*)+\Z\zeta$.
$H^\vee$ is a simplicial cone over $(N+\Z z)^*$ in $W^*$ with $\dim H^\vee=1$.
$H^\vee\cap(N+\Z z)^*=\Z_0\zeta$.
$H^\vee=\R_0\zeta$.
$V^\vee =H^\vee\cup(-H^\vee)=\Vect(H^\vee) =\Vect(-H^\vee)$.
$\{0\}=H^\vee\cap(-H^\vee)=\partial H^\vee =\partial(-H^\vee)$.
$-H^\vee=(-H)^\vee$.
\end{enumerate}
\end{lemma}

Putting $\sigma(a)=\rho(a)/\rho'(a)\in V$ for any $a\in H^\circ$, we define a mapping $\sigma:H^\circ\rightarrow V$.
If $a\in H^\circ$, $b\in V$, $t\in\R$ and $a=b+tz$, then $t>0$ and $\sigma(a)=b/t$.

\begin{lemma}
\label{first correspondence}
\begin{enumerate}
\item
$\sigma$ is surjective.
For any $a\in H^\circ$, $\{\sigma(a)+z\}=(\R_+a)\cap(V+\{z\})=(\R a)\cap(V+\{z\})$ and
$\sigma^{-1}(\sigma(a))=\R_+a$.
For any $b\in V$, $\sigma^{-1}(b)=\R_+(b+z)$.
\item
Consider any convex polyhedral cone $\Delta$ in $W$ satisfying $\Delta\subset H$ and $\Delta\cap H^\circ\neq\emptyset$ and any finite subset $Z$ of $H$ satisfying $\Delta=\Convcone(Z)$.
\begin{enumerate}
\item
$\sigma(\Delta\cap H^\circ)=\Conv(\sigma(Z\cap H^\circ))+\Convcone(Z\cap V)$.

$\sigma(\Delta\cap H^\circ)$ is a convex pseudo polytope in $V$.

If $\Delta$ is rational over $N+\Z z$, then $\sigma(\Delta\cap H^\circ)$ is rational over $N$.

$\dim\Delta=\dim \sigma(\Delta\cap H^\circ)+1$.
$\Delta\cap V=\Stab(\sigma(\Delta\cap H^\circ))$.
\item
$\Delta\cap H^\circ=\sigma^{-1}(\sigma(\Delta\cap H^\circ))$.
$\Delta=\Clos(\sigma^{-1}(\sigma(\Delta\cap H^\circ)))$.
\end{enumerate}
\item
Consider any convex pseudo polytope $S$ in $V$ and any finite subsets $X$, $Y$ of $V$ satisfying $S=\Conv(X)+\Convcone(Y)$ and $X\neq\emptyset$.
\begin{enumerate}
\item
$\Clos(\sigma^{-1}(S))=\Convcone((X+\{z\})\cup Y)$.

$\Clos(\sigma^{-1}(S))$ is a convex polyhedral cone in $W$.
$\Clos(\sigma^{-1}(S))\subset H$.
$\Clos(\sigma^{-1}(S))\cap H^\circ\neq\emptyset$.
$\Clos(\sigma^{-1}(S))\cap V=\Stab(S)$.

If $S$ is rational over $N$, then $\Clos(\sigma^{-1}(S))$ is rational over $N+\Z z$.
\item
$\Clos(\sigma^{-1}(S))\cap H^\circ=\sigma^{-1}(S)$.
$\sigma(\Clos(\sigma^{-1}(S))\cap H^\circ)=S$.
\end{enumerate}
\item
For any subsets $\Delta$, $\Lambda$ of $H$ satisfying $\Delta\cap H^\circ\neq\emptyset$, $\Lambda\cap H^\circ\neq\emptyset$ and $\Delta\subset\Lambda$, $\sigma(\Delta\cap H^\circ)\subset\sigma(\Lambda\cap H^\circ)$.

For any non-empty subsets $S$, $T$ of $V$ satisfying $S\subset T$, $\Clos(\sigma^{-1}(S))\subset\Clos(\sigma^{-1}(T))$.

For any non-empty closed subsets $S$, $T$ of $V$, $S\cap T=\emptyset$, if and only if, $\Clos(\sigma^{-1}(S))\cap\Clos(\sigma^{-1}(T))\cap H^\circ=\emptyset$.
\item
For any convex polyhedral cone $\Delta$ in $W$ satisfying $\Delta\subset H$ and $\Delta\cap H^\circ\neq\emptyset$, $\sigma(\Delta\cap H^\circ)$ is a convex pseudo polytope in $V$.

For any convex pseudo polytope $S$ in $V$, $\Clos(\sigma^{-1}(S))$ is a convex polyhedral cone in $W$, $\Clos(\sigma^{-1}(S))\subset H$ and $\Clos(\sigma^{-1}(S))\cap H^\circ\neq\emptyset$.

The mapping sending any convex polyhedral cone $\Delta$ in $W$ satisfying $\Delta\subset H$ and $\Delta\cap H^\circ\neq\emptyset$ to $\sigma(\Delta\cap H^\circ)$ and the mapping sending any convex pseudo polytope $S$ in $V$ to $\Clos(\sigma^{-1}(S))$ are bijective mappings preserving the inclusion relation between the set of all convex polyhedral cones $\Delta$ in $W$ satisfying $\Delta\subset H$ and $\Delta\cap H^\circ\neq\emptyset$ and the set of all convex pseudo polytopes in $V$, and they are the inverse mappings of each other.

Furthermore, if a convex polyhedral cone $\Delta$ in $W$ satisfying $\Delta\subset H$ and $\Delta\cap H^\circ\neq\emptyset$ and a convex pseudo polytope $S$ in $V$ correspond to each other by them, then $\Delta=\Clos(\sigma^{-1}(S))$, $S=\sigma(\Delta\cap H^\circ)$ and the following claims holds:

$\Delta$ is rational over $N+\Z z$, if and only if, $S$ is rational over $N$.
\begin{equation*}
\begin{split}
\dim\Delta&\geq 1.\\
\dim\Delta&=\dim S+1.\\
\Delta\cap V&=\Stab(S).\\
\Delta\cap H^\circ&=\sigma^{-1}(S).\\
\Delta\cap(V+\{z\})&= S+\{z\}.\\
\Vect(\Delta)\cap(V+\{z\})&=\Affi(S)+\{z\}.\\
\Vect(\Delta)\cap H^\circ&=\sigma^{-1}(\Affi(S)).\\
\sigma(\Vect(\Delta)\cap H^\circ)&= \Affi(S).\\
\Vect(\Delta)&=\Vect(\Affi(S)+\{z\}).\\
\partial\Delta\cap H^\circ&=\sigma^{-1}(\partial S).\\
\Delta^\circ&=\sigma^{-1}(S^\circ).\\
\sigma(\Delta^\circ)&=S^\circ.
\end{split}
\end{equation*}
\end{enumerate}
\end{lemma}

\begin{cor}
\label{intersection}
Consider any convex pseudo polytopes $S$, $T$ in $V$.

If $S\cap T\neq\emptyset$, then $S\cap T$ is a convex pseudo polytope in $V$ and $\Stab(S\cap T)=\Stab(S)\cap\Stab(T)$.

If $S$ and $T$ are rational over $N$ and $S\cap T\neq\emptyset$, then $S\cap T$ is also rational over $N$.

If $S$ and $T$ are convex polytopes and $S\cap T\neq\emptyset$, then $S\cap T$ is also a convex polytope.

If $S$ and $T$ are convex polyhedral cones, then $S\cap T$ is also a convex polyhedral cone.

\end{cor}

\begin{prop}
\label{miracle}
Let $\Delta$ be any convex polyhedral cone in $W$ satisfying $\Delta\subset H$ and $\Delta\cap H^\circ\neq\emptyset$.
We denote
\begin{equation*}
\begin{split}
L&=\Delta\cap(-\Delta)\subset\Delta,\\
\ell&=\dim L\in\Z_0,\\
S&=\sigma(\Delta\cap H^\circ)\subset V,\\
H^\vee&=H^\vee|W\subset W^*,\\
\Delta^\vee&=\Delta^\vee|W\subset W^*,\\
\partial_-\Delta^\vee&=\{\omega\in\Delta^\vee|
(\{\omega\}+\Vect(H^\vee))\cap\Delta^\vee\subset\{\omega\}+H^\vee\}\subset \Delta^\vee,\\
\mathcal{F}(\Delta)_*&=\{\Lambda\in\mathcal{F}(\Delta)|\Lambda\cap H^\circ\neq\emptyset\}\subset\mathcal{F}(\Delta),\\
\mathcal{F}(\Delta^\vee)^*&=\mathcal{F}(\Delta^\vee)\backslash\partial_-\Delta^\vee\subset\mathcal{F}(\Delta^\vee).
\end{split}
\end{equation*}
\begin{enumerate}
\item
$L=(\Delta\cap V)\cap(-(\Delta\cap V))=\Stab(S)\cap(-\Stab(S))\subset V$.
$L\in\mathcal{F}(\Delta)-\mathcal{F}(\Delta)_*$.
$\Delta\in\mathcal{F}(\Delta)_*\neq\emptyset$.
If $\Lambda\in\mathcal{F}(\Delta)_*$, $\Gamma\in\mathcal{F}(\Delta)$ and $\Lambda\subset\Gamma$, then $\Gamma\in\mathcal{F}(\Delta)_*$.
For any $\Lambda\in\mathcal{F}(\Delta)_*$ there exists $\Gamma\in\mathcal{F}(\Delta)_*$ such that $\Lambda\supset\Gamma$ and $\dim\Gamma=\ell+1$.
\item
$H^\vee\subset\Delta^\vee$.
$-H^\vee\not\subset\Delta^\vee$.
$\Delta^\vee+\pi^{*-1}(0)=(\Delta\cap V)^\vee|W=\pi^{*-1}((\Delta\cap V)^\vee|V)$.
$\pi^*(\Delta^\vee)=(\Delta\cap V)^\vee|V=\Stab(S)^\vee|V$.
\item
$\emptyset\neq\partial_-\Delta^\vee\subset\Delta^\vee=\partial_-\Delta^\vee+H^\vee$.
$\pi^*(\partial_-\Delta^\vee)=\Stab(S)^\vee|V$.
The mapping $\pi^*:\partial_-\Delta^\vee\rightarrow\Stab(S)^\vee|V$ induced by $\pi^*$ is bijective.
\item
For any face $\Lambda$ of $\Delta$ the following three conditions are equivalent:
\begin{enumerate}
\item
$\Lambda\in\mathcal{F}(\Delta)_*$.
\item
$\Delta(\Lambda,\Delta|W)\in\mathcal{F}(\Delta^\vee)^*$.
\item
$\Delta^\circ(\Lambda,\Delta|W)\cap\partial_-\Delta^\vee\neq\emptyset$.
\end{enumerate}
\item
$|\mathcal{F}(\Delta^\vee)^*|=\partial_-\Delta^\vee$.
$\Delta^\vee\cap(-\Delta^\vee)\in\mathcal{F}(\Delta^\vee)^*\neq\emptyset$.
$\dim\mathcal{F}(\Delta^\vee)^* =\dim\Delta^\vee-1=\dim\Stab(S)^\vee=\dim V-\ell$.
$(\mathcal{F}(\Delta^\vee)^*)\Mx=(\mathcal{F}(\Delta^\vee)^*)^0$.
\item
For any $\omega\in \partial_-\Delta^\vee$ the following claims hold:
\begin{enumerate}
\item
$\Delta(\omega,\Delta|W)\in\mathcal{F}(\Delta)_*$.
\item
$\Ord(\pi^*(\omega),S|V)=-\langle\omega,z\rangle$.
\item
$\Delta(\omega,\Delta|W)\cap H^\circ=\sigma^{-1}(\Delta(\pi^*(\omega),S|V))$.
\item
$\sigma(\Delta(\omega,\Delta|W)\cap H^\circ)=\Delta(\pi^*(\omega),S|V)$.
\end{enumerate}
\item
For any $\Lambda\in\mathcal{F}(\Delta)_*$, $\sigma(\Lambda\cap H^\circ)\in\mathcal{F}(S)$.

For any $F\in\mathcal{F}(S)$, $\Clos(\sigma^{-1}(F))\in \mathcal{F}(\Delta)_*$.

The mapping from $\mathcal{F}(\Delta)_*$ to $\mathcal{F}(S)$ sending $\Lambda\in\mathcal{F}(\Delta)_*$ to $\sigma(\Lambda\cap H^\circ)\in\mathcal{F}(S)$ and the mapping from $\mathcal{F}(S)$ to $\mathcal{F}(\Delta)_*$ sending $F\in\mathcal{F}(S)$ to $\Clos(\sigma^{-1}(F))\in \mathcal{F}(\Delta)_*$ are bijective mappings preserving the inclusion relation between $\mathcal{F}(\Delta)_*$ and $\mathcal{F}(S)$, and they are the inverse mappings of each other.
\item
Assume $\Lambda\in\mathcal{F}(\Delta)_*$, $F\in\mathcal{F}(S)$, $F=\sigma(\Lambda\cap H^\circ)$ and $\Lambda=\Clos(\sigma^{-1}(F))$.
The following claims hold:
\begin{enumerate}
\item
$F$ is a convex pseudo polytope in $V$.
If $S$ is rational over $N$, then $F$ is also rational over $N$.
$\Stab(F)$ is a face of $\Stab(S)$.
\item
$\dim\Lambda=\dim F+1$.
\item
$\Delta^\circ(\Lambda,\Delta|W)\subset\Delta(\Lambda,\Delta|W)\subset\partial_-\Delta^\vee$.
\item
$\Delta(F,S|V)=\pi^*(\Delta(\Lambda,\Delta|W))$.
$\pi^{*-1}(\Delta(F,S|V))\cap\partial_-\Delta^\vee= \Delta(\Lambda,\Delta|W)$.
$\Delta(F,S|V)$ is a convex polyhedral cone in $V^*$.
\item
$\Delta^\circ(F,S|V)=\pi^*(\Delta^\circ(\Lambda,\Delta|W))$.
$\pi^{*-1}(\Delta^\circ(F,S|V))\cap\partial_-\Delta^\vee=\Delta^\circ(\Lambda,\Delta|W)$.
$\Delta^\circ(F,S|V)=\Delta(F,S|V) ^\circ$.
\item
$\Vect(\Delta(\Lambda,\Delta|W))\cap\pi^{*-1}(0)=\{0\}$.
$\Vect(\Delta(F,S|V))=$\hfill\break$\pi^*(\Vect(\Delta(\Lambda,\Delta|W)))=\Stab(\Affi(F))^\vee|V$.
\end{enumerate}
\item
For any $\Lambda^*\in\mathcal{F}(\Delta^\vee)^*$, $\pi^*(\Lambda^*)\in\Sigma(S|V)$.

For any $\bar{\Lambda}^*\in\Sigma(S|V)$, $\pi^{*-1}(\bar{\Lambda}^*)\cap \partial_-\Delta^\vee\in \mathcal{F}(\Delta^\vee)^*$.

The mapping from $\mathcal{F}(\Delta^\vee)^*$ to $\Sigma(S|V)$ sending $\Lambda^*\in\mathcal{F}(\Delta^\vee)^*$ to $\pi^*(\Lambda^*)\in\Sigma(S|V)$ and the mapping from $\Sigma(S|V)$ to $\mathcal{F}(\Delta^\vee)^*$ sending $\bar{\Lambda}^*\in\Sigma(S|V)$ to $\pi^{*-1}(\bar{\Lambda}^*)\cap \partial_-\Delta^\vee\in \mathcal{F}(\Delta^\vee)^*$ are bijective mappings preserving the dimension and the inclusion relation between $\mathcal{F}(\Delta^\vee)^*$ and $\Sigma(S|V)$, and they are the inverse mappings of each other.
\item
The normal fan $\Sigma(S|V)$ of $S$ in $V$ is a fan in $V^*$.
$|\Sigma(S|V)|=\Stab(S)^\vee|V$.
If $S$ is rational over $N$, then $\Sigma(S|V)$ is rational over $N^*$.
\item
For any $F\in\mathcal{F}(S)$, $\Vect(\Delta(F,S|V))=\Stab(\Affi(F))^\vee|V$ and $\dim F+\dim \Delta(F,$\hfill\break$S|V)=\dim V$.

For any $F\in\mathcal{F}(S)$ and $G\in\mathcal{F}(S)$, $F\subset G$, if and only if, $\Delta(F,S|V)\supset\Delta(G,S|V)$.

The mapping from $\mathcal{F}(S)$ to $\Sigma(S|V)$ sending $F\in\mathcal{F}(S)$ to $\Delta(F,S|V)\in\Sigma(S|V)$ is a bijective mapping.
\item
The function $\Ord(\;,S|V):\Stab(S)^\vee|V\rightarrow\R$ sending $\bar{\omega}\in\Stab(S)^\vee|V$ to $\Ord(\bar{\omega},S|V)\in\R$ is a piecewise linear convex function over $\Sigma(S|V)$.

If $S$ is rational over $N$, then this function $\Ord(\;,S|V)$ is rational over $N^*$.
\item
Denote
\begin{equation*}\begin{split}
\Xi(\Sigma(S|V), \Ord(\;,S|V))=
\{\xi \in W^*|&\xi=\rho(\bar{\omega})+t\zeta\text{ for some }\bar{\omega}\in |\Sigma(S|V)|\\
&\quad\text{and some }t\in\R
\text{ with }t\geq -\Ord(\bar{\omega},S|V)\}.
\end{split}\end{equation*}
Then, $\Xi(\Sigma(S|V), \Ord(\;,S|V))=\Delta^\vee$,
$\Xi(\Sigma(S|V), \Ord(\;,S|V))^\vee|W^*=\Delta$, and
$\sigma((\Xi(\Sigma(S|V), \Ord(\;,S|V))^\vee|W^*)\cap H^\circ)=S$.
\end{enumerate}
\end{prop}

\begin{proof}
We give only the proof of 4.

Since $L\subset V=\{x\in W|\langle\zeta,x\rangle=0\}$, $\langle\zeta,b\rangle=0$ for any $b\in L$.

Consider any face $\Gamma$ of $\Delta$ with $\dim\Gamma=\ell+1$.
We take any point $e_\Gamma\in\Gamma- L$.
We have $\Gamma=\R_0e_\Gamma+L$.
Since $e_\Gamma\in\Gamma\subset\Delta\subset H=\{x\in W|\langle\zeta,x\rangle\geq 0\}$, $\langle\zeta,e_\Gamma\rangle\geq 0$.

We know $\Delta=\sum_{\Gamma\in\mathcal{F}(\Delta)_{\ell+1}}\R_0e_\Gamma+L$.
Take any point $a\in\Delta\cap H^\circ\neq\emptyset$.
$a\in\Delta=\sum_{\Gamma\in\mathcal{F}(\Delta)_{\ell+1}}\R_0e_\Gamma+L$.
Take any function $t: \mathcal{F}(\Delta)_{\ell+1}\rightarrow\R_0$ and $b\in L$ with $a=\sum_{\Gamma\in\mathcal{F}(\Delta)_{\ell+1}}t(\Gamma) e_\Gamma+b$.
Since $a\in H^\circ=\{x\in W|\langle\zeta,x\rangle>0\}$,
$0<\langle\zeta,a\rangle=\sum_{\Gamma\in\mathcal{F}(\Delta)_{\ell+1}}t(\Gamma)\langle\zeta, e_\Gamma\rangle+\langle\zeta,b\rangle=\sum_{\Gamma\in\mathcal{F}(\Delta)_{\ell+1}}t(\Gamma)\langle\zeta, e_\Gamma\rangle$.
We know that there exists $\Gamma\in\mathcal{F}(\Delta)_{\ell+1}$ with $\langle\zeta, e_\Gamma\rangle>0$.

Consider any face $\Lambda$ of $\Delta$.
We know $\Lambda=\sum_{\Gamma\in\mathcal{F}(\Delta)_{\ell+1}\backslash\Lambda}\R_0e_\Gamma+L$.

\noindent $(a)\Rightarrow (b)$.
Assume $(a)$. $\Lambda\cap H^\circ\neq\emptyset$. Take any point $a\in\Lambda\cap H^\circ$. Since $a\in H^\circ=\{x\in W|\langle\zeta,x\rangle>0\}$, $\langle\zeta,a\rangle>0$.
$a\in\Lambda=\sum_{\Gamma\in\mathcal{F}(\Delta)_{\ell+1}\backslash\Lambda}\R_0e_\Gamma+L$. Take any function $t: \mathcal{F}(\Delta)_{\ell+1}\backslash\Lambda\rightarrow\R_0$ and $b\in L$ with $a=\sum_{\Gamma\in\mathcal{F}(\Delta)_{\ell+1}\backslash\Lambda}t(\Gamma)e_\Gamma+b$.
$0 <\langle\zeta,a\rangle=\sum_{\Gamma\in\mathcal{F}(\Delta)_{\ell+1}\backslash\Lambda}t(\Gamma)\langle\zeta, e_\Gamma\rangle+\langle\zeta,b\rangle=\sum_{\Gamma\in\mathcal{F}(\Delta)_{\ell+1}\backslash\Lambda}t(\Gamma)\langle\zeta, e_\Gamma\rangle$.
We know that there exists $\Gamma\in\mathcal{F}(\Delta)_{\ell+1}\backslash\Lambda$ with $\langle\zeta,e_\Gamma\rangle>0$.
We take any $\Gamma\in\mathcal{F}(\Delta)_{\ell+1}\backslash\Lambda$ with $\langle\zeta,e_\Gamma\rangle>0$.

Note that $H^\vee=\R_0\zeta$ and $\Vect(H)=\R\zeta$.

Consider any $\omega\in\Delta(\Lambda,\Delta|W)$.
$e_\Gamma\in\Gamma\subset\Lambda\subset\Delta(\omega,\Delta|W)$.
$\langle\omega, e_\Gamma\rangle=0$.

Consider any $\chi\in(\{\omega\}+\Vect(H))\cap\Delta^\vee$.
Take $t\in\R$ with $\chi=\omega+t\zeta$.
Since $\chi\in\Delta^\vee$ and $e_\Gamma\in\Gamma\subset\Lambda\subset\Delta$, we have $0\leq\langle\chi, e_\Gamma\rangle=\langle\omega+t\zeta, e_\Gamma\rangle=\langle\omega, e_\Gamma\rangle+t\langle\zeta, e_\Gamma\rangle= t\langle\zeta, e_\Gamma\rangle$.
Since $\langle\zeta, e_\Gamma\rangle>0$, we know $t\geq 0$ and $\chi=\omega+t\zeta\in\{\omega\}+\R_0\zeta=\{\omega\}+H$.

We know $(\{\omega\}+\Vect(H))\cap\Delta^\vee\subset\{\omega\}+H$ and $\omega\in\partial_-\Delta^\vee$.

We know $\Delta(\Lambda,\Delta|W)\subset\partial_-\Delta^\vee$ and $\Delta(\Lambda,\Delta|W)\in\mathcal{F}(\Delta^\vee)^*$.

\noindent $(c)\Rightarrow (a)$.
Assume $(a)$ does not hold.
$\Lambda\cap H^\circ=\emptyset$.
Since $\Lambda\subset\Delta\subset H$, $\Lambda\subset H- H^\circ=\partial H=V=\{x\in W|\langle\zeta,a\rangle=0\}$, and $\langle\zeta,a\rangle=0$ for any $a\in\Lambda$.
In particular, $\langle\zeta,e_\Gamma\rangle=0$ for any $\Gamma\in\mathcal{F}(\Delta)_{\ell+1}\backslash\Lambda$.
If $\Gamma\in\mathcal{F}(\Delta)_{\ell+1}$ and $\langle\zeta,e_\Gamma\rangle>0$, then $\Gamma\not\subset\Lambda$.

Consider any $\omega\in\Delta^\circ(\Lambda,\Delta|W)$.
$\omega\in\Delta^\vee$.
Since $L\subset\Lambda=\Delta(\omega,\Delta|W)$, $\langle\omega, b\rangle=0$ for any $b\in L$.
Consider any $\Gamma\in\mathcal{F}(\Delta)_{\ell+1}$.
If $\Gamma\subset\Lambda$, then
$e_\Gamma\in\Gamma\subset\Lambda=\Delta(\omega,\Delta|W)$ and $\langle\omega, e_\Gamma\rangle=0$.
It is easy to see that if $\Gamma\not\subset\Lambda$, then $\langle\omega, e_\Gamma\rangle>0$

Consider any $t\in\R$.
For any $b\in L$ we have $\langle\omega+t\zeta, b\rangle=\langle\omega, b\rangle+t\langle\zeta, b\rangle=0+t0=0$.
For any $\Gamma\in\mathcal{F}(\Delta)_{\ell+1}$, $\langle\omega+t\zeta, e_\Gamma\rangle=\langle\omega, e_\Gamma\rangle+t\langle\zeta, e_\Gamma\rangle$.
$\omega+t\zeta\in\Delta^\vee$, if and only if, $\langle\omega, e_\Gamma\rangle+t\langle\zeta, e_\Gamma\rangle\geq 0$ for any $\Gamma\in\mathcal{F}(\Delta)_{\ell+1}$.

Consider any $\Gamma\in\mathcal{F}(\Delta)_{\ell+1}$.
If $\Gamma\subset\Lambda$, then
$\langle\omega, e_\Gamma\rangle+t\langle\zeta, e_\Gamma\rangle=0+t0=0$.
If $\Gamma\not\subset\Lambda$ and $\langle\zeta, e_\Gamma\rangle=0$, then
$\langle\omega, e_\Gamma\rangle+t\langle\zeta, e_\Gamma\rangle=\langle\omega, e_\Gamma\rangle>0$.
We consider the case $\Gamma\not\subset\Lambda$ and $\langle\zeta, e_\Gamma\rangle>0$.
We have $\langle\omega, e_\Gamma\rangle>0$, $-\langle\omega, e_\Gamma\rangle/\langle\zeta, e_\Gamma\rangle<0$ and $\langle\omega, e_\Gamma\rangle+t\langle\zeta, e_\Gamma\rangle\geq0$, if and only if, $t\geq-\langle\omega, e_\Gamma\rangle/\langle\zeta, e_\Gamma\rangle$.
Put
$$t_0=\max\{-\frac{\langle\omega, e_\Gamma\rangle}{\langle\zeta, e_\Gamma\rangle}|\Gamma\in\mathcal{F}(\Delta)_{\ell+1}, \langle\zeta, e_\Gamma\rangle>0\}\in\R.$$
$t_0<0$ and $\omega+t\zeta\in\Delta^\vee$, if and only if, $t\geq t_0$ for any $t\in\R$.

Since $\Vect(H)=\R\zeta$ and $H=\R_0\zeta$, we know
$\{\omega+t\zeta|t\in\R, t\geq t_0\}=(\{\omega\}+\Vect(H))\cap\Delta^\vee\not\subset\{\omega\}+\R_0\zeta=\{\omega\}+H$, and $\omega\not\in\partial_-\Delta^\vee$.

We know $\Delta^\circ(\Lambda,\Delta|W)\cap\partial_-\Delta^\vee=\emptyset$.
Claim $(c)$ does not hold.

\noindent $(b)\Rightarrow (c)$. Trivial.

\end{proof}

\begin{lemma}
\label{another correspondence}
Denote
\begin{equation*}
\begin{split}
\mathcal{X}(V)&=\text{the set of all convex pseudo polytopes in }V,\\
\mathcal{X}(V,N)&=\text{the set of all rational convex pseudo polytopes over }N\text{ in }V,\\
\mathcal{Y}(V)&=\{(\Sigma,\phi)|
\Sigma\text{ is a fan in }V^*\text{ such that the support }|\Sigma|\text{ of }\Sigma\text{ is a convex}\\
&\qquad\quad\text{polyhedral cone in }V^*\text{ and there exists a piecewise linear convex}\\
&\qquad\quad\text{function }|\Sigma|\rightarrow\R\text{ over }\Sigma, \phi: |\Sigma|\rightarrow\R \text{ is a piecewise linear convex}\\
&\qquad\quad 
\text{function over }\Sigma\},\\
\mathcal{Y}(V, N)&=\{(\Sigma,\phi)|
\Sigma\text{ is a rational fan over }N^*\text{ in }V^*\text{ such that the support }|\Sigma|\text{ of}\\
&\qquad\quad
\Sigma \text{ is a convex polyhedral cone in } V^*\text{ and there exists a piecewise}\\
&\qquad\quad\text{linear function }|\Sigma|\rightarrow\R\text{ which is convex over }\Sigma\text{ and rational over } \\
&\qquad\quad N^*, \phi: |\Sigma|\rightarrow\R \text{ is a piecewise linear function which is convex over}\\
&\qquad\quad\Sigma\text{ and rational over }N^*\}.
\end{split}
\end{equation*}
$\mathcal{X}(V,N)\subset\mathcal{X}(V)$.
$\mathcal{Y}(V,N)\subset\mathcal{Y}(V)$.
Putting $\Phi(S)=(\Sigma(S|V),\Ord(\;, S|V))\in\mathcal{Y}(V)$ for any $S\in\mathcal{X}$, we define a mapping $\Phi:\mathcal{X}(V)\rightarrow\mathcal{Y}(V)$.
$\Phi$ induces a mapping $\Phi': \mathcal{X}(V, N)\rightarrow\mathcal{Y}(V, N)$.

For any $(\Sigma,\phi)\in\mathcal{Y}$ we denote
\begin{equation*}
\begin{split}
\Xi(\Sigma,\phi)=\{\xi\in W^*|&\xi=\rho^*(\bar{\omega})+t\zeta
\text{ for some }\bar{\omega}\in |\Sigma|\text{ and some }t\in\R\\
&\quad\text{with } t\geq-\phi(\bar{\omega})\}.
\end{split}
\end{equation*}
\begin{enumerate}
\item
Consider any $(\Sigma,\phi)\in\mathcal{Y}(V)$.

$\Xi(\Sigma,\phi)$ is a convex polyhedral cone in $W^*$.
$H^\vee\subset\Xi(\Sigma,\phi)$.
$-H^\vee\not\subset\Xi(\Sigma,\phi)$.

$\Xi(\Sigma,\phi)^\vee|W^*$ is a convex polyhedral cone in $W$.
$\Xi(\Sigma,\phi)^\vee|W^*\subset H$.\hfill\break
$(\Xi(\Sigma,\phi)^\vee|W^*)\cap H^\circ\neq\emptyset$.
$\sigma((\Xi(\Sigma,\phi)^\vee|W^*)\cap H^\circ)\in\mathcal{X}(V)$.

If $(\Sigma,\phi)\in\mathcal{Y}(V, N)$, then $\Xi(\Sigma,\phi)$ is rational over $(N+\Z z)^*$, $\Xi(\Sigma,\phi)^\vee|W^*$ is rational over $N+\Z z$ and $\sigma((\Xi(\Sigma,\phi)^\vee|W^*)\cap H^\circ)\in\mathcal{X}(V, N)$.
\end{enumerate}
Putting
$\Psi(\Sigma,\phi)= \sigma((\Xi(\Sigma,\phi)^\vee|W^*)\cap H^\circ)\in\mathcal{X}(V)$ for any $(\Sigma,\phi)\in\mathcal{Y}(V)$, we define a mapping $\Psi: \mathcal{Y}(V)\rightarrow \mathcal{X}(V)$.
$\Psi$ induces a mapping $\Psi': \mathcal{Y}(V, N)\rightarrow \mathcal{X}(V, N)$.
\begin{enumerate}
\setcounter{enumi}{1}
\item
$\Phi$ and $\Psi$ are bijective mappings, and they are the inverse mappings of each other.

$\Phi'$ and $\Psi'$ are bijective mappings, and they are the inverse mappings of each other.
\end{enumerate}
\end{lemma}

\begin{remark}
Assume $\dim V=3$. 
There exists a rational fan $\Sigma$ over $N^*$ in $V^*$ such that $|\Sigma|=V^*$ and there does \emph{not} exist a piecewise linear convex function $|\Sigma|\rightarrow \R$ over $\Sigma$.

See Fulton~\cite{F93} page 71, Cox~\cite{C10} page 183.
\end{remark}

\begin{theorem}
\label{property of polyhedrons}
Let $S$ be any convex pseudo polytope in $V$, and let $X, Y$ be any finite subsets of $V$ satisfying $S=\Conv(X)+\Convcone(Y)$ and $X\neq\emptyset$.
We consider the dual cone $\Stab(S)^\vee=\Stab(S)^\vee|V\subset V^*$ of $\Stab(S)$.
For simplicity we denote
$s=\dim S\in\Z_0$, $L=\Stab(S)\cap(-\Stab(S))\subset \Stab(S)$, $\ell=\dim L\in\Z_0$.
\begin{enumerate}
\item
We consider any vector space $U$ of finite dimension over $\R$ with $\dim S\leq\dim U\leq\dim V$, any injective homomorphism $\nu:U\rightarrow V$ of vector spaces over $\R$, any point $a\in V$ such that $S\subset\nu(U)+\{a\}$, and any subset $F$ of $S$.
Putting $\bar{\nu}(x)=\nu(x)+a\in V$ for any $x\in U$ we define an injective mapping $\bar{\nu}:U\rightarrow V$.
$S\subset\bar{\nu}(U)$.
The inverse image $\bar{\nu}^{-1}(S)$ is a convex polyhedral cone in $U$.
The set $F$ is a face of $S$, if and only if $\bar{\nu}^{-1}(F)$ is a face of $\bar{\nu}^{-1}(S)$.
\item
We consider any vector space $W$ of finite dimension over $\R$ with $\dim V\leq\dim W$, any injective homomorphism $\pi:V\rightarrow W$ of vector spaces over $\R$, any point $b\in W$ and any subset $F$ of $S$.
Putting $\bar{\pi}(x)=\pi(x)+b$ for any $x\in V$ we define an injective mapping $\bar{\pi}:V\rightarrow W$.
The image $\bar{\pi}(S)$ is a convex polyhedral cone in $W$.
The set $F$ is a face of $S$, if and only if $\bar{\pi}(F)$ is a face of $\bar{\pi}(S)$.
\item
Assume that $S$ is a convex polyhedral cone. For any subset $F$ of $S$, $F$ is a face of the convex pseudo polytope $S$, if and only if, $F$ is a face of the convex polyhedral cone $S$.

\item
$\ell\leq s$. $\ell=s\Leftrightarrow L+\{a\}=S$ for some $a\in V\Leftrightarrow S=\Affi(S)$.
\item
Let $F$ be any face of $S$.
\begin{enumerate}
\item
$F$ is a convex pseudo polyhedron in $V$.
$\Stab(F)$ is a face of $\Stab(S)$.
\item
If $\omega\in\Stab(S)^\vee$ and $F=\Delta(\omega, S)$, then $\Stab(F)=\Delta(\omega,\Stab(S))$.
\item
$\Stab(F)=\Convcone(Y\cap\Stab(F))$.
$\Vect(\Stab(F))=\Vect(Y\cap\Stab(F))$.
\item
$F=\Conv(X\cap F)+\Stab(F)=S\cap\Affi(F)$. $\Affi(F)=\Affi(X\cap F)+\Vect(\Stab(F))$.
\item
If $S$ is rational over $N$, then $F$ is also rational over $N$.
\item
$L=\Stab(F)\cap(-\Stab(F))\subset \Stab(F)\subset\Vect(\Stab(F))\subset\Stab(\Affi(F))$.
$\ell\leq\dim \Stab(F)\leq\dim F\leq s$.
\item
Let $G$ be any face of $S$ with $G\subset F$. We have $\dim G\leq\dim F$.  $\dim G=\dim F$, if and only if, $G=F$.
\item
Let $G$ be any subset of $F$. $G$ is a face of the convex pseudo polyhedron $F$, if and only if, $G$ is a face of $S$ with $G\subset F$.
\end{enumerate}
\item
$\mathcal{F}(S)$ is a finite set.
$S\in\mathcal{F}(S)_s$ and $\mathcal{F}(S)_s=\{S\}$.
$S$ contains any face of $S$.
For any $i\in\Z_0$, $\mathcal{F}(S)_i\neq\emptyset$ if and only if $\ell\leq i\leq s$.
The characteristic number $c(S)$ of $S$ is equal to $\sharp\mathcal{F}(S)_\ell$.
$c(S)$ is a positive integer.
\item
$L=\Convcone(Y\cap L)=\Vect(Y\cap L)$.

Any face $G$ of $S$ with $\dim G=\ell$ is an affine space in $V$ with $\Stab(G)=L$.

For any face $F$ of $S$ and any face $G$ of $S$ with $\dim G=\ell$, $F\supset G$, if and only if, $F\cap G\neq\emptyset$.

For any faces $F$, $G$ of $S$ with $\dim F=\dim G=\ell$, $F=G$, if and only if, $F\cap G\neq\emptyset$.

For any face $F$ of $S$, there exists a face $G$ of $S$ such that $\dim G=\ell$ and $F\supset G$.

Consider any face $G$ of $S$ with $\dim G=\ell$ and any point $\omega\in\Vect(\Stab(S)^\vee)$. The function $\langle\omega,\;\rangle:G\rightarrow\R$ sending $x\in G$ to $\langle\omega,x\rangle\in\R$ is a constant function on $G$.
\item
The skeleton $\mathcal{V}(S)$ of $S$ is a non-empty closed subset of $S$ with finite connected components.
Any connected component of $\mathcal{V}(S)$ is an affine space $G$ in $V$ with $\Stab(G)=L$.
The set of connected components of $\mathcal{V}(S)$ is equal to $\mathcal{F}(S)_\ell$.
The number of connected components of $\mathcal{V}(S)$ is equal to $c(S)$.

For any point $\omega\in\Vect(\Stab(S)^\vee)$, the function $\langle\omega,\;\rangle: \mathcal{V}(S)\rightarrow\R$ sending $x\in \mathcal{V}(S)$ to $\langle\omega,x\rangle\in\R$ is constant on each connected component of $\mathcal{V}(S)$, and this function has only finite number of values.

For any face $F$ of $S$, the intersection $F\cap\mathcal{V}(S)$ is non-empty and union of some connected components of $\mathcal{V}(S)$.
\item
Let $F$ and $G$ be any face of $S$ with $F\subset G$.
We denote $f=\dim F$ and $g=\dim G$.
$\ell\leq f\leq g\leq s$.
There exist $(s-\ell+1)$ of faces $F(\ell), F(\ell+1),\ldots, F(s)$ satisfying the following three conditions:
\begin{enumerate}
\item
For any $i\in\{\ell,\ell+1,\ldots,s-1\}$, $F(i)\subset F(i+1)$.
\item
For any $i\in\{\ell,\ell+1,\ldots,s\}$, $\dim F(i)=i$.
\item
$F(f)=F, F(g)=G, F(s)=S$.
\end{enumerate}
\item
Let $F$ be any face of $S$.
\begin{enumerate}
\item
$F=\partial F\cup F^\circ$. $\partial F\cap F^\circ=\emptyset$.
\item
$F^\circ=F\Leftrightarrow\partial F=\emptyset\Leftrightarrow\dim F=\ell$.
\item
$$\partial F=\bigcup_{G\in\mathcal{F}(F)-\{F\}}G.$$
\item
$F^\circ $ is a non-empty open subset of $\Affi(F)$.
For any $a\in F^\circ$ and any $b\in F$ $\Conv(\{a,b\})-\{b\}\subset F^\circ$.
$F^\circ$ is convex.
$\Clos(F^\circ)=F$.
\end{enumerate}
\item
For any face $G$ of $S$ with $\dim G=\ell$ we take any point $a_G\in G$.\begin{enumerate}
\item
$S=\Conv(\{a_G|G\in\mathcal{F}(S)_{\ell}\})+\Stab(S)$.
\item
For any $\omega\in\Vect(\Stab(S)^\vee)$, $\{\langle\omega, x\rangle|x\in\mathcal{V}(S)\}=\{\langle\omega, a_G\rangle| G\in\mathcal{F}(S)_{\ell}\}$.
\item
For any $\omega\in\Stab(S)^\vee$, $\Ord(\omega, S|V)=\min\{\langle\omega, a_G\rangle| G\in\mathcal{F}(S)_{\ell}\}$.
\item
For any face $F$ of $S$,
$F=\Conv(\{a_G| G\in\mathcal{F}(S)_{\ell}, G\subset F\})+\Stab(F)$.
\end{enumerate}
\item
Consider any $m\in\Z_+$ and any mapping $F:\{1,2,\ldots,m\}\rightarrow\mathcal{F}(S)$.
If $\cap_{i\in\{1,2,\ldots,m\}}F(i)\neq\emptyset$, then the intersection $\cap_{i\in\{1,2,\ldots,m\}}F(i)$ is a face of $S$.
\item
Any proper face $F$ of $S$ is the intersection of all facets of $S$ containing $F$.
\item
The normal fan $\Sigma(S|V)$ of $S$ is a fan in $V^*$.
$|\Sigma(S|V)|=\Stab(S)^\vee$.
$\dim \Sigma(S|V)=\dim \Stab(S)^\vee=\dim V-\ell$.
$c(S)=\sharp\Sigma(S|V)^0$.
If $S$ is rational over $N$, then $\Sigma(S|V)$ is rational over $N^*$.
The minimum element of $\Sigma(S|V)$ is $\Delta(S, S|V)$.
$\Delta^\circ(S, S|V)=\Delta(S, S|V)=\Stab(\Affi(S))^\vee|V$.
$\dim\Delta(S, S|V)=\dim V-s$.

For any $i\in\Z$, $\Sigma(S|V)_i\neq\emptyset$ if and only if $\dim V-s\leq i\leq\dim V-\ell$, and $\Sigma(S|V)^i\neq\emptyset$ if and only if $0\leq i\leq s-\ell$.
\item
Let $F$ be any face of $S$.
\begin{enumerate}
\item
$\Delta(F,S|V)\in\Sigma(S|V)$.
\item
$\Vect(\Delta(F,S|V))=\Stab(\Affi(F))^\vee|V$.
\item
$\dim F+\dim\Delta(F,S|V)=\dim V$.
\item
$\Delta^\circ(F, S|V)=\Delta(F,S|V)^\circ$.
\item
$\Delta(F,S|V)=\Clos(\Delta^\circ(F, S|V))$.
\item
$\Delta^\circ(F,S|V)\subset\Delta^\circ(\Stab(F),\Stab(S)|V)\in\mathcal{F}(\Stab(S)^\vee)$.
\end{enumerate}
\item
For any faces $F,G$ of $S$, $F\subset G$, if and only if, $\Delta(F,S|V)\supset\Delta(G,S|V)$.

The mapping from $\mathcal{F}(S)$ to $\Sigma(S|V)$ sending $F\in\mathcal{F}(S)$ to $\Delta(F,S|V)\in\Sigma(S|V)$ is a bijective mapping such that itself and its inverse mappings are reversing the inclusion relation.
\item
Consider any two faces $F, G$ of $S$.
The following four conditions are equivalent:
\begin{enumerate}
\item
$F\subset G$.
\item
$F^\circ\cap G\neq\emptyset$.
\item
$\Delta(F)\supset\Delta(G)$
\item
$\Delta(F)\cap\Delta^\circ(G)\neq\emptyset$.
\end{enumerate}
The following six conditions are also equivalent:
\begin{enumerate}
\setcounter{enumii}{4}
\item
$F=G$.
\item
$F^\circ=G^\circ$.
\item
$F^\circ\cap G^\circ\neq\emptyset$.
\item
$\Delta(F)=\Delta(G)$.
\item
$\Delta^\circ(F)=\Delta^\circ(G)$.
\item
$\Delta^\circ(F)\cap\Delta^\circ(G)\neq\emptyset$.
\end{enumerate}
\item
$c(S)=\sharp\Sigma(S|V)^0=\sharp\mathcal{F}(S)_\ell=$ the number of connected components of $\mathcal{V}(S)\in\Z_+$.
\item
The following three conditions are equivalent:
\begin{enumerate}
\item
$c(S)=1$.
\item
$\Sigma(S|V)=\mathcal{F}(\Stab(S)^\vee)$.
\item
$S=\{a\}+\Stab(S)$ for some $a\in V$.
\end{enumerate}
\item
The family $\{F^\circ|F\in\mathcal{F}(S)\}$ of subsets of $S$ gives the equivalence class decomposition of $S$, in other words, the following three conditions hold:
\begin{enumerate}
\item
$F^\circ\neq\emptyset$ for any $F\in\mathcal{F}(S)$.
\item
$F^\circ=G^\circ$, if and only if, $F^\circ\cap G^\circ\neq\emptyset$ for any $F\in\mathcal{F}(S)$ and any $G\in\mathcal{F}(S)$.
\item
$$S=\bigcup_{ F\in\mathcal{F}(S)}F^\circ.$$
\end{enumerate}
\item
The family $\{\Delta^\circ(F)|F\in\mathcal{F}(S)\}$ of subsets of $\Stab(S)^\vee$ gives the equivalence class decomposition of $\Stab(S)^\vee $, in other words, the following three conditions hold:
\begin{enumerate}
\item
$\Delta^\circ(F)\neq\emptyset$ for any $F\in\mathcal{F}(S)$.
\item
$\Delta^\circ(F)=\Delta^\circ(G)$, if and only if, $\Delta^\circ(F)\cap\Delta^\circ(G)\neq\emptyset$ for any $F\in\mathcal{F}(S)$ and any $G\in\mathcal{F}(S)$.
\item
$$\Stab(S)^\vee =\bigcup_{ F\in\mathcal{F}(S)}\Delta^\circ(F).$$
\end{enumerate}
\item
The function $\Ord(\;,S|V):\Stab(S)^\vee\rightarrow\R$ sending $\omega\in\Stab(S)^\vee$ to $\Ord(\omega,S|$\hfill\break$V)\in\R$ is a piecewise linear convex function over $\Sigma(S|V)$.

If $S$ is rational over $N$, then this function $\Ord(\;,S|V)$ is rational over $N^*$.
\item
Let $m=\dim\Delta(S)\in\Z_0$.
$m=\dim V-s$.
For any $\Lambda\in\Sigma(S|V)_{m+1}$ we take any point $\omega_\Lambda\in\Lambda-\Delta(S)$.
Then, 
\begin{equation*}\begin{split}
S&=\bigcap_{\omega\in\Stab(S)^\vee}\{x\in V|\langle\omega,x\rangle\geq\Ord(\omega,S|V)\}\cap\Affi(S)\\
&=\bigcap_{\Lambda\in\Sigma(S|V)_{m+1}}\{x\in V|\langle\omega_\Lambda,x\rangle\geq\Ord(\omega_\Lambda,S|V)\}\cap\Affi(S).
\end{split}\end{equation*}
\item
Consider any vector space $W$ of finite dimension over $\R$ and any homomorphism $\pi:V\rightarrow W$ of vector spaces over $\R$.
The image $\pi(S)$ is a convex pseudo polytope in $W$, and $\pi(S)^\circ=\pi(S^\circ)$.

If $S$ is a convex polytope in $V$, then $\pi(S)$ is a convex polytope in $W$.
If $S$ is a convex polyhedral cone in $V$, then $\pi(S)$ is a convex polyhedral cone in $W$.
\end{enumerate}
\end{theorem}

\begin{lemma}
\label{sum}
Consider any convex pseudo polytopes $S$, $T$ in $V$.

$S+T$ is a convex pseudo polytope in $V$, $\Stab(S+T)=\Stab(S)+\Stab(T)$,
$\Stab(S+T)^\vee=\Stab(S)^\vee\cap\Stab(T)^\vee$, and
$\Sigma(S+T|V)=\Sigma(S|V)\hat{\cap}\Sigma(T|V)$.

If $S$ and $T$ are rational over $N$, then $S+T$ is also rational over $N$.

If $S$ and $T$ are convex polytopes, then $S+T$ is also a convex polytope.

If $S$ and $T$ are convex polyhedral cones, then $S+T$ is also a convex polyhedral cone.
\end{lemma}

\begin{cor}
\label{sum2}
Consider any convex pseudo polytope $S$ in $V$ and any convex polyhedral cone $\Delta$ in $V$.

$S+\Delta$ is a convex pseudo polytope in $V$, $\Stab(S+\Delta)=\Stab(S)+\Delta$, $\Stab(S+\Delta)^\vee=\Stab(S)^\vee\cap\Delta^\vee$ and
$\Sigma(S+\Delta|V)=\Sigma(S|V)\hat{\cap}\mathcal{F}(\Delta^\vee|V)$.

If $S$ and $\Delta$ are rational over $N$, then $S+\Delta$ is also rational over $N$.

$S\subset S+\Delta$.
For any $\omega\in\Stab(S+\Delta)^\vee$, we have $\omega\in\Stab(S)^\vee$, $\Ord(\omega, S)=\Ord(\omega, S+\Delta)$, and
$\Delta(\omega, S)=\Delta(\omega, S+\Delta)\cap S$.
\end{cor}

\begin{lemma}
\label{Newton1}
Consider any rational simplicial cone $\Delta$ over $N$ in $V$ with $\dim\Delta=\dim V$ and any subset $X$ of $V$ such that $X\subset(\{a\}+\Delta)\cap(1/m)N$ for some $a\in V$ and some $m\in\Z_+$.

There exists a \emph{finite} subset $Y$ of $X$ satisfying $\Conv(X)+\Delta=\Conv(Y)+\Delta$, and $\Conv(X)+\Delta$ is a rational convex pseudo polytope over $N$ in $V$.
\end{lemma}

\begin{remark}
The subset $X$ of $V$ above is not necessarily finite.
\end{remark}

\begin{lemma}
\label{Newton2}
Let $k$ be any field.
Let $A$ be any complete regular local ring such that $\dim A\geq 1$, $A$ contains $k$ as a subring, and the residue field  $A/M(A)$ is isomorphic to $k$ as algebras over $k$.
Let $P$ be any parameter system of $A$.
Let $\phi$ be any non-zero element of $A$.

\begin{enumerate}
\item
The Newton polyhedron $\Gamma_+(P,\phi)$ of $\phi$ over $P$ is a rational convex pseudo polytope over $\Map(P,\Z)$ in $\Map(P,\R)$.
$\Stab(\Gamma_+(P,\phi))=\Map(P,\R_0)$. $\Gamma_+(P,\phi)\subset\Map(P,\R_0)$.
\item
The normal fan $\Sigma(\Gamma_+(P,\phi)| \Map(P,\R))$ of $\Gamma_+(P,\phi)$ is a rational fan over $\Map(P,\Z)^*$ in $\Map(P,\R)^*$.
$|\Sigma(\Gamma_+(P,\phi)| \Map(P,\R))|=\Map(P,\R_0)^\vee|\Map(P,\R)$.
\item
The Newton polyhedron $\Gamma_+(P,\phi)$ has a vertex.
The skeleton $\mathcal{V}(\Gamma_+(P,\phi))$ of $\Gamma_+(P,\phi)$ is a non-empty finite subset of $\Map(P,\Z_0)$, and $\mathcal{V}(\Gamma_+(P,\phi))$ is the union of all vertices of $\Gamma_+(P,\phi)$.
\begin{equation*}\begin{split}
\mathcal{V}(\Gamma_+(P, \phi))= &\{a\in\Gamma_+(P, \phi)|
\text{ There exists }\omega\in\Map(P,\R_0)^\vee|\Map(P,\R)\text{ such that}\\
&\qquad\text{for any } b\in\Gamma_+(P, \phi)\text{ with }\langle\omega,b\rangle=\langle\omega,a\rangle\text{, we have }b=a\}.\\
\end{split}\end{equation*}
\item
$c(\Gamma_+(P,\phi)) =\sharp \Sigma(\Gamma_+(P,\phi)|\Map(P,\R))^0=\sharp\mathcal{V}(\Gamma_+(P,\phi))=\sharp\mathcal{F}(\Gamma_+(P,\phi))_0$\hfill\break$=\text{the number of vertices of }\Gamma_+(P,\phi)$.
\item
For any $\omega\in\Map(P,\R_0)^\vee|\Map(P,\R)$, we have
\begin{equation*}
\begin{split}
\Ord(P,\omega,\phi)&=\Ord(\omega, \Gamma_+(P,\phi)|\Map(P,\R)), \text{ and}\\
\In(P,\omega,\phi)&=\Ps(P,\Delta(\omega, \Gamma_+(P,\phi)|\Map(P,\R)),\phi).
\end{split}
\end{equation*}
\item
$c(\Gamma_+(P,\phi))=1$, if and only if, $\phi$ has normal crossings over $P$.
\item
If $\dim A=1$, then $c(\Gamma_+(P,\phi))=1$.
\item
If $\psi\in A$, $\omega\in A$ and $\phi=\psi\omega$, then
\begin{equation*}
\begin{split}
\Gamma_+(P,\phi)=&\Gamma_+(P,\psi)+\Gamma_+(P,\omega),\\
\Sigma(\Gamma_+(P,\phi)|\Map(P,\R))=&
\Sigma(\Gamma_+(P,\psi)|\Map(P,\R))\hat{\cap}\Sigma(\Gamma_+(P, \omega)|\Map(P,\R)).
\end{split}
\end{equation*}

\end{enumerate}

Let $z\in P$ be any element.

Let $A'$ denote the completion of $k[P - \{z\}]$ by the maximal ideal $k[P - \{z\}]\cap M(A)=(P - \{z\})k[P - \{z\}]$. The ring $A'$ is a complete regular local subring of $A$ and $M(A')=M(A)\cap A' =(P -\{z\})A'$. The set $P-\{z\}$ is a parameter system of $A'$.

\begin{enumerate}
\setcounter{enumi}{8}
\item
Assume that $\Gamma_+(P, \phi)$ is of $z$-Weierstrass type.
\begin{enumerate}
\item
If $\psi\in A$, $\omega\in A$ and $\phi=\psi\omega$, then both $\Gamma_+(P,\psi)$ and $\Gamma_+(P,\omega)$ are of $z$-Weierstrass type.
\item
$\Ht(z,\Gamma_+(P, \phi))=0\Leftrightarrow \Gamma_+(P,\phi)$ has only one vertex $\Leftrightarrow c(\Gamma_+(P, \phi))=1\Leftrightarrow \phi$ has normal crossings over $P$.
\item 
The Newton polyhedron $\Gamma_+(P,\phi)$ has a unique $z$-top vertex.
\end{enumerate}

Below, by $\{a_1\}$ we denote the unique $z$-top vertex of $\Gamma_+(P,\phi)$. Let $b=\Ord(P, f^{P\vee}_z,\phi)\in\Z_0$ and let $h=\Ht(z,\Gamma_+(P, \phi)) \in\Z_0$.

\begin{enumerate}
\setcounter{enumii}{3}
\item
Consider any $a\in \Gamma_+(P, \phi)$. The equality $\langle f^{P\vee}_x, a\rangle=\Ord(P, f^{P\vee}_x,\phi)$ holds for any $x\in P - \{z\}\Leftrightarrow a-a_1\in\R_0f^P_z$.
\item
$\langle f^{P\vee}_z, a_1\rangle=b+h$.
\item
There exist uniquely an invertible element $u\in A^{\times}$ and a mapping $\phi':\{0, 1,\ldots, h-1\}\rightarrow M(A')$ satisfying
$$\phi=u (\prod_{x\in P -\{z\}}x^{\langle f^{P\vee}_x, a_1\rangle}) z^b (z^h+\sum_{i=0}^{h-1} \phi'(i) z^i),$$
and $\phi'(0)\neq 0$ if $h>0$.
\end{enumerate}
\item
The following two conditions are equivalent:
\begin{enumerate}
\item The Newton polyhedron $\Gamma_+(P,\phi)$ is of $z$-Weierstrass type.
\item There exist uniquely an invertible element $u\in A^{\times}$, a mapping $c:P-\{z\}\rightarrow\Z_0$, and a $z$-Weierstrass polynomial $\psi\in A$
over $P$ satisfying
$$\phi=u \prod_{x\in P-\{z\}}x^{c(x)} \psi$$.
\end{enumerate}
\item
If $\dim A=2$, then $\Gamma_+(P,\phi)$ is $z$-simple.
\item
If $\Gamma_+(P,\phi)$ is $z$-simple, then $\Gamma_+(P,\phi)$ is of $z$-Weierstrass type.
\item
Let $r=c(\Gamma_+(P, \phi))\in\Z_+$. The Newton polyhedron $\Gamma_+(P,\phi)$ is $z$-simple, if and only if, the following three conditions are satisfied:
\begin{enumerate}
\item
For any $a\in \mathcal{V}(\Gamma_+(P, \phi))$ and any $b\in \mathcal{V}(\Gamma_+(P, \phi))$, $\langle f^{P\vee}_z, a\rangle = \langle f^{P\vee}_z, b\rangle$, if and only if, $a=b$.
\end{enumerate}

Below we take the unique bijective mapping $a:\{1,2,\ldots,r\}\rightarrow \mathcal{V}(\Gamma_+(P, \phi))$ satisfying $\langle f^{P\vee}_z, a(i)\rangle > \langle f^{P\vee}_z, a(i+1)\rangle$ for any $i\in\{1,2,\ldots,r-1\}$, if $r\geq 2$.
\begin{enumerate}
\setcounter{enumii}{1}
\item For any $x\in P- \{z\}$, $\langle f^{P\vee}_x, a(2)-a(1)\rangle \geq 0$, if $r\geq 2$.
\item For any $i\in\{1,2,\ldots,r-2\}$ and any  $x\in P - \{z\}$,
$$\frac{\langle f^{P\vee}_x, a(i+1)-a(i)\rangle}{\langle f^{P\vee}_z, a(i)-a(i+1)\rangle} \leq
\frac{\langle f^{P\vee}_x, a(i+2)-a(i+1)\rangle}{\langle f^{P\vee}_z, a(i+1)-a(i+2)\rangle},$$
if $r\geq 3$.
\end{enumerate}

Furthermore, if the above equivalent conditions are satisfied, then the following claims hold:
\begin{enumerate}
\setcounter{enumii}{3}
\item There exists $x\in P- \{z\}$ with  $\langle f^{P\vee}_x, a(2)-a(1)\rangle > 0$.
\item For any $i\in\{1,2,\ldots,r-2\}$, there exists  $x\in P - \{z\}$ with
$$\frac{\langle f^{P\vee}_x, a(i+1)-a(i)\rangle}{\langle f^{P\vee}_z, a(i)-a(i+1)\rangle} < 
\frac{\langle f^{P\vee}_x, a(i+2)-a(i+1)\rangle}{\langle f^{P\vee}_z, a(i+1)-a(i+2)\rangle},$$
if $r\geq 3$.
\end{enumerate}
\item
Assume that $\Gamma_+(P,\psi)$ is $z$-simple. If $\psi\in A$, $\omega\in A$ and $\phi=\psi\omega$, then both $\Gamma_+(P,\psi)$ and $\Gamma_+(P,\omega)$ are $z$-simple.
\item
Assume that $\Gamma_+(P,\phi)$ is of $z$-Weierstrass type. Let $\{a_1\}$ denote the unique $z$-top vertex of $\Gamma_+(P,\phi)$. We take an invertible element $u\in A^{\times}$ and a mapping $\phi':\{0,1,\ldots,h-1\}\rightarrow M(A')$ satisfying
$$\phi=u(\prod_{x\in P -\{z\}}x^{\langle f^{P\vee}_x, a_1\rangle}) z^b (z^h+\sum_{i=0}^{h-1} \phi'(i) z^i),$$
and $\phi'(0)\neq 0$ if $h>0$.
Then, $\Gamma_+(P,\phi)$ is $z$-simple, if and only if, there exist positive integer $r$, and a mapping $c:\{1,2,\ldots,r\}\rightarrow \Map(P,\Z_0)$ satisfying the following conditions:
\begin{enumerate}
\item $1\leq r\leq h+1$. $r=1\Leftrightarrow h=0$.
\item $c(1)=hf_z^P$. $\langle f_z^{P\vee}, c(r)\rangle=0$.
\item For any $i\in\{1,2,\ldots,r-1\}$, we have $\langle f_z^{P\vee}, c(i)-c(i+1)\rangle>0$, if $r\geq 2$.
\item For any $x\in P-\{z\}$, we have
$$\langle f_x^{P\vee},c(2)-c(1)\rangle\geq 0,$$ if $r\geq 2$.
\item There exists $x\in P-\{z\}$ with
$$\langle f_x^{P\vee}, c(2)-c(1)\rangle> 0,$$ if $r\geq 2$.
\item For any $i\in\{1,2,\ldots,r-2\}$ and any $x\in P-\{z\}$, we have
$$\frac{\langle f^{P\vee}_x, c(i+1)-c(i)\rangle}{\langle f^{P\vee}_z, c(i)-c(i+1)\rangle} \leq 
\frac{\langle f^{P\vee}_x, c(i+2)-c(i+1)\rangle}{\langle f^{P\vee}_z, c(i+1)-c(i+2)\rangle},$$
if $r\geq 3$.
\item For any $i\in\{1,2,\ldots,r-2\}$, there exists $x\in P-\{z\}$ with
$$\frac{\langle f^{P\vee}_x, c(i+1)-c(i)\rangle}{\langle f^{P\vee}_z, c(i)-c(i+1)\rangle} < 
\frac{\langle f^{P\vee}_x, c(i+2)-c(i+1)\rangle}{\langle f^{P\vee}_z, c(i+1)-c(i+2)\rangle},$$
if $r\geq 3$.
\item For any $i\in\{2,3,\ldots,r\}$ and any $x\in P-\{z\}$, we have
$$\Ord(P, f^{P\vee}_x, \phi'(\langle f^{P\vee}_z, c(i)\rangle))=\langle f^{P\vee}_x,c(i)\rangle,$$
if $r\geq 2$.
\item For any $i\in\{1,2,\ldots,r-1\}$, any $j\in\Z$
with
$$\langle f^{P\vee}_z, c(i+1)\rangle<j<\langle f^{P\vee}_z, c(i)\rangle$$
and any $x\in P-\{z\}$, we have
\begin{equation*}
\begin{split}
&\Ord(P, f^{P\vee}_x, \phi'(j))\geq\\
&\frac{\langle f^{P\vee}_z, c(i)\rangle-j}{\langle f^{P\vee}_z, c(i)-c(i+1)\rangle} \langle f^{P\vee}_x, c(i+1)\rangle
 +\frac{j-\langle f^{P\vee}_z, c(i+1)\rangle}{\langle f^{P\vee}_z, c(i)-c(i+1)\rangle} \langle f^{P\vee}_x, c(i)\rangle,\\
\end{split}
\end{equation*}
if $r\geq 2$.
\end{enumerate}
\end{enumerate}
\end{lemma}

\section{Star subdivisions}
\label{std}
We study star subdivisions of regular fans.

Let $V$ be any vector space of finite dimension over $\R$, let $N$ be any lattice of $V$, let $\Sigma$ be any regular fan over $N$ in $V$ with $\dim\Sigma\geq 1$, and let $F\in\Sigma$ be any element with $\dim F\geq 1$. For simplicity we denote the barycenter $b_{F/N}$ of $F$ over $N$ by $b$.
$b\in F^\circ\cap N$.

\begin{lemma}
\label{prep1}
Consider any element $\Lambda\in\Sigma$ satisfying $\Lambda+F\in\Sigma$ and $F\not\subset\Lambda$.
\begin{enumerate}
\item
$\Lambda+\R_0b$ is a regular cone over $N$ in $V$.
$\R_0b\in\mathcal{F}(\Lambda+\R_0b)_1$. $\Lambda\in\mathcal{F}(\Lambda+\R_0b)^1$.
$\R_0b\cap\Lambda=\{0\}$.
$\R_0b=\Lambda\Op|(\Lambda+\R_0b)$.
$\Lambda=(\R_0b)\Op|(\Lambda+\R_0b)$.
\item
$\Lambda+\R_0b\subset\Lambda+F\in\Sigma/F$.
$(\Lambda+\R_0b)^\circ\subset(\Lambda+F)^\circ$.
If $\dim F=1$, then $\R_0b=F$, and $\Lambda+\R_0b=\Lambda+F$.
If $\dim F\geq 2$, then $\Lambda+\R_0b\neq\Lambda+F$.
\item
$\dim(\Lambda+\R_0b)=\dim\Lambda+1\leq\dim(\Lambda+F)$.
$\Lambda\Op|(\Lambda+F)\in\mathcal{F}(F)\subset\mathcal{F}(\Lambda+F)$.
$\dim(\Lambda\Op|(\Lambda+F))\geq 1$.
$\dim(\Lambda\Op|(\Lambda+F))=1\Leftrightarrow\dim(\Lambda+\R_0b)= \dim(\Lambda+F)$.
\item
For any $\Lambda'\in\mathcal{F}(\Lambda)$, we have $\Lambda'+F\in\Sigma$, and $F\not\subset\Lambda'$. $\mathcal{F}(\Lambda)\subset(\Sigma/F)\Fc-(\Sigma/F)$.
\item
$\{\Lambda'+\R_0b|\Lambda'\in\mathcal{F}(\Lambda)\}=\mathcal{F}(\Lambda+\R_0b)/\R_0b\subset\mathcal{F}(\Lambda+\R_0b)$.
$\mathcal{F}(\Lambda)= \mathcal{F}(\Lambda+\R_0b)-(\mathcal{F}(\Lambda+\R_0b)/\R_0b)\subset\mathcal{F}(\Lambda+\R_0b)$.
\item
For any $\Delta'\in\Sigma-(\Sigma/F)$, we have
$(\Lambda+\R_0b)\cap\Delta'=\Lambda\cap\Delta'\in\mathcal{F}(\Delta')$, and
$\Lambda\cap\Delta'\in\mathcal{F}(\Lambda)\subset\mathcal{F}(\Lambda+\R_0b)$.
\end{enumerate}
\end{lemma}

\begin{definition}
\label{defstd}
We denote
\begin{equation*}
\begin{split}
\Sigma*F=
(\Sigma-(\Sigma/F))\cup
\{\Delta\in 2^V|& \Delta=\Lambda+\R_0b\text{ for some }\Lambda\in\Sigma\\
&\quad\text{satisfying }\Lambda+F\in\Sigma\text{ and } F\not\subset\Lambda\}\subset 2^V,
\end{split}
\end{equation*}
and we call $\Sigma*F$ the \emph{star subdivision} of $\Sigma$ with \emph{center} in $F$, or the \emph{star subdivision} of $\Sigma$ \emph{along} $F$.
\end{definition}

\begin{lemma}
\label{probcd}
\begin{enumerate}
\item
$\Sigma*F$ is a regular fan over $N$ in $V$.
$\Sigma*F$ is a subdivision of $\Sigma$.
$|\Sigma*F|=|\Sigma|$.
$\dim \Sigma*F=\dim \Sigma$.
If $\Sigma\Mx=\Sigma^0$, then $(\Sigma*F)\Mx=(\Sigma*F)^0$.
If $\Sigma$ is flat, them $\Sigma*F$ is also flat.
\item
$\R_0b\in(\Sigma*F)_1$.
$|\Sigma*F/\R_0b|^\circ=|\Sigma/F|^\circ$.
\item
$(\Sigma*F)-(\Sigma*F/\R_0b)= \Sigma-(\Sigma/F)$.
$\Sigma*F/\R_0b=\{\Delta\in 2^V|\Delta=\Lambda+\R_0b\text{ for some }\Lambda\in\Sigma
\text{ satisfying }\Lambda+F\in\Sigma\text{ and }F\not\subset\Lambda\}$.
\item
If $\dim F=1$, then $\R_0b=F\in\Sigma_1$ and $\Sigma*F=\Sigma$.
If $\dim F\geq 2$, then $\R_0b\not\in\Sigma$, $\Sigma*F\neq\Sigma$,
$(\Sigma*F)_1=\Sigma_1\cup\{\R_0b\}$, and $\sharp(\Sigma*F)_1=\sharp\Sigma_1+1$.
\item
Consider any $\Delta\in \Sigma*F/\R_0b$.
We denote $\Lambda=(\R_0b)\Op|\Delta\in\mathcal{F}(\Delta)$.
\begin{enumerate}
\item
$\Delta=\Lambda+\R_0b$. $\Lambda\cap\R_0b=\{0\}$.
$\R_0b\in\mathcal{F}(\Delta)_1$. $\Lambda\in\mathcal{F}(\Delta)^1$.
\item
$\Lambda+F=\Delta+F\in\Sigma$.
$\Delta^\circ\subset(\Lambda+F)^\circ=(\Delta+F )^\circ$.
\end{enumerate}
\item
$(\Sigma*F)\Mx-(\Sigma*F/\R_0b)= \Sigma\Mx-(\Sigma/F)$.
$(\Sigma*F)\Mx\cap(\Sigma*F/\R_0b)=\{\Delta\in 2^V|\Delta=(E\Op|\Lambda)+ \R_0b\text{ for some }
E\in\mathcal{F}(F)_1\text{ and some }\Lambda\in\Sigma\Mx/F\}$.
\end{enumerate}
\end{lemma}

\begin{example}
Assume $\dim V\geq 3$.
Consider any regular cone $S$ over $N$ in $V$ with $\dim S=3$.
Let $E(1), E(2), E(3)$ denote the three edges of $S$.
We denote $b(i)=b_{E(i)/N}\in E(i)^\circ\cap N$ for any $i\in\{1,2,3\}$ for simplicity.
$E(i)=\R_0b(i)$ for any $i\in\{1,2,3\}$.

Put
\begin{equation*}\begin{split}
T(1)&=\R_0b(1)+\R_0(b(1)+b(2))+\R_0(b(1)+b(3)), \\
T(2)&=\R_0b(2)+\R_0(b(2)+b(3))+\R_0(b(2)+b(1)), \\
T(3)&=\R_0b(3)+\R_0(b(3)+b(1))+\R_0(b(3)+b(2)), \\
T(4)&=\R_0(b(1)+b(2)+b(3))+\R_0(b(1)+b(2))+\R_0(b(1)+b(3)), \\
T(5)&=\R_0(b(1)+b(2)+b(3))+\R_0(b(2)+b(3))+\R_0(b(2)+b(1)), \\
T(6)&=\R_0(b(1)+b(2)+b(3))+\R_0(b(3)+b(1))+\R_0(b(3)+b(2)).
\end{split}\end{equation*}
For any $i\in\{1,2,\ldots,6\}$, $T(i)$ is a regular cone over $N$ in $V$ with $\dim T(i)=3$.

Let $\Phi=\cup_{i\in\{1,2,\ldots,6\}}\mathcal{F}(T(i))\subset 2^V$.
$\Phi$ is a regular fan over $N$ in $V$.
$\dim \Phi=3$ and $|\Phi|=S$.
$\Phi$ is a subdivision of the regular fan $\mathcal{F}(S)$ with $|\Phi|=|\mathcal{F}(S)|$.

For any $\bar{F}\in \mathcal{F}(S)_2$, $\R_0b_{\bar{F}/N}\in\Phi$, and $\Phi$ is not a subdivision of $\mathcal{F}(S)*\bar{F}$.

$\R_0b_{S/N}\in\Phi$, and $\Phi$ is not a subdivision of $\mathcal{F}(S)*S$.
\end{example}

\section{Iterated star subdivisions}
\label{istd}
We study iterated star subdivisions of regular fans.

Let $V$ be any vector space of finite dimension over $\R$, let $N$ be any lattice of $V$, and let $\Sigma$ be any regular fan over $N$ in $V$ with $\dim\Sigma\geq 1$.

\begin{definition}
\label{defistd}
Let $m\in\Z_0$ be any non-negative integer.
We call a mapping $F$ from $\{1,2,\ldots,m\}$ to the set $2^V$ of all subsets of $V$ satisfying the following two conditions a \emph{center sequence} of $\Sigma$ of \emph{length} $m$:
\begin{enumerate}
\item
$F(i)$ is a regular cone over $N$ in $V$ and $\dim F(i)\geq 2$ for any $i\in\{1,2,\ldots,m\}$.
\item
There exists uniquely a mapping $\bar{\Sigma}$ from $\{0,1,\ldots,m\}$ to the set of all regular fans over $N$ in $V$ satisfying the following two conditions:
\begin{enumerate}
\item
$\bar{\Sigma}(0)=\Sigma$.
\item
$F(i)\in\bar{\Sigma}(i-1)$ and $\bar{\Sigma}(i)=\bar{\Sigma}(i-1)*F(i)$ for any $i\in\{1,2,\ldots,m\}$.
\end{enumerate}
\end{enumerate}

Consider any $m\in\Z_0$ and any center sequence $F$ of $\Sigma$ of length $m$.
There exists uniquely a mapping $\bar{\Sigma}$ from $\{0,1,\ldots,m\}$ to the set of all regular fans over $N$ in $V$ satisfying the above two conditions $(a)$ and $(b)$.
Since regular fan $\bar{\Sigma}(m)$ is uniquely determined by $\Sigma$ and the center sequence $F$ of $\Sigma$, we denote $\bar{\Sigma}(m)$ by the symbol
$$\Sigma*F(1)*F(2)*\cdots*F(m),$$
and we call $\Sigma*F(1)*F(2)*\cdots*F(m)$ the \emph{iterated star subdivision} of $\Sigma$ \emph{along} the center sequence $F$ of $\Sigma$.

Consider any regular fan $\Phi$ over $N$ in $V$.
If $\Phi=\Sigma*F(1)*F(2)*\cdots*F(m)$ for some $m\in\Z_0$ and some center sequence $F$ of $\Sigma$ of length $m$, then we call $\Phi$ an \emph{iterated star subdivision} of $\Sigma$.
\end{definition}

\begin{lemma}
\label{propistd1}
Consider any $m\in\Z_0$ and any center sequence $F$ of $\Sigma$ of length $m$.

\begin{enumerate}
\item
$\Sigma*F(1)*F(2)*\cdots*F(m)$ is a regular fan over $N$ in $V$.
$\dim\Sigma*F(1)*F(2)*\cdots*F(m)=\dim\Sigma$.
$\Sigma*F(1)*F(2)*\cdots*F(m)$ is a subdivision of $\Sigma$.
$|\Sigma*F(1)*F(2)*\cdots*F(m)|=|\Sigma|$.
If $\Sigma\Mx=\Sigma^0$, then $(\Sigma*F(1)*F(2)*\cdots*F(m))\Mx=(\Sigma*F(1)*F(2)*\cdots*F(m))^0$.
If $\Sigma$ is flat, then $\Sigma*F(1)*F(2)*\cdots*F(m)$ is also flat.
\item
$\Sigma*F(1)*F(2)*\cdots*F(m)=\Sigma$, if $m=0$.

If $m=1$, then $\Sigma*F(1)*F(2)*\cdots*F(m)$ is equal to the star subdivision $\Sigma*F(1)$ of $\Sigma$ with center in $F(1)$.
\item
If $\dim\Sigma=1$, then $m=0$ and $\Sigma*F(1)*F(2)*\cdots*F(m)=\Sigma$.
\item
For any $i\in\{0,1,\ldots,m\}$, the composition of the inclusion mapping\hfill\break $\{1,2,\ldots,i\}\rightarrow \{1,2,\ldots,m\}$ and $F: \{1,2,\ldots,m\}\rightarrow 2^V$ is a center sequence of $\Sigma$ of lemgth $i$.
\item
Assume $m\geq 1$ and consider any $i\in\{1,2,\ldots,m\}$. 
$\dim F(i)\geq 2$.
$F(i)\in\Sigma*F(1)*F(2)*\cdots*F(i-1)$.
$F(i)\subset|\Sigma|$.
$(\Sigma*F(1)*F(2)*\cdots*F(i-1))*F(i)= \Sigma*F(1)*F(2)*\cdots*F(i)$.
\item
Assume $m\geq 2$ and consider any $i\in\{1,2,\ldots,m\}$.
The mapping $G:\{1,2,\ldots, m-i\}\rightarrow 2^V$ satisfying $G(j)=F(i+j)$ for any $j\in\{1,2,\ldots, m-i\}$ is a center sequence of $\Sigma*F(1)*F(2)*\cdots*F(i)$ of length $m-i$, and
$(\Sigma*F(1)*F(2)*\cdots*F(i))*F(i+1)*F(i+2)*\cdots*F(m)= \Sigma*F(1)*F(2)*\cdots*F(m)$.
\item
$(\Sigma*F(1)*F(2)*\cdots*F(m))_1=\Sigma_1\cup\{\R_0 b_{F(i)/N}|i\in\{1,2,\ldots,m\}\}$.
$\Sigma_1\cap\{\R_0 b_{F(i)/N}|i\in\{1,2,\ldots,m\}\}=\emptyset$.
For any $i\in\{1,2,\ldots,m\}$ and any $j\in\{1,2,\ldots,m\}$, $\R_0 b_{F(i)/N}=\R_0 b_{F(j)/N}$, if and only if, $i=j$.

$\sharp(\Sigma*F(1)*F(2)*\cdots*F(m))_1=\sharp\Sigma_1+m$.
\item
For any $\ell\in\Z_0$ and any center sequence $G$ of $\Sigma*F(1)*F(2)*\cdots*F(m)$ of length $\ell$,
the mapping $H:\{1,2,\ldots, m+\ell\}\rightarrow 2^V$ satisfying $H(i)=F(i)$ for any $i\in\{1,2,\ldots,m\}$ and $H(i)=G(i-m)$ for any $i\in\{m+1,m+2,\ldots,m+\ell\}$
is a center sequence of $\Sigma$ of length $m+\ell$ and
$\Sigma*F(1)*F(2)*\cdots*F(m)*G(1)*G(2)*\cdots*G(\ell)
=(\Sigma*F(1)*F(2)*\cdots*F(m))*G(1)*G(2)*\cdots*G(\ell)$.
\end{enumerate}

Consider any non-empty subset $\Phi$ of $\Sigma$ satisfying $\Phi\Fc=\Phi$.
$\Phi$ is a regular fan over $N$ in $V$.
$|\Phi|\subset|\Sigma|$.

\begin{enumerate}
\setcounter{enumi}{8}
\item
Let $\ell=\sharp\{i\in\{1,2,\ldots,m\}|F(i)\subset|\Phi|\}\in\Z_0$ and
let $\nu:\{1,2,\cdots,\ell\}\rightarrow\{1,2,\ldots,m\}$ be the unique injective mapping preserving the order and satisfying
$\nu(\{1,2,\cdots,\ell\})=\{i\in\{1,2,\ldots,m\}|F(i)\subset|\Phi|\}$.

The composition $F\nu$ is a center sequence of $\Phi$ of length $\ell$, and
$\Phi* F\nu(1)*F\nu(2)*\cdots*F\nu(\ell)
=(\Sigma*F(1)*F(2)*\cdots*F(m))\backslash |\Phi|
\subset \Sigma*F(1)*F(2)*\cdots*F(m)$.
\item
If $F(i)\subset |\Phi|$ for any $i\in\{1,2,\ldots,m\}$, then the sequence $F$ is a center sequence of $\Phi$ of length $m$, and
$\Phi* F(1)*F(2)*\cdots*F(m)
=(\Sigma*F(1)*F(2)*\cdots*F(m))\backslash |\Phi|
\subset \Sigma*F(1)*F(2)*\cdots*F(m)$.
\item
For any $n\in\Z_0$ and any center sequence $G$ of $\Phi$ of length $n$, the sequence $G$ is a center sequence of $\Sigma$ of length $n$, and
$\Phi* G(1)*G(2)*\cdots*G(n)
=(\Sigma*G(1)*G(2)*\cdots*G(n))\backslash |\Phi|
\subset \Sigma*G(1)*G(2)*\cdots*G(n)$.
\end{enumerate}
\end{lemma}

\begin{example}
Assume $\dim V\geq 3$.
Consider any regular cone $S$ over $N$ in $V$ with $\dim S=3$.
Let $E(1), E(2), E(3)$ denote the three edges of $S$.
We denote $b(i)=b_{E(i)/N}\in E(i)^\circ\cap N$ for any $i\in\{1,2,3\}$ for simplicity.
$E(i)=\R_0b(i)$ for any $i\in\{1,2,3\}$.

Let
\begin{gather*}\begin{matrix}
F(1)=\R_0b(1)+\R_0b(3), & F(2)=\R_0b(2)+\R_0b(3), \\
G(1)=\R_0(b(1)+b(3))+\R_0b(2), & G(2)=\R_0(b(2)+b(3))+\R_0b(1).
\end{matrix}\end{gather*}

$\dim F(1)=\dim F(2)=\dim G(1)=\dim G(2)=2$.

We consider two mappings $H$ and $\bar{H}$ from $\{1,2,3\}$ to $2^V$ satisfying $(H(1),H(2),$\hfill\break$H(3))=(F(1), F(2), G(1))$ and $(\bar{H}(1),\bar{H}(2),\bar{H}(3))=(F(2), F(1), G(2))$.
Mappings $H$ and $\bar{H}$ are center sequences of the regular fan $\mathcal{F}(S)$ of length $3$.

$$\mathcal{F}(S)*F(1)*F(2)*G(1)=\mathcal{F}(S)*F(2)*F(1)*G(2).$$
\end{example}

\section{Simpleness and semisimpleness}
\label{simple}

Simpleness and semisimpleness of fans or convex pseudo polytopes are important concepts.

Let $V$ be any vector space of finite dimension over $\R$ with $\dim V\geq 1$, let $N$ be any lattice of $V$, let $\Sigma$ be any fan in the dual vector space $V^*$ of $V$ such that the support $|\Sigma|$ of $\Sigma$ is a regular cone over the dual lattice $N^*$ of $N$ in $V^*$ and $\dim |\Sigma|\geq 1$, let $H\in\mathcal{F}(|\Sigma|)_1$ be any edge of the regular cone $|\Sigma|$, let $S$ be any convex pseudo polytope in $V$ such that $|\Sigma(S|V)|=\Stab(S)^\vee|V$ is a regular cone over $N^*$ in $V^*$ and $\dim |\Sigma(S|V)|\geq 1$, and let $G\in\mathcal{F}(|\Sigma(S|V)|)_1$ be any edge of the regular cone $|\Sigma(S|V)|$, where $\Sigma(S|V)$ denotes the normal fan of $S$.

We denote $L=\Stab(S)\cap(-\Stab(S))=\Vect(|\Sigma(S|V)|)^\vee|V^*\subset V$ and $\ell=\dim L\in\Z_0$. $L$ is the maximum vector subspace over $\R$ in $V$ contained in $\Stab(S)$. Recall that $\mathcal{V}(S)=|\mathcal{F}(S)_\ell|$ denotes the skeleton of $S$.

Note that for any $E\in\mathcal{F}(|\Sigma(S|V)|)_1$, any $F\in\mathcal{F}(S)_\ell$, any $a\in F$ and any $b\in F$,
$E\subset|\Sigma(S|V)|\subset\Vect(|\Sigma(S|V)|)$, and we have $\langle b_{E/N^*},a\rangle=\langle b_{E/N^*},b\rangle$. 
(Theorem~\ref{property of polyhedrons}.7.)

\begin{definition}
\label{defsimple}
\begin{enumerate}
\item
We say that $\Sigma$ is \emph{semisimple}, if $\dim\Delta\geq \dim \Sigma-1$ for any  $\Delta\in\Sigma$ with $\Delta^\circ\subset|\Sigma|^\circ$.
\item
We say that $\Sigma$ is of $H$-\emph{Weierstrass type}, if $ H\Op||\Sigma|\in\Sigma$.
\item
We say that $\Sigma$ is $H$-\emph{simple}, if $\Sigma$ is semisimple and $\Sigma$ is of $H$-Weierstrass type.

\item
We say that $S$ is \emph{semisimple}, if $\dim F\leq \ell+1$ for any face $F$ of $S$ with $\Stab(F)=L$.

\item
We say that $S$ is of $G$-\emph{Weierstrass type}, if there exists a face $F$ of $S$ with $\Delta(F,S|V)= G\Op||\Sigma(S|V)|$.

\item
We say that $S$ is $G$-\emph{simple}, if $S$ is semisimple and $S$ is of $G$-Weierstrass type.

\item
Let $F\in\mathcal{F}(S)_\ell$ be any minimal face of $S$.

We say that $F$ is $G$-\emph{top}, if $\langle b_{G/N^*}, a\rangle=\max\{\langle b_{G/N^*}, c\rangle|c\in\mathcal{V}(S)\}$ for some $a\in F$.

We say that $F$ is $G$-\emph{bottom }, if $\langle b_{G/N^*}, a\rangle=\min\{\langle b_{G/N^*}, c\rangle|c\in\mathcal{V}(S)\}$ for some $a\in F$.
\item
We define
$$\Ht(G,S)= \max\{\langle b_{G/N^*}, c\rangle|c\in\mathcal{V}(S)\}- \min\{\langle b_{G/N^*}, c\rangle|c\in\mathcal{V}(S)\}\in\R_0,$$
and we call $\Ht(G,S)$ $G$-\emph{height} of $S$.
\end{enumerate}
\end{definition}

\begin{lemma}
\label{propsimple}
\begin{enumerate}
\item
$S$ is semisimple, if and only if, $\Sigma(S|V)$ is semisimple.
\item
$S$ is of $G$-Weierstrass type, if and only if, $\Sigma(S|V)$ is of $G$-Weierstrass type.
\item
Note that $L\subset\Vect(G\Op||\Sigma(S|V)|)^\vee|V^*\subset V$ and
$\dim \Vect(G\Op||\Sigma(S|V)|)^\vee$\hfill\break$|V^*=\ell+1$.
Let $W=V/\Vect(G\Op||\Sigma(S|V)|)^\vee|V^*$ denote the residue vector space, and let $\rho:V\rightarrow W$ denote the canonical surjective homomorphism of vector spaces over $\R$ to $W$.

$\rho(N)$ is a lattice in $W$, $\rho(S)$ is a convex pseudo polytope in $W$, $\Stab(\rho(S))$ is a regular cone over $\rho(N)$ in $W$, and $\dim \Stab(\rho(S))=\dim W=\dim V-\ell-1$.

$S$ is of $G$-Weierstrass type $\Leftrightarrow c(\rho(S))=1\Leftrightarrow$ there exists $c\in S$ satisfying $\Ord(b_{E/N^*},S|V)=\langle b_{E/N^*},c\rangle$ for any $E\in\mathcal{F}(|\Sigma(S|V)|)_1-\{G\}$.
\item
Assume that $S$ is of $G$-Weierstrass type.
$S$ has a unique $G$-top minimal face.
$c(S)=1$, if and only if, $\Ht(G,S)=0$.
\item
$S$ is $G$-simple, if and only if, $\Sigma(S|V)$ is $G$-simple.
\item
If $\dim |\Sigma|=1$, then  $\Sigma=\mathcal{F}(|\Sigma|)$.

If $\dim |\Sigma|\leq 2$, then $\Sigma$ is $H$-simple.

If $\dim|\Sigma(S|V)|=1$, then $S=\{a\}+\Stab(S)$ for some $a\in S$.

If $\dim|\Sigma(S|V)|\leq 2$, then $S$ is $G$-simple.
\item
Let $k$ be any field. Let $A$ be any complete regular local ring such that $\dim A\geq 1$, $A$ contains $k$ as a subring, and the residue field  $A/M(A)$ is isomorphic to $k$ as algebras over $k$. Let $P$ be any parameter system of $A$. Let $z\in P$ be any element. Let $\phi\in A$ be any non-zero element. We consider the Newton polyhedron $\Gamma_+(P,\phi)$ in the vector space $\Map(P,\R)$. (See Section~\ref{concept}.)
Let 
$G_z=\R_0f^{P\vee}_z \in\mathcal{F}(\Map(P,\R_0)^\vee|\Map(P,\R))_1$.

$\Gamma_+(P,\phi)$ is of $z$-Weierstrass type, if and only if, it is of $G_z$-Weierstrass type.

$\Gamma_+(P,\phi)$ is $z$-simple, if and only if, it is $G_z$-simple.

$\Ht(z,\Gamma_+(P, \phi))= \Ht(G_z,\Gamma_+(P,\phi))$.
\item
Let $F\in\mathcal{F}(S)_\ell$ be any mimimal face of $S$. 

$F$ is $G$-bottom $\Leftrightarrow F\subset\Delta(b_{G/N^*},S|V)\Leftrightarrow G\subset\Delta(F,S|V)$.

If $F$ is $G$-top, then $\dim(\Delta(F,S|V)\cap(G\Op||\Sigma(S|V)|))=\dim |\Sigma(S|V)|-1$ and $\Delta(F,S|V)\subset (\Delta(F,S|V)\cap(G\Op||\Sigma(S|V)|))+G$.
\item
If $\Sigma$ is semisimple, then $\sharp\Sigma^0=\sharp\{\Delta\in\Sigma^1|\Delta^\circ\subset|\Sigma|^\circ\}+1$.
\item
Let $\Lambda$ be any regular cone over $N^*$ in $V^*$ satisfying $\Lambda\subset|\Sigma|$ and $\dim\Lambda\geq 1$.
If $\Sigma$ is semisimple, then $\Sigma\hat{\cap}\mathcal{F}(\Lambda)$ is also semisimple.
\item
If $\Sigma$ is semisimple, then $\Sigma\backslash\Lambda$ is also semisimple for any $\Lambda\in\mathcal{F}(|\Sigma|)$ with $\dim\Lambda\geq 1$.
\item
Assume that $\Sigma$ is of $H$-Weierstrass type.
For any $\Delta\in\Sigma^0$,\hfill\break $\dim(\Delta\cap(H\Op||\Sigma|))=\dim |\Sigma|-1$, if and only if, $\Delta\supset H\Op||\Sigma|$, and there exists only one element $\Delta\in\Sigma^0$ satisfying these equivalent conditions.
\item
If $\Sigma$ is of $H$-Weierstrass type, then $\Sigma\backslash\Lambda$ is also of $H$-Weierstrass type for any $\Lambda\in\mathcal{F}(|\Sigma|)/H$.
\item
If $\Sigma$ is $H$-simple, then $\Sigma\backslash\Lambda$ is also $H$-simple for any $\Lambda\in\mathcal{F}(|\Sigma|)/H$.
\item
Assume that $\Sigma$ is $H$-simple.
We denote $\bar{\Sigma}^1=\{\Delta\in\Sigma^1|\Delta^\circ\subset|\Sigma|^\circ\}\cup\{ H\Op||\Sigma|\}$.
\begin{enumerate}
\item
$ H\Op||\Sigma|\in \bar{\Sigma}^1\subset\Sigma^1$.
$\sharp\Sigma^0=\sharp\bar{\Sigma}^1$.
\item
We denote $\Delta\leq\Lambda$, if $\Delta+H\supset\Lambda+H$ for any $\Delta\in\Sigma^0$ and any $\Lambda\in\Sigma^0$.
Then, the relation $\leq$ is a total order on $\Sigma^0$.
\item
We denote $\Delta\leq\Lambda$, if $\Delta+H\supset\Lambda+H$ for any $\Delta\in\bar{\Sigma}^1$ and any $\Lambda\in\bar{\Sigma}^1$.
Then, the relation $\leq$ is a total order on $\bar{\Sigma}^1$.
\end{enumerate}

Let $r=\sharp\Sigma^0=\sharp\bar{\Sigma}^1\in\Z_+$.

We consider the total order on $\Sigma^0$ described in $(b)$. Let $\Delta:\{1,2,\ldots,r\}\rightarrow \Sigma^0$ denote the unique bijective mapping preserving the order. 

We consider the total order on $\bar{\Sigma}^1$ described in $(c)$.
Let $\bar{\Delta}:\{1,2,\ldots,r\}\rightarrow \bar{\Sigma}^1$ denote the unique bijective mapping preserving the order.
\begin{enumerate}
\setcounter{enumii}{3}
\item
Consider any $i\in\{1,2,\ldots,r\}$ and any $E\in\mathcal{F}(|\Sigma|)_1-\{H\}$.
There exists a unique real number $c(\Sigma,i, E)\in\R$ depending on the pair $(i, E)$ satisfying $b_{E/N^*}+ c(\Sigma,i,E)b_{H/N^*}\in\Vect(\bar{\Delta}(i))$.
\end{enumerate}

Below we assume $c(\Sigma,i,E)\in\R$ and $b_{E/N^*}+ c(\Sigma,i,E)b_{H/N^*}\in\Vect(\bar{\Delta}(i))$ for any $i\in\{1,2,\ldots,r\}$ and any $E\in\mathcal{F}(|\Sigma|)_1-\{H\}$.

\begin{enumerate}
\setcounter{enumii}{4}
\item
For any $i\in\{1,2,\ldots,r\}$,
$\bar{\Delta}(i)=\Convcone(\{b_{E/N^*}+ c(\Sigma,i,E)b_{H/N^*}|E\in\mathcal{F}(|\Sigma|)_1-\{H\}\})$.
\item
For any $\Gamma\in\bar{\Sigma}^1$, $\Gamma=\Vect(\Gamma)\cap|\Sigma|$.
\item
For any $E\in\mathcal{F}(|\Sigma|)_1-\{H\}$, $ c(\Sigma,1,E) =0$.
\item
If $r\geq 2$, then $c(\Sigma,i,E)\leq c(\Sigma,i+1,E)$ for any $i\in\{1,2,\ldots,r-1\}$ and any $E\in\mathcal{F}(|\Sigma|)_1-\{H\}$.
\item
If $r\geq 2$ and $i\in\{1,2,\ldots,r-1\}$, then $c(\Sigma,i,E)< c(\Sigma,i+1,E)$ for some $E\in\mathcal{F}(|\Sigma|)_1-\{H\}$.
\item
$\Sigma$ is rational over $N^*$, if and only if, $c(\Sigma,i,E)\in\Q$
for any $i\in\{2,3,\ldots,r\}$ and any $E\in\mathcal{F}(|\Sigma|)_1-\{H\}$.
\item
If $r\geq 2$, then $\Delta(i)= \bar{\Delta}(i)+ \bar{\Delta}(i+1)$ for any $i\in\{1,2,\ldots,r-1\}$.
\item
$\Delta(r)= \bar{\Delta}(r)+H$.
\item
$\{\Lambda\in\Sigma|\Lambda\not\subset|\mathcal{F}(|\Sigma|)/H|\}=
\Sigma^0\cup\smash{\bar{\Sigma}}^1$.

\item
$\bar{\Delta}(1)=H\Op||\Sigma|\subset\partial|\Sigma|$.
$\Sigma^0/\bar{\Delta}(1)=\{\Delta(1)\}$.

For any $i\in\{2,3,\ldots,r\}$,
$\bar{\Delta}(i)^\circ\subset|\Sigma|^\circ$,
$\bar{\Delta}(i)\not\subset\partial|\Sigma|$, and
$\Sigma^0/\bar{\Delta}(i)=\{\Delta(i-1), \Delta(i)\}$.
\item
$H\subset \Delta(r)$.
$\smash{\bar{\Sigma}}^1\backslash\Delta(r)=\{\bar{\Delta}(r)\}$.

For any $i\in\{1,2,\ldots,r-1\}$,
$H\not\subset\Delta(i)$,
$\smash{\bar{\Sigma}}^1\backslash\Delta(i)=\{\bar{\Delta}(i), \bar{\Delta}(i+1)\}$.
\item
If $r\geq 2$, then $\Delta(i)\cap\Delta(j)=\bar{\Delta}(i+1)\cap\bar{\Delta}(j)$ for any $i\in\{1,2,\ldots,$\hfill\break$r-1\}$ and any $j\in\{2,3,\ldots,r\}$ with $i<j$.
\item
Consider any $\omega\in\Vect(H\Op||\Sigma|)$.

Take the unique function $\bar{\omega}:\mathcal{F}(H\Op||\Sigma|)_1\rightarrow \R$ satisfying
$\omega=$\hfill\break$\sum_{E\in\mathcal{F}(H\Op||\Sigma|)_1}\bar{\omega}(E)b_{E/N^*}$.
For any $i\in\{1,2,\ldots,r\}$, put $t(i)=$\hfill\break$ \sum_{E\in\mathcal{F}(H\Op||\Sigma|)_1}\bar{\omega}(E)c(\Sigma,i, E)\in\R$.
\begin{enumerate}
\item
$t(1)=0$.
\item
For any $i\in\{1,2,\ldots,r\}$, the following claims hold:
\begin{enumerate}
\item
$(\{\omega\}+\R b_{H/N^*})\cap\Vect(\bar{\Delta}(i))=\{\omega+t(i)b_{H/N^*}\}$.
\item
$\omega+t(i)b_{H/N^*}\in\bar{\Delta}(i)\Leftrightarrow \omega\in H\Op||\Sigma|$.
\item
$\omega+t(i)b_{H/N^*}\in\bar{\Delta}(i)^\circ\Leftrightarrow \omega\in (H\Op||\Sigma|)^\circ$.
\end{enumerate}
\item
The following claims hold for any $i\in\{1,2,\ldots,r-1\}$, if $r\geq 2$.
\begin{enumerate}
\item
If $\omega\in H\Op||\Sigma|$, then $t(i)\leq t(i+1)$.
\item
$(\{\omega\}+\R b_{H/N^*})\cap\Delta(i)$\hfill\break
\begin{equation*}
=
\begin{cases}
\{\omega+t b_{H/N^*}|t\in\R, t(i)\leq t\leq t(i+1)\}&\text{if $\omega\in H\Op||\Sigma|$},\\
\emptyset&\text{if $\omega\not\in H\Op||\Sigma|$}.
\end{cases}
\end{equation*}
\item
If $\omega\in (H\Op||\Sigma|)^\circ$, then $t(i)< t(i+1)$.
\item
$(\{\omega\}+\R b_{H/N^*})\cap\Delta(i)^\circ$\hfill\break
\begin{equation*}
=
\begin{cases}
\{\omega+t b_{H/N^*}|t\in\R, t(i)< t< t(i+1)\}&\text{if $\omega\in (H\Op||\Sigma|)^\circ$},\\
\emptyset&\text{if $\omega\not\in (H\Op||\Sigma|)^\circ$}.
\end{cases}
\end{equation*}
\end{enumerate}
\item
$(\{\omega\}+\R b_{H/N^*})\cap\Delta(r)$\hfill\break
\begin{equation*}
=
\begin{cases}
\{\omega+t b_{H/N^*}|t\in\R, t(r)\leq t\}&\text{if $\omega\in H\Op||\Sigma|$},\\
\emptyset&\text{if $\omega\not\in H\Op||\Sigma|$}.
\end{cases}
\end{equation*}
\item
$(\{\omega\}+\R b_{H/N^*})\cap\Delta(r)^\circ$\hfill\break
\begin{equation*}
=
\begin{cases}
\{\omega+t b_{H/N^*}|t\in\R, t(r)< t\}&\text{if $\omega\in (H\Op||\Sigma|)^\circ$},\\
\emptyset&\text{if $\omega\not\in (H\Op||\Sigma|)^\circ$}.
\end{cases}\end{equation*}
\end{enumerate}
\item
Let $\pi:\Vect(|\Sigma|)\rightarrow \Vect(H\Op||\Sigma|)$ denote the unique surjective homomorphism of vector spaces over $\R$ satisfying $\pi^{-1}(0)=\Vect(H)$ and $\pi(x)=x$ for any $x\in\Vect(H\Op||\Sigma|)$.

$\pi(|\Sigma|)= H\Op||\Sigma|$.
For any $\Delta\in\Sigma$, $\pi(\Delta)\in\mathcal{F}(H\Op||\Sigma|)$,
$\pi^{-1}(\pi(\Delta))\cap|\Sigma|\in\mathcal{F}(|\Sigma|)/H$.

For any $\Delta\in\Sigma$ with $\Vect(H)\subset\Vect(\Delta)$, $\dim \pi(\Delta)=\dim\Delta-1$.

For any $\Delta\in\Sigma$ with $\Vect(H)\not\subset\Vect(\Delta)$, $\dim \pi(\Delta)=\dim\Delta$.

For any $\Delta\in\Sigma$ with $\Delta\not\subset H\Op||\Sigma|$,
$\Delta^\circ\subset (\pi^{-1}(\pi(\Delta)) \cap|\Sigma|)^\circ$.
\item
If $\Delta\in\Sigma$, $\Lambda\in\mathcal{F}(|\Sigma|)$, $\Delta^\circ\subset\Lambda^\circ$ and $\Delta\not\subset H\Op||\Sigma|$, then $H\subset\Lambda$, and $\dim \Delta=\dim\Lambda$ or $\dim \Delta=\dim\Lambda-1$.

If $\Delta\in\Sigma$, $\Lambda\in\mathcal{F}(|\Sigma|)$, $\Delta^\circ\subset\Lambda^\circ$ and $\Delta\subset H\Op||\Sigma|$, then $\Delta=\Lambda$ and $\dim\Delta=\dim\Lambda$.
\end{enumerate}
\item
Assume that $\Sigma$ is $H$-simple.
Consider any $\Lambda\in\mathcal{F}(|\Sigma|)/H$.

We use the same notations $\smash{\bar{\Sigma}}^1$, $r$, $\Delta$ and $\bar{\Delta}$ as above.
We denote $\smash{\overline{\Sigma\backslash\Lambda}}^1=\{\Gamma\in\Sigma\backslash\Lambda|\Gamma^\circ\subset\Lambda^\circ\}\cup\{H\Op|\Lambda\}$.
\begin{enumerate}
\item
$r\in\{i\in\{1,2,\ldots,r\}|\dim(\Delta(i)\cap\Lambda)=\dim\Lambda\}\neq\emptyset$.
\end{enumerate}

Put $\bar{r}=\sharp\{i\in\{1,2,\ldots,r\}|\dim(\Delta(i)\cap\Lambda)=\dim\Lambda\}\in\Z_+$.
Let $\nu:\{1,2,\ldots,\bar{r}\}\rightarrow\{1,2,\ldots,r\}$ be the unique injective mapping preserving the order and satisfying $\nu(\{1,2,\ldots,\bar{r}\})=\{i\in\{1,2,\ldots,r\}|\dim(\Delta(i)\cap\Lambda)=\dim\Lambda\}$.
\begin{enumerate}
\setcounter{enumii}{1}
\item
$1\leq\bar{r}\leq r$.
$\nu(\bar{r})=r$.
\item
$\sharp(\Sigma\backslash\Lambda)^0=\sharp\smash{\overline{\Sigma\backslash\Lambda}}^1=\bar{r}$.
\item
$(\Sigma\backslash\Lambda)^0=\{\Delta\nu(i)\cap\Lambda|i\in\{1,2,\ldots,\bar{r}\}\}$.

We consider the total order on $(\Sigma\backslash\Lambda)^0$ described in $15.(b)$.

The bijective mapping $\{1,2,\ldots,\bar{r}\}\rightarrow(\Sigma\backslash\Lambda)^0$ sending $ i\in\{1,2,\ldots,\bar{r}\}$ to $\Delta\nu(i)\cap\Lambda\in(\Sigma\backslash\Lambda)^0$ preserves the order.
\item
$\smash{\overline{\Sigma\backslash\Lambda}}^1=\{\bar{\Delta}\nu(i)\cap\Lambda|i\in\{1,2,\ldots,\bar{r}\}\}$.

We consider the total order on $\smash{\overline{\Sigma\backslash\Lambda}}^1$ described in $15.(c)$.
The bijective mapping $\{1,2,\ldots,\bar{r}\}\rightarrow\smash{\overline{\Sigma\backslash\Lambda}}^1$ sending $ i\in\{1,2,\ldots,\bar{r}\}$ to $\bar{\Delta}\nu(i)\cap\Lambda\in\smash{\overline{\Sigma\backslash\Lambda}}^1$ preserves the order.
\item
For any $j\in\Z$ with $1\leq j<\nu(1)$, $\Delta(j)\cap\Lambda=\bar{\Delta}\nu(1)\cap\Lambda$.

For any $i\in\{2,3,\ldots, \bar{r}\}$ and any $j\in\Z$ with $\nu(i-1)<j<\nu(i)$,
$\Delta(j)\cap\Lambda=\bar{\Delta}\nu(i)\cap\Lambda$.
\item
For any $j\in\Z$ with $1\leq j\leq\nu(1)$, $\bar{\Delta}(j)\cap\Lambda=\bar{\Delta}\nu(1)\cap\Lambda$.

For any $i\in\{2,3,\ldots, \bar{r}\}$ and any $j\in\Z$ with $\nu(i-1)<j\leq\nu(i)$,
$\bar{\Delta}(j)\cap\Lambda=\bar{\Delta}\nu(i)\cap\Lambda$.
\end{enumerate}
\end{enumerate}

Consider any $A\in\mathcal{F}(S)_\ell$ and any $E\in \mathcal{F}(|\Sigma(S|V)|)_1$.
For any $a\in A$, the real number $\langle b_{E/N^*}, a\rangle$ does not depend on the choice of $a\in A$ and it depends only on $A$ and $E$.
We take any $a\in A$ and we define $\langle b_{E/N^*},A\rangle=\langle b_{E/N^*},a\rangle$.

\begin{enumerate}
\setcounter{enumi}{16}
\item
$S$ is $G$-simple, if and only if, the following three conditions are satisfied:
\begin{enumerate}
\item
For any $A\in \mathcal{F}(S)_\ell$ and any $\bar{A}\in \mathcal{F}(S)_\ell$, $\langle b_{G/N^*}, A\rangle=\langle b_{G/N^*}, \bar{A}\rangle$, if and only if, $A=\bar{A}$.
\end{enumerate}

We assume that the first condition is satisfied. Let $r=c(S)$. Let $A:\{1, 2,\ldots, r\}\rightarrow \mathcal{F}(S)_\ell$ be the unique bijective mapping satisfying $\langle b_{G/N^*}, A(i-1)\rangle > \langle b_{G/N^*}, A(i)\rangle$ for any $i\in\{2,3,\ldots,r\}$, if $r\geq 2$. 
\begin{enumerate}
\setcounter{enumii}{1}
\item $\langle b_{E/N^*}, A(2)\rangle \geq \langle b_{E/N^*},A(1)\rangle$ for any  $E\in \mathcal{F}(|\Sigma(S|V)|)_1-\{G\}$, if $r\geq 2$.
\item
$$\frac{\langle b_{E/N^*}, A(i)\rangle- \langle b_{E/N^*}, A(i-1)\rangle}{\langle b_{G/N^*}, A(i-1)\rangle- \langle b_{G/N^*}, A(i)\rangle} \leq 
\frac{\langle b_{E/N^*}, A(i+1)\rangle- \langle b_{E/N^*},A(i)\rangle}{\langle b_{G/N^*}, A(i)\rangle- \langle b_{G/N^*},A(i+1)\rangle},$$
for any $i\in\{2,3,\ldots,r-1\}$ and any  $E\in \mathcal{F}(|\Sigma(S|V)|)_1-\{G\}$, if $r\geq 3$.
\end{enumerate}
\item
Assume that $S$ is $G$-simple.
Let $r=c(S)\in\Z_+$.
Let $A:\{1,2,\ldots,r\}\rightarrow \mathcal{F}(S)_\ell$ be the unique bijective mapping satisfying  $\langle b_{G/N^*}, A(i-1)\rangle > \langle b_{G/N^*}, A(i)\rangle$ for any $i\in\{2,3,\ldots,r\}$, if $r\geq 2$. 

We denote $\bar{\Sigma}(S|V)^1=\{\Delta\in\Sigma(S|V)^1|\Delta^\circ\subset|\Sigma(S|V)|^\circ\}\cup\{G\Op||\Sigma(S|V)|\}\subset\Sigma(S|V)^1$, and $\Delta_G=\Delta(G\Op||\Sigma(S|V)|, |\Sigma(S|V)||V^*)\in\mathcal{F}(\Stab(S))_{\ell+1}$.
The following claims hold:
\begin{enumerate}
\item There exists $E\in \mathcal{F}(|\Sigma(S|V)|)_1-\{G\}$ with $\langle b_{E/N^*}, A(2)\rangle > \langle b_{E/N^*},$\hfill\break$A(1)\rangle$, if $r\geq 2$.
\item For any $i\in\{2,3,\ldots,r-1\}$ there exists  $E\in \mathcal{F}(\Stab(S)^\vee|V)_1-\{G\}$ with
$$\frac{\langle b_{E/N^*}, A(i)\rangle-\langle b_{E/N^*}, A(i-1)\rangle}{\langle b_{G/N^*}, A(i-1)\rangle- \langle b_{G/N^*},A(i)\rangle} < 
\frac{\langle b_{E/N^*}, A(i+1)\rangle-\langle b_{E/N^*},A(i)\rangle}{\langle b_{G/N^*}, A(i)\rangle- \langle b_{G/N^*},A(i+1)\rangle},$$
if $r\geq 3$.
\item
$A(1)+\Delta_G\in\mathcal{F}(S)_ {\ell+1}$.
$\Stab(A(1)+\Delta_G)= \Delta_G$.
$\Delta(A(1)+\Delta_G, S|V)=G\Op||\Sigma(S|V)|\in\Sigma(S|V)^1$.

If $\hat{A}\in\mathcal{F}(S)_{\ell+1}$ and $\Stab(\hat{A})=\Delta_G$, then $\hat{A}=A(1)+\Delta_G$.
\item
$\Conv(A(i-1)\cup A(i))\in\mathcal{F}(S)_{\ell+1}$,
$\Stab(\Conv(A(i-1)\cup A(i)))=L$,
$\Delta(\Conv(A(i-1)\cup A(i)),S|V)\in\Sigma(S|V)^1$, and
$\Delta(\Conv(A(i-1)\cup A(i)),$\hfill\break$S|V)^\circ\subset|\Sigma(S|V)|^\circ$
for any $i\in\{2,3,\ldots,r\}$, if $r\geq 2$.

If $\hat{A}\in\mathcal{F}(S)_{\ell+1}$ and $\Stab(\hat{A})=L$, then $r\geq 2$ and $\hat{A}=\Conv(A(i-1)\cup A(i))$ for some $i\in\{2,3,\ldots,r\}$.
\item
$\bar{\Sigma}(S|V)^1=\{\Delta(A(1)+\Delta_G, S|V)\}\cup\{\Delta(\Conv(A(i-1)\cup A(i)),S|V)| i\in\{2,3,\ldots,r\}\}$.

\item
We define a bijective mapping $\bar{\Delta}:\{1,2,\ldots,r\}\rightarrow\bar{\Sigma}(S|V)^1$ by putting 
$\bar{\Delta}(1)= \Delta(A(1)+\Delta_G, S|V)$ and
$\bar{\Delta}(i)= \Delta(\Conv(A(i-1)\cup A(i)),S|V)$ for any $i\in\{2,3,\ldots,r\}$.
For any $i\in\{1,2,\ldots,r\}$ and any $E\in$\hfill\break$\mathcal{F}(|\Sigma(S|V)|)_1-\{G\}$, we take a unique real number $c(\Sigma(S|V),i,E)\in\R$ satisfying $b_{E/N^*}+ c(\Sigma(S|V),i,E)b_{G/N^*}\in\Vect(\bar{\Delta}(i))$.

For any $i\in\{2,3,\ldots,r\}$, $\bar{\Delta}(i-1)+G\supset\bar{\Delta}(i)+G$, if $r\geq 2$.

If we define a total order described in $15.(c)$ on $\bar{\Sigma}(S|V)^1$, the mapping $\bar{\Delta}$ preserves the order.

For any $i\in\{2,3,\ldots,r\}$ and any $E\in\mathcal{F}(|\Sigma(S|V)|)_1-\{G\}$,
$$c(\Sigma(S|V),i,E)= \frac{\langle b_{E/N^*}, A(i)\rangle- \langle b_{E/N^*}, A(i-1)\rangle}{\langle b_{G/N^*}, A(i-1)\rangle- \langle b_{G/N^*}, A(i)\rangle},$$
if $r\geq 2$.
\end{enumerate}
\item
Assume that $S$ is of $G$-Weierstrass type. Let $A$ be the unique $G$-top minimal face and let $a_0\in A$ be any point of $A$. We denote $U=\{a\in V|\langle b_{G/N},a\rangle<\langle b_{G/N},a_0\rangle\}$ and $W=\{a\in V|\langle b_{G/N},a\rangle=0\}=\Vect(G)^\vee|V^*$. Let $f\in V$ be any point of $V$ with $\langle b_{G/N},f\rangle=1$. We put
$$\rho(a)=\frac{a-\langle b_{G/N},a\rangle f}{\langle b_{G/N},a_0\rangle-\langle b_{G/N},a\rangle}\in W$$
for any $a\in U$ and we define a mapping $\rho:U\rightarrow W$.
\begin{enumerate}
\item $\rho(S\cap U)=\emptyset$, if and only if, $c(S)=1$.
\end{enumerate}
Below, we assume $c(S)\geq 2$. Note that the normal fan $\Sigma(S|V)$ of $S$ is of $G$-Weierstrass type and there exists uniquely $\Delta\in \Sigma(S|V)^0$ with $\Delta\supset G\Op||\Sigma(S|V)|$. We take the unique $\Delta\in \Sigma(S|V)^0$ with $\Delta\supset G\Op||\Sigma(S|V)|$.
\begin{enumerate}
\setcounter{enumii}{1}
\item
$\rho(S\cap U)$ is a convex pseudo polytope in $W$ with $\Stab(\rho(S\cap U))=\Stab(S)\cap W$. If $S$ is rational over $N$, then $\rho(S\cap U)$ is rational over $N\cap W$.
\item
$c(\rho(S\cap U))=1$, if and only if, there exists a vector subspace $X$ of $V^*$ such that $\dim X=\dim V-1$, $X\cap |\Sigma(S|V)|^\circ\neq\emptyset$, $X\cap(G\Op||\Sigma(S|V)|)^\circ=\emptyset$ and $\Delta=(G\Op||\Sigma(S|V)|)+( X\cap |\Sigma(S|V)|)$.
\item
If $S$ is $G$-simple, then $c(\rho(S\cap U))=1$.
\end{enumerate}
\end{enumerate}
\end{lemma}

\begin{definition}
\label{defsimple}
Assume that $\Sigma$ is $H$-simple.
We denote $\bar{\Sigma}^1=\{\Delta\in\Sigma^1|\Delta^\circ\subset|\Sigma|^\circ\}\cup\{ H\Op||\Sigma|\}$.
\begin{enumerate}
\item
We call $\bar{\Sigma}^1$ the $H$-\emph{skeleton} of $\Sigma$.
\item
We call the total order on $\Sigma^0$ described in Lemma~\ref{propsimple}$.15.$(b) the $H$-\emph{order}.
\item
We call the total order on $\bar{\Sigma}^1$ described in Lemma~\ref{propsimple}$.15.$(c) the $H$-\emph{order}.
\item
We consider the $H$-order on $\bar{\Sigma}^1$.
Let $r=\sharp\bar{\Sigma}^1\in\Z_+$.
Let $\bar{\Delta}:\{1,2,\ldots,r\}\rightarrow \bar{\Sigma}^1$ be the unique bijective mapping preserving the order.
Consider any $i\in\{1,2,\ldots, r\}$ and any $E\in\mathcal{F}(|\Sigma|)_1-\{H\}$.
By Lemma~\ref{propsimple}$.15.$(d) there exists uniquely a real number $c(\Sigma,i,E)\in\R$ depending on the pair $(i, E)$ satisfying $b_{E/N^*}+ c(\Sigma,i,E)b_{H/N^*}\in\Vect(\bar{\Delta}(i))$.

We call $c(\Sigma,i,E)$ the \emph{structure constant} of $\Sigma$ corresponding to the pair $(i,E)$.
\end{enumerate}
\end{definition}

\section{Basic subdivisions}
\label{basic}
We define the concept of basic subdivisions. The procedure of a basic subdivision will be used as a unit for constructing an upward subdivision of a normal fan of a Newton polyhedron in Section~\ref{upward}.

Let $V$ be any vector space of finite dimension over $\R$ with $\dim V\geq 2$ and let $N$ be any lattice of $V$.

Let $H$ be any regular cone of dimension one over $N$ in $V$ and let $\Phi$ be any flat regular fan over $N$ in $V$ such that $\dim \Phi\geq 2$, $H\in\Phi_1$ and $\Phi$ is starry with center in $H$. Note that $\dim(\Phi-(\Phi/H))=\dim \Phi-1\geq 1$.

In addition, we consider any non-negative integer $m\in\Z_0$ and any mapping
$$E:\{1,2,\ldots,m\}\rightarrow (\Phi-(\Phi/H))_1.$$

For any $i\in\{0,1,\ldots,m\}$ and any $\bar{E}\in(\Phi-(\Phi/H))_1$, putting
$$s(i, \bar{E})=\sharp(\{1,2,\ldots,i\}\cap E^{-1}(\bar{E}))\in\Z_0,$$
we define a mapping
$$s:\{0,1,\ldots,m\}\times(\Phi-(\Phi/H))_1\rightarrow \Z_0.$$
The mapping $s$ is uniquely determined depending on the mapping $E$.

For any $i\in\{1,2,\ldots,m\}$, we put
\begin{equation*}\begin{split}
F(i)=&\R_0(b_{E(i)/N}+s(i-1, E(i))b_{H/N})+H\subset V,\\
G(i)=&\R_0(b_{E(i)/N}+s(i-1, E(i))b_{H/N})\subset V,\\
H(i)=&\R_0(b_{E(i)/N}+s(i, E(i))b_{H/N})\subset V.
\end{split}\end{equation*}
We put
$$H(m+1)=H\subset V.$$

Three mappings are defined.
$$F, G:\{1,2,\ldots,m\}\rightarrow 2^V.$$
$$H:\{1,2,\ldots,m, m+1\}\rightarrow 2^V.$$
They are uniquely determined depending on the mapping $E$.
\begin{lemma}
\label{basic1}
\begin{enumerate}
\item
For any $\bar{E}\in(\Phi-(\Phi/H))_1$, $s(0,\bar{E})=0$.

For any $i\in\{1,2,\ldots,m\}$ and any $\bar{E}\in(\Phi-(\Phi/H))_1$,
\begin{equation*}
s(i, \bar{E})=
\begin{cases}
s(i-1,\bar{E})&\text{if $\bar{E}\neq E(i)$},\\
s(i-1,\bar{E})+1&\text{if $\bar{E}=E(i)$}.
\end{cases}\end{equation*}
\item
For any $i\in\{1,2,\ldots,m\}$, $E(i), F(i), G(i), H(i), G(i)+H(i)$ and $E(i)+H$ are regular cones over $N$ in $V$,
$\dim E(i)=\dim G(i)=\dim H(i)=1$, $\dim F(i)=\dim(G(i)+H(i))=\dim(E(i)+H)=2$, and
$G(i)+H(i)\subset F(i)\subset E(i)+H\in\Phi_2/H$.
\item
For any $i\in\{1,2,\ldots,m\}$, $\mathcal{F}(G(i)+H(i))_1=\{G(i), H(i)\}$,
$\mathcal{F}(F(i))_1=\{G(i), H\}$,
$\mathcal{F}(E(i)+H)_1=\{E(i), H\}$,
$b_{G(i)/N}=b_{E(i)/N}+s(i-1,E)b_{H/N}$,
$b_{H(i)/N}=b_{E(i)/N}+s(i,E)b_{H/N}= b_{G(i)/N}+ b_{H/N}=b_{F(i)/N}$, and
$H(i)=\R_0 b_{F(i)/N}\subset F(i)$.
\item
The mapping $F$ is a center sequence of $\Phi$ of length $m$ such that $\dim F(i)=2$
and $F(i)^\circ \subset|\Phi/H|^\circ$
for any $i\in\{1,2,\ldots,m\}$, and it is determined by the sextuplet $(V, N, H, \Phi, m, E)$ uniquely.

The iterated star subdivision $\Phi*F(1)*F(2)*\cdots*F(m)$ of $\Phi$ along $F$ is a regular fan determined by the sextuplet $(V, N, H, \Phi, m, E)$ uniquely.
\end{enumerate}
\end{lemma}

Below, for simplicity we denote
$$\Omega= \Phi*F(1)*F(2)*\cdots*F(m)\subset 2^V.$$

For any $i\in\{1,2,\ldots,m\}$, we put
$$\Omega(i)=(\Omega/(G(i)+H(i)))\Fc\subset \Omega\subset 2^V.$$

For $m+1$, we put
$$\Omega(m+1)=(\Omega/H(m+1))\Fc\subset \Omega\subset 2^V.$$

We obtain a mapping
$$\Omega:\{1,2,\ldots,m+1\}\rightarrow 2^{2^V}.$$

\begin{remark}
We denote two different objects by the same symbol $\Omega$.
One satisfies $\Omega\in 2^{2^V}$ and the other $\Omega$ is a mapping from $\{1,2,\ldots,m+1\}$ to $2^{2^V}$.
It is easy to distinguish them.
\end{remark}

\begin{lemma}
\label{property of basic subdivisions}
\begin{enumerate}
\item
The iterated star subdivision $\Omega$ of $\Phi$ is a flat regular fan over $N$ in $V$.
$|\Omega|=|\Phi|$.
$\dim \Omega=\dim \Phi$. 
$\Vect(|\Omega|)=\Vect(|\Phi|)$.
\item
For any $i\in\{1,2,\ldots,m\}$, $G(i)\in\Omega_1$, 
$G(i)+H(i)\in\Omega_2$, and $F(i)^\circ\subset|\Phi/H|^\circ$.
For any $i\in\{1,2,\ldots,m+1\}$, $H(i)\in\Omega_1$.
\item
For any $\Delta\in\Phi/H$, $|\Omega\backslash\Delta|=\Delta$ and $\Omega\backslash\Delta$ is $H$-simple.
\item
$\Omega\backslash|\Phi-(\Phi/H)|= \Phi-(\Phi/H)$.
\item
For any $i\in\{1,2,\ldots,m+1\}$, $\Omega(i)$ is a flat regular fan over $N$ in $V$,
$\dim \Omega(i)=\dim\Phi$, 
$\Vect(|\Omega(i)|)=\Vect(|\Phi|)$,
$H(i)\in\Omega(i)_1$, and $\Omega(i)$ is starry with center in $H(i)$.
\item
For any $i\in\{1,2,\ldots,m\}$,
$G(i)+H(i)=|\Omega(i)|\cap(E(i)+H)\in\Omega(i)_2$, 
$G(i)=| \Omega(i)-(\Omega(i)/H(i))|\cap(G(i)+H(i))\in(\Omega(i)-(\Omega(i)/H(i)))_1\subset\Omega(i)_1$,
$\Omega(i)=(\Omega(i)/(G(i)+H(i)))\Fc=(\Omega(i)/G(i))\Fc$, and
$H(i)=| \Omega(i)-(\Omega(i)/G(i))|\cap(G(i)+H(i))\in(\Omega(i)-(\Omega(i)/G(i)))_1\subset\Omega(i)_1$
\item
For any $i\in\{1,2,\ldots, m\}$ any $j\in\{2,3,\ldots, m+1\}$ with $i<j$,
$\Omega(i)\cap\Omega(j)=(\Omega(i)-(\Omega(i)/G(i)))\cap(\Omega(j)-(\Omega(j)/H(j)))$.
\item
$$|\Phi-(\Phi/H)|\cup(\bigcup_{i\in\{1,2,\ldots,m+1\}}|\Omega(i)/H(i)|^\circ)=|\Phi|$$
\item
For any $i\in\{0,1,\ldots,m+1\}$, $|\Phi-(\Phi/H)|\cap|\Omega(i)/H(i)|^\circ=\emptyset$.

For any $i\in\{0,1,\ldots,m+1\}$ and any $j\in\{0,1,\ldots,m+1\}$ with $i\neq j$,
$|\Omega(i)/H(i)|^\circ\cap|\Omega(j)/H(j)|^\circ=\emptyset$.
\item
$$(\Phi-(\Phi/H))\cup(\bigcup_{i\in\{1,2,\ldots,m+1\}}(\Omega(i)/H(i)))=\Omega$$
\item
For any $i\in\{0,1,\ldots,m+1\}$, $(\Phi-(\Phi/H))\cap(\Omega(i)/H(i))=\emptyset$.

For any $i\in\{0,1,\ldots,m+1\}$ and any $j\in\{0,1,\ldots,m+1\}$ with $i\neq j$,
$(\Omega(i)/H(i))\cap(\Omega(j)/H(j))=\emptyset$.
\item
$$\bigcup_{i\in\{1,2,\ldots,m+1\}}\Omega(i)\Mx=\Omega\Mx$$
\item
For any $i\in\{0,1,\ldots,m+1\}$ and any $j\in\{0,1,\ldots,m+1\}$ with $i\neq j$,
$\Omega(i)\Mx\cap\Omega(j)\Mx=\emptyset$.
\end{enumerate}

For any $i\in\{0,1,\ldots,m+1\}$, we denote
$$X(i)=|\Phi-(\Phi/H)|\cup(\bigcup_{j\in\{1,2,\ldots,i\}}|\Omega(j)/H(j)|^\circ)\subset|\Phi|.$$
\begin{enumerate}
\setcounter{enumi}{13}
\item
$X(0)= |\Phi-(\Phi/H)|$. $X(m+1)=|\Phi|$.

For any $i\in\{1,2,\ldots, m+1\}$, $X(i-1)\subset X(i)$, $X(i-1)\neq X(i)$, and $|\Omega(i)|\cap X(i-1)=|\Omega(i)-(\Omega(i)/H(i))|$.
\item
For any $i\in\{0,1,\ldots,m+1\}$,
$$X(i)= |\Phi-(\Phi/H)|\cup(\bigcup_{j\in\{1,2,\ldots,i\}}|\Omega(j)|),$$
$$|\Omega\backslash X(i)|=X(i),$$

and $X(i)$ is a closed subset of $|\Phi|$.
\end{enumerate}
\end{lemma}

For any $i\in\{1,2,\ldots,m\}$, we put
$$\Omega^\circ(i)=\Omega(i)/G(i) \subset \Omega(i)\subset 2^V.$$

For $m+1$, we put
$$\Omega^\circ(m+1)=\Omega(m+1) \subset 2^V.$$

We obtain a mapping
$$\Omega^\circ:\{1,2,\ldots,m+1\}\rightarrow 2^{2^V}.$$

\begin{lemma}
\label{property of basic subdivisions2}
\begin{enumerate}
\item
For any $i\in\{1,2,\ldots,m+1\}$, the following claims hold:
\begin{enumerate}
\item
$\Omega(i)\Mx\subset\Omega^\circ(i)\subset\Omega(i)$.
$(\Omega(i)\Mx)\Fc=\Omega^\circ(i)\Fc=\Omega(i)$.
$|\Omega(i)\Mx|=|\Omega^\circ(i)|=|\Omega(i)|\supset|\Omega^\circ(i)|^\circ$.
\item
$\Omega^\circ(i)=\Omega(i)\Leftrightarrow|\Omega^\circ(i)|^\circ=|\Omega(i)|\Leftrightarrow i=m+1$.
\item
If $i\neq m+1$, then
$\Omega(i)=\Omega^\circ(i)\cup(\Omega(i)-(\Omega(i)/G(i)))$, $\Omega^\circ(i)\cap(\Omega(i)-(\Omega(i)/G(i)))=\emptyset$,
$|\Omega(i)|=|\Omega^\circ(i)|^\circ\cup|\Omega(i)-(\Omega(i)/G(i))|$, and $|\Omega^\circ(i)|^\circ\cap|\Omega(i)-(\Omega(i)/G(i))|=\emptyset$.
\item
For any $\Theta\in \Omega^\circ(i)/H(i)$, $H(i)\Op|\Theta\in\Omega^\circ(i)-(\Omega^\circ(i)/H(i))$.
For any $\Lambda\in\Omega^\circ(i)-(\Omega^\circ(i)/H(i))$, $\Lambda+H(i)\in \Omega^\circ(i)/H(i)$.

The mapping from $\Omega^\circ(i)/H(i)$ to $\Omega^\circ(i)-(\Omega^\circ(i)/H(i))$ sending\hfill\break$\Theta\in \Omega^\circ(i)/H(i)$ to $H(i)\Op|\Theta\in\Omega^\circ(i)-(\Omega^\circ(i)/H(i))$ and the mapping from $\Omega^\circ(i)-(\Omega^\circ(i)/H(i))$ to $\Omega^\circ(i)/H(i)$ sending $\Lambda\in\Omega^\circ(i)-(\Omega^\circ(i)/H(i))$ to $\Lambda+H(i)\in \Omega^\circ(i)/H(i)$ are bijective mapping preserving the inclusion relation between $\Omega^\circ(i)/H(i)$ and $\Omega^\circ(i)-(\Omega^\circ(i)/H(i))$, and they are the inverse mappings of each other.

Furthermore, if $\Theta\in \Omega^\circ(i)/H(i)$ and $\Lambda\in\Omega^\circ(i)-(\Omega^\circ(i)/H(i))$ correspond to each other by them, then $\dim\Theta=\dim\Lambda+1$.
\end{enumerate}
\item
$$\bigcup_{i\in\{1,2,\ldots,m+1\}}|\Omega^\circ(i)|^\circ=|\Phi|$$
\item
For any $i\in\{0,1,\ldots,m+1\}$ and any $j\in\{0,1,\ldots,m+1\}$ with $i\neq j$,$|\Omega^\circ(i) |^\circ\cap|\Omega^\circ(j)|^\circ=\emptyset$.
\item
$$\bigcup_{i\in\{1,2,\ldots,m+1\}}\Omega^\circ(i)=\Omega$$
\item
For any $i\in\{0,1,\ldots,m+1\}$ and any $j\in\{0,1,\ldots,m+1\}$ with $i\neq j$,
$\Omega^\circ(i)\cap\Omega^\circ(j)=\emptyset$.
\end{enumerate}

For any $i\in\{0,1,\ldots,m+1\}$, we denote
$$Y(i)= \bigcup_{j\in\{i+1,i+2,\ldots,m+1\}}|\Omega^\circ(j)|^\circ\subset|\Phi|.$$
\begin{enumerate}
\setcounter{enumi}{5}
\item
$Y(0)= |\Phi|$. $Y(m+1)=\emptyset$.

For any $i\in\{1,2,\ldots, m+1\}$, $Y(i-1)\supset Y(i)$, $Y(i-1)\neq Y(i)$. For any $i\in\{1,2,\ldots, m\}$, $|\Omega(i)|\cap Y(i)=|\Omega(i)-(\Omega(i)/G(i))|$.
\item
For any $i\in\{0,1,\ldots,m+1\}$,
$$Y(i)= \bigcup_{j\in\{i+1,i+2,\ldots,m+1\}}|\Omega(j)|,$$
$$|\Omega\backslash Y(i)|=Y(i),$$

and $Y(i)$ is a closed subset of $|\Phi|$.

\end{enumerate}
\end{lemma}

\begin{definition}
\label{defbasic}
We denote $\Omega=\mathcal{E}*F(1)*F(2)*\cdots*F(m)$ above by the symbol
$$\Omega(V, N, H, \Phi, m, E)\subset 2^V,$$
and we call it the \emph{basic subdivision} associated with the sextuplet $(V, N, H, \Phi, m, E)$, because $\Omega$ is uniquely determined depending on the sextuplet $(V, N, H, \Phi, m, E)$.

$\Omega(V, N, H, \Phi, m, E)$ is an iterated star subdivision of $\Phi$, it is a flat regular fan over $N$ in $V$, and for any $\Delta\in\Phi/H$, $\Omega(V, N, H, \Phi, m, E) \backslash\Delta$ is $H$-simple.
$|\Omega(V, N, H, \Phi, m, E)|=|\Phi|$.

Note that depending on the sextuplet $(V, N, H, \Phi, m, E)$, six mappings
$$s: \{0,1,\ldots,m\}\times(\Phi-(\Phi/H))_1\rightarrow \Z_0,$$
$$F, G:\{1,2,\ldots,m\}\rightarrow 2^V,$$
$$H:\{1,2,\ldots,m, m+1\}\rightarrow 2^V,$$
$$\Omega:\{1,2,\ldots,m, m+1\}\rightarrow 2^{2^V},$$
$$\Omega^\circ:\{1,2,\ldots,m, m+1\}\rightarrow 2^{2^V},$$
are defined above. 

We denote these six mappings by $s(V, N, H, \Phi, m, E)$, $F(V, N, H, \Phi, m, E)$, \hfill\break$G(V, N, H, \Phi, m, E)$, $H(V, N, H, \Phi, m, E)$, $\Omega(V, N, H, \Phi, m, E)$ and\hfill\break
$\Omega^\circ(V, N, H, \Phi, m, E)$ respectively, and we express the dependence explicitly.
\end{definition}

\begin{remark}
We denote two different objects by the same symbol $\Omega(V, N, H, \Phi, m, E)$.
One satisfies $\Omega(V, N, H, \Phi, m, E)\in 2^{2^V}$ and the other $\Omega(V, N, H, \Phi, m, E)$ is a mapping from $\{1,2,\ldots,m+1\}$ to $2^{2^V}$.
It is easy to distinguish them.
\end{remark}

\begin{lemma}
\label{basic2}
Consider any subset $\hat{\Phi}$ of $\Phi$ satisfying $\dim \hat{\Phi}=\dim \Vect(|\hat{\Phi}|)\geq 2$, $\hat{\Phi}\Mx=\hat{\Phi}^0$, $H\in\hat{\Phi}_1$ and $\hat{\Phi}=(\hat{\Phi}/H)\Fc$.

We know that $\emptyset\neq\hat{\Phi}=\hat{\Phi}\Fc$, $\hat{\Phi}$ is flat regular fan over $N$ in $V$ and it is starry with center in $H$.

Note that $\emptyset\neq(\hat{\Phi}-(\hat{\Phi}/H))_1\subset(\Phi-(\Phi/H))_1$ and $E^{-1}((\hat{\Phi}-(\hat{\Phi}/H))_1)\subset\{1,2,\ldots,m\}$.

Let $\hat{m}=\sharp E^{-1}((\hat{\Phi}-(\hat{\Phi}/H))_1)\in\Z_0$.
Let $\hat{\tau}:\{1,2,\ldots,\hat{m}\}\rightarrow\{1,2,\ldots,m\}$ be the unique injective mapping preserving the order and satisfying $\hat{\tau}(\{1,2,\ldots,\hat{m}\})= E^{-1}((\hat{\Phi}-(\hat{\Phi}/H))_1)$.
Putting $\hat{\tau}(0)=0$ and  $\hat{\tau}(\hat{m}+1)=m+1$, we define an extension
$\hat{\tau}:\{0,1,2,\ldots,\hat{m},\hat{m}+1\}\rightarrow\{0,1,2,\ldots,m,m+1\}$ of $\hat{\tau}:\{1,2,\ldots,\hat{m}\}\rightarrow\{1,2,\ldots,m\}$.
Let $\hat{E}:\{1,2,\ldots,\hat{m}\}\rightarrow(\hat{\Phi}-(\hat{\Phi}/H))_1$ be the unique mapping satisfying $\iota\hat{E}=E\hat{\tau}$, where $\iota: (\hat{\Phi}-(\hat{\Phi}/H))_1\rightarrow(\Phi-(\Phi/H))_1$ denotes the inclusion mapping.

\begin{enumerate}
\item
$\Omega(V, N, H, \Phi, m, E)\backslash|\hat{\Phi}|=
\Omega(V, N, H, \hat{\Phi}, \hat{m}, \hat{E})$.
\item
Let $s=s(V, N, H, \Phi, m, E)$ and $\hat{s}=s(V, N H, \hat{\Phi}, \hat{m}, \hat{E})$.

For any $\bar{E}\in(\hat{\Phi}-(\hat{\Phi}/H))_1$, any $i\in\{0,1,\ldots,\hat{m}\}$ and any $j\in\{0,1,\ldots,m\}$ with $\hat{\tau}(i)\leq j<\hat{\tau}(i+1)$,
$\hat{s}(i, \bar{E})=s(j,\bar{E})=s(\hat{\tau}(i), \bar{E})=s(\hat{\tau}(i+1)-1, \bar{E})$.
\item
Let $F=F(V, N, H, \Phi, m, E)$, $G=G(V, N, H, \Phi, m, E)$,\hfill\break $\hat{F}=F(V, N, H, \hat{\Phi}, \hat{m}, \hat{E})$, and
$\hat{G}=G(V, N, H, \hat{\Phi}, \hat{m}, \hat{E})$.

$\hat{F}=F\hat{\tau}$, and $\hat{G}=G\hat{\tau}$.
\item
Let $H=H(V, N, H, \Phi, m, E)$, $\Omega=\Omega(V, N, H, \Phi, m, E)$,\hfill\break$\Omega^\circ=\Omega^\circ(V, N, H, \Phi, m, E)$, $\hat{H}=H(V, N, H, \hat{\Phi}, \hat{m}, \hat{E})$, 
$\hat{\Omega}=\Omega(V, N, H, \hat{\Phi}, \hat{m}, \hat{E})$, and $\smash{\hat{\Omega}}^\circ=\Omega^\circ(V, N, H, \hat{\Phi}, \hat{m}, \hat{E})$.
\begin{enumerate}
\item
$\hat{H}=H\hat{\tau}$.
\item
$\hat{\Omega}(\hat{m}+1) =\Omega\hat{\tau} (\hat{m}+1)\backslash|\hat{\Phi}|=\Omega(m+1)\backslash|\hat{\Phi}|$.
\item
For any $i\in\{1,2,\ldots,\hat{m}\}$,
\begin{equation*}\begin{split}
&\hat{\Omega}(i)\subset\Omega\hat{\tau}(i)\backslash|\hat{\Phi}|, \text{ and}\\
&(\Omega\hat{\tau}(i)\backslash|\hat{\Phi}|)-\hat{\Omega}(i) \\
&\quad\subset(\Omega\hat{\tau}(i)-(\Omega\hat{\tau}(i)/G\hat{\tau}(i))\cap(\Omega\hat{\tau}(i)-(\Omega\hat{\tau}(i)/H\hat{\tau}(i)).
\end{split}\end{equation*}
\item
For any $j\in\{1,2,\ldots, m\}-\hat{\tau}(\{1,2,\ldots,\hat{m}\})$,
$$\Omega(j)\backslash|\hat{\Phi}|\subset(\Omega(j)-(\Omega(j)/G(j))\cap(\Omega(j))-(\Omega(j)/H(j)).$$
\item
For any $i\in\{1,2,\ldots,\hat{m}, \hat{m}+1\}$, $\smash{\hat{\Omega}}^\circ(i)=\Omega^\circ\hat{\tau}(i)\backslash|\hat{\Phi}|$, and $|\smash{\hat{\Omega}}^\circ(i)|^\circ=|\Omega^\circ\hat{\tau}(i)|^\circ\cap|\hat{\Phi}|$.
\end{enumerate}
\item
For any $j\in\{1,2,\ldots,m\}$, $j\in\hat{\tau}(\{1,2,\ldots,\hat{m}\})\Leftrightarrow E(j)\in(\hat{\Phi}-(\hat{\Phi}/H))_1\Leftrightarrow
F(j)\subset|\hat{\Phi}|\Leftrightarrow
G(j)\subset|\hat{\Phi}|\Leftrightarrow
H(j)\subset|\hat{\Phi}|\Leftrightarrow
G(j)+H(j)\subset|\hat{\Phi}|$.

$\{\hat{F}(i)|i\in\{1,2,\ldots,\hat{m}\}\}=\{F\hat{\tau}(i)|i\in\{1,2,\ldots,\hat{m}\}\}=$\hfill\break$\{F(j)|j\in\{1,2,\ldots,m\}, F(j)\subset|\hat{\Phi}|\}$.

$\{\hat{G}(i)|i\in\{1,2,\ldots,\hat{m}\}\}=\{G\hat{\tau}(i)|i\in\{1,2,\ldots,\hat{m}\}\}=$\hfill\break$\{G(j)|j\in\{1,2,\ldots,m\}, G(j)\subset|\hat{\Phi}|\}$.

$\{\hat{H}(i)|i\in\{1,2,\ldots,\hat{m}, \hat{m}+1\}\}=\{H\hat{\tau}(i)|i\in\{1,2,\ldots,\hat{m}, \hat{m}+1\}\}=\{H(j)|j\in\{1,2,\ldots,m, m+1\}, H(j)\subset|\hat{\Phi}|\}$.
\item
For any $j\in\{1,2,\ldots,m+1\}$, $j\in\hat{\tau}(\{1,2,\ldots,\hat{m}+1\})
\Leftrightarrow
|\Omega^\circ(j)|^\circ\cap|\hat{\Phi}|\neq\emptyset$.

For any $i\in\{1,2,\ldots,\hat{m}+1\}$ and any $\Theta\in\Omega^\circ\hat{\tau}(i)/H\hat{\tau}(i)$, $\Theta\subset|\hat{\Phi}|$, if and only if, $H\hat{\tau}(i)\Op|\Theta\subset|\hat{\Phi}|$.
\item
If $i\in\{1,2,\ldots,\hat{m}\}$, $j\in\{1,2,\ldots,m\}$ and $\hat{\tau}(i)\leq j<\hat{\tau}(i+1)$, then
\begin{equation*}\begin{split}
X(j)\cap|\hat{\Phi}|&=|\hat{\Phi}-(\hat{\Phi}/H)|\cup
(\bigcup_{k\in\{1,2,\ldots, i\}}|\hat{\Omega}(k)|),\\
Y(j)\cap|\hat{\Phi}|&=
\bigcup_{k\in\{i+1,i+2,\ldots, \hat{m}+1\}}|\hat{\Omega}(k)|),
\end{split}\end{equation*}
where $X(j)$ and $Y(j)$ are subsets of $|\Phi|$ defined in Lemma~\ref{property of basic subdivisions}.14 and in\hfill\break Lemma~\ref{property of basic subdivisions2}.6 respectively.
\end{enumerate}
\end{lemma}

\begin{lemma}
\label{basic7}
Consider any $\Delta\in \Phi/H$ with $\dim\Delta\geq 2$.

Note that $\mathcal{F}(\Delta)$ is a flat regular fan over $N$ in $V$,
$\mathcal{F}(\Delta)\subset \Phi$, 
$|\mathcal{F}(\Delta)|=\Delta$,
$\dim \mathcal{F}(\Delta)=\dim\Delta\geq 2$, $H\in\mathcal{F}(\Delta)_1$, and $\mathcal{F}(\Delta)$ is starry with center in $H$.

Note that $\emptyset\neq(\mathcal{F}(\Delta)-(\mathcal{F}(\Delta)/H))_1\subset(\Phi-(\Phi/H))_1$ and \hfill\break
$E^{-1}((\mathcal{F}(\Delta)-(\mathcal{F}(\Delta)/H))_1)
\subset\{1,2,\ldots,m\}$.

Let $\hat{m}=\sharp E^{-1}((\mathcal{F}(\Delta)-(\mathcal{F}(\Delta)/H))_1)\in\Z_0$.
Let $\hat{\tau}:\{1,2,\ldots,\hat{m}\}\rightarrow\{1,2,\ldots,m\}$ be the unique injective mapping preserving the order and satisfying $\hat{\tau}(\{1,2,\ldots,$\hfill\break$\hat{m}\})= E^{-1}((\mathcal{F}(\Delta)-(\mathcal{F}(\Delta)/H))_1)$.
Putting $\hat{\tau}(0)=0$ and  $\hat{\tau}(\hat{m}+1)=m+1$, we define an extension
$\hat{\tau}:\{0,1,2,\ldots,\hat{m},\hat{m}+1\}\rightarrow\{0,1,2,\ldots,m,m+1\}$ of $\hat{\tau}:\{1,2,\ldots,\hat{m}\}\rightarrow\{1,2,\ldots,m\}$.
Let $\hat{E}:\{1,2,\ldots,\hat{m}\}\rightarrow(\mathcal{F}(\Delta)-(\mathcal{F}(\Delta)/H))_1$ be the unique mapping satisfying $\iota\hat{E}=E\hat{\tau}$, where $\iota: (\mathcal{F}(\Delta)-(\mathcal{F}(\Delta)/H))_1\rightarrow(\Phi-(\Phi/H))_1$ denotes the inclusion mapping.

Let
$H\Op=H\Op|\Delta\in\mathcal{F}(\Delta)^1$,
$\hat{\Omega}=\Omega(V, N, H, \mathcal{F}(\Delta), \hat{m}, \hat{E})$,\hfill\break
$\hat{s}=s(V, N, H, \mathcal{F}(\Delta), \hat{m}, \hat{E})$,
$\hat{F}=F(V, N, H, \mathcal{F}(\Delta), \hat{m}, \hat{E})$,\hfill\break
$\hat{G}=G(V, N, H, \mathcal{F}(\Delta), \hat{m}, \hat{E})$,
$\hat{H}=H(V, N, H, \mathcal{F}(\Delta), \hat{m}, \hat{E})$,\hfill\break and
$\hat{\Omega}=\Omega(V, N, H, \mathcal{F}(\Delta), \hat{m}, \hat{E})$.

For any $i\in\{1,2,\ldots,\hat{m}+1\}$, we denote
$$\Theta(i)=|\hat{\Omega}(i)|\subset\Delta,\text{ and}$$
$$\bar{\Theta}(i)= |\hat{\Omega}(i)-(\hat{\Omega}(i)/\hat{H}(i))|\subset\Theta(i).$$

\begin{enumerate}
\item
If $\dim\Delta=2$, then $\sharp\mathcal{F}(H\Op)_1=1$, and for the unique element $\bar{E}\in\mathcal{F}(H\Op)_1$ and any $i\in\{0,1,\ldots,\hat{m}\}$, $\hat{s}(i, \bar{E})=i$.
\item
For any $i\in\{1,2,\ldots,\hat{m}, \hat{m}+1\}$, $\bar{\Theta}(i)$, $\bar{\Theta}(i)+H$ and $\Theta(i)$ are regular cones over $N$ in $V$,
$\dim \bar{\Theta}(i)=\dim \Delta-1$, $\dim(\bar{\Theta}(i)+H)=\dim \Theta(i)=\dim\Delta$,
$\bar{\Theta}(i)\subset \Theta(i)\subset\bar{\Theta}(i)+H\subset\Delta$,
$\Theta(i)+H =\bar{\Theta}(i)+H$,
$\hat{H}(i)\in\mathcal{F}(\Theta(i))_1$,
$\bar{\Theta}(i)=\hat{H}(i)\Op|\Theta(i)\in\mathcal{F}(\Theta(i))^1$,
$\Theta(i)=\bar{\Theta}(i)+\hat{H}(i)\subset\Delta$,
$\Vect(\Theta(i))=\Vect(\Delta)$,
$\Vect(H)\cap\Vect(\bar{\Theta}(i))=\{0\}$,
$\Vect(H)+\Vect(\bar{\Theta}(i))=\Vect(\Delta)$,
$\bar{\Theta}(i)=\Vect(\bar{\Theta}(i))\cap\Delta$, 
$\Theta(i)+\Vect(H) =\bar{\Theta}(i)+\Vect(H)=\Delta+\Vect(H)$,
$\hat{\Omega}(i)=\mathcal{F}(\Theta(i))$,
$\hat{\Omega}(i)-(\hat{\Omega}(i)/\hat{H}(i)) =\mathcal{F}(\bar{\Theta}(i))$,
$\bar{\Theta}(i)=\Convcone(\{b_{E/N^*}+\hat{s}(i-1, E)b_{H/N^*}| E\in\mathcal{F}(H\Op)_1\})\subset\Delta$, and
$(\mathcal{F}(\Delta)*\hat{F}(1)*\hat{F}(2)*\cdots*\hat{F}(i-1))\Mx=
\{\Theta(j)|j\in\{1,2,\ldots,i-1\}\}\cup\{\bar{\Theta}(i)+H\}$.
\item
Consider any $i\in\{1,2,\ldots,\hat{m}, \hat{m}+1\}$.

For any $\bar{E}\in\mathcal{F}(H\Op)_1$, $\bar{E}+H\in\mathcal{F}(\Delta)_2/H$,
$\bar{\Theta}(i)\cap(\bar{E}+H)\in\mathcal{F}(\bar{\Theta}(i))_1$.
The mapping from $\mathcal{F}(H\Op)_1$ to $\mathcal{F}(\bar{\Theta}(i))_1$ sending  $\bar{E}\in\mathcal{F}(H\Op)_1$ to $\bar{\Theta}(i)\cap(\bar{E}+H)\in\mathcal{F}(\bar{\Theta}(i))_1$ is a bijective mapping.

$H\subset\Theta(i)\Leftrightarrow i=\hat{m}+1$.

$\bar{\Theta}(i)^\circ\not\subset\Delta^\circ\Leftrightarrow\bar{\Theta}(i)\subset\partial\Delta\Leftrightarrow\bar{\Theta}(i)=H\Op\Leftrightarrow i=1$.
\item
For any $i\in\{1,2,\ldots,\hat{m}\}$ and any $j\in\{2,3,\ldots,\hat{m}+1\}$ with $i<j$, 
$\Theta(i)\cap\Theta(j)=\Theta(i)\cap(\bar{\Theta}(j)+H)=\bar{\Theta}(i+1)\cap\bar{\Theta}(j)$.
\item
Consider any $i\in\{1,2,\ldots,\hat{m}\}$.
$\hat{F}(i)=(\bar{\Theta}(i)+H)\cap(\hat{E}(i)+H)\in\mathcal{F}(\bar{\Theta}(i)+H)_2/H$.
$\bar{\Theta}(i+1)=\hat{G}(i)\Op|\Theta(i)\in\mathcal{F}(\Theta(i))^1$.
$\hat{H}(i)=\bar{\Theta}(i+1)\cap(\hat{G}(i)+\hat{H}(i)) \in\mathcal{F}(\Theta(i))_1$.
$\bar{\Theta}(i)\cap\bar{\Theta}(i+1)=(\hat{G}(i)+\hat{H}(i))\Op|\Theta(i)$.
$\{\Lambda\in\mathcal{F}(\Theta(i))|$\hfill\break$\Lambda^\circ\subset\Delta^\circ\cup (H\Op)^\circ \}=\{\Theta(i), \bar{\Theta}(i), \bar{\Theta}(i+1)\}$.
For any $\bar{E}\in\mathcal{F}(H\Op)_1-\{\hat{E}(i)\}$.
$\Theta(i)\cap(\bar{E}+H)= \bar{\Theta}(i)\cap\bar{\Theta}(i+1)\cap(\bar{E}+H)\in\mathcal{F}(\bar{\Theta}(i)\cap\bar{\Theta}(i+1))_1=\mathcal{F}(\Theta(i))_1-\{\hat{G}(i), \hat{H}(i)\}$.
The mapping from $\mathcal{F}(H\Op)_1-\{\hat{E}(i)\}$ to $\mathcal{F}(\Theta(i))_1-\{\hat{G}(i), \hat{H}(i)\}$ sending $\bar{E}\in\mathcal{F}(H\Op)_1-\{\hat{E}(i)\}$ to $\Theta(i)\cap(\bar{E}+H)\in\mathcal{F}(\Theta(i))_1-\{\hat{G}(i), \hat{H}(i)\}$ is a bijective mapping.
$\Theta(i)\cap(\hat{E}(i)+H)= \hat{G}(i)+\hat{H}(i) \in\mathcal{F}(\Theta(i))_2$.
\item
$H=\hat{H}(\hat{m}+1)\in\mathcal{F}(\Theta(\hat{m}+1))_1$.
$\bar{\Theta}(\hat{m}+1)=H\Op|\Theta(\hat{m}+1)\in\mathcal{F}(\Theta(\hat{m}+1))^1$.
$\{\Lambda\in\mathcal{F}(\Theta(\hat{m}+1))|\Lambda^\circ\subset\Delta^\circ\cup (H\Op)^\circ \}=\{\Theta(\hat{m}+1), \bar{\Theta}(\hat{m}+1)\}$.
For any $\bar{E}\in\mathcal{F}(H\Op)_1$, $\Theta(\hat{m}+1)\cap(\bar{E}+H)\in\mathcal{F}(\Theta(\hat{m}+1))_2/H$.
The mapping from $\mathcal{F}(H\Op)_1$ to $\mathcal{F}(\Theta(\hat{m}+1))_2/H$ sending $\bar{E}\in\mathcal{F}(H\Op)_1$ to $\Theta(\hat{m}+1)\cap(\bar{E}+H)\in\mathcal{F}(\Theta(\hat{m}+1))_2/H$ is a bijective mapping.
\item
$\hat{\Omega}$ is the iterated star subdivision of $\mathcal{F}(\Delta)$, it is a flat regular fan, and it is determined by the sextuplet $(V, N, H, \mathcal{F}(\Delta), \hat{m}, \hat{E})$ uniquely.
$|\hat{\Omega}|=\Delta$.

$\hat{\Omega}$ is $H$-simple.
Let $\smash{\bar{\hat{\Omega}}}^1=\{\Lambda\in\hat{\Omega}^1|\Lambda^\circ\subset\Delta^\circ\}\cup\{H\Op\}$ denote the $H$-skeleton of $\hat{\Omega}$.

$\sharp\hat{\Omega}^0=\sharp\smash{\bar{\hat{\Omega}}}^1=\hat{m}+1$.
$\hat{\Omega}^0=\{\Theta(i)|i\in\{1,2,\ldots,\hat{m},\hat{m}+1\}\}$.
$\smash{\bar{\hat{\Omega}}}^1=\{\bar{\Theta}(i)|i\in\{1,2,\ldots,\hat{m},\hat{m}+1\}\}$.

We consider the $H$-order on $\hat{\Omega}^0$. The bijective mapping from $\{1,2,\ldots,\hat{m},$\hfill\break$\hat{m}+1\}$ to $\hat{\Omega}^0$ sending $i\in\{1,2,\ldots,\hat{m},\hat{m}+1\}$ to $\Theta(i)\in \hat{\Omega}^0$ preserves the $H$-order.

We consider the $H$-order on $\smash{\bar{\hat{\Omega}}}^1$. The bijective mapping from $\{1,2,\ldots,\hat{m},$\hfill\break$\hat{m}+1\}$ to $\smash{\bar{\hat{\Omega}}}^1$ sending $i\in\{1,2,\ldots,\hat{m},\hat{m}+1\}$ to $\bar{\Theta}(i)\in\smash{\bar{\hat{\Omega}}}^1$ preserves the $H$-order.

Consider any $\bar{E}\in\mathcal{F}(H\Op)_1$ and any $i\in\{1,2,\ldots,\hat{m},\hat{m}+1\}$.
The structure constant of $\hat{\Omega}$ corresponding to the pair $(i,\bar{E})$ is equal to $\hat{s}(i-1, \bar{E})\in\Z_0$.
\end{enumerate}
\end{lemma}

\section{Upper boundaries and lower boundaries}
\label{upper}
Upper boundaries and lower boundaries are applied to define characteristic functions in the next section.

Let $V$ be any vector space of finite dimension over $\R$ with $\dim V\geq 1$, let $N$ be any lattice of $V$, let $H$ be any regular cone over $N$ in $V$ with $\dim H=1$ and let $\pi_H:V\rightarrow V/\Vect(H)$ denote the canonical surjective homomorphism of vector spaces over $\R$ to the residue vector space $V/\Vect(H)$.
 
\begin{definition}
Let $X$ be any subset of $V$. We denote
\begin{equation*}\begin{split}
\partial^H_+X&=\{a\in X|(\{a\}+\Vect(H))\cap X\subset\{a\}+(-H)\},\\
\partial^H_-X&=\{a\in X|(\{a\}+\Vect(H))\cap X\subset\{a\}+H\},
\end{split}\end{equation*}
and we call $\partial^H_+X$ and $\partial^H_-X$ the $H$-\emph{upper boundary} of $X$ and the $H$-\emph{lower boundary} of $X$ respectively.
\end{definition}

\begin{lemma} Let $\Delta$ be any convex polyhedral cone in $V$ satisfying $\Vect(H)\subset\Vect(\Delta)$.

\label{upper and lower}
\begin{enumerate}
\item
\begin{equation*}\begin{split}
\partial^H_+\Delta&=|\{\Lambda\in\mathcal{F}(\Delta)|H\not\subset\Delta+\Vect(\Lambda)\}|\\
&=|\{\Lambda\in\mathcal{F}(\Delta)|H\not\subset\Delta+\Vect(\Lambda), \dim\Lambda=\dim\Delta-1\}|,\\
\partial^H_-\Delta&=|\{\Lambda\in\mathcal{F}(\Delta)|-H\not\subset\Delta+\Vect(\Lambda)\}|\\
&=|\{\Lambda\in\mathcal{F}(\Delta)|-H\not\subset\Delta+\Vect(\Lambda), \dim\Lambda=\dim\Delta-1\}|,\\
\partial^H_+\Delta\cup\partial^H_-\Delta&=|\{\Lambda\in\mathcal{F}(\Delta)|\Vect(H)\not\subset\Vect(\Lambda)\}|\\
&=|\{\Lambda\in\mathcal{F}(\Delta)|\Vect(H)\not\subset\Vect(\Lambda), \dim\Lambda=\dim\Delta-1\}|\\
&\subset\partial\Delta.
\end{split}\end{equation*}
\item
$\partial^H_+\Delta=\emptyset\Leftrightarrow H\subset\Delta$. 
$\partial^H_-\Delta=\emptyset\Leftrightarrow -H\subset\Delta$.
\item
If $\partial^H_+\Delta\neq\emptyset$, then $\pi_H(\partial^H_+\Delta)=\pi_H(\Delta)$ and the mapping $\pi_H: \partial^H_+\Delta\rightarrow\pi_H(\Delta)$ induced by $\pi_H$ is a homeomorphism.

If $\partial^H_-\Delta\neq\emptyset$, then $\pi_H(\partial^H_-\Delta)=\pi_H(\Delta)$ and the mapping $\pi_H: \partial^H_-\Delta\rightarrow\pi_H(\Delta)$ induced by $\pi_H$ is a homeomorphism.
\item
$$\partial^H_+\Delta\cap\partial^H_-\Delta=
|\{\Lambda\in\mathcal{F}(\Delta)|\Lambda\subset\partial^H_+\Delta\cap\partial^H_-\Delta, \dim\Lambda\leq\dim\Delta-2\}|.$$
\item
Consider any $a\in V$.
Note that $\pi_H(a)\in V/\Vect(H)$ and $\pi_H^{-1}(\pi_H(a))=\{a\}+\R b_{H/N}$.
\begin{enumerate}
\item
If $\pi_H(a)\in\pi_H(\Delta)$ and $\partial^H_+\Delta\neq\emptyset$, then there exists uniquely a real number $t_+\in\R$ satisfying $a+t_+b_{H/N}\in\partial^H_+\Delta$.
\end{enumerate}

In the case where $\pi_H(a)\in\pi_H(\Delta)$ and $\partial^H_+\Delta\neq\emptyset$, we take the unique $t_+\in\R$ satisfying $a+t_+b_{H/N}\in\partial^H_+\Delta$.

\begin{enumerate}
\setcounter{enumii}{1}
\item
If $\pi_H(a)\in\pi_H(\Delta)$ and $\partial^H_-\Delta\neq\emptyset$, then there exists uniquely a real number $t_-\in\R$ satisfying $a+t_-b_{H/N}\in\partial^H_-\Delta$.
\end{enumerate}

In the case where $\pi_H(a)\in\pi_H(\Delta)$ and $\partial^H_-\Delta\neq\emptyset$, we take the unique $t_-\in\R$ satisfying $a+t_-b_{H/N}\in\partial^H_-\Delta$.

\begin{enumerate}
\setcounter{enumii}{2}
\item
If $\pi_H(a)\in\pi_H(\Delta)$, $\partial^H_+\Delta\neq\emptyset$ and $\partial^H_-\Delta\neq\emptyset$, then $t_-\leq t_+$.
\item
$(\{a\}+\R b_{H/N})\cap\Delta$
\begin{equation*}
=
\begin{cases}
\emptyset&\text{if $\pi_H(a)\not\in\pi_H(\Delta)$},\\
\{a+t b_{H/N}|t\in\R, t_-\leq t\leq t_+\}&\text{if $\pi_H(a)\in\pi_H(\Delta)$, $\partial^H_+\Delta\neq\emptyset$ and $\partial^H_-\Delta\neq\emptyset$},\\
\{a+t b_{H/N}|t\in\R, t\leq t_+\}&\text{if $\pi_H(a)\in\pi_H(\Delta)$, $\partial^H_+\Delta\neq\emptyset$ and $\partial^H_-\Delta=\emptyset$},\\
\{a+t b_{H/N}|t\in\R, t_-\leq t\}&\text{if $\pi_H(a)\in\pi_H(\Delta)$, $\partial^H_+\Delta=\emptyset$ and $\partial^H_-\Delta\neq\emptyset$},\\
\{a+t b_{H/N}|t\in\R \}&\text{if $\pi_H(a)\in\pi_H(\Delta)$, $\partial^H_+\Delta=\emptyset$ and $\partial^H_-\Delta=\emptyset$}.
\end{cases}
\end{equation*}
\item
If $\pi_H(a)\in\pi_H(\Delta)^\circ$, $\partial^H_+\Delta\neq\emptyset$ and $\partial^H_-\Delta\neq\emptyset$, then $t_-< t_+$.
\item
$(\{a\}+\R b_{H/N})\cap\Delta^\circ$
\begin{equation*}
=
\begin{cases}
\emptyset&\text{if $\pi_H(a)\not\in\pi_H(\Delta)^\circ$},\\
\{a+t b_{H/N}|t\in\R, t_-< t< t_+\}&\text{if $\pi_H(a)\in\pi_H(\Delta)^\circ$, $\partial^H_+\Delta\neq\emptyset$ and $\partial^H_-\Delta\neq\emptyset$},\\
\{a+t b_{H/N}|t\in\R, t< t_+\}&\text{if $\pi_H(a)\in\pi_H(\Delta)^\circ$, $\partial^H_+\Delta\neq\emptyset$ and $\partial^H_-\Delta=\emptyset$},\\
\{a+t b_{H/N}|t\in\R, t_-< t\}&\text{if $\pi_H(a)\in\pi_H(\Delta)^\circ$, $\partial^H_+\Delta=\emptyset$ and $\partial^H_-\Delta\neq\emptyset$},\\
\{a+t b_{H/N}|t\in\R \}&\text{if $\pi_H(a)\in\pi_H(\Delta)^\circ$, $\partial^H_+\Delta=\emptyset$ and $\partial^H_-\Delta=\emptyset$}.
\end{cases}\end{equation*}
\end{enumerate}
\end{enumerate}
\end{lemma}

\section{Height, characteristic functions and compatible mappings}
\label{compatible}
We define the height of a convex pseudo polytope in general situation. Characteristic functions and compatible mappings are applied for constructing a special basic subdivision associated with a given convex pseudo polytope.

In this section we consider the following objects: Let $V$ be any vector space of finite dimension over $\R$ with $\dim V\geq 2$, let $N$ be any lattice of $V$, let $H$ be any regular cone of dimension one over the dual lattice $N^*$ of $N$ in the dual vector space $V^*$ of $V$, let $\Phi$ be any flat regular fan over $N^*$ in $V^*$ satisfying 
$\dim \Phi\geq 2$, $H\in\Phi_1$ and $\Phi$ is starry with center in $H$, let $S$ be any rational convex pseudo polytope over $N$ in $V$ satisfying $\dim(|\Sigma(S|V)|)\geq 2$ and $|\Phi|\subset |\Sigma(S|V)|$, where $\Sigma(S|V)$ denotes the normal fan of $S$, and let $T$ be any convex pseudo polytope in $V$ satisfying $\dim(|\Sigma(T|V)|)\geq 1$ and $H\subset\Vect(|\Sigma(T|V)|)$, where $\Sigma(T|V)$ denotes the normal fan of $T$.

Note that
$b_{H/N^*}\in H\subset|\Phi|\subset|\Sigma(S|V)|\subset
\Vect(|\Sigma(S|V)|)$.

\begin{lemma}
\label{finiteness rationality}
\begin{enumerate}
\item
$\{\langle b_{H/N^*}, a\rangle|a\in\mathcal{V}(T)\}$ is a non-empty finite subset of $\R$, where $\mathcal{V}(T)$ denotes the skeleton of $T$.
\item
$0<\sharp\{\langle b_{H/N^*}, a\rangle|a\in\mathcal{V}(T)\}\leq c(T)$, where $c(T)$ denotes the characteristic number of $T$.
\item
If $T$ is rational over $N$, then the subset
$$\{m\in\Z|ma\in N+(\Vect(|\Sigma(T|V)|)^\vee|V^*)\text{ for any }a\in\mathcal{V}(T)\}$$
of $\Z$ is an ideal of the ring $\Z$ containing a positite integer.
\end{enumerate}
\end{lemma}

\begin{proof}
It follows from Theorem~\ref{property of polyhedrons}.7.
\end{proof}

We denote $\ell=\dim\Vect(|\Sigma(S|V)|)^\vee|V^*=\dim(\Stab(S)\cap(-\Stab(S)))\in\Z_0$.

\begin{definition}
\label{relative}
\begin{enumerate}
\item
We define functions
$$\lceil\ \rceil,\  \lfloor\ \rfloor:\R\rightarrow\Z,$$
by putting
$$\lceil r\rceil=\min\{i\in\Z|r\leq i\},\quad
\lfloor r\rfloor=\max\{i\in\Z|r\geq i\},$$
for any $r\in\R$.
\item
We define
$$\mathcal{H}(H, T)=\{\langle b_{H/N^*},a\rangle|a\in \mathcal{V}(T)\}\subset\R,$$
$$\Ht(H,T)=\max\mathcal{H}(H, T)-\min\mathcal{H}(H, T)\in\R_0,$$
and we call
$\mathcal{H}(H, T)$ and
$\Ht(H,T)$,
the $H$-\emph{height set} of $T$ and
the $H$-\emph{height} of $T$ respectively.
\item
Assume that $T$ is rational over $N$. 
By $\Den(T/N)$ we denote the minimum positive integer in the ideal
$$\{m\in\Z|ma\in N+(\Vect(|\Sigma(T|V)|)^\vee|V^*)\text{ for any }a\in\mathcal{V}(T)\},$$
and we call $\Den(T/N)\in\Z_+$ the \emph{denominator} of $T$ over $N$.
\item
We define
\begin{equation*}\begin{split}
(\Phi, \mathcal{F}(S)_\ell)=
\{F\in\mathcal{F}(S)_\ell|&\dim(\Delta(F,S|V)\cap\Delta)=\dim \Delta\text{ for some }\Delta\in\Phi\Mx.\}\\
&\qquad\qquad\subset \mathcal{F}(S)_\ell,
\end{split}\end{equation*}
$$\mathcal{V}(\Phi, S)=
\bigcup_{F\in(\Phi, \mathcal{F}(S)_\ell)}F\subset\mathcal{V}(S),$$
$$\mathcal{H}(H, \Phi, S)=\{\langle b_{H/N^*},a\rangle|a\in \mathcal{V}(\Phi, S)\}\subset\R,$$
$$\Ht(H,\Phi,S)=\max\mathcal{H}(H, \Phi, S)-\min\mathcal{H}(H, \Phi, S)\in\R_0,$$
and we call $(\Phi, \mathcal{F}(S)_\ell)$,
$\mathcal{V}(\Phi, S)$,
$\mathcal{H}(H, \Phi, S)$ and
$\Ht(H,\Phi,S)$,
the \emph{set of minimal faces} of the pair $(\Phi, S)$,
the \emph{skeleton} of the pair $(\Phi, S)$,
the $H$-\emph{height set} of the pair $(\Phi, S)$ and
the $H$-\emph{height} of the pair $(\Phi, S)$ respectively.
\item
For any $h\in \mathcal{H}(H, \Phi, S)$ we denote
\begin{equation*}\begin{split}
\Pi(h)&=
\{\Delta(F, S|V)\cap\Delta|F\in(\Phi, \mathcal{F}(S)_\ell),
\Delta\in\Phi\Mx,
\dim(\Delta(F, S|V)\cap\Delta)=\dim \Delta,\\
&\qquad\quad\langle b_{H/N^*}, a\rangle=h\text{ for some }a\in F\}\Fc\subset\Sigma(S|V)\hat{\cap}\Phi,\\
\Sigma(h)&=
\{\Delta(F, S|V)\cap\Delta|F\in(\Phi, \mathcal{F}(S)_\ell),
\Delta\in\Phi\Mx,
\dim(\Delta(F, S|V)\cap\Delta)=\dim \Delta,\\
&\qquad\quad\langle b_{H/N^*}, a\rangle\geq h\text{ for some }a\in F\}\Fc\subset\Sigma(S|V)\hat{\cap}\Phi.
\end{split}\end{equation*}
\end{enumerate}
\end{definition}

Consider any $\Delta\in\Phi/H$.
$S+(\Delta^\vee|V^*)$ is a rational conex pseudo polytope over $N$ in $V$.
$\Stab(S+(\Delta^\vee|V^*))=\Delta^\vee|V^*$.
$H\subset |\Sigma(S+(\Delta^\vee|V^*)|V)|=\Stab(S+(\Delta^\vee|V^*))^\vee|V=\Delta\subset|\Phi|\subset|\Sigma(S|V)|\subset\Vect(|\Sigma(S|V)|)$.
$\Sigma(S+(\Delta^\vee|V^*)|V)= \Sigma(S|V)\hat{\cap}\mathcal{F}(\Delta)$, where $\Sigma(S+(\Delta^\vee|V^*)|V)$ denotes the normal fan of $S+(\Delta^\vee|V^*)$.

\begin{lemma}
\label{compatible1}
\begin{enumerate}
\item
The set $\mathcal{H}(H, T)$ is a non-empty finite subset of $\R$.
$0<\sharp\mathcal{H}(H, T)\leq c(T)$.

If $T$ is rational over $N$, then $\mathcal{H}(H, T)\subset(1/\Den(T/N))\Z\subset\Q$ and\hfill\break
$\Ht(H,T)\in(1/\Den(T/N))\Z_0\subset\Q_0$.

If $T$ is rational over $N$ and $\mathcal{V}(T)\subset N+(\Vect(|\Sigma(T|V)|)^\vee|V^*)$, then $\Den(T/N)=1$, $\mathcal{H}(H, T)\subset\Z$ and $\Ht(H,T)\in\Z_0$.
\item
$\emptyset\neq(\Phi,\mathcal{F}(S)_\ell)\subset\mathcal{F}(S)_\ell$.

$\emptyset\neq\mathcal{V}(\Phi, S)\subset\mathcal{V}(S)$. The set $\mathcal{V}(\Phi, S)$ is the union of some connected components of $\mathcal{V}(S)$.

$\emptyset\neq\mathcal{H}(H, \Phi, S)\subset\mathcal{H}(H, S)\subset 
(1/\Den(S/N))\Z\subset\Q$. The set $\mathcal{H}(H, \Phi, S)$ is a non-empty finite subset of $\Q$.

$\min \mathcal{H}(H, \Phi, S)=\min\mathcal{H}(H, S)$.
$\Ht(H,\Phi,S)\leq \Ht(H,S)$.

$\Ht(H,S)\in(1/\Den(S/N))\Z_0\subset\Q_0$.

$\Ht(H,\Phi,S)\in(1/\Den(S/N))\Z_0\subset\Q_0$.

If $\mathcal{V}(S)\subset N+(\Vect(|\Sigma(S|V)|)^\vee|V^*)$, then
$\Den(S/N)=1$,
$\mathcal{H}(H, \Phi, S)\subset\mathcal{H}(H, S)\subset \Z$, $\Ht(H,S)\in\Z_0$ and $\Ht(H,\Phi,S)\in\Z_0$.
\item
Consider any $\Delta\in\Phi/H$. 

$\mathcal{H}(H, S+(\Delta^\vee|V^*))\subset \mathcal{H}(H, \Phi, S)\subset\Q$.

$\min\mathcal{H}(H, S+(\Delta^\vee|V^*))=\min\mathcal{H}(H, \Phi, S)$.

$\Ht(H,S+(\Delta^\vee|V^*))\leq\Ht(H,\Phi,S)$.

$\Den(S/N)$ is a multiple of $\Den(S+(\Delta^\vee|V^*)/N)$.

$\mathcal{H}(H, S+(\Delta^\vee|V^*))\subset(1/\Den(S+(\Delta^\vee|V^*)/N))\Z\subset(1/\Den(S/N))\Z\subset \Q$.

$\Ht(H,S+(\Delta^\vee|V^*))\in(1/\Den(S+(\Delta^\vee|V^*)/N))\Z_0\subset (1/\Den(S/N))\Z_0$\break$\subset \Q_0$.

If $\mathcal{V}(S)\subset N+(\Vect(|\Sigma(S|V)|)^\vee|V^*)$, then 
$\Den(S+(\Delta^\vee|V^*)/N)=1$,
$\mathcal{H}(H, S+(\Delta^\vee|V^*))\subset\Z$ and $\Ht(H,S+(\Delta^\vee|V^*))\in\Z_0$.
\item
$\max\mathcal{H}(H, \Phi, S)=\max\{\max\mathcal{H}(H, S+(\Delta^\vee|V^*))|\Delta\in \Phi\Mx\}=\max\{\max\mathcal{H}(H,$\break$ S+(\Delta^\vee|V^*))|\Delta\in \Phi/H\}$.

$\Ht(H,\Phi,S)= \max\{\Ht(H,S+(\Delta^\vee|V^*))|\Delta\in \Phi\Mx\}= \max\{\Ht($\break$H,S+(\Delta^\vee|V^*))|\Delta\in \Phi/H\}$.
\item
Consider any subset $\hat{\Phi}$ of $\Phi$ satisfying 
$\dim\hat{\Phi}=\dim \Vect(|\hat{\Phi}|)\geq 2$, $\hat{\Phi}\Mx=\hat{\Phi}^0$,
$H\in\hat{\Phi}_1$
and $\hat{\Phi}=(\hat{\Phi}/H)\Fc$.

Note that $\hat{\Phi}$ is a flat regular fan over $N^*$ in $V^*$ and it is starry with center in $H$, and
$|\hat{\Phi}|\subset|\Phi|\subset|\Sigma(S|V)|$.

$\mathcal{V}(\hat{\Phi}, S)\subset\mathcal{V}(\Phi, S)$.
$\mathcal{H}(H, \hat{\Phi},S)\subset\mathcal{H}(H, \Phi, S)$.
$\min\mathcal{H}(H, \hat{\Phi},S)=\min\mathcal{H}(H, $\break$\Phi, S)$.
$\Ht(H, \hat{\Phi}, S)\leq\Ht(H,\Phi,S)$.

\item
$\Sigma(S|V)\hat{\cap}\Phi$ is a flat rational fan over $N^*$ in $V^*$.
$\dim(\Sigma(S|V)\hat{\cap}\Phi)=
\dim \Phi$.
$\Vect(|\Sigma(S|V)\hat{\cap}\Phi|)=\Vect(|\Phi|)$.
\item
Consider any $h\in\mathcal{H}(H, \Phi, S)$.
$\Pi(h)$ and $\Sigma(h)$ are flat rational fans over $N^*$ in $V^*$.
$\dim \Pi(h)=\dim\Phi$.
$\Vect(|\Pi(h)|)=\Vect(|\Phi|)$.
$\dim \Sigma(h)=\dim\Phi$.
$\Vect(|\Sigma(h)|)= \Vect(|\Phi|)$.
$\Pi(h)\subset\Sigma(h)\subset\Sigma(S|V)\hat{\cap}\Phi$.

\item
Consider any $g\in\mathcal{H}(H, \Phi, S)$ and any $h\in\mathcal{H}(H, \Phi, S)$.
$g\leq h\Leftrightarrow \Sigma(g)\supset\Sigma(h)\Leftrightarrow |\Sigma(g)|\supset|\Sigma(h)|$.
$g=h\Leftrightarrow \Sigma(g)=\Sigma(h)\Leftrightarrow |\Sigma(g)|=|\Sigma(h)|$.
\item
We consider any $F\in\mathcal{F}(S)_\ell$ and any $G\in\mathcal{F}(S)_\ell$
such that $\langle b_{H/N^*}, a\rangle=\langle b_{H/N^*}, b\rangle$ for some $a\in F$ and some $b\in G$.
\begin{equation*}\begin{split}
\partial^H_+(\Delta(F,S|V)\cup\Delta(G,S|V))
&=\partial^H_+\Delta(F,S|V)\cup\partial^H_+\Delta(G,S|V),\text{ and}\\
\partial^H_-(\Delta(F,S|V)\cup\Delta(G,S|V))
&=\partial^H_-\Delta(F,S|V)\cup\partial^H_-\Delta(G,S|V).
\end{split}\end{equation*}
\end{enumerate}
\end{lemma}

Below in this section, we assume moreover that $\Sigma(S|V)\hat{\cap}\mathcal{F}(\Delta)$ is $H$-simple for any $\Delta\in\Phi\Mx$.
We denote
$$\max=\max\mathcal{H}(H, \Phi, S)\in(1/\Den(S/N))\Z, \ \min=\min\mathcal{H}(H, \Phi, S)\in(1/\Den(S/N))\Z.$$
By $\pi_H:V^*\rightarrow V^*/\Vect(H)$ we denote the canonical surjective homomorphism of vector spaces over $\R$ to the residue vector space $V^*/\Vect(H)$.
\begin{lemma}
\label{compatible2}
Consider any $h\in\mathcal{H}(H, \Phi, S)$.
\begin{enumerate}
\item
$\max\geq\min$. $\Ht(H,\Phi,S)=\max-\min$.
$\Ht(H,\Phi,S)=0\Leftrightarrow \max=\min\Leftrightarrow \Phi$ is a subdivision of $\Sigma(S|V)$.
\item
$\partial^H_+|\Pi(h)|\neq\emptyset\Leftrightarrow h\neq\min$.
$\partial^H_-|\Pi(h)|\neq\emptyset$.
\item
If $h\neq\min$, then $\pi_H(\partial^H_+|\Pi(h)|)=\pi_H(|\Pi(h)|)$ and the mapping $\pi_H:$\hfill\break $\partial^H_+|\Pi(h)|\rightarrow \pi_H(|\Pi(h)|)$ induced by $\pi_H$ is a homeomorphism.

$\pi_H(\partial^H_-|\Pi(h)|)=\pi_H(|\Pi(h)|)$ and the mapping $\pi_H: \partial^H_-|\Pi(h)|\rightarrow \pi_H(|\Pi(h)|)$ induced by $\pi_H$ is a homeomorphism.
\item
Consider any $a\in V^*$.
Note that $\pi_H(a)\in V^*/\Vect(H)$ and $\pi_H^{-1}(\pi_H(a))=\{a\}+\R b_{H/N^*}$.
\begin{enumerate}
\item
If $\pi_H(a)\in\pi_H(|\Pi(h)|)$ and $h\neq\min$, then there exists uniquely a real number $t_+\in\R$ satisfying $a+t_+b_{H/N^*}\in\partial^H_+|\Pi(h)|$.
\end{enumerate}

In the case where $\pi_H(a)\in\pi_H(|\Pi(h)|)$ and $h\neq\min$, we take the unique $t_+\in\R$ satisfying $a+t_+b_{H/N^*}\in\partial^H_+|\Pi(h)|$.

\begin{enumerate}
\setcounter{enumii}{1}
\item
If $\pi_H(a)\in\pi_H(|\Pi(h)|)$, then there exists uniquely a real number $t_-\in\R$ satisfying $a+t_-b_{H/N^*}\in\partial^H_-|\Pi(h)|$.
\end{enumerate}

In the case where $\pi_H(a)\in\pi_H(|\Pi(h)|)$, we take the unique $t_-\in\R$ satisfying $a+t_-b_{H/N^*}\in\partial^H_-|\Pi(h)|$.

\begin{enumerate}
\setcounter{enumii}{2}
\item
If $\pi_H(a)\in\pi_H(|\Pi(h)|)$ and $h\neq\min$, then $t_-\leq t_+$.
\item
$(\{a\}+\R b_{H/N^*})\cap|\Pi(h)|$
\begin{equation*}
=
\begin{cases}
\emptyset&\text{if $\pi_H(a)\not\in\pi_H(|\Pi(h)|)$},\\
\{a+t b_{H/N^*}|t\in\R, t_-\leq t\leq t_+\}&\text{if $\pi_H(a)\in\pi_H(|\Pi(h)|)$ and $h\neq\min$},\\
\{a+t b_{H/N^*}|t\in\R, t_-\leq t\}&\text{if $\pi_H(a)\in\pi_H(|\Pi(h)|)$ and $h=\min$}.
\end{cases}
\end{equation*}
\end{enumerate}
\item
$\partial^H_+|\Sigma(h)|\neq\emptyset\Leftrightarrow h\neq\min$.
$\partial^H_-|\Sigma(h)|= |\Sigma(h)|\cap|\Phi-(\Phi/H)|\neq\emptyset$.
\item
If $h\neq\min$, then $\pi_H(\partial^H_+|\Sigma(h)|)=\pi_H(|\Sigma(h)|)$ and the mapping $\pi_H:$\hfill\break $\partial^H_+|\Sigma(h)|\rightarrow \pi_H(|\Sigma(h)|)$ induced by $\pi_H$ is a homeomorphism.

$\pi_H(\partial^H_-|\Sigma(h)|)=\pi_H(|\Sigma(h)|)$ and the mapping $\pi_H: \partial^H_-|\Sigma(h)|\rightarrow \pi_H(|\Sigma(h)|)$ induced by $\pi_H$ is a homeomorphism.
\item
Consider any $a\in V^*$.
Note that $\pi_H(a)\in V^*/\Vect(H)$ and $\pi_H^{-1}(\pi_H(a))=\{a\}+\R b_{H/N^*}$.
\begin{enumerate}
\item
If $\pi_H(a)\in\pi_H(|\Sigma(h)|)$ and $h\neq\min$, then there exists uniquely a real number $t_+\in\R$ satisfying $a+t_+b_{H/N^*}\in\partial^H_+|\Sigma(h)|$.
\end{enumerate}

In the case where $\pi_H(a)\in\pi_H(|\Sigma(h)|)$ and $h\neq\min$, we take the unique $t_+\in\R$ satisfying $a+t_+b_{H/N^*}\in\partial^H_+|\Sigma(h)|$.

\begin{enumerate}
\setcounter{enumii}{1}
\item
If $\pi_H(a)\in\pi_H(|\Sigma(h)|)$, then there exists uniquely a real number $t_-\in\R$ satisfying $a+t_-b_{H/N^*}\in\partial^H_-|\Sigma(h)|$.
\end{enumerate}

In the case where $\pi_H(a)\in\pi_H(|\Sigma(h)|)$, we take the unique $t_-\in\R$ satisfying $a+t_-b_{H/N^*}\in\partial^H_-|\Sigma(h)|$.

\begin{enumerate}
\setcounter{enumii}{2}
\item
If $\pi_H(a)\in\pi_H(|\Sigma(h)|)$ and $h\neq\min$, then $t_-\leq t_+$.
\item
$(\{a\}+\R b_{H/N^*})\cap|\Sigma(h)|$
\begin{equation*}
=
\begin{cases}
\emptyset&\text{if $\pi_H(a)\not\in\pi_H(|\Sigma(h)|)$},\\
\{a+t b_{H/N^*}|t\in\R, t_-\leq t\leq t_+\}&\text{if $\pi_H(a)\in\pi_H(|\Sigma(h)|)$ and $h\neq\min$},\\
\{a+t b_{H/N^*}|t\in\R, t_-\leq t\}&\text{if $\pi_H(a)\in\pi_H(|\Sigma(h)|)$ and $h=\min$}.
\end{cases}
\end{equation*}
\end{enumerate}
\item
If $h=\max$, then $\Sigma(h)= \Pi(h)$.
\item
Assume $h\neq\max$.
Put $g=\min\{f\in\mathcal{H}(H, \Phi, S)|f>h\}\in\mathcal{H}(H, \Phi, S)$.

$\Sigma(h)=\Sigma(g)\cup\Pi(h)$.
$|\Sigma(h)|=|\Sigma(g)|\cup|\Pi(h)|$.
$\partial^H_+|\Sigma(g)|\neq\emptyset$.
$\partial^H_-|\Pi(h)|\neq\emptyset$.
$|\Sigma(g)|\cap|\Pi(h)|= \partial^H_+|\Sigma(g)|\cap\partial^H_-|\Pi(h)|$.
$\partial^H_-|\Pi(h)|\subset\partial^H_+|\Sigma(g)|\cup|\Phi-(\Phi/H)|$.
$\partial^H_+|\Sigma(h)|=( \partial^H_+|\Sigma(g)|- \partial^H_-|\Pi(h)|)\cup\partial^H_+|\Pi(h)|$.
\item
If $h=\min$, then $\Sigma(h)=\Sigma(S|V)\hat{\cap}\Phi$ and
$|\Sigma(h)|=| \Phi|$.
\item
$$\pi_H(|\Sigma(h)|)\subset|\pi_{H*}\Phi|.$$
$$\pi_H(|\Sigma(h)|)=
|\{\bar{\Delta}\in(\pi_{H*}\Phi)\Mx|\bar{\Delta}\subset\pi_H(|\Sigma(h)|)\}|.$$
$$\Clos(|\pi_{H*}\Phi|-\pi_H(|\Sigma(h)|))=
|\{\bar{\Delta}\in(\pi_{H*}\Phi)\Mx|\bar{\Delta}\subset\Clos(|\pi_{H*}\Phi|-\pi_H(|\Sigma(h)|))\}|.$$

\item
Consider any $\bar{E}\in(\Phi-(\Phi/H))_1$.
If $\Ht(H,\Phi,S)>0$ and $\pi_H(\bar{E})\subset\pi_H(\Sigma(\max))$, then there exists uniquely a real number $\gamma(\bar{E})\in\R$ satisfying 
$b_{\bar{E}/N^*}+\gamma(\bar{E})b_{H/N^*}\in\partial^H_+\Sigma(\max)$.
\end{enumerate}
\end{lemma}

Below we assume moreover that $\Ht(H,\Phi,S)>0$.

\begin{definition}
\label{defchar}
We define a function $\gamma: (\Phi-(\Phi/H))_1\rightarrow \R$.

Consider any $\bar{E}\in(\Phi-(\Phi/H))_1$.
If $\pi_H(\bar{E})\subset\pi_H(\Sigma(\max))$, then we take the unique real number $\gamma(\bar{E})\in\R$ satisfying 
$b_{\bar{E}/N^*}+\gamma(\bar{E})b_{H/N^*}\in\partial^H_+\Sigma(\max)$.
If $\pi_H(\bar{E})\not\subset\pi_H(\Sigma(\max))$, then we put $\gamma(\bar{E})=0\in\R$.

We call $\gamma$ the \emph{characteristic function} of the triplet $(H,\Phi,S)$.
\end{definition}

\begin{lemma}
\label{propchar}
Let $\gamma: (\Phi-(\Phi/H))_1\rightarrow \R$ denote the characteristic function of $(H, \Phi, S)$.
\begin{enumerate}
\item
For any $\bar{E}\in(\Phi-(\Phi/H))_1$, $\gamma(\bar{E})\in\Q_0$.
\item
Consider any $\Delta\in\Phi\Mx$. 
There exists $\bar{E}\in\mathcal{F}(H\Op|\Delta)_1$ satisfying $\gamma(\bar{E})>0\Leftrightarrow\pi_H(\Delta)\subset\pi_H(\Sigma(\max))\Leftrightarrow \pi_H(\Delta)\not\subset\Clos(|\pi_{H*}\Phi|-\pi_H(|\Sigma(\max)|))$.
\item
Consider any $\Delta\in\Phi\Mx$ satisfying $\pi_H(\Delta)\subset\pi_H(\Sigma(\max))$.
Note that $H\in\mathcal{F}(\Delta)_1$ and $\mathcal{F}(H\Op|\Delta)_1\subset (\Phi-(\Phi/H))_1$.

$\Sigma(S+(\Delta^\vee|V^*)|V)=\Sigma(S|V)\hat{\cap}\mathcal{F}(\Delta)$ is $H$-simple, and
$c(S+(\Delta^\vee|V^*))\geq 2$.

For any $\bar{E}\in\mathcal{F}(H\Op|\Delta)_1$, $\gamma(\bar{E})=c(\Sigma(S+(\Delta^\vee|V^*)|V), 2, \bar{E})$,
where \hfill\break$c(\Sigma(S+(\Delta^\vee|V^*)|V), 2, \bar{E})$ denotes the structure constant of $\Sigma(S+(\Delta^\vee|V^*)|V)$ corresponding to $(2, \bar{E})$.
\item
There exists $\bar{E}\in(\Phi-(\Phi/H))_1$ with $\gamma(\bar{E})>0$.
\item
For any $\bar{E}\in(\Phi-(\Phi/H))_1$ satisfying
$\pi_H(\bar{E})\subset\Clos(|\pi_{H*}\Phi|-\pi_H(|\Sigma(\max)|))$,
$\gamma(\bar{E})=0$.
\item
Consider any $\bar{E}\in(\Phi-(\Phi/H))_1$ satisfying $\gamma(\bar{E})\not\in\Z$.
Then, there exists uniquely $h(\bar{E})\in\mathcal{H}(H, \Phi, S)$ satisfying
\begin{equation*}\begin{split}
\{h\in\mathcal{H}(H, \Phi, S)&|
b_{\bar{E}/N^*}+\lceil\gamma(\bar{E})\rceil b_{H/N^*}\in\Sigma(h)-\partial^H_+\Sigma(h)\}\\
=&
\{h\in\mathcal{H}(H, \Phi, S)|h\leq h(\bar{E})\},\\
\end{split}\end{equation*}
and if $h(\bar{E})\in\mathcal{H}(H, \Phi, S)$ satisfies this equality, then $h(\bar{E})\neq\max$.
\end{enumerate}
\end{lemma}

\begin{definition}
\label{defcomp}
Let $\gamma: (\Phi-(\Phi/H))_1\rightarrow \R$ denote the characteristic function of $(H, \Phi, S)$.
For any $\bar{E}\in(\Phi-(\Phi/H))_1$ satisfying $\gamma(\bar{E})\not\in\Z$, we take the unique element  $h(\bar{E})\in\mathcal{H}(H, \Phi, S)$ satisfying
\begin{equation*}\begin{split}
\{h\in\mathcal{H}(H, \Phi, S)&|
b_{\bar{E}/N^*}+\lceil\gamma(\bar{E})\rceil b_{H/N^*}\in\Sigma(h)-\partial^H_+\Sigma(h)\}\\
=&
\{h\in\mathcal{H}(H, \Phi, S)|h\leq h(\bar{E})\}.\\
\end{split}\end{equation*}

Let
\begin{equation*}\begin{split}
m&=\sum_{\bar{E}\in(\Phi-(\Phi/H))_1}\lceil\gamma(\bar{E})\rceil,\\
\bar{m}&=\sum_{\bar{E}\in(\Phi-(\Phi/H))_1}\lfloor\gamma(\bar{E})\rfloor,\text{ and}\\
\mathcal{R}&=\{ \bar{E}\in(\Phi-(\Phi/H))_1|\gamma(\bar{E})\not\in\Z\}.
\end{split}\end{equation*}
Note that $m\in\Z_+$, $\bar{m}\in\Z_0$, $\bar{m}\leq m$, and
$m-\bar{m}=\sharp\mathcal{R}$.

Consider any mapping $E:\{1,2,\ldots,m\}\rightarrow (\Phi-(\Phi/H))_1$.

We say that the mapping $ E$ is \emph{compatible} with $S$, if the following three conditions are satisfied:
\begin{enumerate}
\item
For any $\bar{E}\in(\Phi-(\Phi/H))_1$,
$\sharp(\{1,2,\ldots,\bar{m}\}\cap E^{-1}(\bar{E}))= \lfloor\gamma(\bar{E})\rfloor$.
\item
$E(\{\bar{m}+1,\bar{m}+2,\ldots,m\})=\mathcal{R}$.
\item
If $m-\bar{m}\geq 2$, then $h(E(i))\geq h(E(i+1))$ for any $i\in\{\bar{m}+1,\bar{m}+2,\ldots,m-1\}$.
\end{enumerate}
\end{definition}

\begin{lemma}
\label{existcomp}
Let $\gamma: (\Phi-(\Phi/H))_1\rightarrow\R$ denote the characteristic function of $(H, \Phi, S)$.
Let $m=\sum_{\bar{E}\in(\Phi-(\Phi/H))_1}\lceil\gamma(\bar{E})\rceil\in\Z_+$,
$\bar{m}=\sum_{\bar{E}\in(\Phi-(\Phi/H))_1}\lfloor\gamma(\bar{E})\rfloor\in\Z_0$, and
$\mathcal{R}=\{\bar{E}\in(\Phi-(\Phi/H))_1|\gamma(\bar{E})\not\in\Z\}\subset (\Phi-(\Phi/H))_1$.
\begin{enumerate}
\item
There exists a compatible mapping $E:\{1,2,\ldots,m\}\rightarrow (\Phi-(\Phi/H))_1$
with $S$.
\item
Assume that a mapping $E:\{1,2,\ldots,m\}\rightarrow (\Phi-(\Phi/H))_1$
is compatible with $S$.
\begin{enumerate}
\item
For any $\bar{E}\in(\Phi-(\Phi/H))_1$,
$\sharp E^{-1}(\bar{E})=\lceil \gamma(\bar{E})\rceil$, and
$\sharp(\{1,2,\ldots, \bar{m}\}\cap E^{-1}(\bar{E}))=\lfloor \gamma(\bar{E})\rfloor$.

\item
$E(\{\bar{m}+1, \bar{m}+2,\ldots,m\})=\mathcal{R}$, and the mapping $E: \{\bar{m}+1, \bar{m}+2,\ldots,m\}\rightarrow \mathcal{R}$ induced by $E$ is bijective.
\item
$E(\{1,2,\ldots,m\})=\{\bar{E}\in(\Phi-(\Phi/H))_1|\gamma(\bar{E})>0\}$.
\item
For any $i\in\{1,2,\ldots,m\}$, $\pi_H(E(i))\subset\pi_H(\Sigma(\max))$ and
$\pi_H(E(i))\not\subset\Clos(|\pi_{H*}\Phi|-\pi_H(|\Sigma(\max)|))$.
\end{enumerate}
\end{enumerate}

Consider any subset $\hat{\Phi}$ of $\Phi$ satisfying $\dim\hat{\Phi}=\dim \Vect(|\hat{\Phi}|)\geq 2$, $\hat{\Phi}\Mx=\hat{\Phi}^0$,
$H\in\hat{\Phi}_1$, and $\hat{\Phi}=(\hat{\Phi}/H)\Fc$.

Note that
$\hat{\Phi}$ is a flat regular fan over $N^*$ in $V^*$ and it is starry with center in $H$,
$\emptyset\neq(\hat{\Phi}-(\hat{\Phi}/H))_1\subset(\Phi-(\Phi/H))_1$ and
$|\hat{\Phi}|\subset|\Phi|\subset|\Sigma(S|V)|$.
\begin{enumerate}
\setcounter{enumi}{2}
\item
$\Ht(H, \hat{\Phi}, S)\leq \Ht(H,\Phi,S)$.

$\Ht(H, \hat{\Phi}, S)=\Ht(H,\Phi,S)$, if and only if,
$\gamma(\bar{E})>0$ for some $\bar{E}\in (\hat{\Phi}-(\hat{\Phi}/H))_1$
\item
If $\gamma(\bar{E})>0$ for some $\bar{E}\in (\hat{\Phi}-(\hat{\Phi}/H))_1$, then $\Ht(H, \hat{\Phi}, S)>0$ and the composition $(\hat{\Phi}-(\hat{\Phi}/H))_1\rightarrow\R$ of the inclusion mapping $(\hat{\Phi}-(\hat{\Phi}/H))_1\rightarrow(\Phi-(\Phi/H))_1$ and $\gamma: (\Phi-(\Phi/H))_1\rightarrow\R$ coincides with the characteristic function of $(H,\hat{\Phi},S)$.
\item
Consider any compatible mapping: $E:\{1,2,\ldots,m\}\rightarrow (\Phi-(\Phi/H))_1$ with $S$.

Let $\hat{m}=\sharp E^{-1}((\hat{\Phi}-(\hat{\Phi}/H))_1)\in\Z_0$
and $\hat{\bar{m}}=\sharp(\{1,2,\ldots,\bar{m}\}\cap E^{-1}((\hat{\Phi}-(\hat{\Phi}/H))_1))\in\Z_0$.
Let $\hat{\tau}:\{1,2,\ldots,\hat{m}\}\rightarrow\{1,2,\ldots,m\}$ be the unique injective mapping preserving the order and satisfying $\hat{\tau}(\{1,2,\ldots,\hat{m}\})=$\hfill\break$E^{-1}((\hat{\Phi}-(\hat{\Phi}/H))_1)$.
Let $\hat{E}:\{1,2,\ldots,\hat{m}\}\rightarrow (\hat{\Phi}-(\hat{\Phi}/H))_1$ be the unique mapping satisfying $\iota\hat{E}=E\hat{\tau}$, where $\iota: (\hat{\Phi}-(\hat{\Phi}/H))_1\rightarrow(\Phi-(\Phi/H))_1$ denotes the inclusion mapping.

If $\gamma(\bar{E})>0$ for some $\bar{E}\in (\hat{\Phi}-(\hat{\Phi}/H))_1$, then $\hat{m}=\sum_{\bar{E}\in(\hat{\Phi}-(\hat{\Phi}/H))_1}\lceil\gamma(\bar{E})\rceil$, $\hat{\bar{m}}=\sum_{\bar{E}\in(\hat{\Phi}-(\hat{\Phi}/H))_1}\lfloor\gamma(\bar{E})\rfloor$ and $\hat{E}$ is compatible with $S$.

If $\gamma(\bar{E})=0$ for any $\bar{E}\in (\hat{\Phi}-(\hat{\Phi}/H))_1$, then $\hat{m}=0$.
\end{enumerate}
\end{lemma}

\begin{lemma}
\label{equiv2}
Assume $\sharp\Phi\Mx=1$.

Let $\Delta\in\Phi\Mx$ denote the unique element. 
$\Delta$ is a regular cone over $N^*$ in $V^*$, $H\in\mathcal{F}(\Delta)_1$ and $\dim\Delta=\dim \Phi\geq 2$, $\Phi=\mathcal{F}(\Delta)$, 
$\Phi-(\Phi/H)= \mathcal{F}(H\Op|\Delta)$, and
 $\Ht(H,S+(\Delta^\vee|V^*))=\Ht(H,\Phi,S)>0$.

Note that $\Sigma(S+(\Delta^\vee|V^*)|V)=\Sigma(S|V)\hat{\cap}\mathcal{F}(\Delta)$ is $H$-simple and $c(S+(\Delta^\vee|V^*))\geq 2$.

Let $\gamma, m, \bar{m}$ and $\mathcal{R}$ be the same as in above Lemma~\ref{existcomp}.

For any $i\in\{1,2,\ldots, c(S+(\Delta^\vee|V^*))\}$ and any $\bar{E}\in\mathcal{F}(H\Op|\Delta)_1$, we can consider the structure constant
$c(\Sigma(S+(\Delta^\vee|V^*)|V),i, \bar{E})\in\Q_0$ of $\Sigma(S+(\Delta^\vee|V^*)|V)$ corresponding to the pair $(i, \bar{E}).$

Denote
\begin{equation*}\begin{split}
\hat{m}&=\sum_{ \bar{E}\in\mathcal{F}(H\Op|\Delta)_1}\lceil c(\Sigma(S+(\Delta^\vee|V^*)|V),2, \bar{E})\rceil\in\Z_+,\\
\hat{\bar{m}}&=\sum_{ \bar{E}\in\mathcal{F}(H\Op|\Delta)_1}\lfloor c(\Sigma(S+(\Delta^\vee|V^*)|V),2, \bar{E})\rfloor\in\Z_0,\text{ and}\\
\hat{\mathcal{R}}&=\{ \bar{E}\in\mathcal{F}(H\Op|\Delta)_1|
c(\Sigma(S+(\Delta^\vee|V^*)|V),2, \bar{E})\not\in\Z\}
\subset\mathcal{F}(H\Op|\Delta)_1.
\end{split}\end{equation*}
For any $\bar{E}\in \hat{\mathcal{R}}$, denote
\begin{equation*}\begin{split}
\bar{c}(\bar{E})&=\max\{j\in\{2,3,\ldots, c(S+(\Delta^\vee|V^*))\}|\\
&\qquad\quad c(\Sigma(S+(\Delta^\vee|V^*)|V),j, \bar{E})<
\lceil c(\Sigma(S+(\Delta^\vee|V^*)|V),2, \bar{E})\rceil\}\\
&\qquad\quad\quad\in\{2,3,\ldots, c(S+(\Delta^\vee|V^*))\}.
\end{split}\end{equation*}

\begin{enumerate}
\item
$\gamma(\bar{E})=c(\Sigma(S+(\Delta^\vee|V^*)|V),2, \bar{E})$ for any $\bar{E}\in \mathcal{F}(H\Op|\Delta)_1$.
\item
$\hat{m}=m$.
$\hat{\bar{m}}=\bar{m}$.
$\hat{\mathcal{R}}=\mathcal{R}$.
\item
A mapping $E:\{1,2,\ldots,m\}\rightarrow (\Phi-(\Phi/H))_1=\mathcal{F}(H\Op|\Delta)_1$ is compatible with $S$, if and only if, the following three conditions are satisfied:
\begin{enumerate}
\item
For any $\bar{E}\in\mathcal{F}(H\Op|\Delta)_1$,
$\sharp(\{1,2,\ldots, \bar{m}\}\cap E^{-1}(\bar{E}))=$\hfill\break$\lfloor c(\Sigma(S+(\Delta^\vee|V^*)|V),2, \bar{E})\rfloor$.
\item
$E(\{\bar{m}+1, \bar{m}+2, \ldots, m\})=\mathcal{R}$.
\item
If $m-\bar{m}\geq 2$, then $\bar{c}(E(i))\leq\bar{c}(E(i+1))$ for any $i\in\{\bar{m}+1, \bar{m}+2, \ldots, m-1\}$.
\end{enumerate}
\item
Assume that a mapping $E:\{1,2,\ldots,m\}\rightarrow \mathcal{F}(H\Op|\Delta)_1$
is compatible with $S$.
\begin{enumerate}
\item
For any $\bar{E}\in\mathcal{F}(H\Op|\Delta)_1$,
$\sharp E^{-1}(\bar{E})=\lceil c(\Sigma(S+(\Delta^\vee|V^*)|V),2, \bar{E})\rceil$, and
$\sharp(\{1,2,\ldots, \bar{m}\}\cap E^{-1}(\bar{E}))=\lfloor c(\Sigma(S+(\Delta^\vee|V^*)|V),2, \bar{E})\rfloor$.

\item
$E(\{\bar{m}+1, \bar{m}+2,\ldots,m\})=\mathcal{R}$, and the mapping $E: \{\bar{m}+1, \bar{m}+2,\ldots,m\}\rightarrow \mathcal{R}$ induced by $E$ is bijective.
\item
$E(\{1,2,\ldots,m\})=\{\bar{E}\in\mathcal{F}(H\Op|\Delta)_1| c(\Sigma(S+(\Delta^\vee|V^*)|V),2, \bar{E})>0\}$.
\end{enumerate}
\end{enumerate}
\end{lemma}

\section{The height inequalities}
\label{height inequalities}
For a basic subdivision associated with a compatible mapping, the height inequalities hold. These iniqualities play an important role in our theory.

In this section we consider the following objects: Let $V$ be any vector space of finite dimension over $\R$ with $\dim V\geq 2$, let $N$ be any lattice of $V$, let $H$ be any regular cone of dimension one over the dual lattice $N^*$ of $N$ in the dual vector space $V^*$ of $V$, let $\Phi$ be any flat regular fan over $N^*$ in $V^*$ such that 
$\dim \Phi\geq 2$, $H\in\Phi_1$ and $\Phi$ is starry with center in $H$ and let $S$ be any rational convex pseudo polytope over $N$ in $V$ satisfying $\dim(|\Sigma(S|V)|)\geq 2$ and $|\Phi|\subset |\Sigma(S|V)|$, where $\Sigma(S|V)$ denotes the normal fan of $S$.

In this section we assume that $\Sigma(S|V)\hat{\cap}\mathcal{F}(\Delta)$ is $H$-simple for any $\Delta\in\Phi\Mx$ and $\Ht(H,\Phi,S)>0$.

Let $\gamma:(\Phi-(\Phi/H))_1\rightarrow\R$ denote the characteristic function of $(H, \Phi, S)$.
Denote
\begin{equation*}\begin{split}
m&=\sum_{\bar{E}\in(\Phi-(\Phi/H))_1}\lceil \gamma(\bar{E})\rceil\in\Z_+,\\
\bar{m}&=\sum_{\bar{E}\in(\Phi-(\Phi/H))_1}\lfloor \gamma(\bar{E})\rfloor\in\Z_0,\text{ and}\\
\end{split}\end{equation*}\begin{equation*}\begin{split}
\mathcal{R}&=\{ \bar{E}\in(\Phi-(\Phi/H))_1| \gamma(\bar{E})\not\in\Z\}\subset(\Phi-(\Phi/H))_1.
\end{split}\end{equation*}
@
We know $\Ht(H,S)\in(1/\Den(S/N))\Z_+$, $\Ht(H,\Phi,S)\in(1/\Den(S/N))\Z_+$ and $\Ht(H,\Phi,S)\leq\Ht(H,S)$.

We denote $\ell=\dim(\Vect(|\Phi|)^\vee|V^*)\in\Z_0$.

Consider any regular cone $\Theta$ over $N^*$ in $V^*$ satisfying $\Theta\subset|\Phi|$ and any $G\in\mathcal{F}(\Theta)_1$.
$\dim\Theta\leq\dim\Vect(|\Phi|)=\dim V-\ell$.
$S+(\Theta^\vee|V^*)$ is a rational conex pseudo polytope over $N$ in $V$.
$\Stab(S+(\Theta^\vee|V^*))=\Theta^\vee|V^*$.
$G\subset |\Sigma(S+(\Theta^\vee|V^*)|V)|=\Stab(S+(\Theta^\vee|V^*))^\vee|V=\Theta\subset|\Phi|\subset|\Sigma(S|V)|\subset\Vect(|\Sigma(S|V)|)$.
$\Sigma(S+(\Theta^\vee|V^*)|V)= \Sigma(S|V)\hat{\cap}\mathcal{F}(\Theta)$.
$\Stab(S+(\Theta^\vee|V^*))\cap(-\Stab(S+(\Theta^\vee|V^*)))=
\Vect(|\Sigma(S+(\Theta^\vee|V^*)|V)|)^\vee|V^*=\Vect(\Theta)^\vee|V^*\supset\Vect(|\Phi|)^\vee|V^*$.

If $\dim\Theta=\dim\Vect(|\Phi|)$, then
$\Vect(|\Sigma(S+(\Theta^\vee|V^*)|V)|)^\vee|V^*=
\Vect(|\Phi|)^\vee|V^*$ and
$\dim(\Vect(|\Sigma(S+(\Theta^\vee|V^*)|V)|)^\vee|V^*)=\ell$.

\begin{theorem}
\label{heightinequality}
Consider any compatible mapping
$$E:\{1,2,\ldots,m\}\rightarrow (\Phi-(\Phi/H))_1$$
with $S$.

Let
$$\Omega=\Omega(V^*, N^*, H, \Phi, m, E)$$
be the basic subdivision associated with the sextuplet $(V^*, N^*, H, \Phi, m, E)$.

Let $s=s(V^*, N^*, H, \Phi, m, E)$, 
$H=H(V^*, N^*, H, \Phi, m, E)$, \hfill\break and
$\Omega=\Omega(V^*, N^*, H, \Phi, m, E)$.
We have three mappings
$$s:\{0,1,\ldots,m\}\times(\Phi-(\Phi/H))_1\rightarrow \Z_0,$$
$$H:\{1,2,\ldots,m+1\}\rightarrow 2^{V^*},$$
$$\Omega:\{1,2,\ldots,m+1\}\rightarrow 2^{2^{V^*}}.$$

\begin{enumerate}
\item
$\bar{m}\in\Z_0$. $m\in\Z_+$. $\bar{m}\leq m$. $m-\bar{m}=\sharp\mathcal{R}$.
\item
For any $\bar{E}\in(\Phi-(\Phi/H))_1$, $s(\bar{m}, \bar{E})=\lfloor \gamma(\bar{E})\rfloor$ and $s(m, \bar{E})=\lceil \gamma(\bar{E})\rceil$.
\item
$\Omega$ is an iterated star subdivision of $\Phi$, and it is a flat regular fan over $N^*$ in $V^*$. $|\Omega|=|\Phi|\subset\Stab(S)^\vee|V=|\Sigma(S|V)|$. 
$\dim \Omega=\dim\Phi$. 
$\Vect(|\Omega|)=\Vect(|\Phi|)$.
\item
Consider any $i\in\{1,2,\ldots,m+1\}$.
\begin{enumerate}
\item
$\Omega(i)$ is a flat regular fan over $N^*$ in $V^*$.
$\Omega(i)\subset\Omega$. 
$|\Omega(i)|\subset|\Omega|$. 
$\dim\Omega(i)=\dim\Phi$.
$\Vect(|\Omega(i)|)= \Vect(|\Phi|)$.
$H(i)\in\Omega(i)_1$.
$\Omega(i)$ is starry with center in $H(i)$.
\item
For any $\Theta\in\Omega(i)$,  $\Theta^\vee|V^*$ is a rational polyhedral cone over $N$ in $V$, $\dim \Theta^\vee|V^*=\dim V$, $S+(\Theta^\vee|V^*)$ is a rational convex pseudo polytope over $N$ in $V$, $\Stab(S+(\Theta^\vee|V^*))=\Theta^\vee|V^*$, $\Sigma(S+(\Theta^\vee|V^*)|V)=\Sigma(S|V)\hat{\cap}\mathcal{F}(\Theta)$, and $\Sigma(S|V)\hat{\cap}\mathcal{F}(\Theta)$ is semisimple.
\item
\emph{[The height inequality]}
$\Ht(H(i), \Omega(i), S)<\Ht(H,\Phi,S)$.
\end{enumerate}
\item
Consider any $i\in\{1,2,\ldots,\bar{m}\}$.

For any $\Theta\in\Omega(i)$,
$\Sigma(S|V)\hat{\cap}\mathcal{F}(\Theta)= \mathcal{F}(\Theta)$ and
$S+(\Theta^\vee|V^*)=\{a\}+(\Theta^\vee|V^*)$ for some $a\in\mathcal{V}(\Phi, S)$.

$\Ht(H(i), \Omega(i), S)=0$.
\item
For any $\Theta\in\Omega(\bar{m}+1)\Mx$,
$\Sigma(S|V)\hat{\cap}\mathcal{F}(\Theta)$
is $H(\bar{m}+1)$-simple.
\item
Consider any $\Delta\in\Phi\Mx$. 
$\Ht(H,S+(\Delta^\vee|V^*))=\Ht(H,\Phi,S)$, if and only if, $\gamma(E)>0$ for some $E\in\mathcal{F}(H\Op|\Delta)_1$
\item
Consider any $\Delta\in\Phi\Mx$ such that $\gamma(\bar{E})=0$ for any $\bar{E}\in\mathcal{F}(H\Op|\Delta)_1$.

$\Ht(H,S+(\Delta^\vee|V^*))<\Ht(H,\Phi,S)$.

$\Delta\subset|\Omega(m+1)|$.
$\Omega\backslash\Delta=\Omega(m+1)\backslash\Delta=\mathcal{F}(\Delta)$.

For any $i\in\{1,2,\ldots,m\}$,
$|\Omega(i)|\cap\Delta\in \mathcal{F}(H\Op|\Delta)$, $|\Omega(i)|\cap\Delta\neq H\Op|\Delta$, and 
$\Omega(i)\backslash\Delta=\mathcal{F}(|\Omega(i)|\cap\Delta)$.
\end{enumerate}

Below we consider any $\Delta\in\Phi\Mx$ such that $\gamma(E)>0$ for some $E\in\mathcal{F}(H\Op|\Delta)_1$.
By $7$ we know $\Ht(H,S+(\Delta^\vee|V^*))=\Ht(H,\Phi,S)>0$.

Note that the normal fan $\Sigma(S+(\Delta^\vee|V^*)|V)=\Sigma(S|V)\hat{\cap}\mathcal{F}(\Delta)$ of $ S+(\Delta^\vee|V^*)$ is $H$-simple and the characteristic number $c(S+(\Delta^\vee|V^*))$ of $S+(\Delta^\vee|V^*)$ satisfies $c(S+(\Delta^\vee|V^*))\geq 2$.
Let $\bar{\Sigma}( S+(\Delta^\vee|V^*)|V)^1=\{\bar{\Lambda}\in\Sigma( S+(\Delta^\vee|V^*)|V)^1|\bar{\Lambda}^\circ\subset\Delta^\circ\}\cup\{H\Op|\Delta\}$ denote the $H$-skeleton of $\Sigma( S+(\Delta^\vee|V^*)|V)$.
$\sharp\Sigma( S+(\Delta^\vee|V^*)|V)^0=\sharp\bar{\Sigma}( S+(\Delta^\vee|V^*)|V)^1=c(S+(\Delta^\vee|V^*))\geq 2$.

We consider the $H$-order on $\Sigma(S+(\Delta^\vee|V^*)|V)^0$.
Let
$$\Lambda:\{1,2,\ldots, c(S+(\Delta^\vee|V^*))\}\rightarrow \Sigma( S+(\Delta^\vee|V^*)|V)^0,$$
denote the unique bijective mapping preserving the $H$-order.

We consider the $H$-order on $\bar{\Sigma}( S+(\Delta^\vee|V^*)|V)^1$.
Let
$$\bar{\Lambda}:\{1,2,\ldots, c(S+(\Delta^\vee|V^*))\}\rightarrow \bar{\Sigma}(S+(\Delta^\vee|V^*)|V)^1,$$
denote the unique bijective mapping preserving the $H$-order.

There exists uniquely a bijective mapping $A: \{1,2,\ldots, c(S+(\Delta^\vee|V^*))\}\rightarrow\mathcal{F}(S+(\Delta^\vee|V^*)|V)_\ell$ satisfying $\Lambda(i)=\Delta(A(i), S+(\Delta^\vee|V^*)|V)$ for any $i\in\{1,2,\ldots, c(S+(\Delta^\vee|V^*))\}$.
We take the bijective mapping $A: \{1,2,\ldots, c(S+(\Delta^\vee|V^*))\}\rightarrow\mathcal{F}(S+(\Delta^\vee|V^*)|V)_\ell$ satisfying these conditions.

For any $i\in\{1,2,\ldots, c(S+(\Delta^\vee|V^*))\}$ and any $\bar{E}\in\mathcal{F}(H\Op|\Delta)_1$, we can consider the structure constant
$c(\Sigma( S+(\Delta^\vee|V^*)|V),i, \bar{E})\in\Q_0$ of $\Sigma( S+(\Delta^\vee|V^*)|V)$ corresponding to the pair $(i, \bar{E}).$

Denote
\begin{equation*}\begin{split}
\hat{m}&=\sum_{ \bar{E}\in\mathcal{F}(H\Op|\Delta)_1}\lceil c(\Sigma( S+(\Delta^\vee|V^*)|V),2, \bar{E})\rceil\in\Z_+,\\
\hat{\bar{m}}&=\sum_{ \bar{E}\in\mathcal{F}(H\Op|\Delta)_1}\lfloor c(\Sigma( S+(\Delta^\vee|V^*)|V),2, \bar{E})\rfloor\in\Z_0,\text{ and}\\
\hat{\mathcal{R}}&=\{ \bar{E}\in\mathcal{F}(H\Op|\Delta)_1|
c(\Sigma( S+(\Delta^\vee|V^*)|V),2, \bar{E})\not\in\Z\}
\subset\mathcal{F}(H\Op|\Delta)_1.
\end{split}\end{equation*}

\begin{enumerate}
\setcounter{enumi}{9}
\item
$\hat{m}\in\Z_+$.
$\hat{\bar{m}}\in\Z_0$.
$\hat{\bar{m}}\leq \hat{m}$.
$\hat{m}=\sharp E^{-1}(\mathcal{F}(H\Op|\Delta)_1)$,
$\hat{\bar{m}}=\sharp(\{1,2,\ldots,\bar{m}\}\cap E^{-1}(\mathcal{F}(H\Op|\Delta)_1))$, and 
$\hat{\mathcal{R}}=\mathcal{R}\cap \mathcal{F}(H\Op|\Delta)_1$.
\end{enumerate}

Let $\hat{\tau}:\{1,2,\ldots, \hat{m}\}\rightarrow\{1,2,\ldots,m\}$ be the unique injective mapping preserving the order and satisfying 
$\hat{\tau}(\{1,2,\ldots, \hat{m}\})= E^{-1}(\mathcal{F}(H\Op|\Delta)_1)$.
Let $\hat{E}:\{1,2,\ldots, \hat{m}\}\rightarrow \mathcal{F}(H\Op|\Delta)_1$ be the unique mapping satisfying $\iota\hat{E}=E\hat{\tau}$ where $\iota: \mathcal{F}(H\Op|\Delta)_1)\rightarrow (\Phi-(\Phi/H))_1$ denotes the inclusion mapping.

Let
$\hat{\Omega}=\Omega(V^*, N^*, H, \mathcal{F}(\Delta), \hat{m}, \hat{E})\subset 2^{V^*}$, 
$\hat{G}=G(V^*, N^*, H, \mathcal{F}(\Delta), \hat{m}, \hat{E})$,\hfill\break 
$\hat{H}=H(V^*, N^*, H, \mathcal{F}(\Delta), \hat{m}, \hat{E})$, and $\hat{\Omega}=\Omega(V^*, N^*, H, \mathcal{F}(\Delta), \hat{m}, \hat{E})$.

$$\hat{G}:\{1,2,\ldots,\hat{m}\}\rightarrow 2^{V^*},$$
$$\hat{H}:\{1,2,\ldots,\hat{m}+1\}\rightarrow 2^{V^*},$$
$$\hat{\Omega}:\{1,2,\ldots,\hat{m}+1\}\rightarrow 2^{2^{V^*}}.$$

For any $i\in\{1,2,\ldots,\hat{m}+1\}$, we denote
\begin{equation*}\begin{split}
\Theta(i)&=|\hat{\Omega}(i)|\subset\Delta,\text{ and}\\
\bar{\Theta}(i)&= |\hat{\Omega}(i)-(\hat{\Omega}(i)/\hat{H}(i))|\subset\Theta(i).
\end{split}\end{equation*}

Denote
\begin{equation*}\begin{split}
\bar{c}(i)&=\max\{j\in\{2,3,\ldots,c(S+(\Delta^\vee|V^*))\}|\\
&\qquad\quad c(\Sigma( S+(\Delta^\vee|V^*)|V),j, \hat{E}(i))<
\lceil c(\Sigma( S+(\Delta^\vee|V^*)|V),2, \hat{E}(i))\rceil\}\\
&\qquad\quad\quad\in\{2,3,\ldots,c(S+(\Delta^\vee|V^*))\}
\end{split}\end{equation*}
for any $i\in\{\hat{\bar{m}}+1, \hat{\bar{m}}+2,\ldots,\hat{m}\}$.

\begin{enumerate}
\setcounter{enumi}{10}
\item
$\hat{E}$ is compatible with $S$.
\item
If $\hat{m}-\hat{\bar{m}}\geq 1$, then
$\bar{c}(i)\in\{2,3,\ldots,c(S+(\Delta^\vee|V^*))\}$ for any $i\in\{\hat{\bar{m}}+1, \hat{\bar{m}}+2,\ldots,\hat{m}\}$.
If $\hat{m}-\hat{\bar{m}}\geq 2$, then 
$\bar{c}(i)\leq\bar{c}(i+1)$  for any $i\in\{\hat{\bar{m}}+1, \hat{\bar{m}}+2,\ldots, \hat{m}-1\}$.
\item
$(H\Op|\Delta)\cap\bar{\Lambda}(2)\in\mathcal{F}(H\Op|\Delta)$.
$\mathcal{F}((H\Op|\Delta)\cap\bar{\Lambda}(2))_1=\{\bar{E}\in\mathcal{F}(H\Op|\Delta)_1| c(\Sigma(S+(\Delta^\vee|V^*)|V),2, \bar{E})=0\}$.
$\mathcal{F}((H\Op|\Delta)\cap\bar{\Lambda}(2))_1\cup \hat{E}(\{1,2,\ldots,\hat{m}\})=\mathcal{F}(H\Op|\Delta)_1$.
$\mathcal{F}((H\Op|\Delta)\cap\bar{\Lambda}(2))_1\cap \hat{E}(\{1,2,\ldots,\hat{m}\})=\emptyset$.
\item
For any $i\in\{1,2,\ldots, \hat{m}+1\}$,
$\hat{H}(i)\in\mathcal{F}(\Theta(i))_1$,
$\bar{\Theta}(i)=\hat{H}(i)\Op|\Theta(i) \in\mathcal{F}(\Theta(i))^1$,
$\dim\Theta(i)=\dim\Vect(|\Phi|)$, and\hfill\break
$\dim\Vect(|\Sigma(S+(\Theta(i)^\vee|V^*)|V)|)^\vee|V^*=\ell$.

\item
Consider any $i\in\{1,2,\ldots, \hat{\bar{m}}\}$.
$\Sigma(S|V)\hat{\cap}\mathcal{F}(\Theta(i))= \mathcal{F}(\Theta(i))$.
$S+(\Theta(i)^\vee$\hfill\break$|V^*)=A(1)+(\Theta(i)^\vee|V^*)$.
$\Ht(\hat{H}(i),S+(\Theta(i)^\vee|V^*))=0<\Ht(H,S+(\Delta^\vee|V^*))$.
For any $j\in\{2,3,\ldots,c(S+(\Delta^\vee|V^*))\}$, $\bar{\Lambda}(j)\cap(\Theta(i)^\circ\cup\bar{\Theta}(i)^\circ)=\emptyset$.
\item
$\Sigma(S|V)\hat{\cap}\mathcal{F}(\Theta(\hat{\bar{m}}+1))$
is $\hat{H}(\hat{\bar{m}}+1)$-simple.
\item
Assume $\hat{m}\neq\hat{\bar{m}}$.
Consider any $i\in\{\hat{\bar{m}}+1, \hat{\bar{m}}+2, \ldots, \hat{m}\}$.

$\Sigma(S|V)\hat{\cap}\mathcal{F}(\Theta(i))$ is semisimple.

$\mathcal{F}(S+(\Theta(i)^\vee|V^*))_\ell=
\{A(j)|j\in\{1,2,\ldots, \bar{c}(i)\}\}$.

Take any point $a(j)\in A(j)$ for any $j\in\{1,2,\ldots, \bar{c}(i)\}$.
For any $j\in\{2,3,\ldots, \bar{c}(i)\}$,
$0<\langle b_{\hat{H}(i)/N^*},a(j-1)-a(j)\rangle<\langle b_{H/N^*},a(j-1)-a(j)\rangle$.

$0<\Ht(\hat{H}(i), S+(\Theta(i)^\vee|V^*))<\Ht(H, S+(\Delta^\vee|V^*))$.

For any $j\in\{2,3,\ldots,c(S+(\Delta^\vee|V^*))\}$,
$j\in\{2,3,\ldots, \bar{c}(i)\}
\Leftrightarrow
c(\Sigma(S+(\Delta^\vee|V^*)|V), j,\hat{E}(i))<
\lceil c(\Sigma(S+(\Delta^\vee|V^*)|V),2,\hat{E}(i))\rceil
\Leftrightarrow
\lfloor c(\Sigma(S+(\Delta^\vee|V^*)$\hfill\break$|V),2,\hat{E}(i))\rfloor<
c(\Sigma(S+(\Delta^\vee|V^*)|V),j,\hat{E}(i))<
\lceil c(\Sigma(S+(\Delta^\vee|V^*)|V),2,\hat{E}(i))\rceil
$\break$\Leftrightarrow
\bar{\Lambda}(j)\cap(\hat{G}(i)+\hat{H}(i))^\circ\neq\emptyset
\Leftrightarrow
\bar{\Lambda}(j)\cap\Theta(i)^\circ\neq\emptyset
\Leftrightarrow
\bar{\Lambda}(j)\cap(\Theta(i)^\circ\cup\bar{\Theta}(i)^\circ)\neq\emptyset
$.

\item
$\Sigma(S|V)\hat{\cap}\mathcal{F}(\Theta(\hat{m}+1))$ is semisimple.

$\mathcal{F}(S+(\Theta(\hat{m}+1)^\vee|V^*))_\ell=
\{A(j)|j\in\{k,k+1,\ldots, \hat{m}+1\}\}$ for some $k\in\{2,3,\ldots, \hat{m}+1\}$.

If $\hat{m}=\hat{\bar{m}}$, then
$\mathcal{F}(S+(\Theta(\hat{m}+1)^\vee|V^*))_\ell=
\{A(j)|j\in\{2,3,\ldots, \hat{m}+1\}\}$.

$\Ht(\hat{H}(\hat{m}+1), S+(\Theta(\hat{m}+1)^\vee|V^*))<\Ht(H, S+(\Delta^\vee|V^*))$.
\item The following three conditions are equivalent:
\begin{enumerate}
\item
$\hat{\Omega}$ is a subdivision of $\Sigma(S+(\Delta^\vee|V^*)|V)$.
\item
$c(S+(\Delta^\vee|V^*)|V)=2$ and $\hat{m}=\bar{\hat{m}}$.
\item
$c(S+(\Delta^\vee|V^*)|V)=2$ and $c(\Sigma(S+(\Delta^\vee|V^*)|V),2, \bar{E})\in\Z$ for any $\bar{E}\in\mathcal{F}(H\Op|\Delta)_1$.
\end{enumerate}
\end{enumerate}

Below we assume that $\hat{\Omega}$ is a subdivision of $\Sigma(S+(\Delta^\vee|V^*)|V)$ until claim $26$.

\begin{enumerate}
\setcounter{enumi}{19}

\item
$\hat{m}=\hat{\bar{m}}\geq 1$. $c(S+(\Delta^\vee|V^*))=2$. $\bar{\Lambda}(1)=\bar{\Theta}(1)=H\Op|\Delta$.
$\Lambda(1)\cap\Lambda(2)=\bar{\Lambda}(2)=\bar{\Theta}(\hat{m}+1)$.
$\Lambda(1)=\cup_{i\in\{1,2,\ldots,\hat{m}\}}\Theta(i)$.
$\Lambda(2)=\Theta(\hat{m}+1)$.
\item
$\{\Theta\in\hat{\Omega}|\Theta^\circ\subset\Delta^\circ\}$\hfill\break
$=\{\Theta(i)|i\in\{1,2,\ldots,\hat{m}+1\}\}\cup\{\bar{\Theta}(i)|i\in\{2,3,\ldots,\hat{m}\}\}\cup\{\bar{\Theta}(\hat{m}+1)\}$.

\noindent $\{\Theta\in\hat{\Omega}|\Theta^\circ\subset\Delta^\circ,$ 
The unique element $\Lambda\in \Sigma(S+(\Delta^\vee|V^*)|V)$ satisfying $\Theta^\circ\subset\Lambda^\circ$ satisfies $\dim\Lambda=\dim\Delta-1\}
=\{\bar{\Theta}(\hat{m}+1)\}$.
\item
$H=\hat{H}(\hat{m}+1)\in\mathcal{F}(\Theta(\hat{m}+1))_1$.
$H\in\mathcal{F}(\Delta)_1$.
$\bar{\Theta}(\hat{m}+1)\in\mathcal{F}(\Theta(\hat{m}+1))^1$.
$\bar{\Theta}(\hat{m}+1)+ H= \Theta(\hat{m}+1) \in\smash{\hat{\Omega}}\Mx$.
$\bar{\Theta}(\hat{m}+1)\cap H=\{0\}$.

$\bar{\Theta}(\hat{m}+1)\subset\Lambda(1)$.
$\bar{\Theta}(\hat{m}+1)\subset\Lambda(2)$.
$\Theta(\hat{m}+1)\not\subset\Lambda(1)$.
$\Theta(\hat{m}+1)\subset\Lambda(2)$.
\end{enumerate}

Let $a(1)\in A(1)$ and $a(2)\in A(2)$ be any points.
Let $F=A(2)+(\Vect(\bar{\Theta}(\hat{m}+1))^\vee|V^*)$. 
\begin{enumerate}
\setcounter{enumi}{22}
\item
$\bar{\Theta}(\hat{m}+1)\subset\Theta(\hat{m}+1)\subset\Delta\subset|\Sigma(S|V)|$, and $S+(\bar{\Theta}(\hat{m}+1)^\vee|V^*)\supset S+(\Theta(\hat{m}+1)^\vee|V^*)\supset S+(\Delta^\vee|V^*)\supset S$.

$S+(\bar{\Theta}(\hat{m}+1)^\vee|V^*)=A(2)+(\bar{\Theta}(\hat{m}+1)^\vee|V^*)=\{a(2)\}+(\bar{\Theta}(\hat{m}+1)^\vee|V^*)$.
$S+(\Theta(\hat{m}+1)^\vee|V^*)=A(2)+ (\Theta(\hat{m}+1)^\vee|V^*)=\{a(2)\}+ (\Theta(\hat{m}+1)^\vee|V^*)$.
$S+(\Delta^\vee|V^*)=\Conv(A(1)\cup A(2))+(\Delta^\vee|V^*)=\Conv(\{a(1),a(2)\})+(\Delta^\vee|V^*)$.
\item
$F\in\mathcal{F}(S+(\bar{\Theta}(\hat{m}+1)^\vee|V^*))$.
$F\cap(S+(\Theta(\hat{m}+1)^\vee|V^*))\in\mathcal{F}(S+(\Theta(\hat{m}+1)^\vee|V^*))$.
$F\cap(S+(\Delta^\vee|V^*))\in\mathcal{F}(S+(\Delta^\vee|V^*))$.

$\Delta(F, S+(\bar{\Theta}(\hat{m}+1)^\vee|V^*)|V)=
\Delta(F\cap(S+(\Theta(\hat{m}+1)^\vee|V^*)), S+(\Theta(\hat{m}+1)^\vee|V^*)|V)=
\Delta(F\cap(S+(\Delta^\vee|V^*)), S+(\Delta^\vee|V^*)|V)=
\bar{\Lambda}(2)$.

$F=\{a(2)\}+ (\Vect(\bar{\Theta}(\hat{m}+1))^\vee|V^*)=\Affi(\{a(1),a(2)\})+(\Vect(\Delta)^\vee|V^*)$.
$F\cap(S+(\Theta(\hat{m}+1)^\vee|V^*))=A(2)+\Delta(\bar{\Theta}(\hat{m}+1), \Theta(\hat{m}+1)|V^*)=\{a(2)\}+ \R_0(a(1)-a(2))+(\Vect(\Delta)^\vee|V^*)$.
$F\cap(S+(\Delta^\vee|V^*))=\Conv(A(1)\cup A(2))= \Conv(\{a(1), a(2)\})+ (\Vect(\Delta)^\vee|V^*)$.
\item
$\langle b_{H/N^*},a(1)\rangle> \langle b_{H/N^*},a(2)\rangle$.

\noindent$\{\langle b_{H/N^*}, a\rangle|a\in F\}=\R$.

\noindent$\{\langle b_{H/N^*}, a\rangle|a\in F\cap(S+(\Theta(\hat{m}+1)^\vee|V^*))\}=
\{t\in\R|\langle b_{H/N^*},a(2)\rangle\leq t \}$.

\noindent$\{\langle b_{H/N^*}, a\rangle|a\in F\cap (S+(\Delta^\vee|V^*))\}=
\{t\in\R|\langle b_{H/N^*},a(2)\rangle\leq t\leq \langle b_{H/N^*},$
\hfill\break$a(1)\rangle\}$.
\item
\begin{equation*}\begin{split}
&\max\{\langle b_{H/N^*}, a\rangle|a\in F\cap (S+(\Delta^\vee|V^*))\}\\
&\qquad -\min\{\langle b_{H/N^*}, a\rangle|a\in F\cap (S+(\Delta^\vee|V^*))\}\\
=\:&\langle b_{H/N^*},a(1)\rangle- \langle b_{H/N^*},a(2)\rangle\\
=\:& \Ht(H, S+(\Delta^\vee|V^*))\\
\end{split}\end{equation*}
\end{enumerate}

Below we consider any rational convex pseudo polytopes $T$ and $U$ over $N$ in $V$ satisfying $T+U=S$ until claim $32$. 
\begin{enumerate}\setcounter{enumi}{26}\item
$\Sigma(T|V)\hat{\cap}\Sigma(U|V)= \Sigma(S|V)$.
$\Sigma(T+(\Delta^\vee|V^*)|V)= \Sigma(T|V) \hat{\cap}\mathcal{F}(\Delta)$ is $H$-simple.
$\Sigma(U+(\Delta^\vee|V^*)|V)= \Sigma(U|V) \hat{\cap}\mathcal{F}(\Delta)$ is $H$-simple.
$\Ht(H, T+(\Delta^\vee|V^*))+\Ht(H, U+(\Delta^\vee|V^*))=\Ht(H, S+(\Delta^\vee|V^*))$.
\end{enumerate}

Let $\bar{\Sigma}(T+(\Delta^\vee|V^*)|V)^1$ and $\bar{\Sigma}(U+(\Delta^\vee|V^*)|V)^1$ denote the $H$-skeleton of $\Sigma(T+(\Delta^\vee|V^*)|V)$ and $\Sigma(U+(\Delta^\vee|V^*)|V)$ respectively.
\begin{enumerate}\setcounter{enumi}{27}\item
$\bar{\Sigma}(T+(\Delta^\vee|V^*)|V)^1\cup\bar{\Sigma}(U+(\Delta^\vee|V^*)|V)^1=\bar{\Sigma}(S+(\Delta^\vee|V^*)|V)^1$.
\item
For any $i\in\{1,2,\ldots,\hat{\bar{m}}\}$, $\Ht(\hat{H}(i), T+(\Theta(i)^\vee|V^*))=0$.
\item
Assume $\hat{m}\neq\hat{\bar{m}}$. Consider any $i\in\{\hat{\bar{m}}+1, \hat{\bar{m}}+2,\ldots,\hat{m}\}$. If $0<\Ht(\hat{H}(i),$\break$T+(\Theta(i)^\vee|V^*))$, then $\Ht(\hat{H}(i), T+(\Theta(i)^\vee|V^*))<\Ht(H, T+(\Delta^\vee|V^*))$.
\item
$\Ht(\hat{H}(\hat{m}+1), T+(\Theta(\hat{m}+1)^\vee|V^*))\leq\Ht(H, T+(\Delta^\vee|V^*))$.
\item
Let $r_U=c(U+(\Delta^\vee|V^*))\in\Z_+$. Assume that if $r_U\geq 2$, then the structure constant of $\Sigma(U+(\Delta^\vee|V^*)|V)$ corresponding to the pair $(2,\bar{E})$ is an integer for any $\bar{E}\in\mathcal{F}(H\Op|\Delta)_1$. Then, $\Ht(\hat{H}(\hat{m}+1), T+(\Theta(\hat{m}+1)^\vee|V^*))=\Ht(H, T+(\Delta^\vee|V^*))$, if and only if, $\bar{\Lambda}(2)\not\in\bar{\Sigma}(T+(\Delta^\vee|V^*)|V)^1$.
\end{enumerate}
\end{theorem}

\begin{proof}
We show only claim $17$. This is the most important. Claim $4$.(c) follows from claims $7, 8, 15, 17$ and $18$.

Assume $\hat{m}\neq\hat{\bar{m}}$. $\hat{m}-\hat{\bar{m}}\geq 1$.
Consider any $i\in\{\hat{\bar{m}}+1, \hat{\bar{m}}+2,\ldots, \hat{m}\}$.
Since $\Theta(i)\subset\Delta$ and $\Sigma(S|V)\hat{\cap}\mathcal{F}(\Delta)$ is $H$-simple, $\Sigma(S|V)\hat{\cap}\mathcal{F}(\Delta)$ is semisimple and we know $\Sigma(S|V)\hat{\cap}\mathcal{F}(\Theta(i)) =\Sigma(S|V)\hat{\cap}\mathcal{F}(\Delta) \hat{\cap}\mathcal{F}(\Theta(i))$ is semisimple by Lemma~\ref{propsimple}.10.

Let $\hat{s}=s(V^*, N^*,H,\mathcal{F}(\Delta),\hat{m},\hat{E})$.
$\hat{s}:\{0,1,\ldots,\hat{m}\}\times\mathcal{F}(H\Op|\Delta)_1\rightarrow \Z_0$.

Note that $\Sigma(S+(\Theta(i)^\vee|V^*)|V)^0=(\Sigma(S|V)\hat{\cap}\mathcal{F}(\Theta(i)))^0=\{\Lambda(j)\cap\Theta(i)|j\in\{1,2,\ldots,c(S+(\Delta^\vee|V^*))\}, \Lambda(j)^\circ\cap\Theta(i)^\circ\neq\emptyset\}$,

\begin{equation*}\begin{split}
\partial^H_-\Theta(i)&=\bar{\Theta}(i) =\sum_{\bar{E}\in\mathcal{F}(H\Op|\Delta)_1}\R_0(b_{\bar{E}/N^*}+\hat{s}(i-1, \bar{E})b_{H/N^*})\\
&=\sum_{\bar{E}\in\mathcal{F}(H\Op|\Delta)_1-\{\hat{E}(i)\}}\R_0(b_{\bar{E}/N^*}
 +\hat{s}(i-1, \bar{E})b_{H/N^*})\\
&\qquad\qquad\quad +\R_0(b_{\hat{E}(i)/N^*}+\hat{s}(i-1, \hat{E}(i))b_{H/N^*}),\\
\partial^H_+\Theta(i)&= \bar{\Theta}(i+1)=\sum_{\bar{E}\in\mathcal{F}(H\Op|\Delta)_1}\R_0(b_{\bar{E}/N^*}+\hat{s}(i, \bar{E})b_{H/N^*})\\
&=\sum_{\bar{E}\in\mathcal{F}(H\Op|\Delta)_1-\{\hat{E}(i)\}}\R_0(b_{\bar{E}/N^*}
 +\hat{s}(i-1,\bar{E})b_{H/N^*})\\
&\qquad\qquad\quad +\R_0(b_{\hat{E}(i)/N^*}+\hat{s}(i, \hat{E}(i))b_{H/N^*}),\text{ and}
\end{split}\end{equation*}
\begin{equation*}\begin{split}
\hat{s}(i-1, \hat{E}(i))&=\lfloor c(\Sigma(S+(\Delta^\vee|V^*)|V^*),2,\hat{E}(i))\rfloor\\
&< c(\Sigma(S+(\Delta^\vee|V^*)|V^*),2,\hat{E}(i))\leq
c(\Sigma(S+(\Delta^\vee|V^*)|V^*),\bar{c}(i),\hat{E}(i))\\
&< \lceil c(\Sigma(S+(\Delta^\vee|V^*)|V^*),2,\hat{E}(i))\rceil
= \hat{s}(i, \hat{E}(i)).
\end{split}\end{equation*} 

Consider any $j\in\{1,2,\cdots, \bar{c}(i)\}$ with $j\neq c(S+(\Delta^\vee|V^*))$. There exists a real number $r_j$ satisfying $ \hat{s}(i-1, \hat{E}(i))<r_j< \hat{s}(i, \hat{E}(i))$ and $ c(\Sigma(S+(\Delta^\vee|V^*)|V^*),j,$\hfill\break$\hat{E}(i))\leq r_j\leq c(\Sigma(S+(\Delta^\vee|V^*)|V^*),j+1,\hat{E}(i))$, and there exist a mapping
$t_j:\mathcal{F}(H\Op|\Delta)_1\rightarrow\R_+$ and a real number $u_j$ satisfying
$$\sum_{\bar{E}\in\mathcal{F}(H\Op|\Delta)_1}t_j(\bar{E}) \hat{s}(i-1,\bar{E})<u_j
<\sum_{\bar{E}\in\mathcal{F}(H\Op|\Delta)_1}t_j(\bar{E}) \hat{s}(i,\bar{E})\text { and}$$
\begin{equation*}\begin{split}
&\sum_{\bar{E}\in\mathcal{F}(H\Op|\Delta)_1}t_j(\bar{E}) c(\Sigma(S+(\Delta^\vee|V^*)|V^*),j, \bar{E})<u_j\\
&\qquad<\sum_{\bar{E}\in\mathcal{F}(H\Op|\Delta)_1}t_j(\bar{E}) c(\Sigma(S+(\Delta^\vee|V^*)|V^*),j+1, \bar{E}).
\end{split}\end{equation*} 

We take a mapping
$t_j:\mathcal{F}(H\Op|\Delta)_1\rightarrow\R_+$ and a real number $u_j$ satisfying
the above conditions.
We know
$(\sum_{\bar{E}\in\mathcal{F}(H\Op|\Delta)_1}t_j(\bar{E})b_{\bar{E}/N^*})+u_jb_{H/N^*}\in\Lambda(j)^\circ\cap\Theta(i)^\circ \neq\emptyset$.

We know that if $ \bar{c}(i)\neq c(S+(\Delta^\vee|V^*))$, then for any $j\in\{1,2,\cdots, \bar{c}(i)\}$, $\Lambda(j)^\circ\cap\Theta(i)^\circ \neq\emptyset$.

Consider the case $ \bar{c}(i)= c(S+(\Delta^\vee|V^*))$.
We denote $j=\bar{c}(i)= c(S+(\Delta^\vee|V^*))$.
There exists a real number $r_j$ satisfying $ \hat{s}(i-1, \hat{E}(i))<r_j< \hat{s}(i, \hat{E}(i))$ and $ c(\Sigma(S+(\Delta^\vee|V^*)|V^*),j,\hat{E}(i))< r_j$, and there exist a mapping
$t_j:\mathcal{F}(H\Op|\Delta)_1\rightarrow\R_+$ and a real number $u_j$ satisfying
$$\sum_{\bar{E}\in\mathcal{F}(H\Op|\Delta)_1}t_j(\bar{E}) \hat{s}(i-1, \bar{E})<u_j
<\sum_{\bar{E}\in\mathcal{F}(H\Op|\Delta)_1}t_j(\bar{E}) \hat{s}(i,\bar{E})\text { and}$$
$$\sum_{\bar{E}\in\mathcal{F}(H\Op|\Delta)_1}t_j(\bar{E}) c(\Sigma(S+(\Delta^\vee|V^*)|V^*),j, \bar{E})<u_j.$$
We take a mapping
$t_j:\mathcal{F}(H\Op|\Delta)_1\rightarrow\R_+$ and a real number $u_j$ satisfying
the above conditions.
We know
$(\sum_{\bar{E}\in\mathcal{F}(H\Op|\Delta)_1}t_j(\bar{E})b_{\bar{E}/N^*})+u_jb_{H/N^*}\in\Lambda(j)^\circ\cap\Theta(i)^\circ \neq\emptyset$.

We know that if $ \bar{c}(i)= c(S+(\Delta^\vee|V^*))$, then for any $j\in\{1,2,\cdots, \bar{c}(i)\}$, $\Lambda(j)^\circ\cap\Theta(i)^\circ \neq\emptyset$.

We know that for any $j\in\{1,2,\cdots, \bar{c}(i)\}$, $\Lambda(j)^\circ\cap\Theta(i)^\circ \neq\emptyset$.

Consider any $j\in\{1,2,\cdots,c(S+(\Delta^\vee|V^*))\}$ satisfying $\Lambda(j)^\circ\cap\Theta(i)^\circ \neq\emptyset$.
We know that there exists $\bar{E}\in\mathcal{F}(H\Op|\Delta)_1$ satisfying
$ c(\Sigma(S+(\Delta^\vee|V^*)|V^*),j,\bar{E})<\hat{s}(i, \bar{E})$.
We take $\bar{E}\in\mathcal{F}(H\Op|\Delta)_1$ satisfying
$ c(\Sigma(S+(\Delta^\vee|V^*)|V^*),j,\bar{E})<\hat{s}(i,\bar{E})$.

If $\bar{E}\not\in\hat{\mathcal{R}}$, then $ c(\Sigma(S+(\Delta^\vee|V^*)|V^*),j,\bar{E})<\hat{s}(i,\bar{E})= c(\Sigma(S+(\Delta^\vee|V^*)|V^*),2,\bar{E})$ and $j=1\leq\bar{c}(i)$.

If $\bar{E}\in\hat{\mathcal{R}}$, then there exists $k\in\{\hat{\bar{m}}+1, \hat{\bar{m}}+2,\ldots,\hat{m}\}$ satisfying $\bar{E}=\hat{E}(k)$. We take $k\in\{\hat{\bar{m}}+1, \hat{\bar{m}}+2,\ldots,\hat{m}\}$ satisfying $\bar{E}=\hat{E}(k)$.

Consider the case $k\leq i$. 
We have
$c(\Sigma(S+(\Delta^\vee|V^*)|V^*),j, \hat{E}(k))= c(\Sigma(S+(\Delta^\vee|V^*)|V^*),j,\bar{E})<\hat{s}(i, \bar{E})= \hat{s}(i, \hat{E}(k))=\lceil c(\Sigma(S+(\Delta^\vee|V^*)|V^*),2, \hat{E}(k))\rceil $ and
$1\leq j\leq\bar{c}(k)\leq\bar{c}(i)$.

Consider the case $k>i$. We have
$c(\Sigma(S+(\Delta^\vee|V^*)|V^*),j, \hat{E}(k))= c(\Sigma(S+(\Delta^\vee|V^*)|V^*),j,\bar{E})<\hat{s}(i, \bar{E})= \hat{s}(i, \hat{E}(k))=\lfloor c(\Sigma(S+(\Delta^\vee|V^*)|V^*),2,\hat{E}(k))\rfloor<
c(\Sigma(S+(\Delta^\vee|V^*)|V^*),2, \hat{E}(k))$ and
$j=1\leq\bar{c}(i)$.

We know that for any $j\in\{1,2,\cdots,c(S+(\Delta^\vee|V^*))\}$ satisfying $\Lambda(j)^\circ\cap\Theta(i)^\circ \neq\emptyset$, $j\in\{1,2,\cdots, \bar{c}(i)\}$.

Therefore,
$\Sigma(S+(\Theta(i)^\vee|V^*)|V)^0=\{\Lambda(j)\cap\Theta(i)|j\in\{1,2,\ldots,c(S+(\Delta^\vee|V^*))\},$\hfill\break$ \Lambda(j)^\circ\cap\Theta(i)^\circ\neq\emptyset\}
=\{\Lambda(j)\cap\Theta(i)|j\in\{1,2,\ldots,\bar{c}(i)\}\}$.
Since $\dim\Theta(i)=\dim\Delta=\dim\Phi$, we know
$\mathcal{F}( S+(\Theta(i)^\vee|V^*))_\ell=\{A(j) |j\in\{1,2,\ldots,\bar{c}(i)\}\}$.

For any $j\in\{1,2,\ldots,c(S+(\Delta^\vee|V^*))\}$ we take any point $a(j)\in A(j)$.

Consider any $ j\in\{2,3,\ldots,\bar{c}(i)\}$.

$\langle b_{H/N^*}, a(j-1)-a(j) \rangle>0$ and 
$ c(S+(\Delta^\vee|V^*), 2, \hat{E}(i)))\leq
c(S+(\Delta^\vee|V^*), j,$\hfill\break$ \hat{E}(i)))\leq
c(S+(\Delta^\vee|V^*), \bar{c}(i), \hat{E}(i)))<
\lceil c(S+(\Delta^\vee|V^*), 2, \hat{E}(i))\rceil$.
$0<\lceil c(S+(\Delta^\vee|V^*), 2, \hat{E}(i))\rceil- c(S+(\Delta^\vee|V^*), 2, \hat{E}(i))<1$.

By Lemma~\ref{propsimple}.$18$.(f) we have
\begin{equation*}\begin{split}
&\langle b_{\hat{H}(i)/N^*},a(j-1) \rangle-\langle b_{\hat{H}(i)/N^*}, a(j)\rangle
=\langle b_{\hat{H}(i)/N^*},a(j-1)-a(j)\rangle\\
=\:&\langle b_{\hat{E}(i)/N^*}+\hat{s}(i, \hat{E}(i))b_{H/N^*},a(j-1)-a(j)\rangle\\
=\:&\langle b_{\hat{E}(i)/N^*}, a(j-1) -a(j)\rangle
+\hat{s}(i, \hat{E}(i), i)\langle b_{H/N^*}, a(j-1)-a(j) \rangle\\
=\:&-c(S+(\Delta^\vee|V^*), j, \hat{E}(i))\langle b_{H/N^*}, a(j-1)-a(j) \rangle\\
\:&\qquad\qquad+\lceil c(S+(\Delta^\vee|V^*), 2, \hat{E}(i))\rceil\langle b_{H/N^*}, a(j-1)-a(j) \rangle\\
=\:&(\lceil c(S+(\Delta^\vee|V^*), 2, \hat{E}(i))\rceil- c(S+(\Delta^\vee|V^*), j, \hat{E}(i))) \langle b_{H/N^*}, a(j-1)-a(j) \rangle\\
>\:&0\text{, and}
\end{split}\end{equation*}
\begin{equation*}\begin{split}
&\langle b_{\hat{H}(i)/N^*},a(j-1)-a(j)\rangle\\
=\:&(\lceil c(S+(\Delta^\vee|V^*), 2, \hat{E}(i))\rceil- c(S+(\Delta^\vee|V^*), j, \hat{E}(i))) \langle b_{H/N^*}, a(j-1)-a(j) \rangle\\
\leq\:&(\lceil c(S+(\Delta^\vee|V^*), 2, \hat{E}(i))\rceil- c(S+(\Delta^\vee|V^*), 2, \hat{E}(i))) \langle b_{H/N^*}, a(j-1)-a(j) \rangle\\
<\:& \langle b_{H/N^*}, a(j-1)-a(j) \rangle.
\end{split}\end{equation*}

We know
\begin{equation*}\begin{split}
&\max\{\langle b_{\hat{H}(i)/N^*},a(j)\rangle|j\in\{1,2,\ldots,\bar{c}(i)\}\}
=\langle b_{\hat{H}(i)/N^*},a(1)\rangle,\\
&\min\{\langle b_{\hat{H}(i)/N^*},a(j)\rangle|j\in\{1,2,\ldots,\bar{c}(i)\}\}
=\langle b_{\hat{H}(i)/N^*},a(\bar{c}(i))\rangle,\\
&\Ht(\hat{H}(i), S+(\Theta(i)^\vee|V^*))
=\langle b_{\hat{H}(i)/N^*},a(1)\rangle-\langle b_{\hat{H}(i)/N^*}, a(\bar{c}(i))\rangle\\
=\:&\langle b_{\hat{H}(i)/N^*},a(1)- a(\bar{c}(i))\rangle
=\sum_{j=2}^{\bar{c}(i)} \langle b_{\hat{H}(i)/N^*},a(j-1)- a(j)\rangle\\
>\:&0\text{, and}\\
&\Ht(\hat{H}(i), S+(\Theta(i)^\vee|V^*))
=\sum_{j=2}^{\bar{c}(i)} \langle b_{\hat{H}(i)/N^*},a(j-1)- a(j)\rangle\\
<\:&\sum_{j=2}^{\bar{c}(i)} \langle b_{H/N^*},a(j-1)- a(j)\rangle
\leq\sum_{j=2}^{c(S+(\Delta^\vee|V^*))} \langle b_{H/N^*},a(j-1)- a(j)\rangle\\
=\:&\langle b_{H/N^*},a(1)- a(c(S+(\Delta^\vee|V^*)))\rangle
=\langle b_{H/N^*},a(1) \rangle-\langle b_{H/N^*}, a(c(S+(\Delta^\vee|V^*)))\rangle\\
=\:&\Ht(H,S+(\Delta^\vee|V^*)).
\end{split}\end{equation*}

Consider any $j\in\{2,3,\ldots,c(S+(\Delta^\vee|V^*))\}$.

By definition of $\bar{c}(i)$, we know
$j\in\{2,3,\ldots, \bar{c}(i)\}
\Leftrightarrow
c(\Sigma(S+(\Delta^\vee|V^*)|V), j,$\hfill\break$\hat{E}(i))<
\lceil c(\Sigma(S+(\Delta^\vee|V^*)|V),2,\hat{E}(i))\rceil$.
Since $\hat{E}(i)\in\hat{\mathcal{R}}$, 
$\lfloor c(\Sigma(S+(\Delta^\vee|V^*)|V),2,$\hfill\break$\hat{E}(i))\rfloor<
c(\Sigma(S+(\Delta^\vee|V^*)|V),2,\hat{E}(i))\leq
c(\Sigma(S+(\Delta^\vee|V^*)|V),j,\hat{E}(i))$.
Therefore,
$c(\Sigma(S+(\Delta^\vee|V^*)|V), j,\hat{E}(i))<
\lceil c(\Sigma(S+(\Delta^\vee|V^*)|V),2,\hat{E}(i))\rceil
\Leftrightarrow
\lfloor c(\Sigma(S+(\Delta^\vee|V^*)|V),2,\hat{E}(i))\rfloor<
c(\Sigma(S+(\Delta^\vee|V^*)|V),j,\hat{E}(i))<
\lceil c(\Sigma(S+(\Delta^\vee|V^*)|V),2, $\break$\hat{E}(i))\rceil$.
Since 
$\hat{G}(i)+\hat{H}(i)=
\R_0(b_{\hat{E}(i)/N^*}+\lfloor c(\Sigma(S+(\Delta^\vee|V^*)|V),2,\hat{E}(i))\rfloor b_{H/N^*})+
\R_0(b_{\hat{E}(i)/N^*}+\lceil c(\Sigma(S+(\Delta^\vee|V^*)|V),2,\hat{E}(i))\rceil b_{H/N^*})$,
we know
$\lfloor c(\Sigma(S+(\Delta^\vee|V^*)|V),$\hfill\break$2,\hat{E}(i))\rfloor<
c(\Sigma(S+(\Delta^\vee|V^*)|V),j,\hat{E}(i))<
\lceil c(\Sigma(S+(\Delta^\vee|V^*)|V),2,\hat{E}(i))\rceil
\Leftrightarrow
\bar{\Lambda}(j)\cap(\hat{G}(i)+\hat{H}(i))^\circ\neq\emptyset$.
Note that $\{k\in\{1,2,\ldots,c(S+(\Delta^\vee|V^*))\}|\Lambda(k)\cap\Theta(i)^\circ\neq\emptyset\}=\{1,2,\ldots,\bar{c}(i)\}$, since $\mathcal{F}(S+(\Theta(i)^\vee|V^*))_\ell=
\{A(j)|j\in\{1,2,\ldots, \bar{c}(i)\}\}$.
Therefore, $j\in\{2,3,\ldots, \bar{c}(i)\}
\Leftrightarrow
\Lambda(j-1)\cap \Theta(i)^\circ\neq\emptyset$ and $\Lambda(j)\cap \Theta(i)^\circ\neq\emptyset$.
Since $\bar{\Lambda}(j)= \Lambda(j-1)\cap\Lambda(j)$, we know 
$\Lambda(j-1)\cap \Theta(i)^\circ\neq\emptyset$ and $\Lambda(j)\cap \Theta(i)^\circ\neq\emptyset
\Leftrightarrow
\bar{\Lambda}(j)\cap \Theta(i)^\circ\neq\emptyset$.

Assume $\bar{\Lambda}(j)\cap\bar{\Theta}(i)^\circ\neq\emptyset$.
Take any point $\omega\in\bar{\Lambda}(j)\cap\bar{\Theta}(i)^\circ$.
Let $\chi= b_{\hat{E}(i)/N^*}+ c(\Sigma(S+(\Delta^\vee|V^*)|V),j,\hat{E}(i))b_{H/N^*}\in\bar{\Lambda}(j)$.
We know that there exists a real number $\epsilon$ with $0<\epsilon\leq 1$ such that
$(1-t)\omega+t\chi\in\bar{\Lambda}(j)\cap\Theta(i)^\circ$ for any real number $t$ with $0<t<\epsilon$, since 
$\lfloor c(\Sigma(S+(\Delta^\vee|V^*)|V),2,\hat{E}(i))\rfloor<
c(\Sigma(S+(\Delta^\vee|V^*)|V),j,\hat{E}(i))$.
It follows $\bar{\Lambda}(j)\cap\Theta(i)^\circ\neq\emptyset$.
Therefore,
$\bar{\Lambda}(j)\cap\Theta(i)^\circ\neq\emptyset
\Leftrightarrow
\bar{\Lambda}(j)\cap(\Theta(i)^\circ\cup\bar{\Theta}(i)^\circ)\neq\emptyset
$.
\end{proof}

\section{Upward subdivisions and the hard height inequalities}
\label{upward}
We define upward subdivisions of a normal fan of a convex pseudo polytope. We construct an upward subdivision by repeating basic subdivisions associated with a compatible mapping. For upward subdivisions, the hard height inequalities hold. These inequalities are essential in our theory.

Let $V$ be any vector space of finite dimension over $\R$ with $\dim V\geq 2$, let $N$ be any lattice of $V$, and let $S$ be any rational convex pseudo polytope over $N$ in $V$ such that $\dim|\Sigma(S|V)|\geq 2$, where $\Sigma(S|V)$ denotes the normal fan of $S$.

By $\mathcal{SF}(V,N,S)$ we denote the set of all pairs $(H,\Phi)$ of a regular cone $H$ of dimension one over the dual lattice $N^*$ of $N$ in the dual vector space $V^*$ of $V$ and a flat regular fan $\Phi$ over $N^*$ in $V^*$ such that $\dim \Phi\geq 2$, $H\in\Phi_1$, $\Phi$ is starry with center in $H$, $|\Phi|\subset |\Sigma(S|V)|$ and $\Sigma(S|V)\hat{\cap}\mathcal{F}(\Delta)$ is $H$-simple for any $\Delta\in\Phi\Mx$.

In this section we assume $\mathcal{SF}(V,N,S)\neq\emptyset$ below.

Note that $\Ht(H,\Phi,S)\in(1/\Den(S/N))\Z_0$ for any $(H,\Phi)\in\mathcal{SF}(V,N,S)$.
Therefore, for any infinite sequence $(H(i),\Phi(i)), i\in\Z_0$ of elements of $\mathcal{SF}(V,N,S)$ such that $\Ht(H(i),\Phi(i),S)\geq\Ht(H(i+1),\Phi(i+1),S)$ for any $i\in\Z_0$, there exists $i_0\in\Z_0$ such that $\Ht(H(i),\Phi(i),S)= \Ht(H(i_0),\Phi(i_0),S)$ for any $i\in\Z_0$ with $i\geq i_0$.

Consider any $(H,\Phi)\in\mathcal{SF}(V,N,S)$.

By $\mathcal{SD}(H,\Phi,S)$ we denote the set of all pairs $(M, F)$ of a non-negative integer $M\in\Z_0$ and a center sequence $F$ of $\Phi$ of length $M$ such that $\dim F(i)=2$ for any $i\in\{1,2,\ldots,M\}$, $F(i)\not\subset|\Phi-(\Phi/H)|$ for any $i\in\{1,2,\ldots,M\}$, and $\Phi* F(1)*F(2)*\cdots *F(M)$ is a subdivision of $\Sigma(S|V)\hat{\cap}\Phi$.

Below, we use induction on $\Ht(H,\Phi,S)$, we will show that $\mathcal{SD}(H,\Phi,S)\neq\emptyset$, and we will define a non-empty subset $\mathcal{USD}(H,\Phi,S)$ of $\mathcal{SD}(H,\Phi,S)$. 

Consider any $(H,\Phi)\in\mathcal{SF}(V,N,S)$ with $\Ht(H,\Phi,S)=0$.
By $0_{\Phi}$ we denote the unique center sequence of $\Phi$ of length $0$.
By Lemma~\ref{compatible2}.1 we know that $\Phi$ is a subdivision of $\Sigma(S|V)\hat{\cap}\Phi$ and therefore $(0, 0_{\Phi})\in\mathcal{SD}(H,\Phi,S)\neq\emptyset$.
We define $\{(0, 0_{\Phi})\}=\mathcal{USD}(H,\Phi,S)$.
Obviously $\emptyset\neq\mathcal{USD}(H,\Phi,S)\subset\mathcal{SD}(H,\Phi,S)$.

Consider any $(H,\Phi)\in\mathcal{SF}(V,N,S)$ with $\Ht(H,\Phi,S)>0$.

By induction hypothesis we can assume that a non-empty subset $\mathcal{USD}(\bar{H},\bar{\Phi},S)$ of $\mathcal{SD}(\bar{H},\bar{\Phi},S)$ is defined for any $(\bar{H},\bar{\Phi})\in \mathcal{SF}(V,N,S)$ satisfying $\Ht(\bar{H}, \bar{\Phi}, S)<\Ht(H,\Phi,S)$.
Below we assume this claim.

The characteristic function $\gamma:(\Phi-(\Phi/H))_1\rightarrow\Q_0$ of $(H,\Phi, S)$ is defined.
Let $m=\sum_{\bar{E}\in(\Phi-(\Phi/H))_1}\lceil\gamma(\bar{E})\rceil\in\Z_+$ and $\bar{m}=\sum_{\bar{E}\in(\Phi-(\Phi/H))_1}\lfloor\gamma(\bar{E})\rfloor\in\Z_0$.
Obviously $\bar{m}\leq m$.

Let $\mathcal{E}$ denote the set of all compatible mappings $E:\{1,2,\ldots,m\}\rightarrow (\Phi-(\Phi/H))_1$ with $S$.
By Lemma~\ref{existcomp}.1 we know $\mathcal{E}\neq\emptyset$.
Consider any $E\in\mathcal{E}$.

Let $F=F(V^*, N^*, H,\Phi,m,E)$ and  $\Omega=\Omega(V^*, N^*, H,\Phi,m,E)$.
We have a mapping
$$F:\{1,2,\ldots,m\}\rightarrow 2^{V^*}.$$
By Lemma~\ref{basic1}.4 and Lemma~\ref{basic1}.2 we know that $F$ is a center sequence of $\Phi$ of length $m$ such that $\dim F(i)=2$, $F(i)\subset E(i)+H$ and $F(i)\not\subset|\Phi-(\Phi/H)|$ for any $i\in\{1,2,\ldots,m\}$.
By definition $\Omega=\Phi* F(1)* F(2)*\ldots*F(m)$.
By Lemma~\ref{property of basic subdivisions}.1 we know that
$\Omega$ is a flat regular fan over $N$ in $V$,
$\dim \Omega=\dim\Phi$,
$\Vect(|\Omega|)=\Vect(|\Phi|)$,
$\Omega$ is an iterated star subdivision of $\Phi$, and
$|\Omega|=|\Phi|$.

We have mappings
$$H=H(V^*, N^*, H,\Phi,m,E):\{1,2,\ldots,m+1\}\rightarrow 2^{V^*},\text{ and}$$
$$\Omega=\Omega(V^*, N^*, H,\Phi,m,E):\{1,2,\ldots,m+1\}\rightarrow 2^{2^{V^*}}.$$

By Lemma~\ref{basic1}.2, Lemma~\ref{property of basic subdivisions}.5 and Theorem~\ref{heightinequality}.4.$(c)$ we know that $H(i)$ is a regular cone of dimension one over $N^*$ in $V^*$, $\Omega(i)$ is a flat regular fan over $N^*$ in $V^*$,
$\dim \Omega(i)=\dim\Phi$,
$ \Vect(|\Omega(i)|)= \Vect(|\Phi|)$,
$H(i)\in\Omega(i)_1$, $\Omega(i)$ is starry with center in $H(i)$, $\Omega(i)\subset\Omega$, $|\Omega(i)|\subset|\Omega|=|\Phi|\subset|\Sigma(S|V)|$, and $\Ht(H(i), \Omega(i), S)<\Ht(H,\Phi,S)$ for any $i\in\{1,2,\ldots,m+1\}$.

By Theorem~\ref{heightinequality}.5 we know that $\Omega\backslash |\Omega(i)|=\Omega(i)=\Sigma(S|V)\hat{\cap}\Omega(i)$ and $\Omega\backslash |\Omega(i)|$ is a subdivision of $\Sigma(S|V)\hat{\cap}\Omega(i)$ for any $i\in\{1,2,\ldots,\bar{m}\}$.

Consider any element $\mu\in\{\bar{m}, \bar{m}+1,\ldots, m+1\}$.

By induction we will show that there exist mappings
$$M: \{\bar{m},\bar{m}+1,\ldots,\mu\}\rightarrow \Z_+,\text{ and}$$
$$\bar{F}:\{1,2,\ldots, M(\mu)\}\rightarrow 2^{V^*}$$ satisfying the following six conditions:
\begin{enumerate}
\item
$M(\bar{m})=m$ and $M(i-1)\leq M(i)$ for any $i\in\{\bar{m}+1,\bar{m}+2,\ldots, \mu\}$.
\item
$\bar{F}(j)=F(j)$ for any $j\in\{1,2,\ldots, m\}$.
\item
$\bar{F}(j)$ is a regular cone over $N^*$ in $V^*$ and $\dim \bar{F}(j)=2$ for and any $j\in\{1,2,\ldots, M(\mu)\}$, and 
$\bar{F}$ is a center sequence of $\Phi$ of length $M(\mu)$.
\item
$\bar{F}(j)\subset|\Omega(i)|$ and $\bar{F} (j)\not\subset|\Omega(i)-(\Omega(i)/H(i))|$ for any $i\in\{\bar{m}+1,\bar{m}+2,\ldots, \mu\}$ and any $j\in\{M(i-1)+1,M(i-1)+2,\ldots, M(i)\}$.
\item
$\Phi*\bar{F}(1)*\bar{F}(2)*\cdots*\bar{F}(M(\mu))\backslash |\Omega(i)|$ is a subdivision of $\Sigma(S|V)\hat{\cap}\Omega(i)$ for any $i\in\{1,2,\ldots,\mu\}$.
\item
Consider any  $i\in\{\bar{m}+1,\bar{m}+2,\ldots,\mu\}$.

Denote $\bar{\Phi}(i)= \Phi*\bar{F}(1)*\bar{F}(2)*\cdots*\bar{F}(M(i-1))\backslash |\Omega(i)|$.

Let $\bar{\bar{F}}(i): \{1,2,\ldots, M(i)-M(i-1)\}\rightarrow 2^{V^*}$ denote the mapping satisfying $\bar{\bar{F}}(i)(j)= \bar{F}(M(i-1)+j)$ for any $j\in\{1,2,\ldots, M(i)-M(i-1)\}$.

$|\bar{\Phi}(i)|=| \Omega(i)|$,
$(H(i), \bar{\Phi}(i))\in\mathcal{SF}(V,N,S)$, 
$\Ht(H(i), \bar{\Phi}(i) , S)<\Ht($\hfill\break$H, \Phi, S)$ and $(M(i)-M(i-1), \bar{\bar{F}}(i))\in\mathcal{USD}(H(i), \bar{\Phi}(i),S)$.
\end{enumerate}

Consider the case $\mu=\bar{m}$.
We define $M(\bar{m})=m$ and $\bar{F}(j)=F(j)$ for any $j\in\{1,2,\ldots, m\}$.
We know that mappings
$$M: \{\bar{m},\bar{m}+1,\ldots,\mu\}=\{\bar{m}\}\rightarrow\Z_+,$$
$$\bar{F}=F:\{1,2,\ldots, M(\mu)\}=\{1,2,\ldots, m\}\rightarrow 2^{V^*}$$ satisfy the above conditions.

Assume that $\mu>\bar{m}$ and there exist mappings
$$M: \{\bar{m},\bar{m}+1,\ldots,\mu-1\}\rightarrow\Z_+,$$
$$\bar{F}:\{1,2,\ldots, M(\mu-1)\}\rightarrow 2^{V^*}$$ satisfying the above conditions in which $\mu$ is replaced by $\mu-1$.
We take mappings $M: \{\bar{m},\bar{m}+1,\ldots,\mu-1\}\rightarrow\Z_+$, $\bar{F}:\{1,2,\ldots, M(\mu-1)\}\rightarrow 2^{V^*}$ satisfying the above conditions in which $\mu$ is replaced by $\mu-1$.

Let $\bar{\Phi}(\mu)=\Phi*\bar{F}(1)* \bar{F}(2)*\cdots*\bar{F}(M(\mu-1))\backslash|\Omega(\mu)|$.
We know that $\bar{\Phi}(\mu)=\Omega*\bar{F}(M(\bar{m})+1)* \bar{F}( M(\bar{m})+2)*\cdots*\bar{F}(M(\mu-1))\backslash|\Omega(\mu)|$,  $\bar{\Phi}(\mu)$ is a flat regular fan over $N^*$ in $V^*$,
$\dim \bar{\Phi}(\mu)=\dim\Phi$,
$\Vect(|\bar{\Phi}(\mu)|)= \Vect(|\Phi|)$,
$\bar{\Phi}(\mu)$ is a subdivision of $\Omega(\mu)$ and 
$|\bar{\Phi}(\mu)|= |\Omega(\mu)|$.
Since $\Sigma(S|V)\hat{\cap}\mathcal{F}(\Theta)$ is semisimple for any $\Theta\in \Omega(\mu)\Mx$, we know that $\Sigma(S|V)\hat{\cap}\mathcal{F}(\Theta)$ is semisimple for any $\Theta\in \bar{\Phi} (\mu)\Mx$

Note that
$$\bar{F}(j)\subset|\Phi-(\Phi/H)|\cup(\bigcup_{i\in\{1,2,\ldots,\mu-1\}}|\Omega(i)|)$$
for any $j\in\{ M(\bar{m})+1, M(\bar{m})+2,\ldots, M(\mu-1)\}$, 
$$(|\Phi-(\Phi/H)|\cup(\bigcup_{i\in\{1,2,\ldots,\mu-1\}}|\Omega(i)|))\cap|\Omega(\mu)|=| \Omega(\mu)-( \Omega(\mu)/H(\mu))|,$$
and $\Omega*\bar{F}(M(\bar{m})+1)* \bar{F}( M(\bar{m})+2)*\cdots*\bar{F}(M(\mu-1))\backslash(|\Phi-(\Phi/H)|\cup$\hfill\break$(\bigcup_{i\in\{1,2,\ldots,\mu-1\}}|\Omega(i)|))$ is a subdivision of $\Sigma(S|V)$.

Let $\ell=\sharp\{j\in\{ M(\bar{m})+1, M(\bar{m})+2,\ldots, M(\mu-1)\}| \bar{F}(j)\subset| \Omega(\mu)-( \Omega(\mu)/H(\mu))|\}$ and let
$\tau:\{1,2,\ldots,\ell\}\rightarrow\{ M(\bar{m})+1, M(\bar{m})+2,\ldots, M(\mu-1)\}$ denote the unique injective mapping preserving the order and satisfying $\tau(\{1,2,\ldots,\ell\})=\{j\in\{ M(\bar{m})+1, M(\bar{m})+2,\ldots, M(\mu-1)\}| \bar{F}(j)\subset| \Omega(\mu)-( \Omega(\mu)/H(\mu))|\}$.
We know that $\bar{\Phi}(\mu)=\Omega(\mu)*\bar{F}\tau(1)* \bar{F}\tau(2)*\cdots*\bar{F}\tau(\ell)$.
Since $\bar{F}\tau(j)\subset|\Omega(\mu)-(\Omega(\mu)/H(\mu))|$ for any $j\in\{1,2,\ldots,\ell\}$, we know that $H(\mu)\in \bar{\Phi}(\mu)_1$, $\bar{\Phi}(\mu)$ is starry with center in $H(\mu)$, 
$|\bar{\Phi}(\mu)|=|\Omega(\mu)|$,
$|\bar{\Phi}(\mu)- (\bar{\Phi}(\mu)/H(\mu))|=
|\Omega(\mu)- (\Omega(\mu)/H(\mu))|$ and
$\Ht(H(\mu), \bar{\Phi}(\mu), S)=\Ht(H(\mu), \Omega(\mu), S)<\Ht(H,\Phi,S)$.

Consider any $\Theta\in \bar{\Phi}(\mu)\Mx$.
$\Sigma(S|V)\hat{\cap}\mathcal{F}(\Theta)$ is semisimple.
$\Theta\in \bar{\Phi}(\mu)/H(\mu)$, and $H(\mu)\Op|\Theta\in \bar{\Phi}(\mu)-( \bar{\Phi}(\mu)/H(\mu))$.
On the other hand, we know
$\bar{\Phi}(\mu)-( \bar{\Phi}(\mu)/H(\mu))
=\bar{\Phi}(\mu)\backslash|\Omega(\mu)-( \Omega(\mu)/H(\mu))|
=\Omega*\bar{F}(M(\bar{m})+1)* \bar{F}( M(\bar{m})+2)*\cdots*\bar{F}(M(\mu-1))\backslash|\Omega(\mu)-( \Omega(\mu)/H(\mu))|
=(\Omega*\bar{F}(M(\bar{m})+1)* \bar{F}( M(\bar{m})+2)*\cdots*\bar{F}(M(\mu-1))\backslash(|\Phi-(\Phi/H)|\cup(\bigcup_{i\in\{1,2,\ldots,\mu-1\}}|\Omega(i)|)))\backslash|\Omega(\mu)-( \Omega(\mu)/H(\mu))|$, and $\bar{\Phi}(\mu)-( \bar{\Phi}(\mu)/H(\mu))$ is a subdivision of $\Sigma(S|V)$.
Therefore $\Sigma(S|V)\hat{\cap}\mathcal{F}( H(\mu)\Op|\Theta)= \mathcal{F}( H(\mu)\Op|\Theta)$, and we know that $\Sigma(S|V)\hat{\cap}\mathcal{F}(\Theta)$ is $H(\mu)$-simple.

We know that $(H(\mu), \bar{\Phi}(\mu))\in\mathcal{SF}(V,N,S)$, $\Ht(H(\mu), \bar{\Phi}(\mu), S)<\Ht(H, \Phi,$\break$ S)$, and a non-empty subset $\mathcal{USD}( H(\mu), \bar{\Phi}(\mu),S)$ of $\mathcal{SD}( H(\mu), \bar{\Phi}(\mu),S)$ is defined.

Take any element $(L, G)\in \mathcal{USD}( H(\mu), \bar{\Phi}(\mu),S)$. 
Put $M(\mu)=M(\mu-1)+L$, and put $\bar{F}(j)=G(j-M(\mu-1))$ for any $j\in\{M(\mu-1)+1,M(\mu-1)+2,\ldots,M(\mu)\}$.
We obtain extended mappings
$$M: \{\bar{m},\bar{m}+1,\ldots,\mu\}\rightarrow\Z_+,$$
$$\bar{F}:\{1,2,\ldots, M(\mu)\}\rightarrow 2^{V^*}$$ satisfying the above conditions.

By induction on $\mu$ we know that there exist mappings 
$$M: \{\bar{m},\bar{m}+1,\ldots,m+1\}\rightarrow\Z_+,$$
$$\bar{F}:\{1,2,\ldots, M(m+1)\}\rightarrow 2^{V^*}$$ satisfying the following six conditions:
\begin{enumerate}
\item
$M(\bar{m})=m$ and $M(i-1)\leq M(i)$ for any $i\in\{\bar{m}+1,\bar{m}+2,\ldots, m+1\}$.
\item
$\bar{F}(j)=F(j)$ for any $j\in\{1,2,\ldots, m\}$
\item
$\bar{F}(j)$ is a regular cone over $N^*$ in $V^*$ and $\dim \bar{F}(j)=2$ for and any $j\in\{1,2,\ldots, M(m+1)\}$, and 
$\bar{F}$ is a center sequence of $\Phi$ of length $M(m+1)$.
\item
$\bar{F}(j)\subset|\Omega(i)|$ and $\bar{F} (j)\not\subset|\Omega(i)-(\Omega(i)/H(i))|$ for any $i\in\{\bar{m}+1,\bar{m}+2,\ldots, m+1\}$ and any $j\in\{M(i-1)+1,M(i-1)+2,\ldots, M(i)\}$.
\item
$\Phi*\bar{F}(1)*\bar{F}(2)*\cdots*\bar{F}(M(m+1))\backslash |\Omega(i)|$ is a subdivision of $\Sigma(S|V)\hat{\cap}\Omega(i)$ for any $i\in\{1,2,\ldots,m+1\}$.
\item
Consider any $i\in\{\bar{m}+1,\bar{m}+2,\ldots,m+1\}$.

Denote $\bar{\Phi}(i)= \Phi*\bar{F}(1)*\bar{F}(2)*\cdots*\bar{F}(M(i-1))\backslash |\Omega(i)|$.

Let $\bar{\bar{F}}(i): \{1,2,\ldots, M(i)-M(i-1)\}\rightarrow 2^{V^*}$ denote the mapping satisfying $\bar{\bar{F}}(i)(j)= \bar{F}(M(i-1)+j)$ for any $j\in\{1,2,\ldots, M(i)-M(i-1)\}$.

$|\bar{\Phi}(i)|=|\Omega(i)|$,
$(H(i), \bar{\Phi}(i))\in\mathcal{SF}(V,N,S)$, 
$\Ht(H(i), S, \bar{\Phi}(i))<$\hfill\break$\Ht(H,\Phi,S)$ and $(M(i)-M(i-1), \bar{\bar{F}}(i))\in\mathcal{USD}(H(i), \bar{\Phi}(i),S)$.
\end{enumerate}

Put
\begin{equation*}\begin{split}
\mathcal{USD}(E, H,\Phi,S)&=
\{(M(m+1),\bar{F})|\\
&\qquad\quad M: \{\bar{m},\bar{m}+1,\ldots,m+1\}\rightarrow\Z_+,\\
&\qquad\quad \bar{F}:\{1,2,\ldots, M(m+1)\}\rightarrow 2^{V^*},\\
&\qquad\quad M\text{ and }\bar{F}\text{ satisfy the above conditions.}\}.
\end{split}\end{equation*}
We know that $\emptyset\neq\mathcal{USD}(E, H,\Phi,S)\subset\mathcal{SD}(H,\Phi,S)$.

Recall that $E\in\mathcal{E}$ is an arbitrary element.
Define
$$\mathcal{USD}(H,\Phi,S)=\bigcup_{ E\in\mathcal{E}}\mathcal{USD}(E, H,\Phi,S).$$
We know that $\emptyset\neq\mathcal{USD}(H,\Phi,S)\subset\mathcal{SD}(H,\Phi,S)$.

Consider any $(H,\Phi)\in\mathcal{SF}(V,N,S)$ and any $(M, F)\in\mathcal{USD}(H,\Phi,S)$.
We call $F$ an \emph{upward center sequence} of $(H,\Phi, S)$, and
we call $\Phi* F(1)*F(2)*\cdots *F(M)$ an \emph{upward subdivision} of $(H,\Phi, S)$.

\begin{theorem}
\label{usd1}
Assume $\mathcal{SF}(V,N,S)\neq\emptyset$, and consider any $(H,\Phi)\in\mathcal{SF}(V,N,S)$.
\begin{enumerate}
\item
$\mathcal{USD}(H,\Phi,S)\neq\emptyset$.
\item
For any $(M, F)\in\mathcal{USD}(H,\Phi,S)$, 
$F$ is a center sequence of $\Phi$ of length $M$,
$\dim F(i)=2$ and $F(i)\not\subset|\Phi-(\Phi/H)|$ for any $i\in\{1,2,\ldots,M\}$, and $\Phi* F(1)*F(2)*\cdots*F(M)$ is a subdivision of $\Sigma(S|V)\hat{\cap}\Phi$.
\item
By $0_{\Phi}$ we denote the unique center sequence of $\Phi$ of length $0$.
The following three conditions are equivalent:
\begin{enumerate}
\item
$\Ht(H,\Phi,S)=0$.
\item
$(0, 0_{\Phi})\in\mathcal{USD}(H,\Phi,S)$.
\item
$\{(0, 0_{\Phi})\}=\mathcal{USD}(H,\Phi,S)$.
\end{enumerate}
\end{enumerate}

Consider any $(M, F)\in\mathcal{USD}(H,\Phi,S)$.
We denote $\widetilde{\Phi}=\Phi*F(1)*F(2)*\cdots*F(M)$ for simplicity.

\begin{enumerate}
\setcounter{enumi}{3}
\item
For any $\Delta\in\Phi_0\cup\Phi_1$,
$\widetilde{\Phi}\backslash\Delta=\Sigma(S|V)\hat{\cap}\mathcal{F}(\Delta)=\mathcal{F}(\Delta)$.
\item
$\widetilde{\Phi}\backslash|\Phi-(\Phi/H)|=(\Sigma(S|V)\hat{\cap}\Phi)\backslash|\Phi-(\Phi/H)|= \Phi-(\Phi/H)$.
\item
Consider any subset $\hat{\Phi}$ of $\Phi$ satisfying $\dim \hat{\Phi}=\dim \Vect(|\hat{\Phi}|)\geq 2$, $\hat{\Phi}\Mx=\hat{\Phi}^0$, $H\in\hat{\Phi}_1$ and $\hat{\Phi}=(\hat{\Phi}/H)\Fc$.

Let $\hat{M}=\sharp\{i\in\{1,2,\ldots,M\}|F(i)\subset|\hat{\Phi}|\}$, and $\tau:\{1,2,\ldots,\hat{M}\}\rightarrow \{1,2,\ldots,M\}$ denote the unique injective mapping preserving the order and satisfying $\tau(\{1,2,\ldots,\hat{M}\})=\{i\in\{1,2,\ldots,M\}|F(i)\subset|\hat{\Phi}|\}$.

Then, $(H, \hat{\Phi})\in\mathcal{SF}(V,N,S)$, $(\hat{M},F\tau)\in\mathcal{USD}(H, \hat{\Phi},S)$, and
$\widetilde{\Phi}\backslash|\hat{\Phi}|=\hat{\Phi}*F\tau (1)*F\tau(2)*\cdots*F\tau(\hat{M})$.
\item
Consider any $\Delta\in\Phi_2/H$.

$|\widetilde{\Phi}\backslash\Delta|=|\Sigma(S|V)\hat{\cap}\mathcal{F}(\Delta)|=\Delta$.

$\widetilde{\Phi}\backslash\Delta$ is the minimum regular subdivision of $\Sigma(S|V)\hat{\cap}\mathcal{F}(\Delta)$ over $N^*$ in $V^*$, in other words, the following three conditions hold:
\begin{enumerate}
\item
$\widetilde{\Phi}\backslash\Delta$ is a regular fan over $N^*$ in $V^*$.
\item
$\widetilde{\Phi}\backslash\Delta$ is a subdivision of $\Sigma(S|V)\hat{\cap}\mathcal{F}(\Delta)$.
\item
If $\Xi$ is a regular fan over $N^*$ in $V^*$ and $\Xi$ is a subdivision of $\Sigma(S|V)\hat{\cap}\mathcal{F}(\Delta)$, then $\Xi$ is a subdivision of $\widetilde{\Phi}\backslash\Delta$.
\end{enumerate}
\end{enumerate}

Below, we consider the case $\Ht(H,\Phi,S)>0$.
Assume $\Ht(H,\Phi,S)>0$.

The characteristic function $\gamma:(\Phi-(\Phi/H))_1\rightarrow\Q_0$ of $(H,\Phi, S)$ is defined.
Let $m=\sum_{\bar{E}\in(\Phi-(\Phi/H))_1}\lceil\gamma(\bar{E})\rceil\in\Z_+$ and $\bar{m}=\sum_{\bar{E}\in(\Phi-(\Phi/H))_1}\lfloor\gamma(\bar{E})\rfloor\in\Z_0$.
Obviously $\bar{m}\leq m$.
\begin{enumerate}
\setcounter{enumi}{7}
\item
$m\leq M$.
\item
There exists uniquely a pair $(E,M)$ of a compatible mapping
$E:\{1,2,\ldots,$\hfill\break$m\}\rightarrow (\Phi-(\Phi/H))_1$
with $S$ and a mapping 
$M:\{\bar{m},\bar{m}+1,\ldots,m+1\}\rightarrow\Z_+$
satisfying the following three conditions.
We denote 
$$F_E=F(V^*,N^*,H,\Phi,m,E):\{1,2,\ldots,m\}\rightarrow 2^{ V^*}$$ 
$$H_E=H(V^*,N^*,H,\Phi,m,E):\{1,2,\ldots,m+1\}\rightarrow 2^{ V^*},\text{ and}$$
$$\Omega_E=\Omega(V^*,N^*,H,\Phi,m,E):\{1,2,\ldots,m+1\}\rightarrow 2^{2^{ V^*}}:$$
\begin{enumerate}
\item
$F(j)=F_E(j)$ for any $j\in\{1,2,\ldots, m\}$.
\item
$M(\bar{m})=m$, $M(m+1)=M$ and $M(i-1)\leq M(i)$ for any $i\in\{\bar{m}+1,\bar{m}+2,\ldots, m+1\}$.
\item
$F(j)\subset|\Omega_E(i)|$ and $F(j)\not\subset|\Omega_E(i)-(\Omega_E(i)/H_E(i))|$ for any $i\in\{\bar{m}+1,\bar{m}+2,\ldots, m+1\}$ and any $j\in\{M(i-1)+1,M(i-1)+2,\ldots, M(i)\}$.
\end{enumerate}
\end{enumerate}

We take the unique pair $(E,M)$ of a compatible mapping
$E:\{1,2,\ldots,m\}\rightarrow (\Phi-(\Phi/H))_1$
with $S$ and a mapping 
$M:\{\bar{m},\bar{m}+1,\ldots,m+1\}\rightarrow\Z_+$
satisfying the above three conditions, and we denote 
$$F_E=F(V^*,N^*,H,\Phi,m,E):\{1,2,\ldots,m\}\rightarrow 2^{V^*}$$ 
$$H_E=H(V^*,N^*,H,\Phi,m,E):\{1,2,\ldots,m+1\}\rightarrow 2^{V^*},\text{ and}$$
$$\Omega_E=\Omega(V^*,N^*,H,\Phi,m,E):\{1,2,\ldots,m+1\}\rightarrow 2^{2^{V^*}}.$$
We put $M(i)=m$ for any $i\in\{0,1,\ldots,\bar{m}-1\}$. We obtain an extension $M:\{0,1,\ldots,m+1\}\rightarrow\Z_+$ of $M:\{\bar{m},\bar{m}+1,\ldots,m+1\}\rightarrow\Z_+$.
$M(0)=m$, $M(m+1)=M$, $M(i-1)\leq M(i)$ for any $i\in\{1,2,\ldots,m+1\}$.

Consider any $i\in\{1, 2,\ldots,m+1\}$.
We denote $\bar{\Phi}(i)= \Phi*F(1)*F(2)*\cdots*F(M(i-1))\backslash |\Omega_E(i)|$.
Let $\bar{F}(i): \{1,2,\ldots, M(i)-M(i-1)\}\rightarrow 2^{V^*}$ denote the mapping satisfying $\bar{F}(i)(j)= F(M(i-1)+j)$ for any $j\in\{1,2,\ldots, M(i)-M(i-1)\}$.

\begin{enumerate}
\setcounter{enumi}{9}
\item
For any $i\in\{1,2,\ldots,m+1\}$,
$|\bar{\Phi}(i)|=|\Omega_E(i)|$,
$(H_E(i), \bar{\Phi}(i))\in\mathcal{SF}(V,N,S)$, 
$\Ht(H_E(i), \bar{\Phi}(i),S)<\Ht(H,\Phi,S)$ and $(M(i)-M(i-1), \bar{F}(i))\in$\hfill\break$\mathcal{USD}(H_E(i), \bar{\Phi}(i),S)$.
\item
For any $i\in\{1, 2,\ldots,m+1\}$,
$\widetilde{\Phi}\backslash|\Omega_E(i)|=
(\Phi*F(1)*F(2)*\cdots*F(M(i))) $\hfill\break$\backslash|\Omega_E(i)|=
\bar{\Phi}(i)*\bar{F}(1)*\bar{F}(2)*\cdots*\bar{F}(M(i)- M(i-1))$.
\item
For any $i\in\{1,2,\ldots,\bar{m}\}$, $\Ht(H_E(i), \bar{\Phi}(i),S)=0$, $ M(i)-M(i-1)=0$, and $\widetilde{\Phi}\backslash|\Omega_E(i)|=
(\Phi*F(1)*F(2)*\cdots*F(m))\backslash|\Omega_E(i)|=\bar{\Phi}(i)*\bar{F}(1)*\bar{F}(2)*\cdots*\bar{F}(M(i)- M(i-1))
=\bar{\Phi}(i)= \Omega_E(i)$.
\item

$$\widetilde{\Phi}=
(\Phi-(\Phi/H))\cup(
\bigcup_{i\in\{1,2,\ldots,m+1\}}
(
(\widetilde{\Phi}\backslash|\Omega_E(i)|)
-(\widetilde{\Phi}\backslash|\Omega_E(i)-(\Omega_E(i)/H_E(i))|)
)).$$
$$(\Phi-(\Phi/H))\cap(
(\widetilde{\Phi}\backslash|\Omega_E(i)|)
-(\widetilde{\Phi}\backslash|\Omega_E(i)-(\Omega_E(i)/H_E(i))|)
)=\emptyset,$$
for any $i\in\{1,2,\ldots,m+1\}$.

\begin{equation*}\begin{split}(
&(\widetilde{\Phi}\backslash|\Omega_E(i)|)
-(\widetilde{\Phi}\backslash|\Omega_E(i)-(\Omega_E(i)/H_E(i))|)
)\\
&\qquad\qquad\cap(
(\widetilde{\Phi}\backslash|\Omega_E(j)|)
-(\widetilde{\Phi}\backslash|\Omega_E(j)-(\Omega_E(j)/H_E(j))|)
)=\emptyset,
\end{split}\end{equation*}
for any $i\in\{1,2,\ldots,m+1\}$ and any $j\in\{1,2,\ldots,m+1\}$ with $i\neq j$.
\item
If $\Gamma\in\widetilde{\Phi}_1$ and $\Gamma\not\subset|\Phi-(\Phi/H)|$,
then there exists uniquely $i\in\{1,2,\ldots,m+1\}$ satisfying
$\Gamma\subset|\Omega_E(i)|$ and $\Gamma\not\subset|\Omega_E(i)-(\Omega_E(i)/H_E(i))|$.

If $\Gamma\in\widetilde{\Phi}_1$, $i\in\{1,2,\ldots,m+1\}$, $\Gamma\subset|\Omega_E(i)|$ and $\Gamma\not\subset|\Omega_E(i)-(\Omega_E(i)/H_E(i))|$, then
$\Gamma\not\subset|\Phi-(\Phi/H)|$.

For any $i\in\{1,2,\ldots,\bar{m}\}$, 
$$\{\Gamma\in\widetilde{\Phi}_1|\Gamma\subset|\Omega_E(i)|, \Gamma\not\subset|\Omega_E(i)-(\Omega_E(i)/H_E(i))|\}=
\{H_E(i)\}.$$

For any $i\in\{\bar{m}+1, \bar{m}+2,\ldots, m+1\}$, 
\begin{equation*}\begin{split}
&\{\Gamma\in\widetilde{\Phi}_1|\Gamma\subset|\Omega_E(i)|, \Gamma\not\subset|\Omega_E(i)-(\Omega_E(i)/H_E(i))|\}=\\
&\qquad\qquad\{H_E(i)\}\cup\{\R_0b_{F(j)/N^*}|j\in\Z, M(i-1)<j\leq M(i)\}.
\end{split}\end{equation*}

For any $i\in\{1, 2, \ldots, m+1\}$, 
$$\sharp\{\Gamma\in\widetilde{\Phi}_1|\Gamma\subset|\Omega_E(i)|, \Gamma\not\subset|\Omega_E(i)-(\Omega_E(i)/H_E(i))|\}=M(i)-M(i-1)+1.$$
\end{enumerate}
\end{theorem}

Consider any $(H,\Phi)\in\mathcal{SF}(V,N,S)$ and any $(M,F)\in\mathcal{USD}(H,\Phi,S)$.
We denote $\widetilde{\Phi}=\Phi*F(1)*F(2)*\cdots*F(M)$ and
$\widetilde{\Phi}_1^\circ=\{\Gamma\in\widetilde{\Phi}_1|
\Gamma\not\subset|\Phi-(\Phi/H)|\}$ for simplicity.
We know $\sharp\widetilde{\Phi}_1^\circ=M+1$.

By induction on $\Ht(H,\Phi,S)$ we define three mappings
$$I(V,N, H,\Phi,S, M,F):\widetilde{\Phi}_1^\circ\rightarrow
\{1,2,\ldots,M+1\},$$
$$\Psi(V,N, H,\Phi,S, M,F):
\widetilde{\Phi}_1^\circ\rightarrow 2^{2^{V^*}},\text{ and}$$
$$\Psi^\circ(V,N, H,\Phi,S, M,F):
\widetilde{\Phi}_1^\circ\rightarrow 2^{2^{V^*}}$$
such that $\Psi(V,N, H, \Phi,S, M,F)(\Gamma)$ is a regular fan over $N^*$ in $V^*$ and $\Psi^\circ(V,N, H,$\break$\Phi,S, M,F) (\Gamma)\subset\Psi(V,N, H,\Phi,S, M,F)(\Gamma)\subset\widetilde{\Phi}$ for any $\Gamma\in \widetilde{\Phi}_1^\circ$.

Consider the case $\Ht(H,\Phi,S)=0$.

We know $M=0$, $\widetilde{\Phi}=\Phi$, and
$\widetilde{\Phi}_1^\circ=\{H\}$.
We define
\begin{equation*}\begin{split}
I(V,N, H,\Phi,S, M,F)(H)&=1,\\
\Psi(V,N, H,\Phi,S, M,F)(H)&= \Phi,\\
\Psi^\circ(V,N, H,\Phi,S, M,F)(H)&=\Phi.
\end{split}\end{equation*}

Consider the case $\Ht(H,\Phi,S)>0$.

The characteristic function $\gamma:(\Phi-(\Phi/H))_1\rightarrow\Q_0$ of $(H,\Phi, S)$ is defined.
Let $m=\sum_{\bar{E}\in(\Phi-(\Phi/H))_1}\lceil\gamma(\bar{E})\rceil\in\Z_+$ and $\bar{m}=\sum_{\bar{E}\in(\Phi-(\Phi/H))_1}\lfloor\gamma(\bar{E})\rfloor\in\Z_0$.
We know $\bar{m}\leq m\leq M$.

We take the unique pair $(E,M)$ of a compatible mapping
$E:\{1,2,\ldots,m\}\rightarrow (\Phi-(\Phi/H))_1$
with $S$ and a mapping 
$M:\{\bar{m},\bar{m}+1,\ldots,m+1\}\rightarrow\Z_+$
satisfying the following three conditions.
We denote 
$$F_E=F(V^*,N^*,H,\Phi,m,E):\{1,2,\ldots,m\}\rightarrow 2^{ V^*}$$ 
$$H_E=H(V^*,N^*,H,\Phi,m,E):\{1,2,\ldots,m+1\}\rightarrow 2^{ V^*},\text{ and}$$
$$\Omega_E=\Omega(V^*,N^*,H,\Phi,m,E):\{1,2,\ldots,m+1\}\rightarrow 2^{2^{ V^*}}.$$
\begin{enumerate}
\item
$F(j)=F_E(j)$ for any $j\in\{1,2,\ldots, m\}$.
\item
$M(\bar{m})=m$, $M(m+1)=M$ and $M(i-1)\leq M(i)$ for any $i\in\{\bar{m}+1,\bar{m}+2,\ldots, m+1\}$.
\item
$F(j)\subset|\Omega_E(i)|$ and $F(j)\not\subset|\Omega_E(i)-(\Omega_E(i)/H_E(i))|$ for any $i\in\{\bar{m}+1,\bar{m}+2,\ldots, m+1\}$ and any $j\in\{M(i-1)+1,M(i-1)+2,\ldots, M(i)\}$.
\end{enumerate}

For any $i\in\{0,1,\ldots, \bar{m}-1\}$, we put $M(i)=m$.
We obtain an extension $M:\{0,1,\ldots,m+1\}\rightarrow\Z_+$ of $M:\{\bar{m},\bar{m}+1,\ldots,m+1\}\rightarrow\Z_+$.
For any $i\in\{1,2,\ldots,m+1\}$, we denote
$\bar{\Phi}(i)=(\Phi*F(1)*F(2)*\cdots*F(M(i-1)))\backslash|\Omega_E(i)|$, and we take the mapping $\bar{F}(i):\{1,2,\ldots, M(i)-M(i-1)\}\rightarrow 2^{V^*}$ satisfying $\bar{F}(i)(j)=F(M(i-1)+j)$ for any $j\in\{1,2,\ldots, M(i)-M(i-1)\}$.
We know $(H_E(i), \bar{\Phi}(i))\in\mathcal{SF}(V,N,S)$,
$|\bar{\Phi}(i)|=|\Omega_E(i)|$,
$|\bar{\Phi}(i)-( \bar{\Phi}(i)/H_E(i))|=
|\Omega_E(i)-( \Omega_E(i)/H_E(i))|$,
$\Ht(H_E(i), \bar{\Phi}(i),S)<\Ht(H,\Phi,S)$, $(M(i)-M(i-1),\bar{F}(i))\in\mathcal{USD}( H_E(i), \bar{\Phi}(i),S)$, and
$\bar{\Phi}(i)* \bar{F}(1)* \bar{F}(2)*\cdots*\bar{F}(M(i)-M(i-1))=
\widetilde{\Phi}\backslash|\Omega_E(i)|$
for any $i\in\{1,2,\ldots,m+1\}$.
We denote
$$G_E=G(V^*,N^*,H,\Phi,m,E):\{1,2,\ldots,m\}\rightarrow 2^{V^*}.$$

Consider any $\Gamma\in\widetilde{\Phi}_1^\circ$.
Take the unique $i\in\{1,2,\ldots,m+1\}$ satisfying $\Gamma\subset|\Omega_E(i)|$ and $\Gamma\not\subset|\Omega_E(i)-( \Omega_E(i)/H_E(i))|$.
We know $\Gamma\in(\bar{\Phi}(i)* \bar{F}(1)* \bar{F}(2)*\cdots*\bar{F}(M(i)-M(i-1)))_1$ and $\Gamma\not\subset|\bar{\Phi}(i)-( \bar{\Phi}(i)/H_E(i))|$.
By induction hypothesis we know that 
$I(V,N, H_E(i),\bar{\Phi},S, M(i)-M(i-1),\bar{F}(i))(\Gamma)$,
$\Psi(V,N, H_E(i),\bar{\Phi},S, M(i)-M(i-1),\bar{F}(i))(\Gamma)$ and
$\Psi^\circ(V,N, H_E(i),\bar{\Phi},S, M(i)-M(i-1),\bar{F}(i))(\Gamma)$ are defined.

We define
\begin{equation*}\begin{split}
&I(V,N, H,\Phi,S, M,F)(\Gamma)=\\
&\qquad i-1+M(i-1)-m+ I(V,N, H_E(i),\bar{\Phi}(i),S, M(i)-M(i-1),\bar{F}(i))(\Gamma),\\
&\Psi(V,N, H,\Phi,S, M,F)(\Gamma)=\Psi(V,N, H_E(i),\bar{\Phi}(i),S, M(i)-M(i-1),\bar{F}(i))(\Gamma).
\end{split}\end{equation*}
In case $i\neq m+1$, we define
\begin{equation*}\begin{split}
&\Psi^\circ(V,N, H,\Phi,S, M,F)(\Gamma)=\\
&\qquad \{\Theta\in\Psi^\circ(V,N, H_E(i),\bar{\Phi}(i),S, M(i)-M(i-1),\bar{F}(i))(\Gamma)|\Theta^\circ\subset|\Omega_E/G_E(i)|^\circ\}.
\end{split}\end{equation*}
In case $i=m+1$, we define
$$\Psi^\circ(V,N, H,\Phi,S, M,F)(\Gamma)=
\Psi^\circ(V,N, H_E(i),\bar{\Phi}(i),S, M(i)-M(i-1),\bar{F}(i))(\Gamma).$$

We call the mapping $I(V,N, H,\Phi,S, M,F)$ the $H$\emph{-ordered enumeration} of $\widetilde{\Phi}_1^\circ$.
Let $\Gamma\in\widetilde{\Phi}_1^\circ$.
We call the regular fan $\Psi(V,N, H,\Phi,S, M,F)(\Gamma)$ the $H$\emph{-lower part} of  $\widetilde{\Phi}$ below $\Gamma$, and we call the subset $\Psi^\circ(V,N, H,\Phi,S, M,F)(\Gamma)$ the $H$\emph{-lower main part} of  $\widetilde{\Phi}$ below $\Gamma$.

\begin{lemma}
\label{lower parts}
Assume $\mathcal{SF}(V,N,S)\neq\emptyset$.
Consider any $(H,\Phi)\in\mathcal{SF}(V,N,S)$ and any $(M,F)\in\mathcal{USD}(H,\Phi,S)$.
We denote $\widetilde{\Phi}=\Phi*F(1)*F(2)*\cdots*F(M)$,
$\widetilde{\Phi}_1^\circ=\{\Gamma\in\widetilde{\Phi}_1|\Gamma\not\subset|\Phi-(\Phi/H)|\}$,
$$I= I(V,N, H,\Phi,S, M,F): \widetilde{\Phi}_1^\circ
\rightarrow\{1,2,\ldots, M+1\},$$
$$\Psi=\Psi(V,N, H,\Phi,S, M,F): \widetilde{\Phi}_1^\circ
\rightarrow 2^{2^{V^*}},\text{ and}$$
$$\Psi^\circ=\Psi^\circ(V,N, H,\Phi,S, M,F): \widetilde{\Phi}_1^\circ
\rightarrow 2^{2^{V^*}}.$$
\begin{enumerate}
\item
The mapping $I$ is bijective. $H\in\widetilde{\Phi}_1^\circ$.
$I(H)=M+1$.
\item
Consider any $\Gamma\in \widetilde{\Phi}_1^\circ$.

$\Psi(\Gamma)$ is a flat regular fan over $N^*$ in $V^*$.
$\Psi(\Gamma)\subset \widetilde{\Phi}$.
$\dim \Psi(\Gamma)=\dim\Phi\geq 2$.
$\Vect(|\Psi(\Gamma)|)=\Vect(|\Phi|)$.
$\Gamma\in\Psi(\Gamma)_1$.
$\Psi(\Gamma)$ is starry with center in $\Gamma$.
$\Psi(\Gamma)\Mx\subset\Psi^\circ(\Gamma)\subset\Psi(\Gamma)$.

$\Gamma\in\Psi^\circ(\Gamma)\Leftrightarrow\Psi^\circ(\Gamma)=\Psi(\Gamma)\Leftrightarrow\Gamma=H$.
\item
$$|\Phi|=|\Phi-(\Phi/H)|\cup(
\bigcup_{\Gamma\in\widetilde{\Phi}_1^\circ}|\Psi(\Gamma)/\Gamma|^\circ).$$

For any $\Gamma\in\widetilde{\Phi}_1^\circ$,
$|\Phi-(\Phi/H)|\cap|\Psi(\Gamma)/\Gamma|^\circ=\emptyset$.

For any $\Gamma\in\widetilde{\Phi}_1^\circ$ and any $\bar{\Gamma}\in\widetilde{\Phi}_1^\circ$ with $\Gamma\neq\bar{\Gamma}$,
$|\Psi(\Gamma)/\Gamma|^\circ\cap|\Psi(\bar{\Gamma})/\bar{\Gamma}|^\circ=\emptyset$.
\item
$$\widetilde{\Phi}=(\Phi-(\Phi/H))\cup(
\bigcup_{\Gamma\in\widetilde{\Phi}_1^\circ}(\Psi(\Gamma)/\Gamma)).$$

For any $\Gamma\in\widetilde{\Phi}_1^\circ$,
$(\Phi-(\Phi/H))\cap(\Psi(\Gamma)/\Gamma)=\emptyset$.

For any $\Gamma\in\widetilde{\Phi}_1^\circ$ and any $\bar{\Gamma}\in\widetilde{\Phi}_1^\circ$ with $\Gamma\neq\bar{\Gamma}$,
$(\Psi(\Gamma)/\Gamma)\cap(\Psi(\bar{\Gamma})/\bar{\Gamma})=\emptyset$.
\item
$$\widetilde{\Phi}\Mx=
\bigcup_{\Gamma\in\widetilde{\Phi}_1^\circ}\Psi(\Gamma)\Mx.$$

For any $\Gamma\in\widetilde{\Phi}_1^\circ$ and any $\bar{\Gamma}\in\widetilde{\Phi}_1^\circ$ with $\Gamma\neq\bar{\Gamma}$,
$\Psi(\Gamma)\Mx\cap\Psi(\bar{\Gamma})\Mx=\emptyset$.
\item
$$|\Phi|=
\bigcup_{\Gamma\in\widetilde{\Phi}_1^\circ}|\Psi^\circ(\Gamma) |^\circ.$$

For any $\Gamma\in\widetilde{\Phi}_1^\circ$ and any $\bar{\Gamma}\in\widetilde{\Phi}_1^\circ$ with $\Gamma\neq\bar{\Gamma}$,
$|\Psi^\circ(\Gamma) |^\circ \cap|\Psi^\circ(\bar{\Gamma}) |^\circ =\emptyset$.
\item
$$\widetilde{\Phi}=
\bigcup_{\Gamma\in\widetilde{\Phi}_1^\circ}\Psi^\circ(\Gamma).$$

For any $\Gamma\in\widetilde{\Phi}_1^\circ$ and any $\bar{\Gamma}\in\widetilde{\Phi}_1^\circ$ with $\Gamma\neq\bar{\Gamma}$,
$\Psi^\circ(\Gamma)\cap\Psi^\circ(\bar{\Gamma})=\emptyset$.
\end{enumerate}

For any $i\in\{0,1,\ldots,M+1\}$, we denote
\begin{equation*}\begin{split}
X(i)=&|\Phi-(\Phi/H)|\cup(
\bigcup_{\Gamma\in\widetilde{\Phi}_1^\circ, I(\Gamma)\leq i}
|\Psi(\Gamma)/\Gamma)|^\circ)\subset V^*,\\
Y(i)=&\bigcup_{\Gamma\in\widetilde{\Phi}_1^\circ, I(\Gamma)> i}
|\Psi^\circ(\Gamma)|^\circ\subset V^*.
\end{split}\end{equation*}

\begin{enumerate}
\setcounter{enumi}{7}
\item
$X(0)= |\Phi-(\Phi/H)|$.
$X(M+1)=|\Phi|$.
$Y(0)=|\Phi|$.
$Y(M+1)=\emptyset$.

For any $i\in\{1,2,\ldots,M+1\}$, $X(i-1)\subset X(i)$, $X(i-1)\neq X(i)$,\hfill\break
$Y(i-1)\supset Y(i)$, and $Y(i-1)\neq Y(i)$.
\item
For any $i\in\{1,2,\ldots,M+1\}$,
\begin{equation*}\begin{split}
X(i-1)\cap|\Psi I^{-1}(i)|&=
|\Psi I^{-1}(i)|-| \Psi I^{-1}(i)/ I^{-1}(i)|^\circ,\\
Y(i)\cap|\Psi I^{-1}(i)|&= |\Psi I^{-1}(i)|- |\Psi^\circ I^{-1}(i)|^\circ.
\end{split}\end{equation*}
\item
For any $i\in\{0,1,\ldots,M+1\}$,
\begin{equation*}\begin{split}
X(i)&=|\Phi-(\Phi/H)|\cup(
\bigcup_{\Gamma\in\widetilde{\Phi}_1^\circ, I(\Gamma)\leq i}
|\Psi(\Gamma)|),\\
Y(i)&=\bigcup_{\Gamma\in\widetilde{\Phi}_1^\circ, I(\Gamma)> i}
|\Psi(\Gamma)|,
\end{split}\end{equation*}
and $X(i)$ and $Y(i)$ are closed subsets of $V^*$.

\item
Consider any subset $\hat{\Phi}$ of $\Phi$ satisfying $\dim \hat{\Phi}=\dim \Vect(|\hat{\Phi}|)\geq 2$, $\hat{\Phi}\Mx=\hat{\Phi}^0$, $H\in\hat{\Phi}_1$ and $\hat{\Phi}=(\hat{\Phi}/H)\Fc$.

Let $\hat{M}=\sharp\{i\in\{1,2,\ldots,M\}|F(i)\subset|\hat{\Phi}|\}$, and $\tau:\{1,2,\ldots,\hat{M}\}\rightarrow \{i\in\{1,2,\ldots,M\}$ denote the unique injective mapping preserving the order and satisfying $\tau(\{1,2,\ldots,\hat{M}\})=\{i\in\{1,2,\ldots,M\}|F(i)\subset|\hat{\Phi}|\}$.

By Theorem~\ref{usd1}.6 we know that $(H, \hat{\Phi})\in\mathcal{SF}(V,N,S)$, $(\hat{M},F\tau)\in$\hfill\break$\mathcal{USD}(H, \hat{\Phi},S)$, and
$\widetilde{\Phi}\backslash|\hat{\Phi}|=\hat{\Phi}*F\tau(1)*F\tau(2)*\cdots*F\tau(\hat{M})$.
\begin{enumerate}
\item
$\{\Gamma\in(\widetilde{\Phi}\backslash|\hat{\Phi}|)_1|
\Gamma\not\subset|\hat{\Phi}-(\hat{\Phi}/H)|\}
=\widetilde{\Phi}_1^\circ\backslash|\hat{\Phi}|$.
\end{enumerate}
We denote
$$\hat{I}= I(V,N, H,\hat{\Phi},S, \hat{M},F\tau): \widetilde{\Phi}_1^\circ\backslash|\hat{\Phi}|
\rightarrow\{1,2,\ldots, \hat{M}+1\},$$
$$\hat{\Psi}=\Psi(V,N, H,\hat{\Phi},S, \hat{M},F\tau): \widetilde{\Phi}_1^\circ\backslash|\hat{\Phi}|
\rightarrow 2^{2^{V^*}},$$
$$\smash{\hat{\Psi}}^\circ=\Psi^\circ(V,N, H,\hat{\Phi},S, \hat{M},F\tau): \widetilde{\Phi}_1^\circ\backslash|\hat{\Phi}|
\rightarrow 2^{2^{V^*}}.$$
\begin{enumerate}
\setcounter{enumii}{1}
\item
Let $\kappa: \{1,2,\ldots, \hat{M}+1\}\rightarrow\{1,2,\ldots, M+1\}$ denote the composition mapping $I\iota\hat{I}^{-1}$, where $\iota: \widetilde{\Phi}_1^\circ\backslash|\hat{\Phi}|\rightarrow\widetilde{\Phi}_1^\circ$ denotes the inclusion mapping.

The mapping $\kappa$ is injective and preserves the order.
$\kappa(\hat{M}+1)=M+1$.
\item
For any $\Gamma\in\widetilde{\Phi}_1^\circ\backslash|\hat{\Phi}|$,
$\hat{\Psi}(\Gamma)/\Gamma= (\Psi(\Gamma)/\Gamma)\backslash |\hat{\Phi}|$ and
$|\hat{\Psi}(\Gamma)/\Gamma|^\circ= |\Psi(\Gamma)/\Gamma|^\circ\cap |\hat{\Phi}|$.

For any $\Gamma\in\widetilde{\Phi}_1^\circ-
(\widetilde{\Phi}_1^\circ\backslash|\hat{\Phi}|)$,
$(\Psi(\Gamma)/\Gamma)\backslash |\hat{\Phi}|=\emptyset$ and
$|\Psi(\Gamma)/\Gamma|^\circ\cap |\hat{\Phi}|=\emptyset$.
\item
For any $\Gamma\in\widetilde{\Phi}_1^\circ\backslash|\hat{\Phi}|$,
$\smash{\hat{\Psi}}^\circ (\Gamma)= \Psi^\circ(\Gamma) \backslash |\hat{\Phi}|$ and
$|\smash{\hat{\Psi}}^\circ(\Gamma) |^\circ = |\Psi^\circ(\Gamma) |^\circ \cap |\hat{\Phi}|$.

For any $\Gamma\in\widetilde{\Phi}_1^\circ-
(\widetilde{\Phi}_1^\circ\backslash|\hat{\Phi}|)$,
$\Psi^\circ(\Gamma) \backslash |\hat{\Phi}|=\emptyset$ and
$|\Psi^\circ(\Gamma)|^\circ \cap |\hat{\Phi}|=\emptyset$.
\end{enumerate}
\item
Consider any $\Gamma\in \widetilde{\Phi}_1^\circ$.

If $\Theta\in\Psi(\Gamma)/\Gamma$, $\Delta\in\Phi$ and $\Theta \subset\Delta$, then $H\subset\Delta$, and $\Gamma\subset\Delta$.

If $\Theta\in\Psi^\circ(\Gamma)$, $\Delta\in\Phi/H$ and $\Theta \subset\Delta$, then $\Gamma\subset\Delta$.

If $\Theta\in\Psi^\circ(\Gamma) $ and $\Theta\subset|\Phi-(\Phi/H)|$, then
$\Gamma\not\subset\Theta$, $\Theta+H\in\Phi/H$ and
$b_{\Gamma/N^*}-b_{H/N^*}\in N^*\cap\Vect(\Theta)$.

If $\Theta\in\Psi^\circ(\Gamma) $,
then $\Theta+\Gamma\in\Psi^\circ(\Gamma)/\Gamma$.
\end{enumerate}

Below, we consider the case $\Ht(H,\Phi,S)>0$.
Assume $\Ht(H,\Phi,S)>0$.

The characteristic function $\gamma:(\Phi-(\Phi/H))_1\rightarrow\Q_0$ of $(H,\Phi, S)$ is defined.
Let $m=\sum_{\bar{E}\in(\Phi-(\Phi/H))_1}\lceil\gamma(\bar{E})\rceil\in\Z_+$ and $\bar{m}=\sum_{\bar{E}\in(\Phi-(\Phi/H))_1}\lfloor\gamma(\bar{E})\rfloor\in\Z_0$.
We know $\bar{m}\leq m\leq M$.

We take the unique pair $(E,M)$ of a compatible mapping
$E:\{1,2,\ldots,m\}\rightarrow (\Phi-(\Phi/H))_1$
with $S$ and a mapping 
$M:\{\bar{m},\bar{m}+1,\ldots,m+1\}\rightarrow\Z_+$
satisfying the following three conditions.
We denote 
$$F_E=F(V^*,N^*,H,\Phi,m,E):\{1,2,\ldots,m\}\rightarrow 2^{ V^*}$$ 
$$H_E=H(V^*,N^*,H,\Phi,m,E):\{1,2,\ldots,m+1\}\rightarrow 2^{ V^*},\text{ and}$$
$$\Omega_E=\Omega(V^*,N^*,H,\Phi,m,E):\{1,2,\ldots,m+1\}\rightarrow 2^{2^{ V^*}}.$$
\begin{description}
\item[\emph{(a)}]
$F(j)=F_E(j)$ for any $j\in\{1,2,\ldots, m\}$.
\item[\emph{(b)}]
$M(\bar{m})=m$, $M(m+1)=M$ and $M(i-1)\leq M(i)$ for any $i\in\{\bar{m}+1,\bar{m}+2,\ldots, m+1\}$.
\item[\emph{(c)}]
$F(j)\subset|\Omega_E(i)|$ and $F(j)\not\subset|\Omega_E(i)-(\Omega_E(i)/H_E(i))|$ for any $i\in\{\bar{m}+1,\bar{m}+2,\ldots, m+1\}$ and any $j\in\{M(i-1)+1,M(i-1)+2,\ldots, M(i)\}$.
\end{description}

For any $i\in\{0,1,\ldots, \bar{m}-1\}$, we put $M(i)=m$.
We obtain an extension $M:\{0,1,\ldots,m+1\}\rightarrow\Z_+$ of $M:\{\bar{m},\bar{m}+1,\ldots,m+1\}\rightarrow\Z_+$.
For any $i\in\{1,2,\ldots,m+1\}$, we denote
$\bar{\Phi}(i)=(\Phi*F(1)*F(2)*\cdots*F(M(i-1)))\backslash|\Omega_E(i)|$, and we take the mapping $\bar{F}(i):\{1,2,\ldots, M(i)-M(i-1)\}\rightarrow 2^{V^*}$ satisfying $\bar{F}(i)(j)=F(M(i-1)+j)$ for any $j\in\{1,2,\ldots, M(i)-M(i-1)\}$.
We know $(H_E(i), \bar{\Phi}(i))\in\mathcal{SF}(V,N,S)$,
$|\bar{\Phi}(i)|=|\Omega_E(i)|$,
$|\bar{\Phi}(i)-( \bar{\Phi}(i)/H_E(i))|=
|\Omega_E(i)-( \Omega_E(i)/H_E(i))|$,
$\Ht(H_E(i), \bar{\Phi}(i),S)<\Ht(H,\Phi,S)$, $(M(i)-M(i-1),\bar{F}(i))\in\mathcal{USD}( H_E(i), \bar{\Phi}(i),S)$, and
$\bar{\Phi}(i)* \bar{F}(1)* \bar{F}(2)*\cdots*\bar{F}(M(i)-M(i-1))=
\widetilde{\Phi}\backslash|\Omega_E(i)|$
for any $i\in\{1,2,\ldots,m+1\}$.
We denote
$$G_E=G(V^*,N^*,H,\Phi,m,E):\{1,2,\ldots,m\}\rightarrow 2^{V^*}.$$

For any $i\in\{0,1,\ldots,m+1\}$, we put $L(i)=i+M(i)-m$.
We obtain a mapping $L: \{0,1,\ldots,m+1\}\rightarrow\Z_0$.

\begin{enumerate}
\setcounter{enumi}{12}
\item
$L(0)=0$. $L(m+1)=M+1$. $L(i-1)<L(i)$ and $L(i)-L(i-1)=M(i)-M(i-1)+1$ for any $i\in\{1,2,\ldots,m+1\}$.
$L(i)=i$ for any $i\in\{0,1,\ldots,\bar{m}\}$.
\item
For any $i\in\{1,2,\ldots,m+1\}$, $H_E(i)\in\widetilde{\Phi}_1^\circ$ and $I(H_E(i))=L(i)$.
For any $i\in\{1,2,\ldots,\bar{m}\}$, $I(H_E(i))=i$.
\item
Consider any $\Gamma\in\widetilde{\Phi}_1^\circ$ and any $i\in\{1,2,\ldots,m+1\}$.
$\Gamma\subset|\Omega_E(i)|$ and $\Gamma\not\subset|\Omega_E(i)-( \Omega_E(i)/H_E(i))|
\Leftrightarrow L(i-1)<I(\Gamma)\leq L(i)
\Leftrightarrow |\Psi(\Gamma)|\subset|\Omega_E(i)|
\Leftrightarrow |\Psi^\circ(\Gamma)|^\circ\subset|\Omega_E(i)|$.
\item
If $\Gamma\in\widetilde{\Phi}_1^\circ$, $i\in\{1,2,\ldots,m+1\}$,
$\Gamma\subset|\Omega_E(i)|$ and $\Gamma\not\subset|\Omega_E(i)-( \Omega_E(i)/H_E(i))|$,
then 
$I(\Gamma)=L(i-1)+ I(V,N, H_E(i),\bar{\Phi}(i),S, M(i)-M(i-1),\bar{F}(i))(\Gamma)$, and
$\Psi(\Gamma)= \Psi(V,N, H_E(i),\bar{\Phi}(i),S, M(i)-M(i-1),\bar{F}(i))(\Gamma)$.
\item
If $\Gamma\in\widetilde{\Phi}_1^\circ$, $i\in\{1,2,\ldots,m\}$,
$\Gamma\subset|\Omega_E(i)|$ and $\Gamma\not\subset|\Omega_E(i)-( \Omega_E(i)/H_E(i))|$,
then $\Psi^\circ(\Gamma)= \{\Theta\in\Psi^\circ(V,N, H_E(i),\bar{\Phi}(i),S, M(i)-M(i-1),\bar{F}(i))(\Gamma)|\Theta^\circ\subset|\Omega_E/G_E(i)|^\circ\}$.

If $\Gamma\in\widetilde{\Phi}_1^\circ$, $\Gamma\subset|\Omega_E(m+1)|$ and $\Gamma\not\subset|\Omega_E(m+1)-( \Omega_E(m+1)/H_E(m+1))|$,
then $\Psi^\circ(\Gamma)= \Psi^\circ(V,N, H_E(m+1),\bar{\Phi}(m+1),S, M(m+1)-M(m),\bar{F}(m+1))(\Gamma)$.
\item
For any  $i\in\{1,2,\ldots,m\}$,
$\Psi^\circ(H_E(i))= \{\Theta\in\Psi(H_E(i))|\Theta^\circ\subset|\Omega_E/G_E(i)|^\circ\}$, and $\Psi (H_E(i))=((\widetilde{\Phi}\backslash|\Omega_E(i)|)/H_E(i))\Fc$.

$\Psi^\circ(H_E(m+1))=\Psi(H_E(m+1))=(\widetilde{\Phi}/H)\Fc$.
\item
Consider any $\Gamma\in \widetilde{\Phi}_1^\circ$ and any $i\in\{1,2,\ldots,m+1\}$ satisfying $\Gamma\subset|\Omega_E(i)|$ and $\Gamma\not\subset|\Omega_E(i)-(\Omega_E(i)/H_E(i))|$.

If $\Theta\in\Psi(\Gamma)/\Gamma$, $\Delta\in\Phi$
and $\Theta \subset\Delta$, then $H_E(i)\subset\Delta$.

If $\Theta\in\Psi^\circ(\Gamma)$, $\Delta\in\Phi/H$
and $\Theta \subset\Delta$, then $H_E(i)\subset\Delta$.
\end{enumerate}
\end{lemma}

\begin{theorem}
\label{important}
Assume $\mathcal{SF}(V,N,S)\neq\emptyset$.
Consider any $(H,\Phi)\in\mathcal{SF}(V,N,S)$ and any $(M,F)\in\mathcal{USD}(H,\Phi,S)$.
We denote $\widetilde{\Phi}=\Phi*F(1)*F(2)*\cdots*F(M)$,
$\widetilde{\Phi}_1^\circ=\{\Gamma\in\widetilde{\Phi}_1|\Gamma\not\subset|\Phi-(\Phi/H)|\}$ and
$$\Psi^\circ=\Psi^\circ(V,N, H,\Phi,S, M,F): \widetilde{\Phi}_1^\circ
\rightarrow 2^{2^{V^*}}.$$

Consider any $\Theta\in\widetilde{\Phi}$ with $\Theta\not\subset|\Phi-(\Phi/H)|$.
\begin{enumerate}
\item
There exists uniquely an element $\Lambda\in\Sigma(S|V)\hat{\cap}\Phi$ satisfying $\Theta^\circ\subset\Lambda^\circ$.
\item
There exists uniquely an element $\Delta\in\Phi$ satisfying $\Theta^\circ\subset\Delta^\circ$.
\item
There exists uniquely an element $A\in\mathcal{F}(S+(\Theta^\vee|V^*))$ satisfying
$\Delta(A, S+(\Theta^\vee|V^*)|V)=\Theta$.
\item
There exists uniquely an element $\Gamma\in\widetilde{\Phi}_1^\circ$ satisfying
$\Theta\in\Psi^\circ(\Gamma)$.
\end{enumerate}

We take the unique element $\Lambda\in\Sigma(S|V)\hat{\cap}\Phi$ satisfying $\Theta^\circ\subset\Lambda^\circ$, the unique element $\Delta\in\Phi$ satisfying $\Theta^\circ\subset\Delta^\circ$,
the unique element $A\in\mathcal{F}(S+(\Theta^\vee|V^*))$ satisfying
$\Delta(A, S+(\Theta^\vee|V^*)|V)=\Theta$ and the unique element $\Gamma\in\widetilde{\Phi}_1^\circ$ satisfying
$\Theta\in\Psi^\circ(\Gamma)$.
\begin{enumerate}
\setcounter{enumi}{4}
\item
$\Theta^\circ\subset\Lambda^\circ\subset\Delta^\circ$.
$\Delta\in\Phi/H$.
$\dim\Lambda=\dim\Delta$ or $\dim\Lambda=\dim\Delta-1$.
$\Gamma\subset\Delta$.
\item
$\dim A=\dim V-\dim\Theta$.
$S+(\Delta^\vee|V^*)\subset S+(\Theta^\vee|V^*)$.
$A\cap(S+(\Delta^\vee|V^*))\in\mathcal{F}(S+(\Delta^\vee|V^*))$.
$\Delta(A\cap(S+(\Delta^\vee|V^*)), S+(\Delta^\vee|V^*)|V)= \Lambda$.
\item
If $\dim\Lambda=\dim\Delta$, then $\langle \omega,a\rangle=\langle\omega,b\rangle$ for any $\omega\in\Vect(\Delta)$, 
any $a\in A\cap(S+(\Delta^\vee|V^*))$, and any $b\in A\cap(S+(\Delta^\vee|V^*))$.
\item
If $\dim\Lambda=\dim\Delta-1$, then $\Ht(H, S+(\Delta^\vee|V^*))>0$ and $\Gamma\in\widetilde{\Phi}_1$ satisfies the following five conditions:
\begin{enumerate}
\item
$\Gamma\subset\Delta$.
$\Gamma\not\subset\Vect(\Lambda)$.
$\Gamma\not\subset\Lambda$.
$\Gamma\not\subset\Theta$.
$\Theta+\Gamma\in \Psi^\circ(\Gamma)/\Gamma\subset\widetilde{\Phi}$.
\item
$\Vect(\Lambda)+\Gamma=\Vect(\Lambda)+H$.
\item
The subset $\{\langle b_{\Gamma/N^*},a\rangle|a\in A\cap(S+(\Delta^\vee|V^*))\}$ of $\R$ is a non-empty bounded closed interval.
\item
\emph{[The hard height inequality]}
\begin{equation*}\begin{split}
&\max\{\langle b_{\Gamma/N^*},a\rangle|a\in A\cap(S+(\Delta^\vee|V^*))\}\\
&\qquad\qquad-\min\{\langle b_{\Gamma/N^*},a\rangle|a\in A\cap(S+(\Delta^\vee|V^*))\}\\
\leq\:&\Ht(H, S+(\Delta^\vee|V^*)).
\end{split}\end{equation*}
\item
The equality 
\begin{equation*}\begin{split}
&\max\{\langle b_{\Gamma/N^*},a\rangle|a\in A\cap(S+(\Delta^\vee|V^*))\}\\
&\qquad\qquad-\min\{\langle b_{\Gamma/N^*},a\rangle|a\in A\cap(S+(\Delta^\vee|V^*))\}\\
=\:&\Ht(H, S+(\Delta^\vee|V^*)),
\end{split}\end{equation*}
holds, if and only if, $c(S+(\Delta^\vee|V^*))=2$ and
the structure constant of $\Sigma(S|V)\hat{\cap}\mathcal{F}(\Delta)$ corresponding to the pair $(2,\bar{E})$ is an integer for any $\bar{E}\in\mathcal{F}(\Delta)_1-\{H\}$.
\end{enumerate}

\item
If $\dim\Lambda=\dim\Delta-1$ and the equivalent conditions in $8$.\emph{(e)} are satisfied, then $\Theta=\Lambda$ and $\Gamma=H$.
\end{enumerate}

Below we consider any rational convex pseudo polytopes $T$ and $U$ over $N$ in $V$ satisfying $T+U=S$.
\begin{enumerate}\setcounter{enumi}{9}\item
$\Sigma(T|V)\hat{\cap}\Sigma(U|V)= \Sigma(S|V)$.
$(H,\Phi)\in\mathcal{SF}(V,N,T)$.
$(H,\Phi)\in\mathcal{SF}(V,N,U)$.
\item
If $\dim\Lambda=\dim\Delta-1$, then $\Lambda\in\Sigma(T+(\Delta^\vee|V^*)|V)^1$ or $\Lambda\in\Sigma(U+(\Delta^\vee|V^*)|V)^1$.
\item
Assume $\dim\Lambda=\dim\Delta-1$ and $\Lambda\in\Sigma(T+(\Delta^\vee|V^*)|V)^1$.

$\Ht(H, T+(\Delta^\vee|V^*))>0$. There exists uniquely an element $A_T\in\mathcal{F}(T+(\Delta^\vee|V^*))$ satisfying $\Delta(A_T, T+(\Delta^\vee|V^*)|V)=\Lambda$.

We take the unique element $A_T\in\mathcal{F}(T+(\Delta^\vee|V^*))$ satisfying $\Delta(A_T, T+(\Delta^\vee|V^*)|V)=\Lambda$. The element $\Gamma\in\widetilde{\Phi}_1^\circ$ satisfies the following three conditions:
\begin{enumerate}
\item
The subset $\{\langle b_{\Gamma/N^*},a\rangle|a\in A_T\}$ of $\R$ is a non-empty bounded closed interval.
\item
\emph{[The hard height inequality]}
\begin{equation*}\begin{split}
&\max\{\langle b_{\Gamma/N^*},a\rangle|a\in A_T\}\\
&\qquad\qquad-\min\{\langle b_{\Gamma/N^*},a\rangle|a\in A_T\}\\
\leq\:&\Ht(H, T+(\Delta^\vee|V^*)).
\end{split}\end{equation*}
\item
Let $r_U=c(U+(\Delta^\vee|V^*))\in\Z_+$.
If $\max\{\langle b_{\Gamma/N^*},a\rangle|a\in A_T\} 
-\min\{\langle b_{\Gamma/N^*},$\break$a\rangle|a\in A_T\}
=\Ht(H, T+(\Delta^\vee|V^*))$ and the structure constant of $\Sigma(U+(\Delta^\vee|V^*)|V)$ corresponding to the pair $(i,\bar{E})$ is an integer for any $i\in\{1,2,\ldots,r_U\}$ and any $\bar{E}\in\mathcal{F}(H\Op|\Delta)_1$, then
$c(T+(\Delta^\vee|V^*))=2$, the structure constant of $\Sigma(T+(\Delta^\vee|V^*)|V)$ corresponding to the pair $(2,\bar{E})$ is an integer for any $\bar{E}\in\mathcal{F}(H\Op|\Delta)_1$ and $\Theta=\Lambda$.
\end{enumerate}
\end{enumerate}
\end{theorem}

\begin{proof}

We show only claim $8$ and $9$.

Assume $\mathcal{SF}(V,N,S)\neq\emptyset$.
Consider any $(H,\Phi)\in\mathcal{SF}(V,N,S)$ and any $(M,F)\in\mathcal{USD}(H,\Phi,S)$.
We denote $\widetilde{\Phi}=\Phi*F(1)*F(2)*\cdots*F(M)$,
$\widetilde{\Phi}_1^\circ=\{\Gamma\in\widetilde{\Phi}_1|\Gamma\not\subset|\Phi-(\Phi/H)|\}$ and
$\Psi^\circ=\Psi^\circ(V,N, H,\Phi,S, M,F): \widetilde{\Phi}_1^\circ
\rightarrow 2^{2^{V^*}}$.

Consider any $\Theta\in\widetilde{\Phi}$ satisfying $\Theta\not\subset|\Phi-(\Phi/H)|$.

We take the unique element $\Lambda\in\Sigma(S|V)\hat{\cap}\Phi$ satisfying $\Theta^\circ\subset\Lambda^\circ$, the unique element $\Delta\in\Phi$ satisfying $\Theta^\circ\subset\Delta^\circ$,
the unique element $A\in\mathcal{F}(S+(\Theta^\vee|V^*))$ satisfying
$\Delta(A, S+(\Theta^\vee|V^*)|V)=\Theta$ and the unique element $\Gamma\in\widetilde{\Phi}_1^\circ$ satisfying
$\Theta\in\Psi^\circ(\Gamma)$.

By $5$ and $6$ we know that 
$\Theta^\circ\subset\Lambda^\circ\subset\Delta^\circ$,
$\Delta\in\Phi/H$,
$\dim\Lambda=\dim\Delta$ or $\dim\Lambda=\dim\Delta-1$,
$\Gamma\subset\Delta$,
$\dim A=\dim V-\dim\Theta$,
$S+(\Delta^\vee|V^*)\subset S+(\Theta^\vee|V^*)$,
$A\cap(S+(\Delta^\vee|V^*))\in\mathcal{F}(S+(\Delta^\vee|V^*))$, and
$\Delta(A\cap(S+(\Delta^\vee|V^*)), S+(\Delta^\vee|V^*)|V)= \Lambda$.

Furthermore, assume $\dim\Lambda=\dim\Delta-1$.

Note that $H\in\mathcal{F}(\Delta)_1$ and $\Sigma( S+(\Delta^\vee|V^*)|V)$ is $H$-simple, since $\Delta\in\Phi/H$.
Let $\bar{\Sigma}(S+(\Delta^\vee|V^*)|V)^1=\{\bar{\Lambda}\in\Sigma( S+(\Delta^\vee|V^*)|V)^1|\bar{\Lambda}^\circ\subset\Delta^\circ\}\cup\{H\Op|\Delta\}$ denote the $H$-skeleton of $\Sigma( S+(\Delta^\vee|V^*)|V)$.
We know $\Lambda\in\bar{\Sigma}(S+(\Delta^\vee|V^*)|V)^1$, $\Lambda\neq H\Op|\Delta$, $c(S+(\Delta^\vee|V^*))=\sharp \bar{\Sigma}(S+(\Delta^\vee|V^*)|V)^1\geq 2$, and $\Ht(H, S+(\Delta^\vee|V^*))>0$.

We consider the $H$-order on $\Sigma( S+(\Delta^\vee|V^*)|V)^0$.
Let $\hat{\Lambda}:\{1,2,\ldots, c(S+(\Delta^\vee|V^*))\}\rightarrow\Sigma( S+(\Delta^\vee|V^*)|V)^0$ denote the unique bijective mapping preserving the $H$-order.
Let $\ell=\dim\Vect(\Delta)^\vee|V^*\in\Z_0$ and let $\hat{A}:\{1,2,\ldots, c(S+(\Delta^\vee|V^*))\}\rightarrow\mathcal{F}(S+(\Delta^\vee))_\ell$ denote the unique bijective mapping satisfying $\Delta(\hat{A}(i), $\hfill\break$S+(\Delta^\vee|V^*)|V)= \hat{\Lambda}(i)$ for any $i\in\{1,2,\ldots, c(S+(\Delta^\vee|V^*))\}$.
For any $i\in\{1,2,\ldots, c(S+(\Delta^\vee|V^*))\}$, we take any point $\hat{a}(i)\in \hat{A}(i)$.

We consider the $H$-order on $\bar{\Sigma}(S+(\Delta^\vee|V^*)|V)^1$.
Let $\hat{\bar{\Lambda}}:\{1,2,\ldots, c(S+(\Delta^\vee|V^*))\}\rightarrow\bar{\Sigma}(S+(\Delta^\vee|V^*)|V)^1$ denote the unique bijective mapping preserving the $H$-order.
We take the unique element $i_\Lambda\in\{2,3,\ldots, c(S+(\Delta^\vee|V^*))\}$ satisfying $\hat{\bar{\Lambda}}( i_\Lambda)=\Lambda$.

Now, by $5$ we know $\Gamma\subset\Delta$. 

By Lemma~\ref{lower parts}.12 
we know $\Theta+\Gamma\in\Psi^\circ(\Gamma)/\Gamma \subset\Psi^\circ(\Gamma)\subset\widetilde{\Phi}$.

Since $\Delta(A\cap(S+(\Delta^\vee|V^*)), S+(\Delta^\vee|V^*)|V)=\Lambda=\hat{\bar{\Lambda}}( i_\Lambda)=
\hat{\Lambda}(i_\Lambda)\cap\hat{\Lambda}(i_\Lambda-1)$ and $\Lambda^\circ\subset\Delta^\circ$,
we know that
$A\cap(S+(\Delta^\vee|V^*))=\Conv(\hat{A}(i_\Lambda)\cup\hat{A}(i_\Lambda-1))=\Conv(\{\hat{a}(i_\Lambda), \hat{a}(i_\Lambda-1)\})+(\Vect(\Delta)^\vee|V^*)$.
Therefore, if $\Vect(\Lambda)+\Gamma=\Vect(\Lambda)+H$, then
$\langle b_{\Gamma/N^*}, \hat{a}(i_\Lambda)\rangle <\langle b_{\Gamma/N^*}, \hat{a}(i_\Lambda-1)\rangle$,
$\{\langle b_{\Gamma/N^*}, a\rangle|a\in A\cap(S+(\Delta^\vee|V^*))\}
=\{t\in\R|\langle b_{\Gamma/N^*}, \hat{a}(i_\Lambda)\rangle \leq t\leq\langle b_{\Gamma/N^*}, \hat{a}(i_\Lambda-1)\rangle\}$, and the subset
$\{\langle b_{\Gamma/N^*}, a\rangle|a\in A\cap(S+(\Delta^\vee|V^*))$ of $\R$ is a non-empty bounded closed interval.
If $\Gamma\subset\Vect(\Lambda)$, then 
$\langle b_{\Gamma/N^*}, \hat{a}(i_\Lambda)\rangle =\langle b_{\Gamma/N^*}, \hat{a}(i_\Lambda-1)\rangle$,
$\{\langle b_{\Gamma/N^*}, a\rangle|a\in A\cap(S+(\Delta^\vee|V^*))\}
=\{\langle b_{\Gamma/N^*}, \hat{a}(i_\Lambda)\rangle\}$, and the subset
$\{\langle b_{\Gamma/N^*}, a\rangle|a\in A\cap(S+(\Delta^\vee|V^*))$ of $\R$ is a non-empty bounded closed interval.
If $\Gamma\not\subset\Vect(\Lambda)$ and $\Vect(\Lambda)+\Gamma\neq\Vect(\Lambda)+H$, then 
$\langle b_{\Gamma/N^*}, \hat{a}(i_\Lambda)\rangle >\langle b_{\Gamma/N^*}, \hat{a}(i_\Lambda-1)\rangle$,
$\{\langle b_{\Gamma/N^*}, a\rangle|a\in A\cap(S+(\Delta^\vee|V^*))\}
=\{t\in\R|\langle b_{\Gamma/N^*}, \hat{a}(i_\Lambda)\rangle \geq t\geq\langle b_{\Gamma/N^*}, \hat{a}(i_\Lambda-1)\rangle\}$, and the subset
$\{\langle b_{\Gamma/N^*}, a\rangle|a\in A\cap(S+(\Delta^\vee|V^*))$ of $\R$ is a non-empty bounded closed interval.

We know that the subset
$\{\langle b_{\Gamma/N^*}, a\rangle|a\in A\cap(S+(\Delta^\vee|V^*))$ of $\R$ is a non-empty bounded closed interval.

Let $\hat{M}=\sharp\{i\in\{1,2,\ldots,M\}|F(i)\subset\Delta\}\in\Z_0$.
We take the unique injective mapping $\tau:\{1,2,\ldots, \hat{M}\}\rightarrow\{1,2,\ldots,M\}$ preserving the order and satisfying $\tau(\{1,2,\ldots,\hat{M}\})=\{i\in\{1,2,\ldots,M\}|F(i)\subset\Delta\}$.
By Therem~\ref{usd1}.6 we know that $(H, \mathcal{F}(\Delta))\in \mathcal{SF}(V, N, S)$, 
$(\hat{M}, F\tau)\in\mathcal{USD}(H,\mathcal{F}(\Delta),S)$, and
$\widetilde{\Phi}\backslash\Delta=\mathcal{F}(\Delta)*F\tau(1)*F\tau(2)*\cdots*F\tau(\hat{M})$.

We denote $\widetilde{\mathcal{F}}(\Delta)=\mathcal{F}(\Delta)*F\tau(1)*F\tau(2)*\cdots*F\tau(\hat{M})$,
$\smash{\widetilde{\mathcal{F}}(\Delta)}_1^\circ=\{\bar{\Gamma}\in\widetilde{\mathcal{F}}(\Delta)_1|\bar{\Gamma}\not\subset H\Op|\Delta\}$ and
$\smash{\hat{\Psi}}^\circ=\Psi^\circ(V,N, H, \mathcal{F}(\Delta),S, \hat{M},F\tau): \smash{\widetilde{\mathcal{F}}(\Delta)}_1^\circ
\rightarrow 2^{2^{V^*}}$.

Since $\Gamma\subset\Delta$, $\Gamma\in\smash{\widetilde{\mathcal{F}}(\Delta)}_1^\circ$ and 
$\smash{\hat{\Psi}}^\circ(\Gamma)=\Psi^\circ(\Gamma)\backslash\Delta$.
Since $\Theta\subset\Delta$, $\Theta\in \smash{\hat{\Psi}}^\circ(\Gamma)$.

For any $j\in\{1,2,\ldots,c(S+(\Delta^\vee|V^*))\}$ and any $\bar{E}\in\mathcal{F}(H\Op|\Delta)_1$, we denote the structure constant of 
$\Sigma( S+(\Delta^\vee|V^*)|V)$ corresponding to the pair $(j, \bar{E})$ by
$c(\Sigma( S+(\Delta^\vee|V^*)|V),j, \bar{E})$.

Let $\hat{m}=\sum_{\bar{E}\in\mathcal{F}(H\Op|\Delta)_1}\lceil c(\Sigma( S+(\Delta^\vee|V^*)|V),2,\bar{E})\rceil\in\Z_+$ and\hfill\break $\hat{\bar{m}}=\sum_{\bar{E}\in\mathcal{F}(H\Op|\Delta)_1}\lfloor c(\Sigma( S+(\Delta^\vee|V^*)|V),2,\bar{E})\rfloor\in\Z_0$.
We know that $\hat{\bar{m}}\leq\hat{m}\leq\hat{M}$.

We take the unique pair $(\hat{E},\hat{M})$ of a compatible mapping
$\hat{E}:\{1,2,\ldots,\hat{m}\}\rightarrow \mathcal{F}(H\Op|\Delta)_1$
with $S$ and a mapping 
$\hat{M}:\{\hat{\bar{m}},\hat{\bar{m}}+1,\ldots,\hat{m}+1\}\rightarrow\Z_+$
satisfying the following three conditions.
We denote 
$$F_{\hat{E}}=F(V^*,N^*,H,\mathcal{F}(\Delta),\hat{m},\hat{E}):\{1,2,\ldots,\hat{m}\}\rightarrow 2^{ V^*}$$ 
$$H_{\hat{E}}=H(V^*,N^*,H,\mathcal{F}(\Delta),\hat{m},\hat{E}):\{1,2,\ldots,\hat{m}+1\}\rightarrow 2^{ V^*},\text{ and}$$
$$\Omega_{\hat{E}}=\Omega(V^*,N^*,H,\mathcal{F}(\Delta),\hat{m},\hat{E}):\{1,2,\ldots,\hat{m+1}\}\rightarrow 2^{2^{ V^*}}.$$
\begin{description}
\item[(a)]
$F\tau(j)=F_{\hat{E}}(j)$ for any $j\in\{1,2,\ldots, \hat{m}\}$.
\item[(b)]
$\hat{M}(\hat{\bar{m}})=\hat{m}$, $\hat{M}(\hat{m}+1)=\hat{M}$ and $\hat{M}(i-1)\leq \hat{M}(i)$ for any $i\in\{\hat{\bar{m}}+1,\hat{\bar{m}}+2,\ldots, \hat{m}+1\}$.
\item[(c)]
$F\tau(j)\subset|\Omega_{\hat{E}}(i)|$ and $F\tau(j)\not\subset|\Omega_{\hat{E}}(i)-(\Omega_{\hat{E}}(i)/H_{\hat{E}}(i))|$ for any $i\in\{\hat{\bar{m}}+1,\hat{\bar{m}}+2,\ldots, \hat{m}+1\}$ and any $j\in\{\hat{M}(i-1)+1,\hat{M}(i-1)+2,\ldots, \hat{M}(i)\}$.
\end{description}

For any $i\in\{1,2,\ldots,\hat{m}+1\}$, we denote
$\hat{\Delta}(i)=| \Omega_{\hat{E}}(i)|\subset\Delta$ and 
$\hat{\bar{\Delta}}(i)= H_{\hat{E}}(i)\Op|\hat{\Delta}(i)\subset\hat{\Delta}(i)$.
For any $i\in\{0,1,\ldots, \hat{\bar{m}}-1\}$ we put $\hat{M}(i)=\hat{m}$.
For any $i\in\{1,2,\ldots,\hat{m}+1\}$, 
let $\bar{\Phi}(i)=\mathcal{F}(\Delta)*F\tau(1)*F\tau(2)*\cdots*F\tau(\hat{M}(i-1))\backslash \hat{\Delta}(i)$ and 
let $\bar{F}(i):\{1,2,\ldots,\hat{M}(i)-\hat{M}(i-1)\}\rightarrow 2^{V^*}$ denote the mapping satisfying $\bar{F}(i)(j)=F\tau(\hat{M}(i-1)+j)$ for any $j\in \{1,2,\ldots,\hat{M}(i)-\hat{M}(i-1)\}$.
Note that we have 
$|\bar{\Phi}(i)|= \hat{\Delta}(i)$,
$(H_{\hat{E}}(i), \bar{\Phi}(i))\in\mathcal{SF}(V,N,S)$,
$(\hat{M}(i)-\hat{M}(i-1), \bar{F}(i))\in\mathcal{USD}(H_{\hat{E}}(i), \bar{\Phi}(i),S)$,
$\Ht(H_{\hat{E}}(i), \bar{\Phi}(i),S)=\Ht(H_{\hat{E}}(i), S+(\hat{\Delta}(i)^\vee|V^*))<\Ht(H, S+(\Delta^\vee|V^*))$ for any 
$i\in\{1,2,\ldots,\hat{m}+1\}$ by Theorem~\ref{usd1}.10.

There exists uniquely $i_\Gamma\in\{1,2,\ldots,\hat{m}+1\}$ satisfying
$\Gamma\subset\hat{\Delta}(i_\Gamma)$ and $\Gamma\not\subset \hat{\bar{\Delta}}(i_\Gamma)$.
We take the unique  $i_\Gamma\in\{1,2,\ldots,\hat{m}+1\}$ satisfying
$\Gamma\subset\hat{\Delta}(i_\Gamma)$ and $\Gamma\not\subset \hat{\bar{\Delta}}(i_\Gamma)$.

By Lemma~\ref{lower parts}.17 we know 
$\Theta\in\smash{\hat{\Psi}}^\circ(\Gamma)\subset
\Psi^\circ(V,N, H_{\hat{E}}(i_\Gamma), \bar{\Phi}(i_\Gamma), S, 
\hat{M}(i_\Gamma)- \hat{M}(i_\Gamma-i),\bar{F}(i_\Gamma))(\Gamma)$, and
$\Theta\subset|\smash{\hat{\Psi}}^\circ(\Gamma)|
\subset|\Psi^\circ(V,N, H_{\hat{E}}(i_\Gamma), \bar{\Phi}(i_\Gamma), S, 
\hat{M}(i_\Gamma)- \hat{M}(i_\Gamma-i),\bar{F}(i_\Gamma))(\Gamma)|\subset
|\bar{\Phi}(i_\Gamma)|=\hat{\Delta}(i_\Gamma)$.
Consider the case $i_\Gamma=\hat{m}+1$. $\Theta^\circ\subset \hat{\Delta}(i_\Gamma)\cap \Delta^\circ=
\hat{\Delta}(i_\Gamma)^\circ\cup\hat{\bar{\Delta}}(i_\Gamma)^\circ $.
Consider the case $i_\Gamma\neq\hat{m}+1$.
We denote  
$$G_{\hat{E}}=G(V^*,N^*,H,\mathcal{F}(\Delta),\hat{m},\hat{E}):\{1,2,\ldots,\hat{m}\}\rightarrow 2^{ V^*}.$$
By Lemma~\ref{lower parts}.17
$\Theta^\circ\subset|\mathcal{F}(\hat{\Delta}(i_\Gamma))/ G_{\hat{E}}(i_\Gamma)|^\circ=
\hat{\Delta}(i_\Gamma)- (G_{\hat{E}}(i_\Gamma)\Op|\Delta(i_\Gamma))$.
Therefore,
$\Theta^\circ\subset(\hat{\Delta}(i_\Gamma)- (G_{\hat{E}}(i_\Gamma)\Op|\Delta(i_\Gamma)))\cap \Delta^\circ=
\hat{\Delta}(i_\Gamma)^\circ\cup\hat{\bar{\Delta}}(i_\Gamma)^\circ $.

We know $\Theta^\circ\subset \hat{\Delta}(i_\Gamma)^\circ\cup\hat{\bar{\Delta}}(i_\Gamma)^\circ$.

Since $\Theta\in\widetilde{\mathcal{F}}(\Delta)$, $\{\hat{\Delta}(i_\Gamma), \hat{\bar{\Delta}}(i_\Gamma)\}\subset \mathcal{F}(\Delta)*F\tau(1)*F\tau(2)*\cdots*F\tau(\hat{m})$, and $\widetilde{\mathcal{F}}(\Delta)$ is a subdivision of $\mathcal{F}(\Delta)*F\tau(1)*F\tau(2)*\cdots*F\tau(\hat{m})$, we know that $\Theta^\circ\subset \hat{\Delta}(i_\Gamma)^\circ$ or $\Theta^\circ\subset \hat{\bar{\Delta}}(i_\Gamma)^\circ$.

We take the unique element $\Delta_0\in\bar{\Phi}(i_\Gamma)/ H_{\hat{E}}(i_\Gamma)$ satisfying $\Theta^\circ\subset \Delta_0^\circ\cup(H_{\hat{E}}(i_\Gamma)\Op$\hfill\break$|\Delta_0)^\circ$.
$\Theta^\circ\subset \hat{\Delta}(i_\Gamma)^\circ$, if and only if, $\Theta^\circ\subset \Delta_0^\circ$.
$\Theta^\circ\subset \hat{\bar{\Delta}}(i_\Gamma)^\circ$, if and only if,
$\Theta^\circ\subset(H_{\hat{E}}(i_\Gamma)\Op|\Delta_0)^\circ$.

Since $\Theta\in\Psi^\circ(V,N, H_{\hat{E}}(i_\Gamma), \bar{\Phi}(i_\Gamma), S, 
\hat{M}(i_\Gamma)- \hat{M}(i_\Gamma-i),\bar{F}(i_\Gamma))(\Gamma)$ and $\Theta\subset\Delta_0$, we know $\Gamma\subset\Delta_0$.
$\Ht(H_{\hat{E}}( i_\Gamma), S+(\Delta_0^\vee|V^*))\leq\Ht(H_{\hat{E}}( i_\Gamma), \bar{\Phi}( i_\Gamma),S)<\Ht(H,S+(\Delta^\vee|V^*))$.
$\Lambda\cap\Delta_0\in\Sigma(S|V)\hat{\cap}\Phi\hat{\cap}\bar{\Phi}(i_\Gamma)= \Sigma(S|V)\hat{\cap}\bar{\Phi}(i_\Gamma)$.
$\emptyset\neq\Theta^\circ\subset\Lambda\cap(\hat{\Delta}(i_\Gamma)^\circ\cup\hat{\bar{\Delta}}(i_\Gamma)^\circ)$.
$\Theta\subset\Lambda\cap\Delta_0$.

Let $\hat{s}=s(V^*,N^*,H,\mathcal{F}(\Delta),\hat{m},\hat{E}):\{0,1,\ldots,\hat{m}\}\times\mathcal{F}(H\Op|\Delta)_1\rightarrow\Z_0$.

Assume that $i_\Gamma\neq\hat{m}+1$.
By Theorem~\ref{heightinequality}.$15$ we know
$\hat{m}- \hat{\bar{m}}\geq 1$,
$i_\Gamma\in\{\hat{\bar{m}}+1,\hat{\bar{m}}+2,\ldots,\hat{m}\}$, since 
$\Lambda\cap(\hat{\Delta}(i_\Gamma)^\circ\cup\hat{\bar{\Delta}}(i_\Gamma)^\circ)\neq\emptyset$.
Since $\Lambda=\hat{\bar{\Lambda}}(i_\Lambda)$ and
$\Lambda\cap(\hat{\Delta}(i_\Gamma)^\circ\cup\hat{\bar{\Delta}}(i_\Gamma)^\circ)\neq\emptyset$, we know 
$c(\Sigma(S+(\Delta^\vee|V^*)|V),i_\Lambda,\hat{E}(i_\Gamma))<
\lceil c(\Sigma(S+(\Delta^\vee|V^*)|V),2,\hat{E}(i_\Gamma))\rceil=
\hat{s}(i_\Gamma,\hat{E}(i_\Gamma))$ by Theorem~\ref{heightinequality}.$17$.

It follows that if $i_\Gamma\neq\hat{m}+1$, then $H_{\hat{E}}(i_\Gamma)\not\subset\Vect(\Lambda)$ and $\Vect(\Lambda)+ H_{\hat{E}}(i_\Gamma)= \Vect(\Lambda)+H$.
Obviously, if $i_\Gamma=\hat{m}+1$, then $H_{\hat{E}}(i_\Gamma)=H$, $H_{\hat{E}}(i_\Gamma)\not\subset\Vect(\Lambda)$ and $\Vect(\Lambda)+ H_{\hat{E}}(i_\Gamma)= \Vect(\Lambda)+H$.
We know that $H_{\hat{E}}(i_\Gamma)\not\subset\Vect(\Lambda)$ and $\Vect(\Lambda)+ H_{\hat{E}}(i_\Gamma)= \Vect(\Lambda)+H$.
Furthermore, it follows that
$\dim \Lambda\cap\Delta_0\leq\dim \Delta_0-1$.

We have two cases.
\begin{enumerate}
\item
$\Theta^\circ\subset \hat{\Delta}(i_\Gamma)^\circ$.
\item
$\Theta^\circ\subset \hat{\bar{\Delta}}(i_\Gamma)^\circ$.
\end{enumerate}

We consider the case $\Theta^\circ\subset \hat{\Delta}(i_\Gamma)^\circ$.
$\Theta^\circ\subset \Delta_0^\circ$.
$\Lambda\cap\Delta_0\in\Sigma(S|V)\hat{\cap}\mathcal{F}(\Delta_0)$.
$\Theta^\circ\subset\Lambda^\circ\cap\Delta_0^\circ=(\Lambda\cap\Delta_0)^\circ$,
$\emptyset\neq\Theta^\circ\subset(\Lambda\cap\Delta_0)^\circ\cap\Delta_0^\circ$.
Since $\Sigma(S|V)\hat{\cap}\mathcal{F}(\Delta_0)$ is $H_{\hat{E}}(i_\Gamma)$-simple, we know $\dim (\Lambda\cap\Delta_0)\geq\dim\Delta_0-1$, $\dim (\Lambda\cap\Delta_0)=\dim\Delta_0-1$ and $\Vect(\Lambda\cap\Delta_0)=\Vect(\Lambda)\cap\Vect(\Delta_0)$.

Since $\Ht(H_{\hat{E}}(i_\Gamma), S+(\Delta_0^\vee|V^*))<\Ht(H, S+(\Delta^\vee|V^*)$, by induction on $\Ht$ we know that the following claims hold:
\begin{description}
\item[(a)] $\Gamma\subset\Delta_0$. $\Gamma\not\subset\Vect(\Lambda\cap\Delta_0)$.
\item[(b)]
$\Vect(\Lambda\cap\Delta_0)+\Gamma=\Vect(\Lambda\cap\Delta_0)+ H_{\hat{E}}(i_\Gamma)$.
\item[(c)]
The subset $\{\langle b_{\Gamma/N^*},a\rangle|a\in A\cap(S+(\Delta_0^\vee|V^*))\}$ of $\R$ is a non-empty bounded closed interval.
\item[(d)]
\begin{equation*}\begin{split}
&\max\{\langle b_{\Gamma/N^*},a\rangle|a\in A\cap(S+(\Delta_0^\vee|V^*))\}\\
&\qquad\qquad -\min\{\langle b_{\Gamma/N^*},a\rangle|a\in A\cap(S+(\Delta_0^\vee|V^*))\}\\
\leq\:&\Ht(H_{\hat{E}}(i_\Gamma), S+(\Delta_0^\vee|V^*))
\end{split}\end{equation*}
\end{description}

Since $\Gamma\subset\Delta_0\subset\Vect(\Delta_0)$ and $\Gamma\not\subset\Vect(\Lambda\cap\Delta_0)=\Vect(\Lambda)\cap\Vect(\Delta_0)$,
we know $\Gamma\not\subset\Vect(\Lambda)$.
Since $\Theta\subset\Lambda\subset\Vect(\Lambda)$, we know $\Gamma\not\subset\Lambda$ and $\Gamma\not\subset\Theta$.

Note that 
$\Lambda=\Vect(\Lambda)\cap\Delta$.
Since $ H_{\hat{E}}(i_\Gamma)\not\subset\Lambda$ and $ H_{\hat{E}}(i_\Gamma)\subset\Delta$, we know $ H_{\hat{E}}(i_\Gamma)\not\subset\Vect(\Lambda)$.
Since $\Vect(\Lambda\cap\Delta_0)+\Gamma=\Vect(\Lambda\cap\Delta_0)+ H_{\hat{E}}(i_\Gamma)$ and $\Vect(\Lambda\cap\Delta_0)\subset\Vect(\Lambda)$, we know
$\Vect(\Lambda)+\Gamma=\Vect(\Lambda)+ H_{\hat{E}}(i_\Gamma) =\Vect(\Lambda)+H$.

Since $\Delta_0\subset\Delta$, we know that
$A\cap(S+(\Delta^\vee|V^*)\subset A\cap(S+(\Delta_0^\vee|V^*)$,\hfill\break
$\max\{\langle b_{\Gamma/N^*},a\rangle|a\in A\cap(S+(\Delta^\vee|V^*))\}
\leq \max\{\langle b_{\Gamma/N^*},a\rangle|a\in A\cap(S+(\Delta_0^\vee|V^*))\}$, 
\hfill\break
$\min\{\langle b_{\Gamma/N^*},a\rangle|a\in A\cap(S+(\Delta^\vee|V^*))\}
\geq \min\{\langle b_{\Gamma/N^*},a\rangle|a\in A\cap(S+(\Delta_0^\vee|V^*))\}$,
\hfill\break
$\max\{\langle b_{\Gamma/N^*},a\rangle|a\in A\cap(S+(\Delta^\vee|V^*))\}
-\min\{\langle b_{\Gamma/N^*},a\rangle|a\in A\cap(S+(\Delta^\vee|V^*))\}
\leq\max\{\langle b_{\Gamma/N^*},a\rangle|a\in A\cap(S+(\Delta_0^\vee|V^*))\}
-\min\{\langle b_{\Gamma/N^*},a\rangle|a\in A\cap(S+(\Delta_0^\vee|V^*))\}
\leq\Ht(H_{\hat{E}}(i_\Gamma), S+(\Delta_0^\vee|V^*))<\Ht(H,S+(\Delta^\vee|V^*))$, and
\begin{equation*}\begin{split}
&\max\{\langle b_{\Gamma/N^*},a\rangle|a\in A\cap(S+(\Delta^\vee|V^*))\}\\
&\qquad\qquad -\min\{\langle b_{\Gamma/N^*},a\rangle|a\in A\cap(S+(\Delta^\vee|V^*))\}\\
<\:&\Ht(H,S+(\Delta^\vee|V^*))
\end{split}\end{equation*}

We consider the case $\Theta^\circ\subset \hat{\bar{\Delta}}(i_\Gamma)^\circ$.

$\Theta\in\widetilde{\Phi}\backslash\hat{\bar{\Delta}}(i_\Gamma)
=(\widetilde{\Phi}\backslash\hat{\Delta}(i_\Gamma)) \backslash\hat{\bar{\Delta}}(i_\Gamma)
=(\bar{\Phi}(i_\Gamma)*\bar{F}(i_\Gamma)(2)* \bar{F}(i_\Gamma)(1)*\cdots*\bar{F}(i_\Gamma)(\hat{M}(i_\Gamma)-\hat{M}(i_\Gamma-1))) \backslash\hat{\bar{\Delta}}(i_\Gamma)
=\bar{\Phi}(i_\Gamma)-( \bar{\Phi}(i_\Gamma)/ H_{\hat{E}}(i_\Gamma))$ and we know $\Theta+ H_{\hat{E}}(i_\Gamma)\in\bar{\Phi}(i_\Gamma)$.

Since $\Theta= H_{\hat{E}}(i_\Gamma)\Op|(\Theta+ H_{\hat{E}}(i_\Gamma))$, we know $\Delta_0=\Theta+ H_{\hat{E}}(i_\Gamma)$ and $\dim\Theta=\dim\Delta_0-1$.

Since $\Theta\subset\Lambda\cap\Delta_0$, $\dim\Delta_0-1=\dim\Theta\leq\dim\Lambda\cap\Delta_0\leq\dim\Delta_0-1$,
$\dim\Theta=\dim\Lambda\cap\Delta_0=\dim\Delta_0-1$,
$\Vect(\Theta)=\Vect(\Lambda\cap\Delta_0)$ and
$\Lambda\cap\Delta_0\subset\Vect(\Lambda\cap\Delta_0)\cap\Delta_0=
\Vect(\Theta)\cap\Delta_0=\Vect(\Theta)\cap(\Theta+ H_{\hat{E}}(i_\Gamma))=\Theta$.
We know $\Theta=\Lambda\cap\Delta_0$.

Since $\Gamma\not\subset\hat{\bar{\Delta}}(i_\Gamma)$ and $\Theta\subset\hat{\bar{\Delta}}(i_\Gamma)$, we know $\Gamma\not\subset\Theta$.
Since $\Gamma\not\subset\Theta=\Lambda\cap\Delta_0$ and $\Gamma\subset\Delta_0$, we know $\Gamma\not\subset\Lambda$.
Since $\Gamma\not\subset\Lambda=\Vect(\Lambda)\cap\Delta$ and $\Gamma\subset\Delta$, we know $\Gamma\not\subset\Vect(\Lambda)$.

Since $\Gamma^\circ\cup H_{\hat{E}}(i_\Gamma)^\circ\subset\Delta_0-(H_{\hat{E}}(i_\Gamma)\Op|\Delta_0)$,
we know $\Vect(\Lambda)+\Gamma=\Vect(\Lambda)+ H_{\hat{E}}(i_\Gamma) =\Vect(\Lambda)+ H$.

We know $\Vect(\Lambda)+\Gamma=\Vect(\Lambda)+ H$,
$\langle b_{\Gamma/N^*}, \hat{a}(i_\Lambda)\rangle <\langle b_{\Gamma/N^*}, \hat{a}(i_\Lambda-1)\rangle$,
$\{\langle b_{\Gamma/N^*}, a\rangle|a\in A\cap(S+(\Delta^\vee|V^*))\}
=\{t\in\R|\langle b_{\Gamma/N^*}, \hat{a}(i_\Lambda)\rangle \leq t\leq\langle b_{\Gamma/N^*}, \hat{a}(i_\Lambda-1)\rangle\}$, and
$\max\{\langle b_{\Gamma/N^*}, a\rangle|a\in A\cap(S+(\Delta^\vee|V^*))\}
-\min\{\langle b_{\Gamma/N^*}, a\rangle|a\in A\cap(S+(\Delta^\vee|V^*))\}
=\langle b_{\Gamma/N^*}, \hat{a}(i_\Lambda-1)\rangle-\langle b_{\Gamma/N^*}, \hat{a}(i_\Lambda)\rangle$.

By Lemma~\ref{lower parts}.12 we know $b_{ H_{\hat{E}}(i_\Gamma)/N^*}-b_{\Gamma/N^*}\in N^*\cap\Vect(\Theta)
\subset\Vect(\Lambda)=\Vect(\hat{\bar{\Lambda}}(i_\Lambda))$.
Since $\{\hat{a}(i_\Lambda-1), \hat{a}(i_\Lambda)\}\subset\Conv(\hat{A}(i_\Lambda-1)\cup\hat{A}(i_\Lambda))\in\mathcal{F}(S+(\Delta^\vee|V^*))_{\ell+1}$ and $\Delta(\Conv(\hat{A}(i_\Lambda-1)\cup\hat{A}(i_\Lambda)), S+(\Delta^\vee|V^*)|V)= \hat{\bar{\Lambda}}(i_\Lambda)$, we know
$\langle b_{ H_{\hat{E}}(i_\Gamma)/N^*}-b_{\Gamma/N^*}, \hat{a}(i_\Lambda-1)\rangle=
\langle b_{H_{\hat{E}}(i_\Gamma)/N^*}-b_{\Gamma/N^*}, \hat{a}(i_\Lambda)\rangle$.
We know
$\langle b_{\Gamma/N^*}, \hat{a}(i_\Lambda-1)\rangle-\langle b_{\Gamma/N^*}, \hat{a}(i_\Lambda)\rangle=
\langle b_{H_{\hat{E}}(i_\Gamma)/N^*}, \hat{a}(i_\Lambda-1)\rangle-\langle b_{H_{\hat{E}}(i_\Gamma)/N^*}, \hat{a}(i_\Lambda)\rangle$.

We have two cases.
\begin{enumerate}
\item $\Lambda\neq\hat{\bar{\Delta}}(i_\Gamma)$.
\item $\Lambda=\hat{\bar{\Delta}}(i_\Gamma)$.
\end{enumerate}

Consider the case $\Lambda\neq\hat{\bar{\Delta}}(i_\Gamma)$.

Since $\hat{\bar{\Lambda}}(i_\Lambda)=\Lambda\neq \hat{\bar{\Delta}}(i_\Gamma)
= H_{\hat{E}}(i_\Gamma)\Op|\hat{\Delta}(i_\Gamma)$, we know  $\hat{\Lambda}(i_\Lambda)\cap \hat{\Delta}(i_\Gamma)^\circ\neq\emptyset$,  $\hat{\Lambda}(i_\Lambda-1)\cap \hat{\Delta}(i_\Gamma)^\circ\neq\emptyset$,
$\{\hat{A}(i_\Gamma), \hat{A}(i_\Gamma-1)\}\subset\mathcal{F}(S+(\hat{\Delta}(i_\Gamma)^\vee|V^*))_\ell$, and 
$\langle b_{H_{\hat{E}}(i_\Gamma)/N^*}, \hat{a}(i_\Lambda-1)\rangle-\langle b_{H_{\hat{E}}(i_\Gamma)/N^*}, \hat{a}(i_\Lambda)\rangle\leq
\Ht(H_{\hat{E}}(i_\Gamma), S+(\hat{\Delta}(i_\Gamma)^\vee|V^*))
<\Ht(H,S+(\Delta^\vee|$\hfill\break$V^*))$.

We know that 
$\max\{\langle b_{\Gamma/N^*}, a\rangle|a\in A\cap(S+(\Delta^\vee|V^*))\}
-\min\{\langle b_{\Gamma/N^*}, a\rangle|a\in A\cap(S+(\Delta^\vee|V^*))\}
<\Ht(H,S+(\Delta^\vee|V^*))$.

Consider the case  $\Lambda=\hat{\bar{\Delta}}(i_\Gamma)$.

If $i_\Gamma=\hat{m}+1$, then $ H_{\hat{E}}(i_\Gamma)=H$ and
$\langle b_{H_{\hat{E}}(i_\Gamma)/N^*}, \hat{a}(i_\Lambda-1)\rangle-\langle b_{H_{\hat{E}}(i_\Gamma)/N^*}, \hat{a}(i_\Lambda)\rangle=
\langle b_{H/N^*}, \hat{a}(i_\Lambda-1)\rangle-\langle b_{H/N^*}, \hat{a}(i_\Lambda)\rangle$.

We consider the case $i_\Gamma\neq\hat{m}+1$. $b_{H_{\hat{E}}(i_\Gamma)/N^*}-b_{H/N^*}=b_{G_{\hat{E}}(i_\Gamma)/N^*}\in 
\hat{\bar{\Delta}}(i_\Gamma)=\Lambda=\hat{\bar{\Lambda}}(i_\Lambda)$.
Therefore, 
$\langle b_{H_{\hat{E}}(i_\Gamma)/N^*}-b_{H/N^*}, \hat{a}(i_\Lambda-1)\rangle=
\langle b_{H_{\hat{E}}(i_\Gamma)/N^*}-b_{H/N^*}, \hat{a}(i_\Lambda)\rangle$, and 
$\langle b_{H_{\hat{E}}(i_\Gamma)/N^*}, \hat{a}(i_\Lambda-1)\rangle-\langle b_{H_{\hat{E}}(i_\Gamma)/N^*}, \hat{a}(i_\Lambda)\rangle=
\langle b_{H/N^*}, \hat{a}(i_\Lambda-1)\rangle-\langle b_{H/N^*}, \hat{a}(i_\Lambda)\rangle$.

We know that 
$\langle b_{H_{\hat{E}}(i_\Gamma)/N^*}-b_{H/N^*}, \hat{a}(i_\Lambda-1)\rangle=
\langle b_{H_{\hat{E}}(i_\Gamma)/N^*}-b_{H/N^*}, \hat{a}(i_\Lambda)\rangle$, and 
$\langle b_{H_{\hat{E}}(i_\Gamma)/N^*}, \hat{a}(i_\Lambda-1)\rangle-\langle b_{H_{\hat{E}}(i_\Gamma)/N^*}, \hat{a}(i_\Lambda)\rangle=
\langle b_{H/N^*}, \hat{a}(i_\Lambda-1)\rangle-\langle b_{H/N^*}, \hat{a}(i_\Lambda)\rangle$.

Note that 
$\langle b_{H/N^*}, \hat{a}(i_\Lambda-1)\rangle-\langle b_{H/N^*}, \hat{a}(i_\Lambda)\rangle
\leq
\langle b_{H/N^*}, \hat{a}(1)\rangle-\langle b_{H/N^*}, \hat{a}(c(S+(\Delta^\vee|V^*)))\rangle
=\Ht(H, S+(\Delta^\vee|V^*))$ and 
$\langle b_{H/N^*}, \hat{a}(i_\Lambda-1)\rangle-\langle b_{H/N^*}, \hat{a}(i_\Lambda)\rangle
=
\langle b_{H/N^*}, \hat{a}(1)\rangle-\langle b_{H/N^*}, \hat{a}(c(S+(\Delta^\vee|V^*)))\rangle$,
if and only if, $c(S+(\Delta^\vee|V^*))=i_\Lambda =2$.

We know that 
$\max\{\langle b_{\Gamma/N^*}, a\rangle|a\in A\cap(S+(\Delta^\vee|V^*))\}
-\min\{\langle b_{\Gamma/N^*}, a\rangle|a\in A\cap(S+(\Delta^\vee|V^*))\}
\leq\Ht(H,S+(\Delta^\vee|V^*))$, and that 
$\max\{\langle b_{\Gamma/N^*}, a\rangle|a\in A\cap(S+(\Delta^\vee|V^*))\}
-\min\{\langle b_{\Gamma/N^*}, a\rangle|a\in A\cap(S+(\Delta^\vee|V^*))\}
=\Ht(H,S+(\Delta^\vee|V^*))$,
if and only if, $c(S+(\Delta^\vee|V^*))=i_\Lambda =2$.

Now, by the arguments so far we know that the inequality 
$\max\{\langle b_{\Gamma/N^*}, a\rangle|a\in A\cap(S+(\Delta^\vee|V^*))\}
-\min\{\langle b_{\Gamma/N^*}, a\rangle|a\in A\cap(S+(\Delta^\vee|V^*))\}
\leq\Ht(H,S+(\Delta^\vee|V^*))$
always holds, and the equality 
$\max\{\langle b_{\Gamma/N^*}, a\rangle|a\in A\cap(S+(\Delta^\vee|V^*))\}
-\min\{\langle b_{\Gamma/N^*}, a\rangle|a\in A\cap(S+(\Delta^\vee|V^*))\}
=\Ht(H,S+(\Delta^\vee|V^*))$
holds, if and only if,  $\Theta^\circ\subset \hat{\bar{\Delta}}(i_\Gamma)^\circ$, 
$\Lambda=\hat{\bar{\Delta}}(i_\Gamma)$ and $c(S+(\Delta^\vee|V^*))=i_\Lambda =2$.

Assume  $\Theta^\circ\subset \hat{\bar{\Delta}}(i_\Gamma)^\circ$, 
$\Lambda=\hat{\bar{\Delta}}(i_\Gamma)$ and $c(S+(\Delta^\vee|V^*))=i_\Lambda =2$.
We have
$\hat{\bar{\Lambda}}(2)= \hat{\bar{\Lambda}}( i_\Lambda)= \Lambda=\hat{\bar{\Delta}}(i_\Gamma)$, and for any $\bar{E}\in\mathcal{F}(H\Op|\Delta)_1=\mathcal{F}(\Delta)_1-\{H\}$, $ c(\Sigma(S|V)\hat{\cap}\mathcal{F}(\Delta), 2,\bar{E})=c(\Sigma(S+(\Delta^\vee|V^*)|V), 2,\bar{E})=\hat{s}( i_\Gamma-1,\bar{E})\in\Z$.

We know that if 
$\max\{\langle b_{\Gamma/N^*}, a\rangle|a\in A\cap(S+(\Delta^\vee|V^*))\}
-\min\{\langle b_{\Gamma/N^*}, a\rangle|a\in A\cap(S+(\Delta^\vee|V^*))\}
=\Ht(H,S+(\Delta^\vee|V^*))$,
then  $c(S+(\Delta^\vee|V^*))= 2$ and the structure constant of  $\Sigma(S|V)\hat{\cap}\mathcal{F}(\Delta)$ corresponding to the pair $(2,\bar{E})$ is an integer for any $\bar{E}\in\mathcal{F}(\Delta)-\{H\}$.

Convesely, assume that $c(S+(\Delta^\vee|V^*))= 2$ and the structure constant of  $\Sigma(S|V)$\break$\hat{\cap}\mathcal{F}(\Delta)$ corresponding to the pair $(2,\bar{E})$ is an integer for any $\bar{E}\in\mathcal{F}(\Delta)-\{H\}$.
By Theorem~\ref{heightinequality}.20-26 we know that
$\hat{\bar{m}}=\hat{m}=\hat{M}$, $i_\Lambda=2$, $i_\Gamma=\hat{m}+1$,
$\Theta=\Lambda=\hat{\bar{\Lambda}}(2)=\hat{\bar{\Delta}}(\hat{m}+1)$, $\Gamma=H$ and the equality
$\max\{\langle b_{\Gamma/N^*}, a\rangle|a\in A\cap(S+(\Delta^\vee|V^*))\}
-\min\{\langle b_{\Gamma/N^*}, a\rangle|a\in A\cap(S+(\Delta^\vee|V^*))\}
=\Ht(H,S+(\Delta^\vee|V^*))$
holds.

Claim 12 follows from similar arguments as in claim 8 and 9 and Theorem~\ref{heightinequality}.27-32.
\end{proof}

\section{Schemes associated with fans}
\label{toric theory}
We review the toric theory (Kempf et al.~\cite{KKMS}, Fulton~\cite{F93}, Cox~\cite{C10}) and arrange our notations for the toric theory.

Let $k$ be any algebraically closed field, let $V$ be any vector space of finite dimension over $\R$, let $N$ be any lattice in $V$ and let $\Sigma$ be any strongly convex rational fan over $N$ in $V$. Associated with the quadruplet $(k, V, N,\Sigma)$, the toric variety $X=X(k, V, N,\Sigma)$ is defined. It has the following properties:
\begin{enumerate}
\item $X$ is a separated reduced irreducible normal scheme of finite type over $k$.
$\dim X=\dim V$.
\item $X$ is complete$\Leftrightarrow |\Sigma|=V$. $X$ is smooth$\Leftrightarrow \Sigma$ is regular.
\item Associated with any $\Delta\in\Sigma$, an affine open subset $U(\Delta)=U(k, V, N,\Sigma,\Delta)$ of $X$, an locally closed subset $O(\Delta)=O(k, V, N,\Sigma,\Delta)$ of $X$ and a closed subset  $V(\Delta)=V(k, V, N,\Sigma,\Delta)$ of $X$ have been defined.
\item A structure of a group scheme over $k$ isomorphic to $\mathbb{G}_{\mathrm{m}}^{\dim V}$ has been defined on the affine open subset $T=U(\{0\})$ of $X$, where $\mathbb{G}_{\mathrm{m}}=\Spec(k[x,1/x])$ denotes the multiplicative group over $k$.

$X$ has an action $T\times X\rightarrow X$ of $T$ extending the action $T\times T\rightarrow T$ of $T$ on $T$ itself. 
\item For any $\Delta\in\Sigma$ the following holds:
\begin{enumerate}
\item
$U(\Delta)$, $O(\Delta)$ and $V(\Delta)$ are $T$-invariant, non-empty and irreducible. $O(\Delta)$ is a $T$-orbit. $V(\Delta)$ is the closure of $O(\Delta)$ in $X$. $O(\Delta)=V(\Delta)\cap U(\Delta)$.
\item
$\dim U(\Delta)=\dim V$. $\dim O(\Delta)=\dim V(\Delta)=\dim V-\dim\Delta$.
\item
$$U(\Delta)=\bigcup_{\Lambda\in\Sigma, \Lambda\subset\Delta}O(\Lambda),\quad
V(\Delta)=\bigcup_{\Lambda\in\Sigma,\Lambda\supset\Delta}O(\Lambda).$$
\end{enumerate}
\item
$$X=\bigcup_{\Delta\in\Sigma}U(\Delta)=\bigcup_{\Delta\in\Sigma}O(\Delta).$$
\item For any $\Delta\in\Sigma$ and any $\Lambda\in\Sigma$ the following holds:
\begin{enumerate}
\item
$U(\Delta\cap\Lambda)=U(\Delta)\cap U(\Lambda)$.
\item
$\Delta=\Lambda\Leftrightarrow U(\Delta)=U(\Lambda)\Leftrightarrow O(\Delta)=O(\Lambda)
\Leftrightarrow O(\Delta)\cap O(\Lambda)\neq\emptyset\Leftrightarrow V(\Delta)=V(\Lambda)$.
\item
$\Delta\subset\Lambda\Leftrightarrow U(\Delta)\subset U(\Lambda)\Leftrightarrow O(\Delta)\subset U(\Lambda)
\Leftrightarrow O(\Delta)\cap U(\Lambda)\neq\emptyset\Leftrightarrow
V(\Delta)\cap U(\Lambda)\neq\emptyset\Leftrightarrow V(\Delta)\supset V(\Lambda)$.
\end{enumerate}
\item
An isomorphism $\chi: N^*\rightarrow \mathrm{Hom}(T, \mathbb{G}_{\mathrm{m}})$ of groups from the dual lattice $N^*$ of $N$ in the dual vector space $V^*$ to the group $\mathrm{Hom}(T, \mathbb{G}_{\mathrm{m}})$ of group homomorphisms from $T$ to $\mathbb{G}_{\mathrm{m}}$ (the group of characters of $T$)
has been defined.
\item
We have natural inclusions $\mathrm{Hom}(T, \mathbb{G}_{\mathrm{m}})\subset \mathcal{O}_X(T)^\times\subset\mathcal{O}_X(T)$, where $\mathcal{O}_X(T)$ denotes the ring of regular functions over $T$, and $\mathcal{O}_X(T)^\times$ denotes the set of invertible elements in $\mathcal{O}_X(T)$.

For any regular function $\phi\in\mathcal{O}_X(T)$ over $T$, there exists a finite subset $I$ of $N^*$ such that $\phi$ is a linear combination of $\chi(m), m\in I$ over $k$. For any finite subset $I$ of $N^*$, elements $\chi(m), m\in I$ are linearly independent over $k$.
\item
For any $\Delta\in\Sigma$, $T=U(\{0\})\subset U(\Delta)$ and the ring $\mathcal{O}_X(U(\Delta))$ of regular functions over $U(\Delta)$ is the $k$-vector subspace of $\mathcal{O}_X(T)$ generated by $\chi(m), m\in N^*\cap\Delta^\vee$.
\item
We consider the the group $\mathrm{Hom}(\mathbb{G}_{\mathrm{m}}, \mathbb{G}_{\mathrm{m}})$ of group homomorphisms from $\mathbb{G}_{\mathrm{m}}$ to $\mathbb{G}_{\mathrm{m}}$ and the group $\mathrm{Hom}(\mathbb{G}_{\mathrm{m}}, T)$ of group homomorphisms from $\mathbb{G}_{\mathrm{m}}$ to $T$ (the group of one parameter subgroups of $T$). The composition of group homomorphisms defines a map $\mathrm{Hom}(\mathbb{G}_{\mathrm{m}}, T)\times\mathrm{Hom}(T, \mathbb{G}_{\mathrm{m}})\rightarrow \mathrm{Hom}(\mathbb{G}_{\mathrm{m}}, \mathbb{G}_{\mathrm{m}})$. There exists uniquely a group isomorphism $\mathrm{Hom}(\mathbb{G}_{\mathrm{m}}, \mathbb{G}_{\mathrm{m}})\rightarrow \Z$ sending $\Id_{\mathbb{G}_{\mathrm{m}}}\in\mathrm{Hom}(\mathbb{G}_{\mathrm{m}}, \mathbb{G}_{\mathrm{m}})$ to $1\in\Z$, where $\Id_{\mathbb{G}_{\mathrm{m}}}$ denotes the identity morphism of the scheme $\mathbb{G}_{\mathrm{m}}$. The composition of these two mappings defines a mapping
$$\langle \;,\;\rangle: \mathrm{Hom}(\mathbb{G}_{\mathrm{m}}, T)\times\mathrm{Hom}(T, \mathbb{G}_{\mathrm{m}})\rightarrow\Z.$$

An isomorphism $\lambda: N\rightarrow \mathrm{Hom}(\mathbb{G}_{\mathrm{m}},T)$ of groups from the lattice $N$ to the group $\mathrm{Hom}(\mathbb{G}_{\mathrm{m}},T)$
has been defined.

For any $m\in N^*$ and any $n\in N$, the equality $\langle m, n\rangle=\langle \lambda(n),\chi(m)\rangle$ holds.
\end{enumerate}

Below, we consider the case where $\Sigma$ is a flat regular fan with $\dim\Sigma=\dim V$. By $2$, $X$ is smooth.
\begin{enumerate}
\setcounter{enumi}{11}
\item For any $\Delta\in\Sigma^0$, the following holds:
\begin{enumerate}
\item
$\dim\Delta=\dim V$, $\Vect(\Delta)=V$, $\{b_{E/N}|E\in\mathcal{F}(\Delta)_1\}$ is a $\Z$-basis of $N$, it is an $\R$-basis of $V$ and $\Delta=\Convcone(\{b_{E/N}|E\in\mathcal{F}(\Delta)_1\})$.
\end{enumerate}

We denote the dual basis of $\{b_{E/N}|E\in\mathcal{F}(\Delta)_1\}$ by $\{{b_{E/N}}^\vee_\Delta |E\in\mathcal{F}(\Delta)_1\}$. For any $D\in\mathcal{F}(\Delta)_1$ and any $E\in\mathcal{F}(\Delta)_1$.
$$
\langle {b_{D/N}}^\vee_\Delta, b_{E/N} \rangle=
\begin{cases}
1&\text{ if $D=E$},\\
0&\text{ if $D\neq E$}.
\end{cases}$$
\begin{enumerate}
\setcounter{enumii}{1}
\item
$\dim \Delta^\vee=\dim V$, $\Vect(\Delta^\vee)=V^*$,  $\{{b_{E/N}}^\vee_\Delta|E\in\mathcal{F}(\Delta)_1\}$ is a $\Z$-basis of $N^*$, it is an $\R$-basis of $V^*$ and $\Delta^\vee=\Convcone(\{{b_{E/N}}^\vee_\Delta|E\in\mathcal{F}(\Delta)_1\})$.
\item
$\mathcal{O}_X(U(\Delta))$ is a polynomial ring over $k$ with variables $\{\chi({b_{E/N}}^\vee_\Delta)|E\in\mathcal{F}(\Delta)_1\}$. In other words, $\mathcal{O}_X(U(\Delta))$ is a $k$-algebra generated by the finite subset $\{\chi({b_{E/N}}^\vee_\Delta)|E\in\mathcal{F}(\Delta)_1\}$ and $\{\chi({b_{E/N}}^\vee_\Delta)|E\in\mathcal{F}(\Delta)_1\}$ is algebraically independent over $k$.
\item
Consider any $\Lambda\in\Sigma$ with $\Lambda\subset\Delta$. Note that $\Lambda$ is a face of $\Delta$ and $\mathcal{F}(\Lambda)_1\subset\mathcal{F}(\Delta)_1$.

$U(\Lambda)$ is the open subset of $U(\Delta)$ defined by $\chi({b_{E/N}}^\vee_\Delta)\neq 0$ for any $E\in\mathcal{F}(\Delta)_1-\mathcal{F}(\Lambda)_1$.

$V(\Lambda)\cap U(\Delta)$ is the closed subset of $U(\Delta)$ defined by $\chi({b_{E/N}}^\vee_\Delta)= 0$ for any $E\in\mathcal{F}(\Lambda)_1$.

$O(\Lambda)= U(\Lambda)\cap V(\Lambda)$ is the locally closed subset of $U(\Delta)$ defined by $\chi({b_{D/N}}^\vee_\Delta)\neq 0$ for any $D\in\mathcal{F}(\Delta)_1-\mathcal{F}(\Lambda)_1$ and $\chi({b_{E/N}}^\vee_\Delta)= 0$ for any $E\in\mathcal{F}(\Lambda)_1$.
\end{enumerate}
\item
For any $\Delta\in\Sigma$, $V(\Delta)$ is a smooth irreducible closed subset of $X$ with $\dim V(\Delta)=\dim V-\dim \Delta$ and it is a local complete intersection.
\item
Assume $\dim \Sigma=\dim V\geq 1$.

We denote $D=\sum_{\Gamma\in\Sigma_1}V(\Gamma)$.

The pair $(X,D)$ is a normal crossing scheme over $k$.

$\Comp(D)=\{ V(\Gamma)| \Gamma\in\Sigma_1\}$.

$(D)_0=\cup_{\Delta\in\Sigma^0}O(\Delta)$.

If $\A\in (D)_0$, $\Delta\in\Sigma^0$ and $\{\A\}= O(\Delta)$, then
$U(X, D,\A)=U(\Delta)$ and
$\Comp(D)(\A)= \{ V(\Gamma)| \Gamma\in\mathcal{F}(\Delta)_1\}$.

Consider any $\A\in (D)_0$.
We take $\Delta\in\Sigma^0$ with $\{\A\}= O(\Delta)$.
We define a mapping $\xi_\A:\Comp(D)(\A)\rightarrow
\mathcal{O}_X(U(X, D,\A))$ by putting
$\xi_\A(V(\Gamma))=\chi({b_{\Gamma/N}}^\vee_\Delta)\in\mathcal{O}_X(U(\Delta))=\mathcal{O}_X(U(X, D,\A))$ for any $\Gamma\in\mathcal{F}(\Delta)_1$.

We put $\xi=\{\xi_\A|\A\in (D)_0\}$.

The triplet $(\Sigma(\mathcal{D}),D,\xi)$ is a coordinated normal crossing scheme over $S$.
\end{enumerate}

We consider morphisms between toric varieties. 

We consider
any algebraically closed field  $k$, any vector space $V$ of finite dimension over $\R$, any lattice $N$ in $V$, any strongly convex rational fan $\Sigma$ over $N$ in $V$, the toric variety $X=X(k, V, N,\Sigma)$ associated with the quadruplet $(k, V, N,\Sigma)$, any vector space $V'$ of finite dimension over $\R$, any lattice $N'$ in $V'$, any strongly convex rational fan $\Sigma'$ over $N'$ in $V'$, the toric variety $X'=X(k, V', N',\Sigma')$ associated with the quadruplet $(k, V', N',\Sigma')$ and any homomorphism $\mu:V\rightarrow V'$ of vector spaces over $\R$.

If $\mu(N)\subset N'$ and if for any $\Delta\in\Sigma$ there exists $\Delta'\in\Sigma'$ with $\mu(\Delta)\subset \Delta'$, then there exists uniquely a morphism $\sigma:X\rightarrow X'$ of varieties over $k$ with the following properties:
\begin{enumerate}
\item
Consider $T=U(k, V, N,\Sigma,\{0\})\subset X$ and  $T'=U(k, V', N',\Sigma',\{0\})\subset X'$. $\sigma(T)\subset T'$ and the induced morphism $\bar{\sigma}:T\rightarrow T'$ by $\sigma$ is a morphism of group schemes.
\item
Let $\bar{\bar{\sigma}}:\mathrm{Hom}(\mathbb{G}_{\mathrm{m}},T)\rightarrow\mathrm{Hom}(\mathbb{G}_{\mathrm{m}},T')$ denote the group homomorphism induced by $\bar{\sigma}$.

Recall that we have group isomorphisms $\lambda:N\rightarrow \mathrm{Hom}(\mathbb{G}_{\mathrm{m}},T)$ and $\lambda':N'\rightarrow \mathrm{Hom}(\mathbb{G}_{\mathrm{m}},T')$ by $11$ above.

The composition $\lambda^{\prime -1}\bar{\bar{\sigma}}\lambda:N\rightarrow N'$ coincides with $\bar{\mu}: N\rightarrow N'$, where $\bar{\mu}$ denotes the group homomorphism induced by $\mu$.
\end{enumerate}

Below, we assume that $\mu(N)\subset N'$, for any $\Delta\in\Sigma$ there exists $\Delta'\in\Sigma'$ with $\mu(\Delta)\subset \Delta'$ and the morphism $\sigma:X\rightarrow X'$ of varieties over $k$ has the above properties.

The following holds:
\begin{enumerate}
\item
$\sigma$ is proper$\Leftrightarrow\mu^{-1}(|\Sigma'|)=|\Sigma|$.
\item
If $\Delta\in\Sigma$, $\Delta'\in\Sigma'$ and $\mu(\Delta)\subset \Delta'$, then $\sigma(U(k,V,N,\Sigma,\Delta))\subset U(k,V',N',\Sigma',$\break$\Delta')$.
\item
For any  $\Delta\in\Sigma$, $\mu(\Delta^\circ)=\mu(\Delta)^\circ$ and there exists uniquely $\Delta'\in\Sigma'$ with $\mu(\Delta^\circ)\subset\Delta^{\prime\circ}$.
\item
If $\Delta\in\Sigma$, $\Delta'\in\Sigma'$ and $\mu(\Delta^\circ)\subset\Delta^{\prime\circ}$, then $\sigma(O(k,V,N,\Sigma,\Delta))\subset O(k,V',N',$\break$\Sigma',\Delta')$.
\item
For any $\Delta'\in\Sigma'$,
\begin{equation*}\begin{split}
\sigma^{-1}( U(k,V',N',\Sigma',\Delta'))=&
\bigcup_{\Delta\in\Sigma,\Delta\subset\Delta'} U(k,V,N,\Sigma,\Delta)),
\text{ and}\\
\sigma^{-1}( O(k,V',N',\Sigma',\Delta'))=&
\bigcup_{\Delta\in\Sigma,\Delta^\circ\subset\Delta^{\prime\circ}} O(k,V,N,\Sigma,\Delta)).
\end{split}\end{equation*}
\end{enumerate}

Below, we consider the case where $V=V'$, $N=N'$ and $\mu=\Id_V$. 
In this case $\Sigma$ and $\Sigma'$ are any strongly convex rational fans over $N$ in $V$. Obviously $\Id_V(N)\subset N$, and for any $\Delta\in\Sigma$ there exists $\Delta'\in\Sigma'$ with $\Id_V (\Delta)\subset \Delta'$, if and only if, $\Sigma$ is a subdivision of $\Sigma'$. 

In case where $\Sigma$ is a subdivision of $\Sigma'$, we denote the above morphism
$$\sigma:X(k,V,N,\Sigma)\rightarrow X(k,V,N,\Sigma')$$
by the symbol $\sigma(k,V,N,\Sigma',\Sigma)$ or $\sigma(\Sigma',\Sigma)$ and we call it the \emph{subdivision morphism} associated with the pair $(\Sigma,\Sigma')$.

Note the order of $\Sigma$ and $\Sigma'$ in the symbol $\sigma(k,V,N,\Sigma',\Sigma)$.
$$\sigma(k,V,N,\Sigma',\Sigma):X(k,V,N,\Sigma)\rightarrow X(k,V,N,\Sigma').$$

The following holds:
\begin{enumerate}
\item
For any strongly convex rational fan $\Sigma$ over $N$ in $V$, $\sigma(k,V,N,\Sigma,\Sigma)=\Id_{X(k,V,N,\Sigma)}$, where $\Id_{X(k,V,N,\Sigma)}$ denotes the identity morphism of the scheme $X(k,V,N,\Sigma)$.
\item
If $\Sigma$, $\Sigma'$ and $\Sigma''$ are any strongly convex rational fans over $N$ in $V$,  $\Sigma$ is a subdivision of $\Sigma'$ and  $\Sigma'$ is a subdivision of $\Sigma''$, then $\Sigma$ is a subdivision of $\Sigma''$ and $\sigma(k,V,N,\Sigma'',\Sigma')\sigma(k,V,N,\Sigma',\Sigma)= \sigma(k,V,N,\Sigma'',\Sigma)$.
\end{enumerate}

Below, we assume that $\Sigma$ and $\Sigma'$ are any strongly convex rational fans over $N$ in $V$, and  $\Sigma$ is a subdivision of $\Sigma'$. We denote $X=X(k,V,N,\Sigma)$, $X'=X(k,V,N,\Sigma')$, and $\sigma=\sigma(k,V,N,\Sigma',\Sigma):X\rightarrow X'$.

The following holds:

\begin{enumerate}
\item
For any $\Delta\in\Sigma\cap \Sigma'$, $\sigma(U(\Sigma,\Delta))= U(\Sigma',\Delta)$ and the morphism $\sigma: U(\Sigma,\Delta)\rightarrow U(\Sigma',\Delta)$ induced by $\sigma$ is an isomorphism of schemes.
\item
If $\Sigma$ is a subset of $\Sigma'$, then $\sigma$ is an open immersion.
\item
If $\Sigma'$ is regular over $N$, $\Delta'\in\Sigma'$, $\dim \Delta'\geq 1$, and $\Sigma=\Sigma'*\Delta'$, then $\Sigma$ is regular over $N$ and the morphism $\sigma:X\rightarrow X'$ coincides with the blowing-up of $X'$ with center in $V(\Sigma',\Delta')$.
\item
Assume that $\Sigma'$ is a flat regular fan with $\dim\Sigma'=\dim V$, $m\in\Z_0$, $F$ is a center sequence of $\Sigma'$ of length $m$ and $\Sigma=\Sigma'*F(1)*F(2)*\cdots*F(m)$.

Refer to claim $14$ two pages before.

We denote $D=\sum_{\Gamma\in\Sigma_1}V(\Sigma, \Gamma)$ and $D'=\sum_{\Gamma\in\Sigma'_1}V(\Sigma', \Gamma)$.

Consider any $\A\in (D)_0$.
We take $\Delta\in\Sigma^0$ with $\{\A\}= O(\Sigma, \Delta)$.
We define a mapping $\xi_\A:\Comp(D)(\A)\rightarrow
\mathcal{O}_X(U(X, D,\A))$ by putting
$\xi_\A(V(\Sigma, \Gamma))=\chi({b_{\Gamma/N}}^\vee_\Delta)\in\mathcal{O}_X(U(\Sigma, \Delta))=\mathcal{O}_X(U(X, D,\A))$ for any $\Gamma\in\mathcal{F}(\Delta)_1$.

We put $\xi=\{\xi_\A|\A\in (D)_0\}$.

Consider any $\A\in (D')_0$.
We take $\Delta'\in\Sigma^{\prime 0}$ with $\{\A\}= O(\Sigma', \Delta')$.
We define a mapping $\xi'_\A:\Comp(D')(\A)\rightarrow
\mathcal{O}_{X'}(U(X', D',\A))$ by putting
$\xi'_\A(V(\Sigma', \Gamma))=\chi' ({b_{\Gamma/N}}^\vee_{\Delta'})\in\mathcal{O}_{X'}(U(\Sigma', \Delta'))=\mathcal{O}_{X'}(U(X', D',\A))$ for any $\Gamma\in\mathcal{F}(\Delta')_1$.

We put $\xi'=\{\xi'_\A|\A\in (D')_0\}$.
\begin{enumerate}
\item
$\Sigma$ is a flat regular fan with $\dim\Sigma=\dim V$.
\item
The morphism $\sigma:X\rightarrow X'$
is an admissible composition of blowing-ups over $D'$.
If $\dim F(i)=2$ for any $i\in\{1,2,\ldots,m\}$, then $\sigma$ is an admissible composition of blowing-ups with center of codimension two over $D'$.
\item
$\Comp(D)=\Comp(\sigma^*D')$
\item
$\sigma((D)_0)= (D')_0$. If $\A\in (D)_0$, $\Delta\in\Sigma^0$, $\{\A\}=O(\Sigma,\Delta)$, $\Delta'\in\Sigma^{\prime 0}$, $\Delta\subset\Delta'$, then $\{\sigma(\A)\}=O(\Sigma',\Delta')$.
\item
$\xi=\sigma^*\xi'$.
\end{enumerate}
\end{enumerate}

\section{Proof of the main theorem}
\label{main proof}
We give the proof of our main theorem Theorem~\ref{main}.

Let $k$ be any \emph{algebraically closed} field, let $R$ be any regular local ring such that $R$ contains $k$ as a subring, the residue field $R/M(R)$ is isomorphic to $k$ as $k$-algebras and $\dim R\geq 2$, let $P$ be any parameter system of $R$, and let $z\in P$ be any element.

Consider any $\phi\in R$ such that $\phi\neq 0$, $\Gamma_+(P,\phi)$ is $z$-simple, and $\phi$ satisfies one of the following two conditions:
\begin{enumerate}
\item
$\Inv(P,z,\phi)>0$, $\Gamma_+(P, \psi)$ has no $z$-removable faces, where $\psi$ denotes a main factor of $(P,z,\phi)$.
\item
$\Inv(P,z,\phi)=0$, $\Inv 2(P,z,\phi)\geq 2$ and $z$ divides $\phi$.
\end{enumerate}

Let $D=\Spec(R/\prod_{x\in P}x R)$, which is a normal crossing divisor on $\Spec(R)$. We define a coordinate system $\xi:\Comp(D)\rightarrow R$ of the normal crossing scheme $(\Spec(R),D)$ at $M(R)$ by putting $\xi(\Spec(R/x R))=x$ for any $x\in P$. The triplet $(\Spec(R),D, \{\xi\})$ is a coordinated normal crossing scheme over $k$.

We denote $V=\Map(P,\R)$, $N=\Map(P,\Z)$, $S=\Gamma_+(P,\phi)$ and $\Sigma=\Sigma(S|V)$.
$V$ is a vector space of finite dimension over $\R$. $\dim V=\dim R$. $N$ is a lattice of $V$. $S$ is a rational convex pseudo polyhedron over $N$ in $V$ with $\dim S=\dim V$. $\Stab(S)=\Map(P,\R_0)$ is a regular cone over $N$ in $V$ with $\dim\Stab(S)=\dim V$. $\mathcal{V}(S)\subset N\cap\Stab(S)$. $\Den(S/N)=1$.
For any $(H,\Psi)\in\mathcal{SF}(V,N,S)$, $\Ht(H,\Psi, S)\in\Z_0$.
The normal fan $\Sigma$ of $S$ is a flat rational strongly convex fan over $N^*$ in $V^*$ with $\Vect(|\Sigma|)=V^*$.
The support $|\Sigma|$ of $\Sigma$ is a regular cone over $N^*$ in $V^*$. $\dim |\Sigma|=\dim V$. $|\Sigma|=\Stab(S)^\vee|V$.

The set $\{f^P_x|x\in P\}$ is a $\R$-basis of $V$. It is a $\Z$-basis of $N$.
$\Stab(S)=\Convcone(\{f^P_x|x\in P\})$. The dual basis $\{f^{P\vee}_x|x\in P\}$ of $\{f^P_x|x\in P\}$ is a $\R$-basis of $V^*$. It is a $\Z$-basis of $N^*$.
$|\Sigma|=\Convcone(\{f^{P\vee}_x|x\in P\})=\Map(P,\R_0)^\vee|V$.

We denote $H=\R_0 f^{P\vee}_z\in\mathcal{F}(|\Sigma|)_1$. 
$S$ and $\Sigma$ are $H$-simple. $(H,\mathcal{F}(|\Sigma|))\in \mathcal{SF}(V,N,S)\neq\emptyset$
and $\mathcal{USD}(H,\mathcal{F}(|\Sigma|),S)\neq\emptyset$. We denote the $H$-skeleton of $\Sigma$ by $\bar{\Sigma}^1$.

Let $\psi\in R$ denote a main factor of $(P,z,\phi)$. $\psi\neq 0$. We denote $S_\psi= \Gamma_+(P,\psi)$ and $\Sigma_\psi=\Sigma(S_\psi|V)$. $S_\psi$ and $\Sigma_\psi$ are $H$-simple and $\{\Inv(P,z,\phi)f^P_z\}$ is the $H$-top vertex of $S_\psi$. $z$ does not divide $\psi$ and $\Inv(P,z,\phi)=\Ht(H, S_\psi)$.  We denote the $H$-skeleton of $\Sigma_\psi$ by $\bar{\Sigma}_\psi^1$.

We take $u\in R^\times$, a mapping $a:P-\{z\}\rightarrow\Z_0$, a finite subset $\Omega$ of $M(R)$ and a mapping $b:\Omega\rightarrow\Z_+$ satisfying the following three conditions:
\begin{enumerate}
\item $\phi=u(\prod_{x\in P-\{z\}}x^{a(x)})(\prod_{\omega\in\Omega}\omega^{b(\omega)})\psi$.
\item For any $\omega\in\Omega$, $\omega$ is of order one and $\partial\omega/\partial z\in R^\times$.
\item If $\omega\in\Omega$, $\omega'\in\Omega$, $v\in R^\times$ and $\omega=v\omega'$, then $v=1$ and $\omega=\omega'$.
\end{enumerate}

Consider any $\omega\in\Omega$. We denote $S_\omega=\Gamma_+(P, \omega)$ and $\Sigma_\omega=\Sigma(S_\omega|V)$. $S_\omega$ and $\Sigma_\omega$ are $H$-simple, and $\{f^P_z\}$ is the $H$-top vertex of $ S_\omega$. $\Ht(H, S_\omega)\leq 1$. $\Ord(P,f^{P\vee}_z,\omega)\leq 1$. $\Ht(H, S_\omega)+ \Ord(P,f^{P\vee}_z,\omega)=1$. We denote the $H$-skeleton of $\Sigma_\omega$ by $\bar{\Sigma}_\omega^1$. $c(S_\omega)\leq 2$.
$c(S_\omega)=1\Leftrightarrow \Ht(H, S_\omega)=0\Leftrightarrow \Ord(P,f^{P\vee}_z,\omega)=1\Leftrightarrow z$ divides $\omega$.
For any $i\in\{1,2,\ldots,c(S_\omega)\}$ and any $\bar{E}\in\mathcal{F}(H\Op||\Sigma|)_1)$, the structure constant of $\Sigma_\omega$ corresponding to the pair $(i,\bar{E})$ is an integer.

If $\Inv(P,z,\phi)>0$, then $\psi\in M(R)$ and $S_\psi$ has no $z$-removable faces. If $\Inv(P,z,\phi)=0$, then $\psi\in R^\times$, $\sharp\Omega=\Inv 2(P,z,\phi)\geq 2$ and $z$ divides $\omega$ for some $\omega\in \Omega$.

Note that $S=\sum_{x\in P-\{z\}}a(x)\Gamma_+(P, x) +\sum_{\omega\in\Omega}b(\omega)S_\omega
+S_\psi$,
$\Sigma=(\hat{\cap}_{\omega\in\Omega}\Sigma_\omega) \hat{\cap}\Sigma_\psi$,
and $\bar{\Sigma}^1=(\cup_{\omega\in\Omega}\bar{\Sigma}_\omega^1)\cup\bar{\Sigma}_\psi^1$.

Take any $(M, F)\in \mathcal{USD}(H,\mathcal{F}(|\Sigma|),S)$.
Let $\Sigma^*=\mathcal{F}(|\Sigma|)*F(1)*F(2)*\cdots*F(M)$.
$\Sigma^*$ is a flat regular fan over $N^*$ in $V^*$ with $\Vect(|\Sigma^*|)=V^*$. It is an upward subdivision of $(H,\mathcal{F}(|\Sigma|), S)$. It is a subdivision of $\mathcal{F}(|\Sigma|)$, it is a subdivision of $\Sigma$ and $|\Sigma^*|=|\mathcal{F}(|\Sigma|)|=|\Sigma|$.

Let $\hat{X}=X(k,V^*,N^*, \Sigma^*)$, $\hat{Y}= X(k,V^*,N^*, \mathcal{F}(|\Sigma|))$ and $\hat{\sigma}=\sigma(k,V^*,N^*,$\break$\mathcal{F}(|\Sigma|),\Sigma^*):\hat{X}\rightarrow \hat{Y}$. $\hat{X}$ is the toric variety over $k$ associated with the regular fan $\Sigma^*$, $\hat{Y}$ is the toric variety over $k$ associated with the regular fan $\mathcal{F}(|\Sigma|)$ and $\hat{\sigma}$ is the subdivision morphism over $k$ associated with the pair $(\Sigma^*, \mathcal{F}(|\Sigma|))$. Since $\dim F(i)=2$ for any $i\in\{1,2,\ldots,m\}$, $\hat{\sigma}$ is a composition of blowing-ups with center in a closed irreducible smooth subschems of codimension two.

$\hat{Y}=U(\mathcal{F}(|\Sigma|), |\Sigma|)$ is an affine scheme. 
The ring of regular functions $\mathcal{O}_{\hat{Y}}(\hat{Y})=\mathcal{O}_{\hat{Y}}(U(\mathcal{F}(|\Sigma|), |\Sigma|))$ over $\hat{Y}$ is a polynomial ring over $k$ with variables $\{\chi ({b_{E/N^*}}^\vee_{|\Sigma|})|$\break$E\in\mathcal{F}(|\Sigma|)_1\}=
\{\chi(f^P_x)|x\in P\}$.
By the injective homomorphism $\mathcal{O}_{\hat{Y}}(\hat{Y})\rightarrow R$ of $k$-algebras sending $\chi(f^P_x)\in \mathcal{O}_{\hat{Y}}(\hat{Y})$ to $x\in R$ for any $x\in P$, we regard $\mathcal{O}_{\hat{Y}}(\hat{Y})$ as an subring of $R$. 
$\chi(f^P_x)=x$ for any $x\in P$. $\mathcal{O}_{\hat{Y}}(\hat{Y})=k[P]\subset R$. $\hat{Y}=\Spec(k[P])$. The inclusion ring homomorphism $k[P]\rightarrow R$ induces a morphism $\hat{\pi}:\Spec(R)\rightarrow \hat{Y}$ of schemes over $k$. 

Let $\hat{D}=\Spec(k[P]/\prod_{x\in P}xk[P])$, which is a normal crossing divisor on $\hat{Y}$. $D=\hat{\pi}^*\hat{D}$. We define a coordinate system $\hat{\xi}:\Comp(\hat{D})\rightarrow k[P]$ of the normal crossing scheme $(\hat{Y},\hat{D})$ at $k[P]\cap M(R)$ by putting $\hat{\xi}(\Spec(k[P]/x k[P]))=x$ for any $x\in P$. The triplet $(\hat{Y},\hat{D}, \{\hat{\xi}\})$ is a coordinated normal crossing scheme over $k$. For any $\hat{C}\in\Comp(\hat{D})$, $\hat{\pi}^*\hat{C}\in\Comp(D)$ and $\hat{\pi}^*\hat{\xi}(\hat{C})=\xi\hat{\pi}^*(\hat{C})$.

The triplet $(\hat{X},\hat{\sigma}^*\hat{D},\hat{\sigma}^* \{\hat{\xi}\})$ is a coordinated normal crossing scheme over $k$. $(\hat{\sigma}^*\hat{D})_0=\cup_{\Xi\in\Sigma^{*0}}O(\Sigma^*,\Xi)$. 
For any $c\in(\hat{\sigma}^*\hat{D})_0$, if we take the unique $\Xi\in\Sigma^{*0}$ with $\{c\}= O(\Sigma^*,\Xi)$, then
$\Comp(\hat{\sigma}^*\hat{D})(c)=\{V(\Sigma^*,E)|E\in\mathcal{F}(\Xi)_1\}$,
$(\hat{\sigma}^*\{\hat{\xi}\})_c (V(\Sigma^*,$\break$E))= \chi({b_{E/N^*}}^\vee_\Xi)$ for any $E\in\mathcal{F}(\Xi)_1$, and thus 
$\{(\hat{\sigma}^*\{\hat{\xi}\})_c(\hat{C})|
\hat{C}\in\Comp(\hat{\sigma}^*\hat{D})(c)\}=\{\chi({b_{E/N^*}}^\vee_\Xi)|E\in\mathcal{F}(\Xi)_1\}$.

Consider the fiber product scheme $X=\hat{X}\times_{\hat{Y}}\Spec(R)$ of $\hat{X}$ and $\Spec(R)$ over $\hat{Y}$, the projection $\pi:X\rightarrow \hat{X}$ and the projection $\sigma:X\rightarrow\Spec(R)$. $\hat{\sigma}\pi=\hat{\pi}\sigma$. $X$ is the toric variety over $\Spec(R)$ associated with the fan $\Sigma^*$ and $\sigma$ is the toric morphism associated with $\Sigma^*$ and it is an admissible composition of blowing-ups with center of codimension two over $D$. $\pi$ induces an isomorphism $\pi:X\times_{\Spec(R)}\Spec(R/M(R))=\sigma^{-1}(M(R))\rightarrow
\hat{X}\times_{\hat{Y}}\Spec(k[P]/(k[P]\cap M(R))k[P])=\hat{\sigma}^{-1}(k[P]\cap M(R))$. For any affine open subset $\hat{U}$ of $\hat{X}$, $\pi^{-1}(\hat{U})$ is an affine open subset of $X$.

The triplet $(X, \sigma^*D,\sigma^* \{\xi\})$ is a coordinated normal crossing scheme over $k$. $(\sigma^*D)_0=\pi^{-1}((\hat{\sigma}^*\hat{D})_0)$. 
For any $c\in(\sigma^*D)_0$, if we take the unique $\Xi\in\Sigma^{*0}$ with $\{\pi(c)\}= O(\Sigma^*,\Xi)$, then
$\Comp(\sigma^*D)(c)=\{\pi^*V(\Sigma^*,E)|E\in\mathcal{F}(\Xi)_1\}$,
$(\sigma^*\{\xi\})_c($\break$\pi^*V(\Sigma^*,E))=\pi^*(\chi({b_{E/N^*}}^\vee_\Xi))$ for any $E\in\mathcal{F}(\Xi)_1$, and thus 
$\{(\sigma^*\{\xi\})_c(C)|C\in\Comp(\sigma^*D)(c)\}=\{\pi^*(\chi({b_{E/N^*}}^\vee_{\Xi}))|E\in\mathcal{F}(\Xi)_1\}$.

Consider any closed point $a\in X$ with $\sigma(a)=M(R)$.

$\pi(a)\in\hat{X}$ and we have a homomorphism $\pi^*:\mathcal{O}_{\hat{X},\pi(a)}\rightarrow\mathcal{O}_{X,a}$ of local $k$-algebras induced by $\pi$. We have $\pi^*(M(\mathcal{O}_{\hat{X},\pi(a)})) \mathcal{O}_{X,a}=M(\mathcal{O}_{X,a})$. We take the unique $\Theta\in\Sigma^*$ with $\pi(a)\in O(\Sigma^*,\Theta)$, we take the unique $\Lambda\in\Sigma$ with $\Theta^\circ\subset\Lambda^\circ$, and we take any $\mu\in \Theta^\circ$.

Since $\pi(a)\in O(\Sigma^*,\Theta)$, $\sigma(a)=M(R)$ and $\{k[P]\cap M(R)\}=O(\mathcal{F}(|\Sigma|),|\Sigma|)$, we know $\hat{\sigma}\pi(a)= k[P]\cap M(R)$ and $\Theta^\circ\subset\Lambda^\circ\subset|\Sigma|^\circ$. Since $\Sigma$ is $H$-simple, $\dim\Lambda=\dim \Sigma$ or $\dim\Lambda=\dim\Sigma-1$.

Let 
\begin{equation*}\begin{split}
\Sigma^{*\circ}_1&=\{\Gamma\in\Sigma^*_1|\Gamma\not\subset H\Op||\Sigma|\},\\
\Psi&=\Psi(V,N,H,\mathcal{F}(|\Sigma|),S,M,F): \Sigma^{*\circ}_1\rightarrow 2^{2^{V^*}},\text{ and}\\
\Psi^\circ&=\Psi^\circ(V,N,H,\mathcal{F}(|\Sigma|),S,M,F): \Sigma^{*\circ}_1\rightarrow 2^{2^{V^*}}.
\end{split}\end{equation*}
We take the unique $\Gamma\in\Sigma^{*\circ}_1$ such that $\Theta$ belongs to the $H$-lower main part $\Psi^\circ(\Gamma)$ of $\Sigma^*$ below $\Gamma$. $H$-lower part $\Psi(\Gamma)$ of $\Sigma^*$ below $\Gamma$ is a flat regular fan with $\Vect(|\Psi(\Gamma)|)=V^*$. 
$\Gamma\in\Psi(\Gamma)$, $\Psi(\Gamma)$ is starry with center $\Gamma$ and $\Theta\in\Psi^\circ(\Gamma)\subset\Psi(\Gamma)\subset\Sigma^*$.
We take any $\Delta\in \Psi(\Gamma)\Mx=\Psi(\Gamma)^0$ with $\Theta+\Gamma\subset\Delta$. $\Delta\in\Sigma^{*0}$. $\dim\Delta=\dim V$. $\Delta^\circ\subset|\Sigma|^\circ$. $\Gamma\in\mathcal{F}(\Delta)_1$. $\Theta\in\mathcal{F}(\Delta)$. $\mathcal{F}(\Theta)_1\subset\mathcal{F}(\Delta)_1$. $\pi(a)\in O(\Sigma^*,\Theta)\subset U(\Sigma^*,\Delta)$.

The ring of regular functions $\mathcal{O}_{\hat{X}}(U(\Sigma^*,\Delta))$ over $U(\Sigma^*,\Delta)$ is a polynomial ring over $k$ with variables $\{\chi({b_{E/N^*}}^\vee_\Delta)|E\in\mathcal{F}(\Delta)_1\}$. $O(\Sigma^*,\Theta)$ is the locally closed subset of $U(\Sigma^*,\Delta)$ defined by $\chi({b_{E/N^*}}^\vee_\Delta)\neq 0$ for any $E\in\mathcal{F}(\Delta)_1-\mathcal{F}(\Theta)_1$ and $\chi({b_{F/N^*}}^\vee_\Delta)= 0$ for any $F\in\mathcal{F}(\Theta)_1$. 
Denote $c_E=\chi({b_{E/N^*}}^\vee_\Delta)\pi(a)\in k$ for any $E\in\mathcal{F}(\Delta)_1$.
$c_E\neq 0$ for any $E\in\mathcal{F}(\Delta)_1-\mathcal{F}(\Theta)_1$ and $c_F= 0$ for any $F\in\mathcal{F}(\Theta)_1$.

Denote $\hat{P}=\{\chi({b_{E/N^*}}^\vee_\Delta)-c_E| E\in\mathcal{F}(\Delta)_1\}$. $\hat{P}$ is a parameter system of the local ring $\mathcal{O}_{\hat{X},\pi(a)}$ of $\hat{X}$ at $\pi(a)$.
Let $\bar{P}=\{\pi^*\chi({b_{E/N^*}}^\vee_\Delta)-c_E| E\in\mathcal{F}(\Delta)_1\}$. $\bar{P}$ is a parameter system of the local ring $\mathcal{O}_{X,a}$ of $X$ at $a$.

Let $\hat{b}\in\hat{X}$ be the unique closed point in $O(\Sigma^*,\Delta)$. $\hat{b}\in O(\Sigma^*,\Delta)\subset U(\Sigma^*,\Delta)$. $\hat{\sigma}(\hat{b})= k[P]\cap M(R)$. $\hat{b}\in(\hat{\sigma}^*\hat{D})_0$. 
$U(\hat{X},\hat{\sigma}^*\hat{D},\hat{b})=U(\Sigma^*,\Delta)\ni\pi(a)$. Let $\hat{P}_0=\{\chi({b_{E/N^*}}^\vee_\Delta) | E\in\mathcal{F}(\Delta)_1\}$. $\hat{P}_0$ is a parameter system of the local ring $\mathcal{O}_{\hat{X},\hat{b}}$ of $\hat{X}$ at $\hat{b}$ and $\hat{P}_0=\{(\hat{\sigma}^*\{\hat{\xi}\})_{\hat{b}}(\hat{C})|\hat{C}\in\Comp(\hat{\sigma}^*\hat{D})(\hat{b})\}$.

Let $b\in\ X$ be the unique closed point in $X$ with $\pi(b)=\hat{b}$. $\hat{b}\in \pi^{-1}(U(\Sigma^*,\Delta))$. $\sigma(b)=M(R)$. $\pi^*\hat{\sigma}^*\hat{D}=\sigma^*\hat{\pi}^*\hat{D}=\sigma^*D$. $b\in\pi^{-1}((\hat{\sigma}^*\hat{D})_0)=(\sigma^*D)_0$. The number of components of the normal crossing divisor $\sigma^*D$ on $X$ passing through $b$ is equal to $\dim R=\dim X$. $a\in\pi^{-1}(U(\hat{X},\hat{\sigma}^*\hat{D},\hat{b}))=U(X,\sigma^*D,b)$.
The point $a$ belongs to the complement $ U(X,\sigma^*D,b)$ in $X$ of the union of all components of $\sigma^*D$ not passing though $b$.
Let $\bar{P}_0=\{\pi^*\chi({b_{E/N^*}}^\vee_\Delta) | E\in\mathcal{F}(\Delta)_1\}$. Since $\pi^*(M(\mathcal{O}_{\hat{X},\hat{b}}))\mathcal{O}_{X,b}=M(\mathcal{O}_{X,b})$, $\bar{P}_0$ is a parameter system of the local ring $\mathcal{O}_{X,b}$ of $X$ at $b$ and $\bar{P}_0=\{(\sigma^*\{\xi\})_{b}(C)|C\in\Comp(\sigma^*D)(b)\}$.

We know $\bar{P}=\{(\sigma^*\{\xi\})_{b}(C)- (\sigma^*\{\xi\})_{b}(C)(a)|C\in\Comp(\sigma^*D)(b)\}$.

$\Comp(\sigma^*D)(b)=\{\pi^*\hat{C}|\hat{C}\in\Comp(\hat{\sigma}^*\hat{D})(\hat{b})\}=\{\pi^*V(\Sigma^*,E)|E\in\mathcal{F}(\Delta)_1\}$. Recall $\Gamma\in\mathcal{F}(\Delta)_1$. Let $\bar{C}=\pi^*V(\Sigma^*,\Gamma)$. $\bar{C}$ is a component passing through $b$ of the pull-back $\sigma^*D$ of the divisor $D$ by $\sigma$. $\bar{C}\in\Comp(\sigma^*D)(b)$. Let $\bar{z}=(\sigma^*\{\xi\})_{b}(\bar{C})- (\sigma^*\{\xi\})_{b}(\bar{C})(a)=\pi^*\chi({b_{\Gamma/N^*}}^\vee_\Delta)-c_\Gamma\in\bar{P}$.

We consider the $k$-algebra homomorphism $\sigma^*:R\rightarrow\mathcal{O}_{X,a}$ induced by $\sigma$. Since $\sigma$ is a composition of blowing-ups, $\sigma^*$ is injective. We have $\mathcal{O}_{X,a}\ni\sigma^*(\phi)\neq 0$.

By the isomorphism $\Map(\bar{P},\R)\rightarrow V$ of vector spaces over $\R$ sending\hfill\break $f^{\bar{P}}_{\chi({b_{E/N^*}}^\vee_\Delta)-c_E}\in\Map(\bar{P},\R)$ to ${b_{E/N^*}}^\vee_\Delta\in V$ for any $E\in\mathcal{F}(\Delta)_1$, we identify $\Map(\bar{P},\R)$ and $V$. Pairs $\Map(\bar{P},\Z)$ and $N$, $\Map(\bar{P},\R_0)$ and $\Delta^\vee|V^*$ are identified.

Consider any $\theta\in\Theta$ and any $\zeta\in R$. Since $\Theta\subset \Delta=\Map(\bar{P},\R_0)^\vee|\Map(\bar{P},\R) \subset|\Sigma|=\Map(P,\R_0)^\vee|V$, elements $\Ord(\bar{P},\theta,\sigma^*(\zeta))\in\R_0$, $\Ord(P,\theta,\zeta)\in\R_0$, $\In(\bar{P},\theta,$\break$\sigma^*(\zeta))\in\mathcal{O}_{X,a}^c$, and
$\In(P,\theta,\zeta)\in R^c$ are defined.
$c_E=0$ for any $E\in\mathcal{F}(\Theta)_1$.

It follows that $\Ord(\bar{P},\theta,\sigma^*(\zeta))=\Ord(P,\theta,\zeta)$ and  $\In(\bar{P},\theta,\sigma^*(\zeta))=\sigma^*(\In(P,\theta,\zeta))$ for any $\theta\in\Theta$ and any $\zeta\in R$.

We consider the case $\dim\Lambda=\dim\Sigma$.

$\Lambda\in\Sigma^0$. We take the unique $A\in\mathcal{V}(S)$ with $\Delta(\{A\},S|V)=\Lambda$. Note that $ A\in\mathcal{V}(S)\subset\Stab(S)\cap N$ and $\langle \theta, A\rangle\in\Z_0$ for any $\theta\in |\Sigma|\cap N^*$. Recall $\mu\in\Theta^\circ\subset\Lambda^\circ=\Delta^\circ(\{A\},S|V)$. It follows $\Ord(P,\mu,\phi)=\langle \mu, A\rangle$ and $\In(P,\mu,\phi)=r\hat{\pi}^*\chi(A)$ for some unique $r\in k-\{0\}$. We take the unique $r\in k-\{0\}$ with $\In(P,\mu,\phi)=r\hat{\pi}^*\chi(A)$.
It follows $\Ord(\bar{P},\mu,\sigma^*(\phi)) = \Ord(P,\mu,\phi)=\langle \mu, A\rangle$ and $\In(\bar{P},\mu,\sigma^*(\phi)) =\sigma^*(\In(P,\mu,\phi))=\sigma^*( r\hat{\pi}^*\chi(A))=r\pi^*\chi(A)$.
We know that $\Supp(\bar{P},\pi^*\chi(A))=\Supp(\bar{P},r\pi^*\chi(A))\subset\Supp(\bar{P},\sigma^*(\phi))\subset\Gamma_+(\bar{P},\sigma^*(\phi))$.

Now, 
\begin{equation*}\begin{split}
&\chi(A)=\chi(\sum_{E\in\mathcal{F}(\Delta)_1}\langle b_{E/N^*}, A\rangle{b_{E/N^*}}^\vee_\Delta)=\prod_{E\in\mathcal{F}(\Delta)_1}\chi({b_{E/N^*}}^\vee_\Delta)^{\langle b_{E/N^*}, A\rangle}\\
&=\prod_{E\in\mathcal{F}(\Delta)_1-\mathcal{F}(\Theta)_1}((\chi({b_{E/N^*}}^\vee_\Delta)-c_E)+c_E) ^{\langle b_{E/N^*}, A\rangle}\prod_{E\in\mathcal{F}(\Theta)_1}\chi({b_{E/N^*}}^\vee_\Delta)^{\langle b_{E/N^*}, A\rangle}.
\end{split}\end{equation*}
Since $c_E\neq 0$ for any $E\in\mathcal{F}(\Delta)_1-\mathcal{F}(\Theta)_1$ and $\bar{P}=\{\pi^*\chi({b_{E/N^*}}^\vee_\Delta) -c_E| E\in\mathcal{F}(\Delta)_1-\mathcal{F}(\Theta)_1\}\cup\{\pi^*\chi({b_{E/N^*}}^\vee_\Delta) | E\in\mathcal{F}(\Theta)_1\}$, we know
$\sum_{ E\in\mathcal{F}(\Theta)_1}\langle b_{E/N^*}, A\rangle{b_{E/N^*}}^\vee_\Delta\in\Supp(\bar{P},\pi^*\chi(A))$.

Let $\bar{A}=\sum_{ E\in\mathcal{F}(\Theta)_1}\langle b_{E/N^*}, A\rangle{b_{E/N^*}}^\vee_\Delta\in (\Delta^\vee|V^*)\cap N$. $\bar{A}\in \Supp(\bar{P},\pi^*\chi(A)) \subset\Gamma_+(\bar{P},\sigma^*(\phi))$.

$0\leq\Ord(\bar{P},b_{E/N^*},\sigma^*(\phi))\leq \langle b_{E/N^*},\bar{A}\rangle=0$ for any $ E\in\mathcal{F}(\Delta)_1-\mathcal{F}(\Theta)_1$. Thus, $\Ord(\bar{P},b_{E/N^*},\sigma^*(\phi))= \langle b_{E/N^*},\bar{A}\rangle=0$ for any $E\in\mathcal{F}(\Delta)_1-\mathcal{F}(\Theta)_1$.

$\Ord(\bar{P},b_{E/N^*},\sigma^*(\phi))=\Ord(P, b_{E/N^*},\phi)= \langle b_{E/N^*},A\rangle=\langle b_{E/N^*},\bar{A}\rangle$ for any $E\in\mathcal{F}(\Theta)_1$, since $ b_{E/N^*}\in E\subset\Theta\subset\Lambda=\Delta(\{A\},S|V)$.

We know that $\Gamma_+(\bar{P},\sigma^*(\phi))=\{\bar{A}\}+\Map(\bar{P},\R_0)$, $c(\Gamma_+(\bar{P},\sigma^*(\phi)))=1$ and $\sigma^*(\phi)$ has normal crossings over $\bar{P}$.

It follows that $\Gamma_+(\bar{P},\sigma^*(\phi))$ is of $\bar{z}$-Weirstrass type, $\partial\bar{\omega}/\partial\bar{z}\in\mathcal{O}_{X,a}^\times$ for any divisor $\bar{\omega}$ of $\sigma^*(\phi)$ of order one such that any $\bar{x}\in\bar{P}-\{\bar{z}\}$ does not divide $\bar{\omega}$, $\Inv(\bar{P},\bar{z},\sigma^*(\phi))=0$ and $\Inv 2(\bar{P},\bar{z},\sigma^*(\phi))\leq 1$.
If $\Inv(P,z,\phi)>0$, then $\Inv(\bar{P},\bar{z},$\break$\sigma^*(\phi))=0<\Inv(P,z,\phi)$.
If $\Inv(P,z,\phi)=0$, then $\Inv 2(P,z,\phi)\geq 2$ by our assumption and $\Inv 2(\bar{P},\bar{z},\sigma^*(\phi))\leq 1<2\leq \Inv 2(P,z,\phi)$.

We conclude that Theorem~\ref{main} holds, if $\dim\Lambda=\dim \Sigma$.

We consider the case $\dim\Lambda=\dim\Sigma-1$.

By Theorem~\ref{important}.$8.(b)$ we know that $\Vect(\Lambda)+\Gamma=\Vect(\Lambda)+H$. Thus $\Vect(\Vect($\break$\Lambda)+\Gamma))=\Vect(\Vect(\Lambda)+H)=V^*\neq\Vect(\Lambda)$ and $\Gamma\not\subset\Vect(\Lambda)$. Since $\Theta\subset\Lambda\subset\Vect(\Lambda)$, we have $\Gamma\not\subset\Theta$, $\Gamma\not\in\mathcal{F}(\Theta)_1$, $\Gamma\in\mathcal{F}(\Delta)_1-\mathcal{F}(\Theta)_1$ and $c_\Gamma\neq 0$.

Since $\Lambda^\circ\subset|\Sigma|^\circ$, $\Lambda\in\bar{\Sigma}^1$. 
Since $\Sigma$ is $H$-simple, there exist uniquely $\Xi_1\in\Sigma^0$, $\Xi_2\in\Sigma^0$ such that $\Lambda=\Xi_1\cap\Xi_2$, $\Lambda+(-H)=\Xi_1+(-H)$ and $\Lambda+H=\Xi_2+H$. We take $\Xi_1\in\Sigma ^0$ and $\Xi_2\in\Sigma ^0$ satisfying this condition. We take the unique $A_1\in\mathcal{V}(S)$ with $\Xi_1=\Delta(\{A_1\},S|V)$ and we take the unique $A_2\in\mathcal{V}(S)$ with $\Xi_2=\Delta(\{A_2\},S|V)$. $\langle b_{H/N^*}, A_1-A_2\rangle>0$. It follows $\langle b_{\Gamma/N^*}, A_1-A_2\rangle>0$, since $\Vect(\Lambda)+\Gamma=\Vect(\Lambda)+H$ and $\Lambda=\Xi_1\cap\Xi_2$. $\langle b_{E/N^*}, A_1\rangle =\langle b_{E/N^*}, A_2\rangle$ for any $E\in\mathcal{F}(\Theta)_1$, since $b_{E/N^*}\in E\subset\Theta\subset\Lambda=\Xi_1\cap\Xi_2$ for any $E\in\mathcal{F}(\Theta)_1$. Note that $\{A_1,A_2\}\subset\mathcal{V}(S)\subset\Stab(S)\cap N$ and $\langle \theta, A_i\rangle\in\Z_0$ for any $\theta\in|\Sigma|\cap N^*$ and any $i\in\{1,2\}$. Let $F=\Conv(\{A_1,A_2\})$. $F\in\mathcal{F}(S)_1$. Recall $\mu\in\Theta^\circ\subset\Lambda^\circ$. $F=\Delta(\mu, S|V)$ and $\Lambda=\Delta(F,S|V)$. $\Ord(\bar{P},\mu,\sigma^*(\phi))=\Ord(P,\mu,\phi)=\langle \mu, A_1\rangle =\langle \mu, A_2\rangle$. We know that there exists uniquely a mapping $r:F\cap N\rightarrow k$ with $\In(P,\mu,\phi)=\Ps(P,F,\phi)=\sum_{B\in F\cap N}r(B)\hat{\pi}^*\chi(B)$. We take the mapping $r:F\cap N\rightarrow k$ satisfying this equality. We know $r(A_1)\neq 0$ and $r(A_2)\neq 0$. $\In(\bar{P},\mu,\sigma^*(\phi))=\sigma^*(\In(P,\mu,\phi))=\sigma^*(\sum_{B\in F\cap N}r(B)\hat{\pi}^*\chi(B))= \sum_{B\in F\cap N}r(B)\pi^*\chi(B)$. Now, for any $B\in F\cap N$, 
\begin{equation*}\begin{split}
&\chi(B)=\chi(\sum_{E\in\mathcal{F}(\Delta)_1}\langle b_{E/N^*}, B\rangle{b_{E/N^*}}^\vee_\Delta)=\prod_{E\in\mathcal{F}(\Delta)_1}\chi({b_{E/N^*}}^\vee_\Delta)^{\langle b_{E/N^*}, B\rangle}\\
&=\prod_{E\in\mathcal{F}(\Delta)_1-\mathcal{F}(\Theta)_1}((\chi({b_{E/N^*}}^\vee_\Delta)-c_E)+c_E) ^{\langle b_{E/N^*}, B\rangle}\prod_{E\in\mathcal{F}(\Theta)_1}\chi({b_{E/N^*}}^\vee_\Delta)^{\langle b_{E/N^*}, A_1\rangle}.
\end{split}\end{equation*}
Since $c_E\neq 0$ for any $E\in\mathcal{F}(\Delta)_1-\mathcal{F}(\Theta)_1$, $\Gamma\in\mathcal{F}(\Delta)_1-\mathcal{F}(\Theta)_1$ and $\bar{P}=\{\pi^*\chi({b_{E/N^*}}^\vee_\Delta) -c_E| E\in\mathcal{F}(\Delta)_1-\mathcal{F}(\Theta)_1\}\cup\{\pi^*\chi({b_{E/N^*}}^\vee_\Delta) | E\in\mathcal{F}(\Theta)_1\}$ and $\langle b_{\Gamma/N^*},$\break$A_1\rangle>\langle b_{\Gamma/N^*},A_2\rangle$, we know
$\sum_{E\in\mathcal{F}(\Theta)_1}\langle b_{E/N^*}, A_1\rangle{b_{E/N^*}}^\vee_\Delta+\langle b_{\Gamma/N^*}, A_1\rangle{b_{\Gamma/N^*}}^\vee_\Delta \in\Supp(\bar{P},\In(\bar{P},\mu,\sigma^*(\phi)))$.

Let $\bar{A}=\sum_{ E\in\mathcal{F}(\Theta)_1}\langle b_{E/N^*}, A_1\rangle{b_{E/N^*}}^\vee_\Delta+\langle b_{\Gamma/N^*}, A_1\rangle{b_{\Gamma/N^*}}^\vee_\Delta \in (\Delta^\vee|V^*)\cap N$. $\bar{A}\in \Supp(\bar{P},\In(\bar{P},\mu,\sigma^*(\phi)))\subset\Supp(\bar{P},\sigma^*(\phi)) \subset\Gamma_+(\bar{P},\sigma^*(\phi))$.

$0\leq\Ord(\bar{P},b_{E/N^*},\sigma^*(\phi))\leq \langle b_{E/N^*},\bar{A}\rangle=0$ for any $ E\in\mathcal{F}(\Delta)_1-(\mathcal{F}(\Theta)_1\cup\{\Gamma\})$. Thus, $\Ord(\bar{P},b_{E/N^*},\sigma^*(\phi))= \langle b_{E/N^*},\bar{A}\rangle=0$ for any $E\in\mathcal{F}(\Delta)_1-(\mathcal{F}(\Theta)_1\cup\{\Gamma\})$.

$\Ord(\bar{P},b_{E/N^*},\sigma^*(\phi))=\Ord(P, b_{E/N^*},\phi)= \langle b_{E/N^*},A_1\rangle=\langle b_{E/N^*},\bar{A}\rangle$ for any $E\in\mathcal{F}(\Theta)_1$, since $ b_{E/N^*}\in E\subset\Theta\subset\Lambda\subset\Xi_1=\Delta(\{A_1\},S|V)$ for any $E\in\mathcal{F}(\Theta)_1$.

Since $\Ord(\bar{P},b_{E/N^*},\sigma^*(\phi))= \langle b_{E/N^*},\bar{A}\rangle$ for any $E\in\mathcal{F}(\Delta)_1-\{\Gamma\}$, $\bar{A}\in\Gamma_+(\bar{P},$\break$\sigma^*(\phi))$ and $\bar{z}=\pi^*\chi({b_{\Gamma/N^*}}^\vee_\Delta)-c_\Gamma$, it follows that $\Gamma_+(\bar{P},\sigma^*(\phi))$ is of $\bar{z}$-Weirstrass type.

We examine each factor of $\phi$ one by one.

$u\in R^\times$ and $\sigma^*(u)\in\mathcal{O}_{X,a}^\times$.

Consider any $x\in P-\{z\}$. It is easy to see that $\sigma^*(x)\in\mathcal{O}_{X,a}$ has normal crossings over $\bar{P}$ and $\Ord(\bar{P}, b_{E/N^*},\sigma^*(x))=0$ for any $E\in\mathcal{F}(\Delta)_1-\mathcal{F}(\Theta)_1$.

Consider any $\omega\in\Omega$. We have three cases.
\begin{enumerate}
\item $z$ divides $\omega$.
\item $z$ does not divide $\omega$ and $\Lambda\not\in\bar{\Sigma}_\omega^1$.
\item $z$ does not divide $\omega$ and $\Lambda\in\bar{\Sigma}_\omega^1$.
\end{enumerate}

We consider the first case. Assume that $z$ divides $\omega$.
Since $\Ord(\omega)=1$, there exists uniquely an element $v\in R^\times$ with $\omega=vz$. It is easy to see that $\sigma^*(\omega)\in\mathcal{O}_{X,a}$ has normal crossings over $\bar{P}$ and $\Ord(\bar{P}, b_{E/N^*},\sigma^*(\omega))=0$ for any $E\in\mathcal{F}(\Delta)_1-\mathcal{F}(\Theta)_1$.

We consider the second case. Assume that $z$ does not divide $\omega$ and $\Lambda\not\in\bar{\Sigma}_\omega^1$. $\Ht(H,S_\omega)=1$ and $\sharp\bar{\Sigma}_\omega^1=c(S_\omega)=2$. Let $\Lambda_\omega$ be the unique element of $\bar{\Sigma}_\omega^1$ different from $H\Op||\Sigma|$. Since $\Lambda\in\bar{\Sigma}^1\supset \bar{\Sigma}_\omega^1$ and $\Sigma$ and $\Sigma_\omega$ are $H$-simple, we know that $\Lambda^\circ\cap\Lambda_\omega^\circ=\emptyset$ and there exists uniquely $\Xi\in\Sigma_\omega^0$ with $\Lambda^\circ\subset\Xi^\circ$. We take the unique $\Xi\in\Sigma_\omega^0$ with $\Lambda^\circ\subset\Xi^\circ$.
Since $\mu\in\Theta^\circ\subset\Lambda^\circ\subset\Xi^\circ$, we can apply the same reasoning as in the case $\dim\Lambda=\dim\Sigma$ and we know that $\sigma^*(\omega)\in\mathcal{O}_{X,a}$ has normal crossings over $\bar{P}$ and $\Ord(\bar{P}, b_{E/N^*},\sigma^*(\omega))=0$ for any $E\in\mathcal{F}(\Delta)_1-\mathcal{F}(\Theta)_1$.

We consider the third case. Assume that $z$ does not divide $\omega$ and $\Lambda\in\bar{\Sigma}_\omega^1$. $\Ht(H,S_\omega)=1$, $\sharp\bar{\Sigma}_\omega^1=c(S_\omega)=2$ and $\bar{\Sigma}_\omega^1=\{\Lambda, H\Op||\Sigma|\}$.

Let $\Xi_1\in\Sigma_\omega^0$ and $\Xi_2\in\Sigma_\omega^0$ be elements with $H\Op||\Sigma\subset\Xi_1$ and $H\subset\Xi_2$. Since $c(S_\omega)=2$ and $\Sigma_\omega$ is $H$-simple, we know $\Xi_1\neq\Xi_2$, $\Sigma_\omega^0=\{\Xi_1,\Xi_2\}$, $\Lambda=\Xi_1\cap\Xi_2$, $\Lambda+(-H)=\Xi_1+(-H)$ and $\Lambda+H=\Xi_2+H$. We take the unique $A_1\in\mathcal{V}(S_\omega)$ with $\Xi_1=\Delta(\{A_1\},S_\omega|V)$, and we take the unique $A_2\in\mathcal{V}(S_\omega)$ with $\Xi_2=\Delta(\{A_2\},S_\omega|V)$.
We know $\langle b_{H/N^*}, A_1-A_2\rangle>0$. Since $\Vect(\Lambda)+\Gamma=\Vect(\Lambda)+H$ and $\Lambda=\Xi_1\cap\Xi_2$, it follows $\langle b_{\Gamma/N^*}, A_1-A_2\rangle>0$. Since $\partial\omega/\partial z\in R^\times$, we know $A_1=f^P_z$ and $\langle b_{H/N^*}, A_2\rangle=0$. Note that $\{A_1, A_2\}\subset\Stab(S_\omega)\cap N=\Stab(S)\cap N$, and $\langle \theta, A_i\rangle\in\Z_0$ for any $\theta\in |\Sigma|\cap N^*$ and any $i\in\{1,2\}$. Since $\Theta\subset\Lambda=\Xi_1\cap\Xi_2$, $\langle b_{E/N^*}, A_1\rangle=\langle b_{E/N^*}, A_2\rangle$ for any $E\in\mathcal{F}(\Theta)_1$. Since $\Sigma^*$ is a subdivision of $\Sigma_\omega$, $\Gamma\subset\Delta\in\Sigma^{*0}$ and $\Vect(\Lambda)+\Gamma=\Vect(\Lambda)+H$, we know $\Delta\subset\Xi_2$, and $\langle b_{E/N^*}, A_1-A_2\rangle\geq 0$ for any $E\in\mathcal{F}(\Delta)_1-\mathcal{F}(\Theta)_1$. 
Let $F=\Convcone(\{A_1, A_2\})$. $F\in\mathcal{F}(S_\omega)_1$ and $\Lambda=\Delta(F,S_\omega|V)$. Recall $\mu\in\Theta^\circ\subset\Lambda^\circ$.
We know $F=\Delta(\mu,S_\omega|V)$.

Now, by Theorem~\ref{important}.$12(b)$,
$\langle b_{\Gamma/N^*}, A_1-A_2\rangle=\langle b_{\Gamma/N^*}, A_1\rangle-\langle b_{\Gamma/N^*}, A_2\rangle=\max\{\langle b_{\Gamma/N^*}, c\rangle|c\in F\}-\min\{\langle b_{\Gamma/N^*}, c\rangle|c\in F\}\leq\Ht(H,S_\omega)=1$.
Since $\langle b_{\Gamma/N^*}, A_1-A_2\rangle\in\Z$, we conclude $\langle b_{\Gamma/N^*}, A_1-A_2\rangle=1$.

We know that there exist uniquely $r_1\in k-\{0\}$ and $r_2\in k-\{0\}$ with
$\In(P,\mu,\omega)=\Ps(P,F,\omega)=r_1\hat{\pi}^*\chi(A_1)+r_2\hat{\pi}^*\chi(A_2)$. We take $r_1\in k-\{0\}$ and $r_2\in k-\{0\}$ satisfying this equality. We have
\begin{equation*}\begin{split}
&\In(\bar{P},\mu,\sigma^*(\omega))=\sigma^*(\In(P,\mu,\omega))=\sigma^*(r_1\hat{\pi}^*\chi(A_1)+r_2\hat{\pi}^*\chi(A_2))=\\
&r_1\pi^*\chi(A_1)+r_2\pi^*\chi(A_2)=\\
&(r_1(\prod_{E\in\mathcal{F}(\Delta)_1-(\mathcal{F}(\Theta)_1\cup\{\Gamma\})}\pi^*\chi({b_{E/N^*}}^\vee_\Delta)^{\langle b_{E/N^*},A_1-A_2\rangle})\pi^*\chi({b_{\Gamma/N^*}}^\vee_\Delta)+r_2)\\
&\qquad\qquad\qquad\qquad\qquad\qquad (\prod_{E\in\mathcal{F}(\Delta)_1}\pi^*\chi({b_{E/N^*}}^\vee_\Delta)^{\langle b_{E/N^*},A_2\rangle}).
\end{split}\end{equation*}
We put $\bar{A}=(\sum_{E\in\mathcal{F}(\Theta)_1}\langle b_{E/N^*},A_2\rangle{b_{E/N^*}}^\vee_\Delta)+{b_{\Gamma/N^*}}^\vee_\Delta\in (\Delta^\vee|V^*)\cap N$. We know $\bar{A}\in\Gamma_+(\bar{P},\sigma^*(\omega))$, $\Ord(\bar{P},b_{E/N^*},\sigma^*(\omega))=\langle b_{E/N^*},\bar{A}\rangle$ for any $E\in\mathcal{F}(\Delta)_1-\{\Gamma\}$ and $\langle b_{E/N^*},\bar{A}\rangle=0$ for any $E\in\mathcal{F}(\Delta)_1-(\mathcal{F}(\Theta)_1\cup\{\Gamma\})$. 

We know that there exists uniquely an element $\bar{\omega}\in\mathcal{O}_{X,a}$ with $\sigma^*(\omega)=$\hfill\break$ (\prod_{E\in\mathcal{F}(\Theta)_1}\pi^*\chi({b_{E/N^*}}^\vee_\Delta)^{\langle b_{E/N^*},A_2\rangle})\bar{\omega}$. We take the unique element $\bar{\omega}\in\mathcal{O}_{X,a}$ satisfying this equality. We know ${b_{\Gamma/N^*}}^\vee_\Delta\in\Gamma_+(\bar{P},\bar{\omega})$.

Since $\bar{z}=\pi^*\chi({b_{\Gamma/N^*}}^\vee_\Delta)-c_\Gamma$, we know that either $\bar{\omega}\in\mathcal{O}_{X,a}^\times$, or $\bar{\omega}\in M(\mathcal{O}_{X,a})$ and $\partial \bar{\omega}/\partial\bar{z}\in \mathcal{O}_{X,a}^\times$.

By the reasoning so far we know that the following claim is true:
There exists $\bar{u}\in\mathcal{O}_{X,a}$, a mapping $\bar{a}:\bar{P}-\{\bar{z}\}\rightarrow\Z_0$, a finite subset $\bar{\Omega}$ of $\mathcal{O}_{X,a}$ and a mapping $\bar{b}:\bar{\Omega}\rightarrow\Z_0$ satisfying the following conditions:
\begin{enumerate}
\item
$\sigma^*(u(\prod_{x\in P-\{z\}}x^{a(x)})(\prod_{\omega\in\Omega}\omega^{b(\omega)}))=
\bar{u}(\prod_{\bar{x}\in \bar{P}-\{\bar{z}\}}\bar{x}^{\bar{a}(\bar{x})})
(\prod_{\bar{\omega}\in\bar{\Omega}}\bar{\omega}^{\bar{b}(\bar{\omega})})$.
\item
$\bar{\omega}$ is of order one and $\partial\bar{\omega}/\partial\bar{z}\in\mathcal{O}_{X,a}^\times$ for any $\bar{\omega}\in\bar{\Omega}$.
\item
If $\bar{\omega}\in\bar{\Omega}$, $\bar{\omega}'\in\bar{\Omega}$, $\bar{v}\in\mathcal{O}_{X,a}^\times$ and $\bar{\omega}=\bar{v}\bar{\omega}'$, then $\bar{v}=1$ and  $\bar{\omega}=\bar{\omega}'$.
\item
$\sharp\bar{\Omega}\leq\sharp\{\omega\in\Omega|z\text{ does not divide }\omega\text{ and }\Lambda\in\bar{\Sigma}_\omega^1\}$.
\end{enumerate}

We take $\bar{u}, \bar{a},\bar{\Omega}$ and $\bar{b}$ satisfying the above conditions.

We consider the case $\Inv(P,z,\phi)=0$. Recall that $\psi\in R$ is a main factor of $(P,z,\phi)$ and $\phi=u(\prod_{x\in P-\{z\}}x^{a(x)})(\prod_{\omega\in\Omega}\omega^{b(\omega)})\psi$. We have
$\sigma^*(\phi)= \bar{u}(\prod_{\bar{x}\in \bar{P}-\{\bar{z}\}}\bar{x}^{\bar{a}(\bar{x})})
(\prod_{\bar{\omega}\in\bar{\Omega}}\bar{\omega}^{\bar{b}(\bar{\omega})})\sigma^*(\psi)$. Since $\Inv(P,z,\phi)=0$, $\psi\in R^\times$ and $\sigma^*(\psi)\in\mathcal{O}_{X,a}^\times$.
We know that $\partial\bar{\omega}/\partial\bar{z}\in \mathcal{O}_{X,a}^\times$ for any divisor $\bar{\omega}\in\mathcal{O}_{X,a}$ of $\sigma^*(\phi)$ of order one such that any $\bar{x}\in\bar{P}-\{\bar{z}\}$ does not divide $\bar{\omega}$, any main factor of $\sigma^*(\phi)$ is an invertible element of $\mathcal{O}_{X,a}$ and $\Inv(\bar{P},\bar{z},\sigma^*(\phi))=0$.

By our assumption $\Inv 2(P,z,\phi)=\sharp\Omega\geq 2$ and $z$ divides $\phi$. It follows that $z$ divides $\omega$ for some $\omega\in\Omega$.
Thus, $\Inv 2(\bar{P},\bar{z},\sigma^*(\phi))=\sharp \bar{\Omega}\leq \sharp\{\omega\in\Omega|z\text{ does not divide }$\break$\omega\text{ and }\Lambda\in\bar{\Sigma}_\omega^1\}<\sharp\Omega=\Inv 2(P,z,\phi)$. 

We conclude that Theorem~\ref{main} holds, if $\dim\Lambda=\dim \Sigma-1$ and $\Inv(P,z,\phi)=0$.

Below, we assume $\Inv(P,z,\phi)>0$. It follows that $S_\psi$ has no $z$-removable faces from our assumption.

We consider the main factor $\psi$ of $(P,z,\phi)$.

$\phi= u(\prod_{x\in P-\{z\}}x^{a(x)})(\prod_{\omega\in\Omega}\omega^{b(\omega)})\psi$. Let $T=\Gamma_+(P, u(\prod_{x\in P-\{z\}}x^{a(x)})$\break$(\prod_{\omega\in\Omega}\omega^{b(\omega)}))$. We have $S=S_\psi+T$, $\Sigma=\Sigma_\psi\hat{\cap}\Sigma(T|V)$ and $\Sigma(T|V)=\hat{\cap}_{\omega\in\Omega}\Sigma_\omega$, where $\Sigma(T|V)$ denotes the normal fan of $T$. Since $\Sigma$ is $H$-simple, it follows that $\Sigma_\psi$ and $\Sigma(T|V)$ are $H$-simple and that $\Sigma_\omega$ is $H$-simple for any $\omega\in\Omega$. It is easy to see that any structure constant of $\Sigma_\omega$ is an integer for any $\omega\in\Omega$, since $\partial\omega/\partial z\in R^\times$ for any $\omega\in\Omega$. We know that any structure constant of $\Sigma(T|V)$ is an integer, and we can apply Theorem~\ref{important}.$12.(c)$ in our situation under consideration.

$\sigma^*(\phi)= \bar{u}(\prod_{\bar{x}\in \bar{P}-\{\bar{z}\}}\bar{x}^{\bar{a}(\bar{x})})
(\prod_{\bar{\omega}\in\bar{\Omega}}\bar{\omega}^{\bar{b}(\bar{\omega})})\sigma^*(\psi)$. It follows that any main factor of $(\bar{P},\bar{z},\sigma^*(\psi))$ is a main factor of $(\bar{P},\bar{z},\sigma^*(\phi))$.

We have two cases.
\begin{enumerate}
\item $\Lambda\not\in\bar{\Sigma}_\psi^1$.
\item $\Lambda\in\bar{\Sigma}_\psi^1$.
\end{enumerate}

We consider the first case. Assume $\Lambda\not\in\bar{\Sigma}_\psi^1$. Let $\Lambda_\psi$ be any element of $\bar{\Sigma}_\psi^1$. Since $\Lambda\in\bar{\Sigma}^1\supset \bar{\Sigma}_\psi^1$ and $\Sigma$ and $\Sigma_\psi$ are $H$-simple, we know that $\Lambda^\circ\cap\Lambda_\psi^\circ=\emptyset$ and there exists uniquely $\Xi\in\Sigma_\psi^0$ with $\Lambda^\circ\subset\Xi^\circ$. We take the unique $\Xi\in\Sigma_\psi^0$ with $\Lambda^\circ\subset\Xi^\circ$.
Since $\mu\in\Theta^\circ\subset\Lambda^\circ\subset\Xi^\circ$, we can apply the same reasoning as in the case $\dim\Lambda=\dim\Sigma$ and we know that $\sigma^*(\psi)\in\mathcal{O}_{X,a}$ has normal crossings over $\bar{P}$ and $\Ord(\bar{P}, b_{E/N^*},\sigma^*(\psi))=0$ for any $E\in\mathcal{F}(\Delta)_1-\mathcal{F}(\Theta)_1$.

We consider the second case. Assume $\Lambda\in\bar{\Sigma}_\psi^1$.

Since $\Sigma_\psi$ is $H$-simple, there exist uniquely $\Xi_1\in\Sigma_\psi^0$, $\Xi_2\in\Sigma_\psi^0$ such that $\Lambda=\Xi_1\cap\Xi_2$, $\Lambda+(-H)=\Xi_1+(-H)$ and $\Lambda+H=\Xi_2+H$. We take $\Xi_1\in\Sigma_\psi^0$ and $\Xi_2\in\Sigma_\psi^0$ satisfying this condition. We take the unique $A_1\in\mathcal{V}(S_\psi)$ with $\Xi_1=\Delta(\{A_1\},S_\psi|V)$ and we take the unique $A_2\in\mathcal{V}(S_\psi)$ with $\Xi_2=\Delta(\{A_2\},S_\psi|V)$. $\langle b_{H/N^*}, A_1-A_2\rangle>0$. It follows $\langle b_{\Gamma/N^*}, A_1-A_2\rangle>0$, since $\Vect(\Lambda)+\Gamma=\Vect(\Lambda)+H$ and $\Lambda=\Xi_1\cap\Xi_2$. $\langle b_{E/N^*}, A_1\rangle =\langle b_{E/N^*}, A_2\rangle$ for any $E\in\mathcal{F}(\Theta)_1$, since $b_{E/N^*}\in E\subset\Theta\subset\Lambda=\Xi_1\cap\Xi_2$ for any $E\in\mathcal{F}(\Theta)_1$.
Since $\{0\}\neq\Theta\subset \Lambda\cap\Delta$, $\Gamma\subset\Delta\in\Sigma^{*0}$, $\Sigma^*$ is a subdivision of $\Sigma_\psi$ and $\Vect(\Lambda)+\Gamma=\Vect(\Lambda)+H$, we know $\Delta\subset\Xi_2$ and $\langle b_{E/N^*}, A_1-A_2\rangle\geq 0$ for any $E\in\mathcal{F}(\Delta)_1-\mathcal{F}(\Theta)_1$.
Note that $\{A_1,A_2\}\subset\mathcal{V}(S_\psi)\subset\Stab(S_\psi)\cap N=\Stab(S)\cap N $ and $\langle \theta, A_i\rangle\in\Z_0$ for any $\theta\in|\Sigma|\cap N^*$ and any $i\in\{1,2\}$. Let $F=\Conv(\{A_1,A_2\})$. $F\in\mathcal{F}(S_\psi)_1$. Recall $\mu\in\Theta^\circ\subset\Lambda^\circ$. $F=\Delta(\mu, S_\psi|V)$ and $\Lambda=\Delta(F,S_\psi|V)$. $\Ord(\bar{P},\mu,\sigma^*(\psi))=\Ord(P,\mu,\psi)=\langle \mu, A_1\rangle =\langle \mu, A_2\rangle$. We know that there exists uniquely a mapping $r:F\cap N\rightarrow k$ with $\In(P,\mu,\psi)=\Ps(P,F,\psi)=\sum_{B\in F\cap N}r(B)\hat{\pi}^*\chi(B)$. We take the mapping $r:F\cap N\rightarrow k$ satisfying this equality. We know $r(A_1)\neq 0$ and $r(A_2)\neq 0$. $\In(\bar{P},\mu,\sigma^*(\psi))=\sigma^*(\In(P,\mu,\psi))=\sigma^*(\sum_{B\in F\cap N}r(B)\hat{\pi}^*\chi(B))= \sum_{B\in F\cap N}r(B)\pi^*\chi(B)$. Now, for any $B\in F\cap N$, 
\begin{equation*}\begin{split}
&\chi(B)=\chi(\sum_{E\in\mathcal{F}(\Delta)_1}\langle b_{E/N^*}, B\rangle{b_{E/N^*}}^\vee_\Delta)=\prod_{E\in\mathcal{F}(\Delta)_1}\chi({b_{E/N^*}}^\vee_\Delta)^{\langle b_{E/N^*}, B\rangle}\\
&=\prod_{E\in\mathcal{F}(\Delta)_1-\mathcal{F}(\Theta)_1}((\chi({b_{E/N^*}}^\vee_\Delta)-c_E)+c_E) ^{\langle b_{E/N^*}, A_2\rangle}\prod_{E\in\mathcal{F}(\Theta)_1}\chi({b_{E/N^*}}^\vee_\Delta)^{\langle b_{E/N^*}, A_2\rangle}\\
&\qquad\qquad\qquad\qquad\qquad \prod_{E\in\mathcal{F}(\Delta)_1-\mathcal{F}(\Theta)_1}((\chi({b_{E/N^*}}^\vee_\Delta)-c_E)+c_E) ^{\langle b_{E/N^*}, B-A_2\rangle}.
\end{split}\end{equation*}

Since $c_E\neq 0$ for any $E\in\mathcal{F}(\Delta)_1-\mathcal{F}(\Theta)_1$, $\Gamma\in\mathcal{F}(\Delta)_1-\mathcal{F}(\Theta)_1$ and $\bar{P}=\{\pi^*\chi({b_{E/N^*}}^\vee_\Delta) -c_E| E\in\mathcal{F}(\Delta)_1-\mathcal{F}(\Theta)_1\}\cup\{\pi^*\chi({b_{E/N^*}}^\vee_\Delta) | E\in\mathcal{F}(\Theta)_1\}$ and $\langle b_{\Gamma/N^*},$\break$A_1\rangle>\langle b_{\Gamma/N^*},A_2\rangle$, we know
$\sum_{E\in\mathcal{F}(\Theta)_1}\langle b_{E/N^*}, A_1\rangle{b_{E/N^*}}^\vee_\Delta+\langle b_{\Gamma/N^*}, A_1\rangle{b_{\Gamma/N^*}}^\vee_\Delta \in\Supp(\bar{P},\In(\bar{P},\mu,\sigma^*(\psi)))$.

Let $\bar{A}=\sum_{ E\in\mathcal{F}(\Theta)_1}\langle b_{E/N^*}, A_1\rangle{b_{E/N^*}}^\vee_\Delta+\langle b_{\Gamma/N^*}, A_1\rangle{b_{\Gamma/N^*}}^\vee_\Delta \in (\Delta^\vee|V^*)\cap N$. $\bar{A}\in \Supp(\bar{P},\In(\bar{P},\mu,\sigma^*(\psi)))\subset\Supp(\bar{P},\sigma^*(\psi)) \subset\Gamma_+(\bar{P},\sigma^*(\psi))$.

$0\leq\Ord(\bar{P},b_{E/N^*},\sigma^*(\psi))\leq \langle b_{E/N^*},\bar{A}\rangle=0$ for any $E\in\mathcal{F}(\Delta)_1-(\mathcal{F}(\Theta)_1\cup\{\Gamma\})$. Thus, $\Ord(\bar{P},b_{E/N^*},\sigma^*(\psi))= \langle b_{E/N^*},\bar{A}\rangle=0$ for any $E\in\mathcal{F}(\Delta)_1-(\mathcal{F}(\Theta)_1\cup\{\Gamma\})$.

$\Ord(\bar{P},b_{E/N^*},\sigma^*(\psi))=\Ord(P, b_{E/N^*},\psi)= \langle b_{E/N^*},A_1\rangle=\langle b_{E/N^*},\bar{A}\rangle$ for any $E\in\mathcal{F}(\Theta)_1$, since $ b_{E/N^*}\in E\subset\Theta\subset\Lambda\subset\Xi_1=\Delta(\{A_1\},S_\psi|V)$ for any $E\in\mathcal{F}(\Theta)_1$.

Since $\Ord(\bar{P},b_{E/N^*},\sigma^*(\psi))= \langle b_{E/N^*},\bar{A}\rangle$ for any $E\in\mathcal{F}(\Delta)_1-\{\Gamma\}$, $\bar{A}\in\Gamma_+(\bar{P},$\break$\sigma^*(\psi))$ and $\bar{z}=\pi^*\chi({b_{\Gamma/N^*}}^\vee_\Delta)-c_\Gamma$, it follows that $\Gamma_+(\bar{P},\sigma^*(\psi))$ is of $\bar{z}$-Weirstrass type, and there exists uniquely an element $\hat{\psi}\in\mathcal{O}_{X,a}$ with 
\begin{equation*}\begin{split}
&\sigma^*(\psi)=\\
&\prod_{E\in\mathcal{F}(\Delta)_1-\mathcal{F}(\Theta)_1}((\chi({b_{E/N^*}}^\vee_\Delta)-c_E)+c_E) ^{\langle b_{E/N^*}, A_2\rangle}\prod_{E\in\mathcal{F}(\Theta)_1}\chi({b_{E/N^*}}^\vee_\Delta)^{\langle b_{E/N^*}, A_2\rangle}\hat{\psi}.
\end{split}\end{equation*}
We take the unique $\hat{\psi}\in\mathcal{O}_{X,a}$ satisfying this equality.
Note that $\prod_{E\in\mathcal{F}(\Delta)_1-\mathcal{F}(\Theta)_1}(($\break$\chi({b_{E/N^*}}^\vee_\Delta)-c_E)+c_E) ^{\langle b_{E/N^*}, A_2\rangle}\in\mathcal{O}_{X,a}^\times$ and $\Ord(\bar{P}, b_{E/N^*}, \prod_{E\in\mathcal{F}(\Theta)_1}\chi({b_{E/N^*}}^\vee_\Delta$\break$)^{\langle b_{E/N^*}, A_2\rangle})=\langle b_{E/N^*}, A_2\rangle=\langle b_{E/N^*}, A_1\rangle=\langle b_{E/N^*}, \bar{A}\rangle$ for any $E\in\mathcal{F}(\Theta)_1$.
\begin{equation*}\begin{split}
\In(\bar{P},\mu,\sigma^*(\psi))&= \prod_{E\in\mathcal{F}(\Delta)_1-\mathcal{F}(\Theta)_1}((\chi({b_{E/N^*}}^\vee_\Delta)-c_E)+c_E) ^{\langle b_{E/N^*}, A_2\rangle}\\
&\qquad\qquad\qquad\qquad \prod_{E\in\mathcal{F}(\Theta)_1}\chi({b_{E/N^*}}^\vee_\Delta)^{\langle b_{E/N^*}, A_2\rangle}\In(\bar{P},\mu,\hat{\psi}), \\
\In(\bar{P},\mu,\hat{\psi})&=\sum_{B\in F\cap N}r(B) \prod_{E\in\mathcal{F}(\Delta)_1-\mathcal{F}(\Theta)_1}((\chi({b_{E/N^*}}^\vee_\Delta)-c_E)+c_E) ^{\langle b_{E/N^*}, B-A_2\rangle}.
\end{split}\end{equation*}
It follows that $\Gamma_+(\bar{P},\hat{\psi})$ is of $\bar{z}$-Weierstrass type and $\langle b_{\Gamma/N^*}, A_1-A_2\rangle{b_{\Gamma/N^*}}^\vee_\Delta\in \Gamma_+(\bar{P},\In(\bar{P},\mu,\hat{\psi}))\subset\Gamma_+(\bar{P},\hat{\psi})$.

Let $h=\langle b_{\Gamma/N^*}, A_1-A_2\rangle\in\Z_+$, let
$$\hat{\psi}_0=\sum_{B\in F\cap N}r(B) \prod_{E\in\mathcal{F}(\Delta)_1-(\mathcal{F}(\Theta)_1\cup\{\Gamma\})}c_E ^{\langle b_{E/N^*}, B-A_2\rangle}(\bar{z}+c_\Gamma)^{\langle b_{\Gamma/N^*}, B-A_2\rangle }\in\mathcal{O}_{X,a},$$
and let $\delta_0=\sum_{E\in\mathcal{F}(\Delta)_1-\{\Gamma\}} b_{E/N^*}\in\Delta\cap N^*$.

$h=\langle b_{\Gamma/N^*}, A_1\rangle-\langle b_{\Gamma/N^*}, A_2\rangle
=\max\{\langle b_{\Gamma/N^*}, c\rangle|c\in F\}-\min\{\langle b_{\Gamma/N^*}, c\rangle|c\in F\}$. $\hat{\psi}_0$ is a polynomial with coefficients in $k$ with variable $\bar{z}$ of degree $h$. $\hat{\psi}_0=\In(\bar{P},\delta_0,\In(\bar{P},\mu,\hat{\psi}))=\In(\bar{P},\delta_0,\hat{\psi})$.
$0\leq\Ord(\hat{P},b_{\Gamma/N^*},\hat{\psi})\leq\Ord(\hat{P},b_{\Gamma/N^*},\hat{\psi}_0)\leq h$. $\{\Ord(\hat{P},b_{\Gamma/N^*},\hat{\psi}_0){b_{\Gamma/N^*}}^\vee_\Delta\}$ is the $\bar{z}$-top vertex of $\Gamma_+(\bar{P},\hat{\psi})$.
$\Ht(\bar{z},\Gamma_+(\bar{P},$\break$\sigma^*(\psi)))=\Ht(\bar{z},\Gamma_+(\bar{P},\hat{\psi}))=\Ord(\hat{P},b_{\Gamma/N^*},\hat{\psi}_0)- \Ord(\hat{P},b_{\Gamma/N^*},\hat{\psi})\leq h$. $\Ht(\bar{z},$\break$\Gamma_+(\bar{P},\sigma^*(\psi)))=h$, if and only if, $\Ord(\hat{P},b_{\Gamma/N^*},\hat{\psi})=0$ and $\Ord(\hat{P},b_{\Gamma/N^*},\hat{\psi}_0)= h$.
By Theorem~\ref{important}.$12.(b)$ we know $h\leq\Ht(H,S_\psi)=\Inv(P,z,\phi)$.
Let $\bar{\psi}\in\mathcal{O}_{X,a}$ be any main factor of $(\bar{P},\bar{z},\sigma^*(\psi))$. $\bar{\psi}$ is a main factor of $(\bar{P},\bar{z},\sigma^*(\phi))$. It follows $\Inv(\bar{P},\bar{z},\sigma^*(\phi))=\Ht(\bar{z}, \Gamma_+(\bar{P},\bar{\psi}))\leq \Ht(\bar{z},\Gamma_+(\bar{P},\sigma^*(\psi)))$.

It follows $\Inv(\bar{P},\bar{z},\sigma^*(\phi))\leq \Inv(P,z,\phi)$.

Assume that $\Inv(\bar{P},\bar{z},\sigma^*(\phi))= \Inv(P,z,\phi)$. It follows $\Ht(\bar{z},\Gamma_+(\bar{P},\sigma^*(\psi)))= \Ord(\hat{P},b_{\Gamma/N^*},\hat{\psi}_0)=\Ht(H,S_\psi)=h$. By Theorem~\ref{important}.$12.(c)$ we know that $c(S_\psi)=2$ and any structure constant of $\Sigma_\psi$ is an integer. Since $c(S_\psi)=2$, $\{A_1,A_2\}\subset\mathcal{V}(S_\psi)$, $\langle f^{P\vee}_z, A_1-A_2\rangle >0$, $\{A_1,A_2\}=\mathcal{V}(S_\psi)$, $A_1$ is the $z$-top vertex, and $A_2$ is the $z$-bottom vertex. Since $x$ does not divide $\psi$ for any $x\in P$, we know $A_1=hf^P_z$, and $\langle f^{P\vee}_z, A_2\rangle=0$, $\Ht(H,S_\psi)=h=\langle b_{\Gamma/N^*}, A_1-A_2\rangle$. Since any structure constant of $\Sigma_\psi$ is an integer, from Lemma~\ref{propsimple}.$18.(f)$ it follows that $A_2/h=(A_2-A_1)/h+f^P_z\in N\cap \Stab(S)\cap(\Vect(H)^\vee|V^*)$. Since $S_\psi$ is of $H$-Weierstrass type and $\Inv(P,z,\phi)>0$, $A_2 \neq 0$.

Let $B=A_2/h $, let $C=(A_1-A_2)/h\in N$ and let $$c=\prod_{E\in\mathcal{F}(\Delta)_1-(\mathcal{F}(\Theta)_1\cup\{\Gamma\})}c_E^{\langle b_{E/N^*}, C\rangle}\in k-\{0\}.$$ $B\in N\cap \Stab(S)\cap(\Vect(H)^\vee|V^*)-\{0\}$,
$C=-B+f^P_z\in N$.
$F\cap N=\{A_2+iC|i\in\{0,1,\ldots,h\}\}$. 
$A_1=A_2+hC$. $\langle f^{P\vee}_z, C\rangle=\langle b_{\Gamma/N^*}, C\rangle=1$.
$ A_2+iC=(h-i)B+if^P_z$ for any $ i\in\{0,1,\ldots,h\}$.
$$\hat{\psi}_0=\sum_{i=0}^hr(A_2+iC)c^i(\bar{z}+c_\Gamma)^i.$$
Since $\Ord(\hat{P},b_{\Gamma/N^*},\hat{\psi}_0)= h$,  $\hat{\psi}_0=r(A_1)c^h\bar{z}^h$.
$\sum_{i=0}^hr(A_2+iC)c^i(\bar{z}+c_\Gamma)^i
=\hat{\psi}_0=r(A_1)c^h\bar{z}^h=\sum_{i=0}^h r(A_1)c^h \binom{h}{i}(-c_\Gamma)^{h-i}(\bar{z}+c_\Gamma)^i.$
We know that $r(A_2+iC)c^i= r(A_1)c^h \binom{h}{i}(-c_\Gamma)^{h-i}$ and $ r(A_2+iC)= r(A_1) \binom{h}{i}(-c c_\Gamma) ^{h-i}$ for any $i\in\{0,1,\ldots,h\}$.
Thus, 
\begin{equation*}\begin{split}
&\sigma^*(\In(P,\mu,\psi))=\In(\bar{P},\mu,\sigma^*(\psi))=\sum_{B\in F\cap N}r(B)\pi^*\chi(B)=\\
&\sum_{i=0}^hr(A_2+iC) \pi^*\chi(A_2+iC)
=\sum_{i=0}^h r(A_1) \binom{h}{i}(-c c_\Gamma) ^{h-i}\pi^*\chi((h-i)B+if^P_z)=\\
&r(A_1) (\sum_{i=0}^h \binom{h}{i}(-c c_\Gamma) ^{h-i} \pi^*\chi(B)^ {h-i})\pi^*\chi(f^P_z )^i=\\
&\sigma^* (r(A_1) ((\hat{\pi}^*\chi(f^P_z )- c c_\Gamma\hat{\pi}^*\chi(B))^h)
= \sigma^* (r(A_1) (z- c c_\Gamma\hat{\pi}^*\chi(B))^h),
\end{split}\end{equation*}
and $\Ps(P,\Delta(\mu, S_\psi|V),\psi)=\In(P,\mu,\psi)= r(A_1) (z- c c_\Gamma\hat{\pi}^*\chi(B))^h$.
$-c c_\Gamma\hat{\pi}^*\chi(B)\in M(R')$,
where $R'$ denotes the localization of $k[P-\{z\}]$ by the maximal ideal $k[P-\{z\}]\cap M(R)=(P-\{z\}) k[P-\{z\}]$.
 $-c c_\Gamma\hat{\pi}^*\chi(B)\neq 0$.
We know that the face $\Delta(\mu, S_\psi|V)$ of $S_\psi$ is a $z$-removable face, which contradicts that $S_\psi$ has no $z$-removable faces.

We conclude  $\Inv(\bar{P},\bar{z},\sigma^*(\phi))< \Inv(P,z,\phi)$.

We conclude that Theorem~\ref{main} holds, if $\dim\Lambda=\dim \Sigma-1$ and $\Inv(P,s,\phi)>0$.

We conclude that Theorem~\ref{main} holds in all cases.

\section{Proof of the submain theorems}
\label{submain proofs}
We give the proof of our submain theorems Theorem~\ref{erase faces}, Theorem~\ref{make simple}, Lemma~\ref{make Weierstrass type}, Corollary~\ref{resolution game} and Corollary~\ref{local uniformization}.

Let $k$ be any \emph{algebraically closed} field, let $R$ be any regular local ring such that $R$ contains $k$ as a subring, the residue field $R/M(R)$ is isomorphic to $k$ as $k$-algebras, $R$ is a localization of a finitely generated $k$-algebra and $\dim R\geq 2$, let $P$ be any parameter system of $R$, and let $z\in P$ be any element.

Let $R'$ denote the localization of $k[P - \{z\}]$ by the maximal ideal $k[P - \{z\}]\cap M(R)=(P-\{z\})k[P-\{z\}]$. The ring $R'$ is a regular local subring of $R$. The set $P - \{z\}$ is a parameter system of $R'$.

We give the proof of Theorem~\ref{erase faces}.

Assume $\dim R\geq 2$.

Consider any element $w\in M(R^c)$ with $\partial w/\partial z\in R^{c\times}$.
We denote $P_w=\{w\}\cup(P-\{z\})$. (Lemma~\ref{coordinate change}.)

The bijection $P_w\rightarrow P$ sending $w\in P_w$ to $z\in P$ and sending any $x\in P_w-\{w\}=P-\{z\}$ to $x\in P-\{z\}$ itself induces an isomorphism $\Map(P,\R)\rightarrow\Map(P_w,\R)$ of vector spaces over $\R$.
By this isomorphism we identify $\Map(P,\R)$ and $\Map(P_w,\R)$.

Consider any element $\psi\in R$ such that $\psi\neq 0$, $\Gamma_+(P,\psi)$ is of $z$-Weierstrass type and any $x\in P-\{z\}$ does not divide $\psi$. We take the unique non-negative integer $h$ such that $\{hf^P_z\}$ is the unique $z$-top vertex of $\Gamma_+(P,\psi)$.

Recall that $\Gamma_+(P,\psi)\subset\Map(P,\R)$ and $\{f^P_x|x\in P\}$ is an $\R$-basis of the vector space $\Map(P,\R)$. Let $U=\{a\in \Map(P,\R)|\langle f^{P\vee}_z,a\rangle<h\}$ and $V=\{a\in \Map(P,\R)|\langle f^{P\vee}_z,a\rangle=0\}$. We put $\rho(a)=(a-\langle f^{P\vee}_z,a\rangle f^P_z)/(h-\langle f^{P\vee}_z,a\rangle)\in V$ for any $a\in U$ and we define a mapping $\rho:U\rightarrow V$.

Note that $V$ is an $\R$-vector subspace of $\Map(P,\R)$ with $\dim V=\dim\Map(P,$\break$\R)-1$ and the set $\{f^P_x|x\in P-\{z\}\}$ is an $\R$-basis of $V$. Using  the isomorphism $\Map(P-\{z\},\R)\rightarrow V$ of vector spaces over $\R$ sending $f^{P-\{z\}}_x\in\Map(P-\{z\},\R)$ to $f^P_x\in V$ for any $x\in P-\{z\}$, we identify $\Map(P-\{z\},\R)$ and $V$. 

We identify the dual vector space $V^*$ of $V$ with the vector subspace $\{\bar{\omega}\in\Map(P,\R)^*|\langle \bar{\omega}, f^P_z\rangle=0\}$ in the dual vector space $\Map(P,\R)^*$ of $\Map(P,\R)$. Under this identification $(\Map(P,\R_0)\cap V)^\vee|V=(\Map(P,\R_0)^\vee|\Map(P,\R))\cap V^*$.

\noindent $1.$\quad It follows from Lemma~\ref{first correspondence}.$2$.

\noindent $2.$\quad
We take the unique pair of an invertible element $\bar{u}\in R^{c\times}$, and a $z$-Weierstrass polynomial $\bar{\psi}\in R^c$ over $P$ with $\psi=\bar{u}\bar{\psi}$.

Consider any $\chi\in M(R^{\prime c})$. The following lemma holds. See Lemma~\ref{first correspondence} and  Proposition~\ref{miracle}:
\begin{lemma}
\label{vertex projection}
\begin{enumerate}
\item
$\Stab(\Gamma_+(P_{z+\chi},\psi))=\Map(P,\R_0)$.
\item
$\Gamma_+( P_{z+\chi},\psi)$ is of $(z+\chi)$-Weierstrass type. The unique $(z+\chi)$-top vertex of $\Gamma_+(P_{z+\chi},\psi)$ is equal to $\{hf^P_z\}$.
\item
$\rho(\Gamma_+(P_{z+\chi},\psi)\cap U)=\emptyset\Leftrightarrow \Gamma_+(P_{z+\chi},\psi)\cap U=\emptyset\Leftrightarrow \Ord(P_{z+\chi}, f^{P\vee}_z, \psi)\geq h\Leftrightarrow \Ord(P_{z+\chi}, f^{P\vee}_z, \psi)= h\Leftrightarrow \bar{\psi}=(z+\chi)^h$. 
\end{enumerate}

Below, we denote $$\bar{\Gamma}_+( P_{z+\chi},\psi)= \rho(\Gamma_+(P_{z+\chi},\psi) \cap U),$$
and we assume $\bar{\Gamma}_+( P_{z+\chi},\psi)\neq\emptyset$.
\begin{enumerate}
\setcounter{enumi}{3}
\item
$\bar{\Gamma}_+(P_{z+\chi},\psi)=
\rho((\Convcone(\Gamma_+(P_{z+\chi},\psi)+\{-hf^P_z\})+\{ hf^P_z \})\cap U)$ and\hfill\break
$\bar{\Gamma}_+( P_{z+\chi},\psi)$ is a rational convex pseudo polytope over the lattice $\Map(P,\Z)\cap V$ in $V$ with $\Stab(\bar{\Gamma}_+(P_{z+\chi},\psi))= \Map(P,\R_0)\cap V$.

$0\not\in\bar{\Gamma}_+(P_{z+\chi},\psi)\subset\Map(P,\R_0)\cap V$.
\item
Consider any face $F$ of $\Gamma_+(P_{z+\chi},\psi)$ satisfying $hf^P_z\in F$ and
$F\cap U\neq \emptyset$.

$\rho(F\cap U)$ is a face of $\bar{\Gamma}_+( P_{z+\chi},\psi)$.

$\dim F\geq 1$ and $\dim \rho(F\cap U)=\dim F-1$.

If $\dim F=1$, then $\Stab(F)=\{0\}$.

Consider any $\omega\in\Map(P,\R_0)^\vee$ with $F=\Delta(\omega, \Gamma_+(P_{z+\chi},\psi))$.
If we take the unique pair of elements $\bar{\omega}\in \Map(P,\R_0)^\vee\cap V^*$ and $t\in\R_0$ satisfying $\omega=\bar{\omega}+tf^{P\vee}_z$, then
$\rho(F\cap U)=\Delta(\bar{\omega}, \bar{\Gamma}_+(P_{z+\chi},\psi))$ and $t=\Ord(\bar{\omega},\bar{\Gamma}_+(P_{z+\chi},\psi))$.

Consider any $\bar{\omega}\in \Map(P,\R_0)^\vee\cap V^*$ with $\rho(F\cap U)=\Delta(\bar{\omega}, \bar{\Gamma}_+(P_{z+\chi},\psi))$.
If $t=\Ord(\bar{\omega},\bar{\Gamma}_+( P_{z+\chi},\psi)) \in\R_0$, then $\bar{\omega}+tf^{P\vee}_z\in \Map(P,\R_0)$ and $F=\Delta(\bar{\omega}+tf^{P\vee}_z , \Gamma_+( P_{z+\chi},\psi))$.
\item
The mapping from the set of all faces of $\Gamma_+(P_{z+\chi},\psi)$ satisfying $hf^P_z\in F$ and
$F\cap U\neq \emptyset$ to the set of all faces of $\bar{\Gamma}_+( P_{z+\chi},\psi)$ sending $F$ to $\rho(F\cap U)$ is bijective and it preserves the inclusion relation.
\item
Consider any face $F$ of $\Gamma_+(P_{z+\chi},\psi)$.

$F$ is a $(z+\chi)$-removable face, if and only if, $hf^P_z\in F$ and there exists $\bar{\chi}\in M(R^{\prime c})$ satisfying $\Ps(P_{z+\chi}, F,\bar{\psi})= (z+\chi+\bar{\chi})^h$ and $\bar{\chi}\neq 0$.
\item
Any $(z+\chi)$-removable face $F$ of $\Gamma_+(P_{z+\chi},\psi)$ satisfies $hf^P_z\in F$, 
and $F\cap U\supset F\cap V\neq \emptyset$.

\item
If $\Gamma_+(P_{z+\chi},\psi)$ has a $(z+\chi)$-removable face, then it has a $(z+\chi)$-removable face of dimension one.
\item
Assume that $\Gamma_+(P_{z+\chi},\psi)$ has a $(z+\chi)$-removable face and
consider any $(z+\chi)$-removable face $F$ of dimension one of $\Gamma_+(P_{z+\chi},\psi)$.
\begin{enumerate}
\item
$\rho(F\cap U)$ is a vertex of $\bar{\Gamma}_+( P_{z+\chi},\psi)$.
\end{enumerate}

We take the unique point $c(F)\in\bar{\Gamma}_+( P_{z+\chi},\psi)$ with $\{c(F)\}=\rho(F\cap U)$.
\begin{enumerate}
\setcounter{enumii}{1}
\item
$\{c(F)\}$ is a vertex of $\bar{\Gamma}_+( P_{z+\chi},\psi)$. $c(F)\in\Map(P,\Z_0)\cap V -\{0\}$.
\item
There exists uniquely an element $\gamma(F)\in k-\{0\}$ satisfying
$\Ps(P_{z+\chi}, F,$\break$\bar{\psi})= (z+\chi+\gamma(F)\prod_{x\in P-\{z\}}x^{\langle f^{P\vee}_x,c(F)\rangle})^h$.
\end{enumerate}

We take the unique element $\gamma(F)\in k-\{0\}$ satisfying
$\Ps(P_{z+\chi}, F,\bar{\psi})=(z+\chi+\gamma(F)\prod_{x\in P-\{z\}}x^{\langle f^{P\vee}_x,c(F)\rangle})^h$ and we denote
\begin{equation*}\begin{split}
\chi(F)&= \gamma(F)\prod_{x\in P-\{z\}}x^{\langle f^{P\vee}_x,c(F)\rangle}\in M(R^{\prime c})-\{0\}\text{, and}\\
\bar{\delta}_0&=\sum_{x\in P-\{z\}}f^{P\vee}_x=b_{(\Map(P,\R_0)^\vee\cap V^*) /\Map(P,\Z)^*}\\
&\qquad\qquad\in(\Map(P,\R_0)^\vee\cap V^*)^\circ\cap \Map(P,\Z)^*.
\end{split}\end{equation*}

\begin{enumerate}
\setcounter{enumii}{3}
\item
$\Ps(P_{z+\chi}, F,\bar{\psi})= (z+\chi+\chi(F))^h$.
\item

$c(F)\not\in\bar{\Gamma}_+( P_{z+\chi+\chi(F)},\psi)\subset\bar{\Gamma}_+(P_{z+\chi},\psi)$.
Any face $\bar{G}$ of $\bar{\Gamma}_+(P_{z+\chi},\psi)$ with $c(F)\not\in\bar{G}$ is a face of $\bar{\Gamma}_+( P_{z+\chi+\chi(F)},\psi)$.

$\Ord(P-\{z\}, \bar{\delta}_0, \chi(F))=\langle  \bar{\delta}_0, c(F)\rangle$.
$\In(P-\{z\}, \bar{\delta}_0, \chi(F))=\chi(F)$. $\Supp($\break$P-\{z\}, \chi(F))=\{c(F)\}$.
\end{enumerate}
\end{enumerate}
\end{lemma}

Below, we use the above notations $\bar{\Gamma}_+(P_{z+\chi},\psi)= \sigma(\Gamma_+(P_{z+\chi},\psi)\cap U)$, $c(F)\in(\Map(P,\Z_0)\cap V)-\{0\}$, $\gamma(F)\in k-\{0\}$, $\chi(F)\in M(R^{\prime c})-\{0\}$ and $\bar{\delta}_0\in(\Map(P,\R_0)^\vee$\break$\cap V^*)^\circ \cap \Map(P,\Z)^*$. For any $w\in M(R^c)$ with $\partial w/\partial z\in R^{c\times}$ we denote $\bar{\Gamma}_+(P_w,\psi)= \sigma(\Gamma_+(P_w,\psi)\cap U)$.

We consider the following algorithm starting from Step 0.

In Step 0 we put $\bar{\chi}(0)=0\in M(R^{\prime c})$ and proceed to Step 1. 

Consider any positive integer $i$. In Step $i$, if $\Gamma_+(P_{z+\bar{\chi} (i-1)},\psi)$ has no $(z+\bar{\chi} (i-1))$-removable faces, then we finish the algorithm. In Step $i$, if $\Gamma_+( P_{z+\bar{\chi} (i-1)},\psi)$ has $(z+\bar{\chi} (i-1))$-removable faces, then we choose any $(z+\bar{\chi} (i-1))$-removable face $F(i)$ of dimension one of $\Gamma_+( P_{z+\bar{\chi} (i-1)},\psi)$ satisfying
$\langle \bar{\delta}_0,c(F(i))\rangle=\min\{\langle \bar{\delta}_0,c(F)\rangle|
F$ is a $(z+\bar{\chi} (i-1))$-removable face of dimension one of $\Gamma_+( P_{z+\bar{\chi} (i-1)},\psi)\}$, we put $\bar{\chi}(i)= \bar{\chi}(i-1)+\chi(F(i))\in M(R^{\prime c})$ and we proceed to Step $i+1$.

Consider the case where we finish this algorithm in finite steps.
Assume that the algorithm has finished in Step $i$ for some positive integer $i$.
$\bar{\chi}(i-1)\in M(R^{\prime c})$ and $\Gamma_+( P_{z+\bar{\chi} (i-1)},\psi)$ has no $(z+\bar{\chi} (i-1))$-removable faces.
By Lemma~\ref{vertex projection}.10.(e) above we know $\Supp(P-\{z\}, \bar{\chi}(i-1))\subset\{ c(F(j))| j\in\{1,2,\ldots,i-1\}\}\subset\bar{\Gamma}_+(P,\psi)-\bar{\Gamma}_+( P_{z+\bar{\chi} (i-1)},\psi)$.
We conclude that there exists $\chi_0\in M(R^{\prime c})$ such that $\Gamma_+(P_{z+\chi_0},$\break$\psi)$ has no $(z+\chi_0)$-removable faces and $\Supp(P-\{z\}, \chi_0)\subset\bar{\Gamma}_+(P,\psi)-\bar{\Gamma}_+( P_{z+\chi_0},\psi)$.

Consider the case where this algorithm has infinite steps.
By Lemma~\ref{vertex projection}.10.(e) and 10.(b) we know that
$\chi(F(i))\in M(R^{\prime c})$, $\Ord(P-\{z\},\bar{\delta}_0, \chi(F(i)))
=\langle \bar{\delta}_0,c(F(i))\rangle$, 
$c(F(i))\in \Map(P,\Z_0)\cap V-\{0\}$, and
$\bar{\chi}(i)=\sum_{j=1}^i\chi(F(j))\in M(R^{\prime c})$ for any $i\in \Z_+$.
By Lemma~\ref{vertex projection}.10.(e) we know that
$c(F(i))\neq c(F(j))$ for any $i\in \Z_+$ and any $j\in \Z_+$ with $i\neq j$.

Since $\{e\in \Map(P,\Z_0)\cap V |\langle \bar{\delta}_0,e\rangle\leq m\}$ is a finite set for any $m\in\Z_0$, we know that $\lim_{i\rightarrow\infty}\langle \bar{\delta}_0,c(F(i))\rangle=\infty$ and the sequence $\bar{\chi}(i)$, $i\in\Z_+$ converges. We put $\chi_0=\lim_{i\rightarrow\infty}\bar{\chi}(i)=
\sum_{i=1}^\infty\chi(F(i))\in M(R^{\prime c})$.

By Lemma~\ref{vertex projection}.10.(e) we know $\Supp(P-\{z\}, \chi_0)\subset\{c(F(i))| i\in\Z_+\}\subset\bar{\Gamma}_+(P,$\break$\psi)-(\cap_{i\in\Z_+}\bar{\Gamma}_+(P_{z+\bar{\chi}(i)}, \psi))$.

Assume that $\Gamma_+(P_{z+\chi_0},\psi)$ has $(z+\chi_0)$-removable faces.
We will deduce a contradiction. $\bar{\Gamma}_+(P_{z+\chi_0},\psi)\neq\emptyset$. Take any $(z+\chi_0)$-removable face $F$ of dimension one of $\Gamma_+(P_{z+\chi_0},\psi)$. 
Take any $\bar{\omega}\in(\Map(P,\R_0)^\vee\cap V^*)^\circ$ with $\{c(F)\}=\Delta(\bar{\omega},\bar{\Gamma}_+(P_{z+\chi_0},\psi))$.
Put $t=\Ord(\bar{\omega}, \bar{\Gamma}_+(P_{z+\chi_0},\psi))\in\R_0$ and $\omega=\bar{\omega}+tf^{P\vee}_z\in \Map(P,\R_0)^\vee$.
We know that $F=\Delta(\omega, \Gamma_+(P_{z+\chi_0},\psi))$, $t=\langle \bar{\omega},c(F)\rangle$ and $\In(P_{z+\chi_0},$\break$ \omega,\bar{\psi})=\Ps(P_{z+\chi_0}, F,\bar{\psi})=(z+\chi_0+\chi(F))^h$.

Since $\bar{\omega}\in(\Map(P,\R_0)^\vee\cap V^*)^\circ$,
$\{e\in \Map(P,\Z_0)\cap V|\langle \bar{\omega},e\rangle\leq t \}$ is a finite set.
Take any $i\in\Z_+$ such that $\langle \bar{\delta}_0,c(F(j))\rangle>\langle \bar{\delta}_0,c(F)\rangle$ and $\langle \bar{\omega},c(F(j))\rangle>t$ for any $j\in \Z_+$ with $j\geq i$.

$\Ord(P_{z+\chi_0},\omega,z+\chi_0)=\langle \omega, f^P_z\rangle =t$.
$(z+\chi_0)-(z+\bar{\chi}(i))=\sum_{j=i+1}^\infty \chi(F(j))$.
For any $j\in\Z_+$ with $j\geq i+1$, $\Ord(P_{z+\chi_0},\omega, \chi(F(j)))=\Ord(P-\{z\},\bar{\omega},\chi(F(j)))=\langle \bar{\omega},c(F(j))\rangle>t$.
Therefore, $\Ord(P_{z+\chi_0},\omega, (z+\chi_0)-(z+\bar{\chi}(i)))= \Ord(P_{z+\chi_0},\omega, $\break$\sum_{j=i+1}^\infty \chi(F(j)))>t$ and
$\Ord(P_{z+\chi_0},\omega,z+\chi_0)= \Ord(P_{z+\chi_0},\omega,z+\chi(i))$.
We know that $\Ord(P_{z+\chi_0},\omega,\zeta)= \Ord(P_{z+\bar{\chi}(i)},\omega,\zeta)$ and $\lambda(\In(P_{z+\chi_0},\omega,\zeta))= \In(P_{z+\bar{\chi}(i)},\omega,\zeta)$ for any $\zeta\in R^c$, where $\lambda:R^c\rightarrow R^c$ is the isomorphism of $k$-algebras satisfying $\lambda(x)=x$ for any $x\in P-\{z\}$ and $\lambda(z+\chi_0)=\In(P_{z+\bar{\chi}(i)},\omega, z+\chi_0)=z+\bar{\chi}(i)$. Thus, $\In(P_{z+\bar{\chi}(i)}, \omega,\bar{\psi})=(z+\bar{\chi}(i)+\chi(F))^h$,
$F=h\Conv(\{f^P_z,c(F)\})=\Conv(\Supp(P_{z+\bar{\chi}(i)}, \In(P_{z+\bar{\chi}(i)}, \omega,\psi)))=\Delta(\omega, \Gamma_+(P_{z+\bar{\chi}(i)},\psi))\in \mathcal{F}(\Gamma_+(P_{z+\bar{\chi}(i)},\psi))$ and $F$ is a $(z+\bar{\chi}(i))$-removable face of dimension one of $\Gamma_+(P_{z+\bar{\chi}(i)},\psi)$.
Since $\langle \bar{\delta}_0,c(F(i+1))\rangle=\min\{\langle \bar{\delta}_0,c(F)\rangle|
F$ is a $(z+\bar{\chi} (i))$-removable face of dimension one of $\Gamma_+( P_{z+\bar{\chi} (i)},\psi)\}$, we have
$\langle \bar{\delta}_0,c(F(i+1))\rangle\leq\langle \bar{\delta}_0,c(F)\rangle$.

Since $\langle \bar{\delta}_0,c(F(i+1))\rangle>\langle \bar{\delta}_0,c(F)\rangle$, we get a contradiction. 

We conclude that $\Gamma_+(P_{z+\chi_0},\psi)$ has no $(z+\chi_0)$-removable faces.

By similar reasoning we know that for any vertex $\bar{G}$ of $\bar{\Gamma}_+(P_{z+\chi_0},\psi)$, there exists $i\in\Z_+$ such that $\bar{G}$ is a vertex of $\bar{\Gamma}_+(P_{z+\bar{\chi}(j)},\psi)$ for any $j\in\Z_+$ with $j\geq i$ and $\bar{G}$ is a vertex of $\cap_{i\in\Z_+}\bar{\Gamma}_+(P_{z+\bar{\chi}(i)},\psi)$. Thus, $\bar{\Gamma}_+(P_{z+\chi_0},\psi)= \cap_{i\in\Z_+}\bar{\Gamma}_+(P_{z+\bar{\chi}(i)},\psi)$ and 
$\Supp(P-\{z\}, \chi_0)\subset\bar{\Gamma}_+(P,\psi)- \bar{\Gamma}_+(P_{z+\chi_0},\psi)$.

We know that there exists $\chi_0\in M(R^{\prime c})$ such that $\Gamma_+(P_{z+\chi_0},\psi)$ has no $(z+\chi_0)$-removable faces and $\Supp(P-\{z\},\chi_0)\subset\bar{\Gamma}_+(P,\psi)- \bar{\Gamma}_+(P_{z+\chi_0},\psi)$ in all cases.

We will show the uniqueness of $\chi_0$ after the proof of claim $3$ and claim $4$.

\text{}\newline

Below, we assume $\chi_0\in M(R^{\prime c})$, $\Gamma_+(P_{z+\chi_0},\psi)$ has no $(z+\chi_0)$-removable faces, and $\Supp(P-\{z\}, \chi_0)\subset \bar{\Gamma}_+(P,\psi)- \bar{\Gamma}_+(P_{z+\chi_0},\psi)$ and that $v\in R^{c\times}$, $\mu\in M(R^{\prime c})$ and $w=v(z+\mu)$. It follows $\Gamma_+(P_w,\psi)=\Gamma_+(P_{z+\mu},\psi)$.

\noindent $3.$\quad
Consider the case where  $\bar{\Gamma}_+(P_{z+\chi_0},\psi)= \bar{\Gamma}_+(P_{z+\mu},\psi)=\emptyset$. $\bar{\psi}= (z+\chi_0)^h=(z+\mu)^h$ by Lemma~\ref{vertex projection}.3, $\mu-\chi_0=0$ and $\Gamma_+(P-\{z\}, \mu-\chi_0)=\emptyset$. $\bar{\Gamma}_+(P_w,\psi)=\bar{\Gamma}_+(P_{z+\mu_0},\psi)=\emptyset=\Conv(\bar{\Gamma}_+(P_{z+\chi_0},\psi)\cup\Gamma_+(P-\{z\}, \mu-\chi_0))$.

Below, we assume either  $\bar{\Gamma}_+(P_{z+\chi_0},\psi)\neq\emptyset$ or $\bar{\Gamma}_+(P_{z+\mu},\psi)\neq\emptyset$.

Consider any $\bar{\omega}\in\Map(P,\R_0)^\vee\cap V^*$ and put \begin{equation*}\begin{split}
t&=\min\{\Ord(\bar{\omega},\Gamma_+(P-\{z\},\mu-\chi_0)), \Ord(\bar{\omega}, \bar{\Gamma}_+(P_{z+\chi_0},\psi)),\\
&\qquad\qquad\Ord(\bar{\omega}, \bar{\Gamma}_+(P_{z+\mu},\psi))\}\in\R_0,
\end{split}\end{equation*}
and put $\omega=\bar{\omega}+tf^{P}_z\in \Map(P,\R_0)$.

Assume that $t=\Ord(\bar{\omega},\Gamma_+(P-\{z\},\mu-\chi_0))< \min\{\Ord(\bar{\omega}, \bar{\Gamma}_+(P_{z+\chi_0},\psi)), \Ord(\bar{\omega}, $\break$\bar{\Gamma}_+(P_{z+\mu},\psi))\}$. We will deduce a contradiction.

It follows $\Ord(P_{z+\chi_0},\omega,\bar{\psi})=ht$, $\In(P_{z+\chi_0},\omega,\bar{\psi})=(z+\chi_0)^h$,
$\Ord(P_{z+\mu},\omega,\bar{\psi})=ht$, and $\In(P_{z+\mu},\omega,\bar{\psi})=(z+\mu)^h$.

We have $t=\Ord(P_{z+\chi_0},\omega, z+\chi_0)= \Ord(P_{z+\mu},\omega, z+\mu)$. 
On the other hand, $t=\Ord(\bar{\omega},\Gamma_+(P-\{z\},\mu-\chi_0))=\Ord(P-\{z\},\bar{\omega}, \chi_0-\mu)= \Ord(P_{z+\chi_0},\omega, \mu-\chi_0)$. Thus
$\Ord(P_{z+\chi_0},\omega, z+\mu)=\Ord(P_{z+\chi_0},\omega, z+\chi_0+(\mu-\chi_0))=t$ and $\In(P_{z+\chi_0},\omega, z+\mu)= z+\chi_0+\In(P-\{z\},\bar{\omega},\mu-\chi_0)$.
Since $\In(P_{z+\mu},\omega,\bar{\psi})=(z+\mu)^h=(z+\chi_0+(\mu-\chi_0))^h$, we have $(z+\chi_0)^h=\In(P_{z+\chi_0},\omega,\bar{\psi})= ( z+\chi_0+\In(P-\{z\},\bar{\omega},\mu-\chi_0))^h$, and $\In(P-\{z\},\bar{\omega},\mu-\chi_0)=0$. On the other hand, $\R_0\ni t=\Ord(P-\{z\},\bar{\omega}, \mu-\chi_0)$, $\In(P-\{z\},\bar{\omega},\mu-\chi_0)\neq 0$ and  we get a contradiction.

We know $t=\min\{\Ord(\bar{\omega}, \bar{\Gamma}_+(P_{z+\chi_0},\psi)), \Ord(\bar{\omega}, \bar{\Gamma}_+(P_{z+\mu},\psi))\}\leq\Ord(\bar{\omega},$\break$\Gamma_+(P-\{z\},\mu-\chi_0))$.

Assume that $t=\Ord(\bar{\omega}, \bar{\Gamma}_+(P_{z+\chi_0},\psi))< \Ord(\bar{\omega}, \bar{\Gamma}_+(P_{z+\mu},\psi))$. We will deduce a contradiction. 

Let $G=\Delta(\omega, \Gamma_+(P_{z+\chi_0},\psi))$. $G$ is a face of $\Gamma_+(P_{z+\chi_0},\psi)$ which is not $(z+\chi_0)$-removable and $G\cap U\neq\emptyset$. 
Sincet $t<\Ord(\bar{\omega}, \bar{\Gamma}_+(P_{z+\mu},\psi))$, we have $\Ord(P_{z+\mu},\omega, \bar{\psi})=ht$, and $\In(P_{z+\mu},\omega, \bar{\psi})=(z+\mu)^h$. Note that $t\leq \Ord(\bar{\omega},\Gamma_+(P-\{z\},\mu-\chi_0))=\Ord(P_{z+\mu_0},\omega, \mu-\chi_0)$. We put
$$\bar{\mu}=
\begin{cases}
\In(P-\{z\},\bar{\omega}, \mu-\chi_0)& \text{if $t=\Ord(P-\{z\},\bar{\omega},\mu-\chi_0)$},\\
0&\text{if $t\neq\Ord(P-\{z\},\bar{\omega},\mu-\chi_0)$}.
\end{cases}.$$
$\bar{\mu}\in M(R^{\prime c})$.
$\In(P_{z+\chi_0},\omega, z+\mu)= \In(P_{z+\chi_0},\omega, z+\chi_0+(\mu-\chi_0))=z+ \chi_0+\bar{\mu}$.
We have $(z+\chi_0+\bar{\mu})^h=\In(P_{z+\chi_0},\omega, \bar{\psi})=\Ps(P_{z+\chi_0},G, \bar{\psi})$, $\bar{\mu}\neq 0$ and it follows that $G$ is $(z+\chi_0)$-removable, which is a contradiction.

We know $t=\Ord(\bar{\omega}, \bar{\Gamma}_+(P_{z+\mu},\psi))\leq\min\{\Ord(\bar{\omega}, \bar{\Gamma}_+(P_{z+\chi_0},\psi)), \Ord(\bar{\omega},\Gamma_+(P-\{z\},\mu-\chi_0))\}$.

Assume that $t=\Ord(\bar{\omega}, \bar{\Gamma}_+(P_{z+\mu},\psi))< \min\{\Ord(\bar{\omega}, \bar{\Gamma}_+(P_{z+\chi_0},\psi)), \Ord(\bar{\omega},\Gamma_+($\break$P-\{z\},\mu-\chi_0))\}$. 

It follows $\In(P_{z+\chi_0},\omega,\bar{\psi})=(z+\chi_0)^h$, $\In(P_{z+\mu},\omega,\bar{\psi})\neq(z+\mu)^h$, $\In(P_{z+\mu},\omega, z+\chi_0)= \In(P_{z+\mu},\omega, z+\mu-(\mu-\chi_0))=z+\mu$, and
$\In(P_{z+\mu},\omega,\bar{\psi})= (z+\mu)^h$.
We get a contradiction.

We know $t=\Ord(\bar{\omega}, \bar{\Gamma}_+(P_{z+\mu},\psi))=\min\{\Ord(\bar{\omega}, \bar{\Gamma}_+(P_{z+\chi_0},\psi)), \Ord(\bar{\omega},\Gamma_+(P-\{z\},\mu-\chi_0))\}$.

Since $\bar{\omega}\in\Map(P,\R_0)^\vee\cap V^*$ is an arbitrary element, we know 
$\bar{\Gamma}_+(P_w,\psi)=\bar{\Gamma}_+(P_{z+\mu},\psi)= \Conv(\bar{\Gamma}_+(P_{z+\chi_0},\psi)\cup\Gamma_+(P-\{z\},\chi_0-\mu))$.

We continue our reasoning.

By the above we have $t\leq \Ord(\bar{\omega},\Gamma_+(P-\{z\}, \mu-\chi_0)=\Ord(P-\{z\},\bar{\omega},\mu-\chi_0)$. Let $\bar{\mu}\in M(R^{\prime c})$ be the same as above.
$\In(P_{z+\chi_0},\omega, z+\mu)=z+ \chi_0+\bar{\mu}$.
Let $\lambda:R^c\rightarrow R^c$ be the isomorphism of $k$-algebras such that $\lambda(z+\mu)= z+ \chi_0+\bar{\mu}$ and $\lambda(x)=x$ for any $x\in P-\{z\}$.

Assume  $t=\Ord(\bar{\omega}, \bar{\Gamma}_+(P_{z+\mu},\psi))= \Ord(\bar{\omega}, \bar{\Gamma}_+(P_{z+\chi_0},\psi))$. 

We have $\lambda(\Ps(P_{z+\mu},\Delta(\omega, \Gamma_+( P_{z+\mu},\psi)),\psi))=\lambda(\In(P_{z+\mu},\omega,\psi))= \In(P_{z+\chi_0},\omega,$\break$\psi)=\Ps(P_{z+\chi_0},\Delta(\omega,\Gamma_+( P_{z+\chi_0},\psi)),\psi) $. Since the face $\Delta(\omega,\Gamma_+( P_{z+\chi_0},\psi))$ of $\Gamma_+($\break$P_{z+\chi_0},\psi)$ is not $(z+\chi_0)$-removable, it follows that the face $\Delta(\omega,\Gamma_+( P_{z+\mu},\psi))$ of $\Gamma_+( P_{z+\mu},\psi)$ is not $(z+\mu)$-removable.

Assume  $t=\Ord(\bar{\omega}, \bar{\Gamma}_+(P_{z+\mu},\psi))< \Ord(\bar{\omega}, \bar{\Gamma}_+(P_{z+\chi_0},\psi))$. 

We have $\In(P_{z+\chi_0},\omega,\bar{\psi})=(z+\chi_0)^h$,
$\lambda(\Ps(P_{z+\mu},\Delta(\omega, \Gamma_+( P_{z+\mu},\psi)),\bar{\psi}))=\lambda(\In($\break$P_{z+\mu},\omega,\bar{\psi}))= \In(P_{z+\chi_0},\omega,\bar{\psi})=(z+\chi_0)^h$, $\Ps(P_{z+\mu},\Delta(\omega, \Gamma_+( P_{z+\mu},\psi)),\bar{\psi})=(z+\mu-\bar{\mu})^h$, and $\bar{\mu}\neq 0$. It follows that the face $\Delta(\omega,\Gamma_+( P_{z+\mu},\psi))$ of $\Gamma_+( P_{z+\mu},\psi)$ is $(z+\mu)$-removable.

Since $\bar{\omega}\in\Map(P,\R_0)^\vee\cap V^*$ is an arbitrary element, we know 
that for any face $F$ of $\Gamma_+(P_w,\psi)$ with $hf^P_z\in F$ and $F\cap U\neq\emptyset$, $F$ is $w$-removable, if and only if, $\rho(F\cap U)\cap\rho(\Gamma_+(P_{z+\chi_0}, \psi)\cap U)=\emptyset$.

\noindent $4.$\quad

\noindent (a)$\Leftrightarrow$(b) It follows from the former half of 3.

\noindent (a)$\Leftrightarrow$(c) I It follows from the latter half of 3.

We show the uniqueness of the element $\chi_0$ in 2 here.

Assume that $\chi_0\in M(R^{\prime c})$, $\Gamma_+(P_{z+\chi_0},\psi)$ has no $(z+\chi_0)$-removable faces, $\Supp(P-\{z\},\chi_0)\subset\rho(\Gamma_+(P,\psi)\cap U)- \rho(\Gamma_+(P_{z+\chi_0},\psi)\cap U)$,
$\mu_0\in M(R^{\prime c})$, $\Gamma_+(P_{z+\mu_0},$\break$\psi)$ has no $(z+\mu_0)$-removable faces and $\Supp(P-\{z\},\mu_0)\subset\rho(\Gamma_+(P,\psi)\cap U)- \rho(\Gamma_+(P_{z+\mu_0},\psi)\cap U)$.

By 4, (c)$\Rightarrow$(a), we have $\rho(\Gamma_+(P_{z+\chi_0},\psi)\cap U)= \rho(\Gamma_+(P_{z+\mu_0},\psi)\cap U)$. By 4, (c)$\Rightarrow$(b), we have $\Supp(P-\{z\},\mu_0-\chi_0)\subset\rho(\Gamma_+(P_{z+\chi_0},\psi)\cap U)$. By assumption $\Supp(P-\{z\}, \mu_0-\chi_0)\subset\Supp(P-\{z\},\mu_0)\cup\Supp(P-\{z\},\chi_0)\subset \rho(\Gamma_+(P,\psi)\cap U)- \rho(\Gamma_+(P_{z+\chi_0},\psi)\cap U)= \rho(\Gamma_+(P,\psi)\cap U)- \rho(\Gamma_+(P_{z+\mu_0},\psi)\cap U)$. We have 
$\Supp(P-\{z\},\mu_0-\chi_0)=\emptyset$ and $\chi_0=\mu_0$.

\noindent $5.$\quad By Lemma~\ref{vertex projection}.3, $\rho(\Gamma_+(P_{z+\chi_0},\psi)\cap U)=\emptyset$, if and only if, $\bar{\psi}=(z+\chi_0)^h$. Assume $\bar{\psi}=(z+\chi_0)^h$. We have $\psi=\bar{u}(z+\chi_0)^h$ and $\bar{u}\in R^{c\times}$. Since $\psi\in R$ and $R$ is a UFD, we know by analytic unramifiedness (Matsumura~\cite{M}, page 236, THEOREM 70 and page 240, THEOREM 72.) that there exist $u\in R^\times$ and $\lambda\in M(R)$ with $\psi=u\lambda^h$. 

Conversely assume that $u\in R^\times$, $\lambda\in M(R)$ and $\psi=u\lambda^h$.
If $h=0$, then $\psi=u\in R^\times$, and $\rho(\Gamma_+(P_{z+\chi_0},\psi)\cap U)=\emptyset$. We assume $h>0$ below.
We have $\Gamma_+(P_{z+\chi_0},$\break$\psi)=h\Gamma_+(P_{z+\chi_0},\lambda)$. Since $\Gamma_+(P_{z+\chi_0},\psi)$ is of $(z+\chi_0)$-Weierstrass type and $\{hf^P_z\}$ is the unique $(z+\chi_0)$-top vertex, we know that $\Gamma_+(P_{z+\chi_0},\lambda)$ is of $(z+\chi_0)$-Weierstrass type and $\{f^P_z\}$ is the unique $(z+\chi_0)$-top vertex. It follows that there exist $v\in R^{c\times}$ and $\mu\in M(R^{\prime c})$ with $\lambda=v(z+\chi_0+\mu)$. We take $v\in R^{c\times}$ and $\mu\in M(R^{\prime c})$ with $\lambda=v(z+\chi_0+\mu)$. $\psi=uv^h(z+\chi_0+\mu)^h$. If $\mu\neq 0$, then $\Gamma_+(P_{z+\chi_0},\psi)$ is $(z+\chi_0)$-removable face of $\Gamma_+(P_{z+\chi_0},\psi)$. Since $\Gamma_+(P_{z+\chi_0},\psi)$ has no $(z+\chi_0)$-removable faces, we know that $\mu=0$, $\psi= uv^h(z+\chi_0)^h$ and $\rho(\Gamma_+(P_{z+\chi_0},\psi)\cap U)=\emptyset$ by Lemma~\ref{vertex projection}.3. 

\noindent $6.$\quad
We consider the case $\rho(\Gamma_+(P_{z+\chi_0},\psi)\cap U)=\emptyset$ first.
If $h=0$, then $\chi_0=0\in M(R^{\prime h})$. We assume $h>0$ below.
By Lemma~\ref{vertex projection}.3, we have $\psi=\bar{u}(z+\chi_0)^h$. By 5, $\psi=u\lambda^h$ for some $u\in R^\times$ and some $\lambda\in M(R)$. We take $u\in R^\times$ and $\lambda\in M(R)$ with $\psi=u\lambda^h$. Since $R^c$ is a UFD, we know that $\lambda=v(z+\chi_0)$ for some $v\in R^{c\times}$. We take $v\in R^{c\times}$ with $\lambda=v(z+\chi_0)$. Since $\lambda\in R\subset R^h$, by Henselian Weierstrass Theorem (Hironaka~\cite{H772}.), we know $\chi_0\in M(R^{\prime h})$ and $v\in R^{h\times}$.

We consider the case $\rho(\Gamma_+(P_{z+\chi_0},\psi)\cap U)\neq\emptyset$.

Consider any $\chi\in M(R^{\prime h})$ with $\Supp(P-\{z\},\chi)\subset \rho(\Gamma_+(P,\psi)\cap U)-\rho(\Gamma_+(P_{z+\chi_0},$\break$\psi)\cap U)$. By 3, we have $\rho(\Gamma_+(P,\psi)\cap U)\supset\Conv(\rho(\Gamma_+(P_{z+\chi_0},\psi)\cap U)\cup\Gamma_+(P-\{z\},\chi-\chi_0))=\rho(\Gamma_+(P_{z+\chi},\psi)\cap U) \supset\rho(\Gamma_+(P_{z+\chi_0},\psi)\cap U)\neq\emptyset$ and $h>0$.

Consider any $F\in\mathcal{F}(\rho(\Gamma_+(P_{z+\chi_0},\psi)\cap U))^1$. We take the unique $\bar{\mu}_F\in\Map(P,\R_0$\break$)^\vee\cap V^*\cap \Map(P,\Z)^*$ with $F=\Delta(\bar{\mu}_F, \rho(\Gamma_+(P_{z+\chi_0},\psi)\cap U))$ and $\Q_0\bar{\mu}_F\cap\Map(P,\Z)^*=\Z_0\bar{\mu}_F$. 
Let $t_F=\Ord(\bar{\mu}_F, \rho(\Gamma_+(P_{z+\chi},\psi)\cap U))\in\Q_0$ and let $\mu_F=\bar{\mu}_F+t_F f^{P\vee}_z\in\Map(P,\R_0)^\vee$. We have 
$\Ord(\bar{\mu}_F, \rho(\Gamma_+(P,\psi)\cap U))\leq t_F\leq\Ord(\bar{\mu}_F, \rho(\Gamma_+(P_{z+\chi_0},\psi)\cap U))\in\Q_0$ and $\Ord(\mu_F, \Gamma_+(P_{z+\chi},\psi))=ht_F$.

By 3, we know the following:
\begin{enumerate}
\item If $t_F=\Ord(\mu_F, \rho(\Gamma_+(P_{z+\chi_0},\psi)\cap U))$, then the face $\Delta(\mu_F, \Gamma_+(P_{z+\chi},\psi))$ of $\Gamma_+(P_{z+\chi},\bar{\psi})$ is not $(z+\chi)$-removable.
\item If $t_F<\Ord(\mu_F, \rho(\Gamma_+(P_{z+\chi_0},\psi)\cap U))$, then $t_F\in\Z_0$ and the face $\Delta(\mu_F, \Gamma_+($\break$P_{z+\chi},\bar{\psi}))$ of $\Gamma_+(P_{z+\chi},\bar{\psi})$ is $(z+\chi)$-removable.
\end{enumerate}

Note that $\mathcal{F}(\rho(\Gamma_+(P_{z+\chi_0},\psi)\cap U))^1$ is a finite set. Let $$m_\chi=\sharp\{ F\in\mathcal{F}(\rho(\Gamma_+(P_{z+\chi_0},\psi)\cap U))^1|
t_F<\Ord(\mu_F, \rho(\Gamma_+(P_{z+\chi_0},\psi)\cap U))\}\in\Z_0.$$

Consider the case $m_\chi=0$. In this case we have $\rho(\Gamma_+(P_{z+\chi},\psi)\cap U)= \rho(\Gamma_+(P_{z+\chi_0},$\break$\psi)\cap U)$. By 4, (a)$\Rightarrow$(c) we know that $\Gamma_+(P_{z+\chi},\psi)$ has no $z$-removable faces. By the uniqueness in 2, we have $\chi_0=\chi\in M(R^{\prime h})$.

Consider the case $m_\chi>0$. We take any $G\in\mathcal{F}(\rho(\Gamma_+(P_{z+\chi_0},\psi)\cap U))^1$ with $t_G<\Ord(\mu_G, \rho(\Gamma_+(P_{z+\chi_0},\psi)\cap U))$. $t_G\in\Z_0$ and $\hat{G}=\Delta(\mu_G, \Gamma_+(P_{z+\chi},\bar{\psi}))$ is a $(z+\chi)$-removable face.
We take $\chi_1\in M(R^{\prime c})$ with $\Ps(P_{z+\chi},\hat{G},\bar{\psi})=(z+\chi_1)^h$. 
Since $\bar{u}\bar{\psi}=\psi\in M(R)\subset M(R^h)$, by Henselian Weierstrass Theorem, we know $\bar{\psi}\in M(R^h)$ and $(z+\chi_1)^h =\Ps(P_{z+\chi},\hat{G},\bar{\psi})=\In(P_{z+\chi},\mu_G,\bar{\psi})\in M(R^h)$. By analytic unramifiedness we know that $(z+\chi_1)^h=u_2\lambda_2^h$ for some $u_2\in R^{h\times}$ and some $\lambda_2\in M(R^h)$. We take $u_2\in R^{h\times}$ and $\lambda_2\in M(R^h)$ with $(z+\chi_1)^h=u_2\lambda_2^h$. Since $R^h$ is a UFD, we know hat $\lambda_2=u_3(z+\chi_1)$ for some $u_3\in R^{c\times}$. By Henselian Weierstrass Theorem we know $\chi_1\in M(R^{\prime h})$. 

$\Supp(P-\{z\},\chi_1-\chi)\subset\Delta(\bar{\mu}_G, \rho(\Gamma_+(P_{z+\chi},\psi)\cap U))\subset\rho(\Gamma_+(P,\psi)\cap U)-\rho(\Gamma_+(P_{z+\chi_0},\psi)\cap U)$. 
$\Supp(P-\{z\},\chi_1)\subset \Supp(P-\{z\},\chi)\cup\Supp(P-\{z\},\chi_1-\chi)\subset\rho(\Gamma_+(P,\psi)\cap U)-\rho(\Gamma_+(P_{z+\chi_0},\psi)\cap U)$.

$\rho(\Gamma_+(P_{z+\chi},\psi)\cap U)-G\supset\rho(\Gamma_+(P_{z +\chi_1},\psi)\cap U) \supset\rho(\Gamma_+(P_{z+\chi_0},\psi)\cap U)$.

$t_G<\Ord(\mu_G, \rho(\Gamma_+(P_{z+\chi_1},\psi)\cap U))\leq\Ord(\mu_G, \rho(\Gamma_+(P_{z+\chi_0},\psi)\cap U))$.
For any $F\in\mathcal{F}(\rho(\Gamma_+(P_{z+\chi_0},\psi)\cap U))^1$ with $F\neq G$, $t_F=\Ord(\mu_F, \rho(\Gamma_+(P_{z +\chi_1},\psi)\cap U))$.
Therefore, if $\Ord(\mu_G, \rho(\Gamma_+(P_{z +\chi_1},\psi)\cap U))=\Ord(\mu_G, \rho(\Gamma_+(P_{z+\chi_0},\psi)\cap U))$, then $m_{\chi_1}<m_\chi$ and by induction on $m_\chi$ we can conclude $\chi_0\in M(R^{\prime h})$. If  $\Ord(\mu_G, \rho(\Gamma_+(P_{z +\chi_1},$\break$\psi)\cap U))<\Ord(\mu_G, \rho(\Gamma_+(P_{z+\chi_0},\psi)\cap U))$, then $m_{\chi_1}=m_\chi$ and by induction on $\Ord(\mu_G,\rho(\Gamma_+(P_{z+\chi_0},\psi)\cap U))- \Ord(\mu_G, \rho(\Gamma_+(P_{z+\chi},\psi)\cap U))$ we can conclude $\chi_0\in M(R^{\prime h})$.

We conclude $\chi_0\in M(R^{\prime h})$ in all cases.

\noindent $7.$\quad
Assume moreover, that either $\rho(\Gamma_+(P_{z+\chi_0},\psi)\cap U)$ has at most one vertex, or $\dim R\leq 3$.

We consider the case $\rho(\Gamma_+(P_{z+\chi_0},\psi)\cap U)=\emptyset$ first.
If $h=0$, then $\Gamma_+(P,\psi)=\Map(P,\R_0)$ and $\Gamma_+(P,\psi)$ has no $z$-removable faces. We assume $h>0$ below. 
By 5, we know that $\psi=uw_1^h$ for some $u\in R^\times$ and some $w_1\in M(R)$. We take $u\in R^\times$ and $w_1\in M(R)$ with $\psi=uw_1^h$. By the proof of 5, we know $\partial w_1/\partial z\in R^\times$. Since $\Gamma_+(P_{w_1},\psi)= \Gamma_+(P_{w_1}, uw_1^h)=\{hf^P_z\}+\Map(P,\R_0)$, we know that $\Gamma_+(P_{w_1},\psi)$ has no $w_1$-removable faces.

We consider the case where $\rho(\Gamma_+(P_{z+\chi_0},\psi)\cap U)$ has only one vertex. Since $\rho(\Gamma_+(P_{z+\chi_0},\psi)\cap U)\neq\emptyset$, we have $h>0$.

Consider any $w\in M(R)$ with $\partial w/\partial z\in R^\times$, and consider any $x\in P-\{z\}$. We denote $$t_x=\Ord(f^{P\vee}_x, \rho(\Gamma_+(P_w,\psi)\cap U))\in \Q_0\text{, and}$$
$$m_w=\sharp\{x\in P-\{z\}| t_x<\Ord(f^{P\vee}_x, \rho(\Gamma_+(P_{z+\chi_0},\psi)\cap U))\}\in\Z_0.$$

We have $\Ord(f^{P\vee}_x, \rho(\Gamma_+(P,\psi)\cap U))\leq t_x\leq \Ord(f^{P\vee}_x, \rho(\Gamma_+(P_{z+\chi_0},\psi)\cap U))$.

Since $\rho(\Gamma_+(P_{z+\chi_0},\psi)\cap U)$ has only one vertex, by 4, $m_w=0$, if and only if, $\rho(\Gamma_+(P_w,\psi)\cap U)= \rho(\Gamma_+(P_{z+\chi_0},\psi)\cap U)$, if and only if, $\Gamma_+(P_w,\psi)$ has no $w$-removable faces.

Consider the case $m_w>0$. 
We denote $P_0=\{x\in P-\{z\}| t_x=\Ord(f^{P\vee}_x, \rho(\Gamma_+($\break$P_{z+\chi_0},\psi)\cap U))\}$ and $r=\sharp P_0$. $P_0\neq P-\{z\}$. We take any $y\in P-(\{z\}\cup P_0)$. We have $t_y<\Ord(f^{P\vee}_y, \rho(\Gamma_+(P_{z+\chi_0},\psi)\cap U))$ and $t_y\in \Z_0$. Let $\mu_y= f^{P\vee}_y+t_y f^{P\vee}_z\in\Map(P,\R_0)^\vee$ and let $\hat{F}_y=\Delta(\mu_y,\Gamma_+(P_w,\psi))$. $\hat{F}_y$ is a $w$-removable face of $\Gamma_+(P_w,\psi)$, $\rho(\hat{F}_y\cap U)$ is a face of $\rho(\Gamma_+(P_w,\psi)\cap U)$, and any vertex of $\rho(\hat{F}_y\cap U)$ belongs to $\Map(P,\Z_0)\cap V$. We have $\Ord(f^{P\vee}_y, \rho(\hat{F}_y\cap U))=t_y\in\Z_0$ and $\lceil t_x\rceil\leq \Ord(f^{P\vee}_x, \rho(\hat{F}_y\cap U))$ for any $x\in P_0$, since $t_x=\Ord(f^{P\vee}_x, \rho(\Gamma_+(P_w,\psi)\cap U))\leq \Ord(f^{P\vee}_x, \rho(\hat{F}_y\cap U))\in\Z_0$ for any $x\in P_0$. We take the unique pair $u\in R^{c\times}$ and $\lambda\in M(R^{\prime c})$ such that $\In(P_w, \mu_y,\psi)=\Ps(P_w, \hat{F}_y,\psi)=u(w+\lambda)^h$ and $\lambda\neq 0$. $\Ord(P-\{z\}, f^{P\vee}_y,\lambda)=t_y$ and $\Ord(P-\{z\}, f^{P\vee}_x,\lambda)\geq \lceil t_x\rceil$ for any $x\in P_0$.

Let $x_1,x_2,\ldots,x_r$ be all the elements in $P_0$.
Let $m=(\sum_{i=1}^r\lceil t_{x_i}\rceil)+t_y\in\Z_0$. We define a mapping $E:\{1,2,\ldots,m\}\rightarrow \mathcal{F}(\Map(P,\R_0)^\vee\cap V^*)_1$ by putting
$$
E(i)=
\begin{cases}
\R_0f^{P\vee}_{x_j}&\text{if $\sum_{k=1}^{j-1}\lceil t_{x_k}\rceil<i\leq\sum_{k=1}^{j}\lceil t_{x_k}\rceil$},\\
\R_0f^{P\vee}_y&\text{if $\sum_{k=1}^{r}\lceil t_{x_k}\rceil<i\leq m$},
\end{cases}
$$
for any $i\in\{1,2,\ldots,m\}$.

We consider the basic subdivision
$$\Omega=\Omega(\Map(P,\R)^*, \Map(P,\Z)^*, \R_0f^{P\vee}_z,\mathcal{F}(\Map(P,\R_0)^\vee),m,E)$$
of $\mathcal{F}(\Map(P,\R_0)^\vee)$ associated with the pair $(m,E)$. $\Omega$ is an iterated star subdivision of $\mathcal{F}(\Map(P,\R_0)^\vee)$ and is a flat regular fan over $\Map(P,\Z)^*$ in $\Map(P,\R)^*$. $|\Omega|=\Map(P,\R_0)^\vee$.

Let $\hat{X}=X(k, \Map(P,\R)^*, \Map(P,\Z)^*,\Omega)$, $\hat{Y}= X(k, \Map(P,\R)^*, \Map(P,\Z)^*,$\break$\mathcal{F}(\Map(P,\R_0)^\vee))$, and $\hat{\sigma}=\sigma(k, \Map(P,\R)^*, \Map(P,\Z)^*, \mathcal{F}(\Map(P,\R_0)^\vee),\Omega):\hat{X}\rightarrow\hat{Y}$. $\hat{X}$ is the toric variety over $k$ associated with the regular fan $\Omega$, $\hat{Y}$ is the toric variety over $k$ associated with the regular fan $\mathcal{F}(\Map(P,\R_0)^\vee)$, and $\hat{\sigma}$ is the subdivision morphism over $k$ associated with the pair $(\Omega, \mathcal{F}(\Map(P,\R_0)^\vee)$.

$\hat{Y}=U(\mathcal{F}(\Map(P,\R_0)^\vee), \Map(P,\R_0)^\vee)$ is an affine scheme. The ring of regular functions $\mathcal{O}_{\hat{Y}}(\hat{Y})=\mathcal{O}_{\hat{Y}}( U(\mathcal{F}(\Map(P,\R_0)^\vee), \Map(P,\R_0)^\vee))$ over $\hat{Y}$ is a polynomial ring over $k$ with variables $\{\chi(f^P_x)|x\in P\}$. By the injective homomorphism $\mathcal{O}_{\hat{Y}}(\hat{Y})\rightarrow R$ of $k$-algebras sending $\chi(f^P_x)\in \mathcal{O}_{\hat{Y}}(\hat{Y})$ to $x\in R$ for any $x\in P-\{z\}$ and sending $\chi(f^P_z)\in \mathcal{O}_{\hat{Y}}(\hat{Y})$ to $w\in R$, we regard $\mathcal{O}_{\hat{Y}}(\hat{Y})$ as an subring of $R$. $\mathcal{O}_{\hat{Y}}(\hat{Y})=k[P_w]\subset R$. $\hat{Y}=\Spec(k[P_w])$.
The inclusion ring homomorphism $k[P_w]\rightarrow R$ induces a morphism $\hat{\pi}:\Spec(R)\rightarrow \hat{Y}$ of schemes over $k$.

Consider the fiber product scheme $X=\hat{X}\times_{\hat{Y}}\Spec(R)$ of $\hat{X}$ and $\Spec(R)$ over $\hat{Y}$, the projection $\pi:X\rightarrow\hat{X}$ and the projection $\sigma:X\rightarrow\Spec(R)$.
Let $D=\Spec(R/\prod_{x\in P_w}x R)$. $D$ is a normal crossing divisor on $\Spec(R)$ and $\sigma$ is an admissible composition of $m$ of blowing-ups with center of codimension two over $D$.

We take the unique element $\Delta_0\in\Omega^0$ with $\R_0f^{P\vee}_z\subset\Delta_0$, the unique closed point $\hat{a}\in\hat{X}$ with  $\{\hat{a}\}=V(\Omega,\Delta_0)$ and the unique closed point $a\in X$ with $\pi(a)=\hat{a}$. We consider the homomorphism $\sigma^*:R\rightarrow\mathcal{O}_{X,a}$ of local $k$-algebras induced by $\sigma$. We know $P-\{z\}\subset \mathcal{O}_{X,a}$ and $\sigma^*(x)=x$ for any $x\in P-\{z\}$. We denote $\hat{w}=\sigma^*(w)/((\prod_{x\in P_0}x^{\lceil t_x\rceil})y^{t_y})$ and $P_{\hat{w}}=(P-\{z\})\cup\{\hat{w}\}$. We know $\hat{w}\in M(\mathcal{O}_{X,a})$ and $P_{\hat{w}}$ is a parameter system of $\mathcal{O}_{X,a}$.
We consider the homomorphism $\sigma^*:R^c\rightarrow\mathcal{O}_{X,a}^c$ between the completions induced by $\sigma^*$. We denote $\hat{\lambda}=\sigma^*(\lambda)/((\prod_{x\in P_0}x^{\lceil t_x\rceil})y^{t_y})$. We know $\hat{\lambda}\in \mathcal{O}_{X,a}^c-y\mathcal{O}_{X,a}^c$, and $\In(P_{\hat{w}},f^{P_{\hat{w}}\vee}_y, \sigma^*(\psi))=\sigma^*(\In(P_w, \mu_y,\psi))= \sigma^*(u) (\prod_{x\in P_0}x^{h\lceil t_x\rceil})y^{ht_y} (\hat{w}+\hat{\lambda})^h$. Let $\theta: \mathcal{O}_{X,a}\rightarrow \mathcal{O}_{X,a}/y\mathcal{O}_{X,a}$ denote the canonical surjective homomorphism to the residue ring. It induces the homomorphism $\theta: \mathcal{O}_{X,a}^c\rightarrow (\mathcal{O}_{X,a}/y\mathcal{O}_{X,a})^c$ between the completions. We know that there exists uniquely an element $\hat{\psi}\in \mathcal{O}_{X,a}$ with $\sigma^*(\psi)= y^{ht_y}\hat{\psi}$. We take the element $\hat{\psi}\in \mathcal{O}_{X,a}$ with $\sigma^*(\psi)= y^{ht_y}\hat{\psi}$. We know $\hat{\psi}\in \mathcal{O}_{X,a}-y\mathcal{O}_{X,a}$ and $\theta (\hat{\psi})\neq 0$. We know that exists uniquely an element $\tilde{\psi}\in \mathcal{O}_{X,a}/y\mathcal{O}_{X,a}$ with $\theta(\hat{\psi})= (\prod_{x\in P_0}x^{h\lceil t_x\rceil})\tilde{\psi}$. We take the element $\tilde{\psi}\in \mathcal{O}_{X,a}/y\mathcal{O}_{X,a}$ with $\theta(\hat{\psi})= (\prod_{x\in P_0}x^{h\lceil t_x\rceil})\tilde{\psi}$. 
(Note that there does not exist $\hat{\psi}_0\in \mathcal{O}_{X,a}$ with $\hat{\psi}= (\prod_{x\in P_0}x^{h\lceil t_x\rceil})\hat{\psi}_0$, if $t_x\not\in\Z$ for some $x\in P_0$.)
We know $\tilde{\psi}=\theta\sigma^*(u) \theta(\hat{w}+\hat{\lambda})^h$. By analytic unramifiednss, we know there exist $\tilde{v}\in (\mathcal{O}_{X,a}/y\mathcal{O}_{X,a})^\times$ and $\tilde{w}_1\in \mathcal{O}_{X,a}/y\mathcal{O}_{X,a}$ with $\tilde{\psi}=\tilde{v}\tilde{w}_1^h$. We know that there exist $\hat{v}\in \mathcal{O}_{X,a}^\times$ and $\hat{w}_0\in \mathcal{O}_{X,a}$ with $\sigma^*(\psi)-\hat{v}((\prod_{x\in P_0}x^{\lceil t_x\rceil})y^{t_y}\hat{w}_0)^h\in y^{ht_y+1}\mathcal{O}_{X,a}$. We take $\hat{v}\in \mathcal{O}_{X,a}^\times$ and $\hat{w}_0\in \mathcal{O}_{X,a}$ with $\sigma^*(\psi)-\hat{v}((\prod_{x\in P_0}x^{\lceil t_x\rceil})y^{t_y}\hat{w}_0)^h\in y^{ht_y+1}\mathcal{O}_{X,a}$.

Now, by induction on $m$, we know that there exists $v_0\in \mathcal{O}_{X,a}^\times $ and $w_0\in R$ with $(\prod_{x\in P_0}x^{\lceil t_x\rceil})y^{t_y}\hat{w}_0= v_0\sigma^*(w_0)$. We take $v_0\in \mathcal{O}_{X,a}^\times $ and $w_0\in R$ with $(\prod_{x\in P_0}x^{\lceil t_x\rceil})y^{t_y}\hat{w}_0= v_0\sigma^*(w_0)$. We have $\sigma^*(\psi)-\hat{v}v_0^h\sigma^*(w_0)^h\in y^{ht_y+1}\mathcal{O}_{X,a}$. 
Therefore, $\In(P_w,\mu_y,w_0)=u_{00}(w+\lambda)$ for some $u_{00}\in R^{c\times}$, $w_0\in M(R)$, $\partial w_0/\partial x\in R^\times$ and $t_y<\Ord(f^{P\vee}_y,\rho(\Gamma_+(P_{w_0},\psi)\cap U))\leq\Ord(f^{P\vee}_y,\rho(\Gamma_+(P_{z+\chi_0},\psi)\cap U))$.

Since $(\prod_{x\in P_0}x^{\lceil t_x\rceil})y^{t_y}\hat{w}_0= v_0\sigma^*(w_0)$, if we take $u_0\in R^{c\times}$ and $\lambda_0\in M(R^{\prime c})$ with $w_0=u_0(w+\lambda_0)$, then $\Ord(P-\{z\}, f^{P\vee}_x,\lambda_0)\geq \lceil t_x\rceil\geq t_x$ for any $x\in P_0$ and $\Ord(P-\{z\}, f^{P\vee}_y,\lambda_0)= t_y$. By 3, we have $t_x=\Ord(f^{P\vee}_x,\rho(\Gamma_+(P_{w_0},\psi)\cap U))=\Ord(f^{P\vee}_x,\rho(\Gamma_+(P_{z+\chi_0},\psi)\cap U))$ for any $x\in P_0$.

Consider the case $m_{w_0}<m_w$. By induction on $m_w$, we can conclude that there exists $w_1\in M(R)$ such that $\partial w_1/\partial z\in R^\times$ and $\Gamma_+(P_{w_1},\psi)$ has no $w_1$-removable faces.

Consider the case $m_{w_0}\geq m_w$. In this case we have  $m_{w_0}=m_w$ and $\Ord(f^{P\vee}_y,\rho($\break$\Gamma_+(P_{w_0},\psi)\cap U))>\Ord(f^{P\vee}_y,\rho(\Gamma_+(P_w,\psi)\cap U))$.
By induction on $\lceil \Ord(f^{P\vee}_y,\rho(\Gamma_+($\break$P_{z+\chi_0},\psi)\cap U))\rceil-\lceil \Ord(f^{P\vee}_y,\rho(\Gamma_+(P_w,\psi)\cap U))\rceil$, we can conclude that there exists $w_1\in M(R)$ such that $\partial w_1/\partial z\in R^\times$ and $\Gamma_+(P_{w_1},\psi)$ has no $w_1$-removable faces.

We conclude that our claim 7 holds, if $\rho(\Gamma_+(P_{z+\chi},\psi)\cap U)$ has at most one vertex.

Note that if $\dim R=2$, then $\rho(\Gamma_+(P_{z+\chi},\psi)\cap U)$ has at most one vertex.

We consider the case $\dim R=3$ and $\rho(\Gamma_+(P_{z+\chi},\psi)\cap U)\neq\emptyset$.

We consider $\mathcal{F}(\rho(\Gamma_+(P_{z+\chi_0},\psi)\cap U))^1$. Since $\dim \rho(\Gamma_+(P_{z+\chi_0},\psi)\cap U)=2$, we have two of non-compact faces of $\rho(\Gamma_+(P_{z+\chi_0},\psi)\cap U)$ of codimension one. Let $F_0, F_1\in \mathcal{F}(\rho(\Gamma_+(P_{z+\chi_0},\psi)\cap U))^1$ be the two non-compact faces. Any $F\in\mathcal{F}(\rho(\Gamma_+(P_{z+\chi_0},\psi)\cap U))^1-\{F_0, F_1\}$ is compact. $\{\Stab(F_0), \Stab(F_1)\}=\mathcal{F}(\Map(P,$\break$\R_0)\cap V)^1$.

Consider any $w\in M(R)$ such that $\partial w/\partial z\in R^\times$.
Since $\rho(\Gamma_+(P_w,\psi)\cap U)\supset\rho(\Gamma_+(P_{z+\chi},\psi)\cap U)\neq\emptyset$, we have $h>0$.

Consider any $F\in\mathcal{F}(\rho(\Gamma_+(P_{z+\chi_0},\psi)\cap U))^1$. We take the unique $\bar{\mu}_F\in\Map(P,\R_0$\break$)^\vee\cap V^*\cap \Map(P,\Z)^*$ with $F=\Delta(\bar{\mu}_F, \rho(\Gamma_+(P_{z+\chi_0},\psi)\cap U))$ and $\Q_0\bar{\mu}_F\cap\Map(P,\Z)^*=\Z_0\bar{\mu}_F$. 
Let $t_F=\Ord(\bar{\mu}_F, \rho(\Gamma_+(P_w,\psi)\cap U))\in\Q_0$ and let $\mu_F=\bar{\mu}_F+t_F f^{P\vee}_z\in\Map(P,\R_0)^\vee$. We have 
$t_F\leq\Ord(\bar{\mu}_F, \rho(\Gamma_+(P_{z+\chi_0},\psi)\cap U))\in\Q_0$ and $\Ord(\mu_F, \Gamma_+(P_w,$\break$\psi))=ht_F$.

Now, repeating the reasoning in the above case where $\rho(\Gamma_+(P_{z+\chi},\psi)\cap U)$ has only one vertex, we know that there exists $w_1\in M(R)$ such that $\partial w_1/\partial z\in R^\times$ and $\Gamma_+(P_{w_1},\psi)$ has no \emph{non-compact} $w_1$-removable faces. We take  $w_1\in M(R)$ such that $\partial w_1/\partial z\in R^\times$ and $\Gamma_+(P_{w_1},\psi)$ has no non-compact $w_1$-removable faces.
By replacing $P$ by $P_{w_1}$ we can assume $t_{F_i}=\Ord(\mu_{F_i}, \rho(\Gamma_+(P_{z+\chi_0},\psi)\cap U))$ for any $i\in\{0,1\}$. We assume this condition below. 

Consider any $F\in\mathcal{F}(\rho(\Gamma_+(P_{z+\chi_0},\psi)\cap U))^1$ with 
$t_F<\Ord(\bar{\mu}_F, \rho(\Gamma_+(P_{z+\chi_0},\psi)\cap U))$. We know that $F$ is compact, $0<t_F$, the face $\hat{F}=\Delta(\mu_F,\Gamma_+(P_w,\psi))$ of $\Gamma_+(P_w,\psi)$ with $\rho(\hat{F}\cap U)=F$ is compact and $\In(P_w,\mu_F,\psi)\in k[P_w]\subset R$. Since $\hat{F}$ is $w$-removable, by analytic unramifiedness we know that $\In(P_w,\mu_F,\psi)=u_F(w+\lambda_F)^h$ for some $u_F\in k-\{0\}$ and some $\lambda_F\in (P-\{z\})k[P-\{z\}]-\{0\}$. Using this fact and repeating reasoning in the proof of claim 6, we know that there exists $w_1\in M(R)$ such that $\partial w_1/\partial z\in R^\times$ and $\Gamma_+(P_{w_1},\psi)$ has no $w_1$-removable faces.

We know that Theorem~\ref{erase faces} holds.

Note here that $\dim R'=\dim R-1<\dim R$, and any $\phi'\in R'$ with $\phi'\neq 0$ has normal crossings over $P'$ if $\dim R=2$. Therefore, we decide that we use induction on $\dim R$, and we can assume the following claim $(*)$ whenever $\dim R\geq 2$.:

\begin{description}
\item[$(*)$]
For any  $\phi'\in R'$ with $\phi'\neq 0$, there exists a composition $\sigma':X'\rightarrow\Spec(R')$ of finite blowing-ups with center in a closed irreducible smooth subscheme such that the divisor on $X'$ defined by the pull-back $\sigma^{\prime *}(\phi')\in\mathcal{O}_{X'}(X')$ of $\phi'$ by $\sigma'$ has normal crossings.
\end{description}

Claim $(*)$ is true, if $\dim R\leq 2$.

Let $\sigma':X'\rightarrow\Spec(R')$ be any composition of blowing-ups with center in a closed irreducible smooth subscheme. The scheme $X'$ is smooth. We consider a morphism $\Spec(R)\rightarrow\Spec(R')$ induced by the inclusion ring homomorphism $R'\rightarrow R$, the product scheme $X=X'\times_{\Spec(R')}\Spec(R)$, the projection $\sigma:X\rightarrow\Spec(R)$, and the projection $\pi:X\rightarrow X'$. We know the following (Lemma~\ref{pull back blowing-ups}.):
\begin{enumerate}
\item The morphism $\sigma$ is a composition of finite blowing-ups with center in a closed irreducible smooth subscheme. The scheme $X$ is smooth.
\item We consider the prime divisor $\Spec(R/zR)$ on $\Spec(R)$  defined by $z\in R$. The pull-back $\sigma^*\Spec(R/zR)$ of $\Spec(R/zR)$ by $\sigma$ is a smooth prime divisor of $X$, and $\sigma^*\Spec(R/zR)\supset \sigma^{-1}(M(R))$.
\item The projection $\pi:X\rightarrow X'$ induces an isomorphism $\sigma^*\Spec (R/zR)\rightarrow X'$.
\item
For any closed point $a\in X$, any $w\in M(R^h)$ with $\partial w/\partial z\in R^{h\times}$ and any parameter system $Q'$ of the local ring $\mathcal{O}_{X',\pi(a)}$ of $X'$ at $\pi(a)$, $\sigma(a)=M(R)$ and $\{\sigma^*(w)\}\cup\pi^*(Q') $ is a parameter system of the Henselization $\mathcal{O}_{X,a}^h$ of the local ring $\mathcal{O}_{X,a}$ of $X$ at $a$ with $\pi^*(Q')\subset\mathcal{O}_{X,a}$, where $\sigma^*:R^h\rightarrow \mathcal{O}_{X,a}^h$ denotes the homomorphism of local $k$-algebras induced by $\sigma$ on the Henselizations of local rings and $\pi^*:\mathcal{O}_{X',\pi(a)}\rightarrow \mathcal{O}_{X,a}\subset \mathcal{O}_{X,a}^h$ denotes the homomorphism of local $k$-algebras induced by $\pi$.
\end{enumerate}

The lemma below plays the role of a key in our proofs below.

Let $X$ be any finite set.
We define a partial order on $\Map(X,\R)$.
Let $e\in\Map(X,\R)$ and $f\in\Map(X,\R)$ be any elements. We denote $e\leq f$ or $f\geq e$, if $e(x)\leq f(x)$ for any $x\in X$. Obviously the relation $\leq$ is a partial order on $\Map(X,\R)$.
We denote $e< f$ or $f>e$, if $e\leq f$ and $e\neq f$.

\begin{lemma}
\label{BM}
\emph{(Bierstone and Milman~\cite{B88}, p. 25, Lemma 4.7)}

Let $\A\in\Map(P, \Z_0)$, $\B\in\Map(P, \Z_0)$ and $\gamma\in\Map(P, \Z_0)$ be any mappings from $P$ to $\Z_0$ and let $u\in R^{c\times}$, $v\in R^{c\times}$ and $w\in R^{c\times}$ be any elements.

If
$$u\prod_{x\in P}x^{\A(x)}- v\prod_{x\in P}x^{\B(x)}= w\prod_{x\in P}x^{\gamma(x)},$$
then, either $\A\leq\B$, or, $\B\leq\A$ with respect to the partial order $\leq$ on $\Map(P,\R)$.
\end{lemma}

We give the proof of Theorem~\ref{make simple}.

Assume the above $(*)$ and $\dim R\geq 2$.

Consider any element $\phi\in R$ and any $w\in M(R^h)$ such that $\phi\neq 0$, $\partial w/\partial z\in R^{h\times}$, and $\Gamma_+(P_w,\phi)$ is of $w$-Weierstrass type, where $P_w=\{w\}\cup(P-\{z\})$. By $\psi$ we denote a main factor of $(P_w,w,\phi)$.

We take the mapping $a:P-\{z\}\rightarrow\Z_0$ and $h\in\Z_0$ such that $\{\sum_{x\in P-\{z\}}a(x)f^{P_w}_x+hf^{P_w}_z\}$ is the unique $w$-top vertex of $\Gamma_+(P_w,\phi)$. We take $\omega\in R$ with $\phi=\prod_{x\in P-\{z\}}x^{a(x)}\omega$. $\Gamma_+(P_w,\omega)$ is of $w$-Weierstrass type and the unique $w$-top vertex of $\Gamma_+(P_w,\omega)$ is $\{hf^{P_w}_z\}$.

By Henselian Weierstrass Theorem (Hironaka~\cite{H772}.), we know that there exist uniquely $u\in R^{h\times}$ and a mapping
$\omega':\{0,1,\ldots,h-1\}\rightarrow M(R^{\prime h})$ satisfying $\omega=u(w^h+\sum_{i=0}^{h-1}\omega'(i)w^i)$. We take $u\in R^{h\times}$ and a mapping
$\omega':\{0,1,\ldots,h-1\}\rightarrow M(R^{\prime h})$ satisfying $\omega=u(w^h+\sum_{i=0}^{h-1}\omega'(i)w^i)$.

We denote $I=\{i\in\{0,1,\ldots,\hat{h}-1\}|\omega'(i)\neq 0\}$.
$I\subset\{0,1,\ldots,\hat{h}-1\}$.
We put $\omega'(h)=1\in R^{\prime h}$.

Note that for any integers $j$, $k$, $\ell$ with $0\leq k<\ell<j$, $j!/(j-k)\in\Z_+$,
$j!/(j-\ell)\in\Z_+$ and $ j!/(j-k)< j!/(j-\ell)$.

For any $j\in I\cup\{h\}$ we denote
\begin{equation*}\begin{split}
K(j)&=\{(k,\ell)|k\in I, \ell\in I, k<\ell<j,\\
&\qquad\quad\omega'(k)^{j!/(j-k)} \omega'(j)^{(j!/(j-\ell))-(j!/(j-k))}\neq\omega'(\ell)^{j!/(j-\ell)}\}\subset I\times I.
\end{split}\end{equation*}
We put
\begin{equation*}\begin{split}
\phi''&=\prod_{x\in P-\{z\}}x\prod_{i\in I}\omega'(i)\\
&\qquad \prod_{j\in I\cup\{h\}}\prod_{(k,\ell)\in K(j)}
(\omega'(k)^{j!/(j-k)} \omega'(j)^{(j!/(j-\ell))-(j!/(j-k))}-\omega'(\ell)^{j!/(j-\ell)}).
\end{split}\end{equation*}
$\phi''\in M(R^{\prime h})$. $\phi''\neq 0$.

By Lemma~\ref{coordinate change}.8, we know that there exist $\phi'\in M(R')$ such that $\phi'R'=(\phi'' R^{\prime h})\cap R'$ and  $\phi'\neq 0$. We take  $\phi'\in M(R')$ such that $\phi'R'=(\phi'' R^{\prime h})\cap R'$ and  $\phi'\neq 0$.

By $(*)$ we know that 
there exists a composition $\sigma':X'\rightarrow\Spec(R')$ of finite blowing-ups with center in a closed irreducible smooth subscheme such that the divisor on $X'$ defined by the pull-back $\sigma^{\prime *}(\phi')\in\mathcal{O}_{X'}(X')$ of $\phi'$ by $\sigma'$ has normal crossings.
We take a composition $\sigma':X'\rightarrow\Spec(R')$ of finite blowing-ups with center in a closed irreducible smooth subscheme such that the divisor on $X'$ defined by the pull-back $\sigma^{\prime *}(\phi')\in\mathcal{O}_{X'}(X')$ of $\phi'$ by $\sigma'$ has normal crossings.

We consider a morphism $\Spec(R)\rightarrow\Spec(R')$ induced by the inclusion ring homomorphism $R'\rightarrow R$, the product scheme $X=X'\times_{\Spec(R')}\Spec(R)$, the projection $\sigma:X\rightarrow\Spec(R)$, and the projection $\pi:X\rightarrow X'$.
The structure sheaves of the schemes $X$ and $X'$ are denoted by $\mathcal{O}_X$ and $\mathcal{O}_{X'}$ respectively.

Consider any closed point $a\in X$ with $\sigma(a)=M(R)$. 
We have the homomorphism $\sigma^*:R\rightarrow \mathcal{O}_{X,a}$ of local $k$-algebras induced by $\sigma$ from $R$ to the local ring $\mathcal{O}_{X,a}$ of $X$ at $a$, the homomorphism $\pi^*:\mathcal{O}_{X',\pi(a)}\rightarrow \mathcal{O}_{X,a}$ of local $k$-algebras induced by $\pi$ from the local ring $\mathcal{O}_{X',\pi(a)}$ of $X'$ at $\pi(a)$ to $\mathcal{O}_{X,a}$ and the homomorphism $\sigma^{\prime *}:R'\rightarrow \mathcal{O}_{X',\pi(a)}$ of local $k$-algebras induced by $\sigma'$ from $R'$ to $\mathcal{O}_{X',\pi(a)}$, and $\sigma^*$ induces a homomorphism $\sigma^*:R^h\rightarrow \mathcal{O}_{X,a}^h$ of local $k$-algebras from the Henselization $R^h$ of $R$ to the Henselization $\mathcal{O}_{X,a}^h$ of $\mathcal{O}_{X,a}$.

$\sigma^*$ and $\sigma^{\prime *}$ are injective and we have $\sigma^*(\phi)\neq 0$ and $\sigma^{\prime *}(\phi')\neq 0$.

We take any parameter system $\bar{Q}$ of $\mathcal{O}_{X',\pi(a)}$ such that $\sigma^{\prime *}(\phi')\in\mathcal{O}_{X',\pi(a)}$ has normal crossing over $\bar{Q}$. We denote $\bar{P}_w=\{\sigma^*(w)\}\cup \pi^*(\bar{Q})$ and by $\bar{\psi}$ we denote a main factor of $(\bar{P}_w, \sigma^*(w), \sigma^*(\phi))$. 

\noindent $1.$\quad
For any $x\in P_w-\{w\}$, $\sigma^{\prime *}(\phi')\in\sigma^{\prime *}(x) \mathcal{O}_{X',\pi(a)}$ by definition of $\phi'$. Thus, $\sigma^{\prime *}(x)$ has normal crossing over $\bar{Q}$ for any $x\in P_w-\{w\}=P-\{z\}$.

\noindent $2.$\quad
We consider the element $\sigma^*(\omega)\in\mathcal{O}_{X, a} $.

$$\sigma^* (\omega)=\sigma^*(u)(\sigma^*(w)^h+\sum_{i=0}^{h-1}\sigma^* (\omega'(i)) \sigma^*(w)^i).$$
We know that the Newton polyhedron $\Gamma_+(\bar{P}_w, \sigma^*(\omega))$ is of $\sigma^*(w)$-Weierstrass type, and the unique $\sigma^*(w)$-vertex of $\Gamma_+(\bar{P}_w, \sigma^*(\omega))$ is equal to $\{hf^{\bar{P}_w}_{\sigma^*(w)}\}$.

Consider any $i\in I$. 

$0\neq \sigma^{\prime *}(\omega'(i))\in M(\mathcal{O}_{X',\pi(a)}^h)$.
Since $\sigma^{\prime *}(\phi')\in\sigma^{\prime *}(\omega'(i)) \mathcal{O}_{X',\pi(a)}$ by definition of $\phi'$, $\sigma^{\prime *} (\omega'(i))$ has normal crossings over $\bar{P}_w-\{\sigma^*(w)\}$. 
We take the unique pair of an element $\bar{v}(i)\in(\mathcal{O}_{X',\pi(a)}^h)^\times$ and a mapping $\bar{c}(i):\bar{P}_w-\{\sigma^*(w)\}\rightarrow\Z_0$ satisfying 
$\sigma^{\prime *}(\omega'(i))= \bar{v}(i)\prod_{\bar{x}\in \bar{P}_w-\{\sigma^*(w)\}}\bar{x}^{\bar{c}(i)(\bar{x})}$.
We know $\bar{c}(i)\neq 0$, since $\sigma^{\prime *}(\omega'(i))\in M(\mathcal{O}_{ X',\pi(a)}^h)$.

We put $\bar{c}(h)=0\in\Map(\bar{P}_w-\{\sigma^*(w)\},\Z_0)$.
$\sigma^{\prime *} (\omega'(h))=1=$\hfill\break$\prod_{\bar{x}\in \bar{P}_w-\{\sigma^*(w)\}}\bar{x}^{\bar{c}(h)(\bar{x})}$.

Assume that a non-negative integer $m$ with $m<\sharp I$ and a mapping $\nu:\{0,1,\ldots,$\break$m\}\rightarrow I\cup\{h\}$ satisfying the following five conditions are given:
\begin{enumerate}
\item
$\nu$ is injective and it reverses the order.
\item
$\nu(0)=h$.
\item
If $m>0$, then for any $i\in \{1,2,\ldots,m\}$ and for any $j\in I$ with $\nu(i)< j< \nu(i-1)$, 
$((\nu(i-1)-j) /(\nu(i-1)-\nu(i)))\bar{c}\nu(i)
 +((j-\nu(i))/(\nu(i-1)-\nu(i)))\bar{c}\nu(i-1)\leq \bar{c}(j)$.
\item
If $m>0$, then for any $i\in \{1,2,\ldots,m\}$ and for any $j\in I$ with $j<\nu(i)$, 
$((\nu(i-1)-j) /(\nu(i-1)-\nu(i)))\bar{c}\nu(i)
 +((j-\nu(i))/(\nu(i-1)-\nu(i)))\bar{c}\nu(i-1)< \bar{c}(j)$.
\item
$\nu(m)\neq \min (I\cup\{h\})$.
\end{enumerate}

Note that if $m=0$, then there exists  uniquely a mapping $\nu:\{0,1,\ldots,m\}\rightarrow I\cup\{h\}$ satisfying four conditions except the last one of the above five.
If $m=0$, a mapping $\nu:\{0\}=\{0,1,\ldots,m\}\rightarrow I\cup\{h\}$ satisfying $\nu(0)=h$ satisfies four conditions except the last one of the above five. This mapping satisfies the above five conditions, if and only if, $I\neq\emptyset$.

Consider any $(k,\ell)\in K\nu(m)$.
$0\leq k<\ell<\nu(m)$.

$\sigma^{\prime *}(\phi')\in\sigma^{\prime *} (\omega'(k)^{\nu(m)!/(\nu(m)-k)} \omega'\nu(m)^{(\nu(m)!/(\nu(m)-\ell))-(\nu(m)!/(\nu(m)-k))}-$\hfill\break$\omega'(\ell)^{\nu(m)!/(\nu(m)-\ell)})\mathcal{O}_{X',\pi(a)}^h$ by definition of $\phi'$, and the element $\sigma^{\prime *} ($\hfill\break$\omega'(k)^{\nu(m)!/(\nu(m)-k)} \omega'\nu(m)^{(\nu(m)!/(\nu(m)-\ell))-(\nu(m)!/(\nu(m)-k))}-\omega'(\ell)^{\nu(m)!/(\nu(m)-\ell)}) $\hfill\break$\in \mathcal{O}_{X',\pi(a)}^h$ has normal crossings over $\bar{P}_w-\{\sigma^*(w)\}$.
By Lemma~\ref{BM}, we know that either $(\bar{c}(k)-\bar{c}\nu(m))/(\nu(m)-k)\leq (\bar{c}(\ell) -\bar{c}\nu(m))/(\nu(m)-\ell)$ or $(\bar{c}(\ell) -\bar{c}\nu(m))/(\nu(m)-\ell)\leq (\bar{c}(k)-\bar{c}\nu(m))/(\nu(m)-k)$ holds.

Consider any $(k,\ell)\in I\times I$ such that $0\leq k<\ell<\nu(m)$ and $(k,\ell)\not\in K\nu(m)$.

We know $\sigma^{\prime *} (\omega'(k)^{\nu(m)!/(\nu(m)-k)} \omega'\nu(m)^{(\nu(m)!/(\nu(m)-\ell))-(\nu(m)!/(\nu(m)-k))}-$\hfill\break$\omega'(\ell)^{\nu(m)!/(\nu(m)-\ell)})=0$ and $(\bar{c}(k)-\bar{c}\nu(m))/(\nu(m)-k)= (\bar{c}(\ell) -\bar{c}\nu(m))/(\nu(m)-\ell)$.

For any $(k,\ell)\in I\times I$ with $0\leq k<\ell<\nu(m)$, either $(\bar{c}(k)-\bar{c}\nu(m))/(\nu(m)-k)\leq (\bar{c}(\ell) -\bar{c}\nu(m))/(\nu(m)-\ell)$ or $(\bar{c}(\ell) -\bar{c}\nu(m))/(\nu(m)-\ell)\leq (\bar{c}(k)-\bar{c}\nu(m))/(\nu(m)-k)$ holds.

Note that $\{(\bar{c}(k) -\bar{c}\nu(m))/(\nu(m)-k)|k\in I, k<\nu(m)\}\neq\emptyset$ by the last condition of the above five.
We know that the set $\{(\bar{c}(k) -\bar{c}\nu(m))/(\nu(m)-k)|k\in I, k<\nu(m)\}$ has the minimum element $\min \{(\bar{c}(k) -\bar{c}\nu(m))/(\nu(m)-k)|k\in I, k<\nu(m)\}$ with respect to the partial order $\leq$.

Putting \begin{equation*}\begin{split}
\nu(m+1)&=\min\{j\in I|j<\nu(m), (\bar{c}(j) -\bar{c}\nu(m))/(\nu(m)-j)=\\
&\qquad\min \{(\bar{c}(k) -\bar{c}\nu(m))/(\nu(m)-k)|k\in I, k<\nu(m)\}\},\end{split}\end{equation*}
we define an extension $\nu:\{0,1,\ldots,m+1\}\rightarrow I\cup\{h\}$ of $\nu:\{0,1,\ldots,m\}\rightarrow I\cup\{h\}$.
$\nu(m+1)<\nu(m)$. 

Consider any $j\in I$ with $\nu(m+1)<j<\nu(m)$. By how to choose $\nu(m+1)$ we know
$(\bar{c}\nu(m+1) -\bar{c}\nu(m))/(\nu(m)-\nu(m+1))\leq(\bar{c}(j) -\bar{c}\nu(m))/(\nu(m)-j)$.
It follows
$((\nu(m)-j)/(\nu(m)-\nu(m+1))\bar{c}\nu(m+1)+((j-\nu(m+1))/ (\nu(m)-\nu(m+1))\bar{c}\nu(m)\leq \bar{c}(j)$.

Consider any $j\in I$ with $j<\nu(m+1) $. By how to choose $\nu(m+1)$ we know
$(\bar{c}\nu(m+1) -\bar{c}\nu(m))/(\nu(m)-\nu(m+1))<(\bar{c}(j) -\bar{c}\nu(m))/(\nu(m)-j)$.
It follows
$((\nu(m)-j)/(\nu(m)-\nu(m+1))\bar{c}\nu(m+1)+((j-\nu(m+1))/ (\nu(m)-\nu(m+1))\bar{c}\nu(m)< \bar{c}(j)$.

We know that if $\nu(m+1)\neq \min(I\cup\{h\})$, then the extension $\nu:\{0,1,\ldots,m+1\}\rightarrow I\cup\{h\}$ satisfies the above five conditions.
If $\nu(m+1)= \min(I\cup\{h\})$, then the extension satisfies the four conditions except the last one of the above five.

By induction we know that there exists uniquely a pair of a positive integer $m$ and a mapping $\nu:\{0,1,\ldots,m\}\rightarrow I\cup\{h\}$ satisfying the following five conditions:
\begin{enumerate}
\item
$\nu$ is injective and it reverses the order.
\item
$\nu(0)=h$.
\item
For any $i\in \{1,2,\ldots,m\}$ and for any $j\in I$ with $\nu(i)< j< \nu(i-1)$, 
$((\nu(i-1)-j) /(\nu(i-1)-\nu(i)))\bar{c}\nu(i)
 +((j-\nu(i))/(\nu(i-1)-\nu(i)))\bar{c}\nu(i-1)\leq \bar{c}(j)$.
\item
For any $i\in \{1,2,\ldots,m\}$ and for any $j\in I$ with $j<\nu(i)$, 
$((\nu(i-1)-j) /(\nu(i-1)-\nu(i)))\bar{c}\nu(i)
 +((j-\nu(i))/(\nu(i-1)-\nu(i)))\bar{c}\nu(i-1)< \bar{c}(j)$.
\item
$\nu(m)= \min(I\cup\{h\})$.
\end{enumerate}

We take the unique pair of a positive integer $m$ and a mapping $\nu:\{0,1,\ldots,m\}\rightarrow I\cup\{h\}$ satisfying the above five conditions.

We know that $(\bar{c}\nu(i)-\bar{c}\nu(i-1))/(\nu(i-1)-\nu(i))<
(\bar{c}\nu(i+1)-\bar{c}\nu(i))/(\nu(i)-\nu(i+1))$ for any $i\in\{1,2,\ldots,m-1\}$, if $m\geq 2$.

By Lemma~\ref{Newton2}.$15$ we know that the Newton polyhedron $\Gamma_+(\bar{P}_w,\sigma^* (\omega))$ is $\sigma^*(w)$-simple.

By 1, we know that there exists $\bar{u}\in\mathcal{O}_{X,a}^\times$ and a mapping $\bar{a}:\bar{P}_w-\{\sigma^*(w)\}\rightarrow\Z_0$ with $\sigma^*(\phi)=\bar{u}\prod_{\bar{x}\in\bar{P}_w-\{\sigma^*(w)\}}\bar{x}^{\bar{a}(\bar{x})}\sigma^*(\omega)$. We know that the Newton polyhedron $\Gamma_+(\bar{P}_w,\sigma^* (\phi))$ is $\sigma^*(w)$-simple.

\noindent $3.$\quad
We take $v\in R^\times$, a mapping $a:P_w-\{w\}\rightarrow \Z_0$, a finite subset $\Lambda\subset M(R)$ and a mapping $b:\Lambda\rightarrow\Z_+$ satisfying the following three conditions:
\begin{enumerate}
\item $\phi=v(\prod_{x\in P_w-\{w\}}x^{a(x)})(\prod_{\lambda\in\Lambda}\lambda^{b(\lambda)})\psi$.
\item $\partial\lambda/\partial w\in R^{h\times}$ for any $\lambda\in\Lambda$.
\item If $\lambda\in\Lambda$, $\lambda'\in\Lambda$, $s\in R^\times$ and $\lambda=s\lambda'$, then $s=1$ and $\lambda=\lambda'$.
\end{enumerate}
By 1, we know that there exist $\bar{v}\in \mathcal{O}_{X,a}^\times$ and a mapping $\bar{a}:\bar{P}_w-\{\sigma^*(w)\}\rightarrow\Z_0$ with $\sigma^*(v(\prod_{x\in P_w-\{w\}}x^{a(x)}))=\bar{v}(\prod_{\bar{x}\in\bar{P}_w-\{\sigma^*(w)\}}\bar{x}^{\bar{a}(\bar{x})})$. We take $\bar{v}\in \mathcal{O}_{X,a}^\times$ and a mapping $\bar{a}:\bar{P}_w-\{\sigma^*(w)\}\rightarrow\Z_0$ satisfying this equality.

We have
$\sigma^*(\phi)= \bar{v}(\prod_{\bar{x}\in\bar{P}_w-\{\sigma^*(w)\}}\bar{x}^{\bar{a}(\bar{x})}) (\prod_{\lambda\in\Lambda}\sigma^*(\lambda)^{b(\lambda)})\sigma^*(\psi)$.

Since $\partial\sigma^*(\lambda)/\partial\sigma^*(w)=\sigma^*(\partial\lambda/\partial w)\in\mathcal{O}_{X,a}^{h\times}$ for any $\lambda\in\Lambda$, we know that the main factor $\bar{\psi}$ of $(\bar{P}_w, \sigma^*(w),\sigma^*(\phi))$ is a main factor of $(\bar{P}_w, \sigma^*(w),\sigma^*(\psi))$. Thus, we have $\Inv(\bar{P}_w, \sigma^*(w),\sigma^*(\phi))=\deg(\bar{P}_w, \sigma^*(w),\bar{\psi})\leq\deg(\bar{P}_w, \sigma^*(w),\sigma^*(\psi))=\deg(P_w,w,\psi)=\Inv(P_w,w,\phi)$. 

$\Inv(\bar{P}_w, \sigma^*(w),\sigma^*(\phi))= \Inv(P_w,w,\phi)$, if and only if, $\sigma^*(\psi)=s\bar{\psi}$ for some $s\in\mathcal{O}_{X,a}^\times$, if and only if, $\Gamma_+(\bar{P}_w, \sigma^*(\psi))= \Gamma_+(\bar{P}_w, \bar{\psi})$.

\noindent $4.$\quad
Assume that $\Inv(\bar{P}_w, \sigma^*(w),\sigma^*(\phi))= \Inv(P_w,w,\phi)$ and $\Gamma_+(P_w,\psi)$ has no $w$-removable faces.

We have $\Gamma_+(\bar{P}_w, \sigma^*(\psi))= \Gamma_+(\bar{P}_w, \bar{\psi})$.
$\Gamma_+(\bar{P}_w, \sigma^*(\psi))$ has no $\sigma^*(w)$-removable faces, if and only if, $\Gamma_+(\bar{P}_w, \bar{\psi})$ has no $\sigma^*(w)$-removable faces.

We denote $\bar{h}=\Inv(\bar{P}_w, \sigma^*(w),\sigma^*(\phi))\in\Z_0$. $\{\bar{h}f^{\bar{P}_w}_{\sigma^*(w)}\}$ is the unique $\sigma^*(w)$-top vertex of $\Gamma_+(\bar{P}_w, \sigma^*(\psi))$. 

Assume moreover, that $\Gamma_+(\bar{P}_w, \sigma^*(\psi))$ has a $\sigma^*(w)$-removable face. We will deduce a contradiction. 

It follows $\bar{h}>0$. Take any $\sigma^*(w)$-removable face $\bar{F}$ of dimension one of $\Gamma_+(\bar{P}_w, $\break$\sigma^*(\psi))$.
$\Stab(\bar{F})=\{0\}$.
$\Delta^\circ(\bar{F}, \Gamma_+(\bar{P}_w, \sigma^*(\psi)))\subset (\Map(\bar{P}_w,\R_0)^\vee)^\circ$.
Take any $\bar{\lambda}\in $\hfill\break$(\Map(\bar{P}_w,\R_0)^\vee)\cap\sum_{\bar{x}\in\bar{P}_w
-\{\sigma^*(w)\}}\R f^{\bar{P}_w\vee}_{\bar{x}}$ such that $\bar{\lambda}+ f^{\bar{P}_w\vee}_{\sigma^*(w)}\in\Delta^\circ(\bar{F}, \Gamma_+(\bar{P}_w, \sigma^*(\psi$\break$)))$.
$\bar{h}f^{\bar{P}_w}_{\sigma^*(w)}\in \bar{F}$ and $\Ord(\bar{P}_w,\bar{\lambda}+ f^{\bar{P}_w\vee}_{\sigma^*(w)},\sigma^* (\psi))=\langle \bar{\lambda}+ f^{\bar{P}_w\vee}_{\sigma^*(w)}, \bar{h}f^{\bar{P}_w}_{\sigma^*(w)}\rangle=\bar{h}$.
Take $\bar{\B}\in k-\{0\}$, $\bar{\gamma}\in k-\{0\}$ and $\bar{e}\in\Map(\bar{P}_w-\{\sigma^*(w)\},\Z_0)$ satisfying
$\In(\bar{P}_w,\bar{\lambda}+ f^{\bar{P}_w\vee}_{\sigma^*(w)},\sigma^*(\psi))=\Ps(\bar{P}_w,F, \sigma^*(\psi))=\bar{\B}(\sigma^*(w)+\bar{\gamma}\prod_{\bar{x}\in\bar{P}_w-\{\sigma^*(w)\}}\bar{x}^{\bar{e}(\bar{x})})^{\bar{h}}=\bar{\B}(\bar{z}^{\bar{h}}+\sum_{i=1}^{\bar{h}}\binom{\bar{h}}{i}\bar{\gamma}^i
\prod_{\bar{x}\in\bar{P}_w-\{\sigma^*(w)\}}\bar{x}^{i\bar{e}(\bar{x})}\sigma^*(w)^{\bar{h}-i}$.

We define an element $\lambda\in\Map(P_w,\R_0)^\vee\cap\sum_{x\in P_w-\{w\}}\R f^{P_w\vee}_x$ by putting $\langle \lambda, $\break$f^{P_w}_x\rangle=\Ord(\bar{P}_w-\{\sigma^*(w)\},\bar{\lambda},\sigma^*(x))\in\R_0$ for any $x\in P_w-\{w\}$.
$\lambda+f^{P_w\vee}_w\in \Map(P_w,\R_0)^\vee$.
We denote $F=\Delta(\lambda+f^{P_w\vee}_w,\Gamma_+(P_w,\psi))\in\mathcal{F}(\Gamma_+(P_w,\psi))$.
Since $\sigma^{\prime *}(x)$ has normal crossing over $\bar{Q}$ for any $x\in P_w-\{w\}$, $\Ord(P_w, \lambda+f^{P_w\vee}_w,\psi)= \Ord(\bar{P}_w,\bar{\lambda}+ f^{\bar{P}_w\vee}_{\sigma^*(w)},\sigma^*(\psi))=\bar{h}$ and
$\sigma^*(\In(P_w, \lambda+f^{P_w\vee}_w,\psi))= \In(\bar{P}_w,\bar{\lambda}+ f^{\bar{P}_w\vee}_{\sigma^*(w)},\sigma^*(\psi$\break$))= \bar{\B}(\sigma^*(w)^{\bar{h}}+\sum_{i=1}^{\bar{h}}\binom{\bar{h}}{i}\bar{\gamma}^i
\prod_{\bar{x}\in\bar{P}_w-\{\sigma^*(w)\}}\bar{x}^{i\bar{e}(\bar{x})}\sigma^*(w)^{\bar{h}-i})$.

We take $u_0\in R^{c\times}$ and a mapping $\psi':\{0,1,\ldots,\bar{h}-1\}\rightarrow M(R^{\prime c})$ satisfying $\psi=u_0(w^{\bar{h}}+\sum_{i=0}^{\bar{h}-1}\psi'(i)w^i)$.

We know that $\bar{h}f^{P_w}_w\in F$ and $\Ord(P_w-\{w\},\lambda, \psi'(\bar{h}-i))\geq i\langle \bar{\lambda}, \bar{e}\rangle$ for any $i\in\{1,2,\ldots, \bar{h}\}$. We put
\begin{equation*}
\hat{\psi}(\bar{h}-i)=
\begin{cases}
\In(P_w-\{w\},\lambda, \psi'(\bar{h}-i))& \text{if $\Ord(P_w-\{w\},\lambda, \psi'(\bar{h}-i))= i\langle \bar{\lambda}, \bar{e}\rangle$},\\
0&\text{if $\Ord(P_w-\{w\},\lambda, \psi'(\bar{h}-i))> i\langle \bar{\lambda}, \bar{e}\rangle$}.
\end{cases}\end{equation*}
We have $\In(P_w, \lambda+f^{P_w\vee}_w,\psi)= \In(P_w, \lambda+f^{P_w\vee}_w, u_0)(w^{\bar{h}}+\sum_{i=1}^{\bar{h}}\hat{\psi}(\bar{h}-i) w^{\bar{h}-i})$ and
$\sigma^*(\In(P, \lambda+f^{P_w\vee}_w,\psi))= \sigma^*(\In(P, \lambda+f^{P_w\vee}_w,u_0)) (\sigma^*(w)^{\bar{h}}+\sum_{i=1}^{\bar{h}}\sigma^*(\hat{\psi}(\bar{h}-i)) \sigma^*(w)^{\bar{h}-i})$.
It follows that $\binom{\bar{h}}{i}\bar{\gamma}^i
\prod_{\bar{x}\in\bar{P}_w-\{\sigma^*(w)\}}\bar{x}^{i\bar{e}(\bar{x})}=
\sigma^*(\hat{\psi}(\bar{h}-i))$ for any $i\in\{1,2,\ldots,\bar{h}\}$.

By $\bar{1}\in k$ we denote the identity element of the field $k$.

We consider the case where the characteristic number of the field $k$ is equal to $0$.

$\binom{\bar{h}}{1}\bar{1}=\bar{h}\bar{1}\neq 0$. We put $\hat{\chi}=\hat{\psi}(\bar{h}-1)/(\bar{h}\bar{1})\in M(R^{\prime c})$.
We have \hfill\break$\bar{\gamma}\prod_{\bar{x}\in\bar{P}_w-\{\sigma^*(w)\}}\bar{x}^{\bar{e}(\bar{x})}=
\sigma^*(\hat{\chi})$.

We consider the case where the characteristic number of the field $k$ is positive.
By $p$ we denote the characteristic number of the field $k$. The integer $p$ is a prime number.
We take the unique pair of $\delta\in\Z_0$ and $\hat{h}\in\Z_+$ such that $\bar{h}=p^\delta\hat{h}$ and $\hat{h}$ is not a multiple of $p$.
$\binom{\bar{h}}{p^\delta}\bar{1}=\hat{h}\bar{1}\neq 0$.
We have $(\bar{\gamma}\prod_{\bar{x}\in\bar{P}_w-\{\sigma^*(w)\}}\bar{x}^{\bar{e}(\bar{x})})^{p^\delta}=\sigma^*(\hat{\psi}(\bar{h}- p^\delta)/(\hat{h}\bar{1}))$ and 
$\hat{\psi}(\bar{h}- p^\delta)/(\hat{h}\bar{1})\in M(R^{\prime c})$.

For any complete regular local ring $A$ such that $A$ contains $k$ as a subring and the residue field $A/M(A)$ is isomorphic to $k$ as $k$-algebras, we denote $A^{p^\delta}=\{s^{p^\delta}|s\in A\}$. $A^{p^\delta}$ is a local $k$-subalgebra of $A$.

$(\bar{\gamma}\prod_{\bar{x}\in\bar{P}_w-\{\sigma^*(w)\}}\bar{x}^{\bar{e}(\bar{x})})^{p^\delta}\in \mathcal{O}_{X,a}^{c\: p^\delta}$.
$\hat{\psi}(\bar{h}- p^\delta)/(\hat{h}\bar{1})\in\sigma^{* -1}(\mathcal{O}_{X,a}^{c\: p^\delta})$.
Now, since $\sigma$ is a composition of finite blowing-ups, $\sigma^{* -1}(\mathcal{O}_{X,a}^{c\: p^\delta})=R^{c\:p^\delta}$.
We know that there exists uniquely an element $\hat{\chi}\in R^{\prime c}$ with $\hat{\chi}^{p^\delta}=\hat{\psi}(\bar{h}- p^\delta)/(\hat{h}\bar{1})$. 
We take the unique element $\hat{\chi}\in R^{\prime c}$ with $\hat{\chi}^{p^\delta}=\hat{\psi}(\bar{h}- p^\delta)/(\hat{h}\bar{1})$.
Since $\hat{\psi}(h- p^\delta)/(\hat{h}\bar{1})\in M(R^{\prime c})$, we know $\hat{\chi}\in M(R^{\prime c})$.
We know $(\bar{\gamma}\prod_{\bar{x}\in\bar{P}_w-\{\sigma^*(w)\}}\bar{x}^{\bar{e}(\bar{x})})^{p^\delta}=\sigma^*( \hat{\chi}^{p^\delta})$ and 
$\bar{\gamma}\prod_{\bar{x}\in\bar{P}_w-\{\sigma^*(w)\}}\bar{x}^{\bar{e}(\bar{x})}=\sigma^*( \hat{\chi})$.

We conclude that there exists $\hat{\chi}\in M(R^{\prime c})$ satisfying $\bar{\gamma}\prod_{\bar{x}\in\bar{P}_w-\{\sigma^*(w)\}\}}\bar{x}^{\bar{e}(\bar{x})}=\sigma^*( \hat{\chi})$ in all cases. We take any  $\hat{\chi}\in M(R^{\prime c})$ satisfying $\bar{\gamma}\prod_{\bar{x}\in\bar{P}_w-\{\sigma^*(w)\}}\bar{x}^{\bar{e}(\bar{x})}=\sigma^*( \hat{\chi})$.

We have 
$\sigma^*(\In(P_w, \lambda+f^{P\vee}_z,\psi))= \In(\bar{P}_w,\bar{\lambda}+ f^{\bar{P}\vee}_{\bar{z}},\sigma^*(\psi))= (\sigma^*(w)+ $\break$\bar{\gamma}\prod_{\bar{x}\in\bar{P}_w-\{\sigma^*(w)\}} \bar{x}^{\bar{e}(\bar{x})})^{\bar{h}}=(\sigma^*(w)+\sigma^*( \hat{\chi}))^{\bar{h}}=\sigma^*((w+\hat{\chi})^{\bar{h}})$.
Since $\sigma^*$ is injective we have
$\Ps(P_w, F,\psi)= \In(P_w, \lambda+f^{P\vee}_z,\psi)= (w+\hat{\chi})^{\bar{h}}$ and we know that the face $F$ of $\Gamma_+(P_w,\psi)$ is $w$-removable.

Since the Newton polyhedron $\Gamma_+(P_w,\psi)$ has no $w$-removable faces, we obtain a contradiction.

We conclude that the Newton polyhedron $\Gamma_+(\bar{P}_w,\bar{\psi})$ has no $\sigma^*(w)$-removable faces.

\noindent $5.$\quad
Assume that $\Inv(\bar{P}_w,\sigma^*(w),\sigma^*(\phi))=\Inv(P_w, w, \phi)$ and $w=z+\chi_0$ where $\chi_0\in M(R^{\prime h})$ is the unique element in Theorem~\ref{erase faces}.2.

Since $\Gamma_+(\bar{P}_w,\sigma^*(\phi))$ is $\sigma^*(w)$-simple by 2 and $\bar{\psi}$ divides $\sigma^*(\phi)$, $\Gamma_+(\bar{P}_w, \bar{\psi})$ is also $\sigma^*(w)$-simple.

By 4, $\Gamma_+(\bar{P}_w, \bar{\psi})$ has no $\sigma^*(w)$-removable faces, since $\Gamma_+(P_w,\psi)= \Gamma_+(P_{z+\chi_0},\psi)$ has no $w$-removable faces.

We denote $\bar{h}=\Inv(\bar{P}_w, \sigma^*(w),\sigma^*(\phi))\in\Z_0$. $\{\bar{h}f^{\bar{P}_w}_{\sigma^*(w)}\}$ is the unique $\sigma^*(w)$-top vertex of $\Gamma_+(\bar{P}_w, \bar{\psi})$. 

Recall that $\Gamma_+(\bar{P}_w, \bar{\psi})\subset\Map(\bar{P}_w,\R)$ and $\{f^{\bar{P}_w }_{\bar{x}}|\bar{x}\in \bar{P}_w \}$ is an $\R$-basis of the vector space $\Map(\bar{P}_w,\R)$. Let $U=\{a\in \Map(\bar{P}_w,\R)|\langle f^{\bar{P}_w \vee}_{\sigma^*(w)},a\rangle<\bar{h}\}$ and $V=\{a\in \Map(\bar{P}_w,\R)|\langle f^{\bar{P}_w \vee}_{\sigma^*(w)},a\rangle=0\}$. We put $\rho(a)=(a-\langle f^{\bar{P}_w\vee}_{\sigma^*(w)},a\rangle f^{\bar{P}_w }_{\sigma^*(w)})/(\bar{h}-\langle f^{\bar{P}_w \vee}_{\sigma^*(w)},a\rangle)\in V$ for any $a\in U$ and we define a mapping $\rho:U\rightarrow V$.

Since $\Gamma_+(\bar{P}_w, \bar{\psi})$ is $\sigma^*(w)$-simple, the convex pseudo polytope $\rho(\Gamma_+(\bar{P}_w, \bar{\psi})\cap U)$ in $V$ has at most one vertex. 

By Theorem~\ref{erase faces}.7, we know that
there exists an element $\bar{w}\in M(\mathcal{O}_{X,a})$ such that $\partial\bar{w}/\partial \sigma^*(w)\in \mathcal{O}_{X,a}^{h\times}$ and if we denote $\bar{P}_{\bar{w}}=\{\bar{w}\}\cup \pi^*(\bar{Q})$, then $\bar{P}_{\bar{w}}$ is a parameter system of $\mathcal{O}_{X,a}$ and $\Gamma_+(\bar{P}_{\bar{w}},\bar{\psi})$ has no $\bar{w}$-removable faces. 
We take any element $\bar{w}\in M(\mathcal{O}_{X,a})$ satisfying these conditions and we denote $\bar{P}_{\bar{w}}=\{\bar{w}\}\cup \pi^*(\bar{Q})$. We have 
$\Inv(\bar{P}_{\bar{w}},\bar{w},\sigma^*(\phi))=\Inv(\bar{P}_w,\sigma^*(w),\sigma^*(\phi))=\Inv(P_w, w, \phi)$ by Lemma~\ref{main factor}.7.(b).

\noindent $6.$\quad
Assume $\Inv(P_w,w,\phi)=0$. It follows $\psi\in R^\times$.
By 3, we know $\Inv(\bar{P}_w,\sigma^*(w),$\break$\sigma^*(\phi))=0$. 

We take $v\in R^\times$, a mapping $a:P_w-\{w\}\rightarrow \Z_0$, a finite subset $\Lambda\subset M(R)$ and a mapping $b:\Lambda\rightarrow\Z_+$ satisfying the following three conditions:
\begin{enumerate}
\item $\phi=v(\prod_{x\in P_w-\{w\}}x^{a(x)})(\prod_{\lambda\in\Lambda}\lambda^{b(\lambda)})$.
\item $\partial\lambda/\partial w\in R^{h\times}$ for any $\lambda\in\Lambda$.
\item If $\lambda\in\Lambda$, $\lambda'\in\Lambda$, $s\in R^\times$ and $\lambda=s\lambda'$, then $s=1$ and $\lambda=\lambda'$.
\end{enumerate}
By definition $\sharp\Lambda=\Inv 2(P_w,w,\phi)$.

By 1, we know that there exist $\bar{v}\in \mathcal{O}_{X,a}^\times$ and a mapping $\bar{a}:\bar{P}_w-\{\sigma^*(w)\}\rightarrow\Z_0$ with $\sigma^*(v(\prod_{x\in P_w-\{w\}}x^{a(x)}))=\bar{v}(\prod_{\bar{x}\in\bar{P}_w-\{\sigma^*(w)\}}\bar{x}^{\bar{a}(\bar{x})})$. We take $\bar{v}\in \mathcal{O}_{X,a}^\times$ and a mapping $\bar{a}:\bar{P}_w-\{\sigma^*(w)\}\rightarrow\Z_0$ satisfying this equality.

We have
$\sigma^*(\phi)= \bar{v}(\prod_{\bar{x}\in\bar{P}_w-\{\sigma^*(w)\}}\bar{x}^{\bar{a}(\bar{x})}) (\prod_{\lambda\in\Lambda}\sigma^*(\lambda)^{b(\lambda)})$.

Consider any $\lambda\in\Lambda$. $\partial\sigma^*(\lambda)/\sigma^*(w)=\sigma^*(\partial \lambda/\partial w)\in\mathcal{O}_{X,a}^\times$.
By condition 2 above, we know that there exists uniquely an element $v_\lambda\in R^{c\times}$ and $\mu_\lambda\in M(R^{\prime c})$ with $\lambda=v_\lambda(w+\mu_\lambda)$. We take $v_\lambda\in R^{c\times}$ and $\mu_\lambda\in M(R^{\prime c})$ with $\lambda=v_\lambda(w+\mu_\lambda)$ for any $\lambda\in\Lambda$.

It follows from condition 3 above that if $\lambda\in\Lambda$, $\lambda'\in\Lambda$, and $\mu_\lambda=\mu_{\lambda'}$, then $\lambda=\lambda'$.

Assume that $\lambda\in\Lambda$, $\lambda'\in\Lambda$, $\bar{s}\in\mathcal{O}_{X,a}^\times$ and $\sigma^*(\lambda)=\bar{s}\sigma^*(\lambda')$.
We have $\sigma^*(v_\lambda)(\sigma^*(w)+\sigma^*(\mu_\lambda))=\bar{s}\sigma^*(v_{\lambda'})(\sigma^*(w)+\sigma^*(\mu_{\lambda'}))$, $\sigma^*(\mu_\lambda)\in\pi^*(M(\mathcal{O}_{X',\pi(a)}^c))$ and $\sigma^*(\mu_{\lambda'})\in\pi^*(M(\mathcal{O}_{X',\pi(a)}^c))$. By Weierstrass' Preparation Theorem, we have $\sigma^*(\mu_\lambda)= \sigma^*(\mu_{\lambda'})$. Since $\sigma^*:R^c\rightarrow \mathcal{O}_{X,a}^c$ is injective, we have $\mu_\lambda=\mu_{\lambda'}$ and $\lambda=\lambda'$.

It follows $\Inv 2(\bar{P}_w,\sigma^*(w),\sigma^*(\phi))=\sharp\Lambda=\Inv 2(P_w,w,\phi)$.

Obviously, if $w$ divides $\phi$, then $\sigma^*(w)$ divides $\sigma^*(\phi)$.

We know that Theorem~\ref{make simple} holds.

We give the proof of Lemma~\ref{make Weierstrass type}.

Consider any $\phi\in R$ with $\phi\neq 0$. Let $h=\Ord(\phi)\in\Z_0$.

Let $\delta_0=\sum_{x\in P}f^{P\vee}_x\in(\Map(P,\R_0)^\vee)^\circ\cap\Map(P,\Z)^*$ and let $\bar{\phi}=\In(P, \delta_0,\phi)\in k[P]$. $\Ord(P, \delta_0,x)=1$ for any $x\in P$. $\bar{\phi}$ is a homogeneous polynomial of degree $h$. 

Since $k$ is an algebraically closed field and $k$ is an infinite field, there exists a mapping $\A:P-\{z\}\rightarrow k$ such that $\theta_\A(\bar{\phi})\neq 0$, where $\theta_\A: k[P]\rightarrow k$ is a homomorphism of $k$-algebras such that $\theta_\A(z)=1$ and $\theta_\A(x)=\A(x)$ for any $x\in P-\{z\}$.

We take any mapping $\A:P-\{z\}\rightarrow k$ such that $\theta_\A(\bar{\phi})\neq 0$ and we put $c=\theta_\A(\bar{\phi})\in k-\{0\}$ and $\bar{P}=\{z\}\cup\{x-\A(x)z|x\in P-\{z\}\}$. Obviously $\bar{P}$ is a parameter system of $R$ containing $z$. $k[\bar{P}]=k[P]$, $\theta_\A(z)=1$ and $\theta_\A(\bar{x})=0$ for any $\bar{x}\in\bar{P}-\{z\}$. $\bar{\phi}-cz^h\in(\bar{P}-\{z\}) k[\bar{P}]$.

Let $\delta_1=\sum_{\bar{x}\in\bar{P}-\{z\}}f^{\bar{P}\vee}_{\bar{x}}\in\Map(\bar{P},\R_0)^\vee\cap\Map(\bar{P},\Z)^*$. $\Ord(\bar{P}, \delta_1,\bar{x})=1$ for any $\bar{x}\in \bar{P}-\{z\}$ and $\Ord(\bar{P}, \delta_1,z)=0$. We have $\In(\bar{P},\delta_1,\bar{\phi})=cz^h$ and $\In(\bar{P},\delta_1, \phi)=uz^h$ for some $u\in R^{c\times}$ with $u-c\in M(R^c)$. We know that $\Gamma_+(\bar{P},\phi)$ is of $z$-Weierstrass type, the unique $z$-top vertex is $\{hf^{\bar{P}}_z\}$ and $\Inv(\bar{P},z,\phi)\leq h$.

We conclude that Lemma~\ref{make Weierstrass type} holds.

We give the proof of Corollary~\ref{resolution game}.

Consider a mathematical game with two players A and B.
At the start of the game a pair $(R,\phi)$ of any regular local ring $R$ with $\dim R\geq 1$ such that $R$ contains $k$ as a subring, the residue field $R/M(R)$ is isomorphic to $k$ as $k$-algebras and $R$ is a localization of a finitely generated $k$-algebra, and any non-zero element $\phi\in R$ is given. We play our game repeating the following step. Before the first step we put $(S,\psi)=(R,\phi)$: At the start of each step, player A chooses a composition $\sigma:X\rightarrow\Spec(S)$ of finite blowing-ups with center in a closed irreducible smooth subscheme. Then, player B chooses a closed point $a\in X$ with $\sigma(a)=M(S)$. We have a morphism $\sigma^*:S\rightarrow\mathcal{O}_{X,a}$ of local $k$-algebras induced by $\sigma$. If the element $\sigma^*(\psi)\in \mathcal{O}_{X,a}$ has normal crossings, then the palyer A wins. Otherwise we proceed to the next step after replacing the pair $(S,\psi)$ by the pair $(\mathcal{O}_{X,a}, \sigma^*(\psi))$.

We call the above game the \emph{main game}. We would like to show that player A can always win the main game after finite steps. To help consideration, we consider three of similar games. The first is called the \emph{invariant game}, the second is called the \emph{simplifying game}, and the third is called the \emph{removing game}. Three have two players A and B.

At the start of the three games a quadruplet $(R,P,z,\phi)$ of any regular local ring $R$ with $\dim R\geq 1$ such that $R$ contains $k$ as a subring, the residue field $R/M(R)$ is isomorphic to $k$ as $k$-algebras and $R$ is a localization of a finitely generated $k$-algebra, a parameter system $P$ of $R$, an element $z\in P$ and any non-zero element $\phi\in R$ such that $\Gamma_+(P,\phi)$ is of $z$-Weierstrass type is given. We play our game repeating the following step. Before the first step we put $(S,Q,w,\psi)=(R,P,z,\phi)$: At the start of each step, player A chooses a composition $\sigma:X\rightarrow\Spec(S)$ of finite blowing-ups with center in a closed irreducible smooth subscheme. Then, player B chooses a closed point $a\in X$ with $\sigma(a)=M(S)$, and successively player A chooses a parameter system $\bar{P}$ of the local ring $\mathcal{O}_{X,a}$ of $X$ at $a$ and an element $\bar{z}\in\bar{P}$. We have a morphism $\sigma^*:S\rightarrow\mathcal{O}_{X,a}$ of local $k$-algebras induced by $\sigma$. 

If $\Gamma_+(\bar{P},\sigma^*(\psi))$ is not of $\bar{z}$-Weierstrass type, then player B wins in three games. 

In the invariant game, if either $\Inv(\bar{P},\bar{z}, \sigma^*(\psi))<\Inv(P,z,\phi)$, or $\Inv(\bar{P},\bar{z}, \sigma^*(\psi$\break$))=\Inv(P,z,\phi)=0$ and $\Inv 2(\bar{P},\bar{z}, \sigma^*(\psi))<\Inv 2(P,z,\phi)$, then player A wins.

In the simplifying game, if the following three conditions are simultaneously satisfied, then player A wins:
\begin{enumerate}
\item Either $\Inv(\bar{P},\bar{z}, \sigma^*(\psi))\leq\Inv(P,z,\phi)$, or $\Inv(\bar{P},\bar{z}, \sigma^*(\psi))=\Inv(P,z,\phi)=0$ and $\Inv 2(\bar{P},\bar{z}, \sigma^*(\psi))\leq\Inv 2(P,z,\phi)$.
\item $\Gamma_+(\bar{P},\sigma^*(\phi))$ is $\bar{z}$-simple.
\item The Newton polyhedron $\Gamma_+(\bar{P},\bar{\psi})$ over $\bar{P}$ of a main factor $\bar{\psi}$ of $(\bar{P},\bar{z}, \sigma^*(\psi))$ has no $\bar{z}$-removable faces.
\end{enumerate}

In the removing game, if the following two conditions are simultaneously satisfied, then player A wins:
\begin{enumerate}
\item Either $\Inv(\bar{P},\bar{z}, \sigma^*(\psi))\leq\Inv(P,z,\phi)$, or $\Inv(\bar{P},\bar{z}, \sigma^*(\psi))=\Inv(P,z,\phi)=0$ and $\Inv 2(\bar{P},\bar{z}, \sigma^*(\psi))\leq\Inv 2(P,z,\phi)$.
\item The Newton polyhedron $\Gamma_+(\bar{P},\bar{\psi})$ over $\bar{P}$ of a main factor $\bar{\psi}$ of $(\bar{P},\bar{z}, \sigma^*(\psi))$ has no $\bar{z}$-removable faces.
\end{enumerate}

If neither A nor B can win the game, we proceed to the next step after replacing the quadruplet $(S,Q,w,\psi)$ by the quadruplet $(\mathcal{O}_{X,a}, \bar{P},\bar{z},\sigma^*(\psi))$.

Consider the following claim $(m)$ corresponding to a positive integer $m$:
If a pair $(R,\phi)$ at the start of the main game satisfies $\dim R\leq m$, then player A can win the main game starting from the pair $(R,\phi)$ after finite steps.

Consider any pair $(R,\phi)$ at the start of the main game. 

Assume $\dim R\leq 1$. Then, $\dim R=1$ and $\phi$ has normal crossings and player A win the game at the first step. We know that claim $(1)$ holds.

Below, we assume $m\geq 2$ and claim $(m-1)$ holds. We will show that claim $(m)$ holds.

Assume $\dim R\leq m$.
If $\dim R\leq m-1$, then by claim $(m-1)$ we can conclude that player A can win the main game starting from the pair $(R,\phi)$ after finite steps.

We consider the case $\dim R=m$.

By Lemma~\ref{make Weierstrass type}, we know that there exist a parameter system $P$ of $R$ and an element $z\in P$ such that the quadruplet $(R,P,z,\phi)$ satisfies conditions at the start of the invariant game, the simplifying game and the removing game. We take $P$ and $z$ such that the quadruplet $(R,P,z,\phi)$ satisfies conditions at the start of the invariant game, the simplifying game and the removing game.

Let $R'$ be the localization of $k[P-\{z\}]$ by the maximal ideal $k[P-\{z\}]\cap M(R)=( P-\{z\}) k[P-\{z\}]$. $\dim R'=\dim R-1<\dim R=m$. 

By the proof of Theorem~\ref{make simple}, we know that there exists a non-zero $\phi'\in R'$ such that if player A can win the main game starting from $(R',\phi')$, then either player A can win the removing game starting from $(R,P,z,\phi)$, or player A can win the invariant game starting from $(R,P,z,\phi)$. By $(m-1)$, we know that player A can win the main game starting from $(R',\phi')$.  We conclude that either player A can win the removing game starting from $(R,P,z,\phi)$, or player A can win the invariant game starting from $(R,P,z,\phi)$.

We consider the case where player A cannot win the invariant game starting from $(R,P,z,\phi)$ so far.
Since player A can win the removing game starting from $(R,P,z,\phi)$, we can assume that the Newton polyhedron $\Gamma_+(P,\psi)$ over $P$ of a main factor $\psi$ of $(P,z, \phi)$ has no $z$-removable faces. We assume this condition. By the proof of Theorem~\ref{make simple}, we know that there exists a non-zero $\phi''\in R'$ such that if player A can win the main game starting from $(R',\phi'')$, then either player A can win the simplifying game starting from $(R,P,z,\phi)$, or player A can win the invariant game starting from $(R,P,z,\phi)$. By $(m-1)$, we know that player A can win the main game starting from $(R',\phi'')$.  We conclude that either player A can win the simplifying game starting from $(R,P,z,\phi)$, or player A can win the invariant game starting from $(R,P,z,\phi)$.

We consider the case where player A cannot win the invariant game starting from $(R,P,z,\phi)$ so far.
Since player A can win the simplifying game starting from $(R,P,z,\phi)$, we can assume that the Newton polyhedron $\Gamma_+(P,\psi)$ over $P$ of a main factor $\psi$ of $(P,z, \phi)$ has no $z$-removable faces and $\Gamma_+(P,\phi)$ is $z$-simple. We assume these conditions.
By Theorem~\ref{main}, we conclude that player A can win the invariant game starting from $(R,P,z,\phi)$ at the first step.

We know that player A can win the invariant game starting from $(R,P,z,\phi)$ in all cases.

Note that $\phi$ has normal crossings, if $\Inv(P,z,\phi)=0$ and $\Inv 2(P,z,\phi)\leq 1$. Thus, we know that if player A can win the invariant game starting from $(R,P,z,\phi)$, then player A can win the main game starting from $(R, \phi)$.

We conclude that player A can win the main game starting from $(R, \phi)$ after finite steps and claim $(m)$ holds.

By induction, we know that claim $(m)$ holds for any positive integer $m$.
We conclude that Corollary~\ref{resolution game} holds.

From Corollary~\ref{resolution game} and the valuation theory, it follows the following claim. We do not explain the valuation theory in this article. (Zariski et al.~\cite{ZS}):
Given any field $\Sigma$ such that $\Sigma$ contains $k$ as a subfield and $\Sigma$ is finitely generated over $k$, given any projective model $X_0$ of $\Sigma$ and given any valuation $B$ of dimension zero of $\Sigma$ containing $k$ with center $a_0$ on $X_0$, there exists a projective model $X$ of $\Sigma$ on which the center of $B$ is at a smooth point $a$ of $X$ such that the inclusion relation $\mathcal{O}_{X,a}\supset \mathcal{O}_{X_0,a_0}$ of local rings holds.

We conclude that Corollary~\ref{local uniformization} holds.

\bibliographystyle{amsplain}

\end{document}